\documentclass{article}
\pdfoutput=1	%
\usepackage{fullpage}
\usepackage[breaklinks,hidelinks]{hyperref}
\hypersetup{
  pdftitle={G\"odel's Incompleteness Theorems in Grounded Arithmetic},
  pdfauthor={Bryan Ford}}

\usepackage{amsmath}
\usepackage{amssymb}
\usepackage{amsthm}

\newtheorem{thm}{Theorem}[section]
\newtheorem{cor}[thm]{Corollary}	%

\usepackage{times}
\usepackage{xspace}
\usepackage{cite}		%
\usepackage{bussproofs}
\usepackage{stmaryrd}
\usepackage[noabbrev,capitalize]{cleveref}
\usepackage{thmtools}
\usepackage{xstring}
\usepackage{graphicx}
\usepackage{subcaption}
\usepackage{nicematrix}		%
\usepackage{longtable}
\usepackage{pifont}

\usepackage{isabelle}
\usepackage{isabellesym}
\usepackage{pdfpages}
\usepackage{ragged2e}
\newcommand{\SNIPTEXT}[1]{#1} %
\newcommand{\SNIP}[2]{\expandafter\newcommand\csname snippet--#1\endcsname{#2}}
\IfFileExists{snips.tex}{\SNIP{Goedel-Paper-Statements:theory:Goedel-Paper-Statements-0}{%
\isakeywordONE{theory}\isamarkupfalse%
\ Goedel{\isacharunderscore}{\kern0pt}Paper{\isacharunderscore}{\kern0pt}Statements\isanewline%
}%
\SNIP{Goedel-Paper-Statements:theory:Goedel-Paper-Statements-1}{%
\ \ \isakeywordTWO{imports}\isanewline%
}%
\SNIP{Goedel-Paper-Statements:theory:Goedel-Paper-Statements-2}{%
\ \ \ \ {\isachardoublequoteopen}RGA{\isacharunderscore}{\kern0pt}Decider{\isachardot}{\kern0pt}Decider{\isacharunderscore}{\kern0pt}Truth{\isachardoublequoteclose}\ {\isachardoublequoteopen}RGA{\isacharunderscore}{\kern0pt}Checker{\isachardot}{\kern0pt}Chk{\isacharunderscore}{\kern0pt}Faithful{\isachardoublequoteclose}\ {\isachardoublequoteopen}GCalc{\isacharunderscore}{\kern0pt}Core{\isachardot}{\kern0pt}GCalc{\isacharunderscore}{\kern0pt}Core{\isachardoublequoteclose}\ {\isachardoublequoteopen}G{\isadigit{2}}{\isacharunderscore}{\kern0pt}Locale{\isachardot}{\kern0pt}G{\isadigit{2}}{\isacharunderscore}{\kern0pt}Locale{\isachardoublequoteclose}\isanewline%
}%
\SNIP{Goedel-Paper-Statements:theory:Goedel-Paper-Statements-3}{%
\ \ \ \ {\isachardoublequoteopen}RGA{\isacharunderscore}{\kern0pt}Goedel{\isachardot}{\kern0pt}Goedel{\isacharunderscore}{\kern0pt}II{\isachardoublequoteclose}\ {\isachardoublequoteopen}RGA{\isacharunderscore}{\kern0pt}Goedel{\isachardot}{\kern0pt}Goedel{\isacharunderscore}{\kern0pt}II{\isacharunderscore}{\kern0pt}PR{\isachardoublequoteclose}\ {\isachardoublequoteopen}RGA{\isacharunderscore}{\kern0pt}Goedel{\isachardot}{\kern0pt}Goedel{\isacharunderscore}{\kern0pt}D{\isadigit{2}}{\isachardoublequoteclose}\ {\isachardoublequoteopen}RGA{\isacharunderscore}{\kern0pt}Goedel{\isachardot}{\kern0pt}Goedel{\isacharunderscore}{\kern0pt}D{\isadigit{3}}{\isachardoublequoteclose}\isanewline%
}%
\SNIP{Goedel-Paper-Statements:theory:Goedel-Paper-Statements-4}{%
\ \ \ \ {\isachardoublequoteopen}RGA{\isacharunderscore}{\kern0pt}Goedel{\isachardot}{\kern0pt}Goedel{\isacharunderscore}{\kern0pt}D{\isadigit{3}}{\isacharunderscore}{\kern0pt}Uniform{\isachardoublequoteclose}\ {\isachardoublequoteopen}RGA{\isacharunderscore}{\kern0pt}Goedel{\isachardot}{\kern0pt}Goedel{\isacharunderscore}{\kern0pt}Transformers{\isachardoublequoteclose}\ {\isachardoublequoteopen}RGA{\isacharunderscore}{\kern0pt}Goedel{\isachardot}{\kern0pt}Goedel{\isacharunderscore}{\kern0pt}Interp{\isachardoublequoteclose}\isanewline%
}%
\SNIP{Goedel-Paper-Statements:theory:Goedel-Paper-Statements-5}{%
\ \ \ \ {\isachardoublequoteopen}RGA{\isacharunderscore}{\kern0pt}Goedel{\isachardot}{\kern0pt}Entailment{\isacharunderscore}{\kern0pt}Goedel{\isachardoublequoteclose}\ {\isachardoublequoteopen}RGA{\isacharunderscore}{\kern0pt}Goedel{\isachardot}{\kern0pt}Goedel{\isacharunderscore}{\kern0pt}Cut{\isachardoublequoteclose}\ {\isachardoublequoteopen}RGA{\isacharunderscore}{\kern0pt}Goedel{\isachardot}{\kern0pt}Goedel{\isacharunderscore}{\kern0pt}Conditions{\isachardoublequoteclose}\isanewline%
}%
\SNIP{Goedel-Paper-Statements:theory:Goedel-Paper-Statements-6}{%
\ \ \ \ {\isachardoublequoteopen}RGA{\isacharunderscore}{\kern0pt}Liar{\isacharunderscore}{\kern0pt}Div{\isachardot}{\kern0pt}Entailment{\isacharunderscore}{\kern0pt}Liar{\isachardoublequoteclose}\ {\isachardoublequoteopen}RGA{\isacharunderscore}{\kern0pt}Checker{\isacharunderscore}{\kern0pt}Gen{\isachardot}{\kern0pt}Chk{\isacharunderscore}{\kern0pt}Gen{\isacharunderscore}{\kern0pt}Base{\isachardoublequoteclose}\ {\isachardoublequoteopen}RGA{\isacharunderscore}{\kern0pt}Checker{\isacharunderscore}{\kern0pt}Gen{\isachardot}{\kern0pt}Chk{\isacharunderscore}{\kern0pt}Gen{\isacharunderscore}{\kern0pt}Graft{\isachardoublequoteclose}\isanewline%
}%
\SNIP{Goedel-Paper-Statements:theory:Goedel-Paper-Statements-7}{%
\ \ \ \ {\isachardoublequoteopen}IGA{\isacharunderscore}{\kern0pt}Checker{\isachardot}{\kern0pt}Chk{\isacharunderscore}{\kern0pt}IGA{\isacharunderscore}{\kern0pt}Pos{\isachardoublequoteclose}\ {\isachardoublequoteopen}IGA{\isacharunderscore}{\kern0pt}Checker{\isachardot}{\kern0pt}Chk{\isacharunderscore}{\kern0pt}IGA{\isacharunderscore}{\kern0pt}Gfix{\isachardoublequoteclose}\ {\isachardoublequoteopen}IGA{\isacharunderscore}{\kern0pt}Checker{\isachardot}{\kern0pt}Chk{\isacharunderscore}{\kern0pt}IGA{\isacharunderscore}{\kern0pt}Table{\isachardoublequoteclose}\ {\isachardoublequoteopen}RGA{\isacharunderscore}{\kern0pt}Adopt{\isachardot}{\kern0pt}Adopt{\isacharunderscore}{\kern0pt}RGA{\isachardoublequoteclose}\isanewline%
}%
\SNIP{Goedel-Paper-Statements:theory:Goedel-Paper-Statements-8}{%
\ \ \ \ {\isachardoublequoteopen}IGA{\isacharunderscore}{\kern0pt}Adopt{\isachardot}{\kern0pt}Adopt{\isacharunderscore}{\kern0pt}IGA{\isachardoublequoteclose}\ {\isachardoublequoteopen}IGA{\isacharunderscore}{\kern0pt}Goedel{\isachardot}{\kern0pt}GCalc{\isacharunderscore}{\kern0pt}IGA{\isachardoublequoteclose}\ {\isachardoublequoteopen}IGA{\isacharunderscore}{\kern0pt}Goedel{\isachardot}{\kern0pt}Goedel{\isacharunderscore}{\kern0pt}Gen{\isacharunderscore}{\kern0pt}Transformers{\isachardoublequoteclose}\ {\isachardoublequoteopen}IGA{\isacharunderscore}{\kern0pt}Goedel{\isachardot}{\kern0pt}Goedel{\isacharunderscore}{\kern0pt}IGA{\isachardoublequoteclose}\isanewline%
}%
\SNIP{Goedel-Paper-Statements:theory:Goedel-Paper-Statements-9}{%
\ \ \ \ {\isachardoublequoteopen}IGA{\isacharunderscore}{\kern0pt}NoGo{\isachardot}{\kern0pt}IGA{\isacharunderscore}{\kern0pt}NoGo{\isachardoublequoteclose}\ {\isachardoublequoteopen}IGA{\isacharunderscore}{\kern0pt}Gentzen{\isachardot}{\kern0pt}Gentzen{\isacharunderscore}{\kern0pt}IGA{\isachardoublequoteclose}\ {\isachardoublequoteopen}IGA{\isacharunderscore}{\kern0pt}PGSem{\isachardot}{\kern0pt}GSem{\isacharunderscore}{\kern0pt}IGA{\isachardoublequoteclose}\ {\isachardoublequoteopen}IGA{\isacharunderscore}{\kern0pt}RE{\isachardot}{\kern0pt}RE{\isacharunderscore}{\kern0pt}IGA{\isachardoublequoteclose}\ {\isachardoublequoteopen}IGA{\isacharunderscore}{\kern0pt}Closing{\isachardot}{\kern0pt}CO{\isacharunderscore}{\kern0pt}PGBiv{\isachardoublequoteclose}\isanewline%
}%
\SNIP{Goedel-Paper-Statements:theory:Goedel-Paper-Statements-10}{%
\isakeywordTWO{begin}%
\endisatagtheory
{\isafoldtheory}%
\isadelimtheory
\endisadelimtheory
\SNIPTEXT{\begin{isamarkuptext}%
A document-only theory: it proves nothing and defines nothing.  For each Isabelle name in a
  theorem or definition header of the Goedel paper (\isa{doc{\isacharslash}{\kern0pt}goedel}) it has one \isa{thm} command
  naming the fact, followed by the fact's statement as Isabelle prints it (a \isa{thm}
  antiquotation in a \isa{text{\isacharunderscore}{\kern0pt}raw} block): a theorem by its name, a defined constant by its
  defining equation \isa{NAME{\isacharunderscore}{\kern0pt}def}, a locale of the paper's definitions by the assumptions it adds
  (\isa{NAME{\isacharunderscore}{\kern0pt}axioms{\isacharunderscore}{\kern0pt}def}; the base locale \isa{gcalc}, whose parent has none, by
  \isa{gcalc{\isacharunderscore}{\kern0pt}def}).  The theories that state the results are the imports above.  The theory
  follows the conventions of \isa{ConRGA{\isacharunderscore}{\kern0pt}Paper{\isacharunderscore}{\kern0pt}Statements}, the statements theory of the
  internal-reasoning paper.

  The session's document output is cut into named snippets by the cartouches tool
  (\url{https://github.com/astahfrom/cartouches}) into \isa{doc{\isacharslash}{\kern0pt}goedel{\isacharslash}{\kern0pt}snips{\isachardot}{\kern0pt}tex}.
  The snippet of the fact \isa{NAME} is \isa{Goedel{\isacharminus}{\kern0pt}Paper{\isacharminus}{\kern0pt}Statements{\isacharcolon}{\kern0pt}thm{\isacharcolon}{\kern0pt}NAME}, with every
  underscore of \isa{NAME} written as a hyphen; its line 0 is the command
  \isa{thm\ NAME} and its line 1 is the printed statement.  So, in the paper,
  \isa{{\isacharbackslash}{\kern0pt}Snippet{\isacharbraceleft}{\kern0pt}Goedel{\isacharminus}{\kern0pt}Paper{\isacharminus}{\kern0pt}Statements{\isacharcolon}{\kern0pt}thm{\isacharcolon}{\kern0pt}NAME{\isacharbraceright}{\kern0pt}} shows both lines and
  \isa{{\isacharbackslash}{\kern0pt}GetSnip{\isacharbraceleft}{\kern0pt}Goedel{\isacharminus}{\kern0pt}Paper{\isacharminus}{\kern0pt}Statements{\isacharcolon}{\kern0pt}thm{\isacharcolon}{\kern0pt}NAME{\isacharminus}{\kern0pt}{\isadigit{1}}{\isacharbraceright}{\kern0pt}} the statement alone.  A name
  shown twice gets the suffix \isa{{\isacharplus}{\kern0pt}{\isadigit{1}}} on its second showing: the generic
  \isa{sigma{\isacharunderscore}{\kern0pt}deduction{\isacharunderscore}{\kern0pt}theorem} of the locale \isa{gcalc{\isacharunderscore}{\kern0pt}strip} is
  \isa{thm{\isacharcolon}{\kern0pt}sigma{\isacharminus}{\kern0pt}deduction{\isacharminus}{\kern0pt}theorem{\isacharplus}{\kern0pt}{\isadigit{1}}}, the plain name being RGA's own.  A qualified name
  \isa{X{\isachardot}{\kern0pt}simps} is cut under \isa{X}.  The facts stand inside the context of the locale
  that states them, one context block per locale.

  Rules of this theory, so that the snippets stay well formed: every statement is printed inline,
  no prose stands between two facts, and prose names commands with \isa{text} antiquotations
  only; a statement too short to contain a symbol carries a LaTeX comment naming one, since the
  tool drops a line without.  Schematic variables are printed without question marks (the session option
  \isa{show{\isacharunderscore}{\kern0pt}question{\isacharunderscore}{\kern0pt}marks}).  As in the companion theory, the biconditional of the statements is
  given back its infix syntax inside \isa{rga}.%
\end{isamarkuptext}}\isamarkuptrue%
}%
\SNIP{Goedel-Paper-Statements:thm:gcalc-def-0}{%
\isakeywordONE{thm}\isamarkupfalse%
\ gcalc{\isacharunderscore}{\kern0pt}def%
\isanewline%
}%
\SNIP{Goedel-Paper-Statements:thm:gcalc-def-1}{%
\isa{\isaconst{gcalc}\ \isavar{ent}\ {\isasymequiv}\ {\isacharparenleft}{\kern0pt}{\isacharparenleft}{\kern0pt}{\isacharparenleft}{\kern0pt}{\isacharparenleft}{\kern0pt}{\isacharparenleft}{\kern0pt}{\isasymforall}\isabound{{\isasymGamma}}\ \isabound{c}{\isachardot}{\kern0pt}\ \isavar{ent}\ \isabound{{\isasymGamma}}\ \isabound{c}\ {\isasymlongrightarrow}\ \isaconst{gigc{\isachardot}{\kern0pt}wft}\ \isabound{c}{\isacharparenright}{\kern0pt}\ {\isasymand}\ {\isacharparenleft}{\kern0pt}{\isasymforall}\isabound{{\isasymGamma}}\ \isabound{c}\ \isabound{s}\ \isabound{v}{\isachardot}{\kern0pt}\ \isavar{ent}\ \isabound{{\isasymGamma}}\ \isabound{c}\ {\isasymlongrightarrow}\ \isaconst{gigc{\isachardot}{\kern0pt}wfc}\ \isabound{{\isasymGamma}}\ {\isasymlongrightarrow}\ \isaconst{gigc{\isachardot}{\kern0pt}wft}\ \isabound{s}\ {\isasymlongrightarrow}\ \isavar{ent}\ {\isacharparenleft}{\kern0pt}{\isacharbrackleft}{\kern0pt}\isabound{v}{\isasymZpinj}\isabound{s}{\isacharbrackright}{\kern0pt}\isabound{{\isasymGamma}}{\isacharparenright}{\kern0pt}\ {\isacharparenleft}{\kern0pt}{\isacharbrackleft}{\kern0pt}\isabound{v}{\isasymZpinj}\isabound{s}{\isacharbrackright}{\kern0pt}\isabound{c}{\isacharparenright}{\kern0pt}{\isacharparenright}{\kern0pt}{\isacharparenright}{\kern0pt}\ {\isasymand}\ {\isacharparenleft}{\kern0pt}{\isasymforall}\isabound{c}\ \isabound{{\isasymGamma}}{\isachardot}{\kern0pt}\ \isaconst{gigc{\isachardot}{\kern0pt}wft}\ \isabound{c}\ {\isasymlongrightarrow}\ \isavar{ent}\ {\isacharparenleft}{\kern0pt}\isabound{c}\ {\isasymtriangleright}\ \isabound{{\isasymGamma}}{\isacharparenright}{\kern0pt}\ \isabound{c}{\isacharparenright}{\kern0pt}\ {\isasymand}\ {\isacharparenleft}{\kern0pt}{\isasymforall}\isabound{{\isasymGamma}}\ \isabound{c}\ \isabound{a}{\isachardot}{\kern0pt}\ \isavar{ent}\ \isabound{{\isasymGamma}}\ \isabound{c}\ {\isasymlongrightarrow}\ \isavar{ent}\ {\isacharparenleft}{\kern0pt}\isabound{a}\ {\isasymtriangleright}\ \isabound{{\isasymGamma}}{\isacharparenright}{\kern0pt}\ \isabound{c}{\isacharparenright}{\kern0pt}\ {\isasymand}\ {\isacharparenleft}{\kern0pt}{\isasymforall}\isabound{{\isasymGamma}}\ \isabound{a}\ \isabound{c}{\isachardot}{\kern0pt}\ \isavar{ent}\ \isabound{{\isasymGamma}}\ \isabound{a}\ {\isasymlongrightarrow}\ \isavar{ent}\ {\isacharparenleft}{\kern0pt}\isabound{a}\ {\isasymtriangleright}\ \isabound{{\isasymGamma}}{\isacharparenright}{\kern0pt}\ \isabound{c}\ {\isasymlongrightarrow}\ \isavar{ent}\ \isabound{{\isasymGamma}}\ \isabound{c}{\isacharparenright}{\kern0pt}{\isacharparenright}{\kern0pt}\ {\isasymand}\ {\isacharparenleft}{\kern0pt}{\isacharparenleft}{\kern0pt}{\isasymforall}\isabound{{\isasymGamma}}\ \isabound{a}{\isachardot}{\kern0pt}\ \isavar{ent}\ \isabound{{\isasymGamma}}\ \isabound{a}\ {\isasymlongrightarrow}\ \isavar{ent}\ \isabound{{\isasymGamma}}\ {\isacharparenleft}{\kern0pt}\isaconst{gigc{\isachardot}{\kern0pt}trmEq}\ \isabound{a}\ {\isacharparenleft}{\kern0pt}\isaconst{gigc{\isachardot}{\kern0pt}trm{\isacharunderscore}{\kern0pt}num}\ {\isadigit{1}}{\isacharparenright}{\kern0pt}{\isacharparenright}{\kern0pt}{\isacharparenright}{\kern0pt}\ {\isasymand}\ {\isacharparenleft}{\kern0pt}{\isasymforall}\isabound{{\isasymGamma}}\ \isabound{a}{\isachardot}{\kern0pt}\ \isavar{ent}\ \isabound{{\isasymGamma}}\ {\isacharparenleft}{\kern0pt}\isaconst{gigc{\isachardot}{\kern0pt}trmEq}\ \isabound{a}\ {\isacharparenleft}{\kern0pt}\isaconst{gigc{\isachardot}{\kern0pt}trm{\isacharunderscore}{\kern0pt}num}\ {\isadigit{1}}{\isacharparenright}{\kern0pt}{\isacharparenright}{\kern0pt}\ {\isasymlongrightarrow}\ \isavar{ent}\ \isabound{{\isasymGamma}}\ \isabound{a}{\isacharparenright}{\kern0pt}{\isacharparenright}{\kern0pt}\ {\isasymand}\ {\isacharparenleft}{\kern0pt}{\isasymforall}\isabound{{\isasymGamma}}\ \isabound{a}{\isachardot}{\kern0pt}\ \isavar{ent}\ \isabound{{\isasymGamma}}\ {\isacharparenleft}{\kern0pt}\isaconst{gigc{\isachardot}{\kern0pt}trmNeg}\ \isabound{a}{\isacharparenright}{\kern0pt}\ {\isasymlongrightarrow}\ \isavar{ent}\ \isabound{{\isasymGamma}}\ {\isacharparenleft}{\kern0pt}\isaconst{gigc{\isachardot}{\kern0pt}trmEq}\ \isabound{a}\ {\isacharparenleft}{\kern0pt}\isaconst{gigc{\isachardot}{\kern0pt}trm{\isacharunderscore}{\kern0pt}num}\ {\isadigit{0}}{\isacharparenright}{\kern0pt}{\isacharparenright}{\kern0pt}{\isacharparenright}{\kern0pt}\ {\isasymand}\ {\isacharparenleft}{\kern0pt}{\isasymforall}\isabound{{\isasymGamma}}\ \isabound{a}{\isachardot}{\kern0pt}\ \isavar{ent}\ \isabound{{\isasymGamma}}\ {\isacharparenleft}{\kern0pt}\isaconst{gigc{\isachardot}{\kern0pt}trmEq}\ \isabound{a}\ {\isacharparenleft}{\kern0pt}\isaconst{gigc{\isachardot}{\kern0pt}trm{\isacharunderscore}{\kern0pt}num}\ {\isadigit{0}}{\isacharparenright}{\kern0pt}{\isacharparenright}{\kern0pt}\ {\isasymlongrightarrow}\ \isavar{ent}\ \isabound{{\isasymGamma}}\ {\isacharparenleft}{\kern0pt}\isaconst{gigc{\isachardot}{\kern0pt}trmNeg}\ \isabound{a}{\isacharparenright}{\kern0pt}{\isacharparenright}{\kern0pt}\ {\isasymand}\ {\isacharparenleft}{\kern0pt}{\isasymforall}\isabound{{\isasymGamma}}\ \isabound{a}\ \isabound{c}{\isachardot}{\kern0pt}\ \isavar{ent}\ \isabound{{\isasymGamma}}\ {\isacharparenleft}{\kern0pt}\isaconst{gigc{\isachardot}{\kern0pt}trmNeg}\ \isabound{a}{\isacharparenright}{\kern0pt}\ {\isasymlongrightarrow}\ \isavar{ent}\ \isabound{{\isasymGamma}}\ \isabound{a}\ {\isasymlongrightarrow}\ \isaconst{gigc{\isachardot}{\kern0pt}wft}\ \isabound{c}\ {\isasymlongrightarrow}\ \isavar{ent}\ \isabound{{\isasymGamma}}\ \isabound{c}{\isacharparenright}{\kern0pt}{\isacharparenright}{\kern0pt}\ {\isasymand}\ {\isacharparenleft}{\kern0pt}{\isacharparenleft}{\kern0pt}{\isacharparenleft}{\kern0pt}{\isasymforall}\isabound{{\isasymGamma}}\ \isabound{a}{\isachardot}{\kern0pt}\ \isavar{ent}\ \isabound{{\isasymGamma}}\ \isabound{a}\ {\isasymlongrightarrow}\ \isavar{ent}\ \isabound{{\isasymGamma}}\ {\isacharparenleft}{\kern0pt}\isaconst{gigc{\isachardot}{\kern0pt}trmNeg}\ {\isacharparenleft}{\kern0pt}\isaconst{gigc{\isachardot}{\kern0pt}trmNeg}\ \isabound{a}{\isacharparenright}{\kern0pt}{\isacharparenright}{\kern0pt}{\isacharparenright}{\kern0pt}\ {\isasymand}\ {\isacharparenleft}{\kern0pt}{\isasymforall}\isabound{{\isasymGamma}}\ \isabound{a}{\isachardot}{\kern0pt}\ \isavar{ent}\ \isabound{{\isasymGamma}}\ {\isacharparenleft}{\kern0pt}\isaconst{gigc{\isachardot}{\kern0pt}trmNeg}\ {\isacharparenleft}{\kern0pt}\isaconst{gigc{\isachardot}{\kern0pt}trmNeg}\ \isabound{a}{\isacharparenright}{\kern0pt}{\isacharparenright}{\kern0pt}\ {\isasymlongrightarrow}\ \isavar{ent}\ \isabound{{\isasymGamma}}\ \isabound{a}{\isacharparenright}{\kern0pt}{\isacharparenright}{\kern0pt}\ {\isasymand}\ {\isacharparenleft}{\kern0pt}{\isasymforall}\isabound{{\isasymGamma}}\ \isabound{a}\ \isabound{b}{\isachardot}{\kern0pt}\ \isavar{ent}\ \isabound{{\isasymGamma}}\ \isabound{a}\ {\isasymlongrightarrow}\ \isaconst{gigc{\isachardot}{\kern0pt}wft}\ \isabound{b}\ {\isasymlongrightarrow}\ \isavar{ent}\ \isabound{{\isasymGamma}}\ {\isacharparenleft}{\kern0pt}\isaconst{gigc{\isachardot}{\kern0pt}trmOr}\ \isabound{a}\ \isabound{b}{\isacharparenright}{\kern0pt}{\isacharparenright}{\kern0pt}\ {\isasymand}\ {\isacharparenleft}{\kern0pt}{\isasymforall}\isabound{{\isasymGamma}}\ \isabound{b}\ \isabound{a}{\isachardot}{\kern0pt}\ \isavar{ent}\ \isabound{{\isasymGamma}}\ \isabound{b}\ {\isasymlongrightarrow}\ \isaconst{gigc{\isachardot}{\kern0pt}wft}\ \isabound{a}\ {\isasymlongrightarrow}\ \isavar{ent}\ \isabound{{\isasymGamma}}\ {\isacharparenleft}{\kern0pt}\isaconst{gigc{\isachardot}{\kern0pt}trmOr}\ \isabound{a}\ \isabound{b}{\isacharparenright}{\kern0pt}{\isacharparenright}{\kern0pt}\ {\isasymand}\ {\isacharparenleft}{\kern0pt}{\isasymforall}\isabound{{\isasymGamma}}\ \isabound{a}\ \isabound{b}{\isachardot}{\kern0pt}\ \isavar{ent}\ \isabound{{\isasymGamma}}\ {\isacharparenleft}{\kern0pt}\isaconst{gigc{\isachardot}{\kern0pt}trmNeg}\ \isabound{a}{\isacharparenright}{\kern0pt}\ {\isasymlongrightarrow}\ \isavar{ent}\ \isabound{{\isasymGamma}}\ {\isacharparenleft}{\kern0pt}\isaconst{gigc{\isachardot}{\kern0pt}trmNeg}\ \isabound{b}{\isacharparenright}{\kern0pt}\ {\isasymlongrightarrow}\ \isavar{ent}\ \isabound{{\isasymGamma}}\ {\isacharparenleft}{\kern0pt}\isaconst{gigc{\isachardot}{\kern0pt}trmNeg}\ {\isacharparenleft}{\kern0pt}\isaconst{gigc{\isachardot}{\kern0pt}trmOr}\ \isabound{a}\ \isabound{b}{\isacharparenright}{\kern0pt}{\isacharparenright}{\kern0pt}{\isacharparenright}{\kern0pt}{\isacharparenright}{\kern0pt}\ {\isasymand}\ {\isacharparenleft}{\kern0pt}{\isacharparenleft}{\kern0pt}{\isasymforall}\isabound{{\isasymGamma}}\ \isabound{a}\ \isabound{b}\ \isabound{c}{\isachardot}{\kern0pt}\ \isavar{ent}\ \isabound{{\isasymGamma}}\ {\isacharparenleft}{\kern0pt}\isaconst{gigc{\isachardot}{\kern0pt}trmOr}\ \isabound{a}\ \isabound{b}{\isacharparenright}{\kern0pt}\ {\isasymlongrightarrow}\ \isavar{ent}\ {\isacharparenleft}{\kern0pt}\isabound{a}\ {\isasymtriangleright}\ \isabound{{\isasymGamma}}{\isacharparenright}{\kern0pt}\ \isabound{c}\ {\isasymlongrightarrow}\ \isavar{ent}\ {\isacharparenleft}{\kern0pt}\isabound{b}\ {\isasymtriangleright}\ \isabound{{\isasymGamma}}{\isacharparenright}{\kern0pt}\ \isabound{c}\ {\isasymlongrightarrow}\ \isavar{ent}\ \isabound{{\isasymGamma}}\ \isabound{c}{\isacharparenright}{\kern0pt}\ {\isasymand}\ {\isacharparenleft}{\kern0pt}{\isasymforall}\isabound{{\isasymGamma}}\ \isabound{a}\ \isabound{b}{\isachardot}{\kern0pt}\ \isavar{ent}\ \isabound{{\isasymGamma}}\ {\isacharparenleft}{\kern0pt}\isaconst{gigc{\isachardot}{\kern0pt}trmNeg}\ {\isacharparenleft}{\kern0pt}\isaconst{gigc{\isachardot}{\kern0pt}trmOr}\ \isabound{a}\ \isabound{b}{\isacharparenright}{\kern0pt}{\isacharparenright}{\kern0pt}\ {\isasymlongrightarrow}\ \isavar{ent}\ \isabound{{\isasymGamma}}\ {\isacharparenleft}{\kern0pt}\isaconst{gigc{\isachardot}{\kern0pt}trmNeg}\ \isabound{a}{\isacharparenright}{\kern0pt}{\isacharparenright}{\kern0pt}\ {\isasymand}\ {\isacharparenleft}{\kern0pt}{\isasymforall}\isabound{{\isasymGamma}}\ \isabound{a}\ \isabound{b}{\isachardot}{\kern0pt}\ \isavar{ent}\ \isabound{{\isasymGamma}}\ {\isacharparenleft}{\kern0pt}\isaconst{gigc{\isachardot}{\kern0pt}trmNeg}\ {\isacharparenleft}{\kern0pt}\isaconst{gigc{\isachardot}{\kern0pt}trmOr}\ \isabound{a}\ \isabound{b}{\isacharparenright}{\kern0pt}{\isacharparenright}{\kern0pt}\ {\isasymlongrightarrow}\ \isavar{ent}\ \isabound{{\isasymGamma}}\ {\isacharparenleft}{\kern0pt}\isaconst{gigc{\isachardot}{\kern0pt}trmNeg}\ \isabound{b}{\isacharparenright}{\kern0pt}{\isacharparenright}{\kern0pt}{\isacharparenright}{\kern0pt}\ {\isasymand}\ {\isacharparenleft}{\kern0pt}{\isasymforall}\isabound{{\isasymGamma}}\ \isabound{a}{\isachardot}{\kern0pt}\ \isavar{ent}\ \isabound{{\isasymGamma}}\ {\isacharparenleft}{\kern0pt}\isaconst{gigc{\isachardot}{\kern0pt}trmNat}\ \isabound{a}{\isacharparenright}{\kern0pt}\ {\isasymlongrightarrow}\ \isavar{ent}\ \isabound{{\isasymGamma}}\ {\isacharparenleft}{\kern0pt}\isaconst{gigc{\isachardot}{\kern0pt}trmEq}\ \isabound{a}\ \isabound{a}{\isacharparenright}{\kern0pt}{\isacharparenright}{\kern0pt}\ {\isasymand}\ {\isacharparenleft}{\kern0pt}{\isasymforall}\isabound{{\isasymGamma}}\ \isabound{a}{\isachardot}{\kern0pt}\ \isavar{ent}\ \isabound{{\isasymGamma}}\ {\isacharparenleft}{\kern0pt}\isaconst{gigc{\isachardot}{\kern0pt}trmEq}\ \isabound{a}\ \isabound{a}{\isacharparenright}{\kern0pt}\ {\isasymlongrightarrow}\ \isavar{ent}\ \isabound{{\isasymGamma}}\ {\isacharparenleft}{\kern0pt}\isaconst{gigc{\isachardot}{\kern0pt}trmNat}\ \isabound{a}{\isacharparenright}{\kern0pt}{\isacharparenright}{\kern0pt}\ {\isasymand}\ {\isacharparenleft}{\kern0pt}{\isasymforall}\isabound{{\isasymGamma}}\ \isabound{a}\ \isabound{b}{\isachardot}{\kern0pt}\ \isavar{ent}\ \isabound{{\isasymGamma}}\ {\isacharparenleft}{\kern0pt}\isaconst{gigc{\isachardot}{\kern0pt}trmEq}\ \isabound{a}\ \isabound{b}{\isacharparenright}{\kern0pt}\ {\isasymlongrightarrow}\ \isavar{ent}\ \isabound{{\isasymGamma}}\ {\isacharparenleft}{\kern0pt}\isaconst{gigc{\isachardot}{\kern0pt}trmEq}\ \isabound{b}\ \isabound{a}{\isacharparenright}{\kern0pt}{\isacharparenright}{\kern0pt}{\isacharparenright}{\kern0pt}\ {\isasymand}\ {\isacharparenleft}{\kern0pt}{\isacharparenleft}{\kern0pt}{\isacharparenleft}{\kern0pt}{\isacharparenleft}{\kern0pt}{\isasymforall}\isabound{{\isasymGamma}}\ \isabound{a}\ \isabound{b}\ \isabound{i}\ \isabound{p}{\isachardot}{\kern0pt}\ \isavar{ent}\ \isabound{{\isasymGamma}}\ {\isacharparenleft}{\kern0pt}\isaconst{gigc{\isachardot}{\kern0pt}trmEq}\ \isabound{a}\ \isabound{b}{\isacharparenright}{\kern0pt}\ {\isasymlongrightarrow}\ \isavar{ent}\ \isabound{{\isasymGamma}}\ {\isacharparenleft}{\kern0pt}{\isacharbrackleft}{\kern0pt}\isabound{i}{\isasymZpinj}\isabound{a}{\isacharbrackright}{\kern0pt}\isabound{p}{\isacharparenright}{\kern0pt}\ {\isasymlongrightarrow}\ \isaconst{gigc{\isachardot}{\kern0pt}wft}\ \isabound{p}\ {\isasymlongrightarrow}\ \isavar{ent}\ \isabound{{\isasymGamma}}\ {\isacharparenleft}{\kern0pt}{\isacharbrackleft}{\kern0pt}\isabound{i}{\isasymZpinj}\isabound{b}{\isacharbrackright}{\kern0pt}\isabound{p}{\isacharparenright}{\kern0pt}{\isacharparenright}{\kern0pt}\ {\isasymand}\ {\isacharparenleft}{\kern0pt}{\isasymforall}\isabound{{\isasymGamma}}\ \isabound{a}\ \isabound{b}{\isachardot}{\kern0pt}\ \isavar{ent}\ \isabound{{\isasymGamma}}\ {\isacharparenleft}{\kern0pt}\isaconst{gigc{\isachardot}{\kern0pt}trmNe}\ \isabound{a}\ \isabound{b}{\isacharparenright}{\kern0pt}\ {\isasymlongrightarrow}\ \isavar{ent}\ \isabound{{\isasymGamma}}\ {\isacharparenleft}{\kern0pt}\isaconst{gigc{\isachardot}{\kern0pt}trmNe}\ \isabound{b}\ \isabound{a}{\isacharparenright}{\kern0pt}{\isacharparenright}{\kern0pt}{\isacharparenright}{\kern0pt}\ {\isasymand}\ {\isacharparenleft}{\kern0pt}{\isasymforall}\isabound{{\isasymGamma}}\ \isabound{a}\ \isabound{b}{\isachardot}{\kern0pt}\ \isavar{ent}\ \isabound{{\isasymGamma}}\ {\isacharparenleft}{\kern0pt}\isaconst{gigc{\isachardot}{\kern0pt}trmNat}\ \isabound{a}{\isacharparenright}{\kern0pt}\ {\isasymlongrightarrow}\ \isavar{ent}\ \isabound{{\isasymGamma}}\ {\isacharparenleft}{\kern0pt}\isaconst{gigc{\isachardot}{\kern0pt}trmNat}\ \isabound{b}{\isacharparenright}{\kern0pt}\ {\isasymlongrightarrow}\ \isavar{ent}\ \isabound{{\isasymGamma}}\ {\isacharparenleft}{\kern0pt}\isaconst{gigc{\isachardot}{\kern0pt}trmTypeBool}\ {\isacharparenleft}{\kern0pt}\isaconst{gigc{\isachardot}{\kern0pt}trmEq}\ \isabound{a}\ \isabound{b}{\isacharparenright}{\kern0pt}{\isacharparenright}{\kern0pt}{\isacharparenright}{\kern0pt}\ {\isasymand}\ {\isacharparenleft}{\kern0pt}{\isasymforall}\isabound{{\isasymGamma}}\ \isabound{a}\ \isabound{b}{\isachardot}{\kern0pt}\ \isavar{ent}\ \isabound{{\isasymGamma}}\ {\isacharparenleft}{\kern0pt}\isaconst{gigc{\isachardot}{\kern0pt}trmTypeBool}\ {\isacharparenleft}{\kern0pt}\isaconst{gigc{\isachardot}{\kern0pt}trmEq}\ \isabound{a}\ \isabound{b}{\isacharparenright}{\kern0pt}{\isacharparenright}{\kern0pt}\ {\isasymlongrightarrow}\ \isavar{ent}\ \isabound{{\isasymGamma}}\ {\isacharparenleft}{\kern0pt}\isaconst{gigc{\isachardot}{\kern0pt}trmAnd}\ {\isacharparenleft}{\kern0pt}\isaconst{gigc{\isachardot}{\kern0pt}trmNat}\ \isabound{a}{\isacharparenright}{\kern0pt}\ {\isacharparenleft}{\kern0pt}\isaconst{gigc{\isachardot}{\kern0pt}trmNat}\ \isabound{b}{\isacharparenright}{\kern0pt}{\isacharparenright}{\kern0pt}{\isacharparenright}{\kern0pt}\ {\isasymand}\ {\isacharparenleft}{\kern0pt}{\isasymforall}\isabound{{\isasymGamma}}{\isachardot}{\kern0pt}\ \isavar{ent}\ \isabound{{\isasymGamma}}\ {\isacharparenleft}{\kern0pt}\isaconst{gigc{\isachardot}{\kern0pt}trmNat}\ \isaconst{gigc{\isachardot}{\kern0pt}trm{\isadigit{0}}}{\isacharparenright}{\kern0pt}{\isacharparenright}{\kern0pt}{\isacharparenright}{\kern0pt}\ {\isasymand}\ {\isacharparenleft}{\kern0pt}{\isacharparenleft}{\kern0pt}{\isasymforall}\isabound{{\isasymGamma}}\ \isabound{a}{\isachardot}{\kern0pt}\ \isavar{ent}\ \isabound{{\isasymGamma}}\ {\isacharparenleft}{\kern0pt}\isaconst{gigc{\isachardot}{\kern0pt}trmNat}\ \isabound{a}{\isacharparenright}{\kern0pt}\ {\isasymlongrightarrow}\ \isavar{ent}\ \isabound{{\isasymGamma}}\ {\isacharparenleft}{\kern0pt}\isaconst{gigc{\isachardot}{\kern0pt}trmNe}\ {\isacharparenleft}{\kern0pt}\isaconst{gigc{\isachardot}{\kern0pt}trmSuc}\ \isabound{a}{\isacharparenright}{\kern0pt}\ \isaconst{gigc{\isachardot}{\kern0pt}trm{\isadigit{0}}}{\isacharparenright}{\kern0pt}{\isacharparenright}{\kern0pt}\ {\isasymand}\ {\isacharparenleft}{\kern0pt}{\isasymforall}\isabound{{\isasymGamma}}\ \isabound{a}\ \isabound{b}{\isachardot}{\kern0pt}\ \isavar{ent}\ \isabound{{\isasymGamma}}\ {\isacharparenleft}{\kern0pt}\isaconst{gigc{\isachardot}{\kern0pt}trmNe}\ \isabound{a}\ \isabound{b}{\isacharparenright}{\kern0pt}\ {\isasymlongrightarrow}\ \isavar{ent}\ \isabound{{\isasymGamma}}\ {\isacharparenleft}{\kern0pt}\isaconst{gigc{\isachardot}{\kern0pt}trmNe}\ {\isacharparenleft}{\kern0pt}\isaconst{gigc{\isachardot}{\kern0pt}trmSuc}\ \isabound{a}{\isacharparenright}{\kern0pt}\ {\isacharparenleft}{\kern0pt}\isaconst{gigc{\isachardot}{\kern0pt}trmSuc}\ \isabound{b}{\isacharparenright}{\kern0pt}{\isacharparenright}{\kern0pt}{\isacharparenright}{\kern0pt}\ {\isasymand}\ {\isacharparenleft}{\kern0pt}{\isasymforall}\isabound{{\isasymGamma}}\ \isabound{a}{\isachardot}{\kern0pt}\ \isavar{ent}\ \isabound{{\isasymGamma}}\ {\isacharparenleft}{\kern0pt}\isaconst{gigc{\isachardot}{\kern0pt}trmNat}\ \isabound{a}{\isacharparenright}{\kern0pt}\ {\isasymlongrightarrow}\ \isavar{ent}\ \isabound{{\isasymGamma}}\ {\isacharparenleft}{\kern0pt}\isaconst{gigc{\isachardot}{\kern0pt}trmNat}\ {\isacharparenleft}{\kern0pt}\isaconst{gigc{\isachardot}{\kern0pt}trmSuc}\ \isabound{a}{\isacharparenright}{\kern0pt}{\isacharparenright}{\kern0pt}{\isacharparenright}{\kern0pt}{\isacharparenright}{\kern0pt}\ {\isasymand}\ {\isacharparenleft}{\kern0pt}{\isasymforall}\isabound{{\isasymGamma}}\ \isabound{a}\ \isabound{b}{\isachardot}{\kern0pt}\ \isavar{ent}\ \isabound{{\isasymGamma}}\ {\isacharparenleft}{\kern0pt}\isaconst{gigc{\isachardot}{\kern0pt}trmEq}\ {\isacharparenleft}{\kern0pt}\isaconst{gigc{\isachardot}{\kern0pt}trmSuc}\ \isabound{a}{\isacharparenright}{\kern0pt}\ {\isacharparenleft}{\kern0pt}\isaconst{gigc{\isachardot}{\kern0pt}trmSuc}\ \isabound{b}{\isacharparenright}{\kern0pt}{\isacharparenright}{\kern0pt}\ {\isasymlongrightarrow}\ \isavar{ent}\ \isabound{{\isasymGamma}}\ {\isacharparenleft}{\kern0pt}\isaconst{gigc{\isachardot}{\kern0pt}trmEq}\ \isabound{a}\ \isabound{b}{\isacharparenright}{\kern0pt}{\isacharparenright}{\kern0pt}\ {\isasymand}\ {\isacharparenleft}{\kern0pt}{\isasymforall}\isabound{{\isasymGamma}}\ \isabound{a}\ \isabound{b}{\isachardot}{\kern0pt}\ \isavar{ent}\ \isabound{{\isasymGamma}}\ {\isacharparenleft}{\kern0pt}\isaconst{gigc{\isachardot}{\kern0pt}trmNe}\ {\isacharparenleft}{\kern0pt}\isaconst{gigc{\isachardot}{\kern0pt}trmSuc}\ \isabound{a}{\isacharparenright}{\kern0pt}\ {\isacharparenleft}{\kern0pt}\isaconst{gigc{\isachardot}{\kern0pt}trmSuc}\ \isabound{b}{\isacharparenright}{\kern0pt}{\isacharparenright}{\kern0pt}\ {\isasymlongrightarrow}\ \isavar{ent}\ \isabound{{\isasymGamma}}\ {\isacharparenleft}{\kern0pt}\isaconst{gigc{\isachardot}{\kern0pt}trmNe}\ \isabound{a}\ \isabound{b}{\isacharparenright}{\kern0pt}{\isacharparenright}{\kern0pt}\ {\isasymand}\ {\isacharparenleft}{\kern0pt}{\isasymforall}\isabound{{\isasymGamma}}\ \isabound{a}{\isachardot}{\kern0pt}\ \isavar{ent}\ \isabound{{\isasymGamma}}\ {\isacharparenleft}{\kern0pt}\isaconst{gigc{\isachardot}{\kern0pt}trmNat}\ {\isacharparenleft}{\kern0pt}\isaconst{gigc{\isachardot}{\kern0pt}trmSuc}\ \isabound{a}{\isacharparenright}{\kern0pt}{\isacharparenright}{\kern0pt}\ {\isasymlongrightarrow}\ \isavar{ent}\ \isabound{{\isasymGamma}}\ {\isacharparenleft}{\kern0pt}\isaconst{gigc{\isachardot}{\kern0pt}trmNat}\ \isabound{a}{\isacharparenright}{\kern0pt}{\isacharparenright}{\kern0pt}{\isacharparenright}{\kern0pt}\ {\isasymand}\ {\isacharparenleft}{\kern0pt}{\isacharparenleft}{\kern0pt}{\isacharparenleft}{\kern0pt}{\isasymforall}\isabound{{\isasymGamma}}\ \isabound{i}\ \isabound{c}\ \isabound{a}{\isachardot}{\kern0pt}\ \isavar{ent}\ \isabound{{\isasymGamma}}\ {\isacharparenleft}{\kern0pt}{\isacharbrackleft}{\kern0pt}\isabound{i}{\isasymZpinj}\isaconst{gigc{\isachardot}{\kern0pt}trm{\isadigit{0}}}{\isacharbrackright}{\kern0pt}\isabound{c}{\isacharparenright}{\kern0pt}\ {\isasymlongrightarrow}\ \isavar{ent}\ {\isacharparenleft}{\kern0pt}\isaconst{gigc{\isachardot}{\kern0pt}trmNat}\ {\isacharparenleft}{\kern0pt}\isaconst{trmVar}\ \isabound{i}{\isacharparenright}{\kern0pt}\ {\isasymtriangleright}\ \isabound{c}\ {\isasymtriangleright}\ \isabound{{\isasymGamma}}{\isacharparenright}{\kern0pt}\ {\isacharparenleft}{\kern0pt}{\isacharbrackleft}{\kern0pt}\isabound{i}{\isasymZpinj}\isaconst{gigc{\isachardot}{\kern0pt}trmSuc}\ {\isacharparenleft}{\kern0pt}\isaconst{trmVar}\ \isabound{i}{\isacharparenright}{\kern0pt}{\isacharbrackright}{\kern0pt}\isabound{c}{\isacharparenright}{\kern0pt}\ {\isasymlongrightarrow}\ \isavar{ent}\ \isabound{{\isasymGamma}}\ {\isacharparenleft}{\kern0pt}\isaconst{gigc{\isachardot}{\kern0pt}trmNat}\ \isabound{a}{\isacharparenright}{\kern0pt}\ {\isasymlongrightarrow}\ \isaconst{gigc{\isachardot}{\kern0pt}wft}\ \isabound{c}\ {\isasymlongrightarrow}\ {\isasymnot}\ \isaconst{free}\ \isabound{i}\ \isabound{{\isasymGamma}}\ {\isasymlongrightarrow}\ \isavar{ent}\ \isabound{{\isasymGamma}}\ {\isacharparenleft}{\kern0pt}{\isacharbrackleft}{\kern0pt}\isabound{i}{\isasymZpinj}\isabound{a}{\isacharbrackright}{\kern0pt}\isabound{c}{\isacharparenright}{\kern0pt}{\isacharparenright}{\kern0pt}\ {\isasymand}\ {\isacharparenleft}{\kern0pt}{\isasymforall}\isabound{{\isasymGamma}}\ \isabound{p}{\isachardot}{\kern0pt}\ \isavar{ent}\ \isabound{{\isasymGamma}}\ {\isacharparenleft}{\kern0pt}\isaconst{gigc{\isachardot}{\kern0pt}trmNat}\ \isaconst{gigc{\isachardot}{\kern0pt}trmBot}{\isacharparenright}{\kern0pt}\ {\isasymlongrightarrow}\ \isaconst{gigc{\isachardot}{\kern0pt}wft}\ \isabound{p}\ {\isasymlongrightarrow}\ \isavar{ent}\ \isabound{{\isasymGamma}}\ \isabound{p}{\isacharparenright}{\kern0pt}{\isacharparenright}{\kern0pt}\ {\isasymand}\ {\isacharparenleft}{\kern0pt}{\isasymforall}\isabound{{\isasymGamma}}\ \isabound{c{\isacharprime}{\kern0pt}}\ \isabound{c}{\isachardot}{\kern0pt}\ \isavar{ent}\ \isabound{{\isasymGamma}}\ \isabound{c{\isacharprime}{\kern0pt}}\ {\isasymlongrightarrow}\ \isaconst{gigc{\isachardot}{\kern0pt}hstep}\ \isabound{c}\ \isabound{c{\isacharprime}{\kern0pt}}\ {\isasymlongrightarrow}\ \isavar{ent}\ \isabound{{\isasymGamma}}\ \isabound{c}{\isacharparenright}{\kern0pt}\ {\isasymand}\ {\isacharparenleft}{\kern0pt}{\isasymforall}\isabound{{\isasymGamma}}\ \isabound{c}\ \isabound{c{\isacharprime}{\kern0pt}}{\isachardot}{\kern0pt}\ \isavar{ent}\ \isabound{{\isasymGamma}}\ \isabound{c}\ {\isasymlongrightarrow}\ \isaconst{gigc{\isachardot}{\kern0pt}hstep}\ \isabound{c}\ \isabound{c{\isacharprime}{\kern0pt}}\ {\isasymlongrightarrow}\ \isavar{ent}\ \isabound{{\isasymGamma}}\ \isabound{c{\isacharprime}{\kern0pt}}{\isacharparenright}{\kern0pt}\ {\isasymand}\ {\isacharparenleft}{\kern0pt}{\isasymforall}\isabound{{\isasymGamma}}\ \isabound{c{\isacharprime}{\kern0pt}}\ \isabound{c}{\isachardot}{\kern0pt}\ \isavar{ent}\ \isabound{{\isasymGamma}}\ {\isacharparenleft}{\kern0pt}\isaconst{gigc{\isachardot}{\kern0pt}trmNat}\ \isabound{c{\isacharprime}{\kern0pt}}{\isacharparenright}{\kern0pt}\ {\isasymlongrightarrow}\ \isaconst{gigc{\isachardot}{\kern0pt}hstep}\ \isabound{c}\ \isabound{c{\isacharprime}{\kern0pt}}\ {\isasymlongrightarrow}\ \isavar{ent}\ \isabound{{\isasymGamma}}\ {\isacharparenleft}{\kern0pt}\isaconst{gigc{\isachardot}{\kern0pt}trmEq}\ \isabound{c}\ \isabound{c{\isacharprime}{\kern0pt}}{\isacharparenright}{\kern0pt}{\isacharparenright}{\kern0pt}{\isacharparenright}{\kern0pt}\ {\isasymand}\ {\isacharparenleft}{\kern0pt}{\isacharparenleft}{\kern0pt}{\isasymforall}\isabound{{\isasymGamma}}\ \isabound{c}\ \isabound{c{\isacharprime}{\kern0pt}}{\isachardot}{\kern0pt}\ \isavar{ent}\ \isabound{{\isasymGamma}}\ {\isacharparenleft}{\kern0pt}\isaconst{gigc{\isachardot}{\kern0pt}trmNeg}\ \isabound{c}{\isacharparenright}{\kern0pt}\ {\isasymlongrightarrow}\ \isaconst{gigc{\isachardot}{\kern0pt}hstep}\ \isabound{c}\ \isabound{c{\isacharprime}{\kern0pt}}\ {\isasymlongrightarrow}\ \isavar{ent}\ \isabound{{\isasymGamma}}\ {\isacharparenleft}{\kern0pt}\isaconst{gigc{\isachardot}{\kern0pt}trmNeg}\ \isabound{c{\isacharprime}{\kern0pt}}{\isacharparenright}{\kern0pt}{\isacharparenright}{\kern0pt}\ {\isasymand}\ {\isacharparenleft}{\kern0pt}{\isasymforall}\isabound{{\isasymGamma}}\ \isabound{m}\ \isabound{h}\ \isabound{g}{\isachardot}{\kern0pt}\ \isavar{ent}\ \isabound{{\isasymGamma}}\ {\isacharparenleft}{\kern0pt}\isaconst{gigc{\isachardot}{\kern0pt}trmNat}\ \isabound{m}{\isacharparenright}{\kern0pt}\ {\isasymlongrightarrow}\ \isavar{ent}\ \isabound{{\isasymGamma}}\ {\isacharparenleft}{\kern0pt}\isaconst{gigc{\isachardot}{\kern0pt}trmNat}\ {\isacharparenleft}{\kern0pt}\isaconst{gigc{\isachardot}{\kern0pt}trmApp}\ {\isacharparenleft}{\kern0pt}\isaconst{gigc{\isachardot}{\kern0pt}trmApp}\ \isabound{h}\ \isabound{m}{\isacharparenright}{\kern0pt}\ {\isacharparenleft}{\kern0pt}\isaconst{gigc{\isachardot}{\kern0pt}trmApp}\ {\isacharparenleft}{\kern0pt}\isaconst{gigc{\isachardot}{\kern0pt}trmApp}\ {\isacharparenleft}{\kern0pt}\isaconst{gigc{\isachardot}{\kern0pt}trmApp}\ \isactrlbold R\ \isabound{g}{\isacharparenright}{\kern0pt}\ \isabound{h}{\isacharparenright}{\kern0pt}\ \isabound{m}{\isacharparenright}{\kern0pt}{\isacharparenright}{\kern0pt}{\isacharparenright}{\kern0pt}\ {\isasymlongrightarrow}\ \isaconst{gigc{\isachardot}{\kern0pt}wft}\ \isabound{g}\ {\isasymlongrightarrow}\ \isavar{ent}\ \isabound{{\isasymGamma}}\ {\isacharparenleft}{\kern0pt}\isaconst{gigc{\isachardot}{\kern0pt}trmEq}\ {\isacharparenleft}{\kern0pt}\isaconst{gigc{\isachardot}{\kern0pt}trmApp}\ {\isacharparenleft}{\kern0pt}\isaconst{gigc{\isachardot}{\kern0pt}trmApp}\ {\isacharparenleft}{\kern0pt}\isaconst{gigc{\isachardot}{\kern0pt}trmApp}\ \isactrlbold R\ \isabound{g}{\isacharparenright}{\kern0pt}\ \isabound{h}{\isacharparenright}{\kern0pt}\ {\isacharparenleft}{\kern0pt}\isaconst{gigc{\isachardot}{\kern0pt}trmSuc}\ \isabound{m}{\isacharparenright}{\kern0pt}{\isacharparenright}{\kern0pt}\ {\isacharparenleft}{\kern0pt}\isaconst{gigc{\isachardot}{\kern0pt}trmApp}\ {\isacharparenleft}{\kern0pt}\isaconst{gigc{\isachardot}{\kern0pt}trmApp}\ \isabound{h}\ \isabound{m}{\isacharparenright}{\kern0pt}\ {\isacharparenleft}{\kern0pt}\isaconst{gigc{\isachardot}{\kern0pt}trmApp}\ {\isacharparenleft}{\kern0pt}\isaconst{gigc{\isachardot}{\kern0pt}trmApp}\ {\isacharparenleft}{\kern0pt}\isaconst{gigc{\isachardot}{\kern0pt}trmApp}\ \isactrlbold R\ \isabound{g}{\isacharparenright}{\kern0pt}\ \isabound{h}{\isacharparenright}{\kern0pt}\ \isabound{m}{\isacharparenright}{\kern0pt}{\isacharparenright}{\kern0pt}{\isacharparenright}{\kern0pt}{\isacharparenright}{\kern0pt}\ {\isasymand}\ {\isacharparenleft}{\kern0pt}{\isasymforall}\isabound{{\isasymGamma}}\ \isabound{f}{\isachardot}{\kern0pt}\ \isavar{ent}\ {\isacharparenleft}{\kern0pt}\isaconst{gigc{\isachardot}{\kern0pt}trmNat}\ \isactrlbold v\isactrlsub {\isadigit{0}}\ {\isasymtriangleright}\ {\isasymupharpoonleft}\isabound{{\isasymGamma}}{\isacharparenright}{\kern0pt}\ {\isacharparenleft}{\kern0pt}\isaconst{gigc{\isachardot}{\kern0pt}trmApp}\ {\isacharparenleft}{\kern0pt}{\isasymupharpoonleft}\isabound{f}{\isacharparenright}{\kern0pt}\ \isactrlbold v\isactrlsub {\isadigit{0}}{\isacharparenright}{\kern0pt}\ {\isasymlongrightarrow}\ \isavar{ent}\ \isabound{{\isasymGamma}}\ {\isacharparenleft}{\kern0pt}\isaconst{gigc{\isachardot}{\kern0pt}trmApp}\ \isactrlbold A\ \isabound{f}{\isacharparenright}{\kern0pt}{\isacharparenright}{\kern0pt}{\isacharparenright}{\kern0pt}\ {\isasymand}\ {\isacharparenleft}{\kern0pt}{\isasymforall}\isabound{{\isasymGamma}}\ \isabound{f}\ \isabound{a}{\isachardot}{\kern0pt}\ \isavar{ent}\ \isabound{{\isasymGamma}}\ {\isacharparenleft}{\kern0pt}\isaconst{gigc{\isachardot}{\kern0pt}trmApp}\ \isactrlbold A\ \isabound{f}{\isacharparenright}{\kern0pt}\ {\isasymlongrightarrow}\ \isavar{ent}\ \isabound{{\isasymGamma}}\ {\isacharparenleft}{\kern0pt}\isaconst{gigc{\isachardot}{\kern0pt}trmNat}\ \isabound{a}{\isacharparenright}{\kern0pt}\ {\isasymlongrightarrow}\ \isavar{ent}\ \isabound{{\isasymGamma}}\ {\isacharparenleft}{\kern0pt}\isaconst{gigc{\isachardot}{\kern0pt}trmApp}\ \isabound{f}\ \isabound{a}{\isacharparenright}{\kern0pt}{\isacharparenright}{\kern0pt}\ {\isasymand}\ {\isacharparenleft}{\kern0pt}{\isasymforall}\isabound{{\isasymGamma}}\ \isabound{f}\ \isabound{a}{\isachardot}{\kern0pt}\ \isavar{ent}\ \isabound{{\isasymGamma}}\ {\isacharparenleft}{\kern0pt}\isaconst{gigc{\isachardot}{\kern0pt}trmNeg}\ {\isacharparenleft}{\kern0pt}\isaconst{gigc{\isachardot}{\kern0pt}trmApp}\ \isabound{f}\ \isabound{a}{\isacharparenright}{\kern0pt}{\isacharparenright}{\kern0pt}\ {\isasymlongrightarrow}\ \isavar{ent}\ \isabound{{\isasymGamma}}\ {\isacharparenleft}{\kern0pt}\isaconst{gigc{\isachardot}{\kern0pt}trmNat}\ \isabound{a}{\isacharparenright}{\kern0pt}\ {\isasymlongrightarrow}\ \isavar{ent}\ \isabound{{\isasymGamma}}\ {\isacharparenleft}{\kern0pt}\isaconst{gigc{\isachardot}{\kern0pt}trmNeg}\ {\isacharparenleft}{\kern0pt}\isaconst{gigc{\isachardot}{\kern0pt}trmApp}\ \isactrlbold A\ \isabound{f}{\isacharparenright}{\kern0pt}{\isacharparenright}{\kern0pt}{\isacharparenright}{\kern0pt}\ {\isasymand}\ {\isacharparenleft}{\kern0pt}{\isasymforall}\isabound{{\isasymGamma}}\ \isabound{f}\ \isabound{c}{\isachardot}{\kern0pt}\ \isavar{ent}\ \isabound{{\isasymGamma}}\ {\isacharparenleft}{\kern0pt}\isaconst{gigc{\isachardot}{\kern0pt}trmNeg}\ {\isacharparenleft}{\kern0pt}\isaconst{gigc{\isachardot}{\kern0pt}trmApp}\ \isactrlbold A\ \isabound{f}{\isacharparenright}{\kern0pt}{\isacharparenright}{\kern0pt}\ {\isasymlongrightarrow}\ \isavar{ent}\ {\isacharparenleft}{\kern0pt}\isaconst{gigc{\isachardot}{\kern0pt}trmNat}\ \isactrlbold v\isactrlsub {\isadigit{0}}\ {\isasymtriangleright}\ \isaconst{gigc{\isachardot}{\kern0pt}trmNeg}\ {\isacharparenleft}{\kern0pt}\isaconst{gigc{\isachardot}{\kern0pt}trmApp}\ {\isacharparenleft}{\kern0pt}{\isasymupharpoonleft}\isabound{f}{\isacharparenright}{\kern0pt}\ \isactrlbold v\isactrlsub {\isadigit{0}}{\isacharparenright}{\kern0pt}\ {\isasymtriangleright}\ {\isasymupharpoonleft}\isabound{{\isasymGamma}}{\isacharparenright}{\kern0pt}\ {\isacharparenleft}{\kern0pt}{\isasymupharpoonleft}\isabound{c}{\isacharparenright}{\kern0pt}\ {\isasymlongrightarrow}\ \isavar{ent}\ \isabound{{\isasymGamma}}\ \isabound{c}{\isacharparenright}{\kern0pt}}%
}%
\SNIP{Goedel-Paper-Statements:thm:gcalc-strip-axioms-def-0}{%
\isakeywordONE{thm}\isamarkupfalse%
\ gcalc{\isacharunderscore}{\kern0pt}strip{\isacharunderscore}{\kern0pt}axioms{\isacharunderscore}{\kern0pt}def%
\isanewline%
}%
\SNIP{Goedel-Paper-Statements:thm:gcalc-strip-axioms-def-1}{%
\isa{\isaconst{gcalc{\isacharunderscore}{\kern0pt}strip{\isacharunderscore}{\kern0pt}axioms}\ \isavar{ent}\ {\isasymequiv}\ {\isasymforall}\isabound{{\isasymGamma}}\ \isabound{c}{\isachardot}{\kern0pt}\ \isavar{ent}\ \isabound{{\isasymGamma}}\ \isabound{c}\ {\isasymlongrightarrow}\ \isavar{ent}\ {\isacharparenleft}{\kern0pt}\isaconst{ffilter}\ \isaconst{gigc{\isachardot}{\kern0pt}wft}\ \isabound{{\isasymGamma}}{\isacharparenright}{\kern0pt}\ \isabound{c}}%
}%
\SNIP{Goedel-Paper-Statements:thm:g2-syn-axioms-def-0}{%
\isakeywordONE{thm}\isamarkupfalse%
\ g{\isadigit{2}}{\isacharunderscore}{\kern0pt}syn{\isacharunderscore}{\kern0pt}axioms{\isacharunderscore}{\kern0pt}def%
\isanewline%
}%
\SNIP{Goedel-Paper-Statements:thm:g2-syn-axioms-def-1}{%
\isa{\isaconst{g{\isadigit{2}}{\isacharunderscore}{\kern0pt}syn{\isacharunderscore}{\kern0pt}axioms}\ \isavar{ent}\ \isavar{chk}\ {\isasymequiv}\ {\isacharparenleft}{\kern0pt}{\isasymforall}\isabound{f}\ \isabound{p}\ \isabound{c}{\isachardot}{\kern0pt}\ \isabound{f}\ {\isasymdownharpoonleft}\ \isavar{chk}\ \isabound{p}\ \isabound{c}\ {\isacharequal}{\kern0pt}\ \isavar{chk}\ {\isacharparenleft}{\kern0pt}\isabound{f}\ {\isasymdownharpoonleft}\ \isabound{p}{\isacharparenright}{\kern0pt}\ {\isacharparenleft}{\kern0pt}\isabound{f}\ {\isasymdownharpoonleft}\ \isabound{c}{\isacharparenright}{\kern0pt}{\isacharparenright}{\kern0pt}\ {\isasymand}\ {\isacharparenleft}{\kern0pt}{\isasymforall}\isabound{p}\ \isabound{c}{\isachardot}{\kern0pt}\ \isaconst{gigc{\isachardot}{\kern0pt}wft}\ {\isacharparenleft}{\kern0pt}\isavar{chk}\ \isabound{p}\ \isabound{c}{\isacharparenright}{\kern0pt}\ {\isacharequal}{\kern0pt}\ {\isacharparenleft}{\kern0pt}\isaconst{gigc{\isachardot}{\kern0pt}wft}\ \isabound{p}\ {\isasymand}\ \isaconst{gigc{\isachardot}{\kern0pt}wft}\ \isabound{c}{\isacharparenright}{\kern0pt}{\isacharparenright}{\kern0pt}\ {\isasymand}\ {\isacharparenleft}{\kern0pt}{\isasymforall}\isabound{{\isasymGamma}}\ \isabound{p}\ \isabound{c}{\isachardot}{\kern0pt}\ \isavar{ent}\ \isabound{{\isasymGamma}}\ {\isacharparenleft}{\kern0pt}\isaconst{gigc{\isachardot}{\kern0pt}trmNat}\ \isabound{p}{\isacharparenright}{\kern0pt}\ {\isasymlongrightarrow}\ \isavar{ent}\ \isabound{{\isasymGamma}}\ {\isacharparenleft}{\kern0pt}\isaconst{gigc{\isachardot}{\kern0pt}trmNat}\ \isabound{c}{\isacharparenright}{\kern0pt}\ {\isasymlongrightarrow}\ \isavar{ent}\ \isabound{{\isasymGamma}}\ {\isacharparenleft}{\kern0pt}\isaconst{gigc{\isachardot}{\kern0pt}trmTypeBool}\ {\isacharparenleft}{\kern0pt}\isavar{chk}\ \isabound{p}\ \isabound{c}{\isacharparenright}{\kern0pt}{\isacharparenright}{\kern0pt}{\isacharparenright}{\kern0pt}}%
}%
\SNIP{Goedel-Paper-Statements:thm:g2-con-axioms-def-0}{%
\isakeywordONE{thm}\isamarkupfalse%
\ g{\isadigit{2}}{\isacharunderscore}{\kern0pt}con{\isacharunderscore}{\kern0pt}axioms{\isacharunderscore}{\kern0pt}def%
\isanewline%
}%
\SNIP{Goedel-Paper-Statements:thm:g2-con-axioms-def-1}{%
\isa{\isaconst{g{\isadigit{2}}{\isacharunderscore}{\kern0pt}con{\isacharunderscore}{\kern0pt}axioms}\ \isavar{ent}\ {\isasymequiv}\ {\isasymnot}\ \isavar{ent}\ {\isasymemptyset}\ {\isacharparenleft}{\kern0pt}\isaconst{gigc{\isachardot}{\kern0pt}trmEq}\ \isaconst{gigc{\isachardot}{\kern0pt}trm{\isadigit{0}}}\ {\isacharparenleft}{\kern0pt}\isaconst{gigc{\isachardot}{\kern0pt}trm{\isacharunderscore}{\kern0pt}num}\ {\isadigit{1}}{\isacharparenright}{\kern0pt}{\isacharparenright}{\kern0pt}}%
}%
\SNIP{Goedel-Paper-Statements:thm:gig-gcalc-strip-0}{%
\isakeywordONE{thm}\isamarkupfalse%
\ gig{\isacharunderscore}{\kern0pt}gcalc{\isacharunderscore}{\kern0pt}strip%
\isanewline%
}%
\SNIP{Goedel-Paper-Statements:thm:gig-gcalc-strip-1}{%
\isa{\isaconst{gcalc{\isacharunderscore}{\kern0pt}strip}\ {\isacharparenleft}{\kern0pt}\isaconst{prf{\isacharunderscore}{\kern0pt}derives{\isachardot}{\kern0pt}entails}\ \isaconst{gi{\isachardot}{\kern0pt}derives}{\isacharparenright}{\kern0pt}}%
}%
\SNIP{Goedel-Paper-Statements:context:rga-0}{%
\isakeywordONE{context}\isamarkupfalse%
\ rga%
}%
\SNIP{Goedel-Paper-Statements:notation:gprop-0}{%
\isakeywordONE{notation}\isamarkupfalse%
\ gprop{\isachardot}{\kern0pt}trmIff{\isacharprime}{\kern0pt}\ {\isacharparenleft}{\kern0pt}\isakeywordTWO{infixr}\ {\isachardoublequoteopen}\isactrlbold {\isasymlongleftrightarrow}{\isachardoublequoteclose}\ {\isadigit{1}}{\isadigit{2}}{\isadigit{5}}{\isacharparenright}{\kern0pt}%
}%
\SNIP{Goedel-Paper-Statements:thm:deduction-theorem-exact-0}{%
\isakeywordONE{thm}\isamarkupfalse%
\ deduction{\isacharunderscore}{\kern0pt}theorem{\isacharunderscore}{\kern0pt}exact%
\isanewline%
}%
\SNIP{Goedel-Paper-Statements:thm:deduction-theorem-exact-1}{%
\isa{\isaconst{wft}\ \isavar{A}\ {\isasymLongrightarrow}\ {\isacharparenleft}{\kern0pt}{\isasymforall}\isabound{B}{\isachardot}{\kern0pt}\ \isavar{A}\ {\isasymtriangleright}\ \isavar{{\isasymGamma}}\ {\isasymturnstile}\ \isabound{B}\ {\isasymlongrightarrow}\ \isavar{{\isasymGamma}}\ {\isasymturnstile}\ \isavar{A}\ \isactrlbold {\isasymlongrightarrow}\ \isabound{B}{\isacharparenright}{\kern0pt}\ {\isacharequal}{\kern0pt}\ {\isacharparenleft}{\kern0pt}\isavar{{\isasymGamma}}\ {\isasymturnstile}\ \isavar{A}\ \isactrlbold {\isasymor}\ \isactrlbold {\isasymnot}\isavar{A}{\isacharparenright}{\kern0pt}}%
}%
\SNIP{Goedel-Paper-Statements:thm:sigma-deduction-theorem-0}{%
\isakeywordONE{thm}\isamarkupfalse%
\ sigma{\isacharunderscore}{\kern0pt}deduction{\isacharunderscore}{\kern0pt}theorem%
\isanewline%
}%
\SNIP{Goedel-Paper-Statements:thm:sigma-deduction-theorem-1}{%
\isa{\isactrlbold v\isactrlsub {\isadigit{0}}\ \isactrlbold N\ {\isasymtriangleright}\ {\isasymupharpoonleft}\isavar{{\isasymGamma}}\ {\isasymturnstile}\ \isavar{D}\ \isactrlbold B\ {\isasymLongrightarrow}\ {\isacharparenleft}{\kern0pt}\isactrlbold {\isasymnot}\isactrlbold A\ \isactrlbold {\isasymcdot}\ {\isacharparenleft}{\kern0pt}\isactrlbold {\isasymlambda}{\isacharparenleft}{\kern0pt}\isactrlbold {\isasymnot}\isavar{D}{\isacharparenright}{\kern0pt}{\isacharparenright}{\kern0pt}\ {\isasymtriangleright}\ \isavar{{\isasymGamma}}\ {\isasymturnstile}\ \isavar{C}{\isacharparenright}{\kern0pt}\ {\isacharequal}{\kern0pt}\ {\isacharparenleft}{\kern0pt}\isavar{{\isasymGamma}}\ {\isasymturnstile}\ \isactrlbold A\ \isactrlbold {\isasymcdot}\ {\isacharparenleft}{\kern0pt}\isactrlbold {\isasymlambda}{\isacharparenleft}{\kern0pt}\isavar{D}\ \isactrlbold {\isasymlongrightarrow}\ {\isasymupharpoonleft}\isavar{C}{\isacharparenright}{\kern0pt}{\isacharparenright}{\kern0pt}{\isacharparenright}{\kern0pt}}%
}%
\SNIP{Goedel-Paper-Statements:thm:dec-contrapos-0}{%
\isakeywordONE{thm}\isamarkupfalse%
\ dec{\isacharunderscore}{\kern0pt}contrapos%
\isanewline%
}%
\SNIP{Goedel-Paper-Statements:thm:dec-contrapos-1}{%
\isa{{\isasymlbrakk}\isavar{{\isasymGamma}}\ {\isasymturnstile}\ \isavar{A}\ \isactrlbold {\isasymor}\ \isactrlbold {\isasymnot}\isavar{A}{\isacharsemicolon}{\kern0pt}\ \isavar{A}\ {\isasymtriangleright}\ \isavar{{\isasymGamma}}\ {\isasymturnstile}\ \isavar{B}{\isasymrbrakk}\ {\isasymLongrightarrow}\ \isactrlbold {\isasymnot}\isavar{B}\ {\isasymtriangleright}\ \isavar{{\isasymGamma}}\ {\isasymturnstile}\ \isactrlbold {\isasymnot}\isavar{A}}%
}%
\SNIP{Goedel-Paper-Statements:thm:sigma-contrapos-0}{%
\isakeywordONE{thm}\isamarkupfalse%
\ sigma{\isacharunderscore}{\kern0pt}contrapos%
\isanewline%
}%
\SNIP{Goedel-Paper-Statements:thm:sigma-contrapos-1}{%
\isa{{\isasymlbrakk}\isactrlbold v\isactrlsub {\isadigit{0}}\ \isactrlbold N\ {\isasymtriangleright}\ {\isasymupharpoonleft}\isavar{{\isasymGamma}}\ {\isasymturnstile}\ \isavar{D}\ \isactrlbold B{\isacharsemicolon}{\kern0pt}\ \isactrlbold {\isasymnot}\isactrlbold A\ \isactrlbold {\isasymcdot}\ {\isacharparenleft}{\kern0pt}\isactrlbold {\isasymlambda}{\isacharparenleft}{\kern0pt}\isactrlbold {\isasymnot}\isavar{D}{\isacharparenright}{\kern0pt}{\isacharparenright}{\kern0pt}\ {\isasymtriangleright}\ \isavar{{\isasymGamma}}\ {\isasymturnstile}\ \isavar{B}{\isasymrbrakk}\ {\isasymLongrightarrow}\ \isactrlbold {\isasymnot}\isavar{B}\ {\isasymtriangleright}\ \isavar{{\isasymGamma}}\ {\isasymturnstile}\ \isactrlbold {\isasymnot}\isactrlbold {\isasymnot}\isactrlbold A\ \isactrlbold {\isasymcdot}\ {\isacharparenleft}{\kern0pt}\isactrlbold {\isasymlambda}{\isacharparenleft}{\kern0pt}\isactrlbold {\isasymnot}\isavar{D}{\isacharparenright}{\kern0pt}{\isacharparenright}{\kern0pt}}%
}%
\SNIP{Goedel-Paper-Statements:thm:ir-disj-prop-0}{%
\isakeywordONE{thm}\isamarkupfalse%
\ ir{\isacharunderscore}{\kern0pt}disj{\isacharunderscore}{\kern0pt}prop%
\isanewline%
}%
\SNIP{Goedel-Paper-Statements:thm:ir-disj-prop-1}{%
\isa{{\isasymemptyset}\ {\isasymturnstile}\ \isavar{a}\ \isactrlbold {\isasymor}\ \isavar{b}\ {\isasymLongrightarrow}\ {\isasymemptyset}\ {\isasymturnstile}\ \isavar{a}\ {\isasymor}\ {\isasymemptyset}\ {\isasymturnstile}\ \isavar{b}}%
}%
\SNIP{Goedel-Paper-Statements:thm:rga-not-all-witness-0}{%
\isakeywordONE{thm}\isamarkupfalse%
\ rga{\isacharunderscore}{\kern0pt}not{\isacharunderscore}{\kern0pt}all{\isacharunderscore}{\kern0pt}witness%
\isanewline%
}%
\SNIP{Goedel-Paper-Statements:thm:rga-not-all-witness-1}{%
\isa{{\isasymemptyset}\ {\isasymturnstile}\ \isactrlbold {\isasymnot}\isactrlbold A\ \isactrlbold {\isasymcdot}\ \isavar{u}\ {\isasymLongrightarrow}\ {\isasymexists}\isabound{n}{\isachardot}{\kern0pt}\ {\isasymemptyset}\ {\isasymturnstile}\ \isactrlbold {\isasymnot}\isavar{u}\ \isactrlbold {\isasymcdot}\ {\isacharparenleft}{\kern0pt}\isactrlbold n\isabound{n}{\isacharparenright}{\kern0pt}}%
}%
\SNIP{Goedel-Paper-Statements:thm:no-goedel-biconditional-0}{%
\isakeywordONE{thm}\isamarkupfalse%
\ no{\isacharunderscore}{\kern0pt}goedel{\isacharunderscore}{\kern0pt}biconditional%
\isanewline%
}%
\SNIP{Goedel-Paper-Statements:thm:no-goedel-biconditional-1}{%
\isa{{\isasymnot}\ {\isasymemptyset}\ {\isasymturnstile}\ \isavar{d}\ \isactrlbold {\isasymlongleftrightarrow}\ \isactrlbold A\ \isactrlbold {\isasymcdot}\ {\isacharparenleft}{\kern0pt}\isactrlbold {\isasymlambda}{\isacharparenleft}{\kern0pt}\isactrlbold {\isasymnot}\isaconst{iChk}\ \isactrlbold v\isactrlsub {\isadigit{0}}\ {\isacharparenleft}{\kern0pt}\isactrlbold n{\isasymlbrace}\isavar{d}{\isasymrbrace}{\isacharparenright}{\kern0pt}{\isacharparenright}{\kern0pt}{\isacharparenright}{\kern0pt}}%
}%
\SNIP{Goedel-Paper-Statements:thm:no-liar-biconditional-0}{%
\isakeywordONE{thm}\isamarkupfalse%
\ no{\isacharunderscore}{\kern0pt}liar{\isacharunderscore}{\kern0pt}biconditional%
\isanewline%
}%
\SNIP{Goedel-Paper-Statements:thm:no-liar-biconditional-1}{%
\isa{{\isasymnot}\ {\isasymemptyset}\ {\isasymturnstile}\ \isavar{d}\ \isactrlbold {\isasymlongleftrightarrow}\ \isactrlbold {\isasymnot}\isaconst{iWtrue}\ {\isacharparenleft}{\kern0pt}\isactrlbold n{\isasymlbrace}\isavar{d}{\isasymrbrace}{\isacharparenright}{\kern0pt}}%
}%
\SNIP{Goedel-Paper-Statements:thm:goedel-diag-ent-0}{%
\isakeywordONE{thm}\isamarkupfalse%
\ goedel{\isacharunderscore}{\kern0pt}diag{\isacharunderscore}{\kern0pt}ent%
\isanewline%
}%
\SNIP{Goedel-Paper-Statements:thm:goedel-diag-ent-1}{%
\isa{\isaconst{goedelG}\ {\isasymtriangleright}\ {\isasymemptyset}\ {\isasymturnstile}\ \isactrlbold {\isasymnot}\isactrlbold {\isasymnot}\isactrlbold A\ \isactrlbold {\isasymcdot}\ {\isacharparenleft}{\kern0pt}\isactrlbold {\isasymlambda}{\isacharparenleft}{\kern0pt}\isactrlbold {\isasymnot}\isaconst{iChk}\ \isactrlbold v\isactrlsub {\isadigit{0}}\ {\isacharparenleft}{\kern0pt}{\isasymupharpoonleft}{\isacharparenleft}{\kern0pt}\isactrlbold n{\isasymlbrace}\isaconst{goedelG}{\isasymrbrace}{\isacharparenright}{\kern0pt}{\isacharparenright}{\kern0pt}{\isacharparenright}{\kern0pt}{\isacharparenright}{\kern0pt}\ {\isasymand}\ \isactrlbold {\isasymnot}\isactrlbold {\isasymnot}\isactrlbold A\ \isactrlbold {\isasymcdot}\ {\isacharparenleft}{\kern0pt}\isactrlbold {\isasymlambda}{\isacharparenleft}{\kern0pt}\isactrlbold {\isasymnot}\isaconst{iChk}\ \isactrlbold v\isactrlsub {\isadigit{0}}\ {\isacharparenleft}{\kern0pt}{\isasymupharpoonleft}{\isacharparenleft}{\kern0pt}\isactrlbold n{\isasymlbrace}\isaconst{goedelG}{\isasymrbrace}{\isacharparenright}{\kern0pt}{\isacharparenright}{\kern0pt}{\isacharparenright}{\kern0pt}{\isacharparenright}{\kern0pt}\ {\isasymtriangleright}\ {\isasymemptyset}\ {\isasymturnstile}\ \isaconst{goedelG}}%
}%
\SNIP{Goedel-Paper-Statements:thm:liar-diag-ent-0}{%
\isakeywordONE{thm}\isamarkupfalse%
\ liar{\isacharunderscore}{\kern0pt}diag{\isacharunderscore}{\kern0pt}ent%
\isanewline%
}%
\SNIP{Goedel-Paper-Statements:thm:liar-diag-ent-1}{%
\isa{\isaconst{goedelW}\ {\isasymtriangleright}\ {\isasymemptyset}\ {\isasymturnstile}\ \isactrlbold {\isasymnot}\isaconst{iWtrue}\ {\isacharparenleft}{\kern0pt}\isactrlbold n{\isasymlbrace}\isaconst{goedelW}{\isasymrbrace}{\isacharparenright}{\kern0pt}\ {\isasymand}\ \isactrlbold {\isasymnot}\isaconst{iWtrue}\ {\isacharparenleft}{\kern0pt}\isactrlbold n{\isasymlbrace}\isaconst{goedelW}{\isasymrbrace}{\isacharparenright}{\kern0pt}\ {\isasymtriangleright}\ {\isasymemptyset}\ {\isasymturnstile}\ \isaconst{goedelW}}%
}%
\SNIP{Goedel-Paper-Statements:thm:no-formalized-goedel1-0}{%
\isakeywordONE{thm}\isamarkupfalse%
\ no{\isacharunderscore}{\kern0pt}formalized{\isacharunderscore}{\kern0pt}goedel{\isadigit{1}}%
\isanewline%
}%
\SNIP{Goedel-Paper-Statements:thm:no-formalized-goedel1-1}{%
\isa{{\isasymnot}\ {\isasymemptyset}\ {\isasymturnstile}\ \isaconst{ConRGA}\ \isactrlbold {\isasymlongrightarrow}\ \isaconst{goedelG}}%
}%
\SNIP{Goedel-Paper-Statements:thm:no-formalized-goedel1-code-0}{%
\isakeywordONE{thm}\isamarkupfalse%
\ no{\isacharunderscore}{\kern0pt}formalized{\isacharunderscore}{\kern0pt}goedel{\isadigit{1}}{\isacharunderscore}{\kern0pt}code%
\isanewline%
}%
\SNIP{Goedel-Paper-Statements:thm:no-formalized-goedel1-code-1}{%
\isa{{\isasymnot}\ {\isasymemptyset}\ {\isasymturnstile}\ \isaconst{ConRGA}\ \isactrlbold {\isasymlongrightarrow}\ \isactrlbold {\isasymnot}\isactrlbold {\isasymnot}\isactrlbold A\ \isactrlbold {\isasymcdot}\ {\isacharparenleft}{\kern0pt}\isactrlbold {\isasymlambda}{\isacharparenleft}{\kern0pt}\isactrlbold {\isasymnot}\isaconst{iChk}\ \isactrlbold v\isactrlsub {\isadigit{0}}\ {\isacharparenleft}{\kern0pt}\isactrlbold n{\isasymlbrace}\isaconst{goedelG}{\isasymrbrace}{\isacharparenright}{\kern0pt}{\isacharparenright}{\kern0pt}{\isacharparenright}{\kern0pt}}%
}%
\SNIP{Goedel-Paper-Statements:thm:ConRGA-unrefutable-0}{%
\isakeywordONE{thm}\isamarkupfalse%
\ ConRGA{\isacharunderscore}{\kern0pt}unrefutable%
\isanewline%
}%
\SNIP{Goedel-Paper-Statements:thm:ConRGA-unrefutable-1}{%
\isa{{\isasymnot}\ {\isasymemptyset}\ {\isasymturnstile}\ \isactrlbold {\isasymnot}\isaconst{ConRGA}}%
}%
\SNIP{Goedel-Paper-Statements:thm:ConRGA-undecided-0}{%
\isakeywordONE{thm}\isamarkupfalse%
\ ConRGA{\isacharunderscore}{\kern0pt}undecided%
\isanewline%
}%
\SNIP{Goedel-Paper-Statements:thm:ConRGA-undecided-1}{%
\isa{{\isasymnot}\ {\isasymemptyset}\ {\isasymturnstile}\ \isaconst{ConRGA}\ \isactrlbold B}%
}%
\SNIP{Goedel-Paper-Statements:thm:goedelG-under-ConRGA-0}{%
\isakeywordONE{thm}\isamarkupfalse%
\ goedelG{\isacharunderscore}{\kern0pt}under{\isacharunderscore}{\kern0pt}ConRGA%
\isanewline%
}%
\SNIP{Goedel-Paper-Statements:thm:goedelG-under-ConRGA-1}{%
\isa{\isaconst{ConRGA}\ {\isasymtriangleright}\ {\isasymemptyset}\ {\isasymturnstile}\ \isaconst{goedelG}}%
}%
\SNIP{Goedel-Paper-Statements:thm:con-notprc-goedelG-ent-0}{%
\isakeywordONE{thm}\isamarkupfalse%
\ con{\isacharunderscore}{\kern0pt}notprc{\isacharunderscore}{\kern0pt}goedelG{\isacharunderscore}{\kern0pt}ent%
\isanewline%
}%
\SNIP{Goedel-Paper-Statements:thm:con-notprc-goedelG-ent-1}{%
\isa{\isaconst{ConRGA}\ {\isasymtriangleright}\ {\isasymemptyset}\ {\isasymturnstile}\ \isactrlbold {\isasymnot}\isactrlbold {\isasymnot}\isactrlbold A\ \isactrlbold {\isasymcdot}\ {\isacharparenleft}{\kern0pt}\isactrlbold {\isasymlambda}{\isacharparenleft}{\kern0pt}\isactrlbold {\isasymnot}\isaconst{iChk}\ \isactrlbold v\isactrlsub {\isadigit{0}}\ {\isacharparenleft}{\kern0pt}{\isasymupharpoonleft}{\isacharparenleft}{\kern0pt}\isactrlbold n{\isasymlbrace}\isaconst{goedelG}{\isasymrbrace}{\isacharparenright}{\kern0pt}{\isacharparenright}{\kern0pt}{\isacharparenright}{\kern0pt}{\isacharparenright}{\kern0pt}}%
}%
\SNIP{Goedel-Paper-Statements:thm:no-D2-schema-goedelG-0}{%
\isakeywordONE{thm}\isamarkupfalse%
\ no{\isacharunderscore}{\kern0pt}D{\isadigit{2}}{\isacharunderscore}{\kern0pt}schema{\isacharunderscore}{\kern0pt}goedelG%
\isanewline%
}%
\SNIP{Goedel-Paper-Statements:thm:no-D2-schema-goedelG-1}{%
\isa{{\isasymnot}\ {\isasymemptyset}\ {\isasymturnstile}\ \isactrlbold {\isasymnot}\isactrlbold A\ \isactrlbold {\isasymcdot}\ {\isacharparenleft}{\kern0pt}\isactrlbold {\isasymlambda}{\isacharparenleft}{\kern0pt}\isactrlbold {\isasymnot}\isaconst{iChk}\ \isactrlbold v\isactrlsub {\isadigit{0}}\ {\isacharparenleft}{\kern0pt}\isactrlbold n{\isasymlbrace}\isaconst{goedelG}\ \isactrlbold {\isasymlongrightarrow}\ \isactrlbold {\isadigit{0}}\ \isactrlbold {\isacharequal}{\kern0pt}\ \isactrlbold n{\isadigit{1}}{\isasymrbrace}{\isacharparenright}{\kern0pt}{\isacharparenright}{\kern0pt}{\isacharparenright}{\kern0pt}\ \isactrlbold {\isasymlongrightarrow}\ \isactrlbold {\isasymnot}\isactrlbold A\ \isactrlbold {\isasymcdot}\ {\isacharparenleft}{\kern0pt}\isactrlbold {\isasymlambda}{\isacharparenleft}{\kern0pt}\isactrlbold {\isasymnot}\isaconst{iChk}\ \isactrlbold v\isactrlsub {\isadigit{0}}\ {\isacharparenleft}{\kern0pt}\isactrlbold n{\isasymlbrace}\isaconst{goedelG}{\isasymrbrace}{\isacharparenright}{\kern0pt}{\isacharparenright}{\kern0pt}{\isacharparenright}{\kern0pt}\ \isactrlbold {\isasymlongrightarrow}\ \isactrlbold {\isasymnot}\isaconst{ConRGA}}%
}%
\SNIP{Goedel-Paper-Statements:thm:no-D3-schema-goedelG-0}{%
\isakeywordONE{thm}\isamarkupfalse%
\ no{\isacharunderscore}{\kern0pt}D{\isadigit{3}}{\isacharunderscore}{\kern0pt}schema{\isacharunderscore}{\kern0pt}goedelG%
\isanewline%
}%
\SNIP{Goedel-Paper-Statements:thm:no-D3-schema-goedelG-1}{%
\isa{{\isasymnot}\ {\isasymemptyset}\ {\isasymturnstile}\ \isactrlbold {\isasymnot}\isactrlbold A\ \isactrlbold {\isasymcdot}\ {\isacharparenleft}{\kern0pt}\isactrlbold {\isasymlambda}{\isacharparenleft}{\kern0pt}\isactrlbold {\isasymnot}\isaconst{iChk}\ \isactrlbold v\isactrlsub {\isadigit{0}}\ {\isacharparenleft}{\kern0pt}\isactrlbold n{\isasymlbrace}\isaconst{goedelG}{\isasymrbrace}{\isacharparenright}{\kern0pt}{\isacharparenright}{\kern0pt}{\isacharparenright}{\kern0pt}\ \isactrlbold {\isasymlongrightarrow}\ \isactrlbold {\isasymnot}\isactrlbold A\ \isactrlbold {\isasymcdot}\ {\isacharparenleft}{\kern0pt}\isactrlbold {\isasymlambda}{\isacharparenleft}{\kern0pt}\isactrlbold {\isasymnot}\isaconst{iChk}\ \isactrlbold v\isactrlsub {\isadigit{0}}\ {\isacharparenleft}{\kern0pt}\isactrlbold n{\isasymlbrace}\isactrlbold {\isasymnot}\isactrlbold A\ \isactrlbold {\isasymcdot}\ {\isacharparenleft}{\kern0pt}\isactrlbold {\isasymlambda}{\isacharparenleft}{\kern0pt}\isactrlbold {\isasymnot}\isaconst{iChk}\ \isactrlbold v\isactrlsub {\isadigit{0}}\ {\isacharparenleft}{\kern0pt}\isactrlbold n{\isasymlbrace}\isaconst{goedelG}{\isasymrbrace}{\isacharparenright}{\kern0pt}{\isacharparenright}{\kern0pt}{\isacharparenright}{\kern0pt}{\isasymrbrace}{\isacharparenright}{\kern0pt}{\isacharparenright}{\kern0pt}{\isacharparenright}{\kern0pt}}%
}%
\SNIP{Goedel-Paper-Statements:thm:D1-rule-pr-0}{%
\isakeywordONE{thm}\isamarkupfalse%
\ D{\isadigit{1}}{\isacharunderscore}{\kern0pt}rule{\isacharunderscore}{\kern0pt}pr%
\isanewline%
}%
\SNIP{Goedel-Paper-Statements:thm:D1-rule-pr-1}{%
\isa{{\isasymemptyset}\ {\isasymturnstile}\ \isavar{{\isasymphi}}\ {\isasymLongrightarrow}\ {\isasymemptyset}\ {\isasymturnstile}\ \isactrlbold {\isasymnot}\isactrlbold A\ \isactrlbold {\isasymcdot}\ {\isacharparenleft}{\kern0pt}\isactrlbold {\isasymlambda}{\isacharparenleft}{\kern0pt}\isactrlbold {\isasymnot}\isaconst{iChk}\ \isactrlbold v\isactrlsub {\isadigit{0}}\ {\isacharparenleft}{\kern0pt}\isactrlbold n{\isasymlbrace}\isavar{{\isasymphi}}{\isasymrbrace}{\isacharparenright}{\kern0pt}{\isacharparenright}{\kern0pt}{\isacharparenright}{\kern0pt}}%
}%
\SNIP{Goedel-Paper-Statements:thm:rga-pr-reflect-0}{%
\isakeywordONE{thm}\isamarkupfalse%
\ rga{\isacharunderscore}{\kern0pt}pr{\isacharunderscore}{\kern0pt}reflect%
\isanewline%
}%
\SNIP{Goedel-Paper-Statements:thm:rga-pr-reflect-1}{%
\isa{{\isasymemptyset}\ {\isasymturnstile}\ \isactrlbold {\isasymnot}\isactrlbold A\ \isactrlbold {\isasymcdot}\ {\isacharparenleft}{\kern0pt}\isactrlbold {\isasymlambda}{\isacharparenleft}{\kern0pt}\isactrlbold {\isasymnot}\isaconst{iChk}\ \isactrlbold v\isactrlsub {\isadigit{0}}\ {\isacharparenleft}{\kern0pt}\isactrlbold n{\isasymlbrace}\isavar{{\isasympsi}}{\isasymrbrace}{\isacharparenright}{\kern0pt}{\isacharparenright}{\kern0pt}{\isacharparenright}{\kern0pt}\ {\isasymLongrightarrow}\ {\isasymemptyset}\ {\isasymturnstile}\ \isavar{{\isasympsi}}}%
}%
\SNIP{Goedel-Paper-Statements:thm:D2-ent-0}{%
\isakeywordONE{thm}\isamarkupfalse%
\ D{\isadigit{2}}{\isacharunderscore}{\kern0pt}ent%
\isanewline%
}%
\SNIP{Goedel-Paper-Statements:thm:D2-ent-1}{%
\isa{{\isasymlbrakk}\isaconst{wft}\ \isavar{{\isasymphi}}{\isacharsemicolon}{\kern0pt}\ \isaconst{wft}\ \isavar{{\isasympsi}}{\isasymrbrakk}\ {\isasymLongrightarrow}\ \isactrlbold {\isasymnot}\isactrlbold A\ \isactrlbold {\isasymcdot}\ {\isacharparenleft}{\kern0pt}\isactrlbold {\isasymlambda}{\isacharparenleft}{\kern0pt}\isactrlbold {\isasymnot}\isaconst{iChk}\ \isactrlbold v\isactrlsub {\isadigit{0}}\ {\isacharparenleft}{\kern0pt}{\isasymupharpoonleft}{\isacharparenleft}{\kern0pt}\isactrlbold n{\isasymlbrace}\isavar{{\isasymphi}}{\isasymrbrace}{\isacharparenright}{\kern0pt}{\isacharparenright}{\kern0pt}{\isacharparenright}{\kern0pt}{\isacharparenright}{\kern0pt}\ {\isasymtriangleright}\ \isactrlbold {\isasymnot}\isactrlbold A\ \isactrlbold {\isasymcdot}\ {\isacharparenleft}{\kern0pt}\isactrlbold {\isasymlambda}{\isacharparenleft}{\kern0pt}\isactrlbold {\isasymnot}\isaconst{iChk}\ \isactrlbold v\isactrlsub {\isadigit{0}}\ {\isacharparenleft}{\kern0pt}{\isasymupharpoonleft}{\isacharparenleft}{\kern0pt}\isactrlbold n{\isasymlbrace}\isavar{{\isasymphi}}\ \isactrlbold {\isasymlongrightarrow}\ \isavar{{\isasympsi}}{\isasymrbrace}{\isacharparenright}{\kern0pt}{\isacharparenright}{\kern0pt}{\isacharparenright}{\kern0pt}{\isacharparenright}{\kern0pt}\ {\isasymtriangleright}\ \isavar{{\isasymGamma}}\ {\isasymturnstile}\ \isactrlbold {\isasymnot}\isactrlbold A\ \isactrlbold {\isasymcdot}\ {\isacharparenleft}{\kern0pt}\isactrlbold {\isasymlambda}{\isacharparenleft}{\kern0pt}\isactrlbold {\isasymnot}\isaconst{iChk}\ \isactrlbold v\isactrlsub {\isadigit{0}}\ {\isacharparenleft}{\kern0pt}{\isasymupharpoonleft}{\isacharparenleft}{\kern0pt}\isactrlbold n{\isasymlbrace}\isavar{{\isasympsi}}{\isasymrbrace}{\isacharparenright}{\kern0pt}{\isacharparenright}{\kern0pt}{\isacharparenright}{\kern0pt}{\isacharparenright}{\kern0pt}}%
}%
\SNIP{Goedel-Paper-Statements:thm:D3-ent-0}{%
\isakeywordONE{thm}\isamarkupfalse%
\ D{\isadigit{3}}{\isacharunderscore}{\kern0pt}ent%
\isanewline%
}%
\SNIP{Goedel-Paper-Statements:thm:D3-ent-1}{%
\isa{\isactrlbold {\isasymnot}\isactrlbold A\ \isactrlbold {\isasymcdot}\ {\isacharparenleft}{\kern0pt}\isactrlbold {\isasymlambda}{\isacharparenleft}{\kern0pt}\isactrlbold {\isasymnot}\isaconst{iChk}\ \isactrlbold v\isactrlsub {\isadigit{0}}\ {\isacharparenleft}{\kern0pt}{\isasymupharpoonleft}{\isacharparenleft}{\kern0pt}\isactrlbold n{\isasymlbrace}\isavar{{\isasymphi}}{\isasymrbrace}{\isacharparenright}{\kern0pt}{\isacharparenright}{\kern0pt}{\isacharparenright}{\kern0pt}{\isacharparenright}{\kern0pt}\ {\isasymtriangleright}\ \isavar{{\isasymGamma}}\ {\isasymturnstile}\ \isactrlbold {\isasymnot}\isactrlbold A\ \isactrlbold {\isasymcdot}\ {\isacharparenleft}{\kern0pt}\isactrlbold {\isasymlambda}{\isacharparenleft}{\kern0pt}\isactrlbold {\isasymnot}\isaconst{iChk}\ \isactrlbold v\isactrlsub {\isadigit{0}}\ {\isacharparenleft}{\kern0pt}{\isasymupharpoonleft}{\isacharparenleft}{\kern0pt}\isactrlbold n{\isasymlbrace}\isactrlbold {\isasymnot}\isactrlbold A\ \isactrlbold {\isasymcdot}\ {\isacharparenleft}{\kern0pt}\isactrlbold {\isasymlambda}{\isacharparenleft}{\kern0pt}\isactrlbold {\isasymnot}\isaconst{iChk}\ \isactrlbold v\isactrlsub {\isadigit{0}}\ {\isacharparenleft}{\kern0pt}{\isasymupharpoonleft}{\isacharparenleft}{\kern0pt}\isactrlbold n{\isasymlbrace}\isavar{{\isasymphi}}{\isasymrbrace}{\isacharparenright}{\kern0pt}{\isacharparenright}{\kern0pt}{\isacharparenright}{\kern0pt}{\isacharparenright}{\kern0pt}{\isasymrbrace}{\isacharparenright}{\kern0pt}{\isacharparenright}{\kern0pt}{\isacharparenright}{\kern0pt}{\isacharparenright}{\kern0pt}}%
}%
\SNIP{Goedel-Paper-Statements:thm:no-code-unprovable-0}{%
\isakeywordONE{thm}\isamarkupfalse%
\ no{\isacharunderscore}{\kern0pt}code{\isacharunderscore}{\kern0pt}unprovable%
\isanewline%
}%
\SNIP{Goedel-Paper-Statements:thm:no-code-unprovable-1}{%
\isa{\isaconst{wft}\ \isavar{{\isasymphi}}\ {\isasymLongrightarrow}\ {\isasymnot}\ {\isasymemptyset}\ {\isasymturnstile}\ \isactrlbold A\ \isactrlbold {\isasymcdot}\ {\isacharparenleft}{\kern0pt}\isactrlbold {\isasymlambda}{\isacharparenleft}{\kern0pt}\isactrlbold {\isasymnot}\isaconst{iChk}\ \isactrlbold v\isactrlsub {\isadigit{0}}\ {\isacharparenleft}{\kern0pt}\isactrlbold n{\isasymlbrace}\isavar{{\isasymphi}}{\isasymrbrace}{\isacharparenright}{\kern0pt}{\isacharparenright}{\kern0pt}{\isacharparenright}{\kern0pt}}%
}%
\SNIP{Goedel-Paper-Statements:thm:rga-loeb-rule-0}{%
\isakeywordONE{thm}\isamarkupfalse%
\ rga{\isacharunderscore}{\kern0pt}loeb{\isacharunderscore}{\kern0pt}rule%
\isanewline%
}%
\SNIP{Goedel-Paper-Statements:thm:rga-loeb-rule-1}{%
\isa{{\isasymemptyset}\ {\isasymturnstile}\ \isactrlbold {\isasymnot}\isactrlbold A\ \isactrlbold {\isasymcdot}\ {\isacharparenleft}{\kern0pt}\isactrlbold {\isasymlambda}{\isacharparenleft}{\kern0pt}\isactrlbold {\isasymnot}\isaconst{iChk}\ \isactrlbold v\isactrlsub {\isadigit{0}}\ {\isacharparenleft}{\kern0pt}\isactrlbold n{\isasymlbrace}\isavar{{\isasymphi}}{\isasymrbrace}{\isacharparenright}{\kern0pt}{\isacharparenright}{\kern0pt}{\isacharparenright}{\kern0pt}\ \isactrlbold {\isasymlongrightarrow}\ \isavar{{\isasymphi}}\ {\isasymLongrightarrow}\ {\isasymemptyset}\ {\isasymturnstile}\ \isavar{{\isasymphi}}}%
}%
\SNIP{Goedel-Paper-Statements:thm:formalized-loeb-ent-0}{%
\isakeywordONE{thm}\isamarkupfalse%
\ formalized{\isacharunderscore}{\kern0pt}loeb{\isacharunderscore}{\kern0pt}ent%
\isanewline%
}%
\SNIP{Goedel-Paper-Statements:thm:formalized-loeb-ent-1}{%
\isa{{\isasymlbrakk}\isaconst{Syntax{\isacharunderscore}{\kern0pt}Term{\isachardot}{\kern0pt}closed}\ \isavar{{\isasympsi}}{\isacharsemicolon}{\kern0pt}\ \isaconst{wft}\ \isavar{{\isasympsi}}{\isasymrbrakk}\ {\isasymLongrightarrow}\ \isactrlbold {\isasymnot}\isactrlbold A\ \isactrlbold {\isasymcdot}\ {\isacharparenleft}{\kern0pt}\isactrlbold {\isasymlambda}{\isacharparenleft}{\kern0pt}\isactrlbold {\isasymnot}\isaconst{iChk}\ \isactrlbold v\isactrlsub {\isadigit{0}}\ {\isacharparenleft}{\kern0pt}{\isasymupharpoonleft}{\isacharparenleft}{\kern0pt}\isactrlbold n{\isasymlbrace}\isactrlbold {\isasymnot}\isactrlbold A\ \isactrlbold {\isasymcdot}\ {\isacharparenleft}{\kern0pt}\isactrlbold {\isasymlambda}{\isacharparenleft}{\kern0pt}\isactrlbold {\isasymnot}\isaconst{iChk}\ \isactrlbold v\isactrlsub {\isadigit{0}}\ {\isacharparenleft}{\kern0pt}{\isasymupharpoonleft}{\isacharparenleft}{\kern0pt}\isactrlbold n{\isasymlbrace}\isavar{{\isasympsi}}{\isasymrbrace}{\isacharparenright}{\kern0pt}{\isacharparenright}{\kern0pt}{\isacharparenright}{\kern0pt}{\isacharparenright}{\kern0pt}\ \isactrlbold {\isasymlongrightarrow}\ \isavar{{\isasympsi}}{\isasymrbrace}{\isacharparenright}{\kern0pt}{\isacharparenright}{\kern0pt}{\isacharparenright}{\kern0pt}{\isacharparenright}{\kern0pt}\ {\isasymtriangleright}\ {\isasymemptyset}\ {\isasymturnstile}\ \isactrlbold {\isasymnot}\isactrlbold A\ \isactrlbold {\isasymcdot}\ {\isacharparenleft}{\kern0pt}\isactrlbold {\isasymlambda}{\isacharparenleft}{\kern0pt}\isactrlbold {\isasymnot}\isaconst{iChk}\ \isactrlbold v\isactrlsub {\isadigit{0}}\ {\isacharparenleft}{\kern0pt}{\isasymupharpoonleft}{\isacharparenleft}{\kern0pt}\isactrlbold n{\isasymlbrace}\isavar{{\isasympsi}}{\isasymrbrace}{\isacharparenright}{\kern0pt}{\isacharparenright}{\kern0pt}{\isacharparenright}{\kern0pt}{\isacharparenright}{\kern0pt}}%
}%
\SNIP{Goedel-Paper-Statements:thm:reflect-ent-iff-con-ctx-0}{%
\isakeywordONE{thm}\isamarkupfalse%
\ reflect{\isacharunderscore}{\kern0pt}ent{\isacharunderscore}{\kern0pt}iff{\isacharunderscore}{\kern0pt}con{\isacharunderscore}{\kern0pt}ctx%
\isanewline%
}%
\SNIP{Goedel-Paper-Statements:thm:reflect-ent-iff-con-ctx-1}{%
\isa{{\isacharparenleft}{\kern0pt}\isactrlbold {\isasymnot}\isactrlbold A\ \isactrlbold {\isasymcdot}\ {\isacharparenleft}{\kern0pt}\isactrlbold {\isasymlambda}{\isacharparenleft}{\kern0pt}\isactrlbold {\isasymnot}\isaconst{iChk}\ \isactrlbold v\isactrlsub {\isadigit{0}}\ {\isacharparenleft}{\kern0pt}{\isasymupharpoonleft}\isaconst{cEq{\isadigit{0}}{\isadigit{1}}}{\isacharparenright}{\kern0pt}{\isacharparenright}{\kern0pt}{\isacharparenright}{\kern0pt}\ {\isasymtriangleright}\ \isavar{{\isasymGamma}}\ {\isasymturnstile}\ \isactrlbold {\isadigit{0}}\ \isactrlbold {\isacharequal}{\kern0pt}\ \isactrlbold n{\isadigit{1}}{\isacharparenright}{\kern0pt}\ {\isacharequal}{\kern0pt}\ {\isacharparenleft}{\kern0pt}\isavar{{\isasymGamma}}\ {\isasymturnstile}\ \isaconst{ConRGA}{\isacharparenright}{\kern0pt}}%
}%
\SNIP{Goedel-Paper-Statements:thm:goedel-diag-0}{%
\isakeywordONE{thm}\isamarkupfalse%
\ goedel{\isacharunderscore}{\kern0pt}diag%
\isanewline%
}%
\SNIP{Goedel-Paper-Statements:thm:goedel-diag-1}{%
\isa{{\isacharbrackleft}{\kern0pt}{\isadigit{0}}{\isasymZpinj}\isactrlbold n{\isasymlbrace}\isaconst{goedelD}{\isasymrbrace}{\isacharbrackright}{\kern0pt}\isaconst{goedelD}\ {\isacharequal}{\kern0pt}\ \isaconst{goedelG}}%
}%
\SNIP{Goedel-Paper-Statements:thm:goedelS-eq-0}{%
\isakeywordONE{thm}\isamarkupfalse%
\ goedelS{\isacharunderscore}{\kern0pt}eq%
\isanewline%
}%
\SNIP{Goedel-Paper-Statements:thm:goedelS-eq-1}{%
\isa{{\isasymemptyset}\ {\isasymturnstile}\ \isaconst{goedelS}\ \isactrlbold {\isacharequal}{\kern0pt}\ \isactrlbold n{\isasymlbrace}\isaconst{goedelG}{\isasymrbrace}}%
}%
\SNIP{Goedel-Paper-Statements:thm:goedelE-chk-0}{%
\isakeywordONE{thm}\isamarkupfalse%
\ goedelE{\isacharunderscore}{\kern0pt}chk%
\isanewline%
}%
\SNIP{Goedel-Paper-Statements:thm:goedelE-chk-1}{%
\isa{{\isasymlbrakk}\isavar{{\isasymGamma}}\ {\isasymturnstile}\ \isavar{a}\ \isactrlbold N{\isacharsemicolon}{\kern0pt}\ \isavar{{\isasymGamma}}\ {\isasymturnstile}\ \isavar{b}\ \isactrlbold N{\isacharsemicolon}{\kern0pt}\ \isavar{{\isasymGamma}}\ {\isasymturnstile}\ \isaconst{iChk}\ \isavar{a}\ {\isacharparenleft}{\kern0pt}\isactrlbold n{\isasymlbrace}\isavar{g}{\isasymrbrace}{\isacharparenright}{\kern0pt}{\isacharsemicolon}{\kern0pt}\ \isavar{{\isasymGamma}}\ {\isasymturnstile}\ \isaconst{iChk}\ \isavar{b}\ {\isacharparenleft}{\kern0pt}\isactrlbold n{\isasymlbrace}\isactrlbold {\isasymnot}\isavar{g}{\isasymrbrace}{\isacharparenright}{\kern0pt}{\isasymrbrakk}\ {\isasymLongrightarrow}\ \isavar{{\isasymGamma}}\ {\isasymturnstile}\ \isaconst{iChk}\ {\isacharparenleft}{\kern0pt}\isaconst{goedelE}\ \isactrlbold {\isasymcdot}\ \isavar{a}\ \isactrlbold {\isasymcdot}\ \isavar{b}{\isacharparenright}{\kern0pt}\ \isaconst{cEq{\isadigit{0}}{\isadigit{1}}}}%
}%
\SNIP{Goedel-Paper-Statements:thm:EXPL-goedelE-0}{%
\isakeywordONE{thm}\isamarkupfalse%
\ EXPL{\isacharunderscore}{\kern0pt}goedelE%
\isanewline%
}%
\SNIP{Goedel-Paper-Statements:thm:EXPL-goedelE-1}{%
\isa{\isaconst{EXPL}\ \isaconst{goedelE}}%
}%
\SNIP{Goedel-Paper-Statements:thm:D2-decided-0}{%
\isakeywordONE{thm}\isamarkupfalse%
\ D{\isadigit{2}}{\isacharunderscore}{\kern0pt}decided%
\isanewline%
}%
\SNIP{Goedel-Paper-Statements:thm:D2-decided-1}{%
\isa{{\isasymlbrakk}\isaconst{wft}\ \isavar{{\isasymphi}}{\isacharsemicolon}{\kern0pt}\ \isaconst{wft}\ \isavar{{\isasympsi}}{\isasymrbrakk}\ {\isasymLongrightarrow}\ {\isasymemptyset}\ {\isasymturnstile}\ \isaconst{D{\isadigit{2}}F}\ {\isacharparenleft}{\kern0pt}\isaconst{d{\isadigit{2}}M}\ \isavar{{\isasymphi}}\ \isavar{{\isasympsi}}{\isacharparenright}{\kern0pt}\ \isavar{{\isasymphi}}\ \isavar{{\isasympsi}}}%
}%
\SNIP{Goedel-Paper-Statements:thm:d2M-chk-0}{%
\isakeywordONE{thm}\isamarkupfalse%
\ d{\isadigit{2}}M{\isacharunderscore}{\kern0pt}chk%
\isanewline%
}%
\SNIP{Goedel-Paper-Statements:thm:d2M-chk-1}{%
\isa{{\isasymlbrakk}\isaconst{wft}\ \isavar{{\isasymphi}}{\isacharsemicolon}{\kern0pt}\ \isaconst{wft}\ \isavar{{\isasympsi}}{\isacharsemicolon}{\kern0pt}\ \isavar{{\isasymGamma}}\ {\isasymturnstile}\ \isavar{a}\ \isactrlbold N{\isacharsemicolon}{\kern0pt}\ \isavar{{\isasymGamma}}\ {\isasymturnstile}\ \isavar{b}\ \isactrlbold N{\isacharsemicolon}{\kern0pt}\ \isavar{{\isasymGamma}}\ {\isasymturnstile}\ \isaconst{iChk}\ \isavar{a}\ {\isacharparenleft}{\kern0pt}\isactrlbold n{\isasymlbrace}\isavar{{\isasymphi}}\ \isactrlbold {\isasymlongrightarrow}\ \isavar{{\isasympsi}}{\isasymrbrace}{\isacharparenright}{\kern0pt}{\isacharsemicolon}{\kern0pt}\ \isavar{{\isasymGamma}}\ {\isasymturnstile}\ \isaconst{iChk}\ \isavar{b}\ {\isacharparenleft}{\kern0pt}\isactrlbold n{\isasymlbrace}\isavar{{\isasymphi}}{\isasymrbrace}{\isacharparenright}{\kern0pt}{\isasymrbrakk}\ {\isasymLongrightarrow}\ \isavar{{\isasymGamma}}\ {\isasymturnstile}\ \isaconst{iChk}\ {\isacharparenleft}{\kern0pt}\isaconst{d{\isadigit{2}}M}\ \isavar{{\isasymphi}}\ \isavar{{\isasympsi}}\ \isactrlbold {\isasymcdot}\ \isavar{a}\ \isactrlbold {\isasymcdot}\ \isavar{b}{\isacharparenright}{\kern0pt}\ {\isacharparenleft}{\kern0pt}\isactrlbold n{\isasymlbrace}\isavar{{\isasympsi}}{\isasymrbrace}{\isacharparenright}{\kern0pt}}%
}%
\SNIP{Goedel-Paper-Statements:thm:D3G-goedelSigma-0}{%
\isakeywordONE{thm}\isamarkupfalse%
\ D{\isadigit{3}}G{\isacharunderscore}{\kern0pt}goedelSigma%
\isanewline%
}%
\SNIP{Goedel-Paper-Statements:thm:D3G-goedelSigma-1}{%
\isa{\isaconst{D{\isadigit{3}}G}\ {\isacharparenleft}{\kern0pt}\isaconst{goedelSigmaOf}\ \isaconst{premCertIdx}{\isacharparenright}{\kern0pt}}%
}%
\SNIP{Goedel-Paper-Statements:thm:goedelSigma-is-pr-0}{%
\isakeywordONE{thm}\isamarkupfalse%
\ goedelSigma{\isacharunderscore}{\kern0pt}is{\isacharunderscore}{\kern0pt}pr%
\isanewline%
}%
\SNIP{Goedel-Paper-Statements:thm:goedelSigma-is-pr-1}{%
\isa{\isaconst{goedelSigmaOf}\ \isaconst{premCertIdx}\ {\isacharequal}{\kern0pt}\ \isactrlbold {\isasymlambda}\isaconst{cmp}\ {\isacharparenleft}{\kern0pt}\isaconst{goedelSigmaIdx}\ \isaconst{premCertIdx}{\isacharparenright}{\kern0pt}\ \isactrlbold {\isasymcdot}\ \isactrlbold v\isactrlsub {\isadigit{0}}\ \isactrlbold {\isasymcdot}\ \isactrlbold v\isactrlsub {\isadigit{0}}\ {\isasymand}\ \isaconst{n{\isadigit{2}}pr}\ {\isacharparenleft}{\kern0pt}\isaconst{pr{\isadigit{2}}fi}\ {\isacharparenleft}{\kern0pt}\isaconst{goedelSigmaIdx}\ \isaconst{premCertIdx}{\isacharparenright}{\kern0pt}{\isacharparenright}{\kern0pt}}%
}%
\SNIP{Goedel-Paper-Statements:thm:D3-decided-any-0}{%
\isakeywordONE{thm}\isamarkupfalse%
\ D{\isadigit{3}}{\isacharunderscore}{\kern0pt}decided{\isacharunderscore}{\kern0pt}any%
\isanewline%
}%
\SNIP{Goedel-Paper-Statements:thm:D3-decided-any-1}{%
\isa{{\isasymemptyset}\ {\isasymturnstile}\ \isactrlbold A\ \isactrlbold {\isasymcdot}\ {\isacharparenleft}{\kern0pt}\isactrlbold {\isasymlambda}{\isacharparenleft}{\kern0pt}\isaconst{iChk}\ \isactrlbold v\isactrlsub {\isadigit{0}}\ {\isacharparenleft}{\kern0pt}\isactrlbold n{\isasymlbrace}\isavar{{\isasymphi}}{\isasymrbrace}{\isacharparenright}{\kern0pt}\ \isactrlbold {\isasymlongrightarrow}\ \isaconst{iChk}\ {\isacharparenleft}{\kern0pt}\isaconst{d{\isadigit{3}}SigmaOf}\ \isavar{{\isasymphi}}\ \isactrlbold {\isasymcdot}\ \isactrlbold v\isactrlsub {\isadigit{0}}{\isacharparenright}{\kern0pt}\ {\isacharparenleft}{\kern0pt}\isactrlbold n{\isasymlbrace}\isactrlbold {\isasymnot}\isactrlbold A\ \isactrlbold {\isasymcdot}\ {\isacharparenleft}{\kern0pt}\isactrlbold {\isasymlambda}{\isacharparenleft}{\kern0pt}\isactrlbold {\isasymnot}\isaconst{iChk}\ \isactrlbold v\isactrlsub {\isadigit{0}}\ {\isacharparenleft}{\kern0pt}\isactrlbold n{\isasymlbrace}\isavar{{\isasymphi}}{\isasymrbrace}{\isacharparenright}{\kern0pt}{\isacharparenright}{\kern0pt}{\isacharparenright}{\kern0pt}{\isasymrbrace}{\isacharparenright}{\kern0pt}{\isacharparenright}{\kern0pt}{\isacharparenright}{\kern0pt}}%
}%
\SNIP{Goedel-Paper-Statements:thm:D3-at-0}{%
\isakeywordONE{thm}\isamarkupfalse%
\ D{\isadigit{3}}{\isacharunderscore}{\kern0pt}at%
\isanewline%
}%
\SNIP{Goedel-Paper-Statements:thm:D3-at-1}{%
\isa{{\isasymlbrakk}\isaconst{d{\isadigit{3}}Src}\ \isavar{S}\ \isavar{SI}\ \isavar{sc}{\isacharsemicolon}{\kern0pt}\ \isaconst{Syntax{\isacharunderscore}{\kern0pt}Term{\isachardot}{\kern0pt}closed}\ \isavar{T}{\isacharsemicolon}{\kern0pt}\ {\isasymemptyset}\ {\isasymturnstile}\ \isavar{T}\ \isactrlbold N{\isacharsemicolon}{\kern0pt}\ {\isasymemptyset}\ {\isasymturnstile}\ \isavar{T}\ \isactrlbold {\isacharequal}{\kern0pt}\ \isavar{S}{\isasymrbrakk}\ {\isasymLongrightarrow}\ {\isasymemptyset}\ {\isasymturnstile}\ \isactrlbold A\ \isactrlbold {\isasymcdot}\ {\isacharparenleft}{\kern0pt}\isactrlbold {\isasymlambda}{\isacharparenleft}{\kern0pt}\isaconst{iChk}\ \isactrlbold v\isactrlsub {\isadigit{0}}\ \isavar{T}\ \isactrlbold {\isasymlongrightarrow}\ \isaconst{iChk}\ {\isacharparenleft}{\kern0pt}\isaconst{d{\isadigit{3}}SigmaAt}\ \isavar{S}\ \isavar{SI}\ \isavar{sc}\ \isactrlbold {\isasymcdot}\ \isactrlbold v\isactrlsub {\isadigit{0}}{\isacharparenright}{\kern0pt}\ {\isacharparenleft}{\kern0pt}\isactrlbold n{\isasymlbrace}\isactrlbold {\isasymnot}\isactrlbold A\ \isactrlbold {\isasymcdot}\ \isaconst{prvF}\ \isavar{S}{\isasymrbrace}{\isacharparenright}{\kern0pt}{\isacharparenright}{\kern0pt}{\isacharparenright}{\kern0pt}}%
}%
\SNIP{Goedel-Paper-Statements:thm:D3-uniform-0}{%
\isakeywordONE{thm}\isamarkupfalse%
\ D{\isadigit{3}}{\isacharunderscore}{\kern0pt}uniform%
\isanewline%
}%
\SNIP{Goedel-Paper-Statements:thm:D3-uniform-1}{%
\isa{{\isasymemptyset}\ {\isasymturnstile}\ \isactrlbold A\ \isactrlbold {\isasymcdot}\ {\isacharparenleft}{\kern0pt}\isactrlbold {\isasymlambda}\isactrlbold A\ \isactrlbold {\isasymcdot}\ {\isacharparenleft}{\kern0pt}\isactrlbold {\isasymlambda}{\isacharparenleft}{\kern0pt}\isaconst{iChk}\ \isactrlbold v\isactrlsub {\isadigit{0}}\ \isactrlbold v\isactrlsub {\isadigit{1}}\ \isactrlbold {\isasymlongrightarrow}\ \isaconst{iChk}\ {\isacharparenleft}{\kern0pt}\isaconst{cmp}\ \isaconst{uSigIdx}\ \isactrlbold {\isasymcdot}\ \isactrlbold v\isactrlsub {\isadigit{1}}\ \isactrlbold {\isasymcdot}\ \isactrlbold v\isactrlsub {\isadigit{0}}{\isacharparenright}{\kern0pt}\ {\isacharparenleft}{\kern0pt}\isaconst{prCU}\ \isactrlbold v\isactrlsub {\isadigit{1}}{\isacharparenright}{\kern0pt}{\isacharparenright}{\kern0pt}{\isacharparenright}{\kern0pt}{\isacharparenright}{\kern0pt}}%
}%
\SNIP{Goedel-Paper-Statements:thm:D3-uniform-num-0}{%
\isakeywordONE{thm}\isamarkupfalse%
\ D{\isadigit{3}}{\isacharunderscore}{\kern0pt}uniform{\isacharunderscore}{\kern0pt}num%
\isanewline%
}%
\SNIP{Goedel-Paper-Statements:thm:D3-uniform-num-1}{%
\isa{{\isasymemptyset}\ {\isasymturnstile}\ \isactrlbold A\ \isactrlbold {\isasymcdot}\ {\isacharparenleft}{\kern0pt}\isactrlbold {\isasymlambda}{\isacharparenleft}{\kern0pt}\isaconst{iChk}\ \isactrlbold v\isactrlsub {\isadigit{0}}\ {\isacharparenleft}{\kern0pt}\isactrlbold n\isavar{k}{\isacharparenright}{\kern0pt}\ \isactrlbold {\isasymlongrightarrow}\ \isaconst{iChk}\ {\isacharparenleft}{\kern0pt}\isaconst{cmp}\ \isaconst{uSigIdx}\ \isactrlbold {\isasymcdot}\ {\isacharparenleft}{\kern0pt}\isactrlbold n\isavar{k}{\isacharparenright}{\kern0pt}\ \isactrlbold {\isasymcdot}\ \isactrlbold v\isactrlsub {\isadigit{0}}{\isacharparenright}{\kern0pt}\ {\isacharparenleft}{\kern0pt}\isactrlbold n{\isasymlbrace}\isactrlbold {\isasymnot}\isactrlbold A\ \isactrlbold {\isasymcdot}\ \isaconst{prvF}\ {\isacharparenleft}{\kern0pt}\isactrlbold n\isavar{k}{\isacharparenright}{\kern0pt}{\isasymrbrace}{\isacharparenright}{\kern0pt}{\isacharparenright}{\kern0pt}{\isacharparenright}{\kern0pt}}%
}%
\SNIP{Goedel-Paper-Statements:thm:goedelG-of-expl-d3g-ctx-0}{%
\isakeywordONE{thm}\isamarkupfalse%
\ goedelG{\isacharunderscore}{\kern0pt}of{\isacharunderscore}{\kern0pt}expl{\isacharunderscore}{\kern0pt}d{\isadigit{3}}g{\isacharunderscore}{\kern0pt}ctx%
\isanewline%
}%
\SNIP{Goedel-Paper-Statements:thm:goedelG-of-expl-d3g-ctx-1}{%
\isa{{\isasymlbrakk}\isaconst{EXPL}\ \isavar{e}{\isacharsemicolon}{\kern0pt}\ \isaconst{D{\isadigit{3}}G}\ \isavar{{\isasymsigma}}{\isacharsemicolon}{\kern0pt}\ \isavar{{\isasymGamma}}\ {\isasymturnstile}\ \isaconst{ConRGA}{\isasymrbrakk}\ {\isasymLongrightarrow}\ \isavar{{\isasymGamma}}\ {\isasymturnstile}\ \isaconst{goedelG}}%
}%
\SNIP{Goedel-Paper-Statements:thm:rga-goedel1-0}{%
\isakeywordONE{thm}\isamarkupfalse%
\ rga{\isacharunderscore}{\kern0pt}goedel{\isadigit{1}}%
\isanewline%
}%
\SNIP{Goedel-Paper-Statements:thm:rga-goedel1-1}{%
\isa{{\isasymnot}\ {\isasymemptyset}\ {\isasymturnstile}\ \isaconst{goedelG}}%
}%
\SNIP{Goedel-Paper-Statements:thm:goedelG-refutes-0}{%
\isakeywordONE{thm}\isamarkupfalse%
\ goedelG{\isacharunderscore}{\kern0pt}refutes%
\isanewline%
}%
\SNIP{Goedel-Paper-Statements:thm:goedelG-refutes-1}{%
\isa{{\isasymemptyset}\ {\isasymturnstile}\ \isaconst{goedelG}\ {\isasymLongrightarrow}\ {\isasymemptyset}\ {\isasymturnstile}\ \isactrlbold {\isadigit{0}}\ \isactrlbold {\isacharequal}{\kern0pt}\ \isactrlbold n{\isadigit{1}}}%
}%
\SNIP{Goedel-Paper-Statements:thm:rga-goedel2-0}{%
\isakeywordONE{thm}\isamarkupfalse%
\ rga{\isacharunderscore}{\kern0pt}goedel{\isadigit{2}}%
\isanewline%
}%
\SNIP{Goedel-Paper-Statements:thm:rga-goedel2-1}{%
\isa{{\isasymnot}\ {\isasymemptyset}\ {\isasymturnstile}\ \isaconst{ConRGA}}%
}%
\SNIP{Goedel-Paper-Statements:thm:cut1M-chk-0}{%
\isakeywordONE{thm}\isamarkupfalse%
\ cut{\isadigit{1}}M{\isacharunderscore}{\kern0pt}chk%
\isanewline%
}%
\SNIP{Goedel-Paper-Statements:thm:cut1M-chk-1}{%
\isa{{\isasymlbrakk}\isavar{Q}\ {\isasymP}{\isasymlbrakk}\isavar{A}\ {\isasymtriangleright}\ {\isasymemptyset}\ {\isasymturnstile}\ \isavar{B}{\isasymrbrakk}{\isacharsemicolon}{\kern0pt}\ \isavar{{\isasymGamma}}\ {\isasymturnstile}\ \isavar{a}\ \isactrlbold N{\isacharsemicolon}{\kern0pt}\ \isavar{{\isasymGamma}}\ {\isasymturnstile}\ \isaconst{iChk}\ \isavar{a}\ {\isacharparenleft}{\kern0pt}\isactrlbold n{\isasymlbrace}\isavar{A}{\isasymrbrace}{\isacharparenright}{\kern0pt}{\isasymrbrakk}\ {\isasymLongrightarrow}\ \isavar{{\isasymGamma}}\ {\isasymturnstile}\ \isaconst{iChk}\ {\isacharparenleft}{\kern0pt}\isaconst{cut{\isadigit{1}}M}\ \isavar{Q}\ \isavar{B}\ \isactrlbold {\isasymcdot}\ \isavar{a}\ \isactrlbold {\isasymcdot}\ \isavar{a}{\isacharparenright}{\kern0pt}\ {\isacharparenleft}{\kern0pt}\isactrlbold n{\isasymlbrace}\isavar{B}{\isasymrbrace}{\isacharparenright}{\kern0pt}}%
}%
\SNIP{Goedel-Paper-Statements:thm:ent-nec-0}{%
\isakeywordONE{thm}\isamarkupfalse%
\ ent{\isacharunderscore}{\kern0pt}nec%
\isanewline%
}%
\SNIP{Goedel-Paper-Statements:thm:ent-nec-1}{%
\isa{\isavar{A}\ {\isasymtriangleright}\ {\isasymemptyset}\ {\isasymturnstile}\ \isavar{B}\ {\isasymLongrightarrow}\ \isactrlbold {\isasymnot}\isactrlbold A\ \isactrlbold {\isasymcdot}\ {\isacharparenleft}{\kern0pt}\isactrlbold {\isasymlambda}{\isacharparenleft}{\kern0pt}\isactrlbold {\isasymnot}\isaconst{iChk}\ \isactrlbold v\isactrlsub {\isadigit{0}}\ {\isacharparenleft}{\kern0pt}{\isasymupharpoonleft}{\isacharparenleft}{\kern0pt}\isactrlbold n{\isasymlbrace}\isavar{A}{\isasymrbrace}{\isacharparenright}{\kern0pt}{\isacharparenright}{\kern0pt}{\isacharparenright}{\kern0pt}{\isacharparenright}{\kern0pt}\ {\isasymtriangleright}\ \isavar{{\isasymGamma}}\ {\isasymturnstile}\ \isactrlbold {\isasymnot}\isactrlbold A\ \isactrlbold {\isasymcdot}\ {\isacharparenleft}{\kern0pt}\isactrlbold {\isasymlambda}{\isacharparenleft}{\kern0pt}\isactrlbold {\isasymnot}\isaconst{iChk}\ \isactrlbold v\isactrlsub {\isadigit{0}}\ {\isacharparenleft}{\kern0pt}{\isasymupharpoonleft}{\isacharparenleft}{\kern0pt}\isactrlbold n{\isasymlbrace}\isavar{B}{\isasymrbrace}{\isacharparenright}{\kern0pt}{\isacharparenright}{\kern0pt}{\isacharparenright}{\kern0pt}{\isacharparenright}{\kern0pt}}%
}%
\SNIP{Goedel-Paper-Statements:thm:rga-goedel2-internal-contra-0}{%
\isakeywordONE{thm}\isamarkupfalse%
\ rga{\isacharunderscore}{\kern0pt}goedel{\isadigit{2}}{\isacharunderscore}{\kern0pt}internal{\isacharunderscore}{\kern0pt}contra%
\isanewline%
}%
\SNIP{Goedel-Paper-Statements:thm:rga-goedel2-internal-contra-1}{%
\isa{\isactrlbold {\isasymnot}\isactrlbold A\ \isactrlbold {\isasymcdot}\ {\isacharparenleft}{\kern0pt}\isactrlbold {\isasymlambda}{\isacharparenleft}{\kern0pt}\isactrlbold {\isasymnot}\isaconst{iChk}\ \isactrlbold v\isactrlsub {\isadigit{0}}\ {\isacharparenleft}{\kern0pt}{\isasymupharpoonleft}{\isacharparenleft}{\kern0pt}\isactrlbold n{\isasymlbrace}\isaconst{ConRGA}{\isasymrbrace}{\isacharparenright}{\kern0pt}{\isacharparenright}{\kern0pt}{\isacharparenright}{\kern0pt}{\isacharparenright}{\kern0pt}\ {\isasymtriangleright}\ {\isasymemptyset}\ {\isasymturnstile}\ \isactrlbold {\isasymnot}\isaconst{ConRGA}}%
}%
\SNIP{Goedel-Paper-Statements:thm:rga-goedel2-internal-0}{%
\isakeywordONE{thm}\isamarkupfalse%
\ rga{\isacharunderscore}{\kern0pt}goedel{\isadigit{2}}{\isacharunderscore}{\kern0pt}internal%
\isanewline%
}%
\SNIP{Goedel-Paper-Statements:thm:rga-goedel2-internal-1}{%
\isa{\isaconst{ConRGA}\ {\isasymtriangleright}\ {\isasymemptyset}\ {\isasymturnstile}\ \isactrlbold {\isasymnot}\isactrlbold {\isasymnot}\isactrlbold A\ \isactrlbold {\isasymcdot}\ {\isacharparenleft}{\kern0pt}\isactrlbold {\isasymlambda}{\isacharparenleft}{\kern0pt}\isactrlbold {\isasymnot}\isaconst{iChk}\ \isactrlbold v\isactrlsub {\isadigit{0}}\ {\isacharparenleft}{\kern0pt}{\isasymupharpoonleft}{\isacharparenleft}{\kern0pt}\isactrlbold n{\isasymlbrace}\isaconst{ConRGA}{\isasymrbrace}{\isacharparenright}{\kern0pt}{\isacharparenright}{\kern0pt}{\isacharparenright}{\kern0pt}{\isacharparenright}{\kern0pt}}%
}%
\SNIP{Goedel-Paper-Statements:thm:rga-goedel2-sentence-0}{%
\isakeywordONE{thm}\isamarkupfalse%
\ rga{\isacharunderscore}{\kern0pt}goedel{\isadigit{2}}{\isacharunderscore}{\kern0pt}sentence%
\isanewline%
}%
\SNIP{Goedel-Paper-Statements:thm:rga-goedel2-sentence-1}{%
\isa{{\isasymemptyset}\ {\isasymturnstile}\ \isactrlbold A\ \isactrlbold {\isasymcdot}\ {\isacharparenleft}{\kern0pt}\isactrlbold {\isasymlambda}{\isacharparenleft}{\kern0pt}\isaconst{iChk}\ \isactrlbold v\isactrlsub {\isadigit{0}}\ {\isacharparenleft}{\kern0pt}\isactrlbold n{\isasymlbrace}\isaconst{ConRGA}{\isasymrbrace}{\isacharparenright}{\kern0pt}\ \isactrlbold {\isasymlongrightarrow}\ \isactrlbold {\isasymnot}\isaconst{ConRGA}{\isacharparenright}{\kern0pt}{\isacharparenright}{\kern0pt}}%
}%
\SNIP{Goedel-Paper-Statements:thm:goedel2T-sentence-0}{%
\isakeywordONE{thm}\isamarkupfalse%
\ goedel{\isadigit{2}}T{\isacharunderscore}{\kern0pt}sentence%
\isanewline%
}%
\SNIP{Goedel-Paper-Statements:thm:goedel2T-sentence-1}{%
\isa{\isavar{Q}\ {\isasymP}{\isasymlbrakk}\isaconst{ConRGA}\ {\isasymtriangleright}\ {\isasymemptyset}\ {\isasymturnstile}\ \isaconst{goedelG}{\isasymrbrakk}\ {\isasymLongrightarrow}\ {\isasymemptyset}\ {\isasymturnstile}\ \isactrlbold A\ \isactrlbold {\isasymcdot}\ {\isacharparenleft}{\kern0pt}\isactrlbold {\isasymlambda}{\isacharparenleft}{\kern0pt}\isaconst{iChk}\ \isactrlbold v\isactrlsub {\isadigit{0}}\ {\isacharparenleft}{\kern0pt}\isactrlbold n{\isasymlbrace}\isaconst{ConRGA}{\isasymrbrace}{\isacharparenright}{\kern0pt}\ \isactrlbold {\isasymlongrightarrow}\ \isaconst{iChk}\ {\isacharparenleft}{\kern0pt}\isaconst{goedel{\isadigit{2}}T}\ \isavar{Q}\ \isactrlbold {\isasymcdot}\ \isactrlbold v\isactrlsub {\isadigit{0}}{\isacharparenright}{\kern0pt}\ \isaconst{cEq{\isadigit{0}}{\isadigit{1}}}{\isacharparenright}{\kern0pt}{\isacharparenright}{\kern0pt}}%
}%
\SNIP{Goedel-Paper-Statements:thm:goedelG-con-interderivable-0}{%
\isakeywordONE{thm}\isamarkupfalse%
\ goedelG{\isacharunderscore}{\kern0pt}con{\isacharunderscore}{\kern0pt}interderivable%
\isanewline%
}%
\SNIP{Goedel-Paper-Statements:thm:goedelG-con-interderivable-1}{%
\isa{\isaconst{goedelG}\ {\isasymtriangleright}\ {\isasymemptyset}\ {\isasymturnstile}\ \isaconst{ConRGA}\ {\isasymand}\ \isaconst{ConRGA}\ {\isasymtriangleright}\ {\isasymemptyset}\ {\isasymturnstile}\ \isaconst{goedelG}}%
}%
\SNIP{Goedel-Paper-Statements:thm:notcon-prc-con-interderivable-0}{%
\isakeywordONE{thm}\isamarkupfalse%
\ notcon{\isacharunderscore}{\kern0pt}prc{\isacharunderscore}{\kern0pt}con{\isacharunderscore}{\kern0pt}interderivable%
\isanewline%
}%
\SNIP{Goedel-Paper-Statements:thm:notcon-prc-con-interderivable-1}{%
\isa{\isactrlbold {\isasymnot}\isaconst{ConRGA}\ {\isasymtriangleright}\ {\isasymemptyset}\ {\isasymturnstile}\ \isactrlbold {\isasymnot}\isactrlbold A\ \isactrlbold {\isasymcdot}\ {\isacharparenleft}{\kern0pt}\isactrlbold {\isasymlambda}{\isacharparenleft}{\kern0pt}\isactrlbold {\isasymnot}\isaconst{iChk}\ \isactrlbold v\isactrlsub {\isadigit{0}}\ {\isacharparenleft}{\kern0pt}{\isasymupharpoonleft}{\isacharparenleft}{\kern0pt}\isactrlbold n{\isasymlbrace}\isaconst{ConRGA}{\isasymrbrace}{\isacharparenright}{\kern0pt}{\isacharparenright}{\kern0pt}{\isacharparenright}{\kern0pt}{\isacharparenright}{\kern0pt}\ {\isasymand}\ \isactrlbold {\isasymnot}\isactrlbold A\ \isactrlbold {\isasymcdot}\ {\isacharparenleft}{\kern0pt}\isactrlbold {\isasymlambda}{\isacharparenleft}{\kern0pt}\isactrlbold {\isasymnot}\isaconst{iChk}\ \isactrlbold v\isactrlsub {\isadigit{0}}\ {\isacharparenleft}{\kern0pt}{\isasymupharpoonleft}{\isacharparenleft}{\kern0pt}\isactrlbold n{\isasymlbrace}\isaconst{ConRGA}{\isasymrbrace}{\isacharparenright}{\kern0pt}{\isacharparenright}{\kern0pt}{\isacharparenright}{\kern0pt}{\isacharparenright}{\kern0pt}\ {\isasymtriangleright}\ {\isasymemptyset}\ {\isasymturnstile}\ \isactrlbold {\isasymnot}\isaconst{ConRGA}}%
}%
\SNIP{Goedel-Paper-Statements:thm:gcalc-strip-rga-0}{%
\isakeywordONE{thm}\isamarkupfalse%
\ gcalc{\isacharunderscore}{\kern0pt}strip{\isacharunderscore}{\kern0pt}rga%
\isanewline%
}%
\SNIP{Goedel-Paper-Statements:thm:gcalc-strip-rga-1}{%
\isa{\isaconst{gcalc{\isacharunderscore}{\kern0pt}strip}\ {\isacharparenleft}{\kern0pt}{\isasymturnstile}{\isacharparenright}{\kern0pt}}%
}%
\SNIP{Goedel-Paper-Statements:thm:g2-con-rga-0}{%
\isakeywordONE{thm}\isamarkupfalse%
\ g{\isadigit{2}}{\isacharunderscore}{\kern0pt}con{\isacharunderscore}{\kern0pt}rga%
\isanewline%
}%
\SNIP{Goedel-Paper-Statements:thm:g2-con-rga-1}{%
\isa{\isaconst{g{\isadigit{2}}{\isacharunderscore}{\kern0pt}con}\ {\isacharparenleft}{\kern0pt}{\isasymturnstile}{\isacharparenright}{\kern0pt}\ \isaconst{iChk}}%
}%
\SNIP{Goedel-Paper-Statements:thm:g2-exf-con-rga-0}{%
\isakeywordONE{thm}\isamarkupfalse%
\ g{\isadigit{2}}{\isacharunderscore}{\kern0pt}exf{\isacharunderscore}{\kern0pt}con{\isacharunderscore}{\kern0pt}rga%
\isanewline%
}%
\SNIP{Goedel-Paper-Statements:thm:g2-exf-con-rga-1}{%
\isa{\isaconst{g{\isadigit{2}}{\isacharunderscore}{\kern0pt}exf{\isacharunderscore}{\kern0pt}con}\ {\isacharparenleft}{\kern0pt}{\isasymturnstile}{\isacharparenright}{\kern0pt}\ \isaconst{iChk}}%
}%
\SNIP{Goedel-Paper-Statements:thm:rga-goedel2-textbook-0}{%
\isakeywordONE{thm}\isamarkupfalse%
\ rga{\isacharunderscore}{\kern0pt}goedel{\isadigit{2}}{\isacharunderscore}{\kern0pt}textbook%
\isanewline%
}%
\SNIP{Goedel-Paper-Statements:thm:rga-goedel2-textbook-1}{%
\isa{{\isasymnot}\ {\isasymemptyset}\ {\isasymturnstile}\ \isaconst{ConRGA}}%
}%
\SNIP{Goedel-Paper-Statements:thm:drvChkIdxG-def-0}{%
\isakeywordONE{thm}\isamarkupfalse%
\ drvChkIdxG{\isacharunderscore}{\kern0pt}def%
\isanewline%
}%
\SNIP{Goedel-Paper-Statements:thm:drvChkIdxG-def-1}{%
\isa{\isaconst{drvChkIdxG}\ \isavar{L}\ {\isasymequiv}\ \isaconst{orNestIdx}\ \isavar{L}}%
}%
\SNIP{Goedel-Paper-Statements:thm:iChkD-def-0}{%
\isakeywordONE{thm}\isamarkupfalse%
\ iChkD{\isacharunderscore}{\kern0pt}def%
\isanewline%
}%
\SNIP{Goedel-Paper-Statements:thm:iChkD-def-1}{%
\isa{\isaconst{iChkD}\ \isavar{Dt}\ \isavar{p}\ \isavar{c}\ {\isasymequiv}\ {\isacharparenleft}{\kern0pt}\isaconst{cmp}\ {\isacharparenleft}{\kern0pt}\isaconst{validChkD}\ \isavar{Dt}{\isacharparenright}{\kern0pt}\ \isactrlbold {\isasymcdot}\ \isavar{p}\ \isactrlbold {\isasymcdot}\ \isavar{p}\ \isactrlbold {\isacharequal}{\kern0pt}\ \isactrlbold n\isaconst{Suc}\ {\isadigit{0}}\ \isactrlbold {\isasymand}\ \isavar{p}\ \isactrlbold {\isasymnoteq}\ \isactrlbold {\isadigit{0}}{\isacharparenright}{\kern0pt}\ \isactrlbold {\isasymand}\ \isaconst{cmp}\ \isaconst{length{\isacharunderscore}{\kern0pt}encIdx}\ \isactrlbold {\isasymcdot}\ \isaconst{cpFst}\ {\isacharparenleft}{\kern0pt}\isaconst{cpFst}\ {\isacharparenleft}{\kern0pt}\isaconst{trmPredR}\ \isavar{p}{\isacharparenright}{\kern0pt}{\isacharparenright}{\kern0pt}\ \isactrlbold {\isasymcdot}\ \isaconst{cpFst}\ {\isacharparenleft}{\kern0pt}\isaconst{cpFst}\ {\isacharparenleft}{\kern0pt}\isaconst{trmPredR}\ \isavar{p}{\isacharparenright}{\kern0pt}{\isacharparenright}{\kern0pt}\ \isactrlbold {\isacharequal}{\kern0pt}\ \isactrlbold {\isadigit{0}}\ \isactrlbold {\isasymand}\ \isaconst{cpSnd}\ {\isacharparenleft}{\kern0pt}\isaconst{cpFst}\ {\isacharparenleft}{\kern0pt}\isaconst{trmPredR}\ \isavar{p}{\isacharparenright}{\kern0pt}{\isacharparenright}{\kern0pt}\ \isactrlbold {\isacharequal}{\kern0pt}\ \isavar{c}}%
}%
\SNIP{Goedel-Paper-Statements:thm:chkG-bridge-drvChkIdx-0}{%
\isakeywordONE{thm}\isamarkupfalse%
\ chkG{\isacharunderscore}{\kern0pt}bridge{\isacharunderscore}{\kern0pt}drvChkIdx%
\isanewline%
}%
\SNIP{Goedel-Paper-Statements:thm:chkG-bridge-drvChkIdx-1}{%
\isa{\isaconst{drvChkIdxG}\ \isaconst{drvChkList}\ {\isacharequal}{\kern0pt}\ \isaconst{drvChkIdx}}%
}%
\SNIP{Goedel-Paper-Statements:thm:chkG-bridge-iChk-0}{%
\isakeywordONE{thm}\isamarkupfalse%
\ chkG{\isacharunderscore}{\kern0pt}bridge{\isacharunderscore}{\kern0pt}iChk%
\isanewline%
}%
\SNIP{Goedel-Paper-Statements:thm:chkG-bridge-iChk-1}{%
\isa{\isaconst{iChkD}\ {\isacharparenleft}{\kern0pt}\isaconst{drvChkIdxG}\ \isaconst{drvChkList}{\isacharparenright}{\kern0pt}\ \isavar{p}\ \isavar{c}\ {\isacharequal}{\kern0pt}\ \isaconst{iChk}\ \isavar{p}\ \isavar{c}}%
}%
\SNIP{Goedel-Paper-Statements:thm:tblShift-def-0}{%
\isakeywordONE{thm}\isamarkupfalse%
\ tblShift{\isacharunderscore}{\kern0pt}def%
\isanewline%
}%
\SNIP{Goedel-Paper-Statements:thm:tblShift-def-1}{%
\isa{\isaconst{tblShift}\ \isavar{E}\ {\isacharequal}{\kern0pt}\ {\isacharparenleft}{\kern0pt}{\isasymforall}\isabound{{\isasymGamma}}\ \isabound{k}\ \isabound{k{\isacharprime}{\kern0pt}}\ \isabound{sh}{\isachardot}{\kern0pt}\ \isabound{{\isasymGamma}}\ {\isasymturnstile}\ \isabound{k}\ \isactrlbold N\ {\isasymlongrightarrow}\ \isabound{{\isasymGamma}}\ {\isasymturnstile}\ \isabound{k{\isacharprime}{\kern0pt}}\ \isactrlbold N\ {\isasymlongrightarrow}\ {\isacharparenleft}{\kern0pt}{\isasymforall}\isabound{{\isasymDelta}}\ \isabound{z}{\isachardot}{\kern0pt}\ {\isacharparenleft}{\kern0pt}{\isasymforall}\isabound{{\isasymphi}}{\isachardot}{\kern0pt}\ \isabound{{\isasymGamma}}\ {\isasymturnstile}\ \isabound{{\isasymphi}}\ {\isasymlongrightarrow}\ \isabound{{\isasymDelta}}\ {\isasymturnstile}\ \isabound{{\isasymphi}}{\isacharparenright}{\kern0pt}\ {\isasymlongrightarrow}\ \isabound{{\isasymDelta}}\ {\isasymturnstile}\ \isabound{z}\ \isactrlbold N\ {\isasymlongrightarrow}\ \isabound{{\isasymDelta}}\ {\isasymturnstile}\ \isabound{sh}\ \isabound{z}\ \isactrlbold N{\isacharparenright}{\kern0pt}\ {\isasymlongrightarrow}\ \isabound{{\isasymGamma}}\ {\isasymturnstile}\ \isaconst{cpFst}\ \isabound{k{\isacharprime}{\kern0pt}}\ \isactrlbold {\isacharequal}{\kern0pt}\ \isaconst{cpFst}\ \isabound{k}\ {\isasymlongrightarrow}\ {\isacharparenleft}{\kern0pt}{\isasymforall}\isabound{{\isasymDelta}}\ \isabound{z}{\isachardot}{\kern0pt}\ {\isacharparenleft}{\kern0pt}{\isasymforall}\isabound{{\isasymphi}}{\isachardot}{\kern0pt}\ \isabound{{\isasymGamma}}\ {\isasymturnstile}\ \isabound{{\isasymphi}}\ {\isasymlongrightarrow}\ \isabound{{\isasymDelta}}\ {\isasymturnstile}\ \isabound{{\isasymphi}}{\isacharparenright}{\kern0pt}\ {\isasymlongrightarrow}\ \isabound{{\isasymDelta}}\ {\isasymturnstile}\ \isabound{z}\ \isactrlbold N\ {\isasymlongrightarrow}\ \isabound{{\isasymDelta}}\ {\isasymturnstile}\ \isaconst{gapSub}\ {\isacharparenleft}{\kern0pt}\isactrlbold S\ \isabound{z}{\isacharparenright}{\kern0pt}\ {\isacharparenleft}{\kern0pt}\isaconst{cmp}\ \isaconst{length{\isacharunderscore}{\kern0pt}encIdx}\ \isactrlbold {\isasymcdot}\ \isaconst{cpSnd}\ \isabound{k}\ \isactrlbold {\isasymcdot}\ \isaconst{cpSnd}\ \isabound{k}{\isacharparenright}{\kern0pt}\ \isactrlbold {\isacharequal}{\kern0pt}\ \isactrlbold {\isadigit{0}}\ {\isasymlongrightarrow}\ \isabound{{\isasymDelta}}\ {\isasymturnstile}\ \isaconst{gapSub}\ {\isacharparenleft}{\kern0pt}\isactrlbold S\ \isabound{sh}\ \isabound{z}{\isacharparenright}{\kern0pt}\ {\isacharparenleft}{\kern0pt}\isaconst{cmp}\ \isaconst{length{\isacharunderscore}{\kern0pt}encIdx}\ \isactrlbold {\isasymcdot}\ \isaconst{cpSnd}\ \isabound{k{\isacharprime}{\kern0pt}}\ \isactrlbold {\isasymcdot}\ \isaconst{cpSnd}\ \isabound{k{\isacharprime}{\kern0pt}}{\isacharparenright}{\kern0pt}\ \isactrlbold {\isacharequal}{\kern0pt}\ \isactrlbold {\isadigit{0}}{\isacharparenright}{\kern0pt}\ {\isasymlongrightarrow}\ {\isacharparenleft}{\kern0pt}{\isasymforall}\isabound{{\isasymDelta}}\ \isabound{z}{\isachardot}{\kern0pt}\ {\isacharparenleft}{\kern0pt}{\isasymforall}\isabound{{\isasymphi}}{\isachardot}{\kern0pt}\ \isabound{{\isasymGamma}}\ {\isasymturnstile}\ \isabound{{\isasymphi}}\ {\isasymlongrightarrow}\ \isabound{{\isasymDelta}}\ {\isasymturnstile}\ \isabound{{\isasymphi}}{\isacharparenright}{\kern0pt}\ {\isasymlongrightarrow}\ \isabound{{\isasymDelta}}\ {\isasymturnstile}\ \isabound{z}\ \isactrlbold N\ {\isasymlongrightarrow}\ \isabound{{\isasymDelta}}\ {\isasymturnstile}\ \isaconst{gapSub}\ {\isacharparenleft}{\kern0pt}\isactrlbold S\ \isabound{z}{\isacharparenright}{\kern0pt}\ {\isacharparenleft}{\kern0pt}\isaconst{cmp}\ \isaconst{length{\isacharunderscore}{\kern0pt}encIdx}\ \isactrlbold {\isasymcdot}\ \isaconst{cpSnd}\ \isabound{k}\ \isactrlbold {\isasymcdot}\ \isaconst{cpSnd}\ \isabound{k}{\isacharparenright}{\kern0pt}\ \isactrlbold {\isacharequal}{\kern0pt}\ \isactrlbold {\isadigit{0}}\ {\isasymlongrightarrow}\ \isabound{{\isasymDelta}}\ {\isasymturnstile}\ \isaconst{trmPredR}\ {\isacharparenleft}{\kern0pt}\isaconst{cmp}\ \isaconst{opt{\isacharunderscore}{\kern0pt}hd{\isacharunderscore}{\kern0pt}encIdx}\ \isactrlbold {\isasymcdot}\ \isaconst{lstDrop}\ {\isacharparenleft}{\kern0pt}\isaconst{cpSnd}\ \isabound{k{\isacharprime}{\kern0pt}}{\isacharparenright}{\kern0pt}\ {\isacharparenleft}{\kern0pt}\isabound{sh}\ \isabound{z}{\isacharparenright}{\kern0pt}\ \isactrlbold {\isasymcdot}\ \isaconst{lstDrop}\ {\isacharparenleft}{\kern0pt}\isaconst{cpSnd}\ \isabound{k{\isacharprime}{\kern0pt}}{\isacharparenright}{\kern0pt}\ {\isacharparenleft}{\kern0pt}\isabound{sh}\ \isabound{z}{\isacharparenright}{\kern0pt}{\isacharparenright}{\kern0pt}\ \isactrlbold {\isacharequal}{\kern0pt}\ \isaconst{trmPredR}\ {\isacharparenleft}{\kern0pt}\isaconst{cmp}\ \isaconst{opt{\isacharunderscore}{\kern0pt}hd{\isacharunderscore}{\kern0pt}encIdx}\ \isactrlbold {\isasymcdot}\ \isaconst{lstDrop}\ {\isacharparenleft}{\kern0pt}\isaconst{cpSnd}\ \isabound{k}{\isacharparenright}{\kern0pt}\ \isabound{z}\ \isactrlbold {\isasymcdot}\ \isaconst{lstDrop}\ {\isacharparenleft}{\kern0pt}\isaconst{cpSnd}\ \isabound{k}{\isacharparenright}{\kern0pt}\ \isabound{z}{\isacharparenright}{\kern0pt}{\isacharparenright}{\kern0pt}\ {\isasymlongrightarrow}\ \isabound{{\isasymGamma}}\ {\isasymturnstile}\ \isaconst{cmp}\ \isavar{E}\ \isactrlbold {\isasymcdot}\ \isaconst{cpPair}\ {\isacharparenleft}{\kern0pt}\isaconst{cpSnd}\ \isabound{k}{\isacharparenright}{\kern0pt}\ {\isacharparenleft}{\kern0pt}\isaconst{cpFst}\ {\isacharparenleft}{\kern0pt}\isaconst{cpFst}\ \isabound{k}{\isacharparenright}{\kern0pt}{\isacharparenright}{\kern0pt}\ \isactrlbold {\isasymcdot}\ \isaconst{cpSnd}\ {\isacharparenleft}{\kern0pt}\isaconst{cpFst}\ \isabound{k}{\isacharparenright}{\kern0pt}\ \isactrlbold {\isasymnoteq}\ \isactrlbold {\isadigit{0}}\ {\isasymlongrightarrow}\ \isabound{{\isasymGamma}}\ {\isasymturnstile}\ \isaconst{cmp}\ \isavar{E}\ \isactrlbold {\isasymcdot}\ \isaconst{cpPair}\ {\isacharparenleft}{\kern0pt}\isaconst{cpSnd}\ \isabound{k{\isacharprime}{\kern0pt}}{\isacharparenright}{\kern0pt}\ {\isacharparenleft}{\kern0pt}\isaconst{cpFst}\ {\isacharparenleft}{\kern0pt}\isaconst{cpFst}\ \isabound{k{\isacharprime}{\kern0pt}}{\isacharparenright}{\kern0pt}{\isacharparenright}{\kern0pt}\ \isactrlbold {\isasymcdot}\ \isaconst{cpSnd}\ {\isacharparenleft}{\kern0pt}\isaconst{cpFst}\ \isabound{k{\isacharprime}{\kern0pt}}{\isacharparenright}{\kern0pt}\ \isactrlbold {\isasymnoteq}\ \isactrlbold {\isadigit{0}}{\isacharparenright}{\kern0pt}}%
}%
\SNIP{Goedel-Paper-Statements:thm:tblShift-orNest-0}{%
\isakeywordONE{thm}\isamarkupfalse%
\ tblShift{\isacharunderscore}{\kern0pt}orNest%
\isanewline%
}%
\SNIP{Goedel-Paper-Statements:thm:tblShift-orNest-1}{%
\isa{{\isacharparenleft}{\kern0pt}{\isasymAnd}\isabound{E}{\isachardot}{\kern0pt}\ \isabound{E}\ {\isasymin}\ \isaconst{set}\ \isavar{L}\ {\isasymLongrightarrow}\ \isaconst{tblShift}\ \isabound{E}{\isacharparenright}{\kern0pt}\ {\isasymLongrightarrow}\ \isaconst{tblShift}\ {\isacharparenleft}{\kern0pt}\isaconst{orNestIdx}\ \isavar{L}{\isacharparenright}{\kern0pt}}%
}%
\SNIP{Goedel-Paper-Statements:thm:tblShift-RGA-0}{%
\isakeywordONE{thm}\isamarkupfalse%
\ tblShift{\isacharunderscore}{\kern0pt}RGA%
\isanewline%
}%
\SNIP{Goedel-Paper-Statements:thm:tblShift-RGA-1}{%
\isa{\isaconst{tblShift}\ {\isacharparenleft}{\kern0pt}\isaconst{drvChkIdxG}\ \isaconst{drvChkList}{\isacharparenright}{\kern0pt}}%
}%
\SNIP{Goedel-Paper-Statements:thm:G-iTreeValid-appC-0}{%
\isakeywordONE{thm}\isamarkupfalse%
\ G{\isacharunderscore}{\kern0pt}iTreeValid{\isacharunderscore}{\kern0pt}appC%
\isanewline%
}%
\SNIP{Goedel-Paper-Statements:thm:G-iTreeValid-appC-1}{%
\isa{{\isasymlbrakk}\isaconst{tblShift}\ \isavar{Dt}{\isacharsemicolon}{\kern0pt}\ \isavar{{\isasymGamma}}\ {\isasymturnstile}\ \isavar{A}\ \isactrlbold N{\isacharsemicolon}{\kern0pt}\ \isavar{{\isasymGamma}}\ {\isasymturnstile}\ \isavar{B}\ \isactrlbold N{\isacharsemicolon}{\kern0pt}\ \isavar{{\isasymGamma}}\ {\isasymturnstile}\ \isaconst{cmp}\ {\isacharparenleft}{\kern0pt}\isaconst{validChkD}\ \isavar{Dt}{\isacharparenright}{\kern0pt}\ \isactrlbold {\isasymcdot}\ \isavar{A}\ \isactrlbold {\isasymcdot}\ \isavar{A}\ \isactrlbold {\isacharequal}{\kern0pt}\ \isactrlbold n\isaconst{Suc}\ {\isadigit{0}}{\isacharsemicolon}{\kern0pt}\ \isavar{{\isasymGamma}}\ {\isasymturnstile}\ \isaconst{cmp}\ {\isacharparenleft}{\kern0pt}\isaconst{validChkD}\ \isavar{Dt}{\isacharparenright}{\kern0pt}\ \isactrlbold {\isasymcdot}\ \isavar{B}\ \isactrlbold {\isasymcdot}\ \isavar{B}\ \isactrlbold {\isacharequal}{\kern0pt}\ \isactrlbold n\isaconst{Suc}\ {\isadigit{0}}{\isasymrbrakk}\ {\isasymLongrightarrow}\ \isavar{{\isasymGamma}}\ {\isasymturnstile}\ \isaconst{cmp}\ {\isacharparenleft}{\kern0pt}\isaconst{validChkD}\ \isavar{Dt}{\isacharparenright}{\kern0pt}\ \isactrlbold {\isasymcdot}\ {\isacharparenleft}{\kern0pt}\isaconst{cmp}\ \isaconst{apCoreIdx}\ \isactrlbold {\isasymcdot}\ \isaconst{cpPair}\ \isavar{A}\ \isavar{B}\ \isactrlbold {\isasymcdot}\ {\isacharparenleft}{\kern0pt}\isaconst{cmp}\ \isaconst{length{\isacharunderscore}{\kern0pt}encIdx}\ \isactrlbold {\isasymcdot}\ \isavar{A}\ \isactrlbold {\isasymcdot}\ \isavar{A}{\isacharparenright}{\kern0pt}{\isacharparenright}{\kern0pt}\ \isactrlbold {\isasymcdot}\ {\isacharparenleft}{\kern0pt}\isaconst{cmp}\ \isaconst{apCoreIdx}\ \isactrlbold {\isasymcdot}\ \isaconst{cpPair}\ \isavar{A}\ \isavar{B}\ \isactrlbold {\isasymcdot}\ {\isacharparenleft}{\kern0pt}\isaconst{cmp}\ \isaconst{length{\isacharunderscore}{\kern0pt}encIdx}\ \isactrlbold {\isasymcdot}\ \isavar{A}\ \isactrlbold {\isasymcdot}\ \isavar{A}{\isacharparenright}{\kern0pt}{\isacharparenright}{\kern0pt}\ \isactrlbold {\isacharequal}{\kern0pt}\ \isactrlbold n\isaconst{Suc}\ {\isadigit{0}}}%
}%
\SNIP{Goedel-Paper-Statements:thm:chkG-iTreeValidD-mono-0}{%
\isakeywordONE{thm}\isamarkupfalse%
\ chkG{\isacharunderscore}{\kern0pt}iTreeValidD{\isacharunderscore}{\kern0pt}mono%
\isanewline%
}%
\SNIP{Goedel-Paper-Statements:thm:chkG-iTreeValidD-mono-1}{%
\isa{{\isasymlbrakk}{\isasymAnd}\isabound{{\isasymDelta}}\ \isabound{u}\ \isabound{y}{\isachardot}{\kern0pt}\ {\isasymlbrakk}\isabound{{\isasymDelta}}\ {\isasymturnstile}\ \isabound{u}\ \isactrlbold N{\isacharsemicolon}{\kern0pt}\ \isabound{{\isasymDelta}}\ {\isasymturnstile}\ \isabound{y}\ \isactrlbold N{\isacharsemicolon}{\kern0pt}\ \isabound{{\isasymDelta}}\ {\isasymturnstile}\ \isaconst{cmp}\ \isavar{D{\isadigit{1}}{\isachardot}{\kern0pt}{\isadigit{0}}}\ \isactrlbold {\isasymcdot}\ \isabound{u}\ \isactrlbold {\isasymcdot}\ \isabound{y}\ \isactrlbold {\isasymnoteq}\ \isactrlbold {\isadigit{0}}{\isasymrbrakk}\ {\isasymLongrightarrow}\ \isabound{{\isasymDelta}}\ {\isasymturnstile}\ \isaconst{cmp}\ \isavar{D{\isadigit{2}}{\isachardot}{\kern0pt}{\isadigit{0}}}\ \isactrlbold {\isasymcdot}\ \isabound{u}\ \isactrlbold {\isasymcdot}\ \isabound{y}\ \isactrlbold {\isasymnoteq}\ \isactrlbold {\isadigit{0}}{\isacharsemicolon}{\kern0pt}\ \isavar{{\isasymGamma}}\ {\isasymturnstile}\ \isavar{m}\ \isactrlbold N{\isacharsemicolon}{\kern0pt}\ \isavar{{\isasymGamma}}\ {\isasymturnstile}\ \isaconst{cmp}\ {\isacharparenleft}{\kern0pt}\isaconst{validChkD}\ \isavar{D{\isadigit{1}}{\isachardot}{\kern0pt}{\isadigit{0}}}{\isacharparenright}{\kern0pt}\ \isactrlbold {\isasymcdot}\ \isavar{m}\ \isactrlbold {\isasymcdot}\ \isavar{m}\ \isactrlbold {\isacharequal}{\kern0pt}\ \isactrlbold n\isaconst{Suc}\ {\isadigit{0}}{\isasymrbrakk}\ {\isasymLongrightarrow}\ \isavar{{\isasymGamma}}\ {\isasymturnstile}\ \isaconst{cmp}\ {\isacharparenleft}{\kern0pt}\isaconst{validChkD}\ \isavar{D{\isadigit{2}}{\isachardot}{\kern0pt}{\isadigit{0}}}{\isacharparenright}{\kern0pt}\ \isactrlbold {\isasymcdot}\ \isavar{m}\ \isactrlbold {\isasymcdot}\ \isavar{m}\ \isactrlbold {\isacharequal}{\kern0pt}\ \isactrlbold n\isaconst{Suc}\ {\isadigit{0}}}%
}%
\SNIP{Goedel-Paper-Statements:thm:chkG-embed-0}{%
\isakeywordONE{thm}\isamarkupfalse%
\ chkG{\isacharunderscore}{\kern0pt}embed%
\isanewline%
}%
\SNIP{Goedel-Paper-Statements:thm:chkG-embed-1}{%
\isa{{\isasymlbrakk}\isaconst{set}\ \isaconst{drvChkList}\ {\isasymsubseteq}\ \isaconst{set}\ \isavar{L}{\isacharsemicolon}{\kern0pt}\ \isavar{{\isasymGamma}}\ {\isasymturnstile}\ \isavar{m}\ \isactrlbold N{\isacharsemicolon}{\kern0pt}\ \isavar{{\isasymGamma}}\ {\isasymturnstile}\ \isaconst{cmp}\ \isaconst{validChkIdx}\ \isactrlbold {\isasymcdot}\ \isavar{m}\ \isactrlbold {\isasymcdot}\ \isavar{m}\ \isactrlbold {\isacharequal}{\kern0pt}\ \isactrlbold n\isaconst{Suc}\ {\isadigit{0}}{\isasymrbrakk}\ {\isasymLongrightarrow}\ \isavar{{\isasymGamma}}\ {\isasymturnstile}\ \isaconst{cmp}\ {\isacharparenleft}{\kern0pt}\isaconst{validChkD}\ {\isacharparenleft}{\kern0pt}\isaconst{drvChkIdxG}\ \isavar{L}{\isacharparenright}{\kern0pt}{\isacharparenright}{\kern0pt}\ \isactrlbold {\isasymcdot}\ \isavar{m}\ \isactrlbold {\isasymcdot}\ \isavar{m}\ \isactrlbold {\isacharequal}{\kern0pt}\ \isactrlbold n\isaconst{Suc}\ {\isadigit{0}}}%
}%
\SNIP{Goedel-Paper-Statements:thm:chkG-embed-cert-0}{%
\isakeywordONE{thm}\isamarkupfalse%
\ chkG{\isacharunderscore}{\kern0pt}embed{\isacharunderscore}{\kern0pt}cert%
\isanewline%
}%
\SNIP{Goedel-Paper-Statements:thm:chkG-embed-cert-1}{%
\isa{{\isasymlbrakk}\isaconst{set}\ \isaconst{drvChkList}\ {\isasymsubseteq}\ \isaconst{set}\ \isavar{L}{\isacharsemicolon}{\kern0pt}\ \isavar{{\isasymGamma}}\ {\isasymturnstile}\ \isavar{HD}\ \isactrlbold N{\isacharsemicolon}{\kern0pt}\ \isavar{{\isasymGamma}}\ {\isasymturnstile}\ \isaconst{iHasCert}\ \isavar{HD}{\isasymrbrakk}\ {\isasymLongrightarrow}\ \isavar{{\isasymGamma}}\ {\isasymturnstile}\ \isaconst{iHasCertD}\ {\isacharparenleft}{\kern0pt}\isaconst{drvChkIdxG}\ \isavar{L}{\isacharparenright}{\kern0pt}\ \isavar{HD}}%
}%
\SNIP{Goedel-Paper-Statements:thm:chkG-iChkD-embed-0}{%
\isakeywordONE{thm}\isamarkupfalse%
\ chkG{\isacharunderscore}{\kern0pt}iChkD{\isacharunderscore}{\kern0pt}embed%
\isanewline%
}%
\SNIP{Goedel-Paper-Statements:thm:chkG-iChkD-embed-1}{%
\isa{{\isasymlbrakk}\isaconst{set}\ \isaconst{drvChkList}\ {\isasymsubseteq}\ \isaconst{set}\ \isavar{L}{\isacharsemicolon}{\kern0pt}\ \isavar{{\isasymGamma}}\ {\isasymturnstile}\ \isavar{p}\ \isactrlbold N{\isacharsemicolon}{\kern0pt}\ \isavar{{\isasymGamma}}\ {\isasymturnstile}\ \isaconst{iChk}\ \isavar{p}\ \isavar{c}{\isasymrbrakk}\ {\isasymLongrightarrow}\ \isavar{{\isasymGamma}}\ {\isasymturnstile}\ \isaconst{iChkD}\ {\isacharparenleft}{\kern0pt}\isaconst{drvChkIdxG}\ \isavar{L}{\isacharparenright}{\kern0pt}\ \isavar{p}\ \isavar{c}}%
}%
\SNIP{Goedel-Paper-Statements:thm:chkG-iChkD-faithful-0}{%
\isakeywordONE{thm}\isamarkupfalse%
\ chkG{\isacharunderscore}{\kern0pt}iChkD{\isacharunderscore}{\kern0pt}faithful%
\isanewline%
}%
\SNIP{Goedel-Paper-Statements:thm:chkG-iChkD-faithful-1}{%
\isa{{\isacharparenleft}{\kern0pt}{\isasymAnd}\isabound{P}\ \isabound{{\isasymGamma}}\ \isabound{c}{\isachardot}{\kern0pt}\ {\isacharparenleft}{\kern0pt}\isaconst{pr{\isadigit{2}}fi}\ \isavar{Dt}\ {\isasymlangle}{\isasymlbrace}\isabound{P}{\isasymrbrace}{\isacharcomma}{\kern0pt}{\isasymlbrace}\isabound{{\isasymGamma}}{\isasymrbrace}{\isasymrangle}\ {\isasymlbrace}\isabound{c}{\isasymrbrace}\ {\isasymnoteq}\ {\isadigit{0}}{\isacharparenright}{\kern0pt}\ {\isacharequal}{\kern0pt}\ \isavar{D}\ \isabound{P}\ \isabound{{\isasymGamma}}\ \isabound{c}{\isacharparenright}{\kern0pt}\ {\isasymLongrightarrow}\ {\isacharparenleft}{\kern0pt}{\isasymemptyset}\ {\isasymturnstile}\ \isaconst{iChkD}\ \isavar{Dt}\ {\isacharparenleft}{\kern0pt}\isactrlbold n{\isasymlbrace}\isavar{P}{\isasymrbrace}{\isacharparenright}{\kern0pt}\ {\isacharparenleft}{\kern0pt}\isactrlbold n{\isasymlbrace}\isavar{c}{\isasymrbrace}{\isacharparenright}{\kern0pt}{\isacharparenright}{\kern0pt}\ {\isacharequal}{\kern0pt}\ \isaconst{prf{\isacharunderscore}{\kern0pt}derives{\isachardot}{\kern0pt}proves}\ \isavar{D}\ \isavar{P}\ {\isasymemptyset}\ \isavar{c}}%
}%
\SNIP{Goedel-Paper-Statements:thm:chkG-iChkD-genuine-proves-0}{%
\isakeywordONE{thm}\isamarkupfalse%
\ chkG{\isacharunderscore}{\kern0pt}iChkD{\isacharunderscore}{\kern0pt}genuine{\isacharunderscore}{\kern0pt}proves%
\isanewline%
}%
\SNIP{Goedel-Paper-Statements:thm:chkG-iChkD-genuine-proves-1}{%
\isa{{\isasymlbrakk}{\isasymAnd}\isabound{P}\ \isabound{{\isasymGamma}}\ \isabound{c}{\isachardot}{\kern0pt}\ {\isacharparenleft}{\kern0pt}\isaconst{pr{\isadigit{2}}fi}\ \isavar{Dt}\ {\isasymlangle}{\isasymlbrace}\isabound{P}{\isasymrbrace}{\isacharcomma}{\kern0pt}{\isasymlbrace}\isabound{{\isasymGamma}}{\isasymrbrace}{\isasymrangle}\ {\isasymlbrace}\isabound{c}{\isasymrbrace}\ {\isasymnoteq}\ {\isadigit{0}}{\isacharparenright}{\kern0pt}\ {\isacharequal}{\kern0pt}\ \isavar{D}\ \isabound{P}\ \isabound{{\isasymGamma}}\ \isabound{c}{\isacharsemicolon}{\kern0pt}\ {\isasymemptyset}\ {\isasymturnstile}\ \isaconst{iChkD}\ \isavar{Dt}\ {\isacharparenleft}{\kern0pt}\isactrlbold n\isavar{p}{\isacharparenright}{\kern0pt}\ {\isacharparenleft}{\kern0pt}\isactrlbold n{\isasymlbrace}\isavar{c}{\isasymrbrace}{\isacharparenright}{\kern0pt}{\isasymrbrakk}\ {\isasymLongrightarrow}\ {\isasymexists}\isabound{P}{\isachardot}{\kern0pt}\ \isavar{p}\ {\isacharequal}{\kern0pt}\ {\isasymlbrace}\isabound{P}{\isasymrbrace}\ {\isasymand}\ \isaconst{prf{\isacharunderscore}{\kern0pt}derives{\isachardot}{\kern0pt}proves}\ \isavar{D}\ \isabound{P}\ {\isasymemptyset}\ \isavar{c}}%
}%
\SNIP{Goedel-Paper-Statements:thm:chkGfixIdx-def-0}{%
\isakeywordONE{thm}\isamarkupfalse%
\ chkGfixIdx{\isacharunderscore}{\kern0pt}def%
\isanewline%
}%
\SNIP{Goedel-Paper-Statements:thm:chkGfixIdx-def-1}{%
\isa{\isaconst{chkGfixIdx}\ {\isasymequiv}\ \isaconst{cite{\isadigit{3}}cIdx}\ \isaconst{payGfix}}%
}%
\SNIP{Goedel-Paper-Statements:thm:pdecF-eq-0}{%
\isakeywordONE{thm}\isamarkupfalse%
\ pdecF{\isacharunderscore}{\kern0pt}eq%
\isanewline%
}%
\SNIP{Goedel-Paper-Statements:thm:pdecF-eq-1}{%
\isa{\isaconst{pdecF}\ \isavar{k}\ \isavar{t}\ {\isacharequal}{\kern0pt}\ {\isacharbrackleft}{\kern0pt}\isaconst{posP}\ \isavar{k}\ \isavar{t}{\isacharcomma}{\kern0pt}\ \isaconst{negP}\ \isavar{k}\ \isavar{t}{\isacharcomma}{\kern0pt}\ \isaconst{pnexp}\ \isavar{k}\ \isavar{t}{\isacharcomma}{\kern0pt}\ {\isasymexists}\isabound{b}{\isachardot}{\kern0pt}\ \isavar{t}\ {\isacharequal}{\kern0pt}\ \isactrlbold {\isasymlambda}\isabound{b}\ {\isasymand}\ \isaconst{posP}\ {\isacharparenleft}{\kern0pt}\isaconst{Suc}\ \isavar{k}{\isacharparenright}{\kern0pt}\ \isabound{b}{\isacharcomma}{\kern0pt}\ {\isasymexists}\isabound{b}{\isachardot}{\kern0pt}\ \isavar{t}\ {\isacharequal}{\kern0pt}\ \isactrlbold {\isasymlambda}\isabound{b}\ {\isasymand}\ \isaconst{negP}\ {\isacharparenleft}{\kern0pt}\isaconst{Suc}\ \isavar{k}{\isacharparenright}{\kern0pt}\ \isabound{b}{\isacharbrackright}{\kern0pt}}%
}%
\SNIP{Goedel-Paper-Statements:thm:pr2fi-pdecPIdx-0}{%
\isakeywordONE{thm}\isamarkupfalse%
\ pr{\isadigit{2}}fi{\isacharunderscore}{\kern0pt}pdecPIdx%
\isanewline%
}%
\SNIP{Goedel-Paper-Statements:thm:pr2fi-pdecPIdx-1}{%
\isa{\isaconst{pr{\isadigit{2}}fi}\ \isaconst{pdecPIdx}\ {\isacharequal}{\kern0pt}\ \isaconst{primrec{\isacharunderscore}{\kern0pt}trmf}\ \isaconst{pdecVar}\ {\isacharparenleft}{\kern0pt}{\isasymlambda}\isabound{n}\ \isabound{v}\ \isabound{t}{\isachardot}{\kern0pt}\ \isabound{v}\ {\isacharplus}{\kern0pt}\ \isaconst{cpx}\ {\isacharparenleft}{\kern0pt}\isaconst{cpy}\ {\isacharparenleft}{\kern0pt}\isaconst{cpy}\ \isabound{t}{\isacharparenright}{\kern0pt}{\isacharparenright}{\kern0pt}{\isacharparenright}{\kern0pt}\ \isaconst{pdecOpr}}%
}%
\SNIP{Goedel-Paper-Statements:thm:chkGfix-fires-iff-0}{%
\isakeywordONE{thm}\isamarkupfalse%
\ chkGfix{\isacharunderscore}{\kern0pt}fires{\isacharunderscore}{\kern0pt}iff%
\isanewline%
}%
\SNIP{Goedel-Paper-Statements:thm:chkGfix-fires-iff-1}{%
\isa{{\isacharparenleft}{\kern0pt}\isaconst{pr{\isadigit{2}}fi}\ \isaconst{chkGfixIdx}\ {\isasymlangle}{\isasymlbrace}\isavar{P}{\isasymrbrace}{\isacharcomma}{\kern0pt}{\isasymlbrace}\isavar{{\isasymGamma}}{\isasymrbrace}{\isasymrangle}\ {\isasymlbrace}\isavar{c}{\isasymrbrace}\ {\isasymnoteq}\ {\isadigit{0}}{\isacharparenright}{\kern0pt}\ {\isacharequal}{\kern0pt}\ \isaconst{gfixRuleI}\ \isavar{P}\ \isavar{{\isasymGamma}}\ \isavar{c}}%
}%
\SNIP{Goedel-Paper-Statements:thm:drvChkI-fires-iff-0}{%
\isakeywordONE{thm}\isamarkupfalse%
\ drvChkI{\isacharunderscore}{\kern0pt}fires{\isacharunderscore}{\kern0pt}iff%
\isanewline%
}%
\SNIP{Goedel-Paper-Statements:thm:drvChkI-fires-iff-1}{%
\isa{{\isacharparenleft}{\kern0pt}\isaconst{pr{\isadigit{2}}fi}\ {\isacharparenleft}{\kern0pt}\isaconst{drvChkIdxG}\ {\isacharparenleft}{\kern0pt}\isaconst{drvChkList}\ {\isacharat}{\kern0pt}\ {\isacharbrackleft}{\kern0pt}\isaconst{chkGfixIdx}{\isacharbrackright}{\kern0pt}{\isacharparenright}{\kern0pt}{\isacharparenright}{\kern0pt}\ {\isasymlangle}{\isasymlbrace}\isavar{P}{\isasymrbrace}{\isacharcomma}{\kern0pt}{\isasymlbrace}\isavar{{\isasymGamma}}{\isasymrbrace}{\isasymrangle}\ {\isasymlbrace}\isavar{c}{\isasymrbrace}\ {\isasymnoteq}\ {\isadigit{0}}{\isacharparenright}{\kern0pt}\ {\isacharequal}{\kern0pt}\ \isaconst{gi{\isachardot}{\kern0pt}derives}\ \isavar{P}\ \isavar{{\isasymGamma}}\ \isavar{c}}%
}%
\SNIP{Goedel-Paper-Statements:thm:tblShift-chkGfix-0}{%
\isakeywordONE{thm}\isamarkupfalse%
\ tblShift{\isacharunderscore}{\kern0pt}chkGfix%
\isanewline%
}%
\SNIP{Goedel-Paper-Statements:thm:tblShift-chkGfix-1}{%
\isa{\isaconst{tblShift}\ \isaconst{chkGfixIdx}}%
}%
\SNIP{Goedel-Paper-Statements:thm:tblShift-IGA-0}{%
\isakeywordONE{thm}\isamarkupfalse%
\ tblShift{\isacharunderscore}{\kern0pt}IGA%
\isanewline%
}%
\SNIP{Goedel-Paper-Statements:thm:tblShift-IGA-1}{%
\isa{\isaconst{tblShift}\ {\isacharparenleft}{\kern0pt}\isaconst{drvChkIdxG}\ {\isacharparenleft}{\kern0pt}\isaconst{drvChkList}\ {\isacharat}{\kern0pt}\ {\isacharbrackleft}{\kern0pt}\isaconst{chkGfixIdx}{\isacharbrackright}{\kern0pt}{\isacharparenright}{\kern0pt}{\isacharparenright}{\kern0pt}}%
}%
\SNIP{Goedel-Paper-Statements:thm:iChkI-faithful-0}{%
\isakeywordONE{thm}\isamarkupfalse%
\ iChkI{\isacharunderscore}{\kern0pt}faithful%
\isanewline%
}%
\SNIP{Goedel-Paper-Statements:thm:iChkI-faithful-1}{%
\isa{{\isacharparenleft}{\kern0pt}{\isasymemptyset}\ {\isasymturnstile}\ \isaconst{iChkD}\ {\isacharparenleft}{\kern0pt}\isaconst{drvChkIdxG}\ {\isacharparenleft}{\kern0pt}\isaconst{drvChkList}\ {\isacharat}{\kern0pt}\ {\isacharbrackleft}{\kern0pt}\isaconst{chkGfixIdx}{\isacharbrackright}{\kern0pt}{\isacharparenright}{\kern0pt}{\isacharparenright}{\kern0pt}\ {\isacharparenleft}{\kern0pt}\isactrlbold n{\isasymlbrace}\isavar{P}{\isasymrbrace}{\isacharparenright}{\kern0pt}\ {\isacharparenleft}{\kern0pt}\isactrlbold n{\isasymlbrace}\isavar{c}{\isasymrbrace}{\isacharparenright}{\kern0pt}{\isacharparenright}{\kern0pt}\ {\isacharequal}{\kern0pt}\ \isaconst{prf{\isacharunderscore}{\kern0pt}derives{\isachardot}{\kern0pt}proves}\ \isaconst{gi{\isachardot}{\kern0pt}derives}\ \isavar{P}\ {\isasymemptyset}\ \isavar{c}}%
}%
\SNIP{Goedel-Paper-Statements:thm:iChkI-genuine-proves-0}{%
\isakeywordONE{thm}\isamarkupfalse%
\ iChkI{\isacharunderscore}{\kern0pt}genuine{\isacharunderscore}{\kern0pt}proves%
\isanewline%
}%
\SNIP{Goedel-Paper-Statements:thm:iChkI-genuine-proves-1}{%
\isa{{\isasymemptyset}\ {\isasymturnstile}\ \isaconst{iChkD}\ {\isacharparenleft}{\kern0pt}\isaconst{drvChkIdxG}\ {\isacharparenleft}{\kern0pt}\isaconst{drvChkList}\ {\isacharat}{\kern0pt}\ {\isacharbrackleft}{\kern0pt}\isaconst{chkGfixIdx}{\isacharbrackright}{\kern0pt}{\isacharparenright}{\kern0pt}{\isacharparenright}{\kern0pt}\ {\isacharparenleft}{\kern0pt}\isactrlbold n\isavar{p}{\isacharparenright}{\kern0pt}\ {\isacharparenleft}{\kern0pt}\isactrlbold n{\isasymlbrace}\isavar{c}{\isasymrbrace}{\isacharparenright}{\kern0pt}\ {\isasymLongrightarrow}\ {\isasymexists}\isabound{P}{\isachardot}{\kern0pt}\ \isavar{p}\ {\isacharequal}{\kern0pt}\ {\isasymlbrace}\isabound{P}{\isasymrbrace}\ {\isasymand}\ \isaconst{prf{\isacharunderscore}{\kern0pt}derives{\isachardot}{\kern0pt}proves}\ \isaconst{gi{\isachardot}{\kern0pt}derives}\ \isabound{P}\ {\isasymemptyset}\ \isavar{c}}%
}%
\SNIP{Goedel-Paper-Statements:thm:G-cut1M-chk-0}{%
\isakeywordONE{thm}\isamarkupfalse%
\ G{\isacharunderscore}{\kern0pt}cut{\isadigit{1}}M{\isacharunderscore}{\kern0pt}chk%
\isanewline%
}%
\SNIP{Goedel-Paper-Statements:thm:G-cut1M-chk-1}{%
\isa{{\isasymlbrakk}{\isasymAnd}\isabound{P}\ \isabound{{\isasymGamma}}\ \isabound{c}{\isachardot}{\kern0pt}\ {\isacharparenleft}{\kern0pt}\isaconst{pr{\isadigit{2}}fi}\ {\isacharparenleft}{\kern0pt}\isaconst{drvChkIdxG}\ \isavar{L}{\isacharparenright}{\kern0pt}\ {\isasymlangle}{\isasymlbrace}\isabound{P}{\isasymrbrace}{\isacharcomma}{\kern0pt}{\isasymlbrace}\isabound{{\isasymGamma}}{\isasymrbrace}{\isasymrangle}\ {\isasymlbrace}\isabound{c}{\isasymrbrace}\ {\isasymnoteq}\ {\isadigit{0}}{\isacharparenright}{\kern0pt}\ {\isacharequal}{\kern0pt}\ \isavar{D}\ \isabound{P}\ \isabound{{\isasymGamma}}\ \isabound{c}{\isacharsemicolon}{\kern0pt}\ \isaconst{tblShift}\ {\isacharparenleft}{\kern0pt}\isaconst{drvChkIdxG}\ \isavar{L}{\isacharparenright}{\kern0pt}{\isacharsemicolon}{\kern0pt}\ \isaconst{chkCutIdx}\ {\isasymin}\ \isaconst{set}\ \isavar{L}{\isacharsemicolon}{\kern0pt}\ \isaconst{prf{\isacharunderscore}{\kern0pt}derives{\isachardot}{\kern0pt}proves}\ \isavar{D}\ \isavar{Q}\ {\isacharparenleft}{\kern0pt}\isavar{A}\ {\isasymtriangleright}\ {\isasymemptyset}{\isacharparenright}{\kern0pt}\ \isavar{B}{\isacharsemicolon}{\kern0pt}\ \isavar{{\isasymGamma}}\ {\isasymturnstile}\ \isavar{a}\ \isactrlbold N{\isacharsemicolon}{\kern0pt}\ \isavar{{\isasymGamma}}\ {\isasymturnstile}\ \isaconst{iChkD}\ {\isacharparenleft}{\kern0pt}\isaconst{drvChkIdxG}\ \isavar{L}{\isacharparenright}{\kern0pt}\ \isavar{a}\ {\isacharparenleft}{\kern0pt}\isactrlbold n{\isasymlbrace}\isavar{A}{\isasymrbrace}{\isacharparenright}{\kern0pt}{\isasymrbrakk}\ {\isasymLongrightarrow}\ \isavar{{\isasymGamma}}\ {\isasymturnstile}\ \isaconst{iChkD}\ {\isacharparenleft}{\kern0pt}\isaconst{drvChkIdxG}\ \isavar{L}{\isacharparenright}{\kern0pt}\ {\isacharparenleft}{\kern0pt}\isaconst{cut{\isadigit{1}}M}\ \isavar{Q}\ \isavar{B}\ \isactrlbold {\isasymcdot}\ \isavar{a}\ \isactrlbold {\isasymcdot}\ \isavar{a}{\isacharparenright}{\kern0pt}\ {\isacharparenleft}{\kern0pt}\isactrlbold n{\isasymlbrace}\isavar{B}{\isasymrbrace}{\isacharparenright}{\kern0pt}}%
}%
\SNIP{Goedel-Paper-Statements:thm:G-goedelE-chk-0}{%
\isakeywordONE{thm}\isamarkupfalse%
\ G{\isacharunderscore}{\kern0pt}goedelE{\isacharunderscore}{\kern0pt}chk%
\isanewline%
}%
\SNIP{Goedel-Paper-Statements:thm:G-goedelE-chk-1}{%
\isa{{\isasymlbrakk}\isaconst{tblShift}\ {\isacharparenleft}{\kern0pt}\isaconst{drvChkIdxG}\ \isavar{L}{\isacharparenright}{\kern0pt}{\isacharsemicolon}{\kern0pt}\ \isaconst{chkNegEIdx}\ {\isasymin}\ \isaconst{set}\ \isavar{L}{\isacharsemicolon}{\kern0pt}\ \isavar{{\isasymGamma}}\ {\isasymturnstile}\ \isavar{a}\ \isactrlbold N{\isacharsemicolon}{\kern0pt}\ \isavar{{\isasymGamma}}\ {\isasymturnstile}\ \isavar{b}\ \isactrlbold N{\isacharsemicolon}{\kern0pt}\ \isavar{{\isasymGamma}}\ {\isasymturnstile}\ \isaconst{iChkD}\ {\isacharparenleft}{\kern0pt}\isaconst{drvChkIdxG}\ \isavar{L}{\isacharparenright}{\kern0pt}\ \isavar{a}\ {\isacharparenleft}{\kern0pt}\isactrlbold n{\isasymlbrace}\isavar{g}{\isasymrbrace}{\isacharparenright}{\kern0pt}{\isacharsemicolon}{\kern0pt}\ \isavar{{\isasymGamma}}\ {\isasymturnstile}\ \isaconst{iChkD}\ {\isacharparenleft}{\kern0pt}\isaconst{drvChkIdxG}\ \isavar{L}{\isacharparenright}{\kern0pt}\ \isavar{b}\ {\isacharparenleft}{\kern0pt}\isactrlbold n{\isasymlbrace}\isactrlbold {\isasymnot}\isavar{g}{\isasymrbrace}{\isacharparenright}{\kern0pt}{\isasymrbrakk}\ {\isasymLongrightarrow}\ \isavar{{\isasymGamma}}\ {\isasymturnstile}\ \isaconst{iChkD}\ {\isacharparenleft}{\kern0pt}\isaconst{drvChkIdxG}\ \isavar{L}{\isacharparenright}{\kern0pt}\ {\isacharparenleft}{\kern0pt}\isaconst{goedelE}\ \isactrlbold {\isasymcdot}\ \isavar{a}\ \isactrlbold {\isasymcdot}\ \isavar{b}{\isacharparenright}{\kern0pt}\ \isaconst{cEq{\isadigit{0}}{\isadigit{1}}}}%
}%
\SNIP{Goedel-Paper-Statements:thm:chkG-iChkD-proves-sf-0}{%
\isakeywordONE{thm}\isamarkupfalse%
\ chkG{\isacharunderscore}{\kern0pt}iChkD{\isacharunderscore}{\kern0pt}proves{\isacharunderscore}{\kern0pt}sf%
\isanewline%
}%
\SNIP{Goedel-Paper-Statements:thm:chkG-iChkD-proves-sf-1}{%
\isa{{\isasymlbrakk}{\isasymAnd}\isabound{P}\ \isabound{{\isasymGamma}}\ \isabound{c}{\isachardot}{\kern0pt}\ {\isacharparenleft}{\kern0pt}\isaconst{pr{\isadigit{2}}fi}\ \isavar{Dt}\ {\isasymlangle}{\isasymlbrace}\isabound{P}{\isasymrbrace}{\isacharcomma}{\kern0pt}{\isasymlbrace}\isabound{{\isasymGamma}}{\isasymrbrace}{\isasymrangle}\ {\isasymlbrace}\isabound{c}{\isasymrbrace}\ {\isasymnoteq}\ {\isadigit{0}}{\isacharparenright}{\kern0pt}\ {\isacharequal}{\kern0pt}\ \isavar{D}\ \isabound{P}\ \isabound{{\isasymGamma}}\ \isabound{c}{\isacharsemicolon}{\kern0pt}\ \isaconst{prf{\isacharunderscore}{\kern0pt}derives{\isachardot}{\kern0pt}proves}\ \isavar{D}\ \isavar{P}\ {\isasymemptyset}\ \isavar{c}{\isasymrbrakk}\ {\isasymLongrightarrow}\ \isavar{{\isasymGamma}}\ {\isasymturnstile}\ \isaconst{iChkD}\ \isavar{Dt}\ {\isacharparenleft}{\kern0pt}\isactrlbold n{\isasymlbrace}\isavar{P}{\isasymrbrace}{\isacharparenright}{\kern0pt}\ {\isacharparenleft}{\kern0pt}\isactrlbold n{\isasymlbrace}\isavar{c}{\isasymrbrace}{\isacharparenright}{\kern0pt}}%
}%
\SNIP{Goedel-Paper-Statements:thm:chkI-proves-sf-0}{%
\isakeywordONE{thm}\isamarkupfalse%
\ chkI{\isacharunderscore}{\kern0pt}proves{\isacharunderscore}{\kern0pt}sf%
\isanewline%
}%
\SNIP{Goedel-Paper-Statements:thm:chkI-proves-sf-1}{%
\isa{\isaconst{prf{\isacharunderscore}{\kern0pt}derives{\isachardot}{\kern0pt}proves}\ \isaconst{gi{\isachardot}{\kern0pt}derives}\ \isavar{P}\ {\isasymemptyset}\ \isavar{c}\ {\isasymLongrightarrow}\ \isavar{{\isasymGamma}}\ {\isasymturnstile}\ \isaconst{iChkD}\ {\isacharparenleft}{\kern0pt}\isaconst{drvChkIdxG}\ {\isacharparenleft}{\kern0pt}\isaconst{drvChkList}\ {\isacharat}{\kern0pt}\ {\isacharbrackleft}{\kern0pt}\isaconst{chkGfixIdx}{\isacharbrackright}{\kern0pt}{\isacharparenright}{\kern0pt}{\isacharparenright}{\kern0pt}\ {\isacharparenleft}{\kern0pt}\isactrlbold n{\isasymlbrace}\isavar{P}{\isasymrbrace}{\isacharparenright}{\kern0pt}\ {\isacharparenleft}{\kern0pt}\isactrlbold n{\isasymlbrace}\isavar{c}{\isasymrbrace}{\isacharparenright}{\kern0pt}}%
}%
\SNIP{Goedel-Paper-Statements:thm:G-iChk-d3SigmaAt-0}{%
\isakeywordONE{thm}\isamarkupfalse%
\ G{\isacharunderscore}{\kern0pt}iChk{\isacharunderscore}{\kern0pt}d{\isadigit{3}}SigmaAt%
\isanewline%
}%
\SNIP{Goedel-Paper-Statements:thm:G-iChk-d3SigmaAt-1}{%
\isa{{\isasymlbrakk}\isaconst{d{\isadigit{3}}Src}\ \isavar{S}\ \isavar{SI}\ \isavar{sc}{\isacharsemicolon}{\kern0pt}\ \isavar{{\isasymGamma}}\ {\isasymturnstile}\ \isavar{v}\ \isactrlbold N{\isacharsemicolon}{\kern0pt}\ \isavar{{\isasymGamma}}\ {\isasymturnstile}\ \isaconst{iChkD}\ \isavar{Dt}\ \isavar{v}\ \isavar{S}{\isasymrbrakk}\ {\isasymLongrightarrow}\ \isavar{{\isasymGamma}}\ {\isasymturnstile}\ \isaconst{iChk}\ {\isacharparenleft}{\kern0pt}\isaconst{d{\isadigit{3}}SigmaAtD}\ \isavar{Dt}\ \isavar{S}\ \isavar{SI}\ \isavar{sc}\ \isactrlbold {\isasymcdot}\ \isavar{v}{\isacharparenright}{\kern0pt}\ {\isacharparenleft}{\kern0pt}\isactrlbold n{\isasymlbrace}\isactrlbold {\isasymnot}\isactrlbold A\ \isactrlbold {\isasymcdot}\ \isaconst{prvFD}\ \isavar{Dt}\ \isavar{S}{\isasymrbrace}{\isacharparenright}{\kern0pt}}%
}%
\SNIP{Goedel-Paper-Statements:thm:chkG-iChkD-d3SigmaAt-0}{%
\isakeywordONE{thm}\isamarkupfalse%
\ chkG{\isacharunderscore}{\kern0pt}iChkD{\isacharunderscore}{\kern0pt}d{\isadigit{3}}SigmaAt%
\isanewline%
}%
\SNIP{Goedel-Paper-Statements:thm:chkG-iChkD-d3SigmaAt-1}{%
\isa{{\isasymlbrakk}\isaconst{set}\ \isaconst{drvChkList}\ {\isasymsubseteq}\ \isaconst{set}\ \isavar{L}{\isacharsemicolon}{\kern0pt}\ \isaconst{d{\isadigit{3}}Src}\ \isavar{S}\ \isavar{SI}\ \isavar{sc}{\isacharsemicolon}{\kern0pt}\ \isavar{{\isasymGamma}}\ {\isasymturnstile}\ \isavar{v}\ \isactrlbold N{\isacharsemicolon}{\kern0pt}\ \isavar{{\isasymGamma}}\ {\isasymturnstile}\ \isaconst{iChkD}\ {\isacharparenleft}{\kern0pt}\isaconst{drvChkIdxG}\ \isavar{L}{\isacharparenright}{\kern0pt}\ \isavar{v}\ \isavar{S}{\isasymrbrakk}\ {\isasymLongrightarrow}\ \isavar{{\isasymGamma}}\ {\isasymturnstile}\ \isaconst{iChkD}\ {\isacharparenleft}{\kern0pt}\isaconst{drvChkIdxG}\ \isavar{L}{\isacharparenright}{\kern0pt}\ {\isacharparenleft}{\kern0pt}\isaconst{d{\isadigit{3}}SigmaAtD}\ {\isacharparenleft}{\kern0pt}\isaconst{drvChkIdxG}\ \isavar{L}{\isacharparenright}{\kern0pt}\ \isavar{S}\ \isavar{SI}\ \isavar{sc}\ \isactrlbold {\isasymcdot}\ \isavar{v}{\isacharparenright}{\kern0pt}\ {\isacharparenleft}{\kern0pt}\isactrlbold n{\isasymlbrace}\isactrlbold {\isasymnot}\isactrlbold A\ \isactrlbold {\isasymcdot}\ \isaconst{prvFD}\ {\isacharparenleft}{\kern0pt}\isaconst{drvChkIdxG}\ \isavar{L}{\isacharparenright}{\kern0pt}\ \isavar{S}{\isasymrbrace}{\isacharparenright}{\kern0pt}}%
}%
\SNIP{Goedel-Paper-Statements:thm:ConIGA-def-0}{%
\isakeywordONE{thm}\isamarkupfalse%
\ ConIGA{\isacharunderscore}{\kern0pt}def%
\isanewline%
}%
\SNIP{Goedel-Paper-Statements:thm:ConIGA-def-1}{%
\isa{\isaconst{ConIGA}\ {\isasymequiv}\ \isactrlbold A\ \isactrlbold {\isasymcdot}\ {\isacharparenleft}{\kern0pt}\isactrlbold {\isasymlambda}{\isacharparenleft}{\kern0pt}\isactrlbold {\isasymnot}\isaconst{iChkD}\ {\isacharparenleft}{\kern0pt}\isaconst{drvChkIdxG}\ {\isacharparenleft}{\kern0pt}\isaconst{drvChkList}\ {\isacharat}{\kern0pt}\ {\isacharbrackleft}{\kern0pt}\isaconst{chkGfixIdx}{\isacharbrackright}{\kern0pt}{\isacharparenright}{\kern0pt}{\isacharparenright}{\kern0pt}\ \isactrlbold v\isactrlsub {\isadigit{0}}\ \isaconst{cEq{\isadigit{0}}{\isadigit{1}}}{\isacharparenright}{\kern0pt}{\isacharparenright}{\kern0pt}}%
}%
\SNIP{Goedel-Paper-Statements:thm:g2-syn-gi-0}{%
\isakeywordONE{thm}\isamarkupfalse%
\ g{\isadigit{2}}{\isacharunderscore}{\kern0pt}syn{\isacharunderscore}{\kern0pt}gi%
\isanewline%
}%
\SNIP{Goedel-Paper-Statements:thm:g2-syn-gi-1}{%
\isa{\isaconst{g{\isadigit{2}}{\isacharunderscore}{\kern0pt}syn}\ {\isacharparenleft}{\kern0pt}\isaconst{prf{\isacharunderscore}{\kern0pt}derives{\isachardot}{\kern0pt}entails}\ \isaconst{gi{\isachardot}{\kern0pt}derives}{\isacharparenright}{\kern0pt}\ {\isacharparenleft}{\kern0pt}\isaconst{iChkD}\ {\isacharparenleft}{\kern0pt}\isaconst{drvChkIdxG}\ {\isacharparenleft}{\kern0pt}\isaconst{drvChkList}\ {\isacharat}{\kern0pt}\ {\isacharbrackleft}{\kern0pt}\isaconst{chkGfixIdx}{\isacharbrackright}{\kern0pt}{\isacharparenright}{\kern0pt}{\isacharparenright}{\kern0pt}{\isacharparenright}{\kern0pt}}%
}%
\SNIP{Goedel-Paper-Statements:thm:g2-con-gi-0}{%
\isakeywordONE{thm}\isamarkupfalse%
\ g{\isadigit{2}}{\isacharunderscore}{\kern0pt}con{\isacharunderscore}{\kern0pt}gi%
\isanewline%
}%
\SNIP{Goedel-Paper-Statements:thm:g2-con-gi-1}{%
\isa{\isaconst{g{\isadigit{2}}{\isacharunderscore}{\kern0pt}con}\ {\isacharparenleft}{\kern0pt}\isaconst{prf{\isacharunderscore}{\kern0pt}derives{\isachardot}{\kern0pt}entails}\ \isaconst{gi{\isachardot}{\kern0pt}derives}{\isacharparenright}{\kern0pt}\ {\isacharparenleft}{\kern0pt}\isaconst{iChkD}\ {\isacharparenleft}{\kern0pt}\isaconst{drvChkIdxG}\ {\isacharparenleft}{\kern0pt}\isaconst{drvChkList}\ {\isacharat}{\kern0pt}\ {\isacharbrackleft}{\kern0pt}\isaconst{chkGfixIdx}{\isacharbrackright}{\kern0pt}{\isacharparenright}{\kern0pt}{\isacharparenright}{\kern0pt}{\isacharparenright}{\kern0pt}}%
}%
\SNIP{Goedel-Paper-Statements:thm:gi-goedel2-internal-contra-0}{%
\isakeywordONE{thm}\isamarkupfalse%
\ gi{\isacharunderscore}{\kern0pt}goedel{\isadigit{2}}{\isacharunderscore}{\kern0pt}internal{\isacharunderscore}{\kern0pt}contra%
\isanewline%
}%
\SNIP{Goedel-Paper-Statements:thm:gi-goedel2-internal-contra-1}{%
\isa{\isaconst{prf{\isacharunderscore}{\kern0pt}derives{\isachardot}{\kern0pt}entails}\ \isaconst{gi{\isachardot}{\kern0pt}derives}\ {\isacharparenleft}{\kern0pt}\isactrlbold {\isasymnot}\isactrlbold A\ \isactrlbold {\isasymcdot}\ {\isacharparenleft}{\kern0pt}\isactrlbold {\isasymlambda}{\isacharparenleft}{\kern0pt}\isactrlbold {\isasymnot}\isaconst{iChkD}\ {\isacharparenleft}{\kern0pt}\isaconst{drvChkIdxG}\ {\isacharparenleft}{\kern0pt}\isaconst{drvChkList}\ {\isacharat}{\kern0pt}\ {\isacharbrackleft}{\kern0pt}\isaconst{chkGfixIdx}{\isacharbrackright}{\kern0pt}{\isacharparenright}{\kern0pt}{\isacharparenright}{\kern0pt}\ \isactrlbold v\isactrlsub {\isadigit{0}}\ {\isacharparenleft}{\kern0pt}{\isasymupharpoonleft}{\isacharparenleft}{\kern0pt}\isactrlbold n{\isasymlbrace}\isaconst{ConIGA}{\isasymrbrace}{\isacharparenright}{\kern0pt}{\isacharparenright}{\kern0pt}{\isacharparenright}{\kern0pt}{\isacharparenright}{\kern0pt}\ {\isasymtriangleright}\ {\isasymemptyset}{\isacharparenright}{\kern0pt}\ {\isacharparenleft}{\kern0pt}\isactrlbold {\isasymnot}\isaconst{ConIGA}{\isacharparenright}{\kern0pt}}%
}%
\SNIP{Goedel-Paper-Statements:thm:gi-goedel2-internal-0}{%
\isakeywordONE{thm}\isamarkupfalse%
\ gi{\isacharunderscore}{\kern0pt}goedel{\isadigit{2}}{\isacharunderscore}{\kern0pt}internal%
\isanewline%
}%
\SNIP{Goedel-Paper-Statements:thm:gi-goedel2-internal-1}{%
\isa{\isaconst{prf{\isacharunderscore}{\kern0pt}derives{\isachardot}{\kern0pt}entails}\ \isaconst{gi{\isachardot}{\kern0pt}derives}\ {\isacharparenleft}{\kern0pt}\isaconst{ConIGA}\ {\isasymtriangleright}\ {\isasymemptyset}{\isacharparenright}{\kern0pt}\ {\isacharparenleft}{\kern0pt}\isactrlbold {\isasymnot}\isactrlbold {\isasymnot}\isactrlbold A\ \isactrlbold {\isasymcdot}\ {\isacharparenleft}{\kern0pt}\isactrlbold {\isasymlambda}{\isacharparenleft}{\kern0pt}\isactrlbold {\isasymnot}\isaconst{iChkD}\ {\isacharparenleft}{\kern0pt}\isaconst{drvChkIdxG}\ {\isacharparenleft}{\kern0pt}\isaconst{drvChkList}\ {\isacharat}{\kern0pt}\ {\isacharbrackleft}{\kern0pt}\isaconst{chkGfixIdx}{\isacharbrackright}{\kern0pt}{\isacharparenright}{\kern0pt}{\isacharparenright}{\kern0pt}\ \isactrlbold v\isactrlsub {\isadigit{0}}\ {\isacharparenleft}{\kern0pt}{\isasymupharpoonleft}{\isacharparenleft}{\kern0pt}\isactrlbold n{\isasymlbrace}\isaconst{ConIGA}{\isasymrbrace}{\isacharparenright}{\kern0pt}{\isacharparenright}{\kern0pt}{\isacharparenright}{\kern0pt}{\isacharparenright}{\kern0pt}{\isacharparenright}{\kern0pt}}%
}%
\SNIP{Goedel-Paper-Statements:thm:gi-goedelG-refutes-0}{%
\isakeywordONE{thm}\isamarkupfalse%
\ gi{\isacharunderscore}{\kern0pt}goedelG{\isacharunderscore}{\kern0pt}refutes%
\isanewline%
}%
\SNIP{Goedel-Paper-Statements:thm:gi-goedelG-refutes-1}{%
\isa{\isaconst{prf{\isacharunderscore}{\kern0pt}derives{\isachardot}{\kern0pt}entails}\ \isaconst{gi{\isachardot}{\kern0pt}derives}\ {\isasymemptyset}\ {\isacharparenleft}{\kern0pt}\isaconst{g{\isadigit{2}}{\isacharunderscore}{\kern0pt}syn{\isachardot}{\kern0pt}goedelGT}\ {\isacharparenleft}{\kern0pt}\isaconst{iChkD}\ {\isacharparenleft}{\kern0pt}\isaconst{drvChkIdxG}\ {\isacharparenleft}{\kern0pt}\isaconst{drvChkList}\ {\isacharat}{\kern0pt}\ {\isacharbrackleft}{\kern0pt}\isaconst{chkGfixIdx}{\isacharbrackright}{\kern0pt}{\isacharparenright}{\kern0pt}{\isacharparenright}{\kern0pt}{\isacharparenright}{\kern0pt}{\isacharparenright}{\kern0pt}\ {\isasymLongrightarrow}\ \isaconst{prf{\isacharunderscore}{\kern0pt}derives{\isachardot}{\kern0pt}entails}\ \isaconst{gi{\isachardot}{\kern0pt}derives}\ {\isasymemptyset}\ {\isacharparenleft}{\kern0pt}\isactrlbold {\isadigit{0}}\ \isactrlbold {\isacharequal}{\kern0pt}\ \isactrlbold n{\isadigit{1}}{\isacharparenright}{\kern0pt}}%
}%
\SNIP{Goedel-Paper-Statements:thm:gi-goedel1-0}{%
\isakeywordONE{thm}\isamarkupfalse%
\ gi{\isacharunderscore}{\kern0pt}goedel{\isadigit{1}}%
\isanewline%
}%
\SNIP{Goedel-Paper-Statements:thm:gi-goedel1-1}{%
\isa{{\isasymnot}\ \isaconst{prf{\isacharunderscore}{\kern0pt}derives{\isachardot}{\kern0pt}entails}\ \isaconst{gi{\isachardot}{\kern0pt}derives}\ {\isasymemptyset}\ {\isacharparenleft}{\kern0pt}\isaconst{g{\isadigit{2}}{\isacharunderscore}{\kern0pt}syn{\isachardot}{\kern0pt}goedelGT}\ {\isacharparenleft}{\kern0pt}\isaconst{iChkD}\ {\isacharparenleft}{\kern0pt}\isaconst{drvChkIdxG}\ {\isacharparenleft}{\kern0pt}\isaconst{drvChkList}\ {\isacharat}{\kern0pt}\ {\isacharbrackleft}{\kern0pt}\isaconst{chkGfixIdx}{\isacharbrackright}{\kern0pt}{\isacharparenright}{\kern0pt}{\isacharparenright}{\kern0pt}{\isacharparenright}{\kern0pt}{\isacharparenright}{\kern0pt}}%
}%
\SNIP{Goedel-Paper-Statements:thm:gi-goedel2-ConIGA-0}{%
\isakeywordONE{thm}\isamarkupfalse%
\ gi{\isacharunderscore}{\kern0pt}goedel{\isadigit{2}}{\isacharunderscore}{\kern0pt}ConIGA%
\isanewline%
}%
\SNIP{Goedel-Paper-Statements:thm:gi-goedel2-ConIGA-1}{%
\isa{{\isasymnot}\ \isaconst{prf{\isacharunderscore}{\kern0pt}derives{\isachardot}{\kern0pt}entails}\ \isaconst{gi{\isachardot}{\kern0pt}derives}\ {\isasymemptyset}\ \isaconst{ConIGA}}%
}%
\SNIP{Goedel-Paper-Statements:thm:gi-separation-IGA-0}{%
\isakeywordONE{thm}\isamarkupfalse%
\ gi{\isacharunderscore}{\kern0pt}separation{\isacharunderscore}{\kern0pt}IGA%
\isanewline%
}%
\SNIP{Goedel-Paper-Statements:thm:gi-separation-IGA-1}{%
\isa{{\isasymnot}\ \isaconst{prf{\isacharunderscore}{\kern0pt}derives{\isachardot}{\kern0pt}entails}\ \isaconst{gi{\isachardot}{\kern0pt}derives}\ {\isasymemptyset}\ \isaconst{ConIGA}\ {\isasymand}\ \isaconst{prf{\isacharunderscore}{\kern0pt}derives{\isachardot}{\kern0pt}entails}\ \isaconst{adi{\isacharunderscore}{\kern0pt}derives}\ {\isasymemptyset}\ \isaconst{ConIGA}\ {\isasymand}\ {\isacharparenleft}{\kern0pt}{\isasymforall}\isabound{{\isasymGamma}}\ \isabound{c}{\isachardot}{\kern0pt}\ \isaconst{prf{\isacharunderscore}{\kern0pt}derives{\isachardot}{\kern0pt}entails}\ \isaconst{gi{\isachardot}{\kern0pt}derives}\ \isabound{{\isasymGamma}}\ \isabound{c}\ {\isasymlongrightarrow}\ \isaconst{prf{\isacharunderscore}{\kern0pt}derives{\isachardot}{\kern0pt}entails}\ \isaconst{adi{\isacharunderscore}{\kern0pt}derives}\ \isabound{{\isasymGamma}}\ \isabound{c}{\isacharparenright}{\kern0pt}\ {\isasymand}\ {\isasymnot}\ \isaconst{prf{\isacharunderscore}{\kern0pt}derives{\isachardot}{\kern0pt}entails}\ \isaconst{adi{\isacharunderscore}{\kern0pt}derives}\ {\isasymemptyset}\ {\isacharparenleft}{\kern0pt}\isactrlbold {\isadigit{0}}\ \isactrlbold {\isacharequal}{\kern0pt}\ \isactrlbold {\isadigit{1}}{\isacharparenright}{\kern0pt}}%
}%
\SNIP{Goedel-Paper-Statements:thm:gi-adiU-unprovable-0}{%
\isakeywordONE{thm}\isamarkupfalse%
\ gi{\isacharunderscore}{\kern0pt}adiU{\isacharunderscore}{\kern0pt}unprovable%
\isanewline%
}%
\SNIP{Goedel-Paper-Statements:thm:gi-adiU-unprovable-1}{%
\isa{{\isasymnot}\ \isaconst{prf{\isacharunderscore}{\kern0pt}derives{\isachardot}{\kern0pt}entails}\ \isaconst{gi{\isachardot}{\kern0pt}derives}\ {\isasymemptyset}\ {\isacharparenleft}{\kern0pt}\isactrlbold A\ \isactrlbold {\isasymcdot}\ \isaconst{adiU}{\isacharparenright}{\kern0pt}}%
}%
\SNIP{Goedel-Paper-Statements:thm:adrW-def-0}{%
\isakeywordONE{thm}\isamarkupfalse%
\ adrW{\isacharunderscore}{\kern0pt}def%
\isanewline%
}%
\SNIP{Goedel-Paper-Statements:thm:adrW-def-1}{%
\isa{\isaconst{adrW}\ \isavar{F}\ \isavar{u}\ {\isasymequiv}\ \isactrlbold A\ \isactrlbold {\isasymcdot}\ {\isacharparenleft}{\kern0pt}\isactrlbold {\isasymlambda}\isaconst{iChk}\ {\isacharparenleft}{\kern0pt}\isaconst{cmp}\ \isavar{F}\ \isactrlbold {\isasymcdot}\ \isactrlbold v\isactrlsub {\isadigit{0}}\ \isactrlbold {\isasymcdot}\ \isactrlbold v\isactrlsub {\isadigit{0}}{\isacharparenright}{\kern0pt}\ {\isacharparenleft}{\kern0pt}\isaconst{adrInst}\ \isavar{u}\ \isactrlbold v\isactrlsub {\isadigit{0}}{\isacharparenright}{\kern0pt}{\isacharparenright}{\kern0pt}}%
}%
\SNIP{Goedel-Paper-Statements:thm:adr-derives-def-0}{%
\isakeywordONE{thm}\isamarkupfalse%
\ adr{\isacharunderscore}{\kern0pt}derives{\isacharunderscore}{\kern0pt}def%
\isanewline%
}%
\SNIP{Goedel-Paper-Statements:thm:adr-derives-def-1}{%
\isa{\isaconst{adr{\isacharunderscore}{\kern0pt}derives}\ \isavar{P}\ \isavar{{\isasymGamma}}\ \isavar{c}\ {\isacharequal}{\kern0pt}\ {\isacharparenleft}{\kern0pt}\isavar{P}\ {\isasymD}{\isasymlbrakk}\isavar{{\isasymGamma}}\ {\isasymturnstile}\ \isavar{c}{\isasymrbrakk}\ {\isasymor}\ \isaconst{adr{\isacharunderscore}{\kern0pt}step}\ \isavar{P}\ \isavar{{\isasymGamma}}\ \isavar{c}{\isacharparenright}{\kern0pt}}%
}%
\SNIP{Goedel-Paper-Statements:thm:adr-adopt-0}{%
\isakeywordONE{thm}\isamarkupfalse%
\ adr{\isacharunderscore}{\kern0pt}adopt%
\isanewline%
}%
\SNIP{Goedel-Paper-Statements:thm:adr-adopt-1}{%
\isa{{\isasymlbrakk}\isaconst{prf{\isacharunderscore}{\kern0pt}derives{\isachardot}{\kern0pt}entails}\ \isaconst{adr{\isacharunderscore}{\kern0pt}derives}\ \isavar{{\isasymGamma}}\ {\isacharparenleft}{\kern0pt}\isaconst{adrW}\ \isavar{F}\ \isavar{u}{\isacharparenright}{\kern0pt}{\isacharsemicolon}{\kern0pt}\ \isaconst{wft}\ \isavar{u}{\isacharsemicolon}{\kern0pt}\ \isaconst{Syntax{\isacharunderscore}{\kern0pt}Term{\isachardot}{\kern0pt}closed}\ \isavar{u}{\isasymrbrakk}\ {\isasymLongrightarrow}\ \isaconst{prf{\isacharunderscore}{\kern0pt}derives{\isachardot}{\kern0pt}entails}\ \isaconst{adr{\isacharunderscore}{\kern0pt}derives}\ \isavar{{\isasymGamma}}\ {\isacharparenleft}{\kern0pt}\isactrlbold A\ \isactrlbold {\isasymcdot}\ \isavar{u}{\isacharparenright}{\kern0pt}}%
}%
\SNIP{Goedel-Paper-Statements:thm:adr-embed-0}{%
\isakeywordONE{thm}\isamarkupfalse%
\ adr{\isacharunderscore}{\kern0pt}embed%
\isanewline%
}%
\SNIP{Goedel-Paper-Statements:thm:adr-embed-1}{%
\isa{\isavar{{\isasymGamma}}\ {\isasymturnstile}\ \isavar{c}\ {\isasymLongrightarrow}\ \isaconst{prf{\isacharunderscore}{\kern0pt}derives{\isachardot}{\kern0pt}entails}\ \isaconst{adr{\isacharunderscore}{\kern0pt}derives}\ \isavar{{\isasymGamma}}\ \isavar{c}}%
}%
\SNIP{Goedel-Paper-Statements:thm:adr-omega-sound-0}{%
\isakeywordONE{thm}\isamarkupfalse%
\ adr{\isacharunderscore}{\kern0pt}omega{\isacharunderscore}{\kern0pt}sound%
\isanewline%
}%
\SNIP{Goedel-Paper-Statements:thm:adr-omega-sound-1}{%
\isa{{\isasymlbrakk}\isaconst{prf{\isacharunderscore}{\kern0pt}derives{\isachardot}{\kern0pt}entails}\ \isaconst{adr{\isacharunderscore}{\kern0pt}derives}\ \isavar{{\isasymGamma}}\ \isavar{c}{\isacharsemicolon}{\kern0pt}\ \isaconst{wfc}\ \isavar{{\isasymGamma}}{\isacharsemicolon}{\kern0pt}\ {\isasymAnd}\isabound{g}{\isachardot}{\kern0pt}\ \isabound{g}\ {\isacharbar}{\kern0pt}{\isasymin}{\isacharbar}{\kern0pt}\ \isavar{{\isasymGamma}}\ {\isasymLongrightarrow}\ \isaconst{omv}\ \isabound{g}\ {\isasymlbrace}\isaconst{True}{\isasymrbrace}{\isasymrbrakk}\ {\isasymLongrightarrow}\ \isaconst{omv}\ \isavar{c}\ {\isasymlbrace}\isaconst{True}{\isasymrbrace}}%
}%
\SNIP{Goedel-Paper-Statements:thm:adr-consistent-0}{%
\isakeywordONE{thm}\isamarkupfalse%
\ adr{\isacharunderscore}{\kern0pt}consistent%
\isanewline%
}%
\SNIP{Goedel-Paper-Statements:thm:adr-consistent-1}{%
\isa{{\isasymnot}\ \isaconst{prf{\isacharunderscore}{\kern0pt}derives{\isachardot}{\kern0pt}entails}\ \isaconst{adr{\isacharunderscore}{\kern0pt}derives}\ {\isasymemptyset}\ {\isacharparenleft}{\kern0pt}\isactrlbold {\isadigit{0}}\ \isactrlbold {\isacharequal}{\kern0pt}\ \isactrlbold {\isadigit{1}}{\isacharparenright}{\kern0pt}}%
}%
\SNIP{Goedel-Paper-Statements:thm:adrU-def-0}{%
\isakeywordONE{thm}\isamarkupfalse%
\ adrU{\isacharunderscore}{\kern0pt}def%
\isanewline%
}%
\SNIP{Goedel-Paper-Statements:thm:adrU-def-1}{%
\isa{\isaconst{adrU}\ {\isasymequiv}\ \isactrlbold {\isasymlambda}\isaconst{adrBody}}%
}%
\SNIP{Goedel-Paper-Statements:thm:adrF-def-0}{%
\isakeywordONE{thm}\isamarkupfalse%
\ adrF{\isacharunderscore}{\kern0pt}def%
\isanewline%
}%
\SNIP{Goedel-Paper-Statements:thm:adrF-def-1}{%
\isa{\isaconst{adrF}\ {\isasymequiv}\ \isaconst{ifIdx}\ \isaconst{adrCIdx}\ \isaconst{adrExfIdx}\ \isaconst{adrGoodIdx}}%
}%
\SNIP{Goedel-Paper-Statements:thm:adr-cert-0}{%
\isakeywordONE{thm}\isamarkupfalse%
\ adr{\isacharunderscore}{\kern0pt}cert%
\isanewline%
}%
\SNIP{Goedel-Paper-Statements:thm:adr-cert-1}{%
\isa{\isavar{{\isasymGamma}}\ {\isasymturnstile}\ \isavar{v}\ \isactrlbold N\ {\isasymLongrightarrow}\ \isavar{{\isasymGamma}}\ {\isasymturnstile}\ \isaconst{iChk}\ {\isacharparenleft}{\kern0pt}\isaconst{cmp}\ \isaconst{adrF}\ \isactrlbold {\isasymcdot}\ \isavar{v}\ \isactrlbold {\isasymcdot}\ \isavar{v}{\isacharparenright}{\kern0pt}\ {\isacharparenleft}{\kern0pt}\isaconst{adrInst}\ \isaconst{adrU}\ \isavar{v}{\isacharparenright}{\kern0pt}}%
}%
\SNIP{Goedel-Paper-Statements:thm:adr-warrant-uC-0}{%
\isakeywordONE{thm}\isamarkupfalse%
\ adr{\isacharunderscore}{\kern0pt}warrant{\isacharunderscore}{\kern0pt}uC%
\isanewline%
}%
\SNIP{Goedel-Paper-Statements:thm:adr-warrant-uC-1}{%
\isa{{\isasymemptyset}\ {\isasymturnstile}\ \isaconst{adrW}\ \isaconst{adrF}\ \isaconst{adrU}}%
}%
\SNIP{Goedel-Paper-Statements:thm:adr-U-ConRGA-0}{%
\isakeywordONE{thm}\isamarkupfalse%
\ adr{\isacharunderscore}{\kern0pt}U{\isacharunderscore}{\kern0pt}ConRGA%
\isanewline%
}%
\SNIP{Goedel-Paper-Statements:thm:adr-U-ConRGA-1}{%
\isa{\isactrlbold A\ \isactrlbold {\isasymcdot}\ \isaconst{adrU}\ {\isasymtriangleright}\ {\isasymemptyset}\ {\isasymturnstile}\ \isaconst{ConRGA}}%
}%
\SNIP{Goedel-Paper-Statements:thm:adr-ConRGA-U-0}{%
\isakeywordONE{thm}\isamarkupfalse%
\ adr{\isacharunderscore}{\kern0pt}ConRGA{\isacharunderscore}{\kern0pt}U%
\isanewline%
}%
\SNIP{Goedel-Paper-Statements:thm:adr-ConRGA-U-1}{%
\isa{\isaconst{ConRGA}\ {\isasymtriangleright}\ {\isasymemptyset}\ {\isasymturnstile}\ \isactrlbold A\ \isactrlbold {\isasymcdot}\ \isaconst{adrU}}%
}%
\SNIP{Goedel-Paper-Statements:thm:adr-ConRGA-0}{%
\isakeywordONE{thm}\isamarkupfalse%
\ adr{\isacharunderscore}{\kern0pt}ConRGA%
\isanewline%
}%
\SNIP{Goedel-Paper-Statements:thm:adr-ConRGA-1}{%
\isa{\isaconst{prf{\isacharunderscore}{\kern0pt}derives{\isachardot}{\kern0pt}entails}\ \isaconst{adr{\isacharunderscore}{\kern0pt}derives}\ {\isasymemptyset}\ \isaconst{ConRGA}}%
}%
\SNIP{Goedel-Paper-Statements:thm:adr-separation-0}{%
\isakeywordONE{thm}\isamarkupfalse%
\ adr{\isacharunderscore}{\kern0pt}separation%
\isanewline%
}%
\SNIP{Goedel-Paper-Statements:thm:adr-separation-1}{%
\isa{{\isasymnot}\ {\isasymemptyset}\ {\isasymturnstile}\ \isaconst{ConRGA}\ {\isasymand}\ \isaconst{prf{\isacharunderscore}{\kern0pt}derives{\isachardot}{\kern0pt}entails}\ \isaconst{adr{\isacharunderscore}{\kern0pt}derives}\ {\isasymemptyset}\ \isaconst{ConRGA}\ {\isasymand}\ {\isacharparenleft}{\kern0pt}{\isasymforall}\isabound{{\isasymGamma}}\ \isabound{c}{\isachardot}{\kern0pt}\ \isabound{{\isasymGamma}}\ {\isasymturnstile}\ \isabound{c}\ {\isasymlongrightarrow}\ \isaconst{prf{\isacharunderscore}{\kern0pt}derives{\isachardot}{\kern0pt}entails}\ \isaconst{adr{\isacharunderscore}{\kern0pt}derives}\ \isabound{{\isasymGamma}}\ \isabound{c}{\isacharparenright}{\kern0pt}\ {\isasymand}\ {\isasymnot}\ \isaconst{prf{\isacharunderscore}{\kern0pt}derives{\isachardot}{\kern0pt}entails}\ \isaconst{adr{\isacharunderscore}{\kern0pt}derives}\ {\isasymemptyset}\ {\isacharparenleft}{\kern0pt}\isactrlbold {\isadigit{0}}\ \isactrlbold {\isacharequal}{\kern0pt}\ \isactrlbold {\isadigit{1}}{\isacharparenright}{\kern0pt}}%
}%
\SNIP{Goedel-Paper-Statements:thm:adi-embed-0}{%
\isakeywordONE{thm}\isamarkupfalse%
\ adi{\isacharunderscore}{\kern0pt}embed%
\isanewline%
}%
\SNIP{Goedel-Paper-Statements:thm:adi-embed-1}{%
\isa{\isaconst{prf{\isacharunderscore}{\kern0pt}derives{\isachardot}{\kern0pt}entails}\ \isaconst{gi{\isachardot}{\kern0pt}derives}\ \isavar{{\isasymGamma}}\ \isavar{c}\ {\isasymLongrightarrow}\ \isaconst{prf{\isacharunderscore}{\kern0pt}derives{\isachardot}{\kern0pt}entails}\ \isaconst{adi{\isacharunderscore}{\kern0pt}derives}\ \isavar{{\isasymGamma}}\ \isavar{c}}%
}%
\SNIP{Goedel-Paper-Statements:thm:adi-omega-sound-0}{%
\isakeywordONE{thm}\isamarkupfalse%
\ adi{\isacharunderscore}{\kern0pt}omega{\isacharunderscore}{\kern0pt}sound%
\isanewline%
}%
\SNIP{Goedel-Paper-Statements:thm:adi-omega-sound-1}{%
\isa{{\isasymlbrakk}\isaconst{prf{\isacharunderscore}{\kern0pt}derives{\isachardot}{\kern0pt}entails}\ \isaconst{adi{\isacharunderscore}{\kern0pt}derives}\ \isavar{{\isasymGamma}}\ \isavar{c}{\isacharsemicolon}{\kern0pt}\ \isaconst{wfc}\ \isavar{{\isasymGamma}}{\isacharsemicolon}{\kern0pt}\ {\isasymAnd}\isabound{g}{\isachardot}{\kern0pt}\ \isabound{g}\ {\isacharbar}{\kern0pt}{\isasymin}{\isacharbar}{\kern0pt}\ \isavar{{\isasymGamma}}\ {\isasymLongrightarrow}\ \isaconst{omv}\ \isabound{g}\ {\isasymlbrace}\isaconst{True}{\isasymrbrace}{\isasymrbrakk}\ {\isasymLongrightarrow}\ \isaconst{omv}\ \isavar{c}\ {\isasymlbrace}\isaconst{True}{\isasymrbrace}}%
}%
\SNIP{Goedel-Paper-Statements:thm:adi-warrant-uC-0}{%
\isakeywordONE{thm}\isamarkupfalse%
\ adi{\isacharunderscore}{\kern0pt}warrant{\isacharunderscore}{\kern0pt}uC%
\isanewline%
}%
\SNIP{Goedel-Paper-Statements:thm:adi-warrant-uC-1}{%
\isa{{\isasymemptyset}\ {\isasymturnstile}\ \isaconst{adiW}\ \isaconst{adiF}\ \isaconst{adiU}}%
}%
\SNIP{Goedel-Paper-Statements:thm:adi-ConRGA-of-ConIGA-0}{%
\isakeywordONE{thm}\isamarkupfalse%
\ adi{\isacharunderscore}{\kern0pt}ConRGA{\isacharunderscore}{\kern0pt}of{\isacharunderscore}{\kern0pt}ConIGA%
\isanewline%
}%
\SNIP{Goedel-Paper-Statements:thm:adi-ConRGA-of-ConIGA-1}{%
\isa{\isaconst{ConIGA}\ {\isasymtriangleright}\ {\isasymemptyset}\ {\isasymturnstile}\ \isaconst{ConRGA}}%
}%
\SNIP{Goedel-Paper-Statements:thm:adi-ConIGA-0}{%
\isakeywordONE{thm}\isamarkupfalse%
\ adi{\isacharunderscore}{\kern0pt}ConIGA%
\isanewline%
}%
\SNIP{Goedel-Paper-Statements:thm:adi-ConIGA-1}{%
\isa{\isaconst{prf{\isacharunderscore}{\kern0pt}derives{\isachardot}{\kern0pt}entails}\ \isaconst{adi{\isacharunderscore}{\kern0pt}derives}\ {\isasymemptyset}\ \isaconst{ConIGA}}%
}%
\SNIP{Goedel-Paper-Statements:thm:adi-ConRGA-0}{%
\isakeywordONE{thm}\isamarkupfalse%
\ adi{\isacharunderscore}{\kern0pt}ConRGA%
\isanewline%
}%
\SNIP{Goedel-Paper-Statements:thm:adi-ConRGA-1}{%
\isa{\isaconst{prf{\isacharunderscore}{\kern0pt}derives{\isachardot}{\kern0pt}entails}\ \isaconst{adi{\isacharunderscore}{\kern0pt}derives}\ {\isasymemptyset}\ \isaconst{ConRGA}}%
}%
\SNIP{Goedel-Paper-Statements:thm:adi-separation-RGA-0}{%
\isakeywordONE{thm}\isamarkupfalse%
\ adi{\isacharunderscore}{\kern0pt}separation{\isacharunderscore}{\kern0pt}RGA%
\isanewline%
}%
\SNIP{Goedel-Paper-Statements:thm:adi-separation-RGA-1}{%
\isa{{\isasymnot}\ {\isasymemptyset}\ {\isasymturnstile}\ \isaconst{ConRGA}\ {\isasymand}\ \isaconst{prf{\isacharunderscore}{\kern0pt}derives{\isachardot}{\kern0pt}entails}\ \isaconst{adi{\isacharunderscore}{\kern0pt}derives}\ {\isasymemptyset}\ \isaconst{ConRGA}\ {\isasymand}\ {\isacharparenleft}{\kern0pt}{\isasymforall}\isabound{{\isasymGamma}}\ \isabound{c}{\isachardot}{\kern0pt}\ \isabound{{\isasymGamma}}\ {\isasymturnstile}\ \isabound{c}\ {\isasymlongrightarrow}\ \isaconst{prf{\isacharunderscore}{\kern0pt}derives{\isachardot}{\kern0pt}entails}\ \isaconst{adi{\isacharunderscore}{\kern0pt}derives}\ \isabound{{\isasymGamma}}\ \isabound{c}{\isacharparenright}{\kern0pt}\ {\isasymand}\ {\isacharparenleft}{\kern0pt}{\isasymforall}\isabound{{\isasymGamma}}\ \isabound{c}{\isachardot}{\kern0pt}\ \isaconst{prf{\isacharunderscore}{\kern0pt}derives{\isachardot}{\kern0pt}entails}\ \isaconst{gi{\isachardot}{\kern0pt}derives}\ \isabound{{\isasymGamma}}\ \isabound{c}\ {\isasymlongrightarrow}\ \isaconst{prf{\isacharunderscore}{\kern0pt}derives{\isachardot}{\kern0pt}entails}\ \isaconst{adi{\isacharunderscore}{\kern0pt}derives}\ \isabound{{\isasymGamma}}\ \isabound{c}{\isacharparenright}{\kern0pt}\ {\isasymand}\ {\isasymnot}\ \isaconst{prf{\isacharunderscore}{\kern0pt}derives{\isachardot}{\kern0pt}entails}\ \isaconst{adi{\isacharunderscore}{\kern0pt}derives}\ {\isasymemptyset}\ {\isacharparenleft}{\kern0pt}\isactrlbold {\isadigit{0}}\ \isactrlbold {\isacharequal}{\kern0pt}\ \isactrlbold {\isadigit{1}}{\isacharparenright}{\kern0pt}}%
}%
\SNIP{Goedel-Paper-Statements:thm:adi-separation-of-goedel2-0}{%
\isakeywordONE{thm}\isamarkupfalse%
\ adi{\isacharunderscore}{\kern0pt}separation{\isacharunderscore}{\kern0pt}of{\isacharunderscore}{\kern0pt}goedel{\isadigit{2}}%
\isanewline%
}%
\SNIP{Goedel-Paper-Statements:thm:adi-separation-of-goedel2-1}{%
\isa{{\isasymnot}\ \isaconst{prf{\isacharunderscore}{\kern0pt}derives{\isachardot}{\kern0pt}entails}\ \isaconst{gi{\isachardot}{\kern0pt}derives}\ {\isasymemptyset}\ \isaconst{ConIGA}\ {\isasymLongrightarrow}\ {\isasymnot}\ \isaconst{prf{\isacharunderscore}{\kern0pt}derives{\isachardot}{\kern0pt}entails}\ \isaconst{gi{\isachardot}{\kern0pt}derives}\ {\isasymemptyset}\ \isaconst{ConIGA}\ {\isasymand}\ \isaconst{prf{\isacharunderscore}{\kern0pt}derives{\isachardot}{\kern0pt}entails}\ \isaconst{adi{\isacharunderscore}{\kern0pt}derives}\ {\isasymemptyset}\ \isaconst{ConIGA}\ {\isasymand}\ {\isacharparenleft}{\kern0pt}{\isasymforall}\isabound{{\isasymGamma}}\ \isabound{c}{\isachardot}{\kern0pt}\ \isaconst{prf{\isacharunderscore}{\kern0pt}derives{\isachardot}{\kern0pt}entails}\ \isaconst{gi{\isachardot}{\kern0pt}derives}\ \isabound{{\isasymGamma}}\ \isabound{c}\ {\isasymlongrightarrow}\ \isaconst{prf{\isacharunderscore}{\kern0pt}derives{\isachardot}{\kern0pt}entails}\ \isaconst{adi{\isacharunderscore}{\kern0pt}derives}\ \isabound{{\isasymGamma}}\ \isabound{c}{\isacharparenright}{\kern0pt}\ {\isasymand}\ {\isasymnot}\ \isaconst{prf{\isacharunderscore}{\kern0pt}derives{\isachardot}{\kern0pt}entails}\ \isaconst{adi{\isacharunderscore}{\kern0pt}derives}\ {\isasymemptyset}\ {\isacharparenleft}{\kern0pt}\isactrlbold {\isadigit{0}}\ \isactrlbold {\isacharequal}{\kern0pt}\ \isactrlbold {\isadigit{1}}{\isacharparenright}{\kern0pt}}%
}%
\SNIP{Goedel-Paper-Statements:thm:iga-not-motive-0}{%
\isakeywordONE{thm}\isamarkupfalse%
\ iga{\isacharunderscore}{\kern0pt}not{\isacharunderscore}{\kern0pt}motive%
\isanewline%
}%
\SNIP{Goedel-Paper-Statements:thm:iga-not-motive-1}{%
\isa{{\isasymlbrakk}\isaconst{wfc}\ \isavar{{\isasymGamma}}{\isacharsemicolon}{\kern0pt}\ {\isasymAnd}\isabound{g}{\isachardot}{\kern0pt}\ \isabound{g}\ {\isacharbar}{\kern0pt}{\isasymin}{\isacharbar}{\kern0pt}\ \isavar{{\isasymGamma}}\ {\isasymLongrightarrow}\ \isaconst{omv}\ \isabound{g}\ {\isasymlbrace}\isaconst{True}{\isasymrbrace}{\isasymrbrakk}\ {\isasymLongrightarrow}\ {\isasymnot}\ \isaconst{prf{\isacharunderscore}{\kern0pt}derives{\isachardot}{\kern0pt}entails}\ \isaconst{gi{\isachardot}{\kern0pt}derives}\ \isavar{{\isasymGamma}}\ {\isacharparenleft}{\kern0pt}\isaconst{fMRb}\ {\isacharparenleft}{\kern0pt}\isactrlbold n{\isasymlbrace}\isaconst{phP}{\isasymrbrace}{\isacharparenright}{\kern0pt}\ {\isacharparenleft}{\kern0pt}\isaconst{iPoRk}\ {\isacharparenleft}{\kern0pt}\isaconst{cmp}\ \isaconst{pomaxIdx}\ \isactrlbold {\isasymcdot}\ \isactrlbold {\isadigit{0}}\ \isactrlbold {\isasymcdot}\ \isaconst{trmPredR}\ {\isacharparenleft}{\kern0pt}\isaconst{iPoHt}\ \isaconst{kidRdIdx}\ {\isacharparenleft}{\kern0pt}\isactrlbold n{\isasymlbrace}\isaconst{phP}{\isasymrbrace}{\isacharparenright}{\kern0pt}\ {\isacharparenleft}{\kern0pt}\isaconst{iPoRk}\ {\isacharparenleft}{\kern0pt}\isaconst{cmp}\ \isaconst{pomaxIdx}\ \isactrlbold {\isasymcdot}\ \isactrlbold {\isadigit{0}}\ \isactrlbold {\isasymcdot}\ \isaconst{trmPredR}\ {\isacharparenleft}{\kern0pt}\isaconst{iPoHt}\ \isaconst{kidRdIdx}\ {\isacharparenleft}{\kern0pt}\isactrlbold n{\isasymlbrace}\isaconst{phP}{\isasymrbrace}{\isacharparenright}{\kern0pt}\ {\isacharparenleft}{\kern0pt}\isaconst{cmp}\ \isaconst{pomaxIdx}\ \isactrlbold {\isasymcdot}\ \isactrlbold {\isadigit{0}}\ \isactrlbold {\isasymcdot}\ {\isacharparenleft}{\kern0pt}\isaconst{cmp}\ \isaconst{length{\isacharunderscore}{\kern0pt}encIdx}\ \isactrlbold {\isasymcdot}\ {\isacharparenleft}{\kern0pt}\isactrlbold n{\isasymlbrace}\isaconst{phP}{\isasymrbrace}{\isacharparenright}{\kern0pt}\ \isactrlbold {\isasymcdot}\ {\isacharparenleft}{\kern0pt}\isactrlbold n{\isasymlbrace}\isaconst{phP}{\isasymrbrace}{\isacharparenright}{\kern0pt}{\isacharparenright}{\kern0pt}{\isacharparenright}{\kern0pt}{\isacharparenright}{\kern0pt}{\isacharparenright}{\kern0pt}\ {\isacharparenleft}{\kern0pt}\isaconst{wBa}\ {\isacharparenleft}{\kern0pt}\isaconst{Suc}\ {\isadigit{0}}{\isacharparenright}{\kern0pt}\ {\isacharparenleft}{\kern0pt}\isaconst{iPoBase}\ \isaconst{kidRdIdx}\ {\isacharparenleft}{\kern0pt}\isactrlbold n{\isasymlbrace}\isaconst{phP}{\isasymrbrace}{\isacharparenright}{\kern0pt}\ {\isacharparenleft}{\kern0pt}\isaconst{cmp}\ \isaconst{pomaxIdx}\ \isactrlbold {\isasymcdot}\ \isactrlbold {\isadigit{0}}\ \isactrlbold {\isasymcdot}\ {\isacharparenleft}{\kern0pt}\isaconst{cmp}\ \isaconst{length{\isacharunderscore}{\kern0pt}encIdx}\ \isactrlbold {\isasymcdot}\ {\isacharparenleft}{\kern0pt}\isactrlbold n{\isasymlbrace}\isaconst{phP}{\isasymrbrace}{\isacharparenright}{\kern0pt}\ \isactrlbold {\isasymcdot}\ {\isacharparenleft}{\kern0pt}\isactrlbold n{\isasymlbrace}\isaconst{phP}{\isasymrbrace}{\isacharparenright}{\kern0pt}{\isacharparenright}{\kern0pt}{\isacharparenright}{\kern0pt}{\isacharparenright}{\kern0pt}\ {\isacharparenleft}{\kern0pt}\isaconst{iPoHt}\ \isaconst{kidRdIdx}\ {\isacharparenleft}{\kern0pt}\isactrlbold n{\isasymlbrace}\isaconst{phP}{\isasymrbrace}{\isacharparenright}{\kern0pt}\ {\isacharparenleft}{\kern0pt}\isaconst{cmp}\ \isaconst{pomaxIdx}\ \isactrlbold {\isasymcdot}\ \isactrlbold {\isadigit{0}}\ \isactrlbold {\isasymcdot}\ {\isacharparenleft}{\kern0pt}\isaconst{cmp}\ \isaconst{length{\isacharunderscore}{\kern0pt}encIdx}\ \isactrlbold {\isasymcdot}\ {\isacharparenleft}{\kern0pt}\isactrlbold n{\isasymlbrace}\isaconst{phP}{\isasymrbrace}{\isacharparenright}{\kern0pt}\ \isactrlbold {\isasymcdot}\ {\isacharparenleft}{\kern0pt}\isactrlbold n{\isasymlbrace}\isaconst{phP}{\isasymrbrace}{\isacharparenright}{\kern0pt}{\isacharparenright}{\kern0pt}{\isacharparenright}{\kern0pt}{\isacharparenright}{\kern0pt}{\isacharparenright}{\kern0pt}{\isacharparenright}{\kern0pt}{\isacharparenright}{\kern0pt}{\isacharparenright}{\kern0pt}\ {\isacharparenleft}{\kern0pt}\isaconst{wBa}\ {\isadigit{0}}\ {\isacharparenleft}{\kern0pt}\isaconst{iPoBase}\ \isaconst{kidRdIdx}\ {\isacharparenleft}{\kern0pt}\isactrlbold n{\isasymlbrace}\isaconst{phP}{\isasymrbrace}{\isacharparenright}{\kern0pt}\ {\isacharparenleft}{\kern0pt}\isaconst{iPoRk}\ {\isacharparenleft}{\kern0pt}\isaconst{cmp}\ \isaconst{pomaxIdx}\ \isactrlbold {\isasymcdot}\ \isactrlbold {\isadigit{0}}\ \isactrlbold {\isasymcdot}\ \isaconst{trmPredR}\ {\isacharparenleft}{\kern0pt}\isaconst{iPoHt}\ \isaconst{kidRdIdx}\ {\isacharparenleft}{\kern0pt}\isactrlbold n{\isasymlbrace}\isaconst{phP}{\isasymrbrace}{\isacharparenright}{\kern0pt}\ {\isacharparenleft}{\kern0pt}\isaconst{cmp}\ \isaconst{pomaxIdx}\ \isactrlbold {\isasymcdot}\ \isactrlbold {\isadigit{0}}\ \isactrlbold {\isasymcdot}\ {\isacharparenleft}{\kern0pt}\isaconst{cmp}\ \isaconst{length{\isacharunderscore}{\kern0pt}encIdx}\ \isactrlbold {\isasymcdot}\ {\isacharparenleft}{\kern0pt}\isactrlbold n{\isasymlbrace}\isaconst{phP}{\isasymrbrace}{\isacharparenright}{\kern0pt}\ \isactrlbold {\isasymcdot}\ {\isacharparenleft}{\kern0pt}\isactrlbold n{\isasymlbrace}\isaconst{phP}{\isasymrbrace}{\isacharparenright}{\kern0pt}{\isacharparenright}{\kern0pt}{\isacharparenright}{\kern0pt}{\isacharparenright}{\kern0pt}{\isacharparenright}{\kern0pt}\ {\isacharparenleft}{\kern0pt}\isaconst{wBa}\ {\isacharparenleft}{\kern0pt}\isaconst{Suc}\ {\isadigit{0}}{\isacharparenright}{\kern0pt}\ {\isacharparenleft}{\kern0pt}\isaconst{iPoBase}\ \isaconst{kidRdIdx}\ {\isacharparenleft}{\kern0pt}\isactrlbold n{\isasymlbrace}\isaconst{phP}{\isasymrbrace}{\isacharparenright}{\kern0pt}\ {\isacharparenleft}{\kern0pt}\isaconst{cmp}\ \isaconst{pomaxIdx}\ \isactrlbold {\isasymcdot}\ \isactrlbold {\isadigit{0}}\ \isactrlbold {\isasymcdot}\ {\isacharparenleft}{\kern0pt}\isaconst{cmp}\ \isaconst{length{\isacharunderscore}{\kern0pt}encIdx}\ \isactrlbold {\isasymcdot}\ {\isacharparenleft}{\kern0pt}\isactrlbold n{\isasymlbrace}\isaconst{phP}{\isasymrbrace}{\isacharparenright}{\kern0pt}\ \isactrlbold {\isasymcdot}\ {\isacharparenleft}{\kern0pt}\isactrlbold n{\isasymlbrace}\isaconst{phP}{\isasymrbrace}{\isacharparenright}{\kern0pt}{\isacharparenright}{\kern0pt}{\isacharparenright}{\kern0pt}{\isacharparenright}{\kern0pt}\ {\isacharparenleft}{\kern0pt}\isaconst{iPoHt}\ \isaconst{kidRdIdx}\ {\isacharparenleft}{\kern0pt}\isactrlbold n{\isasymlbrace}\isaconst{phP}{\isasymrbrace}{\isacharparenright}{\kern0pt}\ {\isacharparenleft}{\kern0pt}\isaconst{cmp}\ \isaconst{pomaxIdx}\ \isactrlbold {\isasymcdot}\ \isactrlbold {\isadigit{0}}\ \isactrlbold {\isasymcdot}\ {\isacharparenleft}{\kern0pt}\isaconst{cmp}\ \isaconst{length{\isacharunderscore}{\kern0pt}encIdx}\ \isactrlbold {\isasymcdot}\ {\isacharparenleft}{\kern0pt}\isactrlbold n{\isasymlbrace}\isaconst{phP}{\isasymrbrace}{\isacharparenright}{\kern0pt}\ \isactrlbold {\isasymcdot}\ {\isacharparenleft}{\kern0pt}\isactrlbold n{\isasymlbrace}\isaconst{phP}{\isasymrbrace}{\isacharparenright}{\kern0pt}{\isacharparenright}{\kern0pt}{\isacharparenright}{\kern0pt}{\isacharparenright}{\kern0pt}{\isacharparenright}{\kern0pt}{\isacharparenright}{\kern0pt}{\isacharparenright}{\kern0pt}\ {\isacharparenleft}{\kern0pt}\isaconst{iPoHt}\ \isaconst{kidRdIdx}\ {\isacharparenleft}{\kern0pt}\isactrlbold n{\isasymlbrace}\isaconst{phP}{\isasymrbrace}{\isacharparenright}{\kern0pt}\ {\isacharparenleft}{\kern0pt}\isaconst{iPoRk}\ {\isacharparenleft}{\kern0pt}\isaconst{cmp}\ \isaconst{pomaxIdx}\ \isactrlbold {\isasymcdot}\ \isactrlbold {\isadigit{0}}\ \isactrlbold {\isasymcdot}\ \isaconst{trmPredR}\ {\isacharparenleft}{\kern0pt}\isaconst{iPoHt}\ \isaconst{kidRdIdx}\ {\isacharparenleft}{\kern0pt}\isactrlbold n{\isasymlbrace}\isaconst{phP}{\isasymrbrace}{\isacharparenright}{\kern0pt}\ {\isacharparenleft}{\kern0pt}\isaconst{cmp}\ \isaconst{pomaxIdx}\ \isactrlbold {\isasymcdot}\ \isactrlbold {\isadigit{0}}\ \isactrlbold {\isasymcdot}\ {\isacharparenleft}{\kern0pt}\isaconst{cmp}\ \isaconst{length{\isacharunderscore}{\kern0pt}encIdx}\ \isactrlbold {\isasymcdot}\ {\isacharparenleft}{\kern0pt}\isactrlbold n{\isasymlbrace}\isaconst{phP}{\isasymrbrace}{\isacharparenright}{\kern0pt}\ \isactrlbold {\isasymcdot}\ {\isacharparenleft}{\kern0pt}\isactrlbold n{\isasymlbrace}\isaconst{phP}{\isasymrbrace}{\isacharparenright}{\kern0pt}{\isacharparenright}{\kern0pt}{\isacharparenright}{\kern0pt}{\isacharparenright}{\kern0pt}{\isacharparenright}{\kern0pt}\ {\isacharparenleft}{\kern0pt}\isaconst{wBa}\ {\isacharparenleft}{\kern0pt}\isaconst{Suc}\ {\isadigit{0}}{\isacharparenright}{\kern0pt}\ {\isacharparenleft}{\kern0pt}\isaconst{iPoBase}\ \isaconst{kidRdIdx}\ {\isacharparenleft}{\kern0pt}\isactrlbold n{\isasymlbrace}\isaconst{phP}{\isasymrbrace}{\isacharparenright}{\kern0pt}\ {\isacharparenleft}{\kern0pt}\isaconst{cmp}\ \isaconst{pomaxIdx}\ \isactrlbold {\isasymcdot}\ \isactrlbold {\isadigit{0}}\ \isactrlbold {\isasymcdot}\ {\isacharparenleft}{\kern0pt}\isaconst{cmp}\ \isaconst{length{\isacharunderscore}{\kern0pt}encIdx}\ \isactrlbold {\isasymcdot}\ {\isacharparenleft}{\kern0pt}\isactrlbold n{\isasymlbrace}\isaconst{phP}{\isasymrbrace}{\isacharparenright}{\kern0pt}\ \isactrlbold {\isasymcdot}\ {\isacharparenleft}{\kern0pt}\isactrlbold n{\isasymlbrace}\isaconst{phP}{\isasymrbrace}{\isacharparenright}{\kern0pt}{\isacharparenright}{\kern0pt}{\isacharparenright}{\kern0pt}{\isacharparenright}{\kern0pt}\ {\isacharparenleft}{\kern0pt}\isaconst{iPoHt}\ \isaconst{kidRdIdx}\ {\isacharparenleft}{\kern0pt}\isactrlbold n{\isasymlbrace}\isaconst{phP}{\isasymrbrace}{\isacharparenright}{\kern0pt}\ {\isacharparenleft}{\kern0pt}\isaconst{cmp}\ \isaconst{pomaxIdx}\ \isactrlbold {\isasymcdot}\ \isactrlbold {\isadigit{0}}\ \isactrlbold {\isasymcdot}\ {\isacharparenleft}{\kern0pt}\isaconst{cmp}\ \isaconst{length{\isacharunderscore}{\kern0pt}encIdx}\ \isactrlbold {\isasymcdot}\ {\isacharparenleft}{\kern0pt}\isactrlbold n{\isasymlbrace}\isaconst{phP}{\isasymrbrace}{\isacharparenright}{\kern0pt}\ \isactrlbold {\isasymcdot}\ {\isacharparenleft}{\kern0pt}\isactrlbold n{\isasymlbrace}\isaconst{phP}{\isasymrbrace}{\isacharparenright}{\kern0pt}{\isacharparenright}{\kern0pt}{\isacharparenright}{\kern0pt}{\isacharparenright}{\kern0pt}{\isacharparenright}{\kern0pt}{\isacharparenright}{\kern0pt}{\isacharparenright}{\kern0pt}{\isacharparenright}{\kern0pt}{\isacharparenright}{\kern0pt}\ {\isacharparenleft}{\kern0pt}\isaconst{lstCons}\ \isactrlbold {\isadigit{0}}\ {\isacharparenleft}{\kern0pt}\isaconst{lstCons}\ \isactrlbold {\isadigit{0}}\ \isactrlbold {\isadigit{0}}{\isacharparenright}{\kern0pt}{\isacharparenright}{\kern0pt}\ {\isacharparenleft}{\kern0pt}\isaconst{lstCons}\ {\isacharparenleft}{\kern0pt}\isactrlbold n{\isadigit{2}}{\isacharparenright}{\kern0pt}\ \isactrlbold {\isadigit{0}}{\isacharparenright}{\kern0pt}{\isacharparenright}{\kern0pt}\isasep\isanewline%
}%
\SNIP{Goedel-Paper-Statements:thm:iga-not-motive-2}{%
{\isasymlbrakk}\isaconst{wfc}\ \isavar{{\isasymGamma}}{\isacharsemicolon}{\kern0pt}\ {\isasymAnd}\isabound{g}{\isachardot}{\kern0pt}\ \isabound{g}\ {\isacharbar}{\kern0pt}{\isasymin}{\isacharbar}{\kern0pt}\ \isavar{{\isasymGamma}}\ {\isasymLongrightarrow}\ \isaconst{omv}\ \isabound{g}\ {\isasymlbrace}\isaconst{True}{\isasymrbrace}{\isasymrbrakk}\ {\isasymLongrightarrow}\ {\isasymnot}\ \isaconst{prf{\isacharunderscore}{\kern0pt}derives{\isachardot}{\kern0pt}entails}\ \isaconst{gi{\isachardot}{\kern0pt}derives}\ \isavar{{\isasymGamma}}\ {\isacharparenleft}{\kern0pt}\isaconst{fMot}\ {\isacharparenleft}{\kern0pt}\isactrlbold n{\isasymlbrace}\isaconst{phP}{\isasymrbrace}{\isacharparenright}{\kern0pt}\ \isactrlbold {\isasymcdot}\ \isaconst{iPoRk}\ {\isacharparenleft}{\kern0pt}\isaconst{cmp}\ \isaconst{pomaxIdx}\ \isactrlbold {\isasymcdot}\ \isactrlbold {\isadigit{0}}\ \isactrlbold {\isasymcdot}\ \isaconst{trmPredR}\ {\isacharparenleft}{\kern0pt}\isaconst{iPoHt}\ \isaconst{kidRdIdx}\ {\isacharparenleft}{\kern0pt}\isactrlbold n{\isasymlbrace}\isaconst{phP}{\isasymrbrace}{\isacharparenright}{\kern0pt}\ {\isacharparenleft}{\kern0pt}\isaconst{iPoRk}\ {\isacharparenleft}{\kern0pt}\isaconst{cmp}\ \isaconst{pomaxIdx}\ \isactrlbold {\isasymcdot}\ \isactrlbold {\isadigit{0}}\ \isactrlbold {\isasymcdot}\ \isaconst{trmPredR}\ {\isacharparenleft}{\kern0pt}\isaconst{iPoHt}\ \isaconst{kidRdIdx}\ {\isacharparenleft}{\kern0pt}\isactrlbold n{\isasymlbrace}\isaconst{phP}{\isasymrbrace}{\isacharparenright}{\kern0pt}\ {\isacharparenleft}{\kern0pt}\isaconst{cmp}\ \isaconst{pomaxIdx}\ \isactrlbold {\isasymcdot}\ \isactrlbold {\isadigit{0}}\ \isactrlbold {\isasymcdot}\ {\isacharparenleft}{\kern0pt}\isaconst{cmp}\ \isaconst{length{\isacharunderscore}{\kern0pt}encIdx}\ \isactrlbold {\isasymcdot}\ {\isacharparenleft}{\kern0pt}\isactrlbold n{\isasymlbrace}\isaconst{phP}{\isasymrbrace}{\isacharparenright}{\kern0pt}\ \isactrlbold {\isasymcdot}\ {\isacharparenleft}{\kern0pt}\isactrlbold n{\isasymlbrace}\isaconst{phP}{\isasymrbrace}{\isacharparenright}{\kern0pt}{\isacharparenright}{\kern0pt}{\isacharparenright}{\kern0pt}{\isacharparenright}{\kern0pt}{\isacharparenright}{\kern0pt}\ {\isacharparenleft}{\kern0pt}\isaconst{wBa}\ {\isacharparenleft}{\kern0pt}\isaconst{Suc}\ {\isadigit{0}}{\isacharparenright}{\kern0pt}\ {\isacharparenleft}{\kern0pt}\isaconst{iPoBase}\ \isaconst{kidRdIdx}\ {\isacharparenleft}{\kern0pt}\isactrlbold n{\isasymlbrace}\isaconst{phP}{\isasymrbrace}{\isacharparenright}{\kern0pt}\ {\isacharparenleft}{\kern0pt}\isaconst{cmp}\ \isaconst{pomaxIdx}\ \isactrlbold {\isasymcdot}\ \isactrlbold {\isadigit{0}}\ \isactrlbold {\isasymcdot}\ {\isacharparenleft}{\kern0pt}\isaconst{cmp}\ \isaconst{length{\isacharunderscore}{\kern0pt}encIdx}\ \isactrlbold {\isasymcdot}\ {\isacharparenleft}{\kern0pt}\isactrlbold n{\isasymlbrace}\isaconst{phP}{\isasymrbrace}{\isacharparenright}{\kern0pt}\ \isactrlbold {\isasymcdot}\ {\isacharparenleft}{\kern0pt}\isactrlbold n{\isasymlbrace}\isaconst{phP}{\isasymrbrace}{\isacharparenright}{\kern0pt}{\isacharparenright}{\kern0pt}{\isacharparenright}{\kern0pt}{\isacharparenright}{\kern0pt}\ {\isacharparenleft}{\kern0pt}\isaconst{iPoHt}\ \isaconst{kidRdIdx}\ {\isacharparenleft}{\kern0pt}\isactrlbold n{\isasymlbrace}\isaconst{phP}{\isasymrbrace}{\isacharparenright}{\kern0pt}\ {\isacharparenleft}{\kern0pt}\isaconst{cmp}\ \isaconst{pomaxIdx}\ \isactrlbold {\isasymcdot}\ \isactrlbold {\isadigit{0}}\ \isactrlbold {\isasymcdot}\ {\isacharparenleft}{\kern0pt}\isaconst{cmp}\ \isaconst{length{\isacharunderscore}{\kern0pt}encIdx}\ \isactrlbold {\isasymcdot}\ {\isacharparenleft}{\kern0pt}\isactrlbold n{\isasymlbrace}\isaconst{phP}{\isasymrbrace}{\isacharparenright}{\kern0pt}\ \isactrlbold {\isasymcdot}\ {\isacharparenleft}{\kern0pt}\isactrlbold n{\isasymlbrace}\isaconst{phP}{\isasymrbrace}{\isacharparenright}{\kern0pt}{\isacharparenright}{\kern0pt}{\isacharparenright}{\kern0pt}{\isacharparenright}{\kern0pt}{\isacharparenright}{\kern0pt}{\isacharparenright}{\kern0pt}{\isacharparenright}{\kern0pt}{\isacharparenright}{\kern0pt}\ {\isacharparenleft}{\kern0pt}\isaconst{wBa}\ {\isadigit{0}}\ {\isacharparenleft}{\kern0pt}\isaconst{iPoBase}\ \isaconst{kidRdIdx}\ {\isacharparenleft}{\kern0pt}\isactrlbold n{\isasymlbrace}\isaconst{phP}{\isasymrbrace}{\isacharparenright}{\kern0pt}\ {\isacharparenleft}{\kern0pt}\isaconst{iPoRk}\ {\isacharparenleft}{\kern0pt}\isaconst{cmp}\ \isaconst{pomaxIdx}\ \isactrlbold {\isasymcdot}\ \isactrlbold {\isadigit{0}}\ \isactrlbold {\isasymcdot}\ \isaconst{trmPredR}\ {\isacharparenleft}{\kern0pt}\isaconst{iPoHt}\ \isaconst{kidRdIdx}\ {\isacharparenleft}{\kern0pt}\isactrlbold n{\isasymlbrace}\isaconst{phP}{\isasymrbrace}{\isacharparenright}{\kern0pt}\ {\isacharparenleft}{\kern0pt}\isaconst{cmp}\ \isaconst{pomaxIdx}\ \isactrlbold {\isasymcdot}\ \isactrlbold {\isadigit{0}}\ \isactrlbold {\isasymcdot}\ {\isacharparenleft}{\kern0pt}\isaconst{cmp}\ \isaconst{length{\isacharunderscore}{\kern0pt}encIdx}\ \isactrlbold {\isasymcdot}\ {\isacharparenleft}{\kern0pt}\isactrlbold n{\isasymlbrace}\isaconst{phP}{\isasymrbrace}{\isacharparenright}{\kern0pt}\ \isactrlbold {\isasymcdot}\ {\isacharparenleft}{\kern0pt}\isactrlbold n{\isasymlbrace}\isaconst{phP}{\isasymrbrace}{\isacharparenright}{\kern0pt}{\isacharparenright}{\kern0pt}{\isacharparenright}{\kern0pt}{\isacharparenright}{\kern0pt}{\isacharparenright}{\kern0pt}\ {\isacharparenleft}{\kern0pt}\isaconst{wBa}\ {\isacharparenleft}{\kern0pt}\isaconst{Suc}\ {\isadigit{0}}{\isacharparenright}{\kern0pt}\ {\isacharparenleft}{\kern0pt}\isaconst{iPoBase}\ \isaconst{kidRdIdx}\ {\isacharparenleft}{\kern0pt}\isactrlbold n{\isasymlbrace}\isaconst{phP}{\isasymrbrace}{\isacharparenright}{\kern0pt}\ {\isacharparenleft}{\kern0pt}\isaconst{cmp}\ \isaconst{pomaxIdx}\ \isactrlbold {\isasymcdot}\ \isactrlbold {\isadigit{0}}\ \isactrlbold {\isasymcdot}\ {\isacharparenleft}{\kern0pt}\isaconst{cmp}\ \isaconst{length{\isacharunderscore}{\kern0pt}encIdx}\ \isactrlbold {\isasymcdot}\ {\isacharparenleft}{\kern0pt}\isactrlbold n{\isasymlbrace}\isaconst{phP}{\isasymrbrace}{\isacharparenright}{\kern0pt}\ \isactrlbold {\isasymcdot}\ {\isacharparenleft}{\kern0pt}\isactrlbold n{\isasymlbrace}\isaconst{phP}{\isasymrbrace}{\isacharparenright}{\kern0pt}{\isacharparenright}{\kern0pt}{\isacharparenright}{\kern0pt}{\isacharparenright}{\kern0pt}\ {\isacharparenleft}{\kern0pt}\isaconst{iPoHt}\ \isaconst{kidRdIdx}\ {\isacharparenleft}{\kern0pt}\isactrlbold n{\isasymlbrace}\isaconst{phP}{\isasymrbrace}{\isacharparenright}{\kern0pt}\ {\isacharparenleft}{\kern0pt}\isaconst{cmp}\ \isaconst{pomaxIdx}\ \isactrlbold {\isasymcdot}\ \isactrlbold {\isadigit{0}}\ \isactrlbold {\isasymcdot}\ {\isacharparenleft}{\kern0pt}\isaconst{cmp}\ \isaconst{length{\isacharunderscore}{\kern0pt}encIdx}\ \isactrlbold {\isasymcdot}\ {\isacharparenleft}{\kern0pt}\isactrlbold n{\isasymlbrace}\isaconst{phP}{\isasymrbrace}{\isacharparenright}{\kern0pt}\ \isactrlbold {\isasymcdot}\ {\isacharparenleft}{\kern0pt}\isactrlbold n{\isasymlbrace}\isaconst{phP}{\isasymrbrace}{\isacharparenright}{\kern0pt}{\isacharparenright}{\kern0pt}{\isacharparenright}{\kern0pt}{\isacharparenright}{\kern0pt}{\isacharparenright}{\kern0pt}{\isacharparenright}{\kern0pt}{\isacharparenright}{\kern0pt}\ {\isacharparenleft}{\kern0pt}\isaconst{iPoHt}\ \isaconst{kidRdIdx}\ {\isacharparenleft}{\kern0pt}\isactrlbold n{\isasymlbrace}\isaconst{phP}{\isasymrbrace}{\isacharparenright}{\kern0pt}\ {\isacharparenleft}{\kern0pt}\isaconst{iPoRk}\ {\isacharparenleft}{\kern0pt}\isaconst{cmp}\ \isaconst{pomaxIdx}\ \isactrlbold {\isasymcdot}\ \isactrlbold {\isadigit{0}}\ \isactrlbold {\isasymcdot}\ \isaconst{trmPredR}\ {\isacharparenleft}{\kern0pt}\isaconst{iPoHt}\ \isaconst{kidRdIdx}\ {\isacharparenleft}{\kern0pt}\isactrlbold n{\isasymlbrace}\isaconst{phP}{\isasymrbrace}{\isacharparenright}{\kern0pt}\ {\isacharparenleft}{\kern0pt}\isaconst{cmp}\ \isaconst{pomaxIdx}\ \isactrlbold {\isasymcdot}\ \isactrlbold {\isadigit{0}}\ \isactrlbold {\isasymcdot}\ {\isacharparenleft}{\kern0pt}\isaconst{cmp}\ \isaconst{length{\isacharunderscore}{\kern0pt}encIdx}\ \isactrlbold {\isasymcdot}\ {\isacharparenleft}{\kern0pt}\isactrlbold n{\isasymlbrace}\isaconst{phP}{\isasymrbrace}{\isacharparenright}{\kern0pt}\ \isactrlbold {\isasymcdot}\ {\isacharparenleft}{\kern0pt}\isactrlbold n{\isasymlbrace}\isaconst{phP}{\isasymrbrace}{\isacharparenright}{\kern0pt}{\isacharparenright}{\kern0pt}{\isacharparenright}{\kern0pt}{\isacharparenright}{\kern0pt}{\isacharparenright}{\kern0pt}\ {\isacharparenleft}{\kern0pt}\isaconst{wBa}\ {\isacharparenleft}{\kern0pt}\isaconst{Suc}\ {\isadigit{0}}{\isacharparenright}{\kern0pt}\ {\isacharparenleft}{\kern0pt}\isaconst{iPoBase}\ \isaconst{kidRdIdx}\ {\isacharparenleft}{\kern0pt}\isactrlbold n{\isasymlbrace}\isaconst{phP}{\isasymrbrace}{\isacharparenright}{\kern0pt}\ {\isacharparenleft}{\kern0pt}\isaconst{cmp}\ \isaconst{pomaxIdx}\ \isactrlbold {\isasymcdot}\ \isactrlbold {\isadigit{0}}\ \isactrlbold {\isasymcdot}\ {\isacharparenleft}{\kern0pt}\isaconst{cmp}\ \isaconst{length{\isacharunderscore}{\kern0pt}encIdx}\ \isactrlbold {\isasymcdot}\ {\isacharparenleft}{\kern0pt}\isactrlbold n{\isasymlbrace}\isaconst{phP}{\isasymrbrace}{\isacharparenright}{\kern0pt}\ \isactrlbold {\isasymcdot}\ {\isacharparenleft}{\kern0pt}\isactrlbold n{\isasymlbrace}\isaconst{phP}{\isasymrbrace}{\isacharparenright}{\kern0pt}{\isacharparenright}{\kern0pt}{\isacharparenright}{\kern0pt}{\isacharparenright}{\kern0pt}\ {\isacharparenleft}{\kern0pt}\isaconst{iPoHt}\ \isaconst{kidRdIdx}\ {\isacharparenleft}{\kern0pt}\isactrlbold n{\isasymlbrace}\isaconst{phP}{\isasymrbrace}{\isacharparenright}{\kern0pt}\ {\isacharparenleft}{\kern0pt}\isaconst{cmp}\ \isaconst{pomaxIdx}\ \isactrlbold {\isasymcdot}\ \isactrlbold {\isadigit{0}}\ \isactrlbold {\isasymcdot}\ {\isacharparenleft}{\kern0pt}\isaconst{cmp}\ \isaconst{length{\isacharunderscore}{\kern0pt}encIdx}\ \isactrlbold {\isasymcdot}\ {\isacharparenleft}{\kern0pt}\isactrlbold n{\isasymlbrace}\isaconst{phP}{\isasymrbrace}{\isacharparenright}{\kern0pt}\ \isactrlbold {\isasymcdot}\ {\isacharparenleft}{\kern0pt}\isactrlbold n{\isasymlbrace}\isaconst{phP}{\isasymrbrace}{\isacharparenright}{\kern0pt}{\isacharparenright}{\kern0pt}{\isacharparenright}{\kern0pt}{\isacharparenright}{\kern0pt}{\isacharparenright}{\kern0pt}{\isacharparenright}{\kern0pt}{\isacharparenright}{\kern0pt}{\isacharparenright}{\kern0pt}{\isacharparenright}{\kern0pt}\isasep\isanewline%
}%
\SNIP{Goedel-Paper-Statements:thm:iga-not-motive-3}{%
{\isasymlbrakk}\isaconst{wfc}\ \isavar{{\isasymGamma}}{\isacharsemicolon}{\kern0pt}\ {\isasymAnd}\isabound{g}{\isachardot}{\kern0pt}\ \isabound{g}\ {\isacharbar}{\kern0pt}{\isasymin}{\isacharbar}{\kern0pt}\ \isavar{{\isasymGamma}}\ {\isasymLongrightarrow}\ \isaconst{omv}\ \isabound{g}\ {\isasymlbrace}\isaconst{True}{\isasymrbrace}{\isasymrbrakk}\ {\isasymLongrightarrow}\ {\isasymnot}\ \isaconst{prf{\isacharunderscore}{\kern0pt}derives{\isachardot}{\kern0pt}entails}\ \isaconst{gi{\isachardot}{\kern0pt}derives}\ \isavar{{\isasymGamma}}\ \isaconst{phMotAll}}%
}%
\SNIP{Goedel-Paper-Statements:thm:iga-not-fold-residual-0}{%
\isakeywordONE{thm}\isamarkupfalse%
\ iga{\isacharunderscore}{\kern0pt}not{\isacharunderscore}{\kern0pt}fold{\isacharunderscore}{\kern0pt}residual%
\isanewline%
}%
\SNIP{Goedel-Paper-Statements:thm:iga-not-fold-residual-1}{%
\isa{{\isasymlbrakk}\isaconst{armsT}\ \isavar{Rs}{\isacharsemicolon}{\kern0pt}\ \isaconst{wft}\ {\isacharparenleft}{\kern0pt}\isavar{Rs}\ {\isacharparenleft}{\kern0pt}\isactrlbold n{\isasymlbrace}\isaconst{phP}{\isasymrbrace}{\isacharparenright}{\kern0pt}{\isacharparenright}{\kern0pt}{\isacharsemicolon}{\kern0pt}\ {\isasymnot}\ \isaconst{free}\ {\isadigit{1}}\ {\isacharparenleft}{\kern0pt}\isavar{Rs}\ {\isacharparenleft}{\kern0pt}\isactrlbold n{\isasymlbrace}\isaconst{phP}{\isasymrbrace}{\isacharparenright}{\kern0pt}{\isacharparenright}{\kern0pt}{\isacharsemicolon}{\kern0pt}\ \isaconst{wfc}\ \isavar{{\isasymGamma}}{\isacharsemicolon}{\kern0pt}\ {\isasymAnd}\isabound{g}{\isachardot}{\kern0pt}\ \isabound{g}\ {\isacharbar}{\kern0pt}{\isasymin}{\isacharbar}{\kern0pt}\ \isavar{{\isasymGamma}}\ {\isasymLongrightarrow}\ \isaconst{omv}\ \isabound{g}\ {\isasymlbrace}\isaconst{True}{\isasymrbrace}{\isasymrbrakk}\ {\isasymLongrightarrow}\ {\isasymnot}\ \isaconst{prf{\isacharunderscore}{\kern0pt}derives{\isachardot}{\kern0pt}entails}\ \isaconst{gi{\isachardot}{\kern0pt}derives}\ \isavar{{\isasymGamma}}\ {\isacharparenleft}{\kern0pt}\isavar{Rs}\ {\isacharparenleft}{\kern0pt}\isactrlbold n{\isasymlbrace}\isaconst{phP}{\isasymrbrace}{\isacharparenright}{\kern0pt}{\isacharparenright}{\kern0pt}}%
}%
\SNIP{Goedel-Paper-Statements:thm:diagSubIdx-def-0}{%
\isakeywordONE{thm}\isamarkupfalse%
\ diagSubIdx{\isacharunderscore}{\kern0pt}def%
\isanewline%
}%
\SNIP{Goedel-Paper-Statements:thm:diagSubIdx-def-1}{%
\isa{\isaconst{diagSubIdx}\ {\isasymequiv}\ \isaconst{prSub}\ \isaconst{submMIdx}\ \isaconst{diagListIdx}\ {\isacharparenleft}{\kern0pt}\isaconst{prSub}\ \isaconst{cpairIdx}\ {\isacharparenleft}{\kern0pt}\isaconst{numIdx}\ {\isadigit{0}}{\isacharparenright}{\kern0pt}\ {\isacharparenleft}{\kern0pt}\isaconst{prInit}\ \isaconst{rInitL}{\isacharparenright}{\kern0pt}{\isacharparenright}{\kern0pt}}%
}%
\SNIP{Goedel-Paper-Statements:thm:iDiagSub-def-0}{%
\isakeywordONE{thm}\isamarkupfalse%
\ iDiagSub{\isacharunderscore}{\kern0pt}def%
\isanewline%
}%
\SNIP{Goedel-Paper-Statements:thm:iDiagSub-def-1}{%
\isa{\isaconst{iDiagSub}\ \isavar{x}\ \isavar{y}\ {\isasymequiv}\ \isaconst{cmp}\ \isaconst{diagSubIdx}\ \isactrlbold {\isasymcdot}\ \isavar{x}\ \isactrlbold {\isasymcdot}\ \isavar{y}}%
}%
\SNIP{Goedel-Paper-Statements:thm:iDiagSub-closed-0}{%
\isakeywordONE{thm}\isamarkupfalse%
\ iDiagSub{\isacharunderscore}{\kern0pt}closed%
\isanewline%
}%
\SNIP{Goedel-Paper-Statements:thm:iDiagSub-closed-1}{%
\isa{\isaconst{wft}\ \isavar{t}\ {\isasymLongrightarrow}\ {\isasymemptyset}\ {\isasymturnstile}\ \isaconst{iDiagSub}\ {\isacharparenleft}{\kern0pt}\isactrlbold n{\isasymlbrace}\isavar{t}{\isasymrbrace}{\isacharparenright}{\kern0pt}\ {\isacharparenleft}{\kern0pt}\isactrlbold n\isavar{k}{\isacharparenright}{\kern0pt}\ \isactrlbold {\isacharequal}{\kern0pt}\ \isactrlbold n{\isasymlbrace}{\isacharbrackleft}{\kern0pt}{\isadigit{0}}{\isasymZpinj}\isactrlbold n\isavar{k}{\isacharbrackright}{\kern0pt}\isavar{t}{\isasymrbrace}}%
}%
\SNIP{Goedel-Paper-Statements:thm:goedelD-def-0}{%
\isakeywordONE{thm}\isamarkupfalse%
\ goedelD{\isacharunderscore}{\kern0pt}def%
\isanewline%
}%
\SNIP{Goedel-Paper-Statements:thm:goedelD-def-1}{%
\isa{\isaconst{goedelD}\ {\isasymequiv}\ \isactrlbold A\ \isactrlbold {\isasymcdot}\ {\isacharparenleft}{\kern0pt}\isactrlbold {\isasymlambda}{\isacharparenleft}{\kern0pt}\isactrlbold {\isasymnot}\isaconst{iChk}\ \isactrlbold v\isactrlsub {\isadigit{0}}\ {\isacharparenleft}{\kern0pt}\isaconst{iDiagSub}\ {\isacharparenleft}{\kern0pt}\isaconst{trmVar}\ {\isadigit{1}}{\isacharparenright}{\kern0pt}\ {\isacharparenleft}{\kern0pt}\isaconst{trmVar}\ {\isadigit{1}}{\isacharparenright}{\kern0pt}{\isacharparenright}{\kern0pt}{\isacharparenright}{\kern0pt}{\isacharparenright}{\kern0pt}}%
}%
\SNIP{Goedel-Paper-Statements:thm:goedelS-def-0}{%
\isakeywordONE{thm}\isamarkupfalse%
\ goedelS{\isacharunderscore}{\kern0pt}def%
\isanewline%
}%
\SNIP{Goedel-Paper-Statements:thm:goedelS-def-1}{%
\isa{\isaconst{goedelS}\ {\isasymequiv}\ \isaconst{iDiagSub}\ {\isacharparenleft}{\kern0pt}\isactrlbold n{\isasymlbrace}\isaconst{goedelD}{\isasymrbrace}{\isacharparenright}{\kern0pt}\ {\isacharparenleft}{\kern0pt}\isactrlbold n{\isasymlbrace}\isaconst{goedelD}{\isasymrbrace}{\isacharparenright}{\kern0pt}}%
}%
\SNIP{Goedel-Paper-Statements:thm:goedelG-def-0}{%
\isakeywordONE{thm}\isamarkupfalse%
\ goedelG{\isacharunderscore}{\kern0pt}def%
\isanewline%
}%
\SNIP{Goedel-Paper-Statements:thm:goedelG-def-1}{%
\isa{\isaconst{goedelG}\ {\isasymequiv}\ \isactrlbold A\ \isactrlbold {\isasymcdot}\ {\isacharparenleft}{\kern0pt}\isactrlbold {\isasymlambda}{\isacharparenleft}{\kern0pt}\isactrlbold {\isasymnot}\isaconst{iChk}\ \isactrlbold v\isactrlsub {\isadigit{0}}\ \isaconst{goedelS}{\isacharparenright}{\kern0pt}{\isacharparenright}{\kern0pt}}%
}%
\SNIP{Goedel-Paper-Statements:thm:goedelF-def-0}{%
\isakeywordONE{thm}\isamarkupfalse%
\ goedelF{\isacharunderscore}{\kern0pt}def%
\isanewline%
}%
\SNIP{Goedel-Paper-Statements:thm:goedelF-def-1}{%
\isa{\isaconst{goedelF}\ {\isasymequiv}\ \isactrlbold {\isasymlambda}{\isacharparenleft}{\kern0pt}\isactrlbold {\isasymnot}\isaconst{iChk}\ \isactrlbold v\isactrlsub {\isadigit{0}}\ \isaconst{goedelS}{\isacharparenright}{\kern0pt}}%
}%
\SNIP{Goedel-Paper-Statements:thm:ConRGA-def-0}{%
\isakeywordONE{thm}\isamarkupfalse%
\ ConRGA{\isacharunderscore}{\kern0pt}def%
\isanewline%
}%
\SNIP{Goedel-Paper-Statements:thm:ConRGA-def-1}{%
\isa{\isaconst{ConRGA}\ {\isasymequiv}\ \isactrlbold A\ \isactrlbold {\isasymcdot}\ {\isacharparenleft}{\kern0pt}\isactrlbold {\isasymlambda}{\isacharparenleft}{\kern0pt}\isactrlbold {\isasymnot}\isaconst{iChk}\ \isactrlbold v\isactrlsub {\isadigit{0}}\ \isaconst{cEq{\isadigit{0}}{\isadigit{1}}}{\isacharparenright}{\kern0pt}{\isacharparenright}{\kern0pt}}%
}%
\SNIP{Goedel-Paper-Statements:thm:iChk-conc-transport-0}{%
\isakeywordONE{thm}\isamarkupfalse%
\ iChk{\isacharunderscore}{\kern0pt}conc{\isacharunderscore}{\kern0pt}transport%
\isanewline%
}%
\SNIP{Goedel-Paper-Statements:thm:iChk-conc-transport-1}{%
\isa{{\isasymlbrakk}\isavar{{\isasymGamma}}\ {\isasymturnstile}\ \isavar{c}\ \isactrlbold {\isacharequal}{\kern0pt}\ \isavar{c{\isacharprime}{\kern0pt}}{\isacharsemicolon}{\kern0pt}\ \isavar{{\isasymGamma}}\ {\isasymturnstile}\ \isaconst{iChk}\ \isavar{p}\ \isavar{c}{\isasymrbrakk}\ {\isasymLongrightarrow}\ \isavar{{\isasymGamma}}\ {\isasymturnstile}\ \isaconst{iChk}\ \isavar{p}\ \isavar{c{\isacharprime}{\kern0pt}}}%
}%
\SNIP{Goedel-Paper-Statements:thm:iChk-proof-transport-0}{%
\isakeywordONE{thm}\isamarkupfalse%
\ iChk{\isacharunderscore}{\kern0pt}proof{\isacharunderscore}{\kern0pt}transport%
\isanewline%
}%
\SNIP{Goedel-Paper-Statements:thm:iChk-proof-transport-1}{%
\isa{{\isasymlbrakk}\isavar{{\isasymGamma}}\ {\isasymturnstile}\ \isavar{p}\ \isactrlbold {\isacharequal}{\kern0pt}\ \isavar{p{\isacharprime}{\kern0pt}}{\isacharsemicolon}{\kern0pt}\ \isavar{{\isasymGamma}}\ {\isasymturnstile}\ \isaconst{iChk}\ \isavar{p}\ \isavar{c}{\isasymrbrakk}\ {\isasymLongrightarrow}\ \isavar{{\isasymGamma}}\ {\isasymturnstile}\ \isaconst{iChk}\ \isavar{p{\isacharprime}{\kern0pt}}\ \isavar{c}}%
}%
\SNIP{Goedel-Paper-Statements:thm:iChk-headI-0}{%
\isakeywordONE{thm}\isamarkupfalse%
\ iChk{\isacharunderscore}{\kern0pt}headI%
\isanewline%
}%
\SNIP{Goedel-Paper-Statements:thm:iChk-headI-1}{%
\isa{{\isasymlbrakk}\isavar{{\isasymGamma}}\ {\isasymturnstile}\ \isavar{P}\ \isactrlbold N{\isacharsemicolon}{\kern0pt}\ \isavar{{\isasymGamma}}\ {\isasymturnstile}\ \isavar{c}\ \isactrlbold N{\isacharsemicolon}{\kern0pt}\ \isavar{{\isasymGamma}}\ {\isasymturnstile}\ \isavar{P}\ \isactrlbold {\isasymnoteq}\ \isactrlbold {\isadigit{0}}{\isacharsemicolon}{\kern0pt}\ \isavar{{\isasymGamma}}\ {\isasymturnstile}\ \isaconst{cmp}\ \isaconst{validChkIdx}\ \isactrlbold {\isasymcdot}\ \isavar{P}\ \isactrlbold {\isasymcdot}\ \isavar{P}\ \isactrlbold {\isacharequal}{\kern0pt}\ \isactrlbold n\isaconst{Suc}\ {\isadigit{0}}{\isacharsemicolon}{\kern0pt}\ \isavar{{\isasymGamma}}\ {\isasymturnstile}\ \isaconst{cpFst}\ {\isacharparenleft}{\kern0pt}\isaconst{trmPredR}\ \isavar{P}{\isacharparenright}{\kern0pt}\ \isactrlbold {\isacharequal}{\kern0pt}\ \isaconst{cpPair}\ {\isacharparenleft}{\kern0pt}\isactrlbold n{\isasymlbrace}{\isasymemptyset}{\isasymrbrace}{\isacharparenright}{\kern0pt}\ \isavar{c}{\isasymrbrakk}\ {\isasymLongrightarrow}\ \isavar{{\isasymGamma}}\ {\isasymturnstile}\ \isaconst{iChk}\ \isavar{P}\ \isavar{c}}%
}%
\SNIP{Goedel-Paper-Statements:thm:iChk-headD-0}{%
\isakeywordONE{thm}\isamarkupfalse%
\ iChk{\isacharunderscore}{\kern0pt}headD%
\isanewline%
}%
\SNIP{Goedel-Paper-Statements:thm:iChk-headD-1}{%
\isa{{\isasymlbrakk}\isavar{{\isasymGamma}}\ {\isasymturnstile}\ \isavar{P}\ \isactrlbold N{\isacharsemicolon}{\kern0pt}\ \isavar{{\isasymGamma}}\ {\isasymturnstile}\ \isaconst{iChk}\ \isavar{P}\ \isavar{c}{\isasymrbrakk}\ {\isasymLongrightarrow}\ \isavar{{\isasymGamma}}\ {\isasymturnstile}\ \isaconst{cpFst}\ {\isacharparenleft}{\kern0pt}\isaconst{trmPredR}\ \isavar{P}{\isacharparenright}{\kern0pt}\ \isactrlbold {\isacharequal}{\kern0pt}\ \isaconst{cpPair}\ {\isacharparenleft}{\kern0pt}\isactrlbold n{\isasymlbrace}{\isasymemptyset}{\isasymrbrace}{\isacharparenright}{\kern0pt}\ \isavar{c}}%
}%
\SNIP{Goedel-Paper-Statements:thm:iChk-appC-0}{%
\isakeywordONE{thm}\isamarkupfalse%
\ iChk{\isacharunderscore}{\kern0pt}appC%
\isanewline%
}%
\SNIP{Goedel-Paper-Statements:thm:iChk-appC-1}{%
\isa{{\isasymlbrakk}\isavar{{\isasymGamma}}\ {\isasymturnstile}\ \isavar{A}\ \isactrlbold N{\isacharsemicolon}{\kern0pt}\ \isavar{{\isasymGamma}}\ {\isasymturnstile}\ \isavar{B}\ \isactrlbold N{\isacharsemicolon}{\kern0pt}\ \isavar{{\isasymGamma}}\ {\isasymturnstile}\ \isaconst{iChk}\ \isavar{A}\ \isavar{a}{\isacharsemicolon}{\kern0pt}\ \isavar{{\isasymGamma}}\ {\isasymturnstile}\ \isaconst{iChk}\ \isavar{B}\ \isavar{b}{\isasymrbrakk}\ {\isasymLongrightarrow}\ \isavar{{\isasymGamma}}\ {\isasymturnstile}\ \isaconst{cmp}\ \isaconst{apCoreIdx}\ \isactrlbold {\isasymcdot}\ \isaconst{cpPair}\ \isavar{A}\ \isavar{B}\ \isactrlbold {\isasymcdot}\ {\isacharparenleft}{\kern0pt}\isaconst{cmp}\ \isaconst{length{\isacharunderscore}{\kern0pt}encIdx}\ \isactrlbold {\isasymcdot}\ \isavar{A}\ \isactrlbold {\isasymcdot}\ \isavar{A}{\isacharparenright}{\kern0pt}\ \isactrlbold N\isasep\isanewline%
}%
\SNIP{Goedel-Paper-Statements:thm:iChk-appC-2}{%
{\isasymlbrakk}\isavar{{\isasymGamma}}\ {\isasymturnstile}\ \isavar{A}\ \isactrlbold N{\isacharsemicolon}{\kern0pt}\ \isavar{{\isasymGamma}}\ {\isasymturnstile}\ \isavar{B}\ \isactrlbold N{\isacharsemicolon}{\kern0pt}\ \isavar{{\isasymGamma}}\ {\isasymturnstile}\ \isaconst{iChk}\ \isavar{A}\ \isavar{a}{\isacharsemicolon}{\kern0pt}\ \isavar{{\isasymGamma}}\ {\isasymturnstile}\ \isaconst{iChk}\ \isavar{B}\ \isavar{b}{\isasymrbrakk}\ {\isasymLongrightarrow}\ \isavar{{\isasymGamma}}\ {\isasymturnstile}\ \isaconst{cmp}\ \isaconst{apCoreIdx}\ \isactrlbold {\isasymcdot}\ \isaconst{cpPair}\ \isavar{A}\ \isavar{B}\ \isactrlbold {\isasymcdot}\ {\isacharparenleft}{\kern0pt}\isaconst{cmp}\ \isaconst{length{\isacharunderscore}{\kern0pt}encIdx}\ \isactrlbold {\isasymcdot}\ \isavar{A}\ \isactrlbold {\isasymcdot}\ \isavar{A}{\isacharparenright}{\kern0pt}\ \isactrlbold {\isasymnoteq}\ \isactrlbold {\isadigit{0}}\isasep\isanewline%
}%
\SNIP{Goedel-Paper-Statements:thm:iChk-appC-3}{%
{\isasymlbrakk}\isavar{{\isasymGamma}}\ {\isasymturnstile}\ \isavar{A}\ \isactrlbold N{\isacharsemicolon}{\kern0pt}\ \isavar{{\isasymGamma}}\ {\isasymturnstile}\ \isavar{B}\ \isactrlbold N{\isacharsemicolon}{\kern0pt}\ \isavar{{\isasymGamma}}\ {\isasymturnstile}\ \isaconst{iChk}\ \isavar{A}\ \isavar{a}{\isacharsemicolon}{\kern0pt}\ \isavar{{\isasymGamma}}\ {\isasymturnstile}\ \isaconst{iChk}\ \isavar{B}\ \isavar{b}{\isasymrbrakk}\ {\isasymLongrightarrow}\ \isavar{{\isasymGamma}}\ {\isasymturnstile}\ \isaconst{cmp}\ \isaconst{validChkIdx}\ \isactrlbold {\isasymcdot}\ {\isacharparenleft}{\kern0pt}\isaconst{cmp}\ \isaconst{apCoreIdx}\ \isactrlbold {\isasymcdot}\ \isaconst{cpPair}\ \isavar{A}\ \isavar{B}\ \isactrlbold {\isasymcdot}\ {\isacharparenleft}{\kern0pt}\isaconst{cmp}\ \isaconst{length{\isacharunderscore}{\kern0pt}encIdx}\ \isactrlbold {\isasymcdot}\ \isavar{A}\ \isactrlbold {\isasymcdot}\ \isavar{A}{\isacharparenright}{\kern0pt}{\isacharparenright}{\kern0pt}\ \isactrlbold {\isasymcdot}\ {\isacharparenleft}{\kern0pt}\isaconst{cmp}\ \isaconst{apCoreIdx}\ \isactrlbold {\isasymcdot}\ \isaconst{cpPair}\ \isavar{A}\ \isavar{B}\ \isactrlbold {\isasymcdot}\ {\isacharparenleft}{\kern0pt}\isaconst{cmp}\ \isaconst{length{\isacharunderscore}{\kern0pt}encIdx}\ \isactrlbold {\isasymcdot}\ \isavar{A}\ \isactrlbold {\isasymcdot}\ \isavar{A}{\isacharparenright}{\kern0pt}{\isacharparenright}{\kern0pt}\ \isactrlbold {\isacharequal}{\kern0pt}\ \isactrlbold n\isaconst{Suc}\ {\isadigit{0}}\isasep\isanewline%
}%
\SNIP{Goedel-Paper-Statements:thm:iChk-appC-4}{%
{\isasymlbrakk}\isavar{{\isasymGamma}}\ {\isasymturnstile}\ \isavar{A}\ \isactrlbold N{\isacharsemicolon}{\kern0pt}\ \isavar{{\isasymGamma}}\ {\isasymturnstile}\ \isavar{B}\ \isactrlbold N{\isacharsemicolon}{\kern0pt}\ \isavar{{\isasymGamma}}\ {\isasymturnstile}\ \isaconst{iChk}\ \isavar{A}\ \isavar{a}{\isacharsemicolon}{\kern0pt}\ \isavar{{\isasymGamma}}\ {\isasymturnstile}\ \isaconst{iChk}\ \isavar{B}\ \isavar{b}{\isasymrbrakk}\ {\isasymLongrightarrow}\ \isavar{{\isasymGamma}}\ {\isasymturnstile}\ \isaconst{cpFst}\ {\isacharparenleft}{\kern0pt}\isaconst{trmPredR}\ {\isacharparenleft}{\kern0pt}\isaconst{cmp}\ \isaconst{apCoreIdx}\ \isactrlbold {\isasymcdot}\ \isaconst{cpPair}\ \isavar{A}\ \isavar{B}\ \isactrlbold {\isasymcdot}\ {\isacharparenleft}{\kern0pt}\isaconst{cmp}\ \isaconst{length{\isacharunderscore}{\kern0pt}encIdx}\ \isactrlbold {\isasymcdot}\ \isavar{A}\ \isactrlbold {\isasymcdot}\ \isavar{A}{\isacharparenright}{\kern0pt}{\isacharparenright}{\kern0pt}{\isacharparenright}{\kern0pt}\ \isactrlbold {\isacharequal}{\kern0pt}\ \isaconst{cpPair}\ {\isacharparenleft}{\kern0pt}\isactrlbold n{\isasymlbrace}{\isasymemptyset}{\isasymrbrace}{\isacharparenright}{\kern0pt}\ \isavar{a}\isasep\isanewline%
}%
\SNIP{Goedel-Paper-Statements:thm:iChk-appC-5}{%
{\isasymlbrakk}\isavar{{\isasymGamma}}\ {\isasymturnstile}\ \isavar{A}\ \isactrlbold N{\isacharsemicolon}{\kern0pt}\ \isavar{{\isasymGamma}}\ {\isasymturnstile}\ \isavar{B}\ \isactrlbold N{\isacharsemicolon}{\kern0pt}\ \isavar{{\isasymGamma}}\ {\isasymturnstile}\ \isaconst{iChk}\ \isavar{A}\ \isavar{a}{\isacharsemicolon}{\kern0pt}\ \isavar{{\isasymGamma}}\ {\isasymturnstile}\ \isaconst{iChk}\ \isavar{B}\ \isavar{b}{\isasymrbrakk}\ {\isasymLongrightarrow}\ \isavar{{\isasymGamma}}\ {\isasymturnstile}\ \isaconst{gapAdd}\ {\isacharparenleft}{\kern0pt}\isaconst{cmp}\ \isaconst{length{\isacharunderscore}{\kern0pt}encIdx}\ \isactrlbold {\isasymcdot}\ \isavar{A}\ \isactrlbold {\isasymcdot}\ \isavar{A}{\isacharparenright}{\kern0pt}\ \isactrlbold {\isadigit{0}}\ \isactrlbold N\isasep\isanewline%
}%
\SNIP{Goedel-Paper-Statements:thm:iChk-appC-6}{%
{\isasymlbrakk}\isavar{{\isasymGamma}}\ {\isasymturnstile}\ \isavar{A}\ \isactrlbold N{\isacharsemicolon}{\kern0pt}\ \isavar{{\isasymGamma}}\ {\isasymturnstile}\ \isavar{B}\ \isactrlbold N{\isacharsemicolon}{\kern0pt}\ \isavar{{\isasymGamma}}\ {\isasymturnstile}\ \isaconst{iChk}\ \isavar{A}\ \isavar{a}{\isacharsemicolon}{\kern0pt}\ \isavar{{\isasymGamma}}\ {\isasymturnstile}\ \isaconst{iChk}\ \isavar{B}\ \isavar{b}{\isasymrbrakk}\ {\isasymLongrightarrow}\ \isavar{{\isasymGamma}}\ {\isasymturnstile}\ \isaconst{gapSub}\ {\isacharparenleft}{\kern0pt}\isactrlbold S\ \isaconst{gapAdd}\ {\isacharparenleft}{\kern0pt}\isaconst{cmp}\ \isaconst{length{\isacharunderscore}{\kern0pt}encIdx}\ \isactrlbold {\isasymcdot}\ \isavar{A}\ \isactrlbold {\isasymcdot}\ \isavar{A}{\isacharparenright}{\kern0pt}\ \isactrlbold {\isadigit{0}}{\isacharparenright}{\kern0pt}\ {\isacharparenleft}{\kern0pt}\isaconst{cmp}\ \isaconst{length{\isacharunderscore}{\kern0pt}encIdx}\ \isactrlbold {\isasymcdot}\ {\isacharparenleft}{\kern0pt}\isaconst{cmp}\ \isaconst{apCoreIdx}\ \isactrlbold {\isasymcdot}\ \isaconst{cpPair}\ \isavar{A}\ \isavar{B}\ \isactrlbold {\isasymcdot}\ {\isacharparenleft}{\kern0pt}\isaconst{cmp}\ \isaconst{length{\isacharunderscore}{\kern0pt}encIdx}\ \isactrlbold {\isasymcdot}\ \isavar{A}\ \isactrlbold {\isasymcdot}\ \isavar{A}{\isacharparenright}{\kern0pt}{\isacharparenright}{\kern0pt}\ \isactrlbold {\isasymcdot}\ {\isacharparenleft}{\kern0pt}\isaconst{cmp}\ \isaconst{apCoreIdx}\ \isactrlbold {\isasymcdot}\ \isaconst{cpPair}\ \isavar{A}\ \isavar{B}\ \isactrlbold {\isasymcdot}\ {\isacharparenleft}{\kern0pt}\isaconst{cmp}\ \isaconst{length{\isacharunderscore}{\kern0pt}encIdx}\ \isactrlbold {\isasymcdot}\ \isavar{A}\ \isactrlbold {\isasymcdot}\ \isavar{A}{\isacharparenright}{\kern0pt}{\isacharparenright}{\kern0pt}{\isacharparenright}{\kern0pt}\ \isactrlbold {\isacharequal}{\kern0pt}\ \isactrlbold {\isadigit{0}}\isasep\isanewline%
}%
\SNIP{Goedel-Paper-Statements:thm:iChk-appC-7}{%
{\isasymlbrakk}\isavar{{\isasymGamma}}\ {\isasymturnstile}\ \isavar{A}\ \isactrlbold N{\isacharsemicolon}{\kern0pt}\ \isavar{{\isasymGamma}}\ {\isasymturnstile}\ \isavar{B}\ \isactrlbold N{\isacharsemicolon}{\kern0pt}\ \isavar{{\isasymGamma}}\ {\isasymturnstile}\ \isaconst{iChk}\ \isavar{A}\ \isavar{a}{\isacharsemicolon}{\kern0pt}\ \isavar{{\isasymGamma}}\ {\isasymturnstile}\ \isaconst{iChk}\ \isavar{B}\ \isavar{b}{\isasymrbrakk}\ {\isasymLongrightarrow}\ \isavar{{\isasymGamma}}\ {\isasymturnstile}\ \isaconst{trmPredR}\ {\isacharparenleft}{\kern0pt}\isaconst{cmp}\ \isaconst{opt{\isacharunderscore}{\kern0pt}hd{\isacharunderscore}{\kern0pt}encIdx}\ \isactrlbold {\isasymcdot}\ \isaconst{lstDrop}\ {\isacharparenleft}{\kern0pt}\isaconst{cmp}\ \isaconst{apCoreIdx}\ \isactrlbold {\isasymcdot}\ \isaconst{cpPair}\ \isavar{A}\ \isavar{B}\ \isactrlbold {\isasymcdot}\ {\isacharparenleft}{\kern0pt}\isaconst{cmp}\ \isaconst{length{\isacharunderscore}{\kern0pt}encIdx}\ \isactrlbold {\isasymcdot}\ \isavar{A}\ \isactrlbold {\isasymcdot}\ \isavar{A}{\isacharparenright}{\kern0pt}{\isacharparenright}{\kern0pt}\ {\isacharparenleft}{\kern0pt}\isaconst{gapAdd}\ {\isacharparenleft}{\kern0pt}\isaconst{cmp}\ \isaconst{length{\isacharunderscore}{\kern0pt}encIdx}\ \isactrlbold {\isasymcdot}\ \isavar{A}\ \isactrlbold {\isasymcdot}\ \isavar{A}{\isacharparenright}{\kern0pt}\ \isactrlbold {\isadigit{0}}{\isacharparenright}{\kern0pt}\ \isactrlbold {\isasymcdot}\ \isaconst{lstDrop}\ {\isacharparenleft}{\kern0pt}\isaconst{cmp}\ \isaconst{apCoreIdx}\ \isactrlbold {\isasymcdot}\ \isaconst{cpPair}\ \isavar{A}\ \isavar{B}\ \isactrlbold {\isasymcdot}\ {\isacharparenleft}{\kern0pt}\isaconst{cmp}\ \isaconst{length{\isacharunderscore}{\kern0pt}encIdx}\ \isactrlbold {\isasymcdot}\ \isavar{A}\ \isactrlbold {\isasymcdot}\ \isavar{A}{\isacharparenright}{\kern0pt}{\isacharparenright}{\kern0pt}\ {\isacharparenleft}{\kern0pt}\isaconst{gapAdd}\ {\isacharparenleft}{\kern0pt}\isaconst{cmp}\ \isaconst{length{\isacharunderscore}{\kern0pt}encIdx}\ \isactrlbold {\isasymcdot}\ \isavar{A}\ \isactrlbold {\isasymcdot}\ \isavar{A}{\isacharparenright}{\kern0pt}\ \isactrlbold {\isadigit{0}}{\isacharparenright}{\kern0pt}{\isacharparenright}{\kern0pt}\ \isactrlbold {\isacharequal}{\kern0pt}\ \isaconst{cpPair}\ {\isacharparenleft}{\kern0pt}\isactrlbold n{\isasymlbrace}{\isasymemptyset}{\isasymrbrace}{\isacharparenright}{\kern0pt}\ \isavar{b}\isasep\isanewline%
}%
\SNIP{Goedel-Paper-Statements:thm:iChk-appC-8}{%
{\isasymlbrakk}\isavar{{\isasymGamma}}\ {\isasymturnstile}\ \isavar{A}\ \isactrlbold N{\isacharsemicolon}{\kern0pt}\ \isavar{{\isasymGamma}}\ {\isasymturnstile}\ \isavar{B}\ \isactrlbold N{\isacharsemicolon}{\kern0pt}\ \isavar{{\isasymGamma}}\ {\isasymturnstile}\ \isaconst{iChk}\ \isavar{A}\ \isavar{a}{\isacharsemicolon}{\kern0pt}\ \isavar{{\isasymGamma}}\ {\isasymturnstile}\ \isaconst{iChk}\ \isavar{B}\ \isavar{b}{\isasymrbrakk}\ {\isasymLongrightarrow}\ \isavar{{\isasymGamma}}\ {\isasymturnstile}\ \isaconst{iChk}\ {\isacharparenleft}{\kern0pt}\isaconst{cmp}\ \isaconst{apCoreIdx}\ \isactrlbold {\isasymcdot}\ \isaconst{cpPair}\ \isavar{A}\ \isavar{B}\ \isactrlbold {\isasymcdot}\ {\isacharparenleft}{\kern0pt}\isaconst{cmp}\ \isaconst{length{\isacharunderscore}{\kern0pt}encIdx}\ \isactrlbold {\isasymcdot}\ \isavar{A}\ \isactrlbold {\isasymcdot}\ \isavar{A}{\isacharparenright}{\kern0pt}{\isacharparenright}{\kern0pt}\ \isavar{a}}%
}%
\SNIP{Goedel-Paper-Statements:thm:iChk-lstCons-0}{%
\isakeywordONE{thm}\isamarkupfalse%
\ iChk{\isacharunderscore}{\kern0pt}lstCons%
\isanewline%
}%
\SNIP{Goedel-Paper-Statements:thm:iChk-lstCons-1}{%
\isa{{\isasymlbrakk}\isavar{{\isasymGamma}}\ {\isasymturnstile}\ \isavar{T}\ \isactrlbold N{\isacharsemicolon}{\kern0pt}\ \isavar{{\isasymGamma}}\ {\isasymturnstile}\ \isaconst{cmp}\ \isaconst{validChkIdx}\ \isactrlbold {\isasymcdot}\ \isavar{T}\ \isactrlbold {\isasymcdot}\ \isavar{T}\ \isactrlbold {\isacharequal}{\kern0pt}\ \isactrlbold n\isaconst{Suc}\ {\isadigit{0}}{\isacharsemicolon}{\kern0pt}\ \isavar{{\isasymGamma}}\ {\isasymturnstile}\ \isavar{c}\ \isactrlbold N{\isacharsemicolon}{\kern0pt}\ \isavar{{\isasymGamma}}\ {\isasymturnstile}\ \isaconst{iIsCode}\ \isavar{c}\ \isactrlbold {\isasymnoteq}\ \isactrlbold {\isadigit{0}}{\isacharsemicolon}{\kern0pt}\ \isavar{E}\ {\isasymin}\ \isaconst{set}\ \isaconst{drvChkList}{\isacharsemicolon}{\kern0pt}\ \isavar{{\isasymGamma}}\ {\isasymturnstile}\ \isaconst{cmp}\ \isavar{E}\ \isactrlbold {\isasymcdot}\ \isaconst{cpPair}\ {\isacharparenleft}{\kern0pt}\isaconst{cpSnd}\ {\isacharparenleft}{\kern0pt}\isaconst{cpPair}\ {\isacharparenleft}{\kern0pt}\isaconst{cpPair}\ {\isacharparenleft}{\kern0pt}\isactrlbold n{\isasymlbrace}{\isasymemptyset}{\isasymrbrace}{\isacharparenright}{\kern0pt}\ \isavar{c}{\isacharparenright}{\kern0pt}\ \isavar{T}{\isacharparenright}{\kern0pt}{\isacharparenright}{\kern0pt}\ {\isacharparenleft}{\kern0pt}\isaconst{cpFst}\ {\isacharparenleft}{\kern0pt}\isaconst{cpFst}\ {\isacharparenleft}{\kern0pt}\isaconst{cpPair}\ {\isacharparenleft}{\kern0pt}\isaconst{cpPair}\ {\isacharparenleft}{\kern0pt}\isactrlbold n{\isasymlbrace}{\isasymemptyset}{\isasymrbrace}{\isacharparenright}{\kern0pt}\ \isavar{c}{\isacharparenright}{\kern0pt}\ \isavar{T}{\isacharparenright}{\kern0pt}{\isacharparenright}{\kern0pt}{\isacharparenright}{\kern0pt}\ \isactrlbold {\isasymcdot}\ \isaconst{cpSnd}\ {\isacharparenleft}{\kern0pt}\isaconst{cpFst}\ {\isacharparenleft}{\kern0pt}\isaconst{cpPair}\ {\isacharparenleft}{\kern0pt}\isaconst{cpPair}\ {\isacharparenleft}{\kern0pt}\isactrlbold n{\isasymlbrace}{\isasymemptyset}{\isasymrbrace}{\isacharparenright}{\kern0pt}\ \isavar{c}{\isacharparenright}{\kern0pt}\ \isavar{T}{\isacharparenright}{\kern0pt}{\isacharparenright}{\kern0pt}\ \isactrlbold {\isasymnoteq}\ \isactrlbold {\isadigit{0}}{\isasymrbrakk}\ {\isasymLongrightarrow}\ \isavar{{\isasymGamma}}\ {\isasymturnstile}\ \isaconst{iChk}\ {\isacharparenleft}{\kern0pt}\isaconst{lstCons}\ {\isacharparenleft}{\kern0pt}\isaconst{cpPair}\ {\isacharparenleft}{\kern0pt}\isactrlbold n{\isasymlbrace}{\isasymemptyset}{\isasymrbrace}{\isacharparenright}{\kern0pt}\ \isavar{c}{\isacharparenright}{\kern0pt}\ \isavar{T}{\isacharparenright}{\kern0pt}\ \isavar{c}}%
}%
\SNIP{Goedel-Paper-Statements:thm:goedelEIdx-def-0}{%
\isakeywordONE{thm}\isamarkupfalse%
\ goedelEIdx{\isacharunderscore}{\kern0pt}def%
\isanewline%
}%
\SNIP{Goedel-Paper-Statements:thm:goedelEIdx-def-1}{%
\isa{\isaconst{goedelEIdx}\ {\isasymequiv}\ \isaconst{sucIdx}\ {\isacharparenleft}{\kern0pt}\isaconst{prSub}\ \isaconst{cpairIdx}\ {\isacharparenleft}{\kern0pt}\isaconst{prSub}\ \isaconst{cpairIdx}\ {\isacharparenleft}{\kern0pt}\isaconst{numIdx}\ {\isasymlbrace}{\isasymemptyset}{\isasymrbrace}{\isacharparenright}{\kern0pt}\ {\isacharparenleft}{\kern0pt}\isaconst{numIdx}\ {\isasymlbrace}\isactrlbold {\isadigit{0}}\ \isactrlbold {\isacharequal}{\kern0pt}\ \isactrlbold n\isaconst{Suc}\ {\isadigit{0}}{\isasymrbrace}{\isacharparenright}{\kern0pt}{\isacharparenright}{\kern0pt}\ {\isacharparenleft}{\kern0pt}\isaconst{prSub}\ \isaconst{apCoreIdx}\ {\isacharparenleft}{\kern0pt}\isaconst{prSub}\ \isaconst{cpairIdx}\ {\isacharparenleft}{\kern0pt}\isaconst{prInit}\ \isaconst{rInitL}{\isacharparenright}{\kern0pt}\ {\isacharparenleft}{\kern0pt}\isaconst{prInit}\ \isaconst{rInitR}{\isacharparenright}{\kern0pt}{\isacharparenright}{\kern0pt}\ {\isacharparenleft}{\kern0pt}\isaconst{prSub}\ \isaconst{length{\isacharunderscore}{\kern0pt}encIdx}\ {\isacharparenleft}{\kern0pt}\isaconst{prInit}\ \isaconst{rInitL}{\isacharparenright}{\kern0pt}\ {\isacharparenleft}{\kern0pt}\isaconst{prInit}\ \isaconst{rInitL}{\isacharparenright}{\kern0pt}{\isacharparenright}{\kern0pt}{\isacharparenright}{\kern0pt}{\isacharparenright}{\kern0pt}}%
}%
\SNIP{Goedel-Paper-Statements:thm:goedelE-def-0}{%
\isakeywordONE{thm}\isamarkupfalse%
\ goedelE{\isacharunderscore}{\kern0pt}def%
\isanewline%
}%
\SNIP{Goedel-Paper-Statements:thm:goedelE-def-1}{%
\isa{\isaconst{goedelE}\ {\isasymequiv}\ \isaconst{cmp}\ \isaconst{goedelEIdx}}%
}%
\SNIP{Goedel-Paper-Statements:thm:EXPLF-def-0}{%
\isakeywordONE{thm}\isamarkupfalse%
\ EXPLF{\isacharunderscore}{\kern0pt}def%
\isanewline%
}%
\SNIP{Goedel-Paper-Statements:thm:EXPLF-def-1}{%
\isa{\isaconst{EXPLF}\ \isavar{e}\ {\isasymequiv}\ \isactrlbold A\ \isactrlbold {\isasymcdot}\ {\isacharparenleft}{\kern0pt}\isactrlbold {\isasymlambda}\isactrlbold A\ \isactrlbold {\isasymcdot}\ {\isacharparenleft}{\kern0pt}\isactrlbold {\isasymlambda}{\isacharparenleft}{\kern0pt}\isaconst{iChk}\ {\isacharparenleft}{\kern0pt}\isaconst{trmVar}\ {\isadigit{1}}{\isacharparenright}{\kern0pt}\ {\isacharparenleft}{\kern0pt}\isactrlbold n{\isasymlbrace}\isaconst{goedelG}{\isasymrbrace}{\isacharparenright}{\kern0pt}\ \isactrlbold {\isasymand}\ \isaconst{iChk}\ \isactrlbold v\isactrlsub {\isadigit{0}}\ {\isacharparenleft}{\kern0pt}\isactrlbold n{\isasymlbrace}\isactrlbold {\isasymnot}\isaconst{goedelG}{\isasymrbrace}{\isacharparenright}{\kern0pt}\ \isactrlbold {\isasymlongrightarrow}\ \isaconst{iChk}\ {\isacharparenleft}{\kern0pt}\isavar{e}\ \isactrlbold {\isasymcdot}\ \isaconst{trmVar}\ {\isadigit{1}}\ \isactrlbold {\isasymcdot}\ \isactrlbold v\isactrlsub {\isadigit{0}}{\isacharparenright}{\kern0pt}\ \isaconst{cEq{\isadigit{0}}{\isadigit{1}}}{\isacharparenright}{\kern0pt}{\isacharparenright}{\kern0pt}{\isacharparenright}{\kern0pt}}%
}%
\SNIP{Goedel-Paper-Statements:thm:EXPL-def-0}{%
\isakeywordONE{thm}\isamarkupfalse%
\ EXPL{\isacharunderscore}{\kern0pt}def%
\isanewline%
}%
\SNIP{Goedel-Paper-Statements:thm:EXPL-def-1}{%
\isa{\isaconst{EXPL}\ \isavar{e}\ {\isasymequiv}\ \isaconst{Syntax{\isacharunderscore}{\kern0pt}Term{\isachardot}{\kern0pt}closed}\ \isavar{e}\ {\isasymand}\ {\isasymemptyset}\ {\isasymturnstile}\ \isactrlbold A\ \isactrlbold {\isasymcdot}\ {\isacharparenleft}{\kern0pt}\isactrlbold {\isasymlambda}\isactrlbold A\ \isactrlbold {\isasymcdot}\ {\isacharparenleft}{\kern0pt}\isactrlbold {\isasymlambda}{\isacharparenleft}{\kern0pt}\isavar{e}\ \isactrlbold {\isasymcdot}\ \isaconst{trmVar}\ {\isadigit{1}}\ \isactrlbold {\isasymcdot}\ \isactrlbold v\isactrlsub {\isadigit{0}}\ \isactrlbold N{\isacharparenright}{\kern0pt}{\isacharparenright}{\kern0pt}{\isacharparenright}{\kern0pt}\ {\isasymand}\ {\isasymemptyset}\ {\isasymturnstile}\ \isaconst{EXPLF}\ \isavar{e}}%
}%
\SNIP{Goedel-Paper-Statements:thm:goedelE-eq-0}{%
\isakeywordONE{thm}\isamarkupfalse%
\ goedelE{\isacharunderscore}{\kern0pt}eq%
\isanewline%
}%
\SNIP{Goedel-Paper-Statements:thm:goedelE-eq-1}{%
\isa{{\isasymlbrakk}\isavar{{\isasymGamma}}\ {\isasymturnstile}\ \isavar{a}\ \isactrlbold N{\isacharsemicolon}{\kern0pt}\ \isavar{{\isasymGamma}}\ {\isasymturnstile}\ \isavar{b}\ \isactrlbold N{\isasymrbrakk}\ {\isasymLongrightarrow}\ \isavar{{\isasymGamma}}\ {\isasymturnstile}\ \isaconst{goedelE}\ \isactrlbold {\isasymcdot}\ \isavar{a}\ \isactrlbold {\isasymcdot}\ \isavar{b}\ \isactrlbold {\isacharequal}{\kern0pt}\ \isaconst{lstCons}\ {\isacharparenleft}{\kern0pt}\isaconst{cpPair}\ {\isacharparenleft}{\kern0pt}\isactrlbold n{\isasymlbrace}{\isasymemptyset}{\isasymrbrace}{\isacharparenright}{\kern0pt}\ \isaconst{cEq{\isadigit{0}}{\isadigit{1}}}{\isacharparenright}{\kern0pt}\ {\isacharparenleft}{\kern0pt}\isaconst{cmp}\ \isaconst{apCoreIdx}\ \isactrlbold {\isasymcdot}\ \isaconst{cpPair}\ \isavar{a}\ \isavar{b}\ \isactrlbold {\isasymcdot}\ {\isacharparenleft}{\kern0pt}\isaconst{cmp}\ \isaconst{length{\isacharunderscore}{\kern0pt}encIdx}\ \isactrlbold {\isasymcdot}\ \isavar{a}\ \isactrlbold {\isasymcdot}\ \isavar{a}{\isacharparenright}{\kern0pt}{\isacharparenright}{\kern0pt}}%
}%
\SNIP{Goedel-Paper-Statements:thm:EXPL-antecedent-B-0}{%
\isakeywordONE{thm}\isamarkupfalse%
\ EXPL{\isacharunderscore}{\kern0pt}antecedent{\isacharunderscore}{\kern0pt}B%
\isanewline%
}%
\SNIP{Goedel-Paper-Statements:thm:EXPL-antecedent-B-1}{%
\isa{{\isasymlbrakk}\isavar{{\isasymGamma}}\ {\isasymturnstile}\ \isavar{a}\ \isactrlbold N{\isacharsemicolon}{\kern0pt}\ \isavar{{\isasymGamma}}\ {\isasymturnstile}\ \isavar{b}\ \isactrlbold N{\isasymrbrakk}\ {\isasymLongrightarrow}\ \isavar{{\isasymGamma}}\ {\isasymturnstile}\ \isaconst{iChk}\ \isavar{a}\ {\isacharparenleft}{\kern0pt}\isactrlbold n{\isasymlbrace}\isaconst{goedelG}{\isasymrbrace}{\isacharparenright}{\kern0pt}\ \isactrlbold {\isasymand}\ \isaconst{iChk}\ \isavar{b}\ {\isacharparenleft}{\kern0pt}\isactrlbold n{\isasymlbrace}\isactrlbold {\isasymnot}\isaconst{goedelG}{\isasymrbrace}{\isacharparenright}{\kern0pt}\ \isactrlbold B}%
}%
\SNIP{Goedel-Paper-Statements:thm:D3GF-def-0}{%
\isakeywordONE{thm}\isamarkupfalse%
\ D{\isadigit{3}}GF{\isacharunderscore}{\kern0pt}def%
\isanewline%
}%
\SNIP{Goedel-Paper-Statements:thm:D3GF-def-1}{%
\isa{\isaconst{D{\isadigit{3}}GF}\ \isavar{{\isasymsigma}}\ {\isasymequiv}\ \isactrlbold A\ \isactrlbold {\isasymcdot}\ {\isacharparenleft}{\kern0pt}\isactrlbold {\isasymlambda}{\isacharparenleft}{\kern0pt}\isaconst{iChk}\ \isactrlbold v\isactrlsub {\isadigit{0}}\ {\isacharparenleft}{\kern0pt}\isactrlbold n{\isasymlbrace}\isaconst{goedelG}{\isasymrbrace}{\isacharparenright}{\kern0pt}\ \isactrlbold {\isasymlongrightarrow}\ \isaconst{iChk}\ {\isacharparenleft}{\kern0pt}\isavar{{\isasymsigma}}\ \isactrlbold {\isasymcdot}\ \isactrlbold v\isactrlsub {\isadigit{0}}{\isacharparenright}{\kern0pt}\ {\isacharparenleft}{\kern0pt}\isactrlbold n{\isasymlbrace}\isactrlbold {\isasymnot}\isaconst{goedelG}{\isasymrbrace}{\isacharparenright}{\kern0pt}{\isacharparenright}{\kern0pt}{\isacharparenright}{\kern0pt}}%
}%
\SNIP{Goedel-Paper-Statements:thm:D3G-def-0}{%
\isakeywordONE{thm}\isamarkupfalse%
\ D{\isadigit{3}}G{\isacharunderscore}{\kern0pt}def%
\isanewline%
}%
\SNIP{Goedel-Paper-Statements:thm:D3G-def-1}{%
\isa{\isaconst{D{\isadigit{3}}G}\ \isavar{{\isasymsigma}}\ {\isasymequiv}\ \isaconst{Syntax{\isacharunderscore}{\kern0pt}Term{\isachardot}{\kern0pt}closed}\ \isavar{{\isasymsigma}}\ {\isasymand}\ {\isasymemptyset}\ {\isasymturnstile}\ \isactrlbold A\ \isactrlbold {\isasymcdot}\ {\isacharparenleft}{\kern0pt}\isactrlbold {\isasymlambda}{\isacharparenleft}{\kern0pt}\isavar{{\isasymsigma}}\ \isactrlbold {\isasymcdot}\ \isactrlbold v\isactrlsub {\isadigit{0}}\ \isactrlbold N{\isacharparenright}{\kern0pt}{\isacharparenright}{\kern0pt}\ {\isasymand}\ {\isasymemptyset}\ {\isasymturnstile}\ \isaconst{D{\isadigit{3}}GF}\ \isavar{{\isasymsigma}}}%
}%
\SNIP{Goedel-Paper-Statements:thm:goedelSigmaIdx-def-0}{%
\isakeywordONE{thm}\isamarkupfalse%
\ goedelSigmaIdx{\isacharunderscore}{\kern0pt}def%
\isanewline%
}%
\SNIP{Goedel-Paper-Statements:thm:goedelSigmaIdx-def-1}{%
\isa{\isaconst{goedelSigmaIdx}\ \isavar{cb}\ {\isasymequiv}\ \isaconst{sucIdx}\ {\isacharparenleft}{\kern0pt}\isaconst{prSub}\ \isaconst{cpairIdx}\ {\isacharparenleft}{\kern0pt}\isaconst{prSub}\ \isaconst{cpairIdx}\ {\isacharparenleft}{\kern0pt}\isaconst{numIdx}\ {\isasymlbrace}{\isasymemptyset}{\isasymrbrace}{\isacharparenright}{\kern0pt}\ {\isacharparenleft}{\kern0pt}\isaconst{numIdx}\ {\isasymlbrace}\isactrlbold {\isasymnot}\isaconst{goedelG}{\isasymrbrace}{\isacharparenright}{\kern0pt}{\isacharparenright}{\kern0pt}\ {\isacharparenleft}{\kern0pt}\isaconst{prSub}\ \isaconst{apCoreIdx}\ {\isacharparenleft}{\kern0pt}\isaconst{prSub}\ \isaconst{cpairIdx}\ \isavar{cb}\ \isaconst{natCertIdx}{\isacharparenright}{\kern0pt}\ {\isacharparenleft}{\kern0pt}\isaconst{prSub}\ \isaconst{length{\isacharunderscore}{\kern0pt}encIdx}\ \isavar{cb}\ \isavar{cb}{\isacharparenright}{\kern0pt}{\isacharparenright}{\kern0pt}{\isacharparenright}{\kern0pt}}%
}%
\SNIP{Goedel-Paper-Statements:thm:goedelSigmaOf-def-0}{%
\isakeywordONE{thm}\isamarkupfalse%
\ goedelSigmaOf{\isacharunderscore}{\kern0pt}def%
\isanewline%
}%
\SNIP{Goedel-Paper-Statements:thm:goedelSigmaOf-def-1}{%
\isa{\isaconst{goedelSigmaOf}\ \isavar{cb}\ {\isasymequiv}\ \isactrlbold {\isasymlambda}\isaconst{cmp}\ {\isacharparenleft}{\kern0pt}\isaconst{goedelSigmaIdx}\ \isavar{cb}{\isacharparenright}{\kern0pt}\ \isactrlbold {\isasymcdot}\ \isactrlbold v\isactrlsub {\isadigit{0}}\ \isactrlbold {\isasymcdot}\ \isactrlbold v\isactrlsub {\isadigit{0}}}%
}%
\SNIP{Goedel-Paper-Statements:thm:D3Gcert-def-0}{%
\isakeywordONE{thm}\isamarkupfalse%
\ D{\isadigit{3}}Gcert{\isacharunderscore}{\kern0pt}def%
\isanewline%
}%
\SNIP{Goedel-Paper-Statements:thm:D3Gcert-def-1}{%
\isa{\isaconst{D{\isadigit{3}}Gcert}\ \isavar{cb}\ {\isasymequiv}\ {\isasymforall}\isabound{{\isasymGamma}}\ \isabound{v}{\isachardot}{\kern0pt}\ \isabound{{\isasymGamma}}\ {\isasymturnstile}\ \isabound{v}\ \isactrlbold N\ {\isasymlongrightarrow}\ \isabound{{\isasymGamma}}\ {\isasymturnstile}\ \isaconst{iChk}\ \isabound{v}\ \isaconst{goedelS}\ {\isasymlongrightarrow}\ \isabound{{\isasymGamma}}\ {\isasymturnstile}\ \isaconst{iChk}\ {\isacharparenleft}{\kern0pt}\isaconst{cmp}\ \isavar{cb}\ \isactrlbold {\isasymcdot}\ \isabound{v}\ \isactrlbold {\isasymcdot}\ \isabound{v}{\isacharparenright}{\kern0pt}\ {\isacharparenleft}{\kern0pt}\isaconst{negInstC}\ {\isacharparenleft}{\kern0pt}\isaconst{cmp}\ \isaconst{num{\isacharunderscore}{\kern0pt}codeIdx}\ \isactrlbold {\isasymcdot}\ \isactrlbold {\isadigit{0}}\ \isactrlbold {\isasymcdot}\ \isabound{v}{\isacharparenright}{\kern0pt}{\isacharparenright}{\kern0pt}}%
}%
\SNIP{Goedel-Paper-Statements:thm:natCert-zero-0}{%
\isakeywordONE{thm}\isamarkupfalse%
\ natCert{\isacharunderscore}{\kern0pt}zero%
\isanewline%
}%
\SNIP{Goedel-Paper-Statements:thm:natCert-zero-1}{%
\isa{\isavar{{\isasymGamma}}\ {\isasymturnstile}\ \isavar{a}\ \isactrlbold N\ {\isasymLongrightarrow}\ \isavar{{\isasymGamma}}\ {\isasymturnstile}\ \isaconst{cmp}\ \isaconst{natCertIdx}\ \isactrlbold {\isasymcdot}\ \isavar{a}\ \isactrlbold {\isasymcdot}\ \isactrlbold {\isadigit{0}}\ \isactrlbold {\isacharequal}{\kern0pt}\ \isaconst{lstCons}\ {\isacharparenleft}{\kern0pt}\isaconst{cpPair}\ {\isacharparenleft}{\kern0pt}\isactrlbold n{\isasymlbrace}{\isasymemptyset}{\isasymrbrace}{\isacharparenright}{\kern0pt}\ {\isacharparenleft}{\kern0pt}\isactrlbold n{\isasymlbrace}\isactrlbold {\isadigit{0}}\ \isactrlbold N{\isasymrbrace}{\isacharparenright}{\kern0pt}{\isacharparenright}{\kern0pt}\ \isactrlbold {\isadigit{0}}}%
}%
\SNIP{Goedel-Paper-Statements:thm:natCert-suc-0}{%
\isakeywordONE{thm}\isamarkupfalse%
\ natCert{\isacharunderscore}{\kern0pt}suc%
\isanewline%
}%
\SNIP{Goedel-Paper-Statements:thm:natCert-suc-1}{%
\isa{{\isasymlbrakk}\isavar{{\isasymGamma}}\ {\isasymturnstile}\ \isavar{a}\ \isactrlbold N{\isacharsemicolon}{\kern0pt}\ \isavar{{\isasymGamma}}\ {\isasymturnstile}\ \isavar{k}\ \isactrlbold N{\isacharsemicolon}{\kern0pt}\ \isaconst{fresh}\ \isavar{i}\ \isavar{a}{\isacharsemicolon}{\kern0pt}\ {\isasymnot}\ \isaconst{free}\ \isavar{i}\ \isavar{{\isasymGamma}}{\isasymrbrakk}\ {\isasymLongrightarrow}\ \isavar{{\isasymGamma}}\ {\isasymturnstile}\ \isaconst{cmp}\ \isaconst{natCertIdx}\ \isactrlbold {\isasymcdot}\ \isavar{a}\ \isactrlbold {\isasymcdot}\ {\isacharparenleft}{\kern0pt}\isactrlbold S\ \isavar{k}{\isacharparenright}{\kern0pt}\ \isactrlbold {\isacharequal}{\kern0pt}\ \isaconst{lstCons}\ {\isacharparenleft}{\kern0pt}\isaconst{cpPair}\ {\isacharparenleft}{\kern0pt}\isactrlbold n{\isasymlbrace}{\isasymemptyset}{\isasymrbrace}{\isacharparenright}{\kern0pt}\ {\isacharparenleft}{\kern0pt}\isaconst{mkEqS}\ {\isacharparenleft}{\kern0pt}\isaconst{cmp}\ \isaconst{num{\isacharunderscore}{\kern0pt}codeIdx}\ \isactrlbold {\isasymcdot}\ \isactrlbold {\isadigit{0}}\ \isactrlbold {\isasymcdot}\ {\isacharparenleft}{\kern0pt}\isactrlbold S\ \isavar{k}{\isacharparenright}{\kern0pt}{\isacharparenright}{\kern0pt}\ {\isacharparenleft}{\kern0pt}\isaconst{cmp}\ \isaconst{num{\isacharunderscore}{\kern0pt}codeIdx}\ \isactrlbold {\isasymcdot}\ \isactrlbold {\isadigit{0}}\ \isactrlbold {\isasymcdot}\ {\isacharparenleft}{\kern0pt}\isactrlbold S\ \isavar{k}{\isacharparenright}{\kern0pt}{\isacharparenright}{\kern0pt}{\isacharparenright}{\kern0pt}{\isacharparenright}{\kern0pt}\ {\isacharparenleft}{\kern0pt}\isaconst{cmp}\ \isaconst{natCertIdx}\ \isactrlbold {\isasymcdot}\ \isavar{a}\ \isactrlbold {\isasymcdot}\ \isavar{k}{\isacharparenright}{\kern0pt}}%
}%
\SNIP{Goedel-Paper-Statements:thm:iChk-natCert-0}{%
\isakeywordONE{thm}\isamarkupfalse%
\ iChk{\isacharunderscore}{\kern0pt}natCert%
\isanewline%
}%
\SNIP{Goedel-Paper-Statements:thm:iChk-natCert-1}{%
\isa{{\isasymlbrakk}\isavar{{\isasymGamma}}\ {\isasymturnstile}\ \isavar{a}\ \isactrlbold N{\isacharsemicolon}{\kern0pt}\ \isavar{{\isasymGamma}}\ {\isasymturnstile}\ \isavar{k}\ \isactrlbold N{\isasymrbrakk}\ {\isasymLongrightarrow}\ \isavar{{\isasymGamma}}\ {\isasymturnstile}\ \isaconst{iChk}\ {\isacharparenleft}{\kern0pt}\isaconst{cmp}\ \isaconst{natCertIdx}\ \isactrlbold {\isasymcdot}\ \isavar{a}\ \isactrlbold {\isasymcdot}\ \isavar{k}{\isacharparenright}{\kern0pt}\ {\isacharparenleft}{\kern0pt}\isaconst{mkEqS}\ {\isacharparenleft}{\kern0pt}\isaconst{cmp}\ \isaconst{num{\isacharunderscore}{\kern0pt}codeIdx}\ \isactrlbold {\isasymcdot}\ \isactrlbold {\isadigit{0}}\ \isactrlbold {\isasymcdot}\ \isavar{k}{\isacharparenright}{\kern0pt}\ {\isacharparenleft}{\kern0pt}\isaconst{cmp}\ \isaconst{num{\isacharunderscore}{\kern0pt}codeIdx}\ \isactrlbold {\isasymcdot}\ \isactrlbold {\isadigit{0}}\ \isactrlbold {\isasymcdot}\ \isavar{k}{\isacharparenright}{\kern0pt}{\isacharparenright}{\kern0pt}}%
}%
\SNIP{Goedel-Paper-Statements:thm:iChk-negG-notAllI-0}{%
\isakeywordONE{thm}\isamarkupfalse%
\ iChk{\isacharunderscore}{\kern0pt}negG{\isacharunderscore}{\kern0pt}notAllI%
\isanewline%
}%
\SNIP{Goedel-Paper-Statements:thm:iChk-negG-notAllI-1}{%
\isa{{\isasymlbrakk}\isavar{{\isasymGamma}}\ {\isasymturnstile}\ \isavar{a}\ \isactrlbold N{\isacharsemicolon}{\kern0pt}\ \isavar{{\isasymGamma}}\ {\isasymturnstile}\ \isavar{P{\isadigit{1}}{\isachardot}{\kern0pt}{\isadigit{0}}}\ \isactrlbold N{\isacharsemicolon}{\kern0pt}\ \isavar{{\isasymGamma}}\ {\isasymturnstile}\ \isavar{P{\isadigit{2}}{\isachardot}{\kern0pt}{\isadigit{0}}}\ \isactrlbold N{\isacharsemicolon}{\kern0pt}\ \isavar{{\isasymGamma}}\ {\isasymturnstile}\ \isaconst{iChk}\ \isavar{P{\isadigit{1}}{\isachardot}{\kern0pt}{\isadigit{0}}}\ {\isacharparenleft}{\kern0pt}\isaconst{negInstC}\ \isavar{a}{\isacharparenright}{\kern0pt}{\isacharsemicolon}{\kern0pt}\ \isavar{{\isasymGamma}}\ {\isasymturnstile}\ \isaconst{iChk}\ \isavar{P{\isadigit{2}}{\isachardot}{\kern0pt}{\isadigit{0}}}\ {\isacharparenleft}{\kern0pt}\isaconst{eqReflC}\ \isavar{a}{\isacharparenright}{\kern0pt}{\isasymrbrakk}\ {\isasymLongrightarrow}\ \isavar{{\isasymGamma}}\ {\isasymturnstile}\ \isaconst{iChk}\ {\isacharparenleft}{\kern0pt}\isaconst{lstCons}\ {\isacharparenleft}{\kern0pt}\isaconst{cpPair}\ {\isacharparenleft}{\kern0pt}\isactrlbold n{\isasymlbrace}{\isasymemptyset}{\isasymrbrace}{\isacharparenright}{\kern0pt}\ {\isacharparenleft}{\kern0pt}\isactrlbold n{\isasymlbrace}\isactrlbold {\isasymnot}\isaconst{goedelG}{\isasymrbrace}{\isacharparenright}{\kern0pt}{\isacharparenright}{\kern0pt}\ {\isacharparenleft}{\kern0pt}\isaconst{cmp}\ \isaconst{apCoreIdx}\ \isactrlbold {\isasymcdot}\ \isaconst{cpPair}\ \isavar{P{\isadigit{1}}{\isachardot}{\kern0pt}{\isadigit{0}}}\ \isavar{P{\isadigit{2}}{\isachardot}{\kern0pt}{\isadigit{0}}}\ \isactrlbold {\isasymcdot}\ {\isacharparenleft}{\kern0pt}\isaconst{cmp}\ \isaconst{length{\isacharunderscore}{\kern0pt}encIdx}\ \isactrlbold {\isasymcdot}\ \isavar{P{\isadigit{1}}{\isachardot}{\kern0pt}{\isadigit{0}}}\ \isactrlbold {\isasymcdot}\ \isavar{P{\isadigit{1}}{\isachardot}{\kern0pt}{\isadigit{0}}}{\isacharparenright}{\kern0pt}{\isacharparenright}{\kern0pt}{\isacharparenright}{\kern0pt}\ {\isacharparenleft}{\kern0pt}\isactrlbold n{\isasymlbrace}\isactrlbold {\isasymnot}\isaconst{goedelG}{\isasymrbrace}{\isacharparenright}{\kern0pt}}%
}%
\SNIP{Goedel-Paper-Statements:thm:goedelSigmaIdx-eq-0}{%
\isakeywordONE{thm}\isamarkupfalse%
\ goedelSigmaIdx{\isacharunderscore}{\kern0pt}eq%
\isanewline%
}%
\SNIP{Goedel-Paper-Statements:thm:goedelSigmaIdx-eq-1}{%
\isa{{\isasymlbrakk}\isavar{{\isasymGamma}}\ {\isasymturnstile}\ \isavar{a}\ \isactrlbold N{\isacharsemicolon}{\kern0pt}\ \isavar{{\isasymGamma}}\ {\isasymturnstile}\ \isavar{b}\ \isactrlbold N{\isasymrbrakk}\ {\isasymLongrightarrow}\ \isavar{{\isasymGamma}}\ {\isasymturnstile}\ \isaconst{cmp}\ {\isacharparenleft}{\kern0pt}\isaconst{goedelSigmaIdx}\ \isavar{cb}{\isacharparenright}{\kern0pt}\ \isactrlbold {\isasymcdot}\ \isavar{a}\ \isactrlbold {\isasymcdot}\ \isavar{b}\ \isactrlbold {\isacharequal}{\kern0pt}\ \isaconst{lstCons}\ {\isacharparenleft}{\kern0pt}\isaconst{cpPair}\ {\isacharparenleft}{\kern0pt}\isactrlbold n{\isasymlbrace}{\isasymemptyset}{\isasymrbrace}{\isacharparenright}{\kern0pt}\ {\isacharparenleft}{\kern0pt}\isactrlbold n{\isasymlbrace}\isactrlbold {\isasymnot}\isaconst{goedelG}{\isasymrbrace}{\isacharparenright}{\kern0pt}{\isacharparenright}{\kern0pt}\ {\isacharparenleft}{\kern0pt}\isaconst{cmp}\ \isaconst{apCoreIdx}\ \isactrlbold {\isasymcdot}\ \isaconst{cpPair}\ {\isacharparenleft}{\kern0pt}\isaconst{cmp}\ \isavar{cb}\ \isactrlbold {\isasymcdot}\ \isavar{a}\ \isactrlbold {\isasymcdot}\ \isavar{b}{\isacharparenright}{\kern0pt}\ {\isacharparenleft}{\kern0pt}\isaconst{cmp}\ \isaconst{natCertIdx}\ \isactrlbold {\isasymcdot}\ \isavar{a}\ \isactrlbold {\isasymcdot}\ \isavar{b}{\isacharparenright}{\kern0pt}\ \isactrlbold {\isasymcdot}\ {\isacharparenleft}{\kern0pt}\isaconst{cmp}\ \isaconst{length{\isacharunderscore}{\kern0pt}encIdx}\ \isactrlbold {\isasymcdot}\ {\isacharparenleft}{\kern0pt}\isaconst{cmp}\ \isavar{cb}\ \isactrlbold {\isasymcdot}\ \isavar{a}\ \isactrlbold {\isasymcdot}\ \isavar{b}{\isacharparenright}{\kern0pt}\ \isactrlbold {\isasymcdot}\ {\isacharparenleft}{\kern0pt}\isaconst{cmp}\ \isavar{cb}\ \isactrlbold {\isasymcdot}\ \isavar{a}\ \isactrlbold {\isasymcdot}\ \isavar{b}{\isacharparenright}{\kern0pt}{\isacharparenright}{\kern0pt}{\isacharparenright}{\kern0pt}}%
}%
\SNIP{Goedel-Paper-Statements:thm:iChk-goedelSigmaOf-0}{%
\isakeywordONE{thm}\isamarkupfalse%
\ iChk{\isacharunderscore}{\kern0pt}goedelSigmaOf%
\isanewline%
}%
\SNIP{Goedel-Paper-Statements:thm:iChk-goedelSigmaOf-1}{%
\isa{{\isasymlbrakk}\isavar{{\isasymGamma}}\ {\isasymturnstile}\ \isavar{v}\ \isactrlbold N{\isacharsemicolon}{\kern0pt}\ \isavar{{\isasymGamma}}\ {\isasymturnstile}\ \isaconst{iChk}\ {\isacharparenleft}{\kern0pt}\isaconst{cmp}\ \isavar{cb}\ \isactrlbold {\isasymcdot}\ \isavar{v}\ \isactrlbold {\isasymcdot}\ \isavar{v}{\isacharparenright}{\kern0pt}\ {\isacharparenleft}{\kern0pt}\isaconst{negInstC}\ {\isacharparenleft}{\kern0pt}\isaconst{cmp}\ \isaconst{num{\isacharunderscore}{\kern0pt}codeIdx}\ \isactrlbold {\isasymcdot}\ \isactrlbold {\isadigit{0}}\ \isactrlbold {\isasymcdot}\ \isavar{v}{\isacharparenright}{\kern0pt}{\isacharparenright}{\kern0pt}{\isasymrbrakk}\ {\isasymLongrightarrow}\ \isavar{{\isasymGamma}}\ {\isasymturnstile}\ \isaconst{iChk}\ {\isacharparenleft}{\kern0pt}\isaconst{goedelSigmaOf}\ \isavar{cb}\ \isactrlbold {\isasymcdot}\ \isavar{v}{\isacharparenright}{\kern0pt}\ {\isacharparenleft}{\kern0pt}\isactrlbold n{\isasymlbrace}\isactrlbold {\isasymnot}\isaconst{goedelG}{\isasymrbrace}{\isacharparenright}{\kern0pt}}%
}%
\SNIP{Goedel-Paper-Statements:thm:D3G-of-d3gcert-0}{%
\isakeywordONE{thm}\isamarkupfalse%
\ D{\isadigit{3}}G{\isacharunderscore}{\kern0pt}of{\isacharunderscore}{\kern0pt}d{\isadigit{3}}gcert%
\isanewline%
}%
\SNIP{Goedel-Paper-Statements:thm:D3G-of-d3gcert-1}{%
\isa{\isaconst{D{\isadigit{3}}Gcert}\ \isavar{cb}\ {\isasymLongrightarrow}\ \isaconst{D{\isadigit{3}}G}\ {\isacharparenleft}{\kern0pt}\isaconst{goedelSigmaOf}\ \isavar{cb}{\isacharparenright}{\kern0pt}}%
}%
\SNIP{Goedel-Paper-Statements:thm:prCertBuilder-0}{%
\isakeywordONE{thm}\isamarkupfalse%
\ prCertBuilder{\isachardot}{\kern0pt}simps%
\isanewline%
}%
\SNIP{Goedel-Paper-Statements:thm:prCertBuilder-1}{%
\isa{\isaconst{prCertBuilder}\ {\isacharparenleft}{\kern0pt}\isaconst{prInit}\ \isavar{i}{\isacharparenright}{\kern0pt}\ {\isacharequal}{\kern0pt}\ \isaconst{bldTop}\ {\isacharparenleft}{\kern0pt}\isaconst{prInit}\ \isavar{i}{\isacharparenright}{\kern0pt}\ {\isacharparenleft}{\kern0pt}\isaconst{bldInit}\ \isavar{i}{\isacharparenright}{\kern0pt}\isasep\isanewline%
}%
\SNIP{Goedel-Paper-Statements:thm:prCertBuilder-2}{%
\isaconst{prCertBuilder}\ {\isacharparenleft}{\kern0pt}\isaconst{prSub}\ \isavar{g}\ \isavar{h{\isadigit{1}}{\isachardot}{\kern0pt}{\isadigit{0}}}\ \isavar{h{\isadigit{2}}{\isachardot}{\kern0pt}{\isadigit{0}}}{\isacharparenright}{\kern0pt}\ {\isacharequal}{\kern0pt}\ \isaconst{bldTop}\ {\isacharparenleft}{\kern0pt}\isaconst{prSub}\ \isavar{g}\ \isavar{h{\isadigit{1}}{\isachardot}{\kern0pt}{\isadigit{0}}}\ \isavar{h{\isadigit{2}}{\isachardot}{\kern0pt}{\isadigit{0}}}{\isacharparenright}{\kern0pt}\ {\isacharparenleft}{\kern0pt}\isaconst{prSub}\ {\isacharparenleft}{\kern0pt}\isaconst{prCertBuilder}\ \isavar{g}{\isacharparenright}{\kern0pt}\ {\isacharparenleft}{\kern0pt}\isaconst{prSub}\ \isaconst{cpairIdx}\ \isaconst{bldRed{\isadigit{2}}}\ {\isacharparenleft}{\kern0pt}\isaconst{prSub}\ \isaconst{cpairIdx}\ {\isacharparenleft}{\kern0pt}\isaconst{prSub}\ \isavar{h{\isadigit{1}}{\isachardot}{\kern0pt}{\isadigit{0}}}\ \isaconst{bldInX}\ \isaconst{bldInY}{\isacharparenright}{\kern0pt}\ {\isacharparenleft}{\kern0pt}\isaconst{prSub}\ \isavar{h{\isadigit{2}}{\isachardot}{\kern0pt}{\isadigit{0}}}\ \isaconst{bldInX}\ \isaconst{bldInY}{\isacharparenright}{\kern0pt}{\isacharparenright}{\kern0pt}{\isacharparenright}{\kern0pt}\ {\isacharparenleft}{\kern0pt}\isaconst{prSub}\ \isaconst{cpairIdx}\ {\isacharparenleft}{\kern0pt}\isaconst{prSub}\ {\isacharparenleft}{\kern0pt}\isaconst{prCertBuilder}\ \isavar{h{\isadigit{1}}{\isachardot}{\kern0pt}{\isadigit{0}}}{\isacharparenright}{\kern0pt}\ {\isacharparenleft}{\kern0pt}\isaconst{prSub}\ \isaconst{cpairIdx}\ {\isacharparenleft}{\kern0pt}\isaconst{kid{\isadigit{2}}Idx}\ {\isacharparenleft}{\kern0pt}\isaconst{kid{\isadigit{1}}Idx}\ \isaconst{bldRed{\isadigit{2}}}{\isacharparenright}{\kern0pt}{\isacharparenright}{\kern0pt}\ \isaconst{bldInXY}{\isacharparenright}{\kern0pt}\ {\isacharparenleft}{\kern0pt}\isaconst{prInit}\ \isaconst{rInitR}{\isacharparenright}{\kern0pt}{\isacharparenright}{\kern0pt}\ {\isacharparenleft}{\kern0pt}\isaconst{prSub}\ {\isacharparenleft}{\kern0pt}\isaconst{prCertBuilder}\ \isavar{h{\isadigit{2}}{\isachardot}{\kern0pt}{\isadigit{0}}}{\isacharparenright}{\kern0pt}\ {\isacharparenleft}{\kern0pt}\isaconst{prSub}\ \isaconst{cpairIdx}\ {\isacharparenleft}{\kern0pt}\isaconst{kid{\isadigit{2}}Idx}\ \isaconst{bldRed{\isadigit{2}}}{\isacharparenright}{\kern0pt}\ \isaconst{bldInXY}{\isacharparenright}{\kern0pt}\ {\isacharparenleft}{\kern0pt}\isaconst{prInit}\ \isaconst{rInitR}{\isacharparenright}{\kern0pt}{\isacharparenright}{\kern0pt}{\isacharparenright}{\kern0pt}{\isacharparenright}{\kern0pt}\isasep\isanewline%
}%
\SNIP{Goedel-Paper-Statements:thm:prCertBuilder-3}{%
\isaconst{prCertBuilder}\ {\isacharparenleft}{\kern0pt}\isaconst{prRec}\ \isavar{g}\ \isavar{h}\ \isavar{h{\isacharprime}{\kern0pt}}{\isacharparenright}{\kern0pt}\ {\isacharequal}{\kern0pt}\ \isaconst{bldTop}\ {\isacharparenleft}{\kern0pt}\isaconst{prRec}\ \isavar{g}\ \isavar{h}\ \isavar{h{\isacharprime}{\kern0pt}}{\isacharparenright}{\kern0pt}\ {\isacharparenleft}{\kern0pt}\isaconst{certSpineIdx}\ {\isacharparenleft}{\kern0pt}\isaconst{kid{\isadigit{2}}Idx}\ {\isacharparenleft}{\kern0pt}\isaconst{kid{\isadigit{1}}Idx}\ {\isacharparenleft}{\kern0pt}\isaconst{kid{\isadigit{1}}Idx}\ \isaconst{bldRed{\isadigit{2}}}{\isacharparenright}{\kern0pt}{\isacharparenright}{\kern0pt}{\isacharparenright}{\kern0pt}\ {\isacharparenleft}{\kern0pt}\isaconst{kid{\isadigit{2}}Idx}\ {\isacharparenleft}{\kern0pt}\isaconst{kid{\isadigit{1}}Idx}\ \isaconst{bldRed{\isadigit{2}}}{\isacharparenright}{\kern0pt}{\isacharparenright}{\kern0pt}\ {\isacharparenleft}{\kern0pt}\isaconst{numCIdx}\ \isaconst{bldInY}{\isacharparenright}{\kern0pt}\ {\isacharparenleft}{\kern0pt}\isaconst{kid{\isadigit{2}}Idx}\ \isaconst{bldRed{\isadigit{2}}}{\isacharparenright}{\kern0pt}\ {\isacharparenleft}{\kern0pt}\isaconst{bldValIdx}\ {\isacharparenleft}{\kern0pt}\isaconst{prRec}\ \isavar{g}\ \isavar{h}\ \isavar{h{\isacharprime}{\kern0pt}}{\isacharparenright}{\kern0pt}{\isacharparenright}{\kern0pt}\ {\isacharparenleft}{\kern0pt}\isaconst{prSub}\ {\isacharparenleft}{\kern0pt}\isaconst{bldRecIdx}\ {\isacharparenleft}{\kern0pt}\isaconst{prRec}\ \isavar{g}\ \isavar{h}\ \isavar{h{\isacharprime}{\kern0pt}}{\isacharparenright}{\kern0pt}\ \isavar{h{\isacharprime}{\kern0pt}}\ {\isacharparenleft}{\kern0pt}\isaconst{prCertBuilder}\ \isavar{g}{\isacharparenright}{\kern0pt}\ {\isacharparenleft}{\kern0pt}\isaconst{prCertBuilder}\ \isavar{h}{\isacharparenright}{\kern0pt}\ {\isacharparenleft}{\kern0pt}\isaconst{prCertBuilder}\ \isavar{h{\isacharprime}{\kern0pt}}{\isacharparenright}{\kern0pt}{\isacharparenright}{\kern0pt}\ {\isacharparenleft}{\kern0pt}\isaconst{prSub}\ \isaconst{cpairIdx}\ {\isacharparenleft}{\kern0pt}\isaconst{prInit}\ \isaconst{rInitL}{\isacharparenright}{\kern0pt}\ {\isacharparenleft}{\kern0pt}\isaconst{prInit}\ \isaconst{rInitR}{\isacharparenright}{\kern0pt}{\isacharparenright}{\kern0pt}\ \isaconst{bldInY}{\isacharparenright}{\kern0pt}\ \isaconst{bldInP{\isadigit{2}}}{\isacharparenright}{\kern0pt}}%
}%
\SNIP{Goedel-Paper-Statements:thm:prCertBuilderOK-def-0}{%
\isakeywordONE{thm}\isamarkupfalse%
\ prCertBuilderOK{\isacharunderscore}{\kern0pt}def%
\isanewline%
}%
\SNIP{Goedel-Paper-Statements:thm:prCertBuilderOK-def-1}{%
\isa{\isaconst{prCertBuilderOK}\ \isavar{f}\ {\isasymequiv}\ {\isasymforall}\isabound{{\isasymGamma}}\ \isabound{z}\ \isabound{T{\isadigit{1}}}\ \isabound{T{\isadigit{2}}}\ \isabound{x}\ \isabound{y}\ \isabound{Pa}\ \isabound{Pb}{\isachardot}{\kern0pt}\ \isabound{{\isasymGamma}}\ {\isasymturnstile}\ \isabound{z}\ \isactrlbold N\ {\isasymlongrightarrow}\ \isabound{{\isasymGamma}}\ {\isasymturnstile}\ \isaconst{sigWfGnChk}\ \isactrlbold {\isasymcdot}\ \isabound{z}\ {\isasymlongrightarrow}\ \isaconst{wft}\ \isabound{T{\isadigit{1}}}\ {\isasymlongrightarrow}\ \isaconst{wft}\ \isabound{T{\isadigit{2}}}\ {\isasymlongrightarrow}\ \isabound{{\isasymGamma}}\ {\isasymturnstile}\ \isabound{x}\ \isactrlbold N\ {\isasymlongrightarrow}\ \isabound{{\isasymGamma}}\ {\isasymturnstile}\ \isabound{y}\ \isactrlbold N\ {\isasymlongrightarrow}\ \isabound{{\isasymGamma}}\ {\isasymturnstile}\ \isabound{Pa}\ \isactrlbold N\ {\isasymlongrightarrow}\ \isabound{{\isasymGamma}}\ {\isasymturnstile}\ \isabound{Pb}\ \isactrlbold N\ {\isasymlongrightarrow}\ \isabound{{\isasymGamma}}\ {\isasymturnstile}\ \isaconst{iChk}\ \isabound{Pa}\ {\isacharparenleft}{\kern0pt}\isaconst{mkEqS}\ {\isacharparenleft}{\kern0pt}\isaconst{iSubM}\ \isabound{z}\ {\isacharparenleft}{\kern0pt}\isactrlbold n{\isasymlbrace}\isabound{T{\isadigit{1}}}{\isasymrbrace}{\isacharparenright}{\kern0pt}{\isacharparenright}{\kern0pt}\ {\isacharparenleft}{\kern0pt}\isaconst{cmp}\ \isaconst{num{\isacharunderscore}{\kern0pt}codeIdx}\ \isactrlbold {\isasymcdot}\ \isactrlbold {\isadigit{0}}\ \isactrlbold {\isasymcdot}\ \isabound{x}{\isacharparenright}{\kern0pt}{\isacharparenright}{\kern0pt}\ {\isasymlongrightarrow}\ \isabound{{\isasymGamma}}\ {\isasymturnstile}\ \isaconst{iChk}\ \isabound{Pb}\ {\isacharparenleft}{\kern0pt}\isaconst{mkEqS}\ {\isacharparenleft}{\kern0pt}\isaconst{iSubM}\ \isabound{z}\ {\isacharparenleft}{\kern0pt}\isactrlbold n{\isasymlbrace}\isabound{T{\isadigit{2}}}{\isasymrbrace}{\isacharparenright}{\kern0pt}{\isacharparenright}{\kern0pt}\ {\isacharparenleft}{\kern0pt}\isaconst{cmp}\ \isaconst{num{\isacharunderscore}{\kern0pt}codeIdx}\ \isactrlbold {\isasymcdot}\ \isactrlbold {\isadigit{0}}\ \isactrlbold {\isasymcdot}\ \isabound{y}{\isacharparenright}{\kern0pt}{\isacharparenright}{\kern0pt}\ {\isasymlongrightarrow}\ \isabound{{\isasymGamma}}\ {\isasymturnstile}\ \isaconst{iChk}\ {\isacharparenleft}{\kern0pt}\isaconst{cmp}\ {\isacharparenleft}{\kern0pt}\isaconst{prCertBuilder}\ \isavar{f}{\isacharparenright}{\kern0pt}\ \isactrlbold {\isasymcdot}\ \isaconst{cpPair}\ {\isacharparenleft}{\kern0pt}\isaconst{iSubM}\ \isabound{z}\ {\isacharparenleft}{\kern0pt}\isactrlbold n{\isasymlbrace}\isaconst{cmp}\ \isavar{f}\ \isactrlbold {\isasymcdot}\ \isabound{T{\isadigit{1}}}\ \isactrlbold {\isasymcdot}\ \isabound{T{\isadigit{2}}}{\isasymrbrace}{\isacharparenright}{\kern0pt}{\isacharparenright}{\kern0pt}\ {\isacharparenleft}{\kern0pt}\isaconst{cpPair}\ \isabound{x}\ \isabound{y}{\isacharparenright}{\kern0pt}\ \isactrlbold {\isasymcdot}\ \isaconst{cpPair}\ \isabound{Pa}\ \isabound{Pb}{\isacharparenright}{\kern0pt}\ {\isacharparenleft}{\kern0pt}\isaconst{mkEqS}\ {\isacharparenleft}{\kern0pt}\isaconst{iSubM}\ \isabound{z}\ {\isacharparenleft}{\kern0pt}\isactrlbold n{\isasymlbrace}\isaconst{cmp}\ \isavar{f}\ \isactrlbold {\isasymcdot}\ \isabound{T{\isadigit{1}}}\ \isactrlbold {\isasymcdot}\ \isabound{T{\isadigit{2}}}{\isasymrbrace}{\isacharparenright}{\kern0pt}{\isacharparenright}{\kern0pt}\ {\isacharparenleft}{\kern0pt}\isaconst{cmp}\ \isaconst{num{\isacharunderscore}{\kern0pt}codeIdx}\ \isactrlbold {\isasymcdot}\ \isactrlbold {\isadigit{0}}\ \isactrlbold {\isasymcdot}\ {\isacharparenleft}{\kern0pt}\isaconst{cmp}\ \isavar{f}\ \isactrlbold {\isasymcdot}\ \isabound{x}\ \isactrlbold {\isasymcdot}\ \isabound{y}{\isacharparenright}{\kern0pt}{\isacharparenright}{\kern0pt}{\isacharparenright}{\kern0pt}}%
}%
\SNIP{Goedel-Paper-Statements:thm:pr-cert-builder-0}{%
\isakeywordONE{thm}\isamarkupfalse%
\ pr{\isacharunderscore}{\kern0pt}cert{\isacharunderscore}{\kern0pt}builder%
\isanewline%
}%
\SNIP{Goedel-Paper-Statements:thm:pr-cert-builder-1}{%
\isa{\isaconst{prCertBuilderOK}\ \isavar{f}}%
}%
\SNIP{Goedel-Paper-Statements:thm:bldRec-ind-0}{%
\isakeywordONE{thm}\isamarkupfalse%
\ bldRec{\isacharunderscore}{\kern0pt}ind%
\isanewline%
}%
\SNIP{Goedel-Paper-Statements:thm:bldRec-ind-1}{%
\isa{{\isasymlbrakk}\isaconst{prCertBuilderOK}\ \isavar{g}{\isacharsemicolon}{\kern0pt}\ \isaconst{prCertBuilderOK}\ \isavar{h}{\isacharsemicolon}{\kern0pt}\ \isaconst{prCertBuilderOK}\ \isavar{h{\isacharprime}{\kern0pt}}{\isacharsemicolon}{\kern0pt}\ \isavar{{\isasymGamma}}\ {\isasymturnstile}\ \isavar{z}\ \isactrlbold N{\isacharsemicolon}{\kern0pt}\ \isavar{{\isasymGamma}}\ {\isasymturnstile}\ \isaconst{sigWfGnChk}\ \isactrlbold {\isasymcdot}\ \isavar{z}{\isacharsemicolon}{\kern0pt}\ \isaconst{wft}\ \isavar{T{\isadigit{1}}{\isachardot}{\kern0pt}{\isadigit{0}}}{\isacharsemicolon}{\kern0pt}\ \isaconst{wft}\ \isavar{T{\isadigit{2}}{\isachardot}{\kern0pt}{\isadigit{0}}}{\isacharsemicolon}{\kern0pt}\ \isavar{{\isasymGamma}}\ {\isasymturnstile}\ \isavar{x}\ \isactrlbold N{\isacharsemicolon}{\kern0pt}\ \isavar{{\isasymGamma}}\ {\isasymturnstile}\ \isavar{y}\ \isactrlbold N{\isacharsemicolon}{\kern0pt}\ \isavar{{\isasymGamma}}\ {\isasymturnstile}\ \isavar{Pa}\ \isactrlbold N{\isacharsemicolon}{\kern0pt}\ \isavar{{\isasymGamma}}\ {\isasymturnstile}\ \isavar{Pb}\ \isactrlbold N{\isacharsemicolon}{\kern0pt}\ \isavar{{\isasymGamma}}\ {\isasymturnstile}\ \isaconst{iChk}\ \isavar{Pa}\ {\isacharparenleft}{\kern0pt}\isaconst{mkEqS}\ {\isacharparenleft}{\kern0pt}\isaconst{iSubM}\ \isavar{z}\ {\isacharparenleft}{\kern0pt}\isactrlbold n{\isasymlbrace}\isavar{T{\isadigit{1}}{\isachardot}{\kern0pt}{\isadigit{0}}}{\isasymrbrace}{\isacharparenright}{\kern0pt}{\isacharparenright}{\kern0pt}\ {\isacharparenleft}{\kern0pt}\isaconst{cmp}\ \isaconst{num{\isacharunderscore}{\kern0pt}codeIdx}\ \isactrlbold {\isasymcdot}\ \isactrlbold {\isadigit{0}}\ \isactrlbold {\isasymcdot}\ \isavar{x}{\isacharparenright}{\kern0pt}{\isacharparenright}{\kern0pt}{\isacharsemicolon}{\kern0pt}\ \isavar{{\isasymGamma}}\ {\isasymturnstile}\ \isavar{q}\ \isactrlbold N{\isasymrbrakk}\ {\isasymLongrightarrow}\ \isavar{{\isasymGamma}}\ {\isasymturnstile}\ \isaconst{iChk}\ {\isacharparenleft}{\kern0pt}\isaconst{cmp}\ {\isacharparenleft}{\kern0pt}\isaconst{bldRecIdx}\ {\isacharparenleft}{\kern0pt}\isaconst{prRec}\ \isavar{g}\ \isavar{h}\ \isavar{h{\isacharprime}{\kern0pt}}{\isacharparenright}{\kern0pt}\ \isavar{h{\isacharprime}{\kern0pt}}\ {\isacharparenleft}{\kern0pt}\isaconst{prCertBuilder}\ \isavar{g}{\isacharparenright}{\kern0pt}\ {\isacharparenleft}{\kern0pt}\isaconst{prCertBuilder}\ \isavar{h}{\isacharparenright}{\kern0pt}\ {\isacharparenleft}{\kern0pt}\isaconst{prCertBuilder}\ \isavar{h{\isacharprime}{\kern0pt}}{\isacharparenright}{\kern0pt}{\isacharparenright}{\kern0pt}\ \isactrlbold {\isasymcdot}\ \isaconst{cpPair}\ {\isacharparenleft}{\kern0pt}\isaconst{cpPair}\ {\isacharparenleft}{\kern0pt}\isaconst{iSubM}\ \isavar{z}\ {\isacharparenleft}{\kern0pt}\isactrlbold n{\isasymlbrace}\isaconst{cmp}\ {\isacharparenleft}{\kern0pt}\isaconst{prRec}\ \isavar{g}\ \isavar{h}\ \isavar{h{\isacharprime}{\kern0pt}}{\isacharparenright}{\kern0pt}\ \isactrlbold {\isasymcdot}\ \isavar{T{\isadigit{1}}{\isachardot}{\kern0pt}{\isadigit{0}}}\ \isactrlbold {\isasymcdot}\ \isavar{T{\isadigit{2}}{\isachardot}{\kern0pt}{\isadigit{0}}}{\isasymrbrace}{\isacharparenright}{\kern0pt}{\isacharparenright}{\kern0pt}\ {\isacharparenleft}{\kern0pt}\isaconst{cpPair}\ \isavar{x}\ \isavar{y}{\isacharparenright}{\kern0pt}{\isacharparenright}{\kern0pt}\ {\isacharparenleft}{\kern0pt}\isaconst{cpPair}\ \isavar{Pa}\ \isavar{Pb}{\isacharparenright}{\kern0pt}\ \isactrlbold {\isasymcdot}\ \isavar{q}{\isacharparenright}{\kern0pt}\ {\isacharparenleft}{\kern0pt}\isaconst{mkEqS}\ {\isacharparenleft}{\kern0pt}\isaconst{mkAppS}\ {\isacharparenleft}{\kern0pt}\isaconst{mkAppS}\ {\isacharparenleft}{\kern0pt}\isaconst{mkAppS}\ \isaconst{mkRS}\ {\isacharparenleft}{\kern0pt}\isaconst{cpFst}\ {\isacharparenleft}{\kern0pt}\isaconst{trmPredR}\ {\isacharparenleft}{\kern0pt}\isaconst{cpSnd}\ {\isacharparenleft}{\kern0pt}\isaconst{trmPredR}\ {\isacharparenleft}{\kern0pt}\isaconst{cpSnd}\ {\isacharparenleft}{\kern0pt}\isaconst{cpSnd}\ {\isacharparenleft}{\kern0pt}\isaconst{cpSnd}\ {\isacharparenleft}{\kern0pt}\isaconst{cpFst}\ {\isacharparenleft}{\kern0pt}\isaconst{trmPredR}\ {\isacharparenleft}{\kern0pt}\isaconst{cpSnd}\ {\isacharparenleft}{\kern0pt}\isaconst{cpSnd}\ {\isacharparenleft}{\kern0pt}\isaconst{cpSnd}\ {\isacharparenleft}{\kern0pt}\isaconst{cpFst}\ {\isacharparenleft}{\kern0pt}\isaconst{trmPredR}\ {\isacharparenleft}{\kern0pt}\isaconst{cpSnd}\ {\isacharparenleft}{\kern0pt}\isaconst{cpSnd}\ {\isacharparenleft}{\kern0pt}\isaconst{cpSnd}\ {\isacharparenleft}{\kern0pt}\isaconst{trmPredR}\ {\isacharparenleft}{\kern0pt}\isaconst{cpFst}\ {\isacharparenleft}{\kern0pt}\isaconst{trmPredR}\ {\isacharparenleft}{\kern0pt}\isaconst{cmp}\ \isaconst{headRedXIdx}\ \isactrlbold {\isasymcdot}\ \isactrlbold {\isadigit{0}}\ \isactrlbold {\isasymcdot}\ \isaconst{cpPair}\ \isactrlbold {\isadigit{0}}\ {\isacharparenleft}{\kern0pt}\isaconst{trmPredR}\ {\isacharparenleft}{\kern0pt}\isaconst{cpFst}\ {\isacharparenleft}{\kern0pt}\isaconst{trmPredR}\ {\isacharparenleft}{\kern0pt}\isaconst{cmp}\ \isaconst{headRedXIdx}\ \isactrlbold {\isasymcdot}\ \isactrlbold {\isadigit{0}}\ \isactrlbold {\isasymcdot}\ \isaconst{cpPair}\ \isactrlbold {\isadigit{0}}\ {\isacharparenleft}{\kern0pt}\isaconst{iSubM}\ \isavar{z}\ {\isacharparenleft}{\kern0pt}\isactrlbold n{\isasymlbrace}\isaconst{cmp}\ {\isacharparenleft}{\kern0pt}\isaconst{prRec}\ \isavar{g}\ \isavar{h}\ \isavar{h{\isacharprime}{\kern0pt}}{\isacharparenright}{\kern0pt}\ \isactrlbold {\isasymcdot}\ \isavar{T{\isadigit{1}}{\isachardot}{\kern0pt}{\isadigit{0}}}\ \isactrlbold {\isasymcdot}\ \isavar{T{\isadigit{2}}{\isachardot}{\kern0pt}{\isadigit{0}}}{\isasymrbrace}{\isacharparenright}{\kern0pt}{\isacharparenright}{\kern0pt}{\isacharparenright}{\kern0pt}{\isacharparenright}{\kern0pt}{\isacharparenright}{\kern0pt}{\isacharparenright}{\kern0pt}{\isacharparenright}{\kern0pt}{\isacharparenright}{\kern0pt}{\isacharparenright}{\kern0pt}{\isacharparenright}{\kern0pt}{\isacharparenright}{\kern0pt}{\isacharparenright}{\kern0pt}{\isacharparenright}{\kern0pt}{\isacharparenright}{\kern0pt}{\isacharparenright}{\kern0pt}{\isacharparenright}{\kern0pt}{\isacharparenright}{\kern0pt}{\isacharparenright}{\kern0pt}{\isacharparenright}{\kern0pt}{\isacharparenright}{\kern0pt}{\isacharparenright}{\kern0pt}{\isacharparenright}{\kern0pt}{\isacharparenright}{\kern0pt}{\isacharparenright}{\kern0pt}{\isacharparenright}{\kern0pt}{\isacharparenright}{\kern0pt}{\isacharparenright}{\kern0pt}{\isacharparenright}{\kern0pt}\ {\isacharparenleft}{\kern0pt}\isaconst{cpFst}\ {\isacharparenleft}{\kern0pt}\isaconst{trmPredR}\ {\isacharparenleft}{\kern0pt}\isaconst{cpSnd}\ {\isacharparenleft}{\kern0pt}\isaconst{trmPredR}\ {\isacharparenleft}{\kern0pt}\isaconst{cpSnd}\ {\isacharparenleft}{\kern0pt}\isaconst{cpSnd}\ {\isacharparenleft}{\kern0pt}\isaconst{cpSnd}\ {\isacharparenleft}{\kern0pt}\isaconst{cpFst}\ {\isacharparenleft}{\kern0pt}\isaconst{trmPredR}\ {\isacharparenleft}{\kern0pt}\isaconst{cpSnd}\ {\isacharparenleft}{\kern0pt}\isaconst{cpSnd}\ {\isacharparenleft}{\kern0pt}\isaconst{cpSnd}\ {\isacharparenleft}{\kern0pt}\isaconst{trmPredR}\ {\isacharparenleft}{\kern0pt}\isaconst{cpFst}\ {\isacharparenleft}{\kern0pt}\isaconst{trmPredR}\ {\isacharparenleft}{\kern0pt}\isaconst{cmp}\ \isaconst{headRedXIdx}\ \isactrlbold {\isasymcdot}\ \isactrlbold {\isadigit{0}}\ \isactrlbold {\isasymcdot}\ \isaconst{cpPair}\ \isactrlbold {\isadigit{0}}\ {\isacharparenleft}{\kern0pt}\isaconst{trmPredR}\ {\isacharparenleft}{\kern0pt}\isaconst{cpFst}\ {\isacharparenleft}{\kern0pt}\isaconst{trmPredR}\ {\isacharparenleft}{\kern0pt}\isaconst{cmp}\ \isaconst{headRedXIdx}\ \isactrlbold {\isasymcdot}\ \isactrlbold {\isadigit{0}}\ \isactrlbold {\isasymcdot}\ \isaconst{cpPair}\ \isactrlbold {\isadigit{0}}\ {\isacharparenleft}{\kern0pt}\isaconst{iSubM}\ \isavar{z}\ {\isacharparenleft}{\kern0pt}\isactrlbold n{\isasymlbrace}\isaconst{cmp}\ {\isacharparenleft}{\kern0pt}\isaconst{prRec}\ \isavar{g}\ \isavar{h}\ \isavar{h{\isacharprime}{\kern0pt}}{\isacharparenright}{\kern0pt}\ \isactrlbold {\isasymcdot}\ \isavar{T{\isadigit{1}}{\isachardot}{\kern0pt}{\isadigit{0}}}\ \isactrlbold {\isasymcdot}\ \isavar{T{\isadigit{2}}{\isachardot}{\kern0pt}{\isadigit{0}}}{\isasymrbrace}{\isacharparenright}{\kern0pt}{\isacharparenright}{\kern0pt}{\isacharparenright}{\kern0pt}{\isacharparenright}{\kern0pt}{\isacharparenright}{\kern0pt}{\isacharparenright}{\kern0pt}{\isacharparenright}{\kern0pt}{\isacharparenright}{\kern0pt}{\isacharparenright}{\kern0pt}{\isacharparenright}{\kern0pt}{\isacharparenright}{\kern0pt}{\isacharparenright}{\kern0pt}{\isacharparenright}{\kern0pt}{\isacharparenright}{\kern0pt}{\isacharparenright}{\kern0pt}{\isacharparenright}{\kern0pt}{\isacharparenright}{\kern0pt}{\isacharparenright}{\kern0pt}{\isacharparenright}{\kern0pt}{\isacharparenright}{\kern0pt}{\isacharparenright}{\kern0pt}{\isacharparenright}{\kern0pt}{\isacharparenright}{\kern0pt}\ {\isacharparenleft}{\kern0pt}\isaconst{cmp}\ \isaconst{num{\isacharunderscore}{\kern0pt}codeIdx}\ \isactrlbold {\isasymcdot}\ \isactrlbold {\isadigit{0}}\ \isactrlbold {\isasymcdot}\ \isavar{q}{\isacharparenright}{\kern0pt}{\isacharparenright}{\kern0pt}\ {\isacharparenleft}{\kern0pt}\isaconst{cmp}\ \isaconst{num{\isacharunderscore}{\kern0pt}codeIdx}\ \isactrlbold {\isasymcdot}\ \isactrlbold {\isadigit{0}}\ \isactrlbold {\isasymcdot}\ {\isacharparenleft}{\kern0pt}\isaconst{cmp}\ {\isacharparenleft}{\kern0pt}\isaconst{prRec}\ \isavar{g}\ \isavar{h}\ \isavar{h{\isacharprime}{\kern0pt}}{\isacharparenright}{\kern0pt}\ \isactrlbold {\isasymcdot}\ \isavar{x}\ \isactrlbold {\isasymcdot}\ \isavar{q}{\isacharparenright}{\kern0pt}{\isacharparenright}{\kern0pt}{\isacharparenright}{\kern0pt}}%
}%
\SNIP{Goedel-Paper-Statements:thm:iChk-premCert-0}{%
\isakeywordONE{thm}\isamarkupfalse%
\ iChk{\isacharunderscore}{\kern0pt}premCert%
\isanewline%
}%
\SNIP{Goedel-Paper-Statements:thm:iChk-premCert-1}{%
\isa{{\isasymlbrakk}\isavar{{\isasymGamma}}\ {\isasymturnstile}\ \isavar{v}\ \isactrlbold N{\isacharsemicolon}{\kern0pt}\ \isavar{{\isasymGamma}}\ {\isasymturnstile}\ \isaconst{iChk}\ \isavar{v}\ \isaconst{goedelS}{\isasymrbrakk}\ {\isasymLongrightarrow}\ \isavar{{\isasymGamma}}\ {\isasymturnstile}\ \isaconst{iChk}\ {\isacharparenleft}{\kern0pt}\isaconst{premCert}\ \isavar{v}{\isacharparenright}{\kern0pt}\ {\isacharparenleft}{\kern0pt}\isaconst{negInstC}\ {\isacharparenleft}{\kern0pt}\isaconst{cmp}\ \isaconst{num{\isacharunderscore}{\kern0pt}codeIdx}\ \isactrlbold {\isasymcdot}\ \isactrlbold {\isadigit{0}}\ \isactrlbold {\isasymcdot}\ \isavar{v}{\isacharparenright}{\kern0pt}{\isacharparenright}{\kern0pt}}%
}%
\SNIP{Goedel-Paper-Statements:thm:premCertIdx-eval-0}{%
\isakeywordONE{thm}\isamarkupfalse%
\ premCertIdx{\isacharunderscore}{\kern0pt}eval%
\isanewline%
}%
\SNIP{Goedel-Paper-Statements:thm:premCertIdx-eval-1}{%
\isa{\isavar{{\isasymGamma}}\ {\isasymturnstile}\ \isavar{v}\ \isactrlbold N\ {\isasymLongrightarrow}\ \isavar{{\isasymGamma}}\ {\isasymturnstile}\ \isaconst{cmp}\ \isaconst{premCertIdx}\ \isactrlbold {\isasymcdot}\ \isavar{v}\ \isactrlbold {\isasymcdot}\ \isavar{v}\ \isactrlbold {\isacharequal}{\kern0pt}\ \isaconst{premCert}\ \isavar{v}}%
}%
\SNIP{Goedel-Paper-Statements:thm:D3Gcert-holds-0}{%
\isakeywordONE{thm}\isamarkupfalse%
\ D{\isadigit{3}}Gcert{\isacharunderscore}{\kern0pt}holds%
\isanewline%
}%
\SNIP{Goedel-Paper-Statements:thm:D3Gcert-holds-1}{%
\isa{\isaconst{D{\isadigit{3}}Gcert}\ \isaconst{premCertIdx}}%
}%
\SNIP{Goedel-Paper-Statements:thm:d3SigmaAt-goedelS-0}{%
\isakeywordONE{thm}\isamarkupfalse%
\ d{\isadigit{3}}SigmaAt{\isacharunderscore}{\kern0pt}goedelS%
\isanewline%
}%
\SNIP{Goedel-Paper-Statements:thm:d3SigmaAt-goedelS-1}{%
\isa{\isaconst{d{\isadigit{3}}SigmaAt}\ \isaconst{goedelS}\ {\isacharparenleft}{\kern0pt}\isaconst{prSub}\ \isaconst{diagSubIdx}\ \isaconst{premDIdx}\ \isaconst{premDIdx}{\isacharparenright}{\kern0pt}\ \isaconst{premCertSIdx}\ {\isacharequal}{\kern0pt}\ \isaconst{goedelSigmaOf}\ \isaconst{premCertIdx}}%
}%
\SNIP{Goedel-Paper-Statements:thm:D3GF-goedelSigma-inst-0}{%
\isakeywordONE{thm}\isamarkupfalse%
\ D{\isadigit{3}}GF{\isacharunderscore}{\kern0pt}goedelSigma{\isacharunderscore}{\kern0pt}inst%
\isanewline%
}%
\SNIP{Goedel-Paper-Statements:thm:D3GF-goedelSigma-inst-1}{%
\isa{{\isasymemptyset}\ {\isasymturnstile}\ \isaconst{D{\isadigit{3}}GF}\ {\isacharparenleft}{\kern0pt}\isaconst{goedelSigmaOf}\ \isaconst{premCertIdx}{\isacharparenright}{\kern0pt}}%
}%
\SNIP{Goedel-Paper-Statements:thm:prCU-num-0}{%
\isakeywordONE{thm}\isamarkupfalse%
\ prCU{\isacharunderscore}{\kern0pt}num%
\isanewline%
}%
\SNIP{Goedel-Paper-Statements:thm:prCU-num-1}{%
\isa{{\isasymemptyset}\ {\isasymturnstile}\ \isaconst{prCU}\ {\isacharparenleft}{\kern0pt}\isactrlbold n\isavar{k}{\isacharparenright}{\kern0pt}\ \isactrlbold {\isacharequal}{\kern0pt}\ \isactrlbold n{\isasymlbrace}\isactrlbold {\isasymnot}\isactrlbold A\ \isactrlbold {\isasymcdot}\ \isaconst{prvF}\ {\isacharparenleft}{\kern0pt}\isactrlbold n\isavar{k}{\isacharparenright}{\kern0pt}{\isasymrbrace}}%
}%
\SNIP{Goedel-Paper-Statements:thm:goedelG-of-expl-d3g-con-0}{%
\isakeywordONE{thm}\isamarkupfalse%
\ goedelG{\isacharunderscore}{\kern0pt}of{\isacharunderscore}{\kern0pt}expl{\isacharunderscore}{\kern0pt}d{\isadigit{3}}g{\isacharunderscore}{\kern0pt}con%
\isanewline%
}%
\SNIP{Goedel-Paper-Statements:thm:goedelG-of-expl-d3g-con-1}{%
\isa{{\isasymlbrakk}\isaconst{EXPL}\ \isavar{e}{\isacharsemicolon}{\kern0pt}\ \isaconst{D{\isadigit{3}}G}\ \isavar{{\isasymsigma}}{\isacharsemicolon}{\kern0pt}\ {\isasymemptyset}\ {\isasymturnstile}\ \isaconst{ConRGA}{\isasymrbrakk}\ {\isasymLongrightarrow}\ {\isasymemptyset}\ {\isasymturnstile}\ \isaconst{goedelG}}%
}%
\SNIP{Goedel-Paper-Statements:thm:rga-no-proof-0eq1-0}{%
\isakeywordONE{thm}\isamarkupfalse%
\ rga{\isacharunderscore}{\kern0pt}no{\isacharunderscore}{\kern0pt}proof{\isacharunderscore}{\kern0pt}{\isadigit{0}}eq{\isadigit{1}}%
\isanewline%
}%
\SNIP{Goedel-Paper-Statements:thm:rga-no-proof-0eq1-1}{%
\isa{{\isasymnot}\ {\isasymemptyset}\ {\isasymturnstile}\ \isactrlbold {\isadigit{0}}\ \isactrlbold {\isacharequal}{\kern0pt}\ \isactrlbold n{\isadigit{1}}}%
}%
\SNIP{Goedel-Paper-Statements:thm:rga-goedel2-of-expl-d3g-0}{%
\isakeywordONE{thm}\isamarkupfalse%
\ rga{\isacharunderscore}{\kern0pt}goedel{\isadigit{2}}{\isacharunderscore}{\kern0pt}of{\isacharunderscore}{\kern0pt}expl{\isacharunderscore}{\kern0pt}d{\isadigit{3}}g%
\isanewline%
}%
\SNIP{Goedel-Paper-Statements:thm:rga-goedel2-of-expl-d3g-1}{%
\isa{{\isasymlbrakk}\isaconst{EXPL}\ \isavar{e}{\isacharsemicolon}{\kern0pt}\ \isaconst{D{\isadigit{3}}G}\ \isavar{{\isasymsigma}}{\isasymrbrakk}\ {\isasymLongrightarrow}\ {\isasymnot}\ {\isasymemptyset}\ {\isasymturnstile}\ \isaconst{ConRGA}}%
}%
\SNIP{Goedel-Paper-Statements:thm:rga-goedel2-of-d3g-0}{%
\isakeywordONE{thm}\isamarkupfalse%
\ rga{\isacharunderscore}{\kern0pt}goedel{\isadigit{2}}{\isacharunderscore}{\kern0pt}of{\isacharunderscore}{\kern0pt}d{\isadigit{3}}g%
\isanewline%
}%
\SNIP{Goedel-Paper-Statements:thm:rga-goedel2-of-d3g-1}{%
\isa{\isaconst{D{\isadigit{3}}G}\ \isavar{{\isasymsigma}}\ {\isasymLongrightarrow}\ {\isasymnot}\ {\isasymemptyset}\ {\isasymturnstile}\ \isaconst{ConRGA}}%
}%
\SNIP{Goedel-Paper-Statements:thm:goedelG-no-accepted-code-0}{%
\isakeywordONE{thm}\isamarkupfalse%
\ goedelG{\isacharunderscore}{\kern0pt}no{\isacharunderscore}{\kern0pt}accepted{\isacharunderscore}{\kern0pt}code%
\isanewline%
}%
\SNIP{Goedel-Paper-Statements:thm:goedelG-no-accepted-code-1}{%
\isa{{\isasymnot}\ {\isasymemptyset}\ {\isasymturnstile}\ \isaconst{iChk}\ {\isacharparenleft}{\kern0pt}\isactrlbold n\isavar{p}{\isacharparenright}{\kern0pt}\ {\isacharparenleft}{\kern0pt}\isactrlbold n{\isasymlbrace}\isaconst{goedelG}{\isasymrbrace}{\isacharparenright}{\kern0pt}}%
}%
\SNIP{Goedel-Paper-Statements:thm:rga-goedel1-unrefutable-0}{%
\isakeywordONE{thm}\isamarkupfalse%
\ rga{\isacharunderscore}{\kern0pt}goedel{\isadigit{1}}{\isacharunderscore}{\kern0pt}unrefutable%
\isanewline%
}%
\SNIP{Goedel-Paper-Statements:thm:rga-goedel1-unrefutable-1}{%
\isa{{\isasymnot}\ {\isasymemptyset}\ {\isasymturnstile}\ \isactrlbold {\isasymnot}\isaconst{goedelG}}%
}%
\SNIP{Goedel-Paper-Statements:thm:sigma-intro-0}{%
\isakeywordONE{thm}\isamarkupfalse%
\ sigma{\isacharunderscore}{\kern0pt}intro%
\isanewline%
}%
\SNIP{Goedel-Paper-Statements:thm:sigma-intro-1}{%
\isa{{\isasymlbrakk}\isavar{{\isasymGamma}}\ {\isasymturnstile}\ \isaconst{bsub}\ \isavar{a}\ \isavar{D}{\isacharsemicolon}{\kern0pt}\ \isavar{{\isasymGamma}}\ {\isasymturnstile}\ \isavar{a}\ \isactrlbold N{\isasymrbrakk}\ {\isasymLongrightarrow}\ \isavar{{\isasymGamma}}\ {\isasymturnstile}\ \isactrlbold {\isasymnot}\isactrlbold A\ \isactrlbold {\isasymcdot}\ {\isacharparenleft}{\kern0pt}\isactrlbold {\isasymlambda}{\isacharparenleft}{\kern0pt}\isactrlbold {\isasymnot}\isavar{D}{\isacharparenright}{\kern0pt}{\isacharparenright}{\kern0pt}}%
}%
\SNIP{Goedel-Paper-Statements:thm:sigma-elim-0}{%
\isakeywordONE{thm}\isamarkupfalse%
\ sigma{\isacharunderscore}{\kern0pt}elim%
\isanewline%
}%
\SNIP{Goedel-Paper-Statements:thm:sigma-elim-1}{%
\isa{{\isasymlbrakk}\isavar{{\isasymGamma}}\ {\isasymturnstile}\ \isactrlbold {\isasymnot}\isactrlbold A\ \isactrlbold {\isasymcdot}\ {\isacharparenleft}{\kern0pt}\isactrlbold {\isasymlambda}{\isacharparenleft}{\kern0pt}\isactrlbold {\isasymnot}\isavar{D}{\isacharparenright}{\kern0pt}{\isacharparenright}{\kern0pt}{\isacharsemicolon}{\kern0pt}\ \isavar{D}\ {\isasymtriangleright}\ \isactrlbold v\isactrlsub {\isadigit{0}}\ \isactrlbold N\ {\isasymtriangleright}\ {\isasymupharpoonleft}\isavar{{\isasymGamma}}\ {\isasymturnstile}\ {\isasymupharpoonleft}\isavar{C}{\isasymrbrakk}\ {\isasymLongrightarrow}\ \isavar{{\isasymGamma}}\ {\isasymturnstile}\ \isavar{C}}%
}%
\SNIP{Goedel-Paper-Statements:thm:sigma-reductio-0}{%
\isakeywordONE{thm}\isamarkupfalse%
\ sigma{\isacharunderscore}{\kern0pt}reductio%
\isanewline%
}%
\SNIP{Goedel-Paper-Statements:thm:sigma-reductio-1}{%
\isa{{\isasymlbrakk}\isactrlbold v\isactrlsub {\isadigit{0}}\ \isactrlbold N\ {\isasymtriangleright}\ {\isasymupharpoonleft}\isavar{{\isasymGamma}}\ {\isasymturnstile}\ \isavar{D}\ \isactrlbold B{\isacharsemicolon}{\kern0pt}\ \isactrlbold {\isasymnot}\isactrlbold A\ \isactrlbold {\isasymcdot}\ {\isacharparenleft}{\kern0pt}\isactrlbold {\isasymlambda}{\isacharparenleft}{\kern0pt}\isactrlbold {\isasymnot}\isavar{D}{\isacharparenright}{\kern0pt}{\isacharparenright}{\kern0pt}\ {\isasymtriangleright}\ \isavar{{\isasymGamma}}\ {\isasymturnstile}\ \isactrlbold {\isasymnot}\isactrlbold {\isasymnot}\isactrlbold A\ \isactrlbold {\isasymcdot}\ {\isacharparenleft}{\kern0pt}\isactrlbold {\isasymlambda}{\isacharparenleft}{\kern0pt}\isactrlbold {\isasymnot}\isavar{D}{\isacharparenright}{\kern0pt}{\isacharparenright}{\kern0pt}{\isasymrbrakk}\ {\isasymLongrightarrow}\ \isavar{{\isasymGamma}}\ {\isasymturnstile}\ \isactrlbold {\isasymnot}\isactrlbold {\isasymnot}\isactrlbold A\ \isactrlbold {\isasymcdot}\ {\isacharparenleft}{\kern0pt}\isactrlbold {\isasymlambda}{\isacharparenleft}{\kern0pt}\isactrlbold {\isasymnot}\isavar{D}{\isacharparenright}{\kern0pt}{\isacharparenright}{\kern0pt}}%
}%
\SNIP{Goedel-Paper-Statements:thm:chk-reflect-num-0}{%
\isakeywordONE{thm}\isamarkupfalse%
\ chk{\isacharunderscore}{\kern0pt}reflect{\isacharunderscore}{\kern0pt}num%
\isanewline%
}%
\SNIP{Goedel-Paper-Statements:thm:chk-reflect-num-1}{%
\isa{\isaconst{wft}\ \isavar{A}\ {\isasymLongrightarrow}\ {\isasymemptyset}\ {\isasymturnstile}\ \isaconst{iChk}\ {\isacharparenleft}{\kern0pt}\isactrlbold n\isavar{n}{\isacharparenright}{\kern0pt}\ {\isacharparenleft}{\kern0pt}\isactrlbold n{\isasymlbrace}\isavar{A}{\isasymrbrace}{\isacharparenright}{\kern0pt}\ \isactrlbold {\isasymlongrightarrow}\ \isavar{A}}%
}%
\SNIP{Goedel-Paper-Statements:thm:chk-reflect-num-ent-0}{%
\isakeywordONE{thm}\isamarkupfalse%
\ chk{\isacharunderscore}{\kern0pt}reflect{\isacharunderscore}{\kern0pt}num{\isacharunderscore}{\kern0pt}ent%
\isanewline%
}%
\SNIP{Goedel-Paper-Statements:thm:chk-reflect-num-ent-1}{%
\isa{\isaconst{wft}\ \isavar{A}\ {\isasymLongrightarrow}\ \isaconst{iChk}\ {\isacharparenleft}{\kern0pt}\isactrlbold n\isavar{n}{\isacharparenright}{\kern0pt}\ {\isacharparenleft}{\kern0pt}\isactrlbold n{\isasymlbrace}\isavar{A}{\isasymrbrace}{\isacharparenright}{\kern0pt}\ {\isasymtriangleright}\ \isavar{{\isasymGamma}}\ {\isasymturnstile}\ \isavar{A}}%
}%
\SNIP{Goedel-Paper-Statements:thm:exfIdx-chk-0}{%
\isakeywordONE{thm}\isamarkupfalse%
\ exfIdx{\isacharunderscore}{\kern0pt}chk%
\isanewline%
}%
\SNIP{Goedel-Paper-Statements:thm:exfIdx-chk-1}{%
\isa{{\isasymlbrakk}\isaconst{wft}\ \isavar{{\isasymphi}}{\isacharsemicolon}{\kern0pt}\ \isavar{{\isasymGamma}}\ {\isasymturnstile}\ \isavar{a}\ \isactrlbold N{\isacharsemicolon}{\kern0pt}\ \isavar{{\isasymGamma}}\ {\isasymturnstile}\ \isaconst{iChk}\ \isavar{a}\ \isaconst{cEq{\isadigit{0}}{\isadigit{1}}}{\isasymrbrakk}\ {\isasymLongrightarrow}\ \isavar{{\isasymGamma}}\ {\isasymturnstile}\ \isaconst{iChk}\ {\isacharparenleft}{\kern0pt}\isaconst{cmp}\ {\isacharparenleft}{\kern0pt}\isaconst{exfIdx}\ \isavar{{\isasymphi}}{\isacharparenright}{\kern0pt}\ \isactrlbold {\isasymcdot}\ \isavar{a}\ \isactrlbold {\isasymcdot}\ \isavar{a}{\isacharparenright}{\kern0pt}\ {\isacharparenleft}{\kern0pt}\isactrlbold n{\isasymlbrace}\isavar{{\isasymphi}}{\isasymrbrace}{\isacharparenright}{\kern0pt}}%
}%
\SNIP{Goedel-Paper-Statements:thm:G-exfIdx-chk-sf-0}{%
\isakeywordONE{thm}\isamarkupfalse%
\ G{\isacharunderscore}{\kern0pt}exfIdx{\isacharunderscore}{\kern0pt}chk{\isacharunderscore}{\kern0pt}sf%
\isanewline%
}%
\SNIP{Goedel-Paper-Statements:thm:G-exfIdx-chk-sf-1}{%
\isa{{\isasymlbrakk}\isaconst{tblShift}\ {\isacharparenleft}{\kern0pt}\isaconst{drvChkIdxG}\ \isavar{L}{\isacharparenright}{\kern0pt}{\isacharsemicolon}{\kern0pt}\ \isaconst{chkNegEIdx}\ {\isasymin}\ \isaconst{set}\ \isavar{L}{\isacharsemicolon}{\kern0pt}\ \isaconst{set}\ \isaconst{drvChkList}\ {\isasymsubseteq}\ \isaconst{set}\ \isavar{L}{\isacharsemicolon}{\kern0pt}\ \isaconst{wft}\ \isavar{{\isasymphi}}{\isacharsemicolon}{\kern0pt}\ \isavar{{\isasymGamma}}\ {\isasymturnstile}\ \isavar{a}\ \isactrlbold N{\isacharsemicolon}{\kern0pt}\ \isavar{{\isasymGamma}}\ {\isasymturnstile}\ \isaconst{iChkD}\ {\isacharparenleft}{\kern0pt}\isaconst{drvChkIdxG}\ \isavar{L}{\isacharparenright}{\kern0pt}\ \isavar{a}\ \isaconst{cEq{\isadigit{0}}{\isadigit{1}}}{\isasymrbrakk}\ {\isasymLongrightarrow}\ \isavar{{\isasymGamma}}\ {\isasymturnstile}\ \isaconst{iChkD}\ {\isacharparenleft}{\kern0pt}\isaconst{drvChkIdxG}\ \isavar{L}{\isacharparenright}{\kern0pt}\ {\isacharparenleft}{\kern0pt}\isaconst{cmp}\ {\isacharparenleft}{\kern0pt}\isaconst{exfIdx}\ \isavar{{\isasymphi}}{\isacharparenright}{\kern0pt}\ \isactrlbold {\isasymcdot}\ \isavar{a}\ \isactrlbold {\isasymcdot}\ \isavar{a}{\isacharparenright}{\kern0pt}\ {\isacharparenleft}{\kern0pt}\isactrlbold n{\isasymlbrace}\isavar{{\isasymphi}}{\isasymrbrace}{\isacharparenright}{\kern0pt}}%
}%
\SNIP{Goedel-Paper-Statements:thm:G-d2M-chk-0}{%
\isakeywordONE{thm}\isamarkupfalse%
\ G{\isacharunderscore}{\kern0pt}d{\isadigit{2}}M{\isacharunderscore}{\kern0pt}chk%
\isanewline%
}%
\SNIP{Goedel-Paper-Statements:thm:G-d2M-chk-1}{%
\isa{{\isasymlbrakk}\isaconst{tblShift}\ {\isacharparenleft}{\kern0pt}\isaconst{drvChkIdxG}\ \isavar{L}{\isacharparenright}{\kern0pt}{\isacharsemicolon}{\kern0pt}\ \isaconst{chkCutIdx}\ {\isasymin}\ \isaconst{set}\ \isavar{L}{\isacharsemicolon}{\kern0pt}\ \isaconst{chkWk{\isadigit{1}}Idx}\ {\isasymin}\ \isaconst{set}\ \isavar{L}{\isacharsemicolon}{\kern0pt}\ \isaconst{set}\ \isaconst{drvChkList}\ {\isasymsubseteq}\ \isaconst{set}\ \isavar{L}{\isacharsemicolon}{\kern0pt}\ \isaconst{wft}\ \isavar{{\isasymphi}}{\isacharsemicolon}{\kern0pt}\ \isaconst{wft}\ \isavar{{\isasympsi}}{\isacharsemicolon}{\kern0pt}\ \isavar{{\isasymGamma}}\ {\isasymturnstile}\ \isavar{a}\ \isactrlbold N{\isacharsemicolon}{\kern0pt}\ \isavar{{\isasymGamma}}\ {\isasymturnstile}\ \isavar{b}\ \isactrlbold N{\isacharsemicolon}{\kern0pt}\ \isavar{{\isasymGamma}}\ {\isasymturnstile}\ \isaconst{iChkD}\ {\isacharparenleft}{\kern0pt}\isaconst{drvChkIdxG}\ \isavar{L}{\isacharparenright}{\kern0pt}\ \isavar{a}\ {\isacharparenleft}{\kern0pt}\isactrlbold n{\isasymlbrace}\isavar{{\isasymphi}}\ \isactrlbold {\isasymlongrightarrow}\ \isavar{{\isasympsi}}{\isasymrbrace}{\isacharparenright}{\kern0pt}{\isacharsemicolon}{\kern0pt}\ \isavar{{\isasymGamma}}\ {\isasymturnstile}\ \isaconst{iChkD}\ {\isacharparenleft}{\kern0pt}\isaconst{drvChkIdxG}\ \isavar{L}{\isacharparenright}{\kern0pt}\ \isavar{b}\ {\isacharparenleft}{\kern0pt}\isactrlbold n{\isasymlbrace}\isavar{{\isasymphi}}{\isasymrbrace}{\isacharparenright}{\kern0pt}{\isasymrbrakk}\ {\isasymLongrightarrow}\ \isavar{{\isasymGamma}}\ {\isasymturnstile}\ \isaconst{iChkD}\ {\isacharparenleft}{\kern0pt}\isaconst{drvChkIdxG}\ \isavar{L}{\isacharparenright}{\kern0pt}\ {\isacharparenleft}{\kern0pt}\isaconst{d{\isadigit{2}}M}\ \isavar{{\isasymphi}}\ \isavar{{\isasympsi}}\ \isactrlbold {\isasymcdot}\ \isavar{a}\ \isactrlbold {\isasymcdot}\ \isavar{b}{\isacharparenright}{\kern0pt}\ {\isacharparenleft}{\kern0pt}\isactrlbold n{\isasymlbrace}\isavar{{\isasympsi}}{\isasymrbrace}{\isacharparenright}{\kern0pt}}%
}%
\SNIP{Goedel-Paper-Statements:thm:G-cut2M-chk-0}{%
\isakeywordONE{thm}\isamarkupfalse%
\ G{\isacharunderscore}{\kern0pt}cut{\isadigit{2}}M{\isacharunderscore}{\kern0pt}chk%
\isanewline%
}%
\SNIP{Goedel-Paper-Statements:thm:G-cut2M-chk-1}{%
\isa{{\isasymlbrakk}{\isasymAnd}\isabound{P}\ \isabound{{\isasymGamma}}\ \isabound{c}{\isachardot}{\kern0pt}\ {\isacharparenleft}{\kern0pt}\isaconst{pr{\isadigit{2}}fi}\ {\isacharparenleft}{\kern0pt}\isaconst{drvChkIdxG}\ \isavar{L}{\isacharparenright}{\kern0pt}\ {\isasymlangle}{\isasymlbrace}\isabound{P}{\isasymrbrace}{\isacharcomma}{\kern0pt}{\isasymlbrace}\isabound{{\isasymGamma}}{\isasymrbrace}{\isasymrangle}\ {\isasymlbrace}\isabound{c}{\isasymrbrace}\ {\isasymnoteq}\ {\isadigit{0}}{\isacharparenright}{\kern0pt}\ {\isacharequal}{\kern0pt}\ \isavar{D}\ \isabound{P}\ \isabound{{\isasymGamma}}\ \isabound{c}{\isacharsemicolon}{\kern0pt}\ \isaconst{tblShift}\ {\isacharparenleft}{\kern0pt}\isaconst{drvChkIdxG}\ \isavar{L}{\isacharparenright}{\kern0pt}{\isacharsemicolon}{\kern0pt}\ \isaconst{chkCutIdx}\ {\isasymin}\ \isaconst{set}\ \isavar{L}{\isacharsemicolon}{\kern0pt}\ \isaconst{chkWk{\isadigit{1}}Idx}\ {\isasymin}\ \isaconst{set}\ \isavar{L}{\isacharsemicolon}{\kern0pt}\ \isaconst{prf{\isacharunderscore}{\kern0pt}derives{\isachardot}{\kern0pt}proves}\ \isavar{D}\ \isavar{Q}\ {\isacharparenleft}{\kern0pt}\isavar{A{\isadigit{2}}{\isachardot}{\kern0pt}{\isadigit{0}}}\ {\isasymtriangleright}\ \isavar{A{\isadigit{1}}{\isachardot}{\kern0pt}{\isadigit{0}}}\ {\isasymtriangleright}\ {\isasymemptyset}{\isacharparenright}{\kern0pt}\ \isavar{B}{\isacharsemicolon}{\kern0pt}\ \isavar{{\isasymGamma}}\ {\isasymturnstile}\ \isavar{a}\ \isactrlbold N{\isacharsemicolon}{\kern0pt}\ \isavar{{\isasymGamma}}\ {\isasymturnstile}\ \isavar{b}\ \isactrlbold N{\isacharsemicolon}{\kern0pt}\ \isavar{{\isasymGamma}}\ {\isasymturnstile}\ \isaconst{iChkD}\ {\isacharparenleft}{\kern0pt}\isaconst{drvChkIdxG}\ \isavar{L}{\isacharparenright}{\kern0pt}\ \isavar{a}\ {\isacharparenleft}{\kern0pt}\isactrlbold n{\isasymlbrace}\isavar{A{\isadigit{1}}{\isachardot}{\kern0pt}{\isadigit{0}}}{\isasymrbrace}{\isacharparenright}{\kern0pt}{\isacharsemicolon}{\kern0pt}\ \isavar{{\isasymGamma}}\ {\isasymturnstile}\ \isaconst{iChkD}\ {\isacharparenleft}{\kern0pt}\isaconst{drvChkIdxG}\ \isavar{L}{\isacharparenright}{\kern0pt}\ \isavar{b}\ {\isacharparenleft}{\kern0pt}\isactrlbold n{\isasymlbrace}\isavar{A{\isadigit{2}}{\isachardot}{\kern0pt}{\isadigit{0}}}{\isasymrbrace}{\isacharparenright}{\kern0pt}{\isasymrbrakk}\ {\isasymLongrightarrow}\ \isavar{{\isasymGamma}}\ {\isasymturnstile}\ \isaconst{iChkD}\ {\isacharparenleft}{\kern0pt}\isaconst{drvChkIdxG}\ \isavar{L}{\isacharparenright}{\kern0pt}\ {\isacharparenleft}{\kern0pt}\isaconst{cut{\isadigit{2}}M}\ \isavar{Q}\ \isavar{A{\isadigit{1}}{\isachardot}{\kern0pt}{\isadigit{0}}}\ \isavar{A{\isadigit{2}}{\isachardot}{\kern0pt}{\isadigit{0}}}\ \isavar{B}\ \isactrlbold {\isasymcdot}\ \isavar{a}\ \isactrlbold {\isasymcdot}\ \isavar{b}{\isacharparenright}{\kern0pt}\ {\isacharparenleft}{\kern0pt}\isactrlbold n{\isasymlbrace}\isavar{B}{\isasymrbrace}{\isacharparenright}{\kern0pt}}%
}%
\SNIP{Goedel-Paper-Statements:thm:adrGood-def-0}{%
\isakeywordONE{thm}\isamarkupfalse%
\ adrGood{\isacharunderscore}{\kern0pt}def%
\isanewline%
}%
\SNIP{Goedel-Paper-Statements:thm:adrGood-def-1}{%
\isa{\isaconst{adrGood}\ \isavar{v}\ {\isasymequiv}\ \isaconst{lstCons}\ {\isacharparenleft}{\kern0pt}\isaconst{cpPair}\ {\isacharparenleft}{\kern0pt}\isactrlbold n{\isasymlbrace}{\isasymemptyset}{\isasymrbrace}{\isacharparenright}{\kern0pt}\ {\isacharparenleft}{\kern0pt}\isaconst{adrInst}\ \isaconst{adrU}\ \isavar{v}{\isacharparenright}{\kern0pt}{\isacharparenright}{\kern0pt}\ {\isacharparenleft}{\kern0pt}\isaconst{certStepBack}\ {\isacharparenleft}{\kern0pt}\isaconst{adrInst}\ \isaconst{adrU}\ \isavar{v}{\isacharparenright}{\kern0pt}\ {\isacharparenleft}{\kern0pt}\isaconst{trmPredR}\ {\isacharparenleft}{\kern0pt}\isaconst{cpFst}\ {\isacharparenleft}{\kern0pt}\isaconst{trmPredR}\ {\isacharparenleft}{\kern0pt}\isaconst{cmp}\ \isaconst{headRedXIdx}\ \isactrlbold {\isasymcdot}\ \isactrlbold {\isadigit{0}}\ \isactrlbold {\isasymcdot}\ \isaconst{cpPair}\ \isactrlbold {\isadigit{0}}\ {\isacharparenleft}{\kern0pt}\isaconst{adrInst}\ \isaconst{adrU}\ \isavar{v}{\isacharparenright}{\kern0pt}{\isacharparenright}{\kern0pt}{\isacharparenright}{\kern0pt}{\isacharparenright}{\kern0pt}{\isacharparenright}{\kern0pt}\ {\isacharparenleft}{\kern0pt}\isactrlbold n{\isasymlbrace}\isactrlbold n{\isadigit{1}}{\isasymrbrace}{\isacharparenright}{\kern0pt}\ {\isacharparenleft}{\kern0pt}\isaconst{lstCons}\ {\isacharparenleft}{\kern0pt}\isaconst{cpPair}\ {\isacharparenleft}{\kern0pt}\isactrlbold n{\isasymlbrace}{\isasymemptyset}{\isasymrbrace}{\isacharparenright}{\kern0pt}\ {\isacharparenleft}{\kern0pt}\isaconst{mkEqS}\ {\isacharparenleft}{\kern0pt}\isaconst{imgV{\isadigit{0}}}\ \isavar{v}\ \isaconst{adrBody}{\isacharparenright}{\kern0pt}\ {\isacharparenleft}{\kern0pt}\isactrlbold n{\isasymlbrace}\isactrlbold n{\isadigit{1}}{\isasymrbrace}{\isacharparenright}{\kern0pt}{\isacharparenright}{\kern0pt}{\isacharparenright}{\kern0pt}\ {\isacharparenleft}{\kern0pt}\isaconst{adrValCert}\ \isavar{v}{\isacharparenright}{\kern0pt}{\isacharparenright}{\kern0pt}{\isacharparenright}{\kern0pt}}%
}%
\SNIP{Goedel-Paper-Statements:thm:adr-exf-code-0}{%
\isakeywordONE{thm}\isamarkupfalse%
\ adr{\isacharunderscore}{\kern0pt}exf{\isacharunderscore}{\kern0pt}code%
\isanewline%
}%
\SNIP{Goedel-Paper-Statements:thm:adr-exf-code-1}{%
\isa{{\isasymlbrakk}\isavar{{\isasymGamma}}\ {\isasymturnstile}\ \isavar{a}\ \isactrlbold N{\isacharsemicolon}{\kern0pt}\ \isavar{{\isasymGamma}}\ {\isasymturnstile}\ \isaconst{iChk}\ \isavar{a}\ \isaconst{cEq{\isadigit{0}}{\isadigit{1}}}{\isacharsemicolon}{\kern0pt}\ \isavar{{\isasymGamma}}\ {\isasymturnstile}\ \isavar{y}\ \isactrlbold N{\isacharsemicolon}{\kern0pt}\ \isavar{{\isasymGamma}}\ {\isasymturnstile}\ \isaconst{iWft}\ \isavar{y}\ \isactrlbold {\isasymnoteq}\ \isactrlbold {\isadigit{0}}{\isacharsemicolon}{\kern0pt}\ \isavar{{\isasymGamma}}\ {\isasymturnstile}\ \isaconst{iIsCode}\ \isavar{y}\ \isactrlbold {\isasymnoteq}\ \isactrlbold {\isadigit{0}}{\isasymrbrakk}\ {\isasymLongrightarrow}\ \isavar{{\isasymGamma}}\ {\isasymturnstile}\ \isaconst{iChk}\ {\isacharparenleft}{\kern0pt}\isaconst{lstCons}\ {\isacharparenleft}{\kern0pt}\isaconst{cpPair}\ {\isacharparenleft}{\kern0pt}\isactrlbold n{\isasymlbrace}{\isasymemptyset}{\isasymrbrace}{\isacharparenright}{\kern0pt}\ \isavar{y}{\isacharparenright}{\kern0pt}\ {\isacharparenleft}{\kern0pt}\isaconst{cmp}\ \isaconst{apCoreIdx}\ \isactrlbold {\isasymcdot}\ \isaconst{cpPair}\ \isavar{a}\ {\isacharparenleft}{\kern0pt}\isactrlbold n{\isasymlbrace}\isaconst{neq{\isadigit{0}}{\isadigit{1}}Prf}{\isasymrbrace}{\isacharparenright}{\kern0pt}\ \isactrlbold {\isasymcdot}\ {\isacharparenleft}{\kern0pt}\isaconst{cmp}\ \isaconst{length{\isacharunderscore}{\kern0pt}encIdx}\ \isactrlbold {\isasymcdot}\ \isavar{a}\ \isactrlbold {\isasymcdot}\ \isavar{a}{\isacharparenright}{\kern0pt}{\isacharparenright}{\kern0pt}{\isacharparenright}{\kern0pt}\ \isavar{y}}%
}%
\SNIP{Goedel-Paper-Statements:thm:adi-exf-code-0}{%
\isakeywordONE{thm}\isamarkupfalse%
\ adi{\isacharunderscore}{\kern0pt}exf{\isacharunderscore}{\kern0pt}code%
\isanewline%
}%
\SNIP{Goedel-Paper-Statements:thm:adi-exf-code-1}{%
\isa{{\isasymlbrakk}\isavar{{\isasymGamma}}\ {\isasymturnstile}\ \isavar{a}\ \isactrlbold N{\isacharsemicolon}{\kern0pt}\ \isavar{{\isasymGamma}}\ {\isasymturnstile}\ \isaconst{iChkD}\ {\isacharparenleft}{\kern0pt}\isaconst{drvChkIdxG}\ {\isacharparenleft}{\kern0pt}\isaconst{drvChkList}\ {\isacharat}{\kern0pt}\ {\isacharbrackleft}{\kern0pt}\isaconst{chkGfixIdx}{\isacharbrackright}{\kern0pt}{\isacharparenright}{\kern0pt}{\isacharparenright}{\kern0pt}\ \isavar{a}\ \isaconst{cEq{\isadigit{0}}{\isadigit{1}}}{\isacharsemicolon}{\kern0pt}\ \isavar{{\isasymGamma}}\ {\isasymturnstile}\ \isavar{y}\ \isactrlbold N{\isacharsemicolon}{\kern0pt}\ \isavar{{\isasymGamma}}\ {\isasymturnstile}\ \isaconst{iWft}\ \isavar{y}\ \isactrlbold {\isasymnoteq}\ \isactrlbold {\isadigit{0}}{\isacharsemicolon}{\kern0pt}\ \isavar{{\isasymGamma}}\ {\isasymturnstile}\ \isaconst{iIsCode}\ \isavar{y}\ \isactrlbold {\isasymnoteq}\ \isactrlbold {\isadigit{0}}{\isasymrbrakk}\ {\isasymLongrightarrow}\ \isavar{{\isasymGamma}}\ {\isasymturnstile}\ \isaconst{iChkD}\ {\isacharparenleft}{\kern0pt}\isaconst{drvChkIdxG}\ {\isacharparenleft}{\kern0pt}\isaconst{drvChkList}\ {\isacharat}{\kern0pt}\ {\isacharbrackleft}{\kern0pt}\isaconst{chkGfixIdx}{\isacharbrackright}{\kern0pt}{\isacharparenright}{\kern0pt}{\isacharparenright}{\kern0pt}\ {\isacharparenleft}{\kern0pt}\isaconst{lstCons}\ {\isacharparenleft}{\kern0pt}\isaconst{cpPair}\ {\isacharparenleft}{\kern0pt}\isactrlbold n{\isasymlbrace}{\isasymemptyset}{\isasymrbrace}{\isacharparenright}{\kern0pt}\ \isavar{y}{\isacharparenright}{\kern0pt}\ {\isacharparenleft}{\kern0pt}\isaconst{cmp}\ \isaconst{apCoreIdx}\ \isactrlbold {\isasymcdot}\ \isaconst{cpPair}\ \isavar{a}\ {\isacharparenleft}{\kern0pt}\isactrlbold n{\isasymlbrace}\isaconst{neq{\isadigit{0}}{\isadigit{1}}Prf}{\isasymrbrace}{\isacharparenright}{\kern0pt}\ \isactrlbold {\isasymcdot}\ {\isacharparenleft}{\kern0pt}\isaconst{cmp}\ \isaconst{length{\isacharunderscore}{\kern0pt}encIdx}\ \isactrlbold {\isasymcdot}\ \isavar{a}\ \isactrlbold {\isasymcdot}\ \isavar{a}{\isacharparenright}{\kern0pt}{\isacharparenright}{\kern0pt}{\isacharparenright}{\kern0pt}\ \isavar{y}}%
}%
\SNIP{Goedel-Paper-Statements:thm:rga-goedel2-sentence-witness-0}{%
\isakeywordONE{thm}\isamarkupfalse%
\ rga{\isacharunderscore}{\kern0pt}goedel{\isadigit{2}}{\isacharunderscore}{\kern0pt}sentence{\isacharunderscore}{\kern0pt}witness%
\isanewline%
}%
\SNIP{Goedel-Paper-Statements:thm:rga-goedel2-sentence-witness-1}{%
\isa{{\isasymexists}\isabound{T}{\isachardot}{\kern0pt}\ \isaconst{Syntax{\isacharunderscore}{\kern0pt}Term{\isachardot}{\kern0pt}closed}\ \isabound{T}\ {\isasymand}\ {\isasymemptyset}\ {\isasymturnstile}\ \isactrlbold A\ \isactrlbold {\isasymcdot}\ {\isacharparenleft}{\kern0pt}\isactrlbold {\isasymlambda}{\isacharparenleft}{\kern0pt}\isaconst{iChk}\ \isactrlbold v\isactrlsub {\isadigit{0}}\ {\isacharparenleft}{\kern0pt}\isactrlbold n{\isasymlbrace}\isaconst{ConRGA}{\isasymrbrace}{\isacharparenright}{\kern0pt}\ \isactrlbold {\isasymlongrightarrow}\ \isaconst{iChk}\ {\isacharparenleft}{\kern0pt}\isabound{T}\ \isactrlbold {\isasymcdot}\ \isactrlbold v\isactrlsub {\isadigit{0}}{\isacharparenright}{\kern0pt}\ \isaconst{cEq{\isadigit{0}}{\isadigit{1}}}{\isacharparenright}{\kern0pt}{\isacharparenright}{\kern0pt}}%
}%
\SNIP{Goedel-Paper-Statements:thm:bs-goedel2-0}{%
\isakeywordONE{thm}\isamarkupfalse%
\ bs{\isacharunderscore}{\kern0pt}goedel{\isadigit{2}}%
\isanewline%
}%
\SNIP{Goedel-Paper-Statements:thm:bs-goedel2-1}{%
\isa{\isactrlbold {\isasymnot}\isaconst{goedelG}\ {\isasymtriangleright}\ {\isasymemptyset}\ {\isasymturnstile}\ \isactrlbold {\isasymnot}\isactrlbold A\ \isactrlbold {\isasymcdot}\ {\isacharparenleft}{\kern0pt}\isactrlbold {\isasymlambda}{\isacharparenleft}{\kern0pt}\isactrlbold {\isasymnot}\isaconst{iChk}\ \isactrlbold v\isactrlsub {\isadigit{0}}\ {\isacharparenleft}{\kern0pt}{\isasymupharpoonleft}{\isacharparenleft}{\kern0pt}\isactrlbold n{\isasymlbrace}\isactrlbold {\isasymnot}\isactrlbold {\isadigit{0}}\ \isactrlbold {\isacharequal}{\kern0pt}\ \isactrlbold {\isadigit{0}}{\isasymrbrace}{\isacharparenright}{\kern0pt}{\isacharparenright}{\kern0pt}{\isacharparenright}{\kern0pt}{\isacharparenright}{\kern0pt}\ {\isasymand}\ \isactrlbold {\isasymnot}\isactrlbold A\ \isactrlbold {\isasymcdot}\ {\isacharparenleft}{\kern0pt}\isactrlbold {\isasymlambda}{\isacharparenleft}{\kern0pt}\isactrlbold {\isasymnot}\isaconst{iChk}\ \isactrlbold v\isactrlsub {\isadigit{0}}\ {\isacharparenleft}{\kern0pt}{\isasymupharpoonleft}{\isacharparenleft}{\kern0pt}\isactrlbold n{\isasymlbrace}\isactrlbold {\isasymnot}\isactrlbold {\isasymnot}\isactrlbold A\ \isactrlbold {\isasymcdot}\ {\isacharparenleft}{\kern0pt}\isactrlbold {\isasymlambda}{\isacharparenleft}{\kern0pt}\isactrlbold {\isasymnot}\isaconst{iChk}\ \isactrlbold v\isactrlsub {\isadigit{0}}\ {\isacharparenleft}{\kern0pt}{\isasymupharpoonleft}{\isacharparenleft}{\kern0pt}\isactrlbold n{\isasymlbrace}\isactrlbold {\isasymnot}\isactrlbold {\isadigit{0}}\ \isactrlbold {\isacharequal}{\kern0pt}\ \isactrlbold {\isadigit{0}}{\isasymrbrace}{\isacharparenright}{\kern0pt}{\isacharparenright}{\kern0pt}{\isacharparenright}{\kern0pt}{\isacharparenright}{\kern0pt}{\isasymrbrace}{\isacharparenright}{\kern0pt}{\isacharparenright}{\kern0pt}{\isacharparenright}{\kern0pt}{\isacharparenright}{\kern0pt}\ {\isasymtriangleright}\ {\isasymemptyset}\ {\isasymturnstile}\ \isactrlbold {\isasymnot}\isactrlbold A\ \isactrlbold {\isasymcdot}\ {\isacharparenleft}{\kern0pt}\isactrlbold {\isasymlambda}{\isacharparenleft}{\kern0pt}\isactrlbold {\isasymnot}\isaconst{iChk}\ \isactrlbold v\isactrlsub {\isadigit{0}}\ {\isacharparenleft}{\kern0pt}{\isasymupharpoonleft}{\isacharparenleft}{\kern0pt}\isactrlbold n{\isasymlbrace}\isactrlbold {\isasymnot}\isactrlbold {\isadigit{0}}\ \isactrlbold {\isacharequal}{\kern0pt}\ \isactrlbold {\isadigit{0}}{\isasymrbrace}{\isacharparenright}{\kern0pt}{\isacharparenright}{\kern0pt}{\isacharparenright}{\kern0pt}{\isacharparenright}{\kern0pt}\ {\isasymand}\ \isactrlbold {\isasymnot}\isactrlbold A\ \isactrlbold {\isasymcdot}\ {\isacharparenleft}{\kern0pt}\isactrlbold {\isasymlambda}{\isacharparenleft}{\kern0pt}\isactrlbold {\isasymnot}\isaconst{iChk}\ \isactrlbold v\isactrlsub {\isadigit{0}}\ {\isacharparenleft}{\kern0pt}{\isasymupharpoonleft}{\isacharparenleft}{\kern0pt}\isactrlbold n{\isasymlbrace}\isactrlbold {\isasymnot}\isactrlbold {\isadigit{0}}\ \isactrlbold {\isacharequal}{\kern0pt}\ \isactrlbold {\isadigit{0}}{\isasymrbrace}{\isacharparenright}{\kern0pt}{\isacharparenright}{\kern0pt}{\isacharparenright}{\kern0pt}{\isacharparenright}{\kern0pt}\ {\isasymtriangleright}\ {\isasymemptyset}\ {\isasymturnstile}\ \isactrlbold {\isasymnot}\isaconst{goedelG}}%
}%
\SNIP{Goedel-Paper-Statements:thm:iChk-def-0}{%
\isakeywordONE{thm}\isamarkupfalse%
\ iChk{\isacharunderscore}{\kern0pt}def%
\isanewline%
}%
\SNIP{Goedel-Paper-Statements:thm:iChk-def-1}{%
\isa{\isaconst{iChk}\ \isavar{p}\ \isavar{c}\ {\isasymequiv}\ {\isacharparenleft}{\kern0pt}\isaconst{cmp}\ \isaconst{validChkIdx}\ \isactrlbold {\isasymcdot}\ \isavar{p}\ \isactrlbold {\isasymcdot}\ \isavar{p}\ \isactrlbold {\isacharequal}{\kern0pt}\ \isactrlbold n\isaconst{Suc}\ {\isadigit{0}}\ \isactrlbold {\isasymand}\ \isavar{p}\ \isactrlbold {\isasymnoteq}\ \isactrlbold {\isadigit{0}}{\isacharparenright}{\kern0pt}\ \isactrlbold {\isasymand}\ \isaconst{cmp}\ \isaconst{length{\isacharunderscore}{\kern0pt}encIdx}\ \isactrlbold {\isasymcdot}\ \isaconst{cpFst}\ {\isacharparenleft}{\kern0pt}\isaconst{cpFst}\ {\isacharparenleft}{\kern0pt}\isaconst{trmPredR}\ \isavar{p}{\isacharparenright}{\kern0pt}{\isacharparenright}{\kern0pt}\ \isactrlbold {\isasymcdot}\ \isaconst{cpFst}\ {\isacharparenleft}{\kern0pt}\isaconst{cpFst}\ {\isacharparenleft}{\kern0pt}\isaconst{trmPredR}\ \isavar{p}{\isacharparenright}{\kern0pt}{\isacharparenright}{\kern0pt}\ \isactrlbold {\isacharequal}{\kern0pt}\ \isactrlbold {\isadigit{0}}\ \isactrlbold {\isasymand}\ \isaconst{cpSnd}\ {\isacharparenleft}{\kern0pt}\isaconst{cpFst}\ {\isacharparenleft}{\kern0pt}\isaconst{trmPredR}\ \isavar{p}{\isacharparenright}{\kern0pt}{\isacharparenright}{\kern0pt}\ \isactrlbold {\isacharequal}{\kern0pt}\ \isavar{c}}%
}%
\SNIP{Goedel-Paper-Statements:thm:cEq01-def-0}{%
\isakeywordONE{thm}\isamarkupfalse%
\ cEq{\isadigit{0}}{\isadigit{1}}{\isacharunderscore}{\kern0pt}def%
\isanewline%
}%
\SNIP{Goedel-Paper-Statements:thm:cEq01-def-1}{%
\isa{\isaconst{cEq{\isadigit{0}}{\isadigit{1}}}\ {\isasymequiv}\ \isactrlbold n{\isasymlbrace}\isactrlbold {\isadigit{0}}\ \isactrlbold {\isacharequal}{\kern0pt}\ \isactrlbold n\isaconst{Suc}\ {\isadigit{0}}{\isasymrbrace}}%
}%
\SNIP{Goedel-Paper-Statements:thm:adiW-def-0}{%
\isakeywordONE{thm}\isamarkupfalse%
\ adiW{\isacharunderscore}{\kern0pt}def%
\isanewline%
}%
\SNIP{Goedel-Paper-Statements:thm:adiW-def-1}{%
\isa{\isaconst{adiW}\ \isavar{F}\ \isavar{u}\ {\isasymequiv}\ \isactrlbold A\ \isactrlbold {\isasymcdot}\ {\isacharparenleft}{\kern0pt}\isactrlbold {\isasymlambda}\isaconst{iChkD}\ {\isacharparenleft}{\kern0pt}\isaconst{drvChkIdxG}\ {\isacharparenleft}{\kern0pt}\isaconst{drvChkList}\ {\isacharat}{\kern0pt}\ {\isacharbrackleft}{\kern0pt}\isaconst{chkGfixIdx}{\isacharbrackright}{\kern0pt}{\isacharparenright}{\kern0pt}{\isacharparenright}{\kern0pt}\ {\isacharparenleft}{\kern0pt}\isaconst{cmp}\ \isavar{F}\ \isactrlbold {\isasymcdot}\ \isactrlbold v\isactrlsub {\isadigit{0}}\ \isactrlbold {\isasymcdot}\ \isactrlbold v\isactrlsub {\isadigit{0}}{\isacharparenright}{\kern0pt}\ {\isacharparenleft}{\kern0pt}\isaconst{adrInst}\ \isavar{u}\ \isactrlbold v\isactrlsub {\isadigit{0}}{\isacharparenright}{\kern0pt}{\isacharparenright}{\kern0pt}}%
}%
\SNIP{Goedel-Paper-Statements:context:gcalc-0}{%
\isakeywordONE{context}\isamarkupfalse%
\ gcalc%
}%
\SNIP{Goedel-Paper-Statements:thm:rga-embed-0}{%
\isakeywordONE{thm}\isamarkupfalse%
\ rga{\isacharunderscore}{\kern0pt}embed%
\isanewline%
}%
\SNIP{Goedel-Paper-Statements:thm:rga-embed-1}{%
\isa{\isaconst{prf{\isacharunderscore}{\kern0pt}derives{\isachardot}{\kern0pt}entails}\ \isaconst{rga{\isachardot}{\kern0pt}derives}\ \isavar{H}\ \isavar{c}\ {\isasymLongrightarrow}\ \isavar{H}\ {\isasymturnstile}\ \isavar{c}}%
}%
\SNIP{Goedel-Paper-Statements:context:gcalc-strip-0}{%
\isakeywordONE{context}\isamarkupfalse%
\ gcalc{\isacharunderscore}{\kern0pt}strip%
}%
\SNIP{Goedel-Paper-Statements:thm:sigma-deduction-theorem+1-0}{%
\isakeywordONE{thm}\isamarkupfalse%
\ sigma{\isacharunderscore}{\kern0pt}deduction{\isacharunderscore}{\kern0pt}theorem%
\isanewline%
}%
\SNIP{Goedel-Paper-Statements:thm:sigma-deduction-theorem+1-1}{%
\isa{\isactrlbold v\isactrlsub {\isadigit{0}}\ \isactrlbold N\ {\isasymtriangleright}\ {\isasymupharpoonleft}\isavar{{\isasymGamma}}\ {\isasymturnstile}\ \isavar{D}\ \isactrlbold B\ {\isasymLongrightarrow}\ {\isacharparenleft}{\kern0pt}\isactrlbold {\isasymnot}\isactrlbold A\ \isactrlbold {\isasymcdot}\ {\isacharparenleft}{\kern0pt}\isactrlbold {\isasymlambda}{\isacharparenleft}{\kern0pt}\isactrlbold {\isasymnot}\isavar{D}{\isacharparenright}{\kern0pt}{\isacharparenright}{\kern0pt}\ {\isasymtriangleright}\ \isavar{{\isasymGamma}}\ {\isasymturnstile}\ \isavar{C}{\isacharparenright}{\kern0pt}\ {\isacharequal}{\kern0pt}\ {\isacharparenleft}{\kern0pt}\isavar{{\isasymGamma}}\ {\isasymturnstile}\ \isactrlbold A\ \isactrlbold {\isasymcdot}\ {\isacharparenleft}{\kern0pt}\isactrlbold {\isasymlambda}{\isacharparenleft}{\kern0pt}\isavar{D}\ \isactrlbold {\isasymlongrightarrow}\ {\isasymupharpoonleft}\isavar{C}{\isacharparenright}{\kern0pt}{\isacharparenright}{\kern0pt}{\isacharparenright}{\kern0pt}}%
}%
\SNIP{Goedel-Paper-Statements:context:gstep-0}{%
\isakeywordONE{context}\isamarkupfalse%
\ gstep%
}%
\SNIP{Goedel-Paper-Statements:thm:is-gcalc-0}{%
\isakeywordONE{thm}\isamarkupfalse%
\ is{\isacharunderscore}{\kern0pt}gcalc%
\isanewline%
}%
\SNIP{Goedel-Paper-Statements:thm:is-gcalc-1}{%
\isa{\isaconst{gcalc}\ {\isacharparenleft}{\kern0pt}\isaconst{prf{\isacharunderscore}{\kern0pt}derives{\isachardot}{\kern0pt}entails}\ \isafree{S}{\isacharparenright}{\kern0pt}}%
}%
\SNIP{Goedel-Paper-Statements:context:g2-syn-0}{%
\isakeywordONE{context}\isamarkupfalse%
\ g{\isadigit{2}}{\isacharunderscore}{\kern0pt}syn%
}%
\SNIP{Goedel-Paper-Statements:thm:ConT-def-0}{%
\isakeywordONE{thm}\isamarkupfalse%
\ ConT{\isacharunderscore}{\kern0pt}def%
\isanewline%
}%
\SNIP{Goedel-Paper-Statements:thm:ConT-def-1}{%
\isa{\isaconst{ConT}\ {\isasymequiv}\ \isactrlbold A\ \isactrlbold {\isasymcdot}\ {\isacharparenleft}{\kern0pt}\isactrlbold {\isasymlambda}{\isacharparenleft}{\kern0pt}\isactrlbold {\isasymnot}\isafree{chk}\ \isactrlbold v\isactrlsub {\isadigit{0}}\ \isaconst{rga{\isachardot}{\kern0pt}cEq{\isadigit{0}}{\isadigit{1}}}{\isacharparenright}{\kern0pt}{\isacharparenright}{\kern0pt}}%
}%
\SNIP{Goedel-Paper-Statements:thm:goedelGT-def-0}{%
\isakeywordONE{thm}\isamarkupfalse%
\ goedelGT{\isacharunderscore}{\kern0pt}def%
\isanewline%
}%
\SNIP{Goedel-Paper-Statements:thm:goedelGT-def-1}{%
\isa{\isaconst{goedelGT}\ {\isasymequiv}\ \isactrlbold A\ \isactrlbold {\isasymcdot}\ {\isacharparenleft}{\kern0pt}\isactrlbold {\isasymlambda}{\isacharparenleft}{\kern0pt}\isactrlbold {\isasymnot}\isafree{chk}\ \isactrlbold v\isactrlsub {\isadigit{0}}\ \isaconst{goedelST}{\isacharparenright}{\kern0pt}{\isacharparenright}{\kern0pt}}%
}%
\SNIP{Goedel-Paper-Statements:thm:goedelST-eq-0}{%
\isakeywordONE{thm}\isamarkupfalse%
\ goedelST{\isacharunderscore}{\kern0pt}eq%
\isanewline%
}%
\SNIP{Goedel-Paper-Statements:thm:goedelST-eq-1}{%
\isa{{\isasymemptyset}\ {\isasymturnstile}\ \isaconst{goedelST}\ \isactrlbold {\isacharequal}{\kern0pt}\ \isactrlbold n{\isasymlbrace}\isaconst{goedelGT}{\isasymrbrace}}%
}%
\SNIP{Goedel-Paper-Statements:thm:goedelT-diag-ent-0}{%
\isakeywordONE{thm}\isamarkupfalse%
\ goedelT{\isacharunderscore}{\kern0pt}diag{\isacharunderscore}{\kern0pt}ent%
\isanewline%
}%
\SNIP{Goedel-Paper-Statements:thm:goedelT-diag-ent-1}{%
\isa{\isaconst{goedelGT}\ {\isasymtriangleright}\ {\isasymemptyset}\ {\isasymturnstile}\ \isactrlbold {\isasymnot}\isactrlbold {\isasymnot}\isactrlbold A\ \isactrlbold {\isasymcdot}\ {\isacharparenleft}{\kern0pt}\isactrlbold {\isasymlambda}{\isacharparenleft}{\kern0pt}\isactrlbold {\isasymnot}\isafree{chk}\ \isactrlbold v\isactrlsub {\isadigit{0}}\ {\isacharparenleft}{\kern0pt}{\isasymupharpoonleft}{\isacharparenleft}{\kern0pt}\isactrlbold n{\isasymlbrace}\isaconst{goedelGT}{\isasymrbrace}{\isacharparenright}{\kern0pt}{\isacharparenright}{\kern0pt}{\isacharparenright}{\kern0pt}{\isacharparenright}{\kern0pt}\ {\isasymand}\ \isactrlbold {\isasymnot}\isactrlbold {\isasymnot}\isactrlbold A\ \isactrlbold {\isasymcdot}\ {\isacharparenleft}{\kern0pt}\isactrlbold {\isasymlambda}{\isacharparenleft}{\kern0pt}\isactrlbold {\isasymnot}\isafree{chk}\ \isactrlbold v\isactrlsub {\isadigit{0}}\ {\isacharparenleft}{\kern0pt}{\isasymupharpoonleft}{\isacharparenleft}{\kern0pt}\isactrlbold n{\isasymlbrace}\isaconst{goedelGT}{\isasymrbrace}{\isacharparenright}{\kern0pt}{\isacharparenright}{\kern0pt}{\isacharparenright}{\kern0pt}{\isacharparenright}{\kern0pt}\ {\isasymtriangleright}\ {\isasymemptyset}\ {\isasymturnstile}\ \isaconst{goedelGT}}%
}%
\SNIP{Goedel-Paper-Statements:context:g2-ent-0}{%
\isakeywordONE{context}\isamarkupfalse%
\ g{\isadigit{2}}{\isacharunderscore}{\kern0pt}ent%
}%
\SNIP{Goedel-Paper-Statements:thm:ent-necT-0}{%
\isakeywordONE{thm}\isamarkupfalse%
\ ent{\isacharunderscore}{\kern0pt}necT%
\isanewline%
}%
\SNIP{Goedel-Paper-Statements:thm:ent-necT-1}{%
\isa{\isavar{A}\ {\isasymtriangleright}\ {\isasymemptyset}\ {\isasymturnstile}\ \isavar{B}\ {\isasymLongrightarrow}\ \isactrlbold {\isasymnot}\isactrlbold A\ \isactrlbold {\isasymcdot}\ {\isacharparenleft}{\kern0pt}\isactrlbold {\isasymlambda}{\isacharparenleft}{\kern0pt}\isactrlbold {\isasymnot}\isafree{chk}\ \isactrlbold v\isactrlsub {\isadigit{0}}\ {\isacharparenleft}{\kern0pt}{\isasymupharpoonleft}{\isacharparenleft}{\kern0pt}\isactrlbold n{\isasymlbrace}\isavar{A}{\isasymrbrace}{\isacharparenright}{\kern0pt}{\isacharparenright}{\kern0pt}{\isacharparenright}{\kern0pt}{\isacharparenright}{\kern0pt}\ {\isasymtriangleright}\ \isavar{{\isasymGamma}}\ {\isasymturnstile}\ \isactrlbold {\isasymnot}\isactrlbold A\ \isactrlbold {\isasymcdot}\ {\isacharparenleft}{\kern0pt}\isactrlbold {\isasymlambda}{\isacharparenleft}{\kern0pt}\isactrlbold {\isasymnot}\isafree{chk}\ \isactrlbold v\isactrlsub {\isadigit{0}}\ {\isacharparenleft}{\kern0pt}{\isasymupharpoonleft}{\isacharparenleft}{\kern0pt}\isactrlbold n{\isasymlbrace}\isavar{B}{\isasymrbrace}{\isacharparenright}{\kern0pt}{\isacharparenright}{\kern0pt}{\isacharparenright}{\kern0pt}{\isacharparenright}{\kern0pt}}%
}%
\SNIP{Goedel-Paper-Statements:thm:D3-entT-0}{%
\isakeywordONE{thm}\isamarkupfalse%
\ D{\isadigit{3}}{\isacharunderscore}{\kern0pt}entT%
\isanewline%
}%
\SNIP{Goedel-Paper-Statements:thm:D3-entT-1}{%
\isa{\isactrlbold {\isasymnot}\isactrlbold A\ \isactrlbold {\isasymcdot}\ {\isacharparenleft}{\kern0pt}\isactrlbold {\isasymlambda}{\isacharparenleft}{\kern0pt}\isactrlbold {\isasymnot}\isafree{chk}\ \isactrlbold v\isactrlsub {\isadigit{0}}\ {\isacharparenleft}{\kern0pt}{\isasymupharpoonleft}{\isacharparenleft}{\kern0pt}\isactrlbold n{\isasymlbrace}\isavar{{\isasymphi}}{\isasymrbrace}{\isacharparenright}{\kern0pt}{\isacharparenright}{\kern0pt}{\isacharparenright}{\kern0pt}{\isacharparenright}{\kern0pt}\ {\isasymtriangleright}\ \isavar{{\isasymGamma}}\ {\isasymturnstile}\ \isactrlbold {\isasymnot}\isactrlbold A\ \isactrlbold {\isasymcdot}\ {\isacharparenleft}{\kern0pt}\isactrlbold {\isasymlambda}{\isacharparenleft}{\kern0pt}\isactrlbold {\isasymnot}\isafree{chk}\ \isactrlbold v\isactrlsub {\isadigit{0}}\ {\isacharparenleft}{\kern0pt}{\isasymupharpoonleft}{\isacharparenleft}{\kern0pt}\isactrlbold n{\isasymlbrace}\isactrlbold {\isasymnot}\isactrlbold A\ \isactrlbold {\isasymcdot}\ {\isacharparenleft}{\kern0pt}\isactrlbold {\isasymlambda}{\isacharparenleft}{\kern0pt}\isactrlbold {\isasymnot}\isafree{chk}\ \isactrlbold v\isactrlsub {\isadigit{0}}\ {\isacharparenleft}{\kern0pt}{\isasymupharpoonleft}{\isacharparenleft}{\kern0pt}\isactrlbold n{\isasymlbrace}\isavar{{\isasymphi}}{\isasymrbrace}{\isacharparenright}{\kern0pt}{\isacharparenright}{\kern0pt}{\isacharparenright}{\kern0pt}{\isacharparenright}{\kern0pt}{\isasymrbrace}{\isacharparenright}{\kern0pt}{\isacharparenright}{\kern0pt}{\isacharparenright}{\kern0pt}{\isacharparenright}{\kern0pt}}%
}%
\SNIP{Goedel-Paper-Statements:thm:goedelGT-under-ConT-textbook-0}{%
\isakeywordONE{thm}\isamarkupfalse%
\ goedelGT{\isacharunderscore}{\kern0pt}under{\isacharunderscore}{\kern0pt}ConT{\isacharunderscore}{\kern0pt}textbook%
\isanewline%
}%
\SNIP{Goedel-Paper-Statements:thm:goedelGT-under-ConT-textbook-1}{%
\isa{\isaconst{ConT}\ {\isasymtriangleright}\ {\isasymemptyset}\ {\isasymturnstile}\ \isaconst{goedelGT}}%
}%
\SNIP{Goedel-Paper-Statements:thm:goedel2-internal-0}{%
\isakeywordONE{thm}\isamarkupfalse%
\ goedel{\isadigit{2}}{\isacharunderscore}{\kern0pt}internal%
\isanewline%
}%
\SNIP{Goedel-Paper-Statements:thm:goedel2-internal-1}{%
\isa{\isaconst{ConT}\ {\isasymtriangleright}\ {\isasymemptyset}\ {\isasymturnstile}\ \isactrlbold {\isasymnot}\isactrlbold {\isasymnot}\isactrlbold A\ \isactrlbold {\isasymcdot}\ {\isacharparenleft}{\kern0pt}\isactrlbold {\isasymlambda}{\isacharparenleft}{\kern0pt}\isactrlbold {\isasymnot}\isafree{chk}\ \isactrlbold v\isactrlsub {\isadigit{0}}\ {\isacharparenleft}{\kern0pt}{\isasymupharpoonleft}{\isacharparenleft}{\kern0pt}\isactrlbold n{\isasymlbrace}\isaconst{ConT}{\isasymrbrace}{\isacharparenright}{\kern0pt}{\isacharparenright}{\kern0pt}{\isacharparenright}{\kern0pt}{\isacharparenright}{\kern0pt}}%
}%
\SNIP{Goedel-Paper-Statements:context:g2-con-0}{%
\isakeywordONE{context}\isamarkupfalse%
\ g{\isadigit{2}}{\isacharunderscore}{\kern0pt}con%
}%
\SNIP{Goedel-Paper-Statements:thm:goedel1-0}{%
\isakeywordONE{thm}\isamarkupfalse%
\ goedel{\isadigit{1}}%
\isanewline%
}%
\SNIP{Goedel-Paper-Statements:thm:goedel1-1}{%
\isa{{\isasymnot}\ {\isasymemptyset}\ {\isasymturnstile}\ \isaconst{goedelGT}}%
}%
\SNIP{Goedel-Paper-Statements:thm:goedel2-textbook-0}{%
\isakeywordONE{thm}\isamarkupfalse%
\ goedel{\isadigit{2}}{\isacharunderscore}{\kern0pt}textbook%
\isanewline%
}%
\SNIP{Goedel-Paper-Statements:thm:goedel2-textbook-1}{%
\isa{{\isasymnot}\ {\isasymemptyset}\ {\isasymturnstile}\ \isaconst{ConT}}%
}%
\SNIP{Goedel-Paper-Statements:thm:goedel2T-of-expl-d3g-0}{%
\isakeywordONE{thm}\isamarkupfalse%
\ goedel{\isadigit{2}}T{\isacharunderscore}{\kern0pt}of{\isacharunderscore}{\kern0pt}expl{\isacharunderscore}{\kern0pt}d{\isadigit{3}}g%
\isanewline%
}%
\SNIP{Goedel-Paper-Statements:thm:goedel2T-of-expl-d3g-1}{%
\isa{{\isasymlbrakk}\isaconst{EXPLT}\ \isavar{e}{\isacharsemicolon}{\kern0pt}\ \isaconst{D{\isadigit{3}}GT}\ \isavar{{\isasymsigma}}{\isasymrbrakk}\ {\isasymLongrightarrow}\ {\isasymnot}\ {\isasymemptyset}\ {\isasymturnstile}\ \isaconst{ConT}}%
}%
\SNIP{Goedel-Paper-Statements:thm:no-loeb-hypothesis-0eq1T-0}{%
\isakeywordONE{thm}\isamarkupfalse%
\ no{\isacharunderscore}{\kern0pt}loeb{\isacharunderscore}{\kern0pt}hypothesis{\isacharunderscore}{\kern0pt}{\isadigit{0}}eq{\isadigit{1}}T%
\isanewline%
}%
\SNIP{Goedel-Paper-Statements:thm:no-loeb-hypothesis-0eq1T-1}{%
\isa{{\isasymnot}\ {\isasymemptyset}\ {\isasymturnstile}\ \isactrlbold A\ \isactrlbold {\isasymcdot}\ {\isacharparenleft}{\kern0pt}\isactrlbold {\isasymlambda}{\isacharparenleft}{\kern0pt}\isafree{chk}\ \isactrlbold v\isactrlsub {\isadigit{0}}\ \isaconst{rga{\isachardot}{\kern0pt}cEq{\isadigit{0}}{\isadigit{1}}}\ \isactrlbold {\isasymlongrightarrow}\ \isactrlbold {\isadigit{0}}\ \isactrlbold {\isacharequal}{\kern0pt}\ \isactrlbold n{\isadigit{1}}{\isacharparenright}{\kern0pt}{\isacharparenright}{\kern0pt}}%
}%
\SNIP{Goedel-Paper-Statements:context:g2-exf-0}{%
\isakeywordONE{context}\isamarkupfalse%
\ g{\isadigit{2}}{\isacharunderscore}{\kern0pt}exf%
}%
\SNIP{Goedel-Paper-Statements:thm:goedelGT-con-interderivable-0}{%
\isakeywordONE{thm}\isamarkupfalse%
\ goedelGT{\isacharunderscore}{\kern0pt}con{\isacharunderscore}{\kern0pt}interderivable%
\isanewline%
}%
\SNIP{Goedel-Paper-Statements:thm:goedelGT-con-interderivable-1}{%
\isa{\isaconst{goedelGT}\ {\isasymtriangleright}\ {\isasymemptyset}\ {\isasymturnstile}\ \isaconst{ConT}\ {\isasymand}\ \isaconst{ConT}\ {\isasymtriangleright}\ {\isasymemptyset}\ {\isasymturnstile}\ \isaconst{goedelGT}}%
}%
\SNIP{Goedel-Paper-Statements:context:g2-d2-0}{%
\isakeywordONE{context}\isamarkupfalse%
\ g{\isadigit{2}}{\isacharunderscore}{\kern0pt}d{\isadigit{2}}%
}%
\SNIP{Goedel-Paper-Statements:thm:D2-entT-0}{%
\isakeywordONE{thm}\isamarkupfalse%
\ D{\isadigit{2}}{\isacharunderscore}{\kern0pt}entT%
\isanewline%
}%
\SNIP{Goedel-Paper-Statements:thm:D2-entT-1}{%
\isa{{\isasymlbrakk}\isaconst{wft}\ \isavar{{\isasymphi}}{\isacharsemicolon}{\kern0pt}\ \isaconst{wft}\ \isavar{{\isasympsi}}{\isasymrbrakk}\ {\isasymLongrightarrow}\ \isactrlbold {\isasymnot}\isactrlbold A\ \isactrlbold {\isasymcdot}\ {\isacharparenleft}{\kern0pt}\isactrlbold {\isasymlambda}{\isacharparenleft}{\kern0pt}\isactrlbold {\isasymnot}\isafree{chk}\ \isactrlbold v\isactrlsub {\isadigit{0}}\ {\isacharparenleft}{\kern0pt}{\isasymupharpoonleft}{\isacharparenleft}{\kern0pt}\isactrlbold n{\isasymlbrace}\isavar{{\isasymphi}}{\isasymrbrace}{\isacharparenright}{\kern0pt}{\isacharparenright}{\kern0pt}{\isacharparenright}{\kern0pt}{\isacharparenright}{\kern0pt}\ {\isasymtriangleright}\ \isactrlbold {\isasymnot}\isactrlbold A\ \isactrlbold {\isasymcdot}\ {\isacharparenleft}{\kern0pt}\isactrlbold {\isasymlambda}{\isacharparenleft}{\kern0pt}\isactrlbold {\isasymnot}\isafree{chk}\ \isactrlbold v\isactrlsub {\isadigit{0}}\ {\isacharparenleft}{\kern0pt}{\isasymupharpoonleft}{\isacharparenleft}{\kern0pt}\isactrlbold n{\isasymlbrace}\isavar{{\isasymphi}}\ \isactrlbold {\isasymlongrightarrow}\ \isavar{{\isasympsi}}{\isasymrbrace}{\isacharparenright}{\kern0pt}{\isacharparenright}{\kern0pt}{\isacharparenright}{\kern0pt}{\isacharparenright}{\kern0pt}\ {\isasymtriangleright}\ \isavar{{\isasymGamma}}\ {\isasymturnstile}\ \isactrlbold {\isasymnot}\isactrlbold A\ \isactrlbold {\isasymcdot}\ {\isacharparenleft}{\kern0pt}\isactrlbold {\isasymlambda}{\isacharparenleft}{\kern0pt}\isactrlbold {\isasymnot}\isafree{chk}\ \isactrlbold v\isactrlsub {\isadigit{0}}\ {\isacharparenleft}{\kern0pt}{\isasymupharpoonleft}{\isacharparenleft}{\kern0pt}\isactrlbold n{\isasymlbrace}\isavar{{\isasympsi}}{\isasymrbrace}{\isacharparenright}{\kern0pt}{\isacharparenright}{\kern0pt}{\isacharparenright}{\kern0pt}{\isacharparenright}{\kern0pt}}%
}%
\SNIP{Goedel-Paper-Statements:thm:formalized-loeb-entT-0}{%
\isakeywordONE{thm}\isamarkupfalse%
\ formalized{\isacharunderscore}{\kern0pt}loeb{\isacharunderscore}{\kern0pt}entT%
\isanewline%
}%
\SNIP{Goedel-Paper-Statements:thm:formalized-loeb-entT-1}{%
\isa{{\isasymlbrakk}\isaconst{Syntax{\isacharunderscore}{\kern0pt}Term{\isachardot}{\kern0pt}closed}\ \isavar{{\isasympsi}}{\isacharsemicolon}{\kern0pt}\ \isaconst{wft}\ \isavar{{\isasympsi}}{\isasymrbrakk}\ {\isasymLongrightarrow}\ \isactrlbold {\isasymnot}\isactrlbold A\ \isactrlbold {\isasymcdot}\ {\isacharparenleft}{\kern0pt}\isactrlbold {\isasymlambda}{\isacharparenleft}{\kern0pt}\isactrlbold {\isasymnot}\isafree{chk}\ \isactrlbold v\isactrlsub {\isadigit{0}}\ {\isacharparenleft}{\kern0pt}{\isasymupharpoonleft}{\isacharparenleft}{\kern0pt}\isactrlbold n{\isasymlbrace}\isactrlbold {\isasymnot}\isactrlbold A\ \isactrlbold {\isasymcdot}\ {\isacharparenleft}{\kern0pt}\isactrlbold {\isasymlambda}{\isacharparenleft}{\kern0pt}\isactrlbold {\isasymnot}\isafree{chk}\ \isactrlbold v\isactrlsub {\isadigit{0}}\ {\isacharparenleft}{\kern0pt}{\isasymupharpoonleft}{\isacharparenleft}{\kern0pt}\isactrlbold n{\isasymlbrace}\isavar{{\isasympsi}}{\isasymrbrace}{\isacharparenright}{\kern0pt}{\isacharparenright}{\kern0pt}{\isacharparenright}{\kern0pt}{\isacharparenright}{\kern0pt}\ \isactrlbold {\isasymlongrightarrow}\ \isavar{{\isasympsi}}{\isasymrbrace}{\isacharparenright}{\kern0pt}{\isacharparenright}{\kern0pt}{\isacharparenright}{\kern0pt}{\isacharparenright}{\kern0pt}\ {\isasymtriangleright}\ {\isasymemptyset}\ {\isasymturnstile}\ \isactrlbold {\isasymnot}\isactrlbold A\ \isactrlbold {\isasymcdot}\ {\isacharparenleft}{\kern0pt}\isactrlbold {\isasymlambda}{\isacharparenleft}{\kern0pt}\isactrlbold {\isasymnot}\isafree{chk}\ \isactrlbold v\isactrlsub {\isadigit{0}}\ {\isacharparenleft}{\kern0pt}{\isasymupharpoonleft}{\isacharparenleft}{\kern0pt}\isactrlbold n{\isasymlbrace}\isavar{{\isasympsi}}{\isasymrbrace}{\isacharparenright}{\kern0pt}{\isacharparenright}{\kern0pt}{\isacharparenright}{\kern0pt}{\isacharparenright}{\kern0pt}}%
}%
\SNIP{Goedel-Paper-Statements:context:g2-exf-acc-0}{%
\isakeywordONE{context}\isamarkupfalse%
\ g{\isadigit{2}}{\isacharunderscore}{\kern0pt}exf{\isacharunderscore}{\kern0pt}acc%
}%
\SNIP{Goedel-Paper-Statements:thm:bsT-structure-conditions-0}{%
\isakeywordONE{thm}\isamarkupfalse%
\ bsT{\isacharunderscore}{\kern0pt}structure{\isacharunderscore}{\kern0pt}conditions%
\isanewline%
}%
\SNIP{Goedel-Paper-Statements:thm:bsT-structure-conditions-1}{%
\isa{{\isacharparenleft}{\kern0pt}{\isasymforall}\isabound{x}{\isachardot}{\kern0pt}\ \isaconst{wft}\ \isabound{x}\ {\isasymlongrightarrow}\ \isabound{x}\ {\isasymtriangleright}\ {\isasymemptyset}\ {\isasymturnstile}\ \isabound{x}{\isacharparenright}{\kern0pt}\ {\isasymand}\ {\isacharparenleft}{\kern0pt}{\isasymforall}\isabound{x}\ \isabound{y}\ \isabound{z}{\isachardot}{\kern0pt}\ \isabound{x}\ {\isasymtriangleright}\ {\isasymemptyset}\ {\isasymturnstile}\ \isabound{y}\ {\isasymlongrightarrow}\ \isabound{y}\ {\isasymtriangleright}\ {\isasymemptyset}\ {\isasymturnstile}\ \isabound{z}\ {\isasymlongrightarrow}\ \isabound{x}\ {\isasymtriangleright}\ {\isasymemptyset}\ {\isasymturnstile}\ \isabound{z}{\isacharparenright}{\kern0pt}\ {\isasymand}\ {\isacharparenleft}{\kern0pt}{\isasymforall}\isabound{x}\ \isabound{y}{\isachardot}{\kern0pt}\ \isabound{x}\ {\isasymtriangleright}\ {\isasymemptyset}\ {\isasymturnstile}\ \isabound{y}\ {\isasymlongrightarrow}\ \isactrlbold {\isasymnot}\isactrlbold A\ \isactrlbold {\isasymcdot}\ {\isacharparenleft}{\kern0pt}\isactrlbold {\isasymlambda}{\isacharparenleft}{\kern0pt}\isactrlbold {\isasymnot}\isafree{chk}\ \isactrlbold v\isactrlsub {\isadigit{0}}\ {\isacharparenleft}{\kern0pt}{\isasymupharpoonleft}{\isacharparenleft}{\kern0pt}\isactrlbold n{\isasymlbrace}\isabound{x}{\isasymrbrace}{\isacharparenright}{\kern0pt}{\isacharparenright}{\kern0pt}{\isacharparenright}{\kern0pt}{\isacharparenright}{\kern0pt}\ {\isasymtriangleright}\ {\isasymemptyset}\ {\isasymturnstile}\ \isactrlbold {\isasymnot}\isactrlbold A\ \isactrlbold {\isasymcdot}\ {\isacharparenleft}{\kern0pt}\isactrlbold {\isasymlambda}{\isacharparenleft}{\kern0pt}\isactrlbold {\isasymnot}\isafree{chk}\ \isactrlbold v\isactrlsub {\isadigit{0}}\ {\isacharparenleft}{\kern0pt}{\isasymupharpoonleft}{\isacharparenleft}{\kern0pt}\isactrlbold n{\isasymlbrace}\isabound{y}{\isasymrbrace}{\isacharparenright}{\kern0pt}{\isacharparenright}{\kern0pt}{\isacharparenright}{\kern0pt}{\isacharparenright}{\kern0pt}{\isacharparenright}{\kern0pt}\ {\isasymand}\ \isactrlbold {\isadigit{0}}\ \isactrlbold {\isacharequal}{\kern0pt}\ \isactrlbold {\isadigit{0}}\ {\isasymtriangleright}\ {\isasymemptyset}\ {\isasymturnstile}\ \isactrlbold {\isasymnot}\isactrlbold A\ \isactrlbold {\isasymcdot}\ {\isacharparenleft}{\kern0pt}\isactrlbold {\isasymlambda}{\isacharparenleft}{\kern0pt}\isactrlbold {\isasymnot}\isafree{chk}\ \isactrlbold v\isactrlsub {\isadigit{0}}\ {\isacharparenleft}{\kern0pt}{\isasymupharpoonleft}{\isacharparenleft}{\kern0pt}\isactrlbold n{\isasymlbrace}\isactrlbold {\isasymnot}\isactrlbold {\isasymnot}\isactrlbold {\isadigit{0}}\ \isactrlbold {\isacharequal}{\kern0pt}\ \isactrlbold {\isadigit{0}}{\isasymrbrace}{\isacharparenright}{\kern0pt}{\isacharparenright}{\kern0pt}{\isacharparenright}{\kern0pt}{\isacharparenright}{\kern0pt}\ {\isasymand}\ {\isacharparenleft}{\kern0pt}{\isasymforall}\isabound{x}\ \isabound{y}{\isachardot}{\kern0pt}\ \isabound{x}\ {\isasymtriangleright}\ {\isasymemptyset}\ {\isasymturnstile}\ \isactrlbold {\isasymnot}\isactrlbold A\ \isactrlbold {\isasymcdot}\ {\isacharparenleft}{\kern0pt}\isactrlbold {\isasymlambda}{\isacharparenleft}{\kern0pt}\isactrlbold {\isasymnot}\isafree{chk}\ \isactrlbold v\isactrlsub {\isadigit{0}}\ {\isacharparenleft}{\kern0pt}{\isasymupharpoonleft}{\isacharparenleft}{\kern0pt}\isactrlbold n{\isasymlbrace}\isabound{y}{\isasymrbrace}{\isacharparenright}{\kern0pt}{\isacharparenright}{\kern0pt}{\isacharparenright}{\kern0pt}{\isacharparenright}{\kern0pt}\ {\isasymlongrightarrow}\ \isabound{x}\ {\isasymtriangleright}\ {\isasymemptyset}\ {\isasymturnstile}\ \isactrlbold {\isasymnot}\isactrlbold A\ \isactrlbold {\isasymcdot}\ {\isacharparenleft}{\kern0pt}\isactrlbold {\isasymlambda}{\isacharparenleft}{\kern0pt}\isactrlbold {\isasymnot}\isafree{chk}\ \isactrlbold v\isactrlsub {\isadigit{0}}\ {\isacharparenleft}{\kern0pt}{\isasymupharpoonleft}{\isacharparenleft}{\kern0pt}\isactrlbold n{\isasymlbrace}\isactrlbold {\isasymnot}\isabound{y}{\isasymrbrace}{\isacharparenright}{\kern0pt}{\isacharparenright}{\kern0pt}{\isacharparenright}{\kern0pt}{\isacharparenright}{\kern0pt}\ {\isasymlongrightarrow}\ \isabound{x}\ {\isasymtriangleright}\ {\isasymemptyset}\ {\isasymturnstile}\ \isactrlbold {\isasymnot}\isactrlbold A\ \isactrlbold {\isasymcdot}\ {\isacharparenleft}{\kern0pt}\isactrlbold {\isasymlambda}{\isacharparenleft}{\kern0pt}\isactrlbold {\isasymnot}\isafree{chk}\ \isactrlbold v\isactrlsub {\isadigit{0}}\ {\isacharparenleft}{\kern0pt}{\isasymupharpoonleft}{\isacharparenleft}{\kern0pt}\isactrlbold n{\isasymlbrace}\isactrlbold {\isasymnot}\isactrlbold {\isadigit{0}}\ \isactrlbold {\isacharequal}{\kern0pt}\ \isactrlbold {\isadigit{0}}{\isasymrbrace}{\isacharparenright}{\kern0pt}{\isacharparenright}{\kern0pt}{\isacharparenright}{\kern0pt}{\isacharparenright}{\kern0pt}{\isacharparenright}{\kern0pt}\ {\isasymand}\ {\isacharparenleft}{\kern0pt}{\isasymforall}\isabound{x}{\isachardot}{\kern0pt}\ \isactrlbold {\isasymnot}\isactrlbold A\ \isactrlbold {\isasymcdot}\ {\isacharparenleft}{\kern0pt}\isactrlbold {\isasymlambda}{\isacharparenleft}{\kern0pt}\isactrlbold {\isasymnot}\isafree{chk}\ \isactrlbold v\isactrlsub {\isadigit{0}}\ {\isacharparenleft}{\kern0pt}{\isasymupharpoonleft}{\isacharparenleft}{\kern0pt}\isactrlbold n{\isasymlbrace}\isactrlbold {\isasymnot}\isabound{x}{\isasymrbrace}{\isacharparenright}{\kern0pt}{\isacharparenright}{\kern0pt}{\isacharparenright}{\kern0pt}{\isacharparenright}{\kern0pt}\ {\isasymtriangleright}\ {\isasymemptyset}\ {\isasymturnstile}\ \isactrlbold {\isasymnot}\isactrlbold A\ \isactrlbold {\isasymcdot}\ {\isacharparenleft}{\kern0pt}\isactrlbold {\isasymlambda}{\isacharparenleft}{\kern0pt}\isactrlbold {\isasymnot}\isafree{chk}\ \isactrlbold v\isactrlsub {\isadigit{0}}\ {\isacharparenleft}{\kern0pt}{\isasymupharpoonleft}{\isacharparenleft}{\kern0pt}\isactrlbold n{\isasymlbrace}\isactrlbold {\isasymnot}\isactrlbold A\ \isactrlbold {\isasymcdot}\ {\isacharparenleft}{\kern0pt}\isactrlbold {\isasymlambda}{\isacharparenleft}{\kern0pt}\isactrlbold {\isasymnot}\isafree{chk}\ \isactrlbold v\isactrlsub {\isadigit{0}}\ {\isacharparenleft}{\kern0pt}{\isasymupharpoonleft}{\isacharparenleft}{\kern0pt}\isactrlbold n{\isasymlbrace}\isactrlbold {\isasymnot}\isabound{x}{\isasymrbrace}{\isacharparenright}{\kern0pt}{\isacharparenright}{\kern0pt}{\isacharparenright}{\kern0pt}{\isacharparenright}{\kern0pt}{\isasymrbrace}{\isacharparenright}{\kern0pt}{\isacharparenright}{\kern0pt}{\isacharparenright}{\kern0pt}{\isacharparenright}{\kern0pt}{\isacharparenright}{\kern0pt}\ {\isasymand}\ {\isacharparenleft}{\kern0pt}{\isasymforall}\isabound{x}{\isachardot}{\kern0pt}\ \isaconst{wft}\ \isabound{x}\ {\isasymlongrightarrow}\ \isabound{x}\ {\isasymtriangleright}\ {\isasymemptyset}\ {\isasymturnstile}\ \isactrlbold {\isadigit{0}}\ \isactrlbold {\isacharequal}{\kern0pt}\ \isactrlbold {\isadigit{0}}{\isacharparenright}{\kern0pt}}%
}%
\SNIP{Goedel-Paper-Statements:context:g2-exf-con-0}{%
\isakeywordONE{context}\isamarkupfalse%
\ g{\isadigit{2}}{\isacharunderscore}{\kern0pt}exf{\isacharunderscore}{\kern0pt}con%
}%
\SNIP{Goedel-Paper-Statements:thm:no-code-unprovableT-0}{%
\isakeywordONE{thm}\isamarkupfalse%
\ no{\isacharunderscore}{\kern0pt}code{\isacharunderscore}{\kern0pt}unprovableT%
\isanewline%
}%
\SNIP{Goedel-Paper-Statements:thm:no-code-unprovableT-1}{%
\isa{\isaconst{wft}\ \isavar{{\isasymphi}}\ {\isasymLongrightarrow}\ {\isasymnot}\ {\isasymemptyset}\ {\isasymturnstile}\ \isactrlbold A\ \isactrlbold {\isasymcdot}\ {\isacharparenleft}{\kern0pt}\isactrlbold {\isasymlambda}{\isacharparenleft}{\kern0pt}\isactrlbold {\isasymnot}\isafree{chk}\ \isactrlbold v\isactrlsub {\isadigit{0}}\ {\isacharparenleft}{\kern0pt}\isactrlbold n{\isasymlbrace}\isavar{{\isasymphi}}{\isasymrbrace}{\isacharparenright}{\kern0pt}{\isacharparenright}{\kern0pt}{\isacharparenright}{\kern0pt}}%
}%
\SNIP{Goedel-Paper-Statements:thm:bsT-goedel2-irrefutable-0}{%
\isakeywordONE{thm}\isamarkupfalse%
\ bsT{\isacharunderscore}{\kern0pt}goedel{\isadigit{2}}{\isacharunderscore}{\kern0pt}irrefutable%
\isanewline%
}%
\SNIP{Goedel-Paper-Statements:thm:bsT-goedel2-irrefutable-1}{%
\isa{{\isasymnot}\ \isactrlbold {\isasymnot}\isactrlbold A\ \isactrlbold {\isasymcdot}\ {\isacharparenleft}{\kern0pt}\isactrlbold {\isasymlambda}{\isacharparenleft}{\kern0pt}\isactrlbold {\isasymnot}\isafree{chk}\ \isactrlbold v\isactrlsub {\isadigit{0}}\ {\isacharparenleft}{\kern0pt}{\isasymupharpoonleft}{\isacharparenleft}{\kern0pt}\isactrlbold n{\isasymlbrace}\isactrlbold {\isasymnot}\isactrlbold {\isadigit{0}}\ \isactrlbold {\isacharequal}{\kern0pt}\ \isactrlbold {\isadigit{0}}{\isasymrbrace}{\isacharparenright}{\kern0pt}{\isacharparenright}{\kern0pt}{\isacharparenright}{\kern0pt}{\isacharparenright}{\kern0pt}\ {\isasymtriangleright}\ {\isasymemptyset}\ {\isasymturnstile}\ \isactrlbold {\isasymnot}\isactrlbold {\isadigit{0}}\ \isactrlbold {\isacharequal}{\kern0pt}\ \isactrlbold {\isadigit{0}}}%
}%
\SNIP{Goedel-Paper-Statements:context:gi-0}{%
\isakeywordONE{context}\isamarkupfalse%
\ gi%
}%
\SNIP{Goedel-Paper-Statements:thm:gig-prf-strip-0}{%
\isakeywordONE{thm}\isamarkupfalse%
\ gig{\isacharunderscore}{\kern0pt}prf{\isacharunderscore}{\kern0pt}strip%
\isanewline%
}%
\SNIP{Goedel-Paper-Statements:thm:gig-prf-strip-1}{%
\isa{{\isasymlbrakk}\isaconst{valid}\ \isavar{P}{\isacharsemicolon}{\kern0pt}\ \isavar{P}\ {\isasymC}{\isasymlbrakk}\isavar{H}\ {\isasymturnstile}\ \isavar{c}{\isasymrbrakk}{\isasymrbrakk}\ {\isasymLongrightarrow}\ \isaconst{ffilter}\ \isaconst{wft}\ \isavar{H}\ {\isasymturnstile}\ \isavar{c}}%
}%
\SNIP{Goedel-Paper-Statements:thm:gi-ConRGA-of-gentzen-0}{%
\isakeywordONE{thm}\isamarkupfalse%
\ gi{\isacharunderscore}{\kern0pt}ConRGA{\isacharunderscore}{\kern0pt}of{\isacharunderscore}{\kern0pt}gentzen%
\isanewline%
}%
\SNIP{Goedel-Paper-Statements:thm:gi-ConRGA-of-gentzen-1}{%
\isa{{\isasymlbrakk}\isaconst{rga{\isachardot}{\kern0pt}GRed}\ \isavar{r}\ \isavar{od}\ \isavar{lt}{\isacharsemicolon}{\kern0pt}\ \isaconst{rga{\isachardot}{\kern0pt}gzIter}\ \isavar{F}\ \isavar{r}{\isacharsemicolon}{\kern0pt}\ \isaconst{prf{\isacharunderscore}{\kern0pt}derives{\isachardot}{\kern0pt}entails}\ \isaconst{rga{\isachardot}{\kern0pt}derives}\ {\isasymemptyset}\ {\isacharparenleft}{\kern0pt}\isactrlbold A\ \isactrlbold {\isasymcdot}\ {\isacharparenleft}{\kern0pt}\isactrlbold {\isasymlambda}{\isacharparenleft}{\kern0pt}\isavar{od}\ \isactrlbold {\isasymcdot}\ \isactrlbold v\isactrlsub {\isadigit{0}}\ \isactrlbold N{\isacharparenright}{\kern0pt}{\isacharparenright}{\kern0pt}{\isacharparenright}{\kern0pt}{\isacharsemicolon}{\kern0pt}\ \isaconst{prf{\isacharunderscore}{\kern0pt}derives{\isachardot}{\kern0pt}entails}\ \isaconst{rga{\isachardot}{\kern0pt}derives}\ {\isasymemptyset}\ {\isacharparenleft}{\kern0pt}\isactrlbold A\ \isactrlbold {\isasymcdot}\ {\isacharparenleft}{\kern0pt}\isactrlbold {\isasymlambda}\isactrlbold A\ \isactrlbold {\isasymcdot}\ {\isacharparenleft}{\kern0pt}\isactrlbold {\isasymlambda}{\isacharparenleft}{\kern0pt}\isactrlbold {\isasymnot}\isavar{lt}\ \isactrlbold {\isasymcdot}\ \isaconst{trmVar}\ {\isadigit{1}}\ \isactrlbold {\isasymcdot}\ \isactrlbold v\isactrlsub {\isadigit{0}}\ \isactrlbold {\isacharequal}{\kern0pt}\ \isactrlbold {\isadigit{0}}\ \isactrlbold {\isasymor}\ \isaconst{rga{\isachardot}{\kern0pt}vOp}\ \isactrlbold {\isasymcdot}\ \isaconst{trmVar}\ {\isadigit{1}}\ \isactrlbold {\isasymcdot}\ \isactrlbold v\isactrlsub {\isadigit{0}}{\isacharparenright}{\kern0pt}{\isacharparenright}{\kern0pt}{\isacharparenright}{\kern0pt}{\isacharparenright}{\kern0pt}{\isacharsemicolon}{\kern0pt}\ \isaconst{prf{\isacharunderscore}{\kern0pt}derives{\isachardot}{\kern0pt}entails}\ \isaconst{rga{\isachardot}{\kern0pt}derives}\ {\isasymemptyset}\ {\isacharparenleft}{\kern0pt}\isactrlbold A\ \isactrlbold {\isasymcdot}\ {\isacharparenleft}{\kern0pt}\isactrlbold {\isasymlambda}\isaconst{rga{\isachardot}{\kern0pt}vOp}\ \isactrlbold {\isasymcdot}\ {\isacharparenleft}{\kern0pt}\isavar{od}\ \isactrlbold {\isasymcdot}\ \isactrlbold v\isactrlsub {\isadigit{0}}{\isacharparenright}{\kern0pt}\ \isactrlbold {\isasymcdot}\ {\isacharparenleft}{\kern0pt}\isactrlbold n\isaconst{veb{\isacharunderscore}{\kern0pt}code}\ \isavar{{\isasymalpha}}{\isacharparenright}{\kern0pt}{\isacharparenright}{\kern0pt}{\isacharparenright}{\kern0pt}{\isacharsemicolon}{\kern0pt}\ \isaconst{veb{\isacharunderscore}{\kern0pt}valid}\ \isavar{{\isasymalpha}}{\isacharsemicolon}{\kern0pt}\ \isaconst{veb{\isacharunderscore}{\kern0pt}lt}\ \isavar{{\isasymalpha}}\ \isaconst{rga{\isachardot}{\kern0pt}veb{\isacharunderscore}{\kern0pt}pw}{\isasymrbrakk}\ {\isasymLongrightarrow}\ {\isasymemptyset}\ {\isasymturnstile}\ \isaconst{rga{\isachardot}{\kern0pt}ConRGA}}%
}%
\SNIP{Goedel-Paper-Statements:thm:gi-GTI-separation-0}{%
\isakeywordONE{thm}\isamarkupfalse%
\ gi{\isacharunderscore}{\kern0pt}GTI{\isacharunderscore}{\kern0pt}separation%
\isanewline%
}%
\SNIP{Goedel-Paper-Statements:thm:gi-GTI-separation-1}{%
\isa{{\isasymlbrakk}\isaconst{rga{\isachardot}{\kern0pt}GRed}\ \isavar{r}\ \isavar{od}\ \isavar{lt}{\isacharsemicolon}{\kern0pt}\ \isaconst{rga{\isachardot}{\kern0pt}gzIter}\ \isavar{F}\ \isavar{r}{\isacharsemicolon}{\kern0pt}\ \isaconst{prf{\isacharunderscore}{\kern0pt}derives{\isachardot}{\kern0pt}entails}\ \isaconst{rga{\isachardot}{\kern0pt}derives}\ {\isasymemptyset}\ {\isacharparenleft}{\kern0pt}\isactrlbold A\ \isactrlbold {\isasymcdot}\ {\isacharparenleft}{\kern0pt}\isactrlbold {\isasymlambda}{\isacharparenleft}{\kern0pt}\isavar{od}\ \isactrlbold {\isasymcdot}\ \isactrlbold v\isactrlsub {\isadigit{0}}\ \isactrlbold N{\isacharparenright}{\kern0pt}{\isacharparenright}{\kern0pt}{\isacharparenright}{\kern0pt}{\isacharsemicolon}{\kern0pt}\ \isaconst{prf{\isacharunderscore}{\kern0pt}derives{\isachardot}{\kern0pt}entails}\ \isaconst{rga{\isachardot}{\kern0pt}derives}\ {\isasymemptyset}\ {\isacharparenleft}{\kern0pt}\isactrlbold A\ \isactrlbold {\isasymcdot}\ {\isacharparenleft}{\kern0pt}\isactrlbold {\isasymlambda}\isactrlbold A\ \isactrlbold {\isasymcdot}\ {\isacharparenleft}{\kern0pt}\isactrlbold {\isasymlambda}{\isacharparenleft}{\kern0pt}\isactrlbold {\isasymnot}\isavar{lt}\ \isactrlbold {\isasymcdot}\ \isaconst{trmVar}\ {\isadigit{1}}\ \isactrlbold {\isasymcdot}\ \isactrlbold v\isactrlsub {\isadigit{0}}\ \isactrlbold {\isacharequal}{\kern0pt}\ \isactrlbold {\isadigit{0}}\ \isactrlbold {\isasymor}\ \isaconst{rga{\isachardot}{\kern0pt}vOp}\ \isactrlbold {\isasymcdot}\ \isaconst{trmVar}\ {\isadigit{1}}\ \isactrlbold {\isasymcdot}\ \isactrlbold v\isactrlsub {\isadigit{0}}{\isacharparenright}{\kern0pt}{\isacharparenright}{\kern0pt}{\isacharparenright}{\kern0pt}{\isacharparenright}{\kern0pt}{\isacharsemicolon}{\kern0pt}\ {\isasymemptyset}\ {\isasymturnstile}\ \isactrlbold A\ \isactrlbold {\isasymcdot}\ {\isacharparenleft}{\kern0pt}\isactrlbold {\isasymlambda}\isaconst{rga{\isachardot}{\kern0pt}vOp}\ \isactrlbold {\isasymcdot}\ {\isacharparenleft}{\kern0pt}\isavar{od}\ \isactrlbold {\isasymcdot}\ \isactrlbold v\isactrlsub {\isadigit{0}}{\isacharparenright}{\kern0pt}\ \isactrlbold {\isasymcdot}\ {\isacharparenleft}{\kern0pt}\isactrlbold n\isaconst{veb{\isacharunderscore}{\kern0pt}code}\ \isavar{{\isasymalpha}}{\isacharparenright}{\kern0pt}{\isacharparenright}{\kern0pt}{\isacharsemicolon}{\kern0pt}\ \isaconst{veb{\isacharunderscore}{\kern0pt}valid}\ \isavar{{\isasymalpha}}{\isacharsemicolon}{\kern0pt}\ \isaconst{veb{\isacharunderscore}{\kern0pt}lt}\ \isavar{{\isasymalpha}}\ \isaconst{rga{\isachardot}{\kern0pt}veb{\isacharunderscore}{\kern0pt}pw}{\isasymrbrakk}\ {\isasymLongrightarrow}\ {\isasymemptyset}\ {\isasymturnstile}\ \isaconst{rga{\isachardot}{\kern0pt}gzTIF}\ \isavar{F}\ \isavar{od}\ \isavar{lt}\isasep\isanewline%
}%
\SNIP{Goedel-Paper-Statements:thm:gi-GTI-separation-2}{%
{\isasymlbrakk}\isaconst{rga{\isachardot}{\kern0pt}GRed}\ \isavar{r}\ \isavar{od}\ \isavar{lt}{\isacharsemicolon}{\kern0pt}\ \isaconst{rga{\isachardot}{\kern0pt}gzIter}\ \isavar{F}\ \isavar{r}{\isacharsemicolon}{\kern0pt}\ \isaconst{prf{\isacharunderscore}{\kern0pt}derives{\isachardot}{\kern0pt}entails}\ \isaconst{rga{\isachardot}{\kern0pt}derives}\ {\isasymemptyset}\ {\isacharparenleft}{\kern0pt}\isactrlbold A\ \isactrlbold {\isasymcdot}\ {\isacharparenleft}{\kern0pt}\isactrlbold {\isasymlambda}{\isacharparenleft}{\kern0pt}\isavar{od}\ \isactrlbold {\isasymcdot}\ \isactrlbold v\isactrlsub {\isadigit{0}}\ \isactrlbold N{\isacharparenright}{\kern0pt}{\isacharparenright}{\kern0pt}{\isacharparenright}{\kern0pt}{\isacharsemicolon}{\kern0pt}\ \isaconst{prf{\isacharunderscore}{\kern0pt}derives{\isachardot}{\kern0pt}entails}\ \isaconst{rga{\isachardot}{\kern0pt}derives}\ {\isasymemptyset}\ {\isacharparenleft}{\kern0pt}\isactrlbold A\ \isactrlbold {\isasymcdot}\ {\isacharparenleft}{\kern0pt}\isactrlbold {\isasymlambda}\isactrlbold A\ \isactrlbold {\isasymcdot}\ {\isacharparenleft}{\kern0pt}\isactrlbold {\isasymlambda}{\isacharparenleft}{\kern0pt}\isactrlbold {\isasymnot}\isavar{lt}\ \isactrlbold {\isasymcdot}\ \isaconst{trmVar}\ {\isadigit{1}}\ \isactrlbold {\isasymcdot}\ \isactrlbold v\isactrlsub {\isadigit{0}}\ \isactrlbold {\isacharequal}{\kern0pt}\ \isactrlbold {\isadigit{0}}\ \isactrlbold {\isasymor}\ \isaconst{rga{\isachardot}{\kern0pt}vOp}\ \isactrlbold {\isasymcdot}\ \isaconst{trmVar}\ {\isadigit{1}}\ \isactrlbold {\isasymcdot}\ \isactrlbold v\isactrlsub {\isadigit{0}}{\isacharparenright}{\kern0pt}{\isacharparenright}{\kern0pt}{\isacharparenright}{\kern0pt}{\isacharparenright}{\kern0pt}{\isacharsemicolon}{\kern0pt}\ {\isasymemptyset}\ {\isasymturnstile}\ \isactrlbold A\ \isactrlbold {\isasymcdot}\ {\isacharparenleft}{\kern0pt}\isactrlbold {\isasymlambda}\isaconst{rga{\isachardot}{\kern0pt}vOp}\ \isactrlbold {\isasymcdot}\ {\isacharparenleft}{\kern0pt}\isavar{od}\ \isactrlbold {\isasymcdot}\ \isactrlbold v\isactrlsub {\isadigit{0}}{\isacharparenright}{\kern0pt}\ \isactrlbold {\isasymcdot}\ {\isacharparenleft}{\kern0pt}\isactrlbold n\isaconst{veb{\isacharunderscore}{\kern0pt}code}\ \isavar{{\isasymalpha}}{\isacharparenright}{\kern0pt}{\isacharparenright}{\kern0pt}{\isacharsemicolon}{\kern0pt}\ \isaconst{veb{\isacharunderscore}{\kern0pt}valid}\ \isavar{{\isasymalpha}}{\isacharsemicolon}{\kern0pt}\ \isaconst{veb{\isacharunderscore}{\kern0pt}lt}\ \isavar{{\isasymalpha}}\ \isaconst{rga{\isachardot}{\kern0pt}veb{\isacharunderscore}{\kern0pt}pw}{\isasymrbrakk}\ {\isasymLongrightarrow}\ {\isasymnot}\ \isaconst{prf{\isacharunderscore}{\kern0pt}derives{\isachardot}{\kern0pt}entails}\ \isaconst{rga{\isachardot}{\kern0pt}derives}\ {\isasymemptyset}\ {\isacharparenleft}{\kern0pt}\isaconst{rga{\isachardot}{\kern0pt}gzTIF}\ \isavar{F}\ \isavar{od}\ \isavar{lt}{\isacharparenright}{\kern0pt}}%
}%
\SNIP{Goedel-Paper-Statements:context:rga+1-0}{%
\isakeywordONE{context}\isamarkupfalse%
\ rga%
}%
\SNIP{Goedel-Paper-Statements:thm:gi-gwosr-open-iff-0}{%
\isakeywordONE{thm}\isamarkupfalse%
\ gi{\isacharunderscore}{\kern0pt}gwosr{\isacharunderscore}{\kern0pt}open{\isacharunderscore}{\kern0pt}iff%
\isanewline%
}%
\SNIP{Goedel-Paper-Statements:thm:gi-gwosr-open-iff-1}{%
\isa{\isaconst{wft}\ \isavar{c}\ {\isasymLongrightarrow}\ \isaconst{prf{\isacharunderscore}{\kern0pt}derives{\isachardot}{\kern0pt}entails}\ \isaconst{gi{\isachardot}{\kern0pt}derives}\ {\isasymemptyset}\ \isavar{c}\ {\isacharequal}{\kern0pt}\ {\isacharparenleft}{\kern0pt}{\isasymexists}\isabound{s}{\isachardot}{\kern0pt}\ \isaconst{gsem{\isachardot}{\kern0pt}gwosr}\ \isaconst{gi{\isachardot}{\kern0pt}derives}\ \isabound{s}\ \isavar{c}\ {\isasymlbrace}\isaconst{True}{\isasymrbrace}{\isacharparenright}{\kern0pt}}%
}%
\SNIP{Goedel-Paper-Statements:thm:gi-consistent-gsem-0}{%
\isakeywordONE{thm}\isamarkupfalse%
\ gi{\isacharunderscore}{\kern0pt}consistent{\isacharunderscore}{\kern0pt}gsem%
\isanewline%
}%
\SNIP{Goedel-Paper-Statements:thm:gi-consistent-gsem-1}{%
\isa{{\isasymnot}\ \isaconst{prf{\isacharunderscore}{\kern0pt}derives{\isachardot}{\kern0pt}entails}\ \isaconst{gi{\isachardot}{\kern0pt}derives}\ {\isasymemptyset}\ {\isacharparenleft}{\kern0pt}\isactrlbold {\isadigit{0}}\ \isactrlbold {\isacharequal}{\kern0pt}\ \isactrlbold {\isadigit{1}}{\isacharparenright}{\kern0pt}}%
}%
\SNIP{Goedel-Paper-Statements:thm:gi-disj-prop-0}{%
\isakeywordONE{thm}\isamarkupfalse%
\ gi{\isacharunderscore}{\kern0pt}disj{\isacharunderscore}{\kern0pt}prop%
\isanewline%
}%
\SNIP{Goedel-Paper-Statements:thm:gi-disj-prop-1}{%
\isa{\isaconst{prf{\isacharunderscore}{\kern0pt}derives{\isachardot}{\kern0pt}entails}\ \isaconst{gi{\isachardot}{\kern0pt}derives}\ {\isasymemptyset}\ {\isacharparenleft}{\kern0pt}\isavar{a}\ \isactrlbold {\isasymor}\ \isavar{b}{\isacharparenright}{\kern0pt}\ {\isasymLongrightarrow}\ \isaconst{prf{\isacharunderscore}{\kern0pt}derives{\isachardot}{\kern0pt}entails}\ \isaconst{gi{\isachardot}{\kern0pt}derives}\ {\isasymemptyset}\ \isavar{a}\ {\isasymor}\ \isaconst{prf{\isacharunderscore}{\kern0pt}derives{\isachardot}{\kern0pt}entails}\ \isaconst{gi{\isachardot}{\kern0pt}derives}\ {\isasymemptyset}\ \isavar{b}}%
}%
\SNIP{Goedel-Paper-Statements:thm:gi-not-all-witness-0}{%
\isakeywordONE{thm}\isamarkupfalse%
\ gi{\isacharunderscore}{\kern0pt}not{\isacharunderscore}{\kern0pt}all{\isacharunderscore}{\kern0pt}witness%
\isanewline%
}%
\SNIP{Goedel-Paper-Statements:thm:gi-not-all-witness-1}{%
\isa{\isaconst{prf{\isacharunderscore}{\kern0pt}derives{\isachardot}{\kern0pt}entails}\ \isaconst{gi{\isachardot}{\kern0pt}derives}\ {\isasymemptyset}\ {\isacharparenleft}{\kern0pt}\isactrlbold {\isasymnot}\isactrlbold A\ \isactrlbold {\isasymcdot}\ \isavar{u}{\isacharparenright}{\kern0pt}\ {\isasymLongrightarrow}\ {\isasymexists}\isabound{n}{\isachardot}{\kern0pt}\ \isaconst{prf{\isacharunderscore}{\kern0pt}derives{\isachardot}{\kern0pt}entails}\ \isaconst{gi{\isachardot}{\kern0pt}derives}\ {\isasymemptyset}\ {\isacharparenleft}{\kern0pt}\isactrlbold {\isasymnot}\isavar{u}\ \isactrlbold {\isasymcdot}\ {\isacharparenleft}{\kern0pt}\isactrlbold n\isabound{n}{\isacharparenright}{\kern0pt}{\isacharparenright}{\kern0pt}}%
}%
\SNIP{Goedel-Paper-Statements:thm:gi-entails-re-0}{%
\isakeywordONE{thm}\isamarkupfalse%
\ gi{\isacharunderscore}{\kern0pt}entails{\isacharunderscore}{\kern0pt}re%
\isanewline%
}%
\SNIP{Goedel-Paper-Statements:thm:gi-entails-re-1}{%
\isa{\isaconst{renp}\ {\isacharparenleft}{\kern0pt}{\isasymlambda}\isabound{z}{\isachardot}{\kern0pt}\ \isaconst{prf{\isacharunderscore}{\kern0pt}derives{\isachardot}{\kern0pt}entails}\ \isaconst{gi{\isachardot}{\kern0pt}derives}\ {\isacharparenleft}{\kern0pt}\isaconst{decode}\ {\isacharparenleft}{\kern0pt}\isaconst{cpx}\ \isabound{z}{\isacharparenright}{\kern0pt}{\isacharparenright}{\kern0pt}\ {\isacharparenleft}{\kern0pt}\isaconst{decode}\ {\isacharparenleft}{\kern0pt}\isaconst{cpy}\ \isabound{z}{\isacharparenright}{\kern0pt}{\isacharparenright}{\kern0pt}{\isacharparenright}{\kern0pt}}%
}%
\SNIP{Goedel-Paper-Statements:context:gi+1-0}{%
\isakeywordONE{context}\isamarkupfalse%
\ gi%
}%
\SNIP{Goedel-Paper-Statements:thm:gi-halting-undecided-0}{%
\isakeywordONE{thm}\isamarkupfalse%
\ gi{\isacharunderscore}{\kern0pt}halting{\isacharunderscore}{\kern0pt}undecided%
\isanewline%
}%
\SNIP{Goedel-Paper-Statements:thm:gi-halting-undecided-1}{%
\isa{{\isasymexists}\isabound{x}{\isachardot}{\kern0pt}\ \isabound{x}\ {\isasymin}\ {\isacharminus}{\kern0pt}\ \isaconst{diagK}\ {\isasymand}\ {\isasymnot}\ {\isasymemptyset}\ {\isasymturnstile}\ \isactrlbold A\ \isactrlbold {\isasymcdot}\ \isaconst{rga{\isachardot}{\kern0pt}ofam}\ \isabound{x}\ \isactrlbold B}%
}%
\SNIP{Goedel-Paper-Statements:context:rga+2-0}{%
\isakeywordONE{context}\isamarkupfalse%
\ rga%
}%
\SNIP{Goedel-Paper-Statements:thm:gi-famT-undecided-0}{%
\isakeywordONE{thm}\isamarkupfalse%
\ gi{\isacharunderscore}{\kern0pt}famT{\isacharunderscore}{\kern0pt}undecided%
\isanewline%
}%
\SNIP{Goedel-Paper-Statements:thm:gi-famT-undecided-1}{%
\isa{{\isasymexists}\isabound{x}{\isachardot}{\kern0pt}\ \isabound{x}\ {\isasymin}\ {\isacharminus}{\kern0pt}\ \isaconst{diagK}\ {\isasymand}\ {\isasymnot}\ \isaconst{prf{\isacharunderscore}{\kern0pt}derives{\isachardot}{\kern0pt}entails}\ \isaconst{gi{\isachardot}{\kern0pt}derives}\ {\isasymemptyset}\ {\isacharparenleft}{\kern0pt}\isaconst{allq}\ {\isacharparenleft}{\kern0pt}\isaconst{RF}\ {\isacharparenleft}{\kern0pt}\isactrlbold {\isasymlambda}{\isacharparenleft}{\kern0pt}\isactrlbold {\isasymlambda}{\isacharparenleft}{\kern0pt}\isactrlbold {\isasymlambda}\isaconst{recBody}{\isacharparenright}{\kern0pt}{\isacharparenright}{\kern0pt}{\isacharparenright}{\kern0pt}\ {\isacharparenleft}{\kern0pt}\isactrlbold n\isaconst{Suc}\ {\isadigit{0}}{\isacharparenright}{\kern0pt}\ {\isacharparenleft}{\kern0pt}\isaconst{mkAppS}\ {\isacharparenleft}{\kern0pt}{\isasymupharpoonleft}{\isacharparenleft}{\kern0pt}\isactrlbold n{\isasymlbrace}\isaconst{ofam}\ \isabound{x}{\isasymrbrace}{\isacharparenright}{\kern0pt}{\isacharparenright}{\kern0pt}\ {\isacharparenleft}{\kern0pt}\isaconst{cmp}\ \isaconst{num{\isacharunderscore}{\kern0pt}codeIdx}\ \isactrlbold {\isasymcdot}\ \isactrlbold {\isadigit{0}}\ \isactrlbold {\isasymcdot}\ \isactrlbold v\isactrlsub {\isadigit{0}}{\isacharparenright}{\kern0pt}{\isacharparenright}{\kern0pt}{\isacharparenright}{\kern0pt}\ \isactrlbold B{\isacharparenright}{\kern0pt}}%
}%
\SNIP{Goedel-Paper-Statements:thm:gi-not-halting-biv-0}{%
\isakeywordONE{thm}\isamarkupfalse%
\ gi{\isacharunderscore}{\kern0pt}not{\isacharunderscore}{\kern0pt}halting{\isacharunderscore}{\kern0pt}biv%
\isanewline%
}%
\SNIP{Goedel-Paper-Statements:thm:gi-not-halting-biv-1}{%
\isa{{\isasymnot}\ \isaconst{prf{\isacharunderscore}{\kern0pt}derives{\isachardot}{\kern0pt}entails}\ \isaconst{gi{\isachardot}{\kern0pt}derives}\ {\isasymemptyset}\ {\isacharparenleft}{\kern0pt}\isactrlbold A\ \isactrlbold {\isasymcdot}\ {\isacharparenleft}{\kern0pt}\isactrlbold {\isasymlambda}{\isacharparenleft}{\kern0pt}\isaconst{allq}\ {\isacharparenleft}{\kern0pt}\isaconst{RF}\ {\isacharparenleft}{\kern0pt}\isactrlbold {\isasymlambda}{\isacharparenleft}{\kern0pt}\isactrlbold {\isasymlambda}{\isacharparenleft}{\kern0pt}\isactrlbold {\isasymlambda}\isaconst{recBody}{\isacharparenright}{\kern0pt}{\isacharparenright}{\kern0pt}{\isacharparenright}{\kern0pt}\ {\isacharparenleft}{\kern0pt}\isactrlbold n\isaconst{Suc}\ {\isadigit{0}}{\isacharparenright}{\kern0pt}\ {\isacharparenleft}{\kern0pt}\isaconst{mkAppS}\ {\isacharparenleft}{\kern0pt}{\isasymupharpoonleft}\isactrlbold v\isactrlsub {\isadigit{0}}{\isacharparenright}{\kern0pt}\ {\isacharparenleft}{\kern0pt}\isaconst{cmp}\ \isaconst{num{\isacharunderscore}{\kern0pt}codeIdx}\ \isactrlbold {\isasymcdot}\ \isactrlbold {\isadigit{0}}\ \isactrlbold {\isasymcdot}\ \isactrlbold v\isactrlsub {\isadigit{0}}{\isacharparenright}{\kern0pt}{\isacharparenright}{\kern0pt}{\isacharparenright}{\kern0pt}\ \isactrlbold B{\isacharparenright}{\kern0pt}{\isacharparenright}{\kern0pt}{\isacharparenright}{\kern0pt}}%
}%
\SNIP{Goedel-Paper-Statements:thm:gi-not-iLFI-0}{%
\isakeywordONE{thm}\isamarkupfalse%
\ gi{\isacharunderscore}{\kern0pt}not{\isacharunderscore}{\kern0pt}iLFI%
\isanewline%
}%
\SNIP{Goedel-Paper-Statements:thm:gi-not-iLFI-1}{%
\isa{{\isasymnot}\ \isaconst{prf{\isacharunderscore}{\kern0pt}derives{\isachardot}{\kern0pt}entails}\ \isaconst{gi{\isachardot}{\kern0pt}derives}\ {\isasymemptyset}\ \isaconst{iLFI}}%
}%
\SNIP{Goedel-Paper-Statements:thm:giv-gsem-sound-0}{%
\isakeywordONE{thm}\isamarkupfalse%
\ giv{\isacharunderscore}{\kern0pt}gsem{\isacharunderscore}{\kern0pt}sound%
\isanewline%
}%
\SNIP{Goedel-Paper-Statements:thm:giv-gsem-sound-1}{%
\isa{\isaconst{gsem{\isacharunderscore}{\kern0pt}sound}\ {\isacharparenleft}{\kern0pt}\isaconst{giv{\isacharunderscore}{\kern0pt}calc{\isachardot}{\kern0pt}givR}\ \isaconst{vOp}{\isacharparenright}{\kern0pt}\ {\isacharparenleft}{\kern0pt}{\isasymlambda}\isabound{P}\ \isabound{H}\ \isabound{c}{\isachardot}{\kern0pt}\ \isaconst{gfixStep}\ \isabound{P}\ \isabound{H}\ \isabound{c}\ {\isasymor}\ \isaconst{giv{\isacharunderscore}{\kern0pt}calc{\isachardot}{\kern0pt}accbivS}\ \isaconst{vOp}\ \isabound{P}\ \isabound{H}\ \isabound{c}{\isacharparenright}{\kern0pt}}%
}%
\SNIP{Goedel-Paper-Statements:thm:giv-disj-prop-0}{%
\isakeywordONE{thm}\isamarkupfalse%
\ giv{\isacharunderscore}{\kern0pt}disj{\isacharunderscore}{\kern0pt}prop%
\isanewline%
}%
\SNIP{Goedel-Paper-Statements:thm:giv-disj-prop-1}{%
\isa{\isaconst{prf{\isacharunderscore}{\kern0pt}derives{\isachardot}{\kern0pt}entails}\ {\isacharparenleft}{\kern0pt}\isaconst{giv{\isacharunderscore}{\kern0pt}calc{\isachardot}{\kern0pt}givR}\ \isaconst{vOp}{\isacharparenright}{\kern0pt}\ {\isasymemptyset}\ {\isacharparenleft}{\kern0pt}\isavar{a}\ \isactrlbold {\isasymor}\ \isavar{b}{\isacharparenright}{\kern0pt}\ {\isasymLongrightarrow}\ \isaconst{prf{\isacharunderscore}{\kern0pt}derives{\isachardot}{\kern0pt}entails}\ {\isacharparenleft}{\kern0pt}\isaconst{giv{\isacharunderscore}{\kern0pt}calc{\isachardot}{\kern0pt}givR}\ \isaconst{vOp}{\isacharparenright}{\kern0pt}\ {\isasymemptyset}\ \isavar{a}\ {\isasymor}\ \isaconst{prf{\isacharunderscore}{\kern0pt}derives{\isachardot}{\kern0pt}entails}\ {\isacharparenleft}{\kern0pt}\isaconst{giv{\isacharunderscore}{\kern0pt}calc{\isachardot}{\kern0pt}givR}\ \isaconst{vOp}{\isacharparenright}{\kern0pt}\ {\isasymemptyset}\ \isavar{b}}%
}%
\SNIP{Goedel-Paper-Statements:context:oiga-0}{%
\isakeywordONE{context}\isamarkupfalse%
\ oiga%
}%
\SNIP{Goedel-Paper-Statements:thm:iga-not-contains-oga-0}{%
\isakeywordONE{thm}\isamarkupfalse%
\ iga{\isacharunderscore}{\kern0pt}not{\isacharunderscore}{\kern0pt}contains{\isacharunderscore}{\kern0pt}oga%
\isanewline%
}%
\SNIP{Goedel-Paper-Statements:thm:iga-not-contains-oga-1}{%
\isa{{\isasymexists}\isabound{c}{\isachardot}{\kern0pt}\ \isaconst{prf{\isacharunderscore}{\kern0pt}derives{\isachardot}{\kern0pt}entails}\ \isaconst{oga{\isachardot}{\kern0pt}derives}\ {\isasymemptyset}\ \isabound{c}\ {\isasymand}\ {\isasymnot}\ \isaconst{prf{\isacharunderscore}{\kern0pt}derives{\isachardot}{\kern0pt}entails}\ \isaconst{gi{\isachardot}{\kern0pt}derives}\ {\isasymemptyset}\ \isabound{c}}%
}%
\SNIP{Goedel-Paper-Statements:thm:iga-not-contains-oiga-0}{%
\isakeywordONE{thm}\isamarkupfalse%
\ iga{\isacharunderscore}{\kern0pt}not{\isacharunderscore}{\kern0pt}contains{\isacharunderscore}{\kern0pt}oiga%
\isanewline%
}%
\SNIP{Goedel-Paper-Statements:thm:iga-not-contains-oiga-1}{%
\isa{{\isasymexists}\isabound{c}{\isachardot}{\kern0pt}\ \isaconst{oigaPR{\isachardot}{\kern0pt}entails}\ {\isasymemptyset}\ \isabound{c}\ {\isasymand}\ {\isasymnot}\ \isaconst{prf{\isacharunderscore}{\kern0pt}derives{\isachardot}{\kern0pt}entails}\ \isaconst{gi{\isachardot}{\kern0pt}derives}\ {\isasymemptyset}\ \isabound{c}}%
}%
\SNIP{Goedel-Paper-Statements:context:gi+2-0}{%
\isakeywordONE{context}\isamarkupfalse%
\ gi%
}%
\SNIP{Goedel-Paper-Statements:thm:gi-desc-pr-contra-0}{%
\isakeywordONE{thm}\isamarkupfalse%
\ gi{\isacharunderscore}{\kern0pt}desc{\isacharunderscore}{\kern0pt}pr{\isacharunderscore}{\kern0pt}contra%
\isanewline%
}%
\SNIP{Goedel-Paper-Statements:thm:gi-desc-pr-contra-1}{%
\isa{{\isasymlbrakk}\isaconst{Syntax{\isacharunderscore}{\kern0pt}Term{\isachardot}{\kern0pt}closed}\ \isavar{R}{\isacharsemicolon}{\kern0pt}\ \isaconst{wft}\ \isavar{R}{\isacharsemicolon}{\kern0pt}\ {\isasymemptyset}\ {\isasymturnstile}\ \isactrlbold A\ \isactrlbold {\isasymcdot}\ {\isacharparenleft}{\kern0pt}\isactrlbold {\isasymlambda}\isavar{R}\ \isactrlbold {\isasymcdot}\ {\isacharparenleft}{\kern0pt}\isaconst{rga{\isachardot}{\kern0pt}cmp}\ \isavar{f}\ \isactrlbold {\isasymcdot}\ \isactrlbold {\isadigit{0}}\ \isactrlbold {\isasymcdot}\ {\isacharparenleft}{\kern0pt}\isactrlbold S\ \isactrlbold v\isactrlsub {\isadigit{0}}{\isacharparenright}{\kern0pt}{\isacharparenright}{\kern0pt}\ \isactrlbold {\isasymcdot}\ {\isacharparenleft}{\kern0pt}\isaconst{rga{\isachardot}{\kern0pt}cmp}\ \isavar{f}\ \isactrlbold {\isasymcdot}\ \isactrlbold {\isadigit{0}}\ \isactrlbold {\isasymcdot}\ \isactrlbold v\isactrlsub {\isadigit{0}}{\isacharparenright}{\kern0pt}{\isacharparenright}{\kern0pt}{\isacharsemicolon}{\kern0pt}\ {\isasymemptyset}\ {\isasymturnstile}\ \isaconst{rga{\isachardot}{\kern0pt}trmAcc}\ \isavar{R}\ \isactrlbold {\isasymcdot}\ {\isacharparenleft}{\kern0pt}\isaconst{rga{\isachardot}{\kern0pt}cmp}\ \isavar{f}\ \isactrlbold {\isasymcdot}\ \isactrlbold {\isadigit{0}}\ \isactrlbold {\isasymcdot}\ \isactrlbold {\isadigit{0}}{\isacharparenright}{\kern0pt}{\isasymrbrakk}\ {\isasymLongrightarrow}\ {\isasymemptyset}\ {\isasymturnstile}\ \isactrlbold {\isadigit{0}}\ \isactrlbold {\isacharequal}{\kern0pt}\ \isactrlbold {\isadigit{1}}}%
}%
\SNIP{Goedel-Paper-Statements:thm:gi-no-desc-below-pw-0}{%
\isakeywordONE{thm}\isamarkupfalse%
\ gi{\isacharunderscore}{\kern0pt}no{\isacharunderscore}{\kern0pt}desc{\isacharunderscore}{\kern0pt}below{\isacharunderscore}{\kern0pt}pw%
\isanewline%
}%
\SNIP{Goedel-Paper-Statements:thm:gi-no-desc-below-pw-1}{%
\isa{{\isasymlbrakk}\isaconst{veb{\isacharunderscore}{\kern0pt}valid}\ \isavar{x}{\isacharsemicolon}{\kern0pt}\ \isaconst{veb{\isacharunderscore}{\kern0pt}lt}\ \isavar{x}\ \isaconst{rga{\isachardot}{\kern0pt}veb{\isacharunderscore}{\kern0pt}pw}{\isacharsemicolon}{\kern0pt}\ {\isasymemptyset}\ {\isasymturnstile}\ \isaconst{rga{\isachardot}{\kern0pt}cmp}\ \isavar{f}\ \isactrlbold {\isasymcdot}\ \isactrlbold {\isadigit{0}}\ \isactrlbold {\isasymcdot}\ \isactrlbold {\isadigit{0}}\ \isactrlbold {\isacharequal}{\kern0pt}\ \isactrlbold n\isaconst{veb{\isacharunderscore}{\kern0pt}code}\ \isavar{x}{\isasymrbrakk}\ {\isasymLongrightarrow}\ {\isasymnot}\ {\isasymemptyset}\ {\isasymturnstile}\ \isactrlbold A\ \isactrlbold {\isasymcdot}\ {\isacharparenleft}{\kern0pt}\isactrlbold {\isasymlambda}\isaconst{rga{\isachardot}{\kern0pt}vOp}\ \isactrlbold {\isasymcdot}\ {\isacharparenleft}{\kern0pt}\isaconst{rga{\isachardot}{\kern0pt}cmp}\ \isavar{f}\ \isactrlbold {\isasymcdot}\ \isactrlbold {\isadigit{0}}\ \isactrlbold {\isasymcdot}\ {\isacharparenleft}{\kern0pt}\isactrlbold S\ \isactrlbold v\isactrlsub {\isadigit{0}}{\isacharparenright}{\kern0pt}{\isacharparenright}{\kern0pt}\ \isactrlbold {\isasymcdot}\ {\isacharparenleft}{\kern0pt}\isaconst{rga{\isachardot}{\kern0pt}cmp}\ \isavar{f}\ \isactrlbold {\isasymcdot}\ \isactrlbold {\isadigit{0}}\ \isactrlbold {\isasymcdot}\ \isactrlbold v\isactrlsub {\isadigit{0}}{\isacharparenright}{\kern0pt}{\isacharparenright}{\kern0pt}}%
}%
}{}

\newcommand{\GetSnip}[1]{%
    \ifcsname snippet--#1\endcsname%
        \csname snippet--#1\endcsname%
    \else%
        \PackageWarning{snips}{Snippet ``#1'' is undefined.}%
        \emph{Warning: Snippet ``#1'' is undefined.}%
    \fi%
}

\newcommand{\Snippet}[1]{{%
  \newcount\i
  \i=0
  \loop
    \GetSnip{#1-\the\i}%
    \advance \i 1
  \ifcsname snippet--#1-\the\i\endcsname
  \repeat
}}

\newcommand{\SnippetPart}[3]{{%
  \newcount\i
  \i=#1
  \loop
    \ifnum \i=#2
      \renewcommand{\isanewline}{}%
    \fi
    \GetSnip{#3-\the\i}%
    \advance \i 1
    \ifnum \i>#2 {}
    \else \repeat
}}

\hypersetup{breaklinks=true}

\newcommand{\com}[1]{}

\long\def\note#1{}	%

\newcommand{\rga}{RGA\xspace}		%

\newcommand{\gdl}{G\"odel\xspace}

\newcommand{\limp}{\rightarrow}		%
\newcommand{\liff}{\leftrightarrow}	%
\newcommand{\ldef}{\equiv}		%

\newcommand{\typestyle}[1]{\textsf{#1}\xspace}
\com{	%
\newcommand{\tbool}{\typestyle{bool}}
\newcommand{\ttrue}{\typestyle{true}}
\newcommand{\tfalse}{\typestyle{false}}

\newcommand{\tnat}{\typestyle{nat}}

}%
\newcommand{\tbool}{\typestyle{B}}	%
\newcommand{\ttrue}{\typestyle{T}}
\newcommand{\tfalse}{\typestyle{F}}

\newcommand{\tnat}{\typestyle{N}}

\newcommand{\ctrue}{\ttrue}
\newcommand{\cfalse}{\tfalse}

\newcommand{\judgment}[1]{\textsf{ #1}}

\newcommand{\jbool}{\judgment{\tbool}}

\newcommand{\jnat}{\judgment{\tnat}}

\makeatletter
\newlength{\doublefracgap}
\DeclareRobustCommand{\doublefrac}[2]{%
  \mathinner{\mathpalette\doublefrac@{{#1}{#2}}}%
}
\newcommand{\doublefrac@}[2]{\doublefrac@@#1#2}
\newcommand{\doublefrac@@}[3]{%
  \ooalign{%
    \raisebox{\doublefracgap}{$\m@th#1\frac{#2}{\phantom{#3}}$}\cr
    \raisebox{-\doublefracgap}{$\m@th#1\frac{\phantom{#2}}{#3}$}\cr
  }%
}
\newcommand{\ddoublefrac}[2]{{\displaystyle\doublefrac{#1}{#2}}}

\makeatother

\newcommand{\irl}[1]{\ensuremath{\mathit{#1}}}		%

\newcommand{\infrule}[3][]{\cfrac{#2}{#3}\IfStrEq{#1}{}{}{\ \irl{#1}}}
\newcommand{\infeqv}[3][]{\ddoublefrac{#2}{#3}\IfStrEq{#1}{}{}{\ \irl{#1}}}
\newcommand{\infceqv}[4][]{\cfrac{#2\qquad}{}\ddoublefrac{#3}{#4}
				\IfStrEq{#1}{}{}{\ \irl{#1}}}

\newcommand{\z}{\mathbf 0}

\newcommand{\suc}{\mathbf S}

\newcommand{\cR}{\mathcal R}

\newcommand{\quo}[1]{\ulcorner{#1}\urcorner}

\newcommand{\tofn}[1]{\underline{#1}}

\SetSymbolFont{stmry}{bold}{U}{stmry}{m}{n}

\theoremstyle{definition}
\newtheorem{defi}[thm]{Definition}
\newtheorem{rem}[thm]{Remark}	%
\theoremstyle{plain}
\newcommand{\cmp}{\mathit{cmp}}		%
\newcommand{\den}[1]{\llbracket #1 \rrbracket}	%
\newcommand{\cA}{\mathcal A}		%
\newcommand{\lft}[1]{{\uparrow}#1}	%
\newcommand{\vz}{v_0}			%
\newcommand{\wf}{\mathrm{wf}}		%
\newcommand{\proves}{\mathrel{:}}	%
\newcommand{\cpair}[2]{\langle #1, #2 \rangle}	%

\newcommand{\iRec}{\mathsf{Rec}}	%
\newcommand{\iSub}{\mathsf{Sub}}	%
\newcommand{\sof}[2]{\hat{#1}({#2})}	%
\newcommand{\ssub}[3]{[{#1}{\mapsto}{#2}]{#3}}	%
\newcommand{\bquo}[1]{\underline{\vrule width0pt depth0.08em\relax\smash[b]{\quo{#1}}}}

\newcommand{\Dec}{\mathit{Dec}}		%
\newcommand{\Aff}{\mathit{Aff}}		%
\newcommand{\Tru}{\mathit{True}}	%
\newcommand{\Chk}{\mathit{Chk}}		%
\newcommand{\Cert}{\mathit{Cert}}	%
\newcommand{\Subst}{\mathit{Subst}}	%
\newcommand{\Num}{\mathit{Num}}		%
\newcommand{\Wfc}{\mathit{IsWf}}		%
\newcommand{\IsCode}{\mathit{IsCode}}	%
\newcommand{\Diag}{\mathit{Diag}}	%
\newcommand{\ConR}{\mathit{Con}}	%
\newcommand{\ChkQf}{\Chk_{\mathrm{qf}}}	%
\newcommand{\Rcd}{\mathit{Val}}		%
\newcommand{\Frag}{\mathit{Frag}}	%
\newcommand{\LvC}{\mathit{Lv}}		%
\newcommand{\LF}{\mathit{LF}}		%
\newcommand{\Res}{\mathit{Res}}		%
\newcommand{\oga}{OGA\xspace}		%
\newcommand{\eqv}{\mathrel{\dashv\vdash}}	%
\newcommand{\PrC}{\mathit{Pr}}		%

\newcommand{\iga}{IGA\xspace}			%
\newcommand{\oiga}{OIGA\xspace}			%
\newcommand{\rgaP}{\mbox{RGA$^{+}$}\xspace}	%
\newcommand{\igaP}{\mbox{IGA$^{+}$}\xspace}	%
\newcommand{\adoptR}{adoption\textsubscript{\rga}\xspace}	%
\newcommand{\adoptI}{adoption\textsubscript{\iga}\xspace}	%
\newcommand{\ChkI}{\Chk_{\mathrm{I}}}	%
\newcommand{\ConI}{\mathit{Con}_{\mathrm{IGA}}}	%
\newcommand{\PrI}{\PrC_{\mathrm{I}}}	%
\newcommand{\vdashI}{\vdash_{\mathrm{I}}}	%
\newcommand{\nvdashI}{\nvdash_{\mathrm{I}}}
\newcommand{\vdashRP}{\vdash_{\mathrm{R}^{+}}}	%
\newcommand{\nvdashRP}{\nvdash_{\mathrm{R}^{+}}}
\newcommand{\vdashIP}{\vdash_{\mathrm{I}^{+}}}	%
\newcommand{\nvdashIP}{\nvdash_{\mathrm{I}^{+}}}
\newcommand{\Inst}{\mathit{Inst}}	%
\newcommand{\WarR}{W_{\mathrm{R}}}	%
\newcommand{\WarI}{W_{\mathrm{I}}}	%
\newcommand{\Mw}{M_{\omega}}		%
\newcommand{\KO}{\mathcal O}			%
\newcommand{\Ia}[1]{\mbox{IGA$[#1]$}}		%
\newcommand{\vdashA}[1]{\vdash_{[#1]}}		%
\newcommand{\nvdashA}[1]{\nvdash_{[#1]}}
\newcommand{\ChkA}[1]{\Chk_{[#1]}}		%
\newcommand{\PrA}[1]{\PrC_{[#1]}}		%
\newcommand{\ConA}[1]{\mathit{Con}_{[#1]}}	%
\newcommand{\WF}{\mathrm{WF}}			%
\newcommand{\Acc}{\mathrm{Acc}}		%
\newcommand{\gchk}{\mathit{chk}}		%
\newcommand{\gPr}{\PrC_{\gchk}}		%
\newcommand{\gCon}{\ConR_{\gchk}}		%
\newcommand{\gG}{G_{\gchk}}			%
\newcommand{\gS}{S_{\gchk}}			%
\newcommand{\gD}{D_{\gchk}}			%

\begin{document}

\title{G\"odel's Incompleteness Theorems in Grounded Arithmetic}
\author{Bryan Ford}
\date{}
\maketitle

\begin{abstract}
Reflective Grounded Arithmetic (RGA) is a paracomplete arithmetic
in which every term is a program and truth is grounded evaluation.
RGA discharges a hypothesis only when it can prove the hypothesis decided.
This paper shows, in a machine-checked Isabelle/HOL development,
how \gdl's incompleteness theorems hold for RGA and for systems built on it.
RGA's implication is material and RGA has the disjunction property,
so the implications that the textbook proof of the second theorem uses fail,
while the entailments and rules it uses hold:
\gdl's proof survives as entailments.
RGA does not prove its consistency sentence $\ConR$,
and RGA derives its own second theorem
as the entailment $\ConR \vdash \neg\PrC(\bquo{\ConR})$
and as the sentence $\forall q.\,(\Chk(q, \bquo{\ConR}) \limp \neg\ConR)$.
The second theorem is then proved once,
for every \emph{G\"odel extension} of RGA:
an entailment relation on RGA's language, closed under RGA's rules,
with a checker that meets three structural conditions,
NEC, D3 and EXPL as proof transformers, D1 as closed acceptance,
and consistency.
D2 is not among the conditions.
RGA with its proof checker is one instance.
Inductive Grounded Arithmetic (IGA), RGA with ground induction, is another,
with a checker that adds one rule detector to RGA's;
IGA does not derive its consistency sentence for that checker.
Adoption, a reflection rule whose premise is a derived warrant
that the base system's checker accepts computed proofs
of every instance of a universal,
gives consistent extensions RGA$^{+}$ of RGA and IGA$^{+}$ of IGA.
RGA$^{+}$ derives RGA's consistency sentence,
and IGA$^{+}$ derives RGA's and IGA's,
so RGA$^{+}$ strictly extends RGA,
and IGA$^{+}$ strictly extends both RGA and IGA.
The second theorem holds at every adopted notation claim:
IGA with the claim that $a$ is a notation does not derive its own consistency sentence,
which adopting Turing's code built from that claim's checker does,
so no claim certifies itself and each needs a larger one.
Whether IGA alone derives RGA's consistency sentence is open,
and the paper names what each of two routes to an answer must supply.
\end{abstract}

\section{Introduction}
\label{sec:intro}

G\"odel's second incompleteness theorem is usually read as the point
where a formal system's knowledge of itself runs out: a consistent
system of sufficient strength cannot prove its own consistency.  The
theorem's textbook proof is valid even in intuitionistic arithmetic,
but the proof relies on two key elements.  The proof asserts a
self-referential sentence to be equivalent to that sentence's own
unprovability, and the proof derives implications by discharging
hypotheses, provability and consistency statements, that the system
can neither prove nor refute.  This paper asks what happens to
G\"odel's theorems, and to their proofs, in a system where neither
element is available and where everything the system can say is a
computation; and then whether the answer holds for one system only or
for a family of systems.

Reflective Grounded Arithmetic, or \rga~\cite{rga-paper}, is such a
system.  \rga's terms are programs: numerals, $\lambda$-terms, a
recursor, and a quantifier combinator $\cA$.  A term is true when its
evaluation yields $1$, with a proof certificate at each universal
quantifier, false when it yields $0$, and \emph{ungrounded} when it
yields neither value.  So \rga is paracomplete: excluded middle fails
at ungrounded terms.  \rga discharges a hypothesis only when \rga can
prove that hypothesis decided, whether the discharge is by implication
introduction or by case analysis.  The RGA paper proves in
Isabelle/HOL that \rga is sound and consistent, and complete for
grounded truth: for well-formed terms, provability coincides with
grounded truth~\cite{rga-paper}.  This completeness is not negation
completeness (every sentence provable or refutable).  The G\"odel
sentence and the consistency sentence are ungrounded, because \rga's
universal quantifier is settled only by a proof certificate, not
instance by instance, so completeness for grounded truth does not
reach them.  The internal-reasoning paper~\cite{conrga-paper} builds,
inside \rga, the machinery that the proofs here run on: codes of terms
and proofs, compiled primitive-recursive machines, and a proof checker
$\Chk$ whose verdict \rga proves decided (Section~\ref{sec:recall}).

Paracomplete logics have a reputation for being too weak for
metamathematics.  Developing a closely related paracomplete theory of
truth, Maudlin concluded that ``banning the use of T-Inferences in
\emph{reductio} arguments eliminates all of the profound results of
metalogic''~\cite[p.~134]{maudlin06truth}, and earlier work on
grounded deduction left open whether grounded reasoning loses
G\"odel's diagonal arguments or gains from them~\cite{ford24reasoning}.
\rga's restriction is of the same kind: \rga discharges a hypothesis
only when the hypothesis is provably decided.  This paper shows that,
for \rga and for the systems built on \rga in the same way, G\"odel's
theorems survive, in a precisely stated form.  G\"odel's first and
second theorems hold as facts of HOL about \rga.  Within \rga,
formalized G\"odel~II holds as an entailment and in
decided-antecedent form (defined below), and L\"ob's theorem as an
entailment and as a rule.  The second theorem's proof, read in these
forms, needs only a short list of conditions on the system and its
checker, and the paper proves the theorem once, from those conditions.

\paragraph{Takeaways.}

\begin{table}[tbp]
\centering
\footnotesize
\setlength{\tabcolsep}{3pt}
\begin{tabular}{@{}>{\raggedright\arraybackslash}p{0.13\linewidth}>{\raggedright\arraybackslash}p{0.22\linewidth}>{\raggedright\arraybackslash}p{0.34\linewidth}>{\raggedright\arraybackslash}p{0.27\linewidth}@{}}
\hline
step & textbook proof (L\"ob's conditions) & \rga: implication, entailment, rule & \rga: decided-antecedent form \\
\hline
self-reference & $\vdash G \liff \neg\PrC(\quo{G})$ & $G$ is literally $\forall q.\ \neg\Chk(q, S_G)$, and $\vdash S_G = \bquo{G}$; $G \eqv \neg\PrC(\bquo{G})$ holds; the biconditional fails & --- \\
provability & the $\Sigma_1$ formula $\exists q.\ \mathit{Prf}(q, x)$ & $\PrC(x) = \exists q.\ \Chk(q, x)$, with a decided checker $\Chk$ and its structural laws & an entailment from $\PrC(c)$ is the sentence $\forall q.\,(\Chk(q, c) \limp B)$ \\
D1 & a rule & a rule & --- \\
D2 & a single sentence for each $\varphi, \psi$ & fails as an implication at $G$; holds as an entailment and as a rule; the proof uses its ex falso form EXPL & decided form for every $\varphi, \psi$, with an explicit PR transformer \\
D3 & $\PrC(\quo\varphi) \limp \PrC(\quo{\PrC(\quo\varphi)})$ & fails as an implication at $G$; holds as an entailment and as a rule & decided form for every $\varphi$, with an explicit PR transformer, and one sentence over the code; D3G is the same construction at $G$'s code term $S_G$ \\
formalized G\"odel~I & $\vdash \ConR \limp \neg\PrC(\quo{G})$ & fails as an implication; the entailment $\ConR \vdash G$ holds, and the proof uses it as a case split on the decided $\Chk(v, S_G)$ under $\forall$-introduction, as in free-variable arithmetic & contraposed: $\forall q.\,(\Chk(q, \bquo{G}) \limp \neg\ConR)$, also with an explicit transformer \\
hypotheses discharged & $\Sigma_1$ provability and $\Pi_1$ consistency statements & only decided ones, as \rga's rules require: one decided hypothesis in each of EXPL and D3G & --- \\
G\"odel~I & $\nvdash G$; $\nvdash \neg G$ from $\Sigma_1$-soundness & both $\nvdash G$ and $\nvdash \neg G$; the witness property, a consequence of soundness, plays the role of $\Sigma_1$-soundness & --- \\
L\"ob & formalized L\"ob provable; G\"odel~II is its instance at $\z = \tofn 1$ & formalized L\"ob fails as an implication (at $\z = \tofn 1$) and holds as an entailment; L\"ob's rule holds & $\forall q.\,(\Chk(q, \bquo{\PrC(\bquo\psi) \limp \psi}) \limp \PrC(\bquo\psi))$, every closed $\psi$ \\
formalized G\"odel~II & $\vdash \ConR \limp \neg\PrC(\quo{\ConR})$ & fails as an implication; holds as an entailment, $\ConR \vdash \neg\PrC(\bquo{\ConR})$ & $\forall q.\,(\Chk(q, \bquo{\ConR}) \limp \neg\ConR)$, also with an explicit transformer to proofs of $\z = \tofn 1$ \\
conditions used & the diagonal lemma and D1--D3 & three structural conditions, NEC, D3 and EXPL as proof transformers, D1 as closed acceptance, and consistency; D2 is not among them & the decided-antecedent form of formalized G\"odel~II uses the structural and transformer conditions only \\
\hline
\end{tabular}
\caption{G\"odel's second theorem: the textbook proof and \rga's.  Here
  $\PrC(x)$ is $\exists q.\ \Chk(q, x)$, written $\neg\forall q.\
  \neg\Chk(q, x)$ in \rga; $S_G$ is a closed term that \rga proves equal
  to $G$'s code; EXPL (``a proof of $G$ and a proof of $\neg G$ give a
  proof of $\z = \tofn 1$'') and D3G (D3 at $S_G$) are the decided
  forms of D2 and D3 that the proof uses (Section~\ref{sec:goedel}).
  The last column gives the decided-antecedent form of
  \S\ref{sec:recall-forms}, which \rga proves where the hypothesis is
  a provability statement.  The implication with $\PrC(c)$ as
  antecedent fails in the rows D2 and D3 (at $G$), L\"ob (at $\z =
  \tofn 1$) and formalized G\"odel~II.  The last row is the G\"odel
  conditions of Section~\ref{sec:locale}, under which the second
  theorem holds for every G\"odel extension of \rga.  Hilbert and
  Bernays already state the third condition with a primitive-recursive
  equation as its antecedent, and free-variable presentations use
  decided forms throughout (Section~\ref{sec:related}).}
\label{tab:goedel-compare}
\end{table}

Six results organize the paper, in the order of the paper's sections,
and all six are mechanized in Isabelle/HOL.  A term is \emph{decided}
when \rga proves the term true or false, $\vdash t \lor \neg t$.

\begin{enumerate}
\item[(a)] \emph{G\"odel's incompleteness proof survives as
  entailments} (Section~\ref{sec:breaks}).  A conditional step has four
  forms in \rga (\S\ref{sec:recall-forms}): an implication $\vdash A
  \limp B$, one sentence; for a hypothesis $\exists s.\,D(s)$ with $D$
  decided, such as a provability statement, the
  \emph{decided-antecedent form} $\vdash \forall s.\,(D(s) \limp B)$; an
  entailment $A \vdash B$, one derivation under a hypothesis, which cut
  turns into a transformation of derivations; and a rule, ``if $\vdash
  A$ then $\vdash B$''.  \rga's conditional is material, $A \limp B$ is
  $\neg A \lor B$, so by the disjunction property
  (Theorem~\ref{thm:br-mech}) \rga proves no implication between two
  sentences that \rga neither proves nor refutes, not even $G \limp G$.
  The sentences that the textbook proof relates are independent, so
  every implication that the textbook proof uses fails, by theorem;
  formalized L\"ob fails too, at $\z = \tofn 1$
  (Theorem~\ref{thm:br-loeb}).  The deduction theorem holds exactly at
  provably decided hypotheses (Theorem~\ref{thm:rc-dt}).  Every
  entailment and every rule that the textbook proof uses holds:
  $G \eqv \neg\PrC(\bquo{G})$, $\ConR \vdash G$, D2 at well-formed
  formulas, D3 at every term, and formalized L\"ob for every closed
  well-formed sentence (Theorems~\ref{thm:br-diag-ent},
  \ref{thm:br-hyp}, \ref{thm:br-dc-ent} and~\ref{thm:br-loeb-ent}).
  Where the hypothesis is a provability statement, the entailment is
  equivalent to a sentence in decided-antecedent form, which \rga
  proves; only the extraction of the quantifier, to $(\exists
  s.\,D(s)) \limp B$, fails in \rga.
\item[(b)] \emph{G\"odel's theorems hold for \rga, and \rga proves the
  second theorem about itself} (Section~\ref{sec:goedel}).  \rga does
  not derive its consistency sentence $\ConR = \forall p.\ \neg\Chk(p,
  \bquo{\z = \tofn 1})$ (Theorem~\ref{thm:g-g2}).  This
  non-derivability is proved in HOL, and the HOL proof also shows \rga
  consistent, as any proof that a system does not derive some sentence
  must.  The proof follows the route through a decided proof predicate
  that free-variable arithmetic takes to the second theorem;
  free-variable arithmetic, like \rga, discharges a hypothesis only
  when the hypothesis is decided.  \rga carries the proof out internally
  up to the proof's last step: decided forms of D2 and D3 with explicit
  primitive-recursive proof transformers, and a case split on the
  decided proof predicate under $\forall$-introduction in place of
  formalized G\"odel~I (Theorem~\ref{thm:g-split}).  The last step,
  that \rga does not derive $\z = \tofn 1$, is a fact of HOL that rests
  on \rga's soundness.  \rga also proves the second theorem about
  itself, as the entailment $\ConR \vdash \neg\PrC(\bquo{\ConR})$
  (Theorem~\ref{thm:g-formal}) and in decided-antecedent form,
  $\forall q.\,(\Chk(q, \bquo{\ConR}) \limp \neg\ConR)$, with an
  explicit transformer from accepted proofs of $\ConR$ to accepted
  proofs of $\z = \tofn 1$ (Theorem~\ref{thm:g-sentence}).  Only the
  implication forms fail.  G\"odel's first theorem holds too, and the
  plain G\"odel sentence is independent (Theorem~\ref{thm:g-g1}):
  \rga's soundness supplies a witness property that plays the role of
  $\Sigma_1$-soundness.  Table~\ref{tab:goedel-compare} sets the
  textbook proof of the second theorem and \rga's proof side by side.
\item[(c)] \emph{The second theorem, proved once}
  (Section~\ref{sec:locale}).  The proof of~(b) uses few facts about
  \rga and its checker, and Section~\ref{sec:locale} states them as
  hypotheses about an abstract entailment relation on \rga's language,
  closed under \rga's rules, and an abstract checker: three structural
  conditions, NEC, D3 and EXPL as proof transformers, D1 as closed
  acceptance, and consistency (Definition~\ref{def:lc-cond}).  The
  conditions are the \emph{G\"odel conditions}, and a system that meets them is a
  \emph{G\"odel extension of \rga}.  Every G\"odel extension fails to
  derive its G\"odel sentence and its consistency sentence for its
  checker (Theorem~\ref{thm:lc-g2}), and derives formalized G\"odel~II
  as an entailment and in decided-antecedent form
  (Theorem~\ref{thm:lc-textbook}).  D2 is not among the conditions, and
  neither is the checker's faithfulness.  \rga with its tree checker is
  a G\"odel extension (Theorem~\ref{thm:lc-rga}), and 62 facts of the
  G\"odel theories for \rga, among them G\"odel's second theorem in a
  second proof (Corollary~\ref{cor:lc-rga-g2}), are defined as
  instances of the locale's facts, with their propositions fixed by
  build-time fingerprints.
\item[(d)] \emph{\iga as a G\"odel extension}
  (Sections~\ref{sec:checker} and~\ref{sec:iga-goedel}).  Inductive
  Grounded Arithmetic~\cite{iga-paper}, \iga, is \rga with one more
  rule, ground induction.  Section~\ref{sec:checker} writes \rga's
  checker once over its table of rule detectors, so that the derivation
  laws hold at every table that satisfies one law, the shift law
  (Theorem~\ref{thm:ck-shift}), and \iga's checker $\ChkI$ is \rga's
  table plus one tag (Definition~\ref{def:ck-iga} and
  Theorem~\ref{thm:ck-iga-faith}).  The certificate
  of the third derivability condition, the substantial piece of the
  proof for \rga, is a tree of \rga's own rules at every table
  (Theorem~\ref{thm:ig-d3}), and closed acceptance holds by one-way
  closed evaluation, without \rga's soundness
  (Theorem~\ref{thm:ig-closed}).  So \iga with $\ChkI$ is a G\"odel
  extension (Theorem~\ref{thm:ig-witness}), and \iga does not derive
  $\ConI$, its consistency sentence for the checker $\ChkI$
  (Theorem~\ref{thm:ig-g2}).
\item[(e)] \emph{Adoption and the consistency family}
  (Section~\ref{sec:adopt}).  Adoption is a reflection rule keyed to a
  base system's checker (Definition~\ref{def:ad-rule}).  The rule's
  premise is a derived warrant: the sentence that, for every $v$, the
  checker accepts the proof that an explicit primitive-recursive tree
  computes for the $v$-th instance of a universal.  From that premise,
  adoption concludes the universal.  Adoption differs from Shoenfield's
  recursive $\omega$-rule: the premise is internal and decided, keyed
  to a fixed base checker, and not iterated.  The adoption extensions
  \rgaP and \igaP are sound for the $\omega$-model $\Mw$ and consistent
  (Theorems~\ref{thm:ad-sound} and~\ref{thm:ad-iga-rule}).  \rga alone
  derives the warrant that yields $\ConR$, so \rgaP derives $\ConR$
  (Theorem~\ref{thm:ad-sep}), and \igaP derives both $\ConI$ and $\ConR$
  (Theorem~\ref{thm:ad-iga-con}).  The separations are strict: \rga is
  a proper part of \rgaP and of \igaP, and \iga is a proper part of
  \igaP (Theorems~\ref{thm:ad-sep}, \ref{thm:ad-iga-sep}
  and~\ref{thm:ig-sep}).  The section ends with a second kind of
  adoption, \iga with a supplied claim that a code lies in Kleene's
  $\mathcal O$.  G\"odel's second theorem holds at every sound claim,
  with a checker that adds one axiom cell to \iga's table, and the
  sound claim built from that checker derives the claim extension's
  consistency sentence (Theorems~\ref{thm:ac-g2}
  and~\ref{thm:ac-larger}).  A claim that \iga never refutes can still
  give an inconsistent extension (Theorem~\ref{thm:ac-quine}).
  Unrefutedness is no closure condition on adopted claims or on
  well-foundedness sentences (Theorem~\ref{thm:ac-direct}), and
  consistency with the system is none on well-foundedness sentences
  (Theorems~\ref{thm:ac-liar} and~\ref{thm:ac-rosser}); that it is
  none on adopted claims either is a prediction.
\item[(f)] \emph{Does \iga alone derive $\ConR$?}
  (Section~\ref{sec:open}).  The question is open, and the paper states
  what is known and what each route to an answer would have to supply.
  The internal-reasoning paper's no-go theorem excludes, at \iga, the
  internal-soundness fold with that paper's motive
  (Theorem~\ref{thm:op-nogo}).  On the Gentzen route, \iga derives
  $\ConR$ from a reduction on \rga's proof codes that \rga proves
  correct, under the further premises of Theorem~\ref{thm:op-gentzen},
  among them an ordinal assignment that \rga proves bounded below
  $\varphi(\omega, 0)$; the missing premise is an ordinal analysis of
  \rga.  On the route through a stage-indexed truth, the fold is a
  theorem one rule up.  \oiga, which adds to \iga the decidedness rule
  ATI of Omega Grounded Arithmetic (\oga)~\cite{oga-paper}, derives
  $\ConR$ from no hypothesis through a truth predicate indexed by
  stages below $\varepsilon_0$, so \oiga strictly exceeds \rga in
  consistency strength (Theorem~\ref{thm:op-oiga-con}).  That \iga
  alone does not derive the bivalence of that predicate, which the
  route needs, is a prediction.
\end{enumerate}

\paragraph{Kinds of claim, and viewpoints.}
Every theorem in this paper is a theorem of HOL about \rga or about a
system built on \rga, and the paper says of each theorem which kind of
claim the theorem makes: a schema of derivations, one internal
sentence, an entailment, or a HOL fact about derivability.  For every
conditional claim the paper also says which of the four forms of a
conditional is meant.  \S\ref{sec:recall-forms} sets out both
distinctions (Table~\ref{tab:rc-forms}).  Where a result is read, the
paper labels the reading by viewpoint.  A rule, ``if $\vdash A$ then
$\vdash B$'', or a non-derivability $\nvdash A$ is \emph{classical}: a
HOL statement about derivability, in HOL's reading.  An entailment or a
sentence that the system derives is labeled by naming the system,
\emph{a derivation in \rga} (or in \iga, in \oga, or in the extended
system at hand): the system itself derives that judgment.  A claim that
the development does not mechanize is a \emph{prediction}.  The paper
does not use the bare word ``grounded'' as a label, because here it
names grounded truth: \iga derives $\ConI \vdash G$, and $G$ is
ungrounded.  That $\ConR$ or $\ConI$ expresses the
consistency of \rga or \iga is the classical reading, through the
checker's faithfulness at closed codes.  Every theorem carries its
Isabelle name in its header.

\paragraph{What is new.}
The arithmetization and the equation form of self-reference are
classical, and so are proofs of the second theorem that avoid the
textbook route: free-variable arithmetic proves the second theorem
with decided proof transformers~\cite{coquand26auto}, building on
Hilbert and Bernays' proof-predicate
formulation~\cite{hilbert39grundlagen,kurahashi20note}; Jeroslow's
proof avoids D2~\cite{jeroslow73redundancies}; and Kritchman and Raz's
proof uses no diagonal sentence~\cite{kritchman10surprise}.  Abstract
second theorems are not new either.  Beklemishev and Shamkanov proved
the second theorem and its formalized version as entailments, for any
consequence relation that meets four conditions and has a Jeroslow
fixed point~\cite{beklemishev16abstract}, and Popescu and Traytel
mechanized abstract versions of both theorems, with the derivability
conditions as assumptions~\cite{popescu21distilling}.
\rga meets Beklemishev and Shamkanov's conditions on \rga's decided and
$\Sigma$ sentences, though not on all \rga sentences, and the
contrapositive form of \rga's formalized second theorem is an instance
of Beklemishev and Shamkanov's Theorem~1 on that fragment
(Section~\ref{sec:related}).  What is new is, first, the internal
decided route: G\"odel's proof read with entailments, carried out
inside a concrete arithmetic with gaps at the arithmetic level, where
the decided forms must be built and proved.  The substantial pieces
there are the builder of D3's premise certificates at abstract
arguments, the cut transformers that internalize entailments, and
formalized G\"odel~II in decided-antecedent form with an explicit
witness.  Second, the G\"odel conditions: conditions for the second
theorem that a system without implication introduction over undecided
hypotheses meets, without D2 and without the checker's faithfulness,
together with two concrete instances, \rga and \iga, whose checkers
are \rga terms.  Third, the checker at a table, under which an
extension of \rga by a rule gets its checker, and its D3 certificate,
from \rga's.  Fourth, adoption and the strict separations it gives,
and, for adopted notation claims, the second theorem at every sound
claim and the theorems on what falsifies a claim.
The author knows of no earlier proof of the second theorem for an
arithmetic whose arithmetic sentences themselves can lack a truth
value.  Heyting arithmetic has no truth-value gaps, and the gaps of
KF, PKF and Field's theory are confined to a truth predicate over
classical arithmetic.  The author also knows of no earlier mechanized
proof of the second theorem for a concrete arithmetic with quantifiers
whose object logic is not classical (Section~\ref{sec:related} lists
the mechanized proofs).  The metatheory is classical HOL, and
$\emptyset \nvdash \ConR$ is proved there, like Paulson's
theorem~\cite{paulson15mechanised}, by a classical argument; what is
new lies in the objects inside \rga and in what \rga and \iga derive
about those objects.

\paragraph{Relation to the earlier version of the internal-reasoning paper.}
The internal-reasoning paper's version~3, arXiv:2608.16140v3, contained
G\"odel's theorems for \rga.  Section~\ref{sec:bg-goedel}, Sections~\ref{sec:breaks}
and~\ref{sec:goedel}, Table~\ref{tab:goedel-compare} and
Appendix~\ref{app:goedel} revise that material, and
Section~\ref{sec:recall} restates, by citation, the machinery of that
paper which the proofs use.  The current version of the
internal-reasoning paper keeps the internal machinery and cites this
paper for G\"odel's theorems.  Sections~\ref{sec:locale}
to~\ref{sec:open} are new.  Version~3 asked whether the second theorem
holds uniformly for consistent extensions of \rga.
Section~\ref{sec:locale} answers that question for rule extensions on
\rga's language with a checker that meets the G\"odel conditions, and
Sections~\ref{sec:checker} and~\ref{sec:iga-goedel} supply \iga as an
instance.  Adoption, the consistency sentences that adoption derives,
and the open question of Section~\ref{sec:open} are new as well.

\paragraph{Outline.}
Section~\ref{sec:bg} recalls \rga and the classical proof of the
second theorem, and Section~\ref{sec:recall} recalls the codes, the
four forms of a conditional and the tree checker from the
internal-reasoning paper.  Section~\ref{sec:breaks} reads G\"odel's
proof in the four forms, and Section~\ref{sec:goedel} proves G\"odel's
theorems for \rga.  Section~\ref{sec:locale} states the G\"odel
conditions and proves the second theorem once, from the conditions.
Section~\ref{sec:checker} builds the checker at a table and \iga's
checker, and Section~\ref{sec:iga-goedel} proves the second theorem at
\iga.  Section~\ref{sec:adopt} on adoption derives the consistency
sentences in the adoption extensions and proves the separations, and
Section~\ref{sec:open} sets out the open question.
Section~\ref{sec:mech} describes the mechanization, its trusted base
and its build-time checks, Section~\ref{sec:related} discusses related
work, and Section~\ref{sec:concl} concludes.  The appendices give a
glossary of the notation (Appendix~\ref{app:notation}), the proof of
G\"odel's second theorem for \rga at the level of the mechanized lemmas
(Appendix~\ref{app:goedel}), the locale layers and their
interpretations at \rga and \iga (Appendix~\ref{app:locale}), the
generic checker's and \iga's tables (Appendix~\ref{app:checker}), and
adoption in detail (Appendix~\ref{app:adopt}).

\paragraph{AI disclosure.}
The machine-checked proofs and the first drafts of this paper's text
were produced with substantial assistance from Claude models from
Anthropic, working under the direction and review of the author.  The
definitions that a reader must trust are listed in
Section~\ref{sec:mech}, which also names the released development.

\section{Background}
\label{sec:bg}

This section recalls \rga from the RGA paper~\cite{rga-paper} and from
the paper on internal reasoning in \rga~\cite{conrga-paper}, the
internal-reasoning paper for short, and then the classical proof of
G\"odel's second theorem.  Section~\ref{sec:recall} recalls the
machinery inside \rga on which this paper's proofs run.

\subsection{Reflective Grounded Arithmetic}
\label{sec:bg-rga}

\paragraph{Grounded truth.}
\rga belongs to the grounded-deduction family of
logics~\cite{ford24reasoning}.  A \emph{sentence} is a closed term, and
a sentence is a computation: it is \emph{grounded-true} when its
evaluation settles on true, \emph{grounded-false} when its evaluation
settles on false, and \emph{ungrounded} when its evaluation never
settles.  Ungroundedness is not a third truth value but the absence
of a verdict, as divergence is the absence of a return value.  Logics
with such gaps are \emph{paracomplete}.  The gaps of \rga sit at the
arithmetic level.  They include the Liar, and also the consistency
sentence, which does not refer to itself.
The internal-reasoning paper compares this notion of grounding with
Kripke's ([\citen{conrga-paper}, \S\,``Reflective Grounded Arithmetic'']).

\paragraph{The language and the semantics.}
\rga's terms form an untyped $\lambda$-calculus with ten operators: the
divergent constant $\bot$, a term with no value, zero $\z$, successor
$\suc$, equality of naturals $=$, negation $\neg$, disjunction $\lor$,
$\lambda$-abstraction, application $\cdot$, a recursor $\cR$, and the
universal quantifier $\cA$.  A formula is a term like any other, with
no separate syntax.  Variables are de Bruijn indices, $\forall x.\,b$
abbreviates $\cA\cdot(\lambda b)$, and $a \limp b$ abbreviates $\neg a
\lor b$.  The ordinary Liar is a term $L$ whose evaluation unfolds to
$\neg L$.  The semantics is a step-indexed evaluation: at fuel $s$ a
term evaluates to a natural number, with $\ctrue = 1$ and $\cfalse =
0$.  The universal $\cA\cdot u$ is true at fuel $s+1$ when a proof in
\rga's own proof system, with code at most $s$, derives the open
instance of $u$, and false when such a proof refutes some numeral
instance of $u$.  Because of this proof certificate, grounded truth is
recursively enumerable, like the set of \rga's theorems.

\paragraph{The proof system.}
The judgment $\Gamma \vdash c$ says that $c$ is derivable from the
finite set of hypotheses $\Gamma$.  Three differences from classical
logic shape every internal argument below.  First, \rga splits cases
on a sentence $a$, and discharges a hypothesis $a$ into an implication
$a \limp b$, only when \rga derives $a \jbool$, which abbreviates $a
\lor \neg a$; so the Liar entails itself, $L \vdash L$, while
$\emptyset \nvdash L \limp L$.  Second, the judgment $t \jnat$, which
abbreviates $t = t$, says that $t$ is a natural number, and universal
elimination and induction carry $\jnat$ premises.  Third, witnesses
are real: a derivable disjunction has a derivable disjunct, and a
derivable negated universal has a derivable counterexample at some
numeral (\S\ref{sec:breaks-mech}).  The internal-reasoning paper lists
the rules that the proofs below cite, with their Isabelle names
([\citen{conrga-paper}, the table of \rga's rules]).

\paragraph{The RGA paper's results.}
The RGA paper proves in Isabelle/HOL that every derivable sentence is
grounded-true (soundness), so that \rga does not derive $\z = \tofn 1$
(consistency), and that for well-formed terms, derivability from no
hypotheses coincides with grounded truth (open
completeness)~\cite{rga-paper}.  Open completeness is not negation
completeness, and it fails under hypotheses.  The G\"odel sentence and
the consistency sentence are not grounded-true, so open completeness
does not reach them.  In this paper, ``in HOL'' refers to statements
about \rga proved in Isabelle/HOL, and ``in \rga'' to derivations in
\rga's own proof system.

\subsection{The classical proof of G\"odel's second theorem}
\label{sec:bg-goedel}

Let $T$ be a consistent, recursively axiomatized theory containing
enough arithmetic, such as Peano arithmetic.  Code formulas and proofs
by numbers, and let $\mathit{Prf}(q, x)$ be a formula, decidable in
$T$, expressing that $q$ codes a $T$-proof of the formula coded by
$x$.  The provability predicate is the $\Sigma_1$ formula
$\PrC(x) \ldef \exists q.\ \mathit{Prf}(q, x)$, and the consistency
sentence is $\ConR_T \ldef \neg\PrC(\quo{\z = \tofn 1})$.  The proof
of $T \nvdash \ConR_T$ has five ingredients~\cite{godel31uber,hilbert39grundlagen,lob55solution,boolos93logic}.

\begin{enumerate}
\item \emph{The diagonal lemma.}  For every formula $\varphi(x)$ there
  is a sentence $G$ with $T \vdash G \liff \varphi(\quo{G})$.  Applied
  to $\neg\PrC(x)$, the lemma gives the G\"odel sentence, with $T \vdash G
  \liff \neg\PrC(\quo{G})$.
\item \emph{D1.}  If $T \vdash \varphi$, then $T \vdash
  \PrC(\quo\varphi)$.
\item \emph{D2.}  $T \vdash \PrC(\quo{\varphi \limp \psi}) \limp
  (\PrC(\quo\varphi) \limp \PrC(\quo\psi))$.
\item \emph{D3.}  $T \vdash \PrC(\quo\varphi) \limp
  \PrC(\quo{\PrC(\quo\varphi)})$, an instance of formalized
  $\Sigma_1$-completeness.
\item \emph{Formalized G\"odel~I.}  From D1--D3 and the diagonal
  equivalence, $T \vdash \ConR_T \limp \neg\PrC(\quo{G})$, and hence
  $T \vdash \ConR_T \limp G$.
\end{enumerate}

\noindent
If $T \vdash \ConR_T$, then $T \vdash G$, so $T \vdash \PrC(\quo{G})$
by D1 and $T \vdash \neg G$ by the diagonal equivalence, and $T$ is
inconsistent.  L\"ob's theorem~\cite{lob55solution} generalizes the
argument: if $T \vdash \PrC(\quo\varphi) \limp \varphi$, then $T
\vdash \varphi$; G\"odel's second theorem is the instance of L\"ob's theorem at $\varphi
= (\z = \tofn 1)$.  G\"odel's first theorem says that a consistent $T$
does not prove $G$.  That $T$ does not refute $G$ needs more than
consistency: $\mathrm{PA} + \neg\ConR_{\mathrm{PA}}$ is consistent and
refutes $G$.  G\"odel assumed $\omega$-consistency (1-consistency
suffices); Rosser's modified sentence is independent of every
consistent $T$~\cite{rosser36extensions}.

The whole argument is valid in intuitionistic arithmetic: Heyting
arithmetic does not prove its own consistency, by the same
proof~\cite{troelstra73metamathematical}.  What the argument uses is
implication introduction over undecided antecedents.  The diagonal
lemma asserts a biconditional, and D2, D3 and formalized G\"odel~I are
implications whose antecedents, $\PrC(\ldots)$ and $\ConR_T$, a
consistent $T$ can in general neither prove nor refute.  In Heyting
arithmetic $(\exists x.\,D(x)) \limp B$ and $\forall x.\,(D(x) \limp
B)$ are equivalent.  In \rga, whose $\limp$ is material, only the
second form is available at an undecided $\exists x.\,D(x)$ with $D$
decided.  Section~\ref{sec:breaks} shows that
in \rga each of these steps fails as an implication, by theorem, while
the entailments behind them hold, and Section~\ref{sec:goedel} proves
G\"odel's second theorem with the entailments in place of the implications.

\section{Codes, conditionals and the checker}
\label{sec:recall}

The proofs of this paper run on machinery that the internal-reasoning
paper~\cite{conrga-paper} builds inside \rga: codes of terms and
proofs, compiled primitive-recursive machines, and a proof checker.
This section recalls the parts of that machinery that the proofs
below use, and states the four forms of a conditional by which this
paper reads G\"odel's proof.  The statements are the
internal-reasoning paper's, cited by number, and so are their proofs.

\subsection{Codes and numerals}
\label{sec:recall-codes}

Every \rga term $t$ has a \emph{code} $\quo{t}$, a natural number built
with the Cantor pairing function $\cpair{x}{y}$: a variable, or an
operator node with the operator's number, the number of variables it
binds and the list of its arguments, each nested as a pair
([\citen{conrga-paper}, \S\,``Codes, numerals, and four kinds of statement'']).  A proof is a finite list of judgments,
and it has a code in the same way.  The coding is injective but not
onto.  The \emph{numeral} $\tofn{n}$ is the term $\suc^n\z$, and the
\emph{code numeral} of $t$ is the numeral of its code, $\bquo{t} =
\tofn{\quo{t}}$.  Plain corners $\quo{t}$ denote a number of the
metalanguage HOL; underlined corners $\bquo{t}$ denote a closed \rga
term, which \rga's own derivations can mention.  Codes are large: the
code of the false equation $\z = \tofn 1$ has 206 decimal digits, so
the paper never writes a code out.

The development describes primitive-recursive functions by
\emph{index trees} of initial functions, composition nodes and
recursion nodes, and the compiler $\cmp$ of the RGA
paper~\cite{rga-paper} maps each index tree $f$ to a closed \rga term
$\cmp(f)$.  At numerals, \rga proves every compiled function's values,
$\emptyset \vdash \cmp(f) \cdot \tofn{x} \cdot \tofn{y} =
\tofn{\den{f}(x,y)}$ ([\citen{conrga-paper}, theorem~\textsf{cmp\_eq\_num}]), whatever the
tree.  At arguments that are not numerals, \rga reasons about
$\cmp(f)$ by unfolding the tree one node at a time
([\citen{conrga-paper}, theorem~\textsf{cmp\_prRec\_zero}]), so that reasoning needs a tree that
the development writes out.

\subsection{Four kinds of statement, and four forms of a conditional}
\label{sec:recall-forms}

\paragraph{Four kinds of statement.}
Every theorem in this paper is a theorem of HOL about \rga.  The
theorems differ in where the quantifier over codes sits, and the
internal-reasoning paper distinguishes four kinds
([\citen{conrga-paper}, \S\,``Codes, numerals, and four kinds of statement'']).
\begin{itemize}
\item A \emph{schematic} theorem quantifies in HOL: for every HOL
  object that meets the theorem's hypotheses, \rga derives one judgment
  about that object's code numeral.  A derived rule, which says that
  \rga derives a conclusion whenever it derives the premises, is
  schematic too.  A \emph{rule in every context} derives its conclusion
  in every context $\Gamma$ in which \rga derives the premises; a
  \emph{rule on theorems} applies only to closed theorems $\emptyset
  \vdash c$.
\item An \emph{internal} theorem says that \rga derives one particular
  sentence, in which the quantifier over codes is \rga's own $\cA$.
\item A \emph{derivation under hypotheses} says that \rga derives one
  judgment $\Gamma \vdash c$ in a particular non-empty context
  $\Gamma$.
\item A \emph{meta} theorem concludes a HOL proposition about \rga
  other than a derivation, such as $\mathit{False}$ or that \rga does
  not derive some sentence.
\end{itemize}
The consistency sentence shows the difference between the first two
kinds.  For each number $p$, \rga derives $\neg\Chk(\tofn{p},
\bquo{\z = \tofn 1})$, one derivation per number: a schema
([\citen{conrga-paper}, theorem~\textsf{goedelW\_omega\_incomplete}]).  The consistency sentence $\ConR =
\forall p.\ \neg\Chk(p, \bquo{\z = \tofn 1})$ is one sentence, and an
internal theorem would derive it.  G\"odel's second theorem for \rga
says that no such derivation exists (Theorem~\ref{thm:g-g2}), and
that theorem is a meta theorem.

\paragraph{Four forms of a conditional.}
A conditional claim, ``$A$ gives $B$'', has four forms.
\begin{itemize}
\item An \emph{implication} $\vdash A \limp B$ is one sentence, which
  abbreviates $\neg A \lor B$.  Read computationally, for closed $A$
  and $B$, the implication is one closed computation that ends with the
  value true.  \rga introduces $A \limp B$ freely from $\neg A$ or from
  $B$, but \rga discharges a hypothesis $A$ into $A \limp B$ only when
  \rga proves $A$ decided.  The implication carries no transformation
  of evidence.
\item An \emph{entailment} $A \vdash B$ is one derivation under the
  hypothesis $A$, a derivation under hypotheses in the sense above.
  Read as a statement about derivations, an entailment always has
  content: by cut, it turns any derivation of $A$ in any context into a
  derivation of $B$ in the same context.  For hypotheses about proof
  codes this transformation also exists inside \rga, as an explicit
  primitive-recursive map of accepted proof codes
  (Theorem~\ref{thm:g-cut}).  Read as a statement about truth, an
  entailment says that $B$ is grounded-true whenever $A$ is, which
  holds by soundness.  At a hypothesis that is never grounded-true,
  such as $\ConR$ or the G\"odel sentence $G$, this truth reading
  holds vacuously.  The entailment's content is then the derivational
  one, together with truth preservation in the $\omega$-semantics of
  [\citen{conrga-paper}, \S\,``The $\omega$-semantics''], in which $\ConR$ and $G$ are true and
  \rga is sound in every context.  Two-way
  entailment, $A \vdash B$ and $B \vdash A$, is written $A \eqv B$.
\item A \emph{rule}, ``if $\vdash A$ then $\vdash B$'', is a closure
  property of \rga's theorems, a fact of HOL.  The rule may hold in
  every context or on theorems only.
\item The \emph{decided-antecedent form} $\vdash \forall s.\,(D(s)
  \limp B)$ applies to a hypothesis of the shape $\exists s.\,D(s)$ with
  $D$ decided.  This form is the subject of the next paragraph.
\end{itemize}
An implication gives the entailment by modus ponens, and an entailment
gives the rule by cut.  The step back from an entailment to an
implication is the deduction theorem, and the internal-reasoning paper
proves its exact boundary.

\begin{thm}[\textsf{deduction\_theorem\_exact}, \textsf{sigma\_deduction\_theorem}]
\label{thm:rc-dt}
Let $\Gamma$ be a context.
\begin{enumerate}
\item For every well-formed $A$: $A, \Gamma \vdash B$ gives $\Gamma
  \vdash A \limp B$ for every $B$ if and only if $\Gamma \vdash A \lor
  \neg A$.  For the ordinary Liar $L$, $L \vdash L$ while $\emptyset
  \nvdash L \limp L$.
\item If $\vz\jnat, \lft\Gamma \vdash D(\vz)\jbool$, then for every
  $C$: $\exists s.\,D(s), \Gamma \vdash C$ if and only if $\Gamma
  \vdash \forall s.\,(D(s) \limp C)$.  In particular, at a closed code
  term $c$ with $\Gamma \vdash c\jnat$ and a closed $C$,
  $\PrC(c), \Gamma \vdash C$ if and only if $\Gamma \vdash \forall
  q.\,(\Chk(q, c) \limp C)$.
\end{enumerate}
\end{thm}

\noindent
Theorem~\ref{thm:rc-dt} is [\citen{conrga-paper}, theorem~\textsf{deduction\_theorem\_exact}].  Here $\vz$ is
the de Bruijn variable~0 and $\lft\Gamma$ lifts the free variables of
$\Gamma$ by one, so $\vz$ is fresh; $\PrC(c) = \exists q.\,\Chk(q, c)$
is the provability formula over the checker of
\S\ref{sec:recall-chk}.  Both parts of the theorem are HOL
equivalences between derivability facts.  Part~2, the $\Sigma$
deduction theorem, yields the two moves that replace excluded middle in
the proofs below.

\begin{cor}[\textsf{dec\_contrapos}, \textsf{sigma\_contrapos}]
\label{cor:rc-contrapos}
Contraposition, from $A, \Gamma \vdash B$ to $\neg B, \Gamma \vdash
\neg A$, holds when $\Gamma \vdash A \lor \neg A$ and when $A$ is
$\exists s.\,D(s)$ with $D$ decided as in Theorem~\ref{thm:rc-dt}(2);
in particular a hypothesis of that shape that entails $\z = \tofn 1$
is refuted.  Contraposition fails in general: $L_{\Tru} \wedge \neg L_{\Tru}
\vdash \z = \tofn 1$, while $\neg(\z = \tofn 1) \nvdash \neg(L_{\Tru}
\wedge \neg L_{\Tru})$.
\end{cor}

\noindent
Corollary~\ref{cor:rc-contrapos} is [\citen{conrga-paper}, corollary~\textsf{dec\_contrapos}].
$L_{\Tru}$ is the Liar for \rga's truth predicate
([\citen{conrga-paper}, \S\,``The truth predicate and the Liar'']), a sentence that \rga neither proves nor
refutes.  The corollary's second case gives reductio and contraposition
at a hypothesis $E = \exists s.\,D(s)$: if $E \vdash \neg E$, then
$\emptyset \vdash \neg E$, and if $E \vdash B$, then $\neg B \vdash
\neg E$.  Both moves split on the decided $D(s)$ at a fresh $s$, never
on $E$.

\paragraph{The decided-antecedent form.}
For a hypothesis $\exists s.\,D(s)$ with $D$ decided, the entailment
$\exists s.\,D(s) \vdash B$ is equivalent to the single sentence
$\forall s.\,(D(s) \limp B)$ (Theorem~\ref{thm:rc-dt}(2)).  Provability
statements have this shape, since $\PrC(c) = \exists q.\,\Chk(q, c)$
and the checker's verdict is decided; so does $\neg\ConR$, which is
literally $\PrC(\bquo{\z = \tofn 1})$.  When the entailment holds, the
sentence in decided-antecedent form is an \rga theorem, and so
grounded-true in its own right, whatever the status of $\exists
s.\,D(s)$.  The form sits between the implication and the entailment:
the implication gives the decided-antecedent sentence, by cases on
the disjunction $\neg\exists s.\,D(s) \lor B$, and that sentence gives
the entailment.  At an undecided $\exists s.\,D(s)$, only the step from
the decided-antecedent sentence to the implication $(\exists s.\,D(s))
\limp B$ fails.  At $\ConR$ that step fails although the sentence holds:
$\vdash \forall q.\,(\Chk(q, \bquo{\ConR}) \limp \neg\ConR)$, while
$\nvdash (\exists q.\,\Chk(q, \bquo{\ConR})) \limp \neg\ConR$
(Theorem~\ref{thm:g-sentence}, parts~1 and~5).  The decided forms of
Section~\ref{sec:goedel} are sentences of this shape, strengthened by
an explicit primitive-recursive transformer that names the
consequent's witness.

At a hypothesis that is neither decided nor of the shape $\exists
s.\,D(s)$, such as $\ConR$ or $G$, entailment and implication come
apart: $\ConR \vdash G$ holds (Theorem~\ref{thm:br-hyp}), while $\ConR
\limp G$ is not a theorem (Theorem~\ref{thm:br-g1}).  The entailments
from $\ConR$ or $G$ still follow from sentences in decided-antecedent
form, by contraposition: $\neg G$ is literally $\PrC(S_G)$, the
sentence $\forall q.\,(\Chk(q, S_G) \limp \neg\ConR)$ holds
(Appendix~\ref{app:goedel}; at the code numeral $\bquo{G}$ it is
Theorem~\ref{thm:g-sentence}(3)), and Corollary~\ref{cor:rc-contrapos}
with double negation turns that sentence into $\ConR \vdash G$.

\paragraph{Why discharge needs decidedness.}
\rga's terms include fixed points, so with an unrestricted deduction
theorem Curry's sentence $C$, equal to $C \limp \z = \tofn 1$, would
derive $\z = \tofn 1$ and with it every
sentence~\cite{curry42inconsistency}.  Since $C \limp \z = \tofn 1$ is
$\neg C \lor \z = \tofn 1$, Curry's sentence is Liar-like, and the
restriction to decided hypotheses blocks both paradoxes.  In
Curry--Howard terms, an entailment $A \vdash B$ is a hypothetical
judgment $x{:}A \vdash M{:}B$, while \rga's implication is the boolean
$\neg A \lor B$, which \rga introduces by discharge exactly when $A$ is
decided, as a proof assistant reflects a decidable proposition into a
boolean.  The decided-antecedent form is what \rga has of the BHK
reading of $(\exists s.\,D(s)) \limp B$: one sentence that discharges
the decided $D(s)$ at every candidate witness $s$.

Table~\ref{tab:rc-forms} follows two conditionals of G\"odel's proof
through the four forms.  Meta facts of HOL, such as $\nvdash \ConR$,
are not a form of the conditional, and the table lists them below the
rules.  Each form can be a HOL schema, with one instance per HOL
object, or one internal sentence over codes; the decided D3 of
Section~\ref{sec:goedel} comes in both versions
(Theorems~\ref{thm:g-d3-any} and~\ref{thm:g-d3-uniform}).

\begin{table}[tbp]
\centering
\small
\setlength{\tabcolsep}{3pt}
\begin{tabular}{@{}>{\raggedright\arraybackslash}p{0.15\linewidth}>{\raggedright\arraybackslash}p{0.17\linewidth}>{\raggedright\arraybackslash}p{0.31\linewidth}>{\raggedright\arraybackslash}p{0.31\linewidth}@{}}
\hline
form & written & from $\ConR$ to $G$ & from $\PrC(\bquo{\ConR})$ to $\neg\ConR$ \\
\hline
implication & $\vdash A \limp B$ & fails (Thm.~\ref{thm:br-g1}) & fails (Thm.~\ref{thm:g-sentence}(5)) \\
decided-antecedent form & $\vdash \forall s.\,(D(s) \limp B)$ & not of this shape ($\ConR$ is not $\exists s.\,D(s)$); the conditional's contrapositive is: $\forall q.\,(\Chk(q, \bquo{G}) \limp \neg\ConR)$ holds (Thm.~\ref{thm:g-sentence}(3)) & holds: $\forall q.\,(\Chk(q, \bquo{\ConR}) \limp \neg\ConR)$ (Thm.~\ref{thm:g-sentence}(1)) \\
entailment & $A \vdash B$ & holds (Thm.~\ref{thm:br-hyp}) & holds (Thm.~\ref{thm:g-formal}) \\
rule in every context & $\Gamma \vdash A \Rightarrow \Gamma \vdash B$ & holds, as the entailment does & holds, as the entailment does \\
rule on theorems & $\vdash A \Rightarrow {}\vdash B$ & holds, vacuously: $\nvdash \ConR$ & holds, vacuously: $\nvdash \PrC(\bquo{\ConR})$ \\
\hline
meta facts & HOL facts & $\nvdash \ConR$, $\nvdash G$ (Section~\ref{sec:goedel}) & $\nvdash \ConR$ \\
\hline
\end{tabular}
\caption{The four forms of a conditional (\S\ref{sec:recall-forms}),
  with two worked conditionals; the internal-reasoning paper's
  [\citen{conrga-paper}, the table of the four forms of a conditional].  A rule holds either in every context
  (a derivable rule) or on theorems only (an admissible rule).  The
  decided-antecedent form applies to a hypothesis $\exists s.\,D(s)$
  with $D$ decided, such as $\PrC(c)$; $\neg G$ is literally
  $\PrC(S_G)$, and $\vdash S_G = \bquo{G}$.  $\nvdash \PrC(\bquo{\ConR})$
  follows from $\nvdash \ConR$ by the reflection rule
  (Theorem~\ref{thm:br-rules}(4)).  Meta facts are not a form of the
  conditional.}
\label{tab:rc-forms}
\end{table}

\subsection{The tree checker}
\label{sec:recall-chk}

The provability predicate of this paper is the internal-reasoning
paper's tree checker $\Chk$ ([\citen{conrga-paper}, definition~\textsf{iChk}]).  The
checker consists of 41 rule detectors, explicit index trees that each
recognize the instances of one of \rga's rules at codes, and a compiled
validity test that requires, at every position of the coded list, a
genuine judgment code on which some detector fires against the
entries after it.  $\Chk(p, c)$ is the conjunction of four decided
tests: validity, $p \neq 0$, a conclusion with no hypotheses, and that
conclusion's code equal to $c$.  The checker is faithful at proof
codes, $\emptyset \vdash \Chk(\bquo{P}, \bquo{c})$ exactly when $P$
is a proof of $\emptyset \vdash c$, and genuine at every number: a
numeral that \rga proves accepted for $\bquo{c}$ is the code of such a
proof ([\citen{conrga-paper}, theorem~\textsf{iChk\_faithful}]).  The verdict is decided at
abstract arguments, $\Gamma \vdash \Chk(p, c)\jbool$ whenever $\Gamma$
derives $p\jnat$ and $c\jnat$ ([\citen{conrga-paper}, theorem~\textsf{hq\_iChk}]), and
that fact licenses the case split on $\Chk(v, S_G)$ under
$\forall$-introduction on which the second theorem's proof turns
(Theorem~\ref{thm:g-split}).  The consistency sentence is $\ConR =
\forall p.\ \neg\Chk(p, \bquo{\z = \tofn 1})$.

\subsection{Transparent witnesses}
\label{sec:recall-transparent}

The second theorem's proof needs a witness for the third derivability
condition at $G$, the closed term $\sigma$ of D3G
(Theorem~\ref{thm:g-d3g}), and the internal proofs reason about
$\sigma$ at an abstract argument.  That reasoning unfolds the compiled
tree of $\sigma$ node by node, so the tree must be written out.  An
index chosen by Hilbert choice, ``some index that denotes the function
$F$'', also compiles and also obeys the closed-numeral laws
([\citen{conrga-paper}, theorem~\textsf{cmp\_eq\_num}]), but its nodes are invisible to the
development's proofs, and the statement that such an index denotes a
primitive-recursive function holds of every index
([\citen{conrga-paper}, \S\,``Compiled primitive-recursive machines'']).  No theorem about an index's value can
therefore say that the index is explicit.  The internal-reasoning
paper makes the requirement a build-time check instead
([\citen{conrga-paper}, \S\,``Machines without opaque indices'']): the check walks the definitional closure
of every machine that the internal-reasoning paper calls transparent
and of every explicit witness and transformer of
Section~\ref{sec:goedel}, $\sigma$ among them.  The check fails the
build if a closure reaches an index chosen by Hilbert choice, or any
other constant defined by choice outside a fixed whitelist of HOL
library constants ([\citen{conrga-paper}, Appendix ``The build-time check''] lists the
roots, the whitelist and the controls).
In this paper, as in the internal-reasoning paper, ``transparent''
means this checked property of machines.

\section{G\"odel's proof in four forms}
\label{sec:breaks}

G\"odel's second theorem has a short classical proof.  A diagonal
sentence $G$ is provably equivalent to its own unprovability.  The
Hilbert--Bernays--L\"ob derivability conditions let the theory reason
about its own provability predicate.  Formalizing the first theorem
gives $\vdash \ConR \limp \neg\PrC(\quo{G})$, and with the diagonal
equivalence, $\vdash \ConR \limp G$.  A proof of $\ConR$ would then
yield a proof of $G$, which the first theorem forbids.  L\"ob's
theorem generalizes the same argument.  This section reads that proof
in the four forms of a conditional and finds which forms \rga has;
Appendix~\ref{app:goedel} gives the proofs of the G\"odel-side results,
and [\citen{conrga-paper}, Appendix ``The traps catalog''] the mechanized traps behind
\S\ref{sec:breaks-mech}.

Each step of that proof is a conditional, and each conditional has the
forms of [\citen{conrga-paper}, \S\,``Codes, numerals, and four kinds of statement'']: an implication, such as $\vdash
\ConR \limp G$; an entailment, such as $\ConR \vdash G$; a rule, such
as ``if $\vdash \ConR$ then $\vdash G$''; and, where the hypothesis is
a provability statement $\PrC(c) = \exists q.\,\Chk(q, c)$, the
decided-antecedent form $\vdash \forall q.\,(\Chk(q, c) \limp B)$.  The classical
proof states its steps as implications.  This section shows that in
\rga every implication the classical proof uses fails, while every
entailment and every rule it uses holds, and the entailments from
provability hypotheses hold in decided-antecedent form.  The deduction
theorem, which would turn the entailments back into implications,
holds only at decided hypotheses, and the hypotheses here are not
decided.  So the proof does not break so much as change form: read
with entailments, the proof goes through, and Section~\ref{sec:goedel}
carries it out inside \rga.  Among the entailments that do fail are D1
at $G$ and reflection at $\z = \tofn 1$, which classical arithmetic
too has only as rules, and the full T-biconditional at the Liar
(Table~\ref{tab:breaks} lists all the failing entailments).

Throughout, $\PrC(x) = \neg\forall q.\ \neg\Chk(q, x)$ is the
provability formula over the tree checker
([\citen{conrga-paper}, \S\,``Internalized completeness'']), that is, $\exists q.\ \Chk(q, x)$.  $\PrC$
is notation, not a constant of the development.  Its matrix $\Chk(q,
x)$ is decided at every number ([\citen{conrga-paper}, theorem~\textsf{hq\_iChk}]), so a
hypothesis $\PrC(c)$ has the shape $\exists s.\, D(s)$ with $D$
decided, the shape of the $\Sigma$ deduction theorem of
[\citen{conrga-paper}, \S\,``Codes, numerals, and four kinds of statement'']; so does $\neg\ConR$, which is literally
$\PrC(\bquo{\z = \tofn 1})$.  For a closed hypothesis $A$, $A \vdash
B$ is the judgment $\{A\} \vdash B$, and $A \nvdash B$ says that \rga
does not derive this judgment.  Every theorem below that says that \rga does not
derive something rests on a soundness theorem proved in HOL, as such a theorem
must, since a claim that \rga does not derive some sentence implies
\rga's consistency.  Most such theorems rest on the companion paper's soundness theorem, through \rga's
consistency, and the $\omega$-model arguments on \rga's soundness for
the $\omega$-semantics of [\citen{conrga-paper}, \S\,``The $\omega$-semantics''].  The entailments
that hold do not rest on soundness ([\citen{conrga-paper}, Appendix ``The build-time check'']);
\S\ref{app:goedel-ent} gives their proofs.

\subsection{One mechanism}
\label{sec:breaks-mech}

Three properties of \rga account for the failures of the implications.
The first property is the deduction theorem's boundary
([\citen{conrga-paper}, theorem~\textsf{deduction\_theorem\_exact}]): \rga's rule for discharging a hypothesis $a$
into an implication $a \limp b$ requires $a$ to be provably decided,
and the deduction theorem holds exactly there; at the Liar sentence
$L$, $L \vdash L$, yet $\emptyset \nvdash L \limp L$.  The other two properties,
and that failure, follow from the soundness theorem; the disjunction
property uses open completeness as well.

\begin{thm}[\textsf{ir\_disj\_prop}, \textsf{rga\_not\_all\_witness}]
\label{thm:br-mech}
\leavevmode
\begin{enumerate}
\item \emph{The disjunction property.}  If $\emptyset \vdash a \lor b$,
  then $\emptyset \vdash a$ or $\emptyset \vdash b$.
\item \emph{The witness property.}  If $\emptyset \vdash
  \neg\forall x.\,\varphi$, then $\emptyset \vdash \neg\varphi(\tofn n)$
  for some numeral $\tofn n$.
\item If $\emptyset \nvdash \delta$ and $\emptyset \nvdash \neg\delta$,
  then $\emptyset \nvdash \delta \limp \delta$ and $\emptyset \nvdash
  \delta \liff \delta$.
\end{enumerate}
\end{thm}

\noindent
Since $a \limp b$ is $\neg a \lor b$, the disjunction property says
that an implication $a \limp b$ is an \rga theorem only if \rga
refutes $a$ or proves $b$.  So \emph{no} implication between two
sentences that \rga neither proves nor refutes is an \rga theorem: not
$\ConR \limp G$, but not $G \limp G$ or $\ConR \limp \ConR$ either
(part~3).  What fails is the implication, a single sentence; the
entailment $\delta \vdash \delta$ holds at every well-formed $\delta$,
by the hypothesis rule.  A general form of this obstruction is
mechanized as a trap: an implication between two sentences that have
no value at any fuel is not provable, and no sentence built with
$\neg$ and $\lor$ from well-formed such sentences is provable or
refutable ([\citen{conrga-paper}, Appendix ``The traps catalog'']).  Excluded middle fails for the
same reason, and it has no entailment form: its place is taken by the
hypothesis $a \lor \neg a$ that $a$ is decided, which at the Liar is
explosive, $L \lor \neg L \vdash c$ for every well-formed $c$
([\citen{conrga-paper}, Appendix ``The traps catalog'']).

The implication failures below are therefore not about
diagonalization or about provability as such.  Each failure combines the
disjunction property with the independence of the sentences involved:
for $G$ and $\PrC(\bquo{G})$, G\"odel's first theorem; for $\ConR$, and
for the unrefutability of $\PrC(\bquo{\varphi})$ at every well-formed
$\varphi$, G\"odel's second theorem (Section~\ref{sec:goedel}).  Several of these failures are
thus corollaries of G\"odel's theorems for \rga, which Section~\ref{sec:goedel} proves
by other means.  What these failures show is that the classical proof, which
derives implications between independent sentences at every step,
cannot be run in \rga as written.

The entailments need two further tools, both from
[\citen{conrga-paper}, \S\,``Codes, numerals, and four kinds of statement''].  At a hypothesis $E = \exists s.\, D(s)$
with $D$ decided, such as $\PrC(c)$, $\neg\ConR$ or $\Tru(c)$, the
$\Sigma$ deduction theorem makes an entailment from $E$ one sentence
([\citen{conrga-paper}, theorem~\textsf{deduction\_theorem\_exact}](2)), and two of its corollaries replace the
classical moves that need excluded middle
([\citen{conrga-paper}, corollary~\textsf{dec\_contrapos}]).  Reductio: if $E \vdash \neg E$,
then $\emptyset \vdash \neg E$.  Contraposition: if $E \vdash B$, then
$\neg B \vdash \neg E$.  Both corollaries split on the decided $D(s)$ at a fresh
$s$, never on $E$.  An entailment that fails is shown not derivable in
one of two ways.  If its hypothesis is true in the $\omega$-semantics
of [\citen{conrga-paper}, \S\,``The $\omega$-semantics''], in which \rga is sound in context, and its
conclusion is not, it fails; $G$, $\ConR$ and the truth-predicate Liar
are such hypotheses ([\citen{conrga-paper}, theorem~\textsf{goedelW\_omega\_incomplete}]).  Otherwise, reductio
or the $\Sigma$ deduction theorem turns the entailment into a sentence
that G\"odel's theorems exclude.

\paragraph{The classical ingredients, one by one.}
\label{sec:breaks-table}
Table~\ref{tab:breaks} lists each ingredient of the classical proof in
its forms.  Every entry is a theorem of the development.

\begin{table}[tbp]
\centering
\footnotesize
\setlength{\tabcolsep}{2.5pt}
\begin{tabular}{@{}>{\raggedright\arraybackslash}p{0.11\linewidth}>{\raggedright\arraybackslash}p{0.21\linewidth}>{\raggedright\arraybackslash}p{0.21\linewidth}>{\raggedright\arraybackslash}p{0.27\linewidth}>{\raggedright\arraybackslash}p{0.15\linewidth}@{}}
\hline
ingredient & implication / internal sentence & decided-antecedent form & entailment & rule \\
\hline
self-reference & $G \liff \forall q.\,\neg\Chk(q, \bquo{G})$ fails (Thm.~\ref{thm:br-bicond}); the code equation $\vdash S_G = \bquo{G}$ holds (Thm.~\ref{thm:g-diag}) & from $\PrC(\bquo{G}) \vdash \neg G$, by [\citen{conrga-paper}, theorem~\textsf{deduction\_theorem\_exact}](2) & holds: $G \eqv \neg\PrC(\bquo{G})$ and $\neg G \eqv \PrC(\bquo{G})$, and for the Liar, $L_{\Tru} \eqv \neg\Tru(\bquo{L_{\Tru}})$ (Thm.~\ref{thm:br-diag-ent}) & --- \\
D1 & fails at $G$: $\nvdash G \limp \PrC(\bquo{G})$ (Thm.~\ref{thm:br-d2d3}) & not of this shape & fails at $G$: $G \nvdash \PrC(\bquo{G})$ (Thm.~\ref{thm:br-dc-ent}) & holds (Thm.~\ref{thm:br-rules}) \\
D2 & fails at $G$ (Thm.~\ref{thm:br-d2d3}) & decided form with explicit transformer (Thm.~\ref{thm:g-d2}) & holds for all well-formed $\varphi$, $\psi$; EXPL for every $g$ (Thm.~\ref{thm:br-dc-ent}) & holds (Thm.~\ref{thm:br-rules}) \\
D3 & fails at $G$ (Thm.~\ref{thm:br-d2d3}) & decided form with explicit transformer, for every $\varphi$ and over the code (Thms.~\ref{thm:g-d3-any}, \ref{thm:g-d3-uniform}) & holds for every $\varphi$, and as D3G at $G$ (Thm.~\ref{thm:br-dc-ent}) & holds (Thm.~\ref{thm:br-rules}) \\
formalized G\"odel~I & $\ConR \limp G$ and $\ConR \limp \neg\PrC(\bquo{G})$ fail (Thm.~\ref{thm:br-g1}) & contraposed: $\forall q.\,(\Chk(q, \bquo{G}) \limp \neg\ConR)$ holds, also with an explicit transformer (Thm.~\ref{thm:g-sentence}) & holds: $\ConR \vdash G$, $\ConR \vdash \neg\PrC(\bquo{G})$, $\neg G \vdash \neg\ConR$ (Thm.~\ref{thm:br-hyp}); in fact $G \eqv \ConR$ (Thm.~\ref{thm:g-interderivable}) & by cut; vacuous, since $\nvdash \ConR$ \\
formalized G\"odel~II & $\ConR \limp \neg\PrC(\bquo{\ConR})$ and $\PrC(\bquo{\ConR}) \limp \neg\ConR$ fail (Thms.~\ref{thm:g-formal}, \ref{thm:g-sentence}) & $\forall q.\,(\Chk(q, \bquo{\ConR}) \limp \neg\ConR)$ holds, also with an explicit transformer (Thm.~\ref{thm:g-sentence}) & holds: $\PrC(\bquo{\ConR}) \vdash \neg\ConR$ and $\ConR \vdash \neg\PrC(\bquo{\ConR})$ (Thm.~\ref{thm:g-formal}); in fact $\neg\ConR \eqv \PrC(\bquo{\ConR})$ (Thm.~\ref{thm:g-interderivable}) & by cut; vacuous, since $\nvdash \ConR$ \\
L\"ob & fails at $\z = \tofn 1$ (Thm.~\ref{thm:br-loeb}) & holds for every closed well-formed $\psi$ (Thm.~\ref{thm:g-sentence}) & holds for every closed well-formed $\psi$ (Thm.~\ref{thm:br-loeb-ent}) & L\"ob's rule holds (Thm.~\ref{thm:br-loeb}) \\
local reflection at $\z = \tofn 1$ & $\PrC(\bquo{\z = \tofn 1}) \limp \z = \tofn 1$ fails (Thm.~\ref{thm:br-loeb}) & $\forall p.\,(\Chk(p, \bquo{\z = \tofn 1}) \limp \z = \tofn 1)$ fails (Thm.~\ref{thm:br-loeb}) & fails: $\PrC(\bquo{\z = \tofn 1}) \nvdash \z = \tofn 1$; holding would give $\vdash \ConR$ (Thm.~\ref{thm:br-loeb-ent}) & holds (Thm.~\ref{thm:br-rules}) \\
T-bicon\-ditionals & fail at every ungrounded sentence ([\citen{conrga-paper}, theorem~\textsf{tschema\_\allowbreak unprovable\_\allowbreak ungrounded}]); hold at every grounded-true sentence and every grounded-false $\forall$-free one ([\citen{conrga-paper}, \S\,``Truth, the T-schema, and the Liars'']) & --- & T-In fails at $L_{\Tru}$, $G$ and $\ConR$; T-Out fails at $L_{\Tru}$ ([\citen{conrga-paper}, theorem~\textsf{tin\_\allowbreak ent\_\allowbreak fails\_\allowbreak liar}]); the full biconditional at $L_{\Tru}$, as a hypothesis, is explosive ([\citen{conrga-paper}, theorem~\textsf{tbicond\_\allowbreak liar\_\allowbreak explosive}]) & NEC and CONEC hold ([\citen{conrga-paper}, \S\,``Internalized completeness'']); other rules: [\citen{conrga-paper}, the table of the truth predicate against the axioms of truth] \\
excluded middle & fails at $L$: $\nvdash L \lor \neg L$ & --- & no entailment form; $L \lor \neg L$ is an explosive hypothesis (\S\ref{sec:breaks-mech}) & as a scheme, inconsistent (\S\ref{sec:breaks-mech}) \\
$\ConR$ decided & the internal sentence $\ConR\jbool$ is not a theorem (Thm.~\ref{thm:br-con}) & --- & settles nothing: $\ConR\jbool \nvdash \ConR$, $\ConR\jbool \nvdash \neg\ConR$ (Thm.~\ref{thm:br-con}) & $\nvdash \ConR$, $\nvdash \neg\ConR$ (Thms.~\ref{thm:g-g2}, \ref{thm:br-con}) \\
\hline
G\"odel~II & \multicolumn{4}{l}{\textbf{holds}: $\emptyset \nvdash \ConR$ (Thm.~\ref{thm:g-g2}), by the entailments above and one cut} \\
\hline
\end{tabular}
\caption{The ingredients of G\"odel's proof in \rga, in the four forms
  of [\citen{conrga-paper}, \S\,``Codes, numerals, and four kinds of statement''].  The decided-antecedent form applies
  where the hypothesis is a provability statement $\PrC(c)$; ``not of
  this shape'' marks a hypothesis that is not.  The row ``$\ConR$
  decided'' has an internal sentence, not an implication, in its
  second column.  The deduction theorem, which would turn each
  entailment into the implication, holds at a hypothesis exactly when
  \rga proves that hypothesis decided, and fails at $L$ (\S\ref{sec:breaks-mech});
  the entailment $a \vdash a$, which the deduction theorem would turn into $a \limp a$, is
  the hypothesis rule.  In a classical theory such as Peano
  arithmetic, implication and entailment coincide by the deduction
  theorem, excluded middle holds, and the textbook proof uses the
  implications.  Even there, the implications of the rows D1 at $G$,
  local reflection at $\z = \tofn 1$ and the T-biconditionals at the Liar are
  not theorems (for a $\Sigma_1$-sound theory, and by Tarski's
  theorem).  In \rga every implication the proof uses fails, every
  entailment it uses holds, and the theorem holds; the entailments
  that fail are ones the proof does not use.}
\label{tab:breaks}
\end{table}

\subsection{Self-reference survives as an equation and as entailments}
\label{sec:breaks-diag}

The diagonal construction works in \rga exactly.  Section~\ref{sec:goedel}
builds a term $D_G$ with one free variable and a closed term $S_G$ such
that substituting the code numeral of $D_G$ into $D_G$ yields the G\"odel
sentence $G = \forall q.\ \neg\Chk(q, S_G)$, as an identity of terms,
and \rga proves $S_G = \bquo{G}$ (Theorem~\ref{thm:g-diag}).  So $G$ is
literally the statement ``no $q$ is an accepted proof of the term whose
code is $S_G$'', and \rga proves that $S_G$ is $G$'s own code.
The internal-reasoning paper built the Liar for the truth predicate the
same way [\citen{conrga-paper}, \S\,``The truth predicate and the Liar''].  The equation form itself is classical: it is the strong form
of the diagonal lemma~\cite{halbach14self,grabmayr20self}.

What fails is the biconditional that the classical diagonal lemma
derives from the equation, read as an implication.

\begin{thm}[\textsf{no\_goedel\_biconditional}, \textsf{no\_liar\_biconditional}]
\label{thm:br-bicond}
For every term $\delta$ (trivially when $\delta$ is ill-formed),
$\emptyset \nvdash \delta \liff \forall q.\ \neg\Chk(q, \bquo{\delta})$
and
$\emptyset \nvdash \delta \liff \neg\Tru(\bquo{\delta})$.
\end{thm}

\noindent
\emph{Proof idea.}  Suppose \rga proves $\delta \liff \forall q.\
\neg\Chk(q, \bquo{\delta})$.  One conjunct is the disjunction
$\neg\forall q.\ \neg\Chk(q, \bquo{\delta}) \lor \delta$, that is,
$\PrC(\bquo{\delta}) \lor \delta$.  By
the disjunction property \rga proves one disjunct.  If \rga proves
$\PrC(\bquo{\delta})$, the witness property gives a numeral that \rga
proves accepted, and genuineness ([\citen{conrga-paper}, theorem~\textsf{iChk\_faithful}]) makes that numeral a
real proof of $\delta$.  So in either case \rga proves $\delta$, and
faithfulness makes the code of that proof provably accepted.  The other
conjunct, $\neg\delta \lor \forall q.\ \neg\Chk(q, \bquo{\delta})$, then forces \rga
to prove $\neg\delta$ or to refute an accepted code; both contradict
consistency.  The truth-predicate form runs the same way, with
[\citen{conrga-paper}, theorem~\textsf{iWtrue\_of\_wosr}](2), open completeness and
[\citen{conrga-paper}, theorem~\textsf{rga\_sound\_sch}] in place of genuineness and faithfulness.
Neither part of the theorem has a well-formedness hypothesis: the checker form does
not need one, and in the truth form a derivable biconditional is
well-formed.

By \S\ref{sec:breaks-mech}, the theorem says less about
self-reference than it appears to: at an undecided $\delta$, even
$\delta \liff \delta$ is unprovable (Theorem~\ref{thm:br-mech}(3)),
since its conjunct $\neg\delta \lor \delta$ would give a proof or a
refutation of $\delta$ by the disjunction property.  Read as
entailments, the two directions of the diagonal biconditional hold.

\begin{thm}[\textsf{goedel\_\allowbreak diag\_\allowbreak ent}, \textsf{liar\_\allowbreak diag\_\allowbreak ent}]
\label{thm:br-diag-ent}
The diagonal pairs hold as entailments: $G \eqv \neg\PrC(\bquo{G})$,
$G \eqv \forall q.\ \neg\Chk(q, \bquo{G})$ and $\neg G \eqv
\PrC(\bquo{G})$; and $L_{\Tru} \eqv \neg\Tru(\bquo{L_{\Tru}})$.
\end{thm}

\noindent
\emph{Proof idea.}  $G$ is $\forall q.\ \neg\Chk(q, S_G)$, and
$\neg\PrC(\bquo{G})$ is $\neg\neg\forall q.\ \neg\Chk(q, \bquo{G})$.
Each direction transports the checker's verdict along $S_G = \bquo{G}$
under the binder, at an eigenvariable, and adds or removes the double
negation.  $\neg G$ is literally $\PrC(S_G)$, so the third pair is the
same transport of an existential.  For the Liar, the transport of $\neg\Tru$ along
$S_{\Tru} = \bquo{L_{\Tru}}$ does the same.  \qed

\medskip
\noindent
The biconditional fails and the two entailments hold, and the gap
between the biconditional and the entailments is the deduction theorem at the undecided $G$ and
$L_{\Tru}$.  As [\citen{conrga-paper}, \S\,``Codes, numerals, and four kinds of statement''] cautions, $G \eqv
\neg\PrC(\bquo{G})$ is interderivability, not intersubstitutivity.
What is specific to \rga is that self-reference holds as an equation
of codes and as a pair of entailments, but not as a sentence.

\subsection{Formalized G\"odel~I and the consistency sentence}
\label{sec:breaks-g1}

Classically the core of the second theorem is formalized G\"odel~I,
$\vdash \ConR \limp \neg\PrC(\quo{G})$.  In \rga that implication is
not a theorem, in either form.

\begin{thm}[\textsf{no\_formalized\_goedel1}, \textsf{no\_formalized\_goedel1\_code}]
\label{thm:br-g1}
$\emptyset \nvdash \ConR \limp G$ and
$\emptyset \nvdash \ConR \limp \neg\PrC(\bquo{G})$.
\end{thm}

\noindent
The reason is that the antecedent $\ConR$ is undecided, and a hypothesis that
$\ConR$ is decided settles nothing:

\begin{thm}[\textsf{ConRGA\_\allowbreak unrefutable}, \textsf{ConRGA\_\allowbreak undecided}]
\label{thm:br-con}
$\emptyset \nvdash \neg\ConR$, and $\emptyset \nvdash \ConR \jbool$.
Moreover $\ConR\jbool \nvdash \ConR$ and $\ConR\jbool \nvdash
\neg\ConR$.
\end{thm}

\noindent
\emph{Proof idea.}  A refutation of $\ConR = \forall p.\ \neg\Chk(p,
\bquo{\z = \tofn 1})$ would, by the witness property, provide a numeral
that \rga proves accepted as a proof of $\z = \tofn 1$, and genuineness
makes that numeral a real proof of a false equation.  So $\ConR$ is unrefutable;
by G\"odel's second theorem (Theorem~\ref{thm:g-g2}) it is also
unprovable, and by the disjunction property it is not provably decided.
For the entailments: $\neg\ConR$ entails $\ConR\jbool$, so
$\ConR\jbool \vdash \ConR$ would give $\neg\ConR \vdash \neg\neg\ConR$,
and reductio at the hypothesis $\neg\ConR = \PrC(\bquo{\z = \tofn 1})$
would give $\vdash \ConR$; and $\ConR\jbool$ is true in the
$\omega$-semantics while $\neg\ConR$ is not.

The entailment $\ConR \vdash G$, however, holds.

\begin{thm}[\textsf{goedelG\_\allowbreak under\_\allowbreak ConRGA}, \textsf{con\_\allowbreak notprc\_\allowbreak goedelG\_\allowbreak ent}]
\label{thm:br-hyp}
$\ConR \vdash G$ and $\ConR \vdash \neg\PrC(\bquo{G})$; and,
contraposed, $\neg G \vdash \neg\ConR$.
\end{thm}

\noindent
The proof of Theorem~\ref{thm:g-split} runs unchanged in the context
$\{\ConR\}$, since $\ConR$ enters that proof only through
$\forall$-elimination; the proof uses only the checker, with the explicit
witnesses of EXPL and D3G.  The second entailment follows by a cut with the first direction of
Theorem~\ref{thm:br-diag-ent}.  The third entailment holds although neither
$\ConR$ nor $G$ is decided, because $\neg G$ and $\neg\ConR$ are
$\PrC(S_G)$ and $\PrC(\bquo{\z = \tofn 1})$, hypotheses of the $\Sigma$
shape.  The third entailment's route is not a generic contraposition, which fails in
general ([\citen{conrga-paper}, \S\,``Codes, numerals, and four kinds of statement'']), but D3G and EXPL as entailments
(Theorem~\ref{thm:br-dc-ent}): from $\PrC(\bquo{G})$ they give
$\PrC(\bquo{\neg G})$ and then $\PrC(\bquo{\z = \tofn 1})$.  So the gap between $\ConR \vdash G$
and $\emptyset \vdash \ConR \limp G$ is exactly the deduction theorem
at the undecided hypothesis $\ConR$.  By Theorem~\ref{thm:br-con}, that gap
cannot be closed by assuming $\ConR$ decided either.  The rule ``if
$\vdash \ConR$ then $\vdash G$'' follows by cut, and is vacuous, since
$\nvdash \ConR$.  Section~\ref{sec:goedel} uses the entailment to
derive formalized G\"odel~II as an entailment
(Theorem~\ref{thm:g-formal}).

$\ConR$ has no value ([\citen{conrga-paper}, \S\,``Consistency through the truth predicate, and its limit'']), although \rga
derives each of its numeral instances $\neg\Chk(\tofn n, \bquo{\z =
\tofn 1})$.  So $\ConR$ is $\omega$-incomplete, like $G$, and the
$\omega$-semantics makes it true ([\citen{conrga-paper}, theorem~\textsf{goedelW\_omega\_incomplete}]).  On the
grounded reading of [\citen{conrga-paper}, \S\,``The reflection pillars''] this lack of a value is the
expected status of a consistency sentence that quantifies over all of
\rga's own proofs with no restriction.

\subsection{The derivability conditions \rga has}
\label{sec:breaks-dc}

The derivability conditions are where the four forms come apart most
clearly, and where \rga's own conditions can be read off.  Classically,
D2 and D3 are implications, one for each $\varphi$ and $\psi$, whose
antecedents are provability statements $\PrC(\ldots)$.  In \rga those
implications fail already at the G\"odel sentence, and so does D1 as an
implication.

\begin{thm}[\textsf{no\_D2\_schema\_goedelG}, \textsf{no\_D3\_schema\_goedelG}]
\label{thm:br-d2d3}
$\emptyset \nvdash \PrC(\bquo{G \limp \z = \tofn 1}) \limp
(\PrC(\bquo{G}) \limp \PrC(\bquo{\z = \tofn 1}))$, 
$\emptyset \nvdash \PrC(\bquo{G}) \limp \PrC(\bquo{\PrC(\bquo{G})})$,
and $\emptyset \nvdash G \limp \PrC(\bquo{G})$.
\end{thm}

\noindent
\emph{Proof idea.}  By the disjunction property, each implication
needs \rga to refute an antecedent, which for $\PrC(\ldots)$
Theorem~\ref{thm:br-loeb}(1) below excludes (for $\PrC(\bquo{G})$
already G\"odel's first theorem) and for $G$ the first theorem
excludes, or to prove a consequent.
For D3, proving $\PrC(\bquo{\PrC(\bquo{G})})$ would, by the witness
property and genuineness, give a proof of $\PrC(\bquo{G})$, and then of
$G$, against Theorem~\ref{thm:g-g1}.  For D2, the consequent
$\PrC(\bquo{\z = \tofn 1})$ is $\neg\ConR$, which
Theorem~\ref{thm:br-con} excludes.  For D1, the consequent
$\PrC(\bquo{G})$ is excluded as for D3.

The rule forms, by contrast, all hold.

\begin{thm}[\textsf{D1\_rule\_pr}, \textsf{rga\_pr\_reflect}]
\label{thm:br-rules}
For all terms $\varphi$ and $\psi$:
\begin{enumerate}
\item if $\emptyset \vdash \varphi$, then $\emptyset \vdash
  \PrC(\bquo{\varphi})$;
\item if $\emptyset \vdash \PrC(\bquo{\varphi \limp \psi})$ and
  $\emptyset \vdash \PrC(\bquo{\varphi})$, then $\emptyset \vdash
  \PrC(\bquo{\psi})$;
\item if $\emptyset \vdash \PrC(\bquo{\varphi})$, then $\emptyset
  \vdash \PrC(\bquo{\PrC(\bquo{\varphi})})$;
\item if $\emptyset \vdash \PrC(\bquo{\psi})$, then $\emptyset \vdash
  \psi$.
\end{enumerate}
\end{thm}

\noindent
Item~3 is item~1 applied to $\PrC(\bquo{\varphi})$.  Item~2 uses
item~4 twice, modus ponens and item~1.  Item~4, the reflection rule,
holds because \rga proves $\PrC(\bquo{\psi})$ only when $\psi$ is a
theorem, by the witness property and genuineness.
These rules are closure properties of \rga's theorems, proved in HOL;
none of these rules is a single \rga sentence.

Between the failing implications and the HOL rules lie the
entailments.  D2 and D3 hold as entailments; D1 at $G$ does not.

\begin{thm}[\textsf{D2\_\allowbreak ent}, \textsf{D3\_\allowbreak ent}]
\label{thm:br-dc-ent}
\leavevmode
\begin{enumerate}
\item $G \nvdash \PrC(\bquo{G})$.
\item For all well-formed $\varphi$ and $\psi$, $\PrC(\bquo{\varphi}),\
  \PrC(\bquo{\varphi \limp \psi}) \vdash \PrC(\bquo{\psi})$.
\item For every term $\varphi$, $\PrC(\bquo{\varphi}) \vdash
  \PrC(\bquo{\PrC(\bquo{\varphi})})$; and $\PrC(\bquo{G}) \vdash
  \PrC(\bquo{\neg G})$.
\item For every term $g$, $\PrC(\bquo{g}),\ \PrC(\bquo{\neg g}) \vdash
  \PrC(\bquo{\z = \tofn 1})$; and for every well-formed $\varphi$,
  $\PrC(\bquo{\z = \tofn 1}) \vdash \PrC(\bquo{\varphi})$.
\end{enumerate}
Items~2 to~4 hold in every context, not only in the empty one.
\end{thm}

\noindent
\emph{Proof idea.}  (1) $G$ is true in the $\omega$-semantics, and
$\PrC(\bquo{G})$ is not, since \rga proves every instance
$\neg\Chk(\tofn n, \bquo{G})$.  (2)--(4) come from the decided
forms below.  Each of these proofs eliminates the existential of a hypothesis
$\PrC(c)$ at a fresh proof variable $v$, applies the decided
implication to $\Chk(v, c)$, and introduces the existential again at
the transformer's output (\S\ref{app:goedel-ent}).  Item~3's second part is D3G
(Theorem~\ref{thm:g-d3g}); item~4's parts are EXPL and ex
falso.  \qed

\medskip
\noindent
D1 at $G$ fails as an entailment, as it fails classically as an
implication: for a $\Sigma_1$-sound theory such as Peano arithmetic,
$G \limp \PrC(\quo{G})$ with the diagonal equivalence would refute
$G$.  The classical internal counterpart of D1 is provable
$\Sigma_1$-completeness, $\sigma \limp \PrC(\quo{\sigma})$ for
$\Sigma_1$ sentences $\sigma$; in \rga the counterparts of provable $\Sigma_1$-completeness are D3 and
internal completeness, both entailments from a hypothesis of the
$\Sigma$ shape ([\citen{conrga-paper}, \S\,``Internal completeness as one sentence'']).

By the $\Sigma$ deduction theorem, each of the entailments in items~2
to~4 is equivalent to one sentence.  \rga proves stronger, decided
forms of items~2 to~4, internal implications over the checker's
verdict with explicit primitive-recursive transformers: decided D2 and
EXPL (Theorems~\ref{thm:g-expl} and~\ref{thm:g-d2}), decided D3 for
every $\varphi$ and over the code (Theorems~\ref{thm:g-d3-any}
and~\ref{thm:g-d3-uniform}), and ex falso
(Appendix~\ref{app:goedel-cons}).

The development collects the checker's conditions in one theorem
(\S\ref{app:goedel-ent}): D1 as a rule, proved without \rga's
soundness by closed evaluation of the checker; the decided D2 for every pair of
well-formed formulas; the decided D3 as one sentence over the code;
ex falso; and EXPL, each with its transformer named.  This collection is what the
intensional criterion of [\citen{conrga-paper}, \S\,``Compiled primitive-recursive machines''] asks for: \rga proves the
derivability conditions for the very checker that $\ConR$ mentions.

These conditions are the ones that the rest of the paper uses.  EXPL and D3G
carry G\"odel's second theorem (Section~\ref{sec:goedel}); D2 and D3 as
entailments, with the cut transformers of \S\ref{sec:goedel-formal},
carry formalized G\"odel~II and L\"ob's theorem as entailments; ex
falso, with Theorem~\ref{thm:g-g2}, carries Theorem~\ref{thm:br-loeb};
the closed-numeral laws carry every evaluation at a closed code.
Where the classical condition is an implication over $\PrC$, \rga has
an entailment, and a decided implication over the checker with the
transformer written out.  The decided form is the one that
free-variable presentations of the second theorem use throughout, and
the one that \rga's discipline admits as a sentence.

\subsection{L\"ob's theorem and reflection}
\label{sec:breaks-loeb}

L\"ob's theorem says that if $\vdash \PrC(\quo{\varphi}) \limp \varphi$
then $\vdash \varphi$.  The classical proof of L\"ob's theorem needs a fixed point as a
provable biconditional, and D2 and D3 as implications.  The
biconditionals of Theorem~\ref{thm:br-bicond} and D2 and D3 at $G$
fail as implications (Theorems~\ref{thm:br-bicond},
\ref{thm:br-d2d3}), and by \S\ref{sec:breaks-mech} so does any
biconditional or implication between independent sentences.  L\"ob's
rule holds anyway, by the decided ex-falso transformer.

\begin{thm}[\textsf{no\_code\_unprovable}, \textsf{rga\_loeb\_rule}]
\label{thm:br-loeb}
\leavevmode
\begin{enumerate}
\item For every well-formed $\varphi$, $\emptyset \nvdash \forall q.\
  \neg\Chk(q, \bquo{\varphi})$.
\item L\"ob's rule: for every term $\varphi$, if $\emptyset \vdash
  \PrC(\bquo{\varphi}) \limp \varphi$, then $\emptyset \vdash \varphi$.
\item Formalized L\"ob fails at $\z = \tofn 1$ as an implication:
  $\emptyset \nvdash \PrC(\bquo{\PrC(\bquo{\z = \tofn 1}) \limp \z =
  \tofn 1}) \limp \PrC(\bquo{\z = \tofn 1})$.
\item Local reflection at $\z = \tofn 1$ fails in decided-antecedent
  form: $\emptyset \nvdash \forall p.\ (\Chk(p, \bquo{\z = \tofn 1})
  \limp \z = \tofn 1)$.
\end{enumerate}
\end{thm}

\noindent
\emph{Proof idea.}  (1) If \rga proved $\forall q.\ \neg\Chk(q,
\bquo{\varphi})$, it would prove $\ConR$: open a fresh $p$ with $p
\jnat$ and split on the decided $\Chk(p, \bquo{\z = \tofn 1})$; in the
accepting case, ex falso turns $p$ into an accepted proof of $\varphi$,
against the hypothesis instantiated at that code.  Theorem~\ref{thm:g-g2}
excludes $\vdash \ConR$.  So \rga proves the unprovability of no
well-formed sentence.  (2) By the disjunction property, \rga proves
$\neg\PrC(\bquo{\varphi})$, which (1) excludes after double-negation
elimination ($\varphi$ is well-formed, since the implication is
derivable), or proves $\varphi$.  (3) By the disjunction property
again, with (1) for the antecedent and Theorem~\ref{thm:br-con} for the
consequent $\PrC(\bquo{\z = \tofn 1}) = \neg\ConR$.  (4) Open a
fresh $p$ and split on the decided $\Chk(p, \bquo{\z = \tofn 1})$: in
the accepting case the assumed sentence gives $\z = \tofn 1$, which with the
refutation of that equation gives anything, so both cases give $\neg\Chk(p, \bquo{\z =
\tofn 1})$, and $\forall$-introduction gives $\ConR$, which
Theorem~\ref{thm:g-g2} excludes.  \qed

So L\"ob's rule holds for every sentence, by a route that is not
L\"ob's, and the formalized principle fails as an implication at its
simplest instance.  At $\varphi = (\z = \tofn 1)$, item~2 and \rga's
consistency also exclude the implication $\PrC(\bquo{\z = \tofn 1})
\limp \z = \tofn 1$, local reflection at $\z = \tofn 1$.  As
entailments, formalized L\"ob holds and local reflection fails.

\begin{thm}[\textsf{formalized\_\allowbreak loeb\_\allowbreak ent}, \textsf{reflect\_\allowbreak ent\_\allowbreak iff\_\allowbreak con\_\allowbreak ctx}]
\label{thm:br-loeb-ent}
\leavevmode
\begin{enumerate}
\item For every closed well-formed $\psi$,
  $\PrC(\bquo{\PrC(\bquo{\psi}) \limp \psi}) \vdash \PrC(\bquo{\psi})$;
  at $\psi = (\z = \tofn 1)$, $\PrC(\bquo{\neg\ConR \limp \z = \tofn
  1}) \vdash \neg\ConR$.
\item Local reflection at $\z = \tofn 1$ holds in a context exactly
  when consistency does: for every context $\Gamma$,
  $\PrC(\bquo{\z = \tofn 1}), \Gamma \vdash \z = \tofn 1$ if and only
  if $\Gamma \vdash \ConR$.  In particular $\ConR,
  \PrC(\bquo{\z = \tofn 1}) \vdash \z = \tofn 1$.  In the empty context
  local reflection fails, since holding would give $\vdash \ConR$: $\PrC(\bquo{\z =
  \tofn 1}) \nvdash \z = \tofn 1$, $\vz\jnat,\ \Chk(\vz, \bquo{\z =
  \tofn 1}) \nvdash \z = \tofn 1$, and $\Cert_\emptyset(\bquo{\z =
  \tofn 1}) \nvdash \z = \tofn 1$.
\end{enumerate}
\end{thm}

\noindent
\emph{Proof idea.}  (1) is L\"ob's proof read with entailments, with a
cut transformer in place of each implication introduction
(\S\ref{sec:goedel-formal}).  (2) If $\PrC(\bquo{\z = \tofn 1})
\vdash \z = \tofn 1$, explosion gives $\PrC(\bquo{\z = \tofn 1})
\vdash \neg\PrC(\bquo{\z = \tofn 1})$, and reductio at this hypothesis
of the $\Sigma$ shape gives $\vdash \neg\neg\ConR$, hence $\vdash
\ConR$; the same argument runs in every context $\Gamma$.  Conversely,
$\Gamma \vdash \ConR$ refutes the hypothesis.  The empty context is
G\"odel's second theorem.  The checker-premised form reduces to the first form by the $\Sigma$ deduction
theorem.  For the certificate form, reductio over the decided test
$\pi$ refutes $\Cert_\emptyset(\bquo{\z = \tofn 1})$, and since
every accepted proof yields a certificate, that refutation gives
$\vdash \ConR$.  \qed

\medskip
\noindent
Classically too, local reflection at $\z = \tofn 1$ is a rule and not
an implication: a theory that proved $\PrC(\quo{\z = \tofn 1}) \limp
\z = \tofn 1$ would prove its consistency.  The decided-antecedent
form fails as well (Theorem~\ref{thm:br-loeb}(4)), and so do the forms
with the truth predicate, as a single soundness sentence
([\citen{conrga-paper}, theorem~\textsf{rga\_no\_uniform\_soundness}], \S\ref{sec:goedel-cons}) and as
entailments (\S\ref{app:goedel-ent}).

\subsection{Truth, briefly}
\label{sec:breaks-truth}

G\"odel's proof uses no truth predicate, and neither does \rga's
version of the proof (\S\ref{sec:goedel-g2}).  The internal-reasoning
paper reads the truth predicate $\Tru$ in the same four forms
([\citen{conrga-paper}, \S\,``Truth, the T-schema, and the Liars'']).  As an implication, the T-biconditional
fails at every sentence that has no value at any fuel, as the
disjunction property predicts ([\citen{conrga-paper}, theorem~\textsf{tschema\_unprovable\_ungrounded}]), and
T-In fails even as an entailment at the truth-predicate Liar, at $G$
and at $\ConR$ ([\citen{conrga-paper}, theorem~\textsf{tin\_ent\_fails\_liar}]).  For this paper the
relevant fact is the one about $G$.  \rga derives every numeral
instance $\neg\Chk(\tofn n, S_G)$ of the universal inside $G$, but not
$G$ itself, and the $\omega$-semantics makes $G$ true
([\citen{conrga-paper}, theorem~\textsf{goedelW\_omega\_incomplete}]).  So $G$ is $\omega$-incomplete, as
$\ConR$ is (\S\ref{sec:breaks-g1}).

\subsection{What survives}
\label{sec:breaks-survive}

What survives is what Section~\ref{sec:goedel}, on G\"odel's theorems
for \rga, uses: the diagonal as an
equation of codes and as a pair of entailments, the checker's decided
verdict at every number ([\citen{conrga-paper}, theorem~\textsf{hq\_iChk}]), D1 as a rule,
D2 and D3 as entailments, and the decided implications of
\S\ref{sec:breaks-dc} with their explicit transformers.  Every
entailment the textbook proof uses holds, and the entailments from provability
hypotheses hold in decided-antecedent form.  The repair consists in
never asking \rga to discharge an undecided hypothesis into an
implication, or to split on an undecided sentence: where the textbook
proof applies modus ponens to an implication from an undecided
antecedent, \rga cuts against the entailment instead.

\section{G\"odel's theorems for \rga}
\label{sec:goedel}

Section~\ref{sec:breaks} read the classical proof in four forms.
This section proves the theorem: \rga does not derive
\[
  \ConR \;=\; \forall p.\ \neg\Chk(p, \bquo{\z = \tofn 1})
\]
(Theorem~\ref{thm:g-g2}).  The statement uses the tree checker of [\citen{conrga-paper}, \S\,``The proof checker''].  What
makes $\ConR$ a consistency statement, rather than a sentence that
merely has the right extension, is that \rga proves the checker's
structural laws at abstract arguments ([\citen{conrga-paper}, \S\,``The proof checker'']).  The
internal derivations below (EXPL, D3G and the case split) use only
those laws; faithfulness at closed codes enters once, in the G\"odel~I
step (Theorem~\ref{thm:g-g1}(2)).  The theorem is proved in HOL,
which also proves \rga consistent, as any proof that \rga does not
derive a sentence must.  Appendix~\ref{app:goedel} gives the whole
proof at the level of the mechanized lemmas.

The proof is the textbook proof read with entailments.  Where the
textbook proof derives an implication by discharging an undecided
hypothesis,
such as $\ConR \limp G$, \rga derives the entailment, here $\ConR
\vdash G$, and where the textbook proof applies modus ponens to that
implication, \rga cuts.  The conditions that the proof applies to
abstract proof codes, EXPL and D3G, have hypotheses of the form
$\PrC(c) = \exists q.\ \Chk(q, c)$, and by the $\Sigma$ deduction
theorem of [\citen{conrga-paper}, \S\,``Codes, numerals, and four kinds of statement''] an entailment from such a
hypothesis is equivalent to one sentence, a decided implication over
the checker's verdict, with an existential consequent.  \rga proves
stronger sentences, with the transformer written out in place of that
existential.  The decided forms of this section are therefore these
entailments written as sentences, with their witnesses made explicit.

The route works with the proof predicate rather than with
provability.  Hilbert and Bernays already restrict the third
derivability condition to primitive-recursive
antecedents~\cite{hilbert39grundlagen,kurahashi20note}; the decided
form used here is the one that free-variable arithmetic gives the
theorem~\cite{coquand26auto}: the derivability conditions are explicit
primitive-recursive transformations of proofs with decided
antecedents, and the G\"odel formula is derived from the consistency
formula by propositional reasoning on decidable formulas followed by
generalization.  That route discharges a hypothesis, by implication
introduction or by case analysis, only when the hypothesis is decided,
so that route is one that \rga's discipline admits
(Section~\ref{sec:breaks}).  The proof below carries that route out
inside \rga, in four mechanized pieces.  Self-reference is an equation
of codes (\S\ref{sec:goedel-sentences}).  Two derivability facts,
EXPL and D3G, are proof transformers with decided antecedents and
explicit primitive-recursive witnesses (\S\ref{sec:goedel-expl},
\S\ref{sec:goedel-d3g}).  A case split on a decided sentence, under
$\forall$-introduction, derives $G$ from $\ConR$ without the
implication $\ConR \limp G$ (\S\ref{sec:goedel-split}).  No step
discharges an undecided hypothesis.  After G\"odel's first and second
theorems (\S\ref{sec:goedel-g1}, \S\ref{sec:goedel-g2}),
\S\ref{sec:goedel-formal} derives the formalized theorems inside \rga,
as entailments and in decided-antecedent form: formalized G\"odel~II and
L\"ob's theorem, with explicit cut transformers, and formalized
G\"odel~II as one sentence with its transformer named.

How this route relates to earlier proofs of the second theorem,
Beklemishev and Shamkanov's abstract theorem among them, and what is
new in this route, is set out in Sections~\ref{sec:intro}
and~\ref{sec:related}.  What is specific to \rga is where the work
lies.  In free-variable arithmetic every formula is decided, so the
restriction to decided antecedents costs nothing.  \rga has quantifiers
and undecided sentences, and its universal elimination and induction
demand that every term used as a number be shown to be one.  The
substantial piece is therefore the premise certificate of D3G
(\S\ref{sec:goedel-d3g}): an explicit builder that, from the checker's
acceptance of a proof of $G$ at an abstract $v$, produces inside \rga
an accepted proof that the corresponding instance of $G$'s matrix
fails.

\subsection{The sentences}
\label{sec:goedel-sentences}

Let $\Diag$ be the compiled diagonal function: at the code of a
well-formed term $t$ and a numeral $\tofn k$, it computes the code of
$t$ with $\tofn k$ substituted for variable $0$.  The development
builds $\Diag$ from an explicit index tree, and \rga proves the values of $\Diag$
(Appendix~\ref{app:goedel-diag}).  Define
\[
  D_G \;\ldef\; \forall q.\ \neg\Chk(q, \Diag(x, x)) , \qquad
  S_G \;\ldef\; \Diag(\bquo{D_G}, \bquo{D_G}) , \qquad
  G \;\ldef\; \forall q.\ \neg\Chk(q, S_G) ,
\]
where $x$ is the one free variable of $D_G$.

\begin{thm}[\textsf{goedel\_diag}, \textsf{goedelS\_eq}]
\label{thm:g-diag}
$\ssub{x}{\bquo{D_G}}{D_G} = G$, as an identity of terms, and
$\emptyset \vdash S_G = \bquo{G}$.
\end{thm}

\noindent
So $\neg G$ is literally $\PrC(S_G)$, and \rga proves that $S_G$ is
$G$'s code.  Isabelle never computes the large numeral: the equation
is an instance of the closed-numeral law of [\citen{conrga-paper}, theorem~\textsf{cmp\_eq\_num}].  A
transport rule moves the checker's verdict along this equation, in the
checker's conclusion slot and in its proof slot (Appendix~\ref{app:goedel}).

\subsection{Explosion as a proof transformer}
\label{sec:goedel-expl}

\begin{thm}[\textsf{goedelE\_chk}, \textsf{EXPL\_goedelE}]
\label{thm:g-expl}
There is a closed term $e$, compiled from an explicit index tree, such
that for every context $\Gamma$, every term $g$ and all terms $a, b$
with $\Gamma \vdash a \jnat$ and $\Gamma \vdash b \jnat$: if $\Gamma
\vdash \Chk(a, \bquo{g})$ and $\Gamma \vdash \Chk(b, \bquo{\neg g})$,
then $\Gamma \vdash \Chk(e\cdot a\cdot b, \bquo{\z = \tofn 1})$.
Moreover \rga proves $e$ total and proves, for every term $g$, the
internal sentence
\[
  \forall a.\ \forall b.\ \bigl(\Chk(a, \bquo{g}) \wedge \Chk(b,
  \bquo{\neg g})\bigr) \limp \Chk(e\cdot a\cdot b, \bquo{\z = \tofn 1}) .
\]
\end{thm}

\noindent
The term $e\cdot a\cdot b$ appends the two proof lists and puts one new
entry in front: the judgment $\emptyset \vdash \z = \tofn 1$, justified
by the negation-elimination rule, which cites the heads of the two old
lists.  The checker accepts the result because each old head sits at a
position the new entry can name.  The internal sentence discharges its
antecedent legitimately, since the checker's verdicts at numbers are
decided.

EXPL is the decided form of D2 at the instance the proof uses,
explosion.  Modus
ponens on accepted proofs has the same decided form at every pair of
formulas, although the second theorem does not need that form.

\begin{thm}[\textsf{D2\_decided}, \textsf{d2M\_chk}]
\label{thm:g-d2}
For all well-formed terms $\varphi$ and $\psi$ there is a closed term
$m_{\varphi\psi}$, compiled from an explicit index tree, such that
\rga proves
\[
  \forall a.\ \forall b.\ \bigl(\Chk(a, \bquo{\varphi \limp \psi})
  \wedge \Chk(b, \bquo{\varphi})\bigr) \limp \Chk(m_{\varphi\psi}\cdot
  a\cdot b, \bquo{\psi}) .
\]
The rule form of this sentence holds in every context: if $\Gamma \vdash a\jnat$,
$\Gamma \vdash b\jnat$, $\Gamma \vdash \Chk(a, \bquo{\varphi \limp
\psi})$ and $\Gamma \vdash \Chk(b, \bquo{\varphi})$, then $\Gamma
\vdash \Chk(m_{\varphi\psi}\cdot a\cdot b, \bquo{\psi})$.
\end{thm}

\noindent
The theorem is a schema: one internal sentence for each pair of formulas.  The
proof list $m_{\varphi\psi}\cdot a\cdot b$ puts three entries above a
fixed six-entry proof of $\psi$ from the hypotheses $\neg\varphi \lor
\psi$ and $\varphi$, and the two given proofs: two cuts and one
weakening.  The fixed proof is accepted by closed evaluation of the
checker, so neither statement uses \rga's soundness.

\subsection{D3G with an explicit primitive-recursive witness}
\label{sec:goedel-d3g}

Classical D3 says that the theory proves, for every $\varphi$, that a
proof of $\varphi$ yields a proof that $\varphi$ is provable.  \rga
proves a decided form of D3 for every $\varphi$
(Theorem~\ref{thm:g-d3-any} below).  The second theorem needs only the
same construction at $G$'s diagonal code term $S_G$, which is D3G.

\begin{thm}[\textsf{D3G\_goedelSigma}, \textsf{goedelSigma\_is\_pr}]
\label{thm:g-d3g}
There is a closed term $\sigma = \lambda v.\ \cmp(\tau_\sigma)\cdot
v\cdot v$, where $\tau_\sigma$ is an explicit index tree, such that
\rga proves $\sigma$ total and proves
\[
  \forall v.\ \Chk(v, \bquo{G}) \limp \Chk(\sigma\cdot v, \bquo{\neg G}) .
\]
\end{thm}

\noindent
The argument $v$ is doubled because index trees denote binary
functions ([\citen{conrga-paper}, \S\,``Compiled primitive-recursive machines'']): $\cmp(\tau_\sigma)\cdot v\cdot v$
feeds the one input $v$ to both argument slots of the binary function
that $\tau_\sigma$ denotes.  The same device recurs in the
transformers below, such as $c_Q\cdot a\cdot a$ of
Theorem~\ref{thm:g-cut}.  That $\tau_\sigma$ denotes a
primitive-recursive function holds of every index tree, including one
chosen by Hilbert choice; what shows that $\sigma$ is built without
choice is the build-time audit of [\citen{conrga-paper}, \S\,``Machines without opaque indices''], which
checks the definitional closure of $\tau_\sigma$.

The witness has a simple shape.  For a number $v$, $\sigma\cdot v$ is
a proof list whose head is $\emptyset \vdash \neg G$, justified by the
rule that refutes a universal through one refuted instance.  Below the
head come two appended proofs: one of $\neg(F\cdot \tofn v)$, where
$F$ is $G$'s matrix, and one of $\tofn v \jnat$.  So $\sigma$ turns a
proof $v$ of $G$ into a proof of $\neg G$ whose counterexample is $v$
itself.

The work lies in the proof of $\neg(F\cdot \tofn v)$.  From $\Gamma \vdash
\Chk(v, S_G)$ at an abstract $v$, \rga must produce an accepted proof
of $\neg(F\cdot\tofn v)$, the refuted instance of $G$'s matrix at
$v$'s numeral, whose content is that $\Chk(\tofn v, S_G)$ holds.  The verdict $\Chk(v, S_G)$ is a conjunction
of four equations and disequations between compiled
primitive-recursive applications.  A
certificate builder, defined by HOL recursion on index trees, produces
for every index $f$ an accepted proof that the compiled application
$\cmp(f)\cdot a\cdot b$ equals the numeral of its value
(Appendix~\ref{app:goedel-d3g}).  Each recursion node of the tree costs one
internal induction.  Four such certificates, combined with conjunction,
double-negation and false-equation steps, give the instance.  This construction is
formalized $\Sigma_1$-completeness for the decided family $\Chk(v,
S_G)$, with its witness written out as an index tree.

The construction depends on $G$ only through the certificate for $S_G$
and through $G$'s matrix, so it works for every sentence.

\begin{thm}[\textsf{D3\_decided\_any}, \textsf{D3\_at}]
\label{thm:g-d3-any}
For every term $\varphi$, well-formed or not, there is a closed term $\sigma_\varphi =
\lambda v.\ \cmp(\tau_\varphi)\cdot v\cdot v$, where $\tau_\varphi$ is
an explicit index tree, such that \rga proves
\[
  \forall v.\ \Chk(v, \bquo{\varphi}) \limp \Chk(\sigma_\varphi\cdot
  v, \bquo{\PrC(\bquo{\varphi})}) .
\]
The same construction at the code term $S_G$ in place of
$\bquo{\varphi}$ yields the witness $\sigma$ of
Theorem~\ref{thm:g-d3g}, and the resulting sentence is D3G.
\end{thm}

\noindent
So D3G is the decided D3 at the diagonal code term $S_G$, with the
second theorem's witness.  D3G is not the instance $\varphi \ldef G$ of
the sentence of Theorem~\ref{thm:g-d3-any}.  The antecedents agree, since $\Chk(v, S_G)$ and
$\Chk(v, \bquo{G})$ are interchangeable by transport along $S_G =
\bquo{G}$.  The consequents differ: the instance at $G$ concludes
$\Chk(\sigma_G\cdot v, \bquo{\PrC(\bquo{G})})$, whose code contains
the code numeral $\bquo{G}$, while D3G concludes $\Chk(\sigma\cdot v,
\bquo{\neg G})$, and $\neg G$ is $\PrC(S_G)$, whose code contains the
diagonal term $S_G$ instead.  The transformers differ accordingly.  The
transformer of Theorem~\ref{thm:g-d3-any} is a HOL function of
$\varphi$, so the theorem is a schema: one internal sentence for each
$\varphi$ (Appendix~\ref{app:goedel-d3any}).  The decided D3 is also
one internal sentence quantified over the code.

\begin{thm}[\textsf{D3\_uniform}, \textsf{D3\_uniform\_num}]
\label{thm:g-d3-uniform}
There is a closed term $\sigma_{\mathrm u} = \cmp(\tau_{\mathrm u})$, where $\tau_{\mathrm u}$
is one explicit index tree, and a compiled function $\mathit{pr}$ such that
\begin{enumerate}
\item $\emptyset \vdash \forall c.\ \forall v.\ \Chk(v, c) \limp
  \Chk(\sigma_{\mathrm u}\cdot c\cdot v, \mathit{pr}(c))$;
\item for every number $k$, $\emptyset \vdash \mathit{pr}(\tofn k) =
  \bquo{\PrC(\tofn k)}$, and hence $\emptyset \vdash \forall v.\
  \Chk(v, \tofn k) \limp \Chk(\sigma_{\mathrm u}\cdot \tofn k\cdot v,
  \bquo{\PrC(\tofn k)})$.
\end{enumerate}
\end{thm}

\noindent
$\mathit{pr}(c)$ computes, like $\Diag$, the code of $\PrC(x)$ with the
numeral of $c$ in place of $x$, as $\mathit{tr}(x)$ of [\citen{conrga-paper}, Appendix ``$\Tru$ against the axioms of truth''] computes the code of $\Tru$ applied to a numeral.  Part~1 is one \rga sentence, and that sentence's
implication is introduced over the decided verdict $\Chk(v, c)$.  The
construction is the premise certificate of
Theorem~\ref{thm:g-d3-any} with the code a second argument of the
transformer instead of a HOL parameter.  Part~2 is a schema over $k$;
at $k = \bquo{\varphi}$ the schema is the decided D3 at $\varphi$, with the
transformer $\sigma_{\mathrm u}\cdot\bquo{\varphi}$ in place of $\sigma_\varphi$.
Neither part uses \rga's soundness;
Appendix~\ref{app:goedel-d3-uniform} gives the construction.

\subsection{The decided case split}
\label{sec:goedel-split}

The classical proof derives $G$ from $\ConR$ through the implication
$\ConR \limp \neg\PrC(\quo{G})$, which is not an \rga theorem
(Theorem~\ref{thm:br-g1}).  \rga derives $G$ from $\ConR$ differently.  Say that
\emph{EXPL holds of} $e$ when $e$ is closed, \rga proves $e$ total, and
\rga proves the internal sentence of Theorem~\ref{thm:g-expl}; and
that \emph{D3G holds of} $\sigma$ likewise, with
Theorem~\ref{thm:g-d3g}.

\begin{thm}[\textsf{goedelG\_of\_expl\_d3g\_ctx}, \textsf{goedelG\_under\_ConRGA}]
\label{thm:g-split}
If EXPL holds of $e$ and D3G holds of $\sigma$, then for every context
$\Gamma$: if $\Gamma \vdash \ConR$, then $\Gamma \vdash G$.  In
particular, with $\Gamma = \{\ConR\}$, $\ConR \vdash G$.
\end{thm}

\noindent
\emph{Proof.}  Open a fresh variable $v$ with $v \jnat$.  The verdict
$\Chk(v, S_G)$ is decided ([\citen{conrga-paper}, theorem~\textsf{hq\_iChk}]), so \rga may
split on that verdict.
\begin{itemize}
\item If $\neg\Chk(v, S_G)$, that is the goal.
\item If $\Chk(v, S_G)$, transport along $S_G = \bquo{G}$ makes $v$ an
  accepted proof of $G$.  D3G at $v$ gives the accepted proof
  $\sigma\cdot v$ of $\neg G$.  EXPL at $v$ and $\sigma\cdot v$ gives
  the accepted proof $e\cdot v\cdot(\sigma\cdot v)$ of $\z = \tofn 1$.
  $\ConR$, instantiated at that code, says that this code is not accepted.  So this
  case is contradictory, and in particular yields $\neg\Chk(v, S_G)$.
\end{itemize}
Both cases give $\neg\Chk(v, S_G)$, and $\forall$-introduction gives
$G$.  \qed

\medskip
\noindent
The proof introduces no implication (Appendix~\ref{app:goedel-assembly}
counts the proof's rules).  The two case hypotheses are discharged by
disjunction elimination, and the typing hypothesis $v \jnat$ by
$\forall$-introduction.  $\ConR$ enters only through
$\forall$-elimination at one term, never as a hypothesis to be
discharged.

\emph{Classical versus \rga.}  In the textbook proof this step is
``$\ConR \limp \neg\PrC(\quo{G})$, then $G$ by the diagonal
biconditional'': an internal implication over the undecided $\ConR$,
followed by an equivalence that \rga does not have
(Theorem~\ref{thm:br-bicond}).  In free-variable arithmetic the step
is already the case split: modus tollens on the decidable proof
predicate at a free variable, then generalization.  In \rga the case split is the only form
available, and that form is available because $G$'s matrix is a
primitive-recursive check, decided at every number.  Read with
entailments, the derivation is the entailment $\ConR \vdash G$ of
Theorem~\ref{thm:br-hyp}, since the derivation runs unchanged in the context
$\{\ConR\}$.  Only the derivation's discharge into an implication
fails.  The instance of Theorem~\ref{thm:g-split} at the empty context, ``if $\emptyset \vdash \ConR$ then
$\emptyset \vdash G$'', is vacuous
as a statement about derivability, since $\nvdash \ConR$; the content
is the entailment.

The decided split under $\forall$-introduction is the pattern of the
$\Sigma$ reductio of [\citen{conrga-paper}, corollary~\textsf{dec\_contrapos}]: split on the decided matrix at a fresh
eigenvariable, never on the undecided hypothesis.  So the decided
route and the textbook chain with $\Sigma$-contraposition are one
argument.  The development also derives $\ConR \vdash G$ literally by
the textbook chain read with entailments
(\S\ref{app:goedel-ent}).  Formalized G\"odel~I
contraposed, $\neg G \vdash \neg\ConR$ (Theorem~\ref{thm:br-hyp}),
has the hypothesis $\neg G = \PrC(S_G)$, of the $\Sigma$ shape, so
$\Sigma$-contraposition gives $\neg\neg\ConR \vdash \neg\neg G$, and
double negation in both places gives $\ConR \vdash G$.  This
derivation is independent of the decided split, and a proof of $\ConR$
cut against this derivation gives a proof of $G$, against G\"odel's
first theorem.  That argument proves Theorem~\ref{thm:g-g2} below
along the textbook chain.

\subsection{G\"odel's first theorem}
\label{sec:goedel-g1}

\begin{thm}[\textsf{rga\_goedel1}, \textsf{goedelG\_refutes}]
\label{thm:g-g1}
\leavevmode
\begin{enumerate}
\item $\emptyset \nvdash G$ and $\emptyset \nvdash \neg G$.
\item If $\emptyset \vdash G$, then $\emptyset \vdash \z = \tofn 1$.
\item \rga is $\omega$-consistent: for no term $u$ does \rga derive
  $\neg\forall x.\,u(x)$ together with $u(\tofn n)$ for every $n$.
\end{enumerate}
\end{thm}

\noindent
\emph{Proof.}  (2) Take a proof $P$ of $G$.  D1 for the checker gives
$\emptyset \vdash \Chk(\bquo{P}, \bquo{G})$; transport along $S_G =
\bquo{G}$, and $G$ instantiated at $\bquo{P}$, then give $\emptyset
\vdash \z = \tofn 1$.  (1) The first half is~(2) with \rga's
consistency.  The second half needs no Rosser sentence.  $\neg G$ is
$\neg\forall q.\ \neg\Chk(q, S_G)$, so a refutation of $G$ would, by
the witness property, give a numeral that \rga proves accepted as a
proof of $G$, and genuineness would make that numeral a real proof of $G$,
contradicting the first half.  (3) The witness property turns a
refutation of $\forall x.\,u(x)$ into a refuted instance $u(\tofn n)$,
which consistency excludes.  \qed

\medskip
\noindent
Part~2 is a transformation of derivations that does not use \rga's
consistency: G\"odel~I as a recipe.  The witness property, a
consequence of \rga's soundness theorem, plays the role of
$\Sigma_1$-soundness, the part of G\"odel's $\omega$-consistency
hypothesis that $\nvdash \neg G$ needs; with consistency the witness property gives
part~3.  So the plain G\"odel sentence is independent.  By
[\citen{conrga-paper}, theorem~\textsf{goedelW\_omega\_incomplete}], $G$ is moreover $\omega$-incomplete, and
true in the classical $\omega$-evaluation.  $G$ has no value at any
fuel (Appendix~\ref{app:goedel-cons}): it is ungrounded.  In \rga's
own semantics $G$ is therefore a gap, neither grounded-true nor
grounded-false, consistent with grounded completeness
([\citen{conrga-paper}, \S\,``Internalized completeness'']): completeness concerns grounded-true
terms, and $G$ is not one, though it is true classically.

\subsection{G\"odel's second theorem}
\label{sec:goedel-g2}

\begin{thm}[\textsf{rga\_goedel2}]
\label{thm:g-g2}
$\emptyset \nvdash \ConR$.
\end{thm}

\noindent
\emph{Proof.}  By the first half of Theorem~\ref{thm:g-g1}, \rga does
not prove $G$.  By Theorem~\ref{thm:g-split}, with $e$ and $\sigma$
from Theorems~\ref{thm:g-expl} and~\ref{thm:g-d3g}, a proof of $\ConR$
would give one of $G$.  So \rga does not prove $\ConR$.  \qed

\medskip
\noindent
In entailment form the proof is two cuts.  A derivation of $\ConR$, cut
against $\ConR \vdash G$, gives a derivation of $G$, and
Theorem~\ref{thm:g-g1}(2) turns that derivation into a derivation of $\z = \tofn
1$.  Only the last step, that \rga does not derive $\z = \tofn 1$, uses
\rga's consistency.

A second proof of the same statement
(Appendix~\ref{app:goedel-search}) differs only in the D3G witness: a
search for an accepted proof of $\neg G$, proved total by an internal
induction from a certificate that the decider's internal completeness
supplies ([\citen{conrga-paper}, theorem~\textsf{rga\_complete\_int2}], in its rule form).  Extensionally
the search is the zero function, but \rga does not prove this extensional fact.  The
headline proof uses neither the decider nor the truth predicate; the
headline proof is built entirely at the checker, and its dependency cone is the smaller
of the two proofs' cones ([\citen{conrga-paper}, Section ``Mechanization'']).  The internal theorems of
\S\ref{sec:goedel-formal} are built at the checker too.

\subsection{The formalized theorems, as entailments}
\label{sec:goedel-formal}

Classically, the second theorem is also a theorem of the theory
itself: Peano arithmetic proves the implication $\ConR \limp
\neg\PrC(\quo{\ConR})$, formalized G\"odel~II, and more generally
formalized L\"ob's theorem.  In \rga these implications fail
(Theorem~\ref{thm:g-formal} below and Theorem~\ref{thm:br-loeb}(3)).
The entailment forms of these implications hold, and so does the
entailment form of formalized G\"odel~I (Theorem~\ref{thm:br-hyp}).  The tool is necessitation for
entailments: every derivation under one or two hypotheses is mirrored
inside \rga by an explicit primitive-recursive transformer of proof
codes, a cut transformer.  Appendix~\ref{app:goedel-ent} gives the
constructions at the level of the mechanized lemmas.

\begin{thm}[\textsf{cut1M\_\allowbreak chk}, \textsf{ent\_\allowbreak nec}]
\label{thm:g-cut}
Let $Q$ be a proof of $\{A\} \vdash B$.  There is a closed term $c_Q$,
compiled from an explicit index tree, such that for every context
$\Gamma$ and every term $a$: if $\Gamma \vdash a\jnat$ and $\Gamma
\vdash \Chk(a, \bquo{A})$, then $\Gamma \vdash \Chk(c_Q\cdot a\cdot a,
\bquo{B})$.  Likewise, for a proof $Q$ of $\{A_1, A_2\} \vdash B$
there is a closed $c_Q$ with $\Gamma \vdash \Chk(c_Q\cdot a\cdot b,
\bquo{B})$ whenever $\Gamma$ derives $a\jnat$, $b\jnat$, $\Chk(a,
\bquo{A_1})$ and $\Chk(b, \bquo{A_2})$.  Consequently, for every
context $\Gamma$:
\begin{enumerate}
\item if $A \vdash B$, then $\PrC(\bquo{A}),\ \Gamma \vdash
  \PrC(\bquo{B})$;
\item if $A_1, A_2 \vdash B$, then $\PrC(\bquo{A_1}),\
  \PrC(\bquo{A_2}),\ \Gamma \vdash \PrC(\bquo{B})$.
\end{enumerate}
\end{thm}

\noindent
The transformer has the shape of EXPL's.  The proof list
$c_Q\cdot a\cdot a$ appends the given proof $a$, whose head is
$\emptyset \vdash A$, to the code of the fixed proof $Q$, whose head is
$\{A\} \vdash B$, and puts one new entry in front: $\emptyset \vdash
B$, justified by cut, which cites the two heads.  Closed evaluation,
without \rga's soundness, shows that the code of $Q$ is a valid tree
whose head is $\{A\} \vdash B$; the checker itself does not accept it,
since that head has a hypothesis.  The new part, relative to EXPL, is reading a head whose hypothesis set is
not empty.  The two-hypothesis transformer adds three entries above
the two given proofs and $Q$: two cuts and one weakening.  For the
consequences, open a fresh variable for the proof of each hypothesis
and introduce the existential again at the transformer's output, as
for the entailments of Theorem~\ref{thm:br-dc-ent}.  The index tree of
$c_Q$ is a HOL function of $Q$, so the theorem is a schema, one
transformer for each proof.  The corollaries obtain $Q$ from the
derivation in HOL; the tree itself contains no choice, which the
build-time check of [\citen{conrga-paper}, \S\,``Machines without opaque indices''] confirms.  ``Explicit'' means choice-free, not small: the
tree of $c_Q$ embeds the codes of $Q$'s formulas as unary numerals
([\citen{conrga-paper}, \S\,``Compiled primitive-recursive machines'']), so the tree is astronomically large.  None of the
four statements uses \rga's soundness.

\begin{thm}[\textsf{rga\_\allowbreak goedel2\_\allowbreak internal\_\allowbreak contra}, \textsf{rga\_\allowbreak goedel2\_\allowbreak internal}]
\label{thm:g-formal}
\leavevmode
\begin{enumerate}
\item $\PrC(\bquo{\ConR}) \vdash \neg\ConR$, that is,
  $\PrC(\bquo{\ConR}) \vdash \PrC(\bquo{\z = \tofn 1})$.
\item $\ConR \vdash \neg\PrC(\bquo{\ConR})$.
\item $\emptyset \nvdash \ConR \limp \neg\PrC(\bquo{\ConR})$.
\item \mbox{\rga{}${}+\ConR$} is sound for the $\omega$-semantics: if
  $\Gamma \vdash c$, $\Gamma$ is well-formed and each member of
  $\Gamma$ is $\ConR$ or $\omega$-true, then $c$ is $\omega$-true.  So
  $\ConR \nvdash \z = \tofn 1$, $G \nvdash \z = \tofn 1$ and $\ConR,
  G \nvdash \z = \tofn 1$.
\end{enumerate}
\end{thm}

\noindent
\emph{Proof.}  (1) Necessitation (Theorem~\ref{thm:g-cut}) at $\ConR
\vdash G$ (Theorem~\ref{thm:br-hyp}) gives $\PrC(\bquo{\ConR}) \vdash
\PrC(\bquo{G})$.  D3G as an entailment gives $\PrC(\bquo{\neg G})$
from the same hypothesis, and EXPL as an entailment gives
$\PrC(\bquo{\z = \tofn 1})$, which is $\neg\ConR$
(Theorem~\ref{thm:br-dc-ent}(3, 4)).  (2) The hypothesis
$\PrC(\bquo{\ConR})$ of~(1) has the $\Sigma$ shape, so contraposition
gives $\neg\neg\ConR \vdash \neg\PrC(\bquo{\ConR})$
([\citen{conrga-paper}, corollary~\textsf{dec\_contrapos}]), and a cut with $\ConR \vdash
\neg\neg\ConR$ gives~(2).  (3) By the disjunction property, \rga would
have to refute $\ConR$, which Theorem~\ref{thm:br-con} excludes, or
prove $\neg\PrC(\bquo{\ConR})$, which Theorem~\ref{thm:br-loeb}(1)
excludes.  (4) \rga's soundness for the $\omega$-semantics, with the
$\omega$-truth of $\ConR$ and of $G$ ([\citen{conrga-paper}, theorem~\textsf{goedelW\_omega\_incomplete}]).  \qed

\medskip
\noindent
Parts~1 and~2, with the decided-antecedent form of
Theorem~\ref{thm:g-sentence} below, are what this paper means by
``\rga derives its own second incompleteness theorem''.  Part~1 is the
primary entailment form, and reads as a transformer: inside \rga, an
accepted proof code of $\ConR$ is turned into an accepted proof code
of $\z = \tofn 1$.  Part~1 is one \rga derivation, from the
hypothesis that some number is an accepted proof of $\ConR$, of the
sentence that some number is an accepted proof of $\z = \tofn 1$; by
cut part~1 turns any derivation of $\PrC(\bquo{\ConR})$, in any context,
into a derivation of $\neg\ConR$.  The derivation composes explicit
transformers: from an accepted proof code of $\ConR$, the cut
transformer of Theorem~\ref{thm:g-cut} builds one of $G$, D3G's
$\sigma$ one of $\neg G$, and EXPL's $e$ one of $\z = \tofn 1$.  Part~2 is the
contrapositive of part~1, the textbook statement: under the hypothesis $\ConR$,
\rga derives $\neg\PrC(\bquo{\ConR})$, ``no number is an accepted
proof of $\ConR$''.  Neither derivation is Theorem~\ref{thm:g-g2},
which is a fact of HOL about \rga; parts~1 and~2 are reasoning inside
\rga, the counterpart of Peano arithmetic's proof of $\ConR \limp
\neg\PrC(\quo{\ConR})$ in the form that \rga admits.  Both parts are built
at the checker: they use neither \rga's soundness, nor the decider,
nor the truth predicate, and their transformers contain no choice
([\citen{conrga-paper}, Appendix ``The build-time check'']).  Part~3 says that the form matters:
the implication form fails, for the reason of \S\ref{sec:breaks-mech}.

What does part~2 say, given that its hypothesis $\ConR$ is never
grounded-true ([\citen{conrga-paper}, \S\,``Consistency through the truth predicate, and its limit''])?  Read as a statement about
truth, the entailment holds vacuously.  The entailment's content is derivational:
the entailment is a theorem of the theory \mbox{\rga{}${}+\ConR$}, and that theory is
consistent and sound for the $\omega$-semantics, in which $\ConR$ is
true (part~4).  Formalized G\"odel~II is thus a theorem of a
consistent, $\omega$-sound extension of \rga, as it is of Peano
arithmetic.

The entailments from provability hypotheses are also single sentences.
By [\citen{conrga-paper}, theorem~\textsf{deduction\_theorem\_exact}](2), an entailment from $\PrC(c)$ to a
closed conclusion is equivalent to a sentence in decided-antecedent
form over proof codes, and \rga proves the sentences of the formalized
theorems, with their transformers named.

\begin{thm}[\textsf{rga\_goedel2\_sentence}, \textsf{goedel2T\_sentence}]
\label{thm:g-sentence}
\leavevmode
\begin{enumerate}
\item $\emptyset \vdash \forall q.\,(\Chk(q, \bquo{\ConR}) \limp
  \neg\ConR)$.
\item For every proof $Q$ of $\{\ConR\} \vdash G$ there is a closed
  term $T_Q$, compiled from an explicit index tree, with $\emptyset
  \vdash \forall q.\,(\Chk(q, \bquo{\ConR}) \limp \Chk(T_Q\cdot q,
  \bquo{\z = \tofn 1}))$; so some closed $T$ satisfies this sentence.
\item $\emptyset \vdash \forall q.\,(\Chk(q, \bquo{G}) \limp
  \neg\ConR)$, and for a closed term $T_1$ compiled from an explicit
  index tree, $\emptyset \vdash \forall q.\,(\Chk(q, \bquo{G}) \limp
  \Chk(T_1\cdot q, \bquo{\z = \tofn 1}))$.
\item For every closed well-formed $\psi$, $\emptyset \vdash \forall
  q.\,(\Chk(q, \bquo{\PrC(\bquo{\psi}) \limp \psi}) \limp
  \PrC(\bquo{\psi}))$.
\item $\emptyset \nvdash \PrC(\bquo{\ConR}) \limp \neg\ConR$, although
  part~1 holds: the quantifier of a sentence in decided-antecedent form
  cannot in general be extracted into an existential antecedent.
\end{enumerate}
\end{thm}

\noindent
\emph{Proof idea.}  Parts~1, 3 and~4 are the decided-antecedent form
([\citen{conrga-paper}, theorem~\textsf{deduction\_theorem\_exact}](2)) of the entailments of
Theorems~\ref{thm:g-formal}(1), \ref{thm:br-hyp} (contraposed, with
the diagonal pair of Theorem~\ref{thm:br-diag-ent}) and
\ref{thm:br-loeb-ent}(1).  Part~2 names the transformer of
part~1.  $T_Q$ composes three explicit trees: on a proof $q$ of
$\ConR$, the cut transformer of Theorem~\ref{thm:g-cut} at $Q$ builds
a proof $a$ of $G$, D3G builds the proof $\sigma\cdot a$ of $\neg G$,
and EXPL builds the proof $e\cdot a\cdot(\sigma\cdot a)$ of $\z = \tofn
1$.  The matrix is decided on both sides, so the implication under
the quantifier is introduced by discharging a decided hypothesis.  The
witness $T_1$ of part~3 is the same composition without the cut, and
the closed $T$ of part~2 takes $Q$ from Theorem~\ref{thm:g-split}.
Part~5 is the disjunction property again: \rga would have to refute $\PrC(\bquo{\ConR})$, which
Theorem~\ref{thm:br-loeb}(1) excludes, or $\ConR$, which
Theorem~\ref{thm:br-con} excludes.  \qed

\medskip
\noindent
Part~2 is formalized G\"odel~II with the transformer that turns a proof of consistency into a proof of $\z =
\tofn 1$ written out as a term.  The trees of $T_Q$ and $T_1$ contain
no choice ([\citen{conrga-paper}, \S\,``Machines without opaque indices'']).  Like $c_Q$, these trees are
astronomically large.  The derivations of parts~1 to~4, and the
transformer laws behind them, use no soundness theorem; part~5 is a
non-derivability fact and rests on soundness.

The transformers of this section differ in how uniform they are.
Table~\ref{tab:g-uniform} sums up the transformers' uniformity.

\begin{table}[tbp]
\centering
\small
\begin{tabular}{@{}lll@{}}
\hline
transformer & uniform in & one internal sentence? \\
\hline
$e$ of EXPL (Thm.~\ref{thm:g-expl}) & the term $g$: one closed term & yes, for each $g$ \\
$\sigma_{\mathrm u}$ of D3 (Thm.~\ref{thm:g-d3-uniform}) & the code $c$ & yes, one sentence over all codes \\
$m_{\varphi\psi}$ of D2 (Thm.~\ref{thm:g-d2}) & nothing: one term per pair $\varphi$, $\psi$ & one sentence per pair \\
$c_Q$, $T_Q$ (Thms.~\ref{thm:g-cut}, \ref{thm:g-sentence}) & nothing: one term per HOL proof $Q$ & one sentence per $Q$ \\
\hline
\end{tabular}
\caption{The uniformity of the explicit transformers.  Every
  transformer is a closed term compiled from an explicit index tree;
  the trees of $c_Q$ and $T_Q$ embed the codes of $Q$'s formulas and
  are astronomically large.  Whether D2 and cut have transformers that
  take the codes as run-time arguments, as $\sigma_{\mathrm u}$ does,
  is open (Section~\ref{sec:concl}).}
\label{tab:g-uniform}
\end{table}

The G\"odel sentence and the consistency sentence are interderivable,
and so are the inconsistency sentence and the provability of
consistency.

\begin{thm}[\textsf{goedelG\_con\_interderivable}, \textsf{notcon\_prc\_con\_interderivable}]
\label{thm:g-interderivable}
$G \eqv \ConR$, and $\neg\ConR \eqv \PrC(\bquo{\ConR})$.
\end{thm}

\noindent
\emph{Proof idea.}  The directions $\ConR \vdash G$ and
$\PrC(\bquo{\ConR}) \vdash \neg\ConR$ are Theorems~\ref{thm:br-hyp}
and~\ref{thm:g-formal}(1).  For $G \vdash \ConR$: ex falso as an
entailment, $\PrC(\bquo{\z = \tofn 1}) \vdash \PrC(\bquo{G})$
(Theorem~\ref{thm:br-dc-ent}(4)), with $\PrC(\bquo{G}) \vdash \neg G$
(Theorem~\ref{thm:br-diag-ent}), gives $\neg\ConR \vdash \neg G$; the
hypothesis $\neg\ConR = \PrC(\bquo{\z = \tofn 1})$ has the $\Sigma$
shape, so contraposition and double negation give $G \vdash \ConR$.
$\neg\ConR \vdash \PrC(\bquo{\ConR})$ is ex falso as an entailment at
$\ConR$.  \qed

Formalized L\"ob's theorem as an entailment, Theorem~\ref{thm:br-loeb-ent}(1),
is L\"ob's proof read the same way.  Fix a closed well-formed $\psi$.
The diagonal construction of \S\ref{sec:goedel-sentences} gives a
closed term $S$ and the sentence $L = \PrC(S) \limp \psi$ with
$\emptyset \vdash S = \bquo{L}$ (\S\ref{app:goedel-ent}).
\begin{enumerate}
\item $L \vdash \PrC(\bquo{L}) \limp \psi$, by contraposition of
  $\PrC(\bquo{L}) \vdash \PrC(S)$ inside the disjunction that $L$ is.
\item Necessitation of step~1 gives $\PrC(\bquo{L}) \vdash
  \PrC(\bquo{\PrC(\bquo{L}) \limp \psi})$; D3 gives $\PrC(\bquo{L})
  \vdash \PrC(\bquo{\PrC(\bquo{L})})$; and D2 combines the two entailments into
  $\PrC(\bquo{L}) \vdash \PrC(\bquo{\psi})$, hence $\PrC(S) \vdash
  \PrC(\bquo{\psi})$.
\item $\PrC(\bquo{\psi}) \limp \psi \vdash L$, by cases on the
  disjunction: from $\psi$ by $\lor$-introduction, and from
  $\neg\PrC(\bquo{\psi})$ by contraposition of step~2, which gives
  $\neg\PrC(S)$.
\item Necessitation of step~3 gives $\PrC(\bquo{\PrC(\bquo{\psi})
  \limp \psi}) \vdash \PrC(\bquo{L})$, and step~2 completes the
  entailment.
\end{enumerate}
Every implication in the argument is a disjunction introduced by
$\lor$-introduction, and every case analysis is on a disjunction that
is a hypothesis; no step discharges an undecided hypothesis.  The
classical proof is the same argument with modus ponens in place of
each cut.  At $\psi = (\z = \tofn 1)$ the entailment is formalized
L\"ob at $\z = \tofn 1$, whose implication form Theorem~\ref{thm:br-loeb}(3)
refutes.  D2 in general, Theorem~\ref{thm:br-dc-ent}(2), comes from
the decided D2 of Theorem~\ref{thm:g-d2}; that entailment is also the
two-hypothesis necessitation of Theorem~\ref{thm:g-cut} at the
derivation of $\psi$ from $\varphi \limp \psi$ and $\varphi$.

So all three formalized theorems hold in \rga as entailments:
G\"odel~I ($\ConR \vdash G$ and $\ConR \vdash \neg\PrC(\bquo{G})$, and,
contraposed, $\neg G \vdash \neg\ConR$),
G\"odel~II ($\ConR \vdash \neg\PrC(\bquo{\ConR})$) and L\"ob's
theorem for every closed well-formed sentence.  As implications, all
three theorems fail; in decided-antecedent form over proof codes, formalized
G\"odel~II, the contraposed formalized G\"odel~I and formalized L\"ob
hold (Theorem~\ref{thm:g-sentence}).

\subsection{Consequences}
\label{sec:goedel-cons}

\paragraph{Soundness is not a single sentence.}
The decided-antecedent form of L\"ob's hypothesis at $\z = \tofn 1$
(Theorem~\ref{thm:br-loeb}(4)) and [\citen{conrga-paper}, theorem~\textsf{rga\_no\_uniform\_soundness}], its
form with the truth predicate in the consequent, are the second theorem
read as L\"ob's theorem at $\z = \tofn 1$.  The form $\PrC(\bquo{\z = \tofn 1})
\limp \z = \tofn 1$ is excluded by L\"ob's rule
(Theorem~\ref{thm:br-loeb}(2)) and $\emptyset \nvdash \z = \tofn 1$.
For [\citen{conrga-paper}, theorem~\textsf{rga\_no\_uniform\_soundness}], the same decided case split, with the internal
refutation of $\Tru(\bquo{\z = \tofn 1})$ ([\citen{conrga-paper}, theorem~\textsf{iPrDecTrueAt\_cEq01\_never}])
in the role that EXPL and D3G play above, would turn the single-sentence
soundness sentence into a proof of $\ConR$.  So \rga proves its
soundness proof by proof ([\citen{conrga-paper}, theorem~\textsf{rga\_sound\_sch}]), and not for all
codes at once.  The entailment forms fail as well:
$\PrC(\bquo{\z = \tofn 1}) \nvdash \Tru(\bquo{\z = \tofn 1})$
(\S\ref{app:goedel-ent}).

\paragraph{Entailments do not give self-consistency.}
Entailments rescue G\"odel's steps, but not \rga's consistency.  The
entailment form of ``\rga is consistent'' is local reflection at
$\z = \tofn 1$,
$\PrC(\bquo{\z = \tofn 1}) \vdash \z = \tofn 1$: from a proof of $\z =
\tofn 1$, derive $\z = \tofn 1$.  The entailment's hypothesis has the $\Sigma$
shape, so by reductio the entailment holds exactly when $\vdash
\ConR$ does (Theorem~\ref{thm:br-loeb-ent}(2)), and G\"odel's second
theorem excludes both the entailment and $\vdash \ConR$.  What holds is ex falso one level up:
$\PrC(\bquo{\z = \tofn 1}) \vdash \PrC(\bquo{\varphi})$ for every
well-formed $\varphi$ (Theorem~\ref{thm:br-dc-ent}(4)).  The consistency that \rga
does derive in these forms is that of [\citen{conrga-paper}, Section ``Consistency, from unconditional to conditional'']:
sentences over decided classes of proof codes, and entailments from
hypotheses that \rga neither proves nor refutes.  That section also
shows that no Gentzen-style consistency proof of a natural shape exists
in \rga, and that the hypothesis that the proof search for $\z = \tofn
1$ never stops entails $\ConR$, $G$ and the unprovability of $\ConR$ and of $G$
([\citen{conrga-paper}, \S\,``What the theorems show'']).

\subsection{The two proofs side by side}
\label{sec:goedel-compare}

Table~\ref{tab:goedel-compare} in the introduction sets G\"odel's
proof and \rga's proof side by side, step by step.  In summary, \rga
keeps the diagonal construction, as an equation of codes and as a
pair of entailments; has D1--D3 as rules (Theorem~\ref{thm:br-rules})
and D2 and D3 as entailments (Theorem~\ref{thm:br-dc-ent}), but not
as implications over provability, which fail at $G$
(Theorem~\ref{thm:br-d2d3}); uses as sentences the decided forms with
explicit witnesses (EXPL, D3G); and, as free-variable arithmetic does,
replaces formalized G\"odel~I as an implication by the decided case
split of Theorem~\ref{thm:g-split}, which is the entailment $\ConR
\vdash G$.  The conclusion is the classical one, and the formalized
theorems hold as entailments and, where the hypothesis is a
provability statement, in decided-antecedent form
(\S\ref{sec:goedel-formal}).

\section{The second theorem once: the G\"odel conditions}
\label{sec:locale}

Section~\ref{sec:goedel} proved G\"odel's theorems for \rga's tree
checker.  The textbook route of that section, read with entailments,
uses few facts about the checker: the checker commutes with
substitution and its verdict is decided, three maps, necessitation
(NEC), D3 and EXPL, carry accepted proofs to accepted proofs, every
theorem has an accepted proof (D1), and \rga is consistent.  This
section states those facts as hypotheses about an abstract entailment
relation and an abstract checker, the \emph{G\"odel conditions}, and
proves once, from the conditions alone, the first half of G\"odel's
first theorem, $\nvdash G$, and the second theorem.  A system that meets
the conditions is a \emph{G\"odel extension of \rga}.  \rga with its
tree checker is one (\S\ref{sec:locale-rga}).  \iga, the extension of
\rga by ground induction~\cite{iga-paper}, with a checker of its own,
is another (\cref{sec:iga-goedel}).

The development states the conditions as Isabelle \emph{locales}.  A
locale is a named list of parameters and assumptions.  A theorem
proved inside a locale holds for every choice of the parameters that
meets the assumptions, and an \emph{interpretation} of the locale, at
particular parameters and with proofs of the assumptions, yields
every theorem of the locale for those parameters.  The G\"odel
conditions form nine locales, each built on the ones below it, and
the development places each theorem in a layer whose assumptions
include those that the theorem's proof uses.

The theorems of this section come in two kinds, which the text labels
by viewpoint.  A rule, ``if $\vdash A$ then $\vdash B$'', and a
non-derivability $\nvdash A$ are \emph{classical}: HOL statements
about derivability, for every G\"odel extension at once, in HOL's
reading.  An entailment or a sentence that the system derives is \emph{a
derivation in the extended system}: for every G\"odel extension, the
extended system itself derives that judgment, by reasoning inside the system.

\subsection{The calculus}
\label{sec:locale-calc}

The calculus fixes \rga's language and rules and leaves the
entailment relation open.

\begin{defi}[\textsf{gcalc}, \textsf{gcalc\_strip}]
\label{def:lc-calc}
An \emph{entailment relation over \rga} is a relation $\Gamma \vdash
c$ between finite sets $\Gamma$ of \rga terms and \rga terms $c$ with
three properties.
\begin{enumerate}
\item Every conclusion is well-formed: $\Gamma \vdash c$ gives
  $\wf(c)$.
\item The relation is closed under each of \rga's 42 rules, read as
  closure conditions, with \rga's own head step in the conversion
  rules.
\item Strip: $\Gamma \vdash c$ gives $\Gamma' \vdash c$, where
  $\Gamma'$ keeps the well-formed members of $\Gamma$.
\end{enumerate}
\end{defi}

\noindent
The term language, substitution, well-formedness and the head step
are \rga's own; only the relation $\vdash$ varies.  The calculus
has 44 named assumptions: well-formed conclusions, the 42 rules and
strip (\S\ref{app:locale-calc}).  The definition does not ask $\vdash$
to be the derivability relation of a proof system, nor to be
recursively enumerable.  Every relation with the three properties is
an instance.

\begin{thm}[\textsf{rga\_embed}, \textsf{sigma\_\allowbreak deduction\_\allowbreak theorem}]
\label{thm:lc-embed}
Let $\vdash$ be an entailment relation over \rga.
\begin{enumerate}
\item If \rga derives $H \vdash c$, then $H \vdash c$.
\item If $\vz\jnat, \lft\Gamma \vdash D(\vz)\jbool$, then for every
  $C$: $\exists s.\,D(s), \Gamma \vdash C$ if and only if $\Gamma
  \vdash \forall s.\,(D(s) \limp C)$.
\end{enumerate}
\end{thm}

\noindent
Part~1 makes every \rga derivation available in every instance; the
proofs below use part~1 for derivations about closed codes, such as
the values of the diagonal function.  Part~1 holds without strip.
Part~2 is the $\Sigma$ deduction theorem of
Theorem~\ref{thm:rc-dt}(2), proved once from the rules and strip, and
$\Sigma$ contraposition (Corollary~\ref{cor:rc-contrapos}) carries
over in the same way.  Both parts are classical.

\begin{thm}[\textsf{gcalc\_strip\_rga}, \textsf{is\_gcalc}]
\label{thm:lc-inst}
\leavevmode
\begin{enumerate}
\item \rga's derivability is an entailment relation over \rga.
\item Let $S$ be a one-step relation on proof lists that contains
  \rga's step relation, is monotone in the list, and gives only
  well-formed conclusions from lists whose conclusions are
  well-formed.  Then the entailment that $S$ defines has properties~1
  and~2 of Definition~\ref{def:lc-calc}.
\end{enumerate}
\end{thm}

\noindent
Part~2 is the route by which a rule extension of \rga becomes an
instance: the extension's step relation adds rules to \rga's, and a
strip theorem for the extension supplies property~3.  \iga takes this
route (\cref{sec:iga-goedel}).

\subsection{The conditions}
\label{sec:locale-cond}

A \emph{checker} is a HOL function $\gchk$ that maps two terms, a
proof slot and a conclusion slot, to a term $\gchk(p, c)$.  The
checker is a term former of HOL, not an object-level term.  Every
statement below is relative to the checker.

\begin{defi}[\textsf{ConT}, \textsf{goedelGT}]
\label{def:lc-sentences}
For a checker $\gchk$,
\[
  \gPr(c) \;\ldef\; \exists q.\ \gchk(q, c) , \qquad
  \gCon \;\ldef\; \forall q.\ \neg\gchk(q, \bquo{\z = \tofn 1}) ,
\]
with $\exists$ written as $\neg\forall\neg$, and, with $x$ free,
\[
  \gD \;\ldef\; \forall q.\ \neg\gchk(q, \Diag(x, x)) , \qquad
  \gS \;\ldef\; \Diag(\bquo{\gD}, \bquo{\gD}) , \qquad
  \gG \;\ldef\; \forall q.\ \neg\gchk(q, \gS) .
\]
\end{defi}

\noindent
$\Diag$ is \rga's compiled diagonal function of
\S\ref{sec:goedel-sentences}.  The locale uses a copy of $\Diag$, the
same explicit index tree under another name, and \rga proves the
copy's values at closed codes as it proves $\Diag$'s.  At \rga's tree
checker, $\gPr$ is literally $\PrC$, and the interpretation rewrites
$\gCon$ and $\gG$ to $\ConR$ and $G$ (\S\ref{sec:locale-rga}).

\begin{defi}[\textsf{g2\_syn}, \textsf{g2\_con}]
\label{def:lc-cond}
An entailment relation $\vdash$ over \rga and a checker $\gchk$ meet
the \emph{G\"odel conditions} when the following hold for all terms
and contexts.
\begin{description}
\item[Structural conditions.]  (S1) The checker commutes with every
  substitution $f$ of terms for variables: $\sof{f}{\gchk(p, c)} =
  \gchk(\sof{f}{p}, \sof{f}{c})$.  (S2) $\gchk(p, c)$ is well-formed
  exactly when $p$ and $c$ are.  (S3) The verdict is decided at
  naturals: $\Gamma \vdash p\jnat$ and $\Gamma \vdash c\jnat$ give
  $\Gamma \vdash \gchk(p, c)\jbool$.
\item[Transformer conditions.]  Each of NEC, D3 and EXPL asserts, for
  each instance, a HOL function $T$ on terms such that in every
  context $\Delta$, at all terms $u, v$ with $\Delta \vdash u\jnat$ and
  $\Delta \vdash v\jnat$, the stated acceptances give $\Delta \vdash
  T(v)\jnat$ (or $T(u, v)\jnat$) and the stated conclusion.
  \begin{itemize}
  \item (NEC) Whenever $\{A\} \vdash B$, $\Delta \vdash
    \gchk(v, \bquo{A})$ gives $\Delta \vdash \gchk(T(v), \bquo{B})$.
  \item (D3) For each term $\varphi$: $\Delta \vdash \gchk(v,
    \bquo{\varphi})$ gives $\Delta \vdash \gchk(T(v),
    \bquo{\gPr(\bquo{\varphi})})$.
  \item (EXPL) For each term $g$: $\Delta \vdash \gchk(u, \bquo{g})$
    and $\Delta \vdash \gchk(v, \bquo{\neg g})$ give $\Delta \vdash
    \gchk(T(u, v), \bquo{\z = \tofn 1})$.
  \end{itemize}
\item[Closed acceptance.]  (D1) If $\emptyset \vdash \varphi$, then
  $\emptyset \vdash \gchk(\tofn k, \bquo{\varphi})$ for some number
  $k$.
\item[Consistency.]  (CON) $\emptyset \nvdash \z = \tofn 1$.
\end{description}
A \emph{G\"odel extension of \rga} is an entailment relation over
\rga together with a checker such that the two meet the G\"odel
conditions.
\end{defi}

\noindent
D1 is not a transformer; it asserts one accepted numeral for each
theorem of the empty context and maps no proof codes.  The
transformers are uniform in the context and in the proof code, which
may be an open term, but not in the sentence, and the HOL functions $T$
for different sentences need not be one internal term.  Only D1 gives
an accepted numeral from a theorem, and no condition goes the other
way, so the checker need not be faithful.  Condition S3 is a
derivation in the system, not a claim of HOL decidability.

Three further conditions are optional.  They are not hypotheses of the
second theorem, and only the theorems of \S\ref{sec:locale-thm} that
name them use them.  Each optional condition is a transformer law of
the same shape.
(EXF), ex falso: for well-formed $\varphi$, $\gchk(v, \bquo{\z =
\tofn 1})$ gives $\gchk(T(v), \bquo{\varphi})$.  (D2): for
well-formed $\varphi$ and $\psi$, $\gchk(u, \bquo{\varphi \limp
\psi})$ and $\gchk(v, \bquo{\varphi})$ give $\gchk(T(u, v),
\bquo{\psi})$.  (NEC2): whenever $\{A_1, A_2\} \vdash
B$, $\gchk(u, \bquo{A_1})$ and $\gchk(v, \bquo{A_2})$ give $\gchk(T(u,
v), \bquo{B})$.

\begin{table}[tbp]
\centering
\small
\begin{tabular}{@{}lll>{\raggedright\arraybackslash}p{0.42\linewidth}@{}}
\hline
locale & built on & adds & role \\
\hline
\textsf{g2\_syn} & the calculus & S1, S2, S3 & the checker's syntax and decidedness; the sentences of Definition~\ref{def:lc-sentences} \\
\textsf{g2\_ent} & \textsf{g2\_syn} & NEC, D3, EXPL & the textbook route as entailments (\S\ref{sec:locale-route}) \\
\textsf{g2\_acc} & \textsf{g2\_syn} & D1 & closed acceptance \\
\textsf{g2\_con} & \textsf{g2\_ent}, \textsf{g2\_acc} & CON & $\nvdash \gG$ and G\"odel's second theorem \\
\textsf{g2\_exf} & \textsf{g2\_ent} & EXF & ex falso; the interderivabilities \\
\textsf{g2\_d2} & \textsf{g2\_ent} & D2 & D2 as an entailment; formalized L\"ob \\
\textsf{g2\_nec2} & \textsf{g2\_ent} & NEC2 & necessitation for two hypotheses \\
\textsf{g2\_exf\_acc} & \textsf{g2\_exf}, \textsf{g2\_acc} & --- & Beklemishev and Shamkanov's structure conditions \\
\textsf{g2\_exf\_con} & \textsf{g2\_exf\_acc}, \textsf{g2\_con} & --- & Beklemishev and Shamkanov's second theorem; no derivable unprovability \\
\hline
\end{tabular}
\caption{The nine locales of the G\"odel conditions.  The second
  theorem is a theorem of \textsf{g2\_con}, whose assumptions are the
  calculus and the eight conditions of Definition~\ref{def:lc-cond}.
  Appendix~\ref{app:locale-cond} names the eleven assumptions of the
  layers and gives the form in which Isabelle states them, and
  \S\ref{app:locale-calc} names the calculus's 44.}
\label{tab:lc-layers}
\end{table}

Table~\ref{tab:lc-layers} lists the layers.  Isabelle gives a theorem
proved in a locale all of that locale's assumptions as premises, so
the layering bounds the conditions that each theorem may use.  A
theorem's proof uses at most the assumptions of its locale.  The
second theorem lives in \textsf{g2\_con}, which is built from
\textsf{g2\_ent}, \textsf{g2\_acc} and consistency, and none of these
three includes \textsf{g2\_d2}.  D2 is not among the conditions of the
second theorem.  The layering does not say that D2 fails, nor that D2 does
not follow from the other conditions; at \rga, D2 holds
(Theorem~\ref{thm:g-d2}).  Under \textsf{g2\_con} the named
assumptions number 52, the calculus's 44 and the eight conditions;
across all nine locales they number 55.

\subsection{The textbook route, as entailments}
\label{sec:locale-route}

Fix an entailment relation over \rga and a checker.  The textbook
proof of \S\ref{sec:bg-goedel}, read with entailments as in
Section~\ref{sec:goedel}, goes through in the locale.  Write $\gPr$,
$\gCon$, $\gG$ and $\gS$ as in Definition~\ref{def:lc-sentences}.

\begin{thm}[\textsf{goedelST\_eq}, \textsf{goedelT\_\allowbreak diag\_\allowbreak ent}]
\label{thm:lc-diag}
If the structural conditions hold, then, for every context $\Gamma$:
\begin{enumerate}
\item $\emptyset \vdash \gS = \bquo{\gG}$;
\item $\gG \vdash \neg\gPr(\bquo{\gG})$ and $\neg\gPr(\bquo{\gG})
  \vdash \gG$;
\item $\neg\gG, \Gamma \vdash \gPr(\bquo{\gG})$ and
  $\gPr(\bquo{\gG}), \Gamma \vdash \neg\gG$;
\item $\gPr(\bquo{\z = \tofn 1}), \Gamma \vdash \z = \tofn 1$ if and
  only if $\Gamma \vdash \gCon$.
\end{enumerate}
\end{thm}

\noindent
As in Theorem~\ref{thm:g-diag}, $\neg\gG$ is literally $\gPr(\gS)$.
The equation of part~1 is \rga's derivation of the diagonal's value,
carried into the instance by Theorem~\ref{thm:lc-embed}(1), and
condition~S1 transports the checker's conclusion along the equation.
Parts~1 to~3 are derivations in the extended system.  Part~2 is two entailments; the theorem does
not give one biconditional sentence, which \rga does not derive
(Theorem~\ref{thm:br-bicond}).  Part~4, local reflection at $\z =
\tofn 1$, is classical, an equivalence of two derivability facts.  Its
proof uses the $\Sigma$ rules of Theorem~\ref{thm:lc-embed}(2) and
Corollary~\ref{cor:rc-contrapos} at the decided verdict, as at \rga
(Theorem~\ref{thm:br-loeb-ent}(2)).

\begin{thm}[\textsf{ent\_necT}, \textsf{D3\_entT}]
\label{thm:lc-ent}
If the structural conditions and NEC, D3 and EXPL hold, then, for
every context $\Gamma$:
\begin{enumerate}
\item if $\{A\} \vdash B$, then $\gPr(\bquo{A}), \Gamma \vdash
  \gPr(\bquo{B})$;
\item $\gPr(\bquo{\varphi}), \Gamma \vdash
  \gPr(\bquo{\gPr(\bquo{\varphi})})$ for every term $\varphi$;
\item $\gPr(\bquo{\gG}), \Gamma \vdash \gPr(\bquo{\neg\gG})$;
\item $\gPr(\bquo{\neg g}), \gPr(\bquo{g}), \Gamma \vdash
  \gPr(\bquo{\z = \tofn 1})$ for every term $g$.
\end{enumerate}
\end{thm}

\noindent
Each part eliminates the existential of a hypothesis $\gPr(c)$ at a
fresh proof variable $v$, applies the transformer, and introduces the
existential again at the transformer's output, as at \rga
(Theorem~\ref{thm:br-dc-ent}).  Part~3, D3G, is D3 at $\gG$ followed
by NEC applied to the entailment $\gPr(\bquo{\gG}) \vdash \neg\gG$ of
Theorem~\ref{thm:lc-diag}(3).  Part~1 is a rule, and so classical;
parts~2 to~4 are derivations in the extended system.

\begin{thm}[\textsf{goedelGT\_\allowbreak under\_\allowbreak ConT\_\allowbreak textbook}, \textsf{goedel2\_\allowbreak internal}]
\label{thm:lc-textbook}
If the structural conditions and NEC, D3 and EXPL hold, then:
\begin{enumerate}
\item $\neg\gG, \Gamma \vdash \neg\gCon$ for every context $\Gamma$;
\item $\gCon \vdash \gG$;
\item $\gCon \vdash \neg\gPr(\bquo{\gG})$;
\item $\gPr(\bquo{\gCon}) \vdash \neg\gCon$;
\item $\gCon \vdash \neg\gPr(\bquo{\gCon})$;
\item $\emptyset \vdash \forall q.\,(\gchk(q, \bquo{\gCon}) \limp
  \neg\gCon)$.
\end{enumerate}
\end{thm}

\noindent
\emph{Proof.}  (1) From $\neg\gG$, Theorem~\ref{thm:lc-diag}(3) gives
$\gPr(\bquo{\gG})$, Theorem~\ref{thm:lc-ent}(3) gives
$\gPr(\bquo{\neg\gG})$, and Theorem~\ref{thm:lc-ent}(4) gives
$\gPr(\bquo{\z = \tofn 1})$, which is literally $\neg\gCon$.  (2)
Since $\neg\gG$ is literally $\gPr(\gS)$, part~1 reads $\gPr(\gS),
\Gamma \vdash \neg\gCon$.  The hypothesis $\gPr(\gS)$ has the
decided matrix $\gchk(v, \gS)$, by condition~S3, so $\Sigma$
contraposition (Corollary~\ref{cor:rc-contrapos}, in the instance)
gives $\neg\neg\gCon \vdash \neg\gPr(\gS)$, which is
$\neg\neg\gG$; double negation, introduced and eliminated by \rga's
rules, gives part~2.  (3) Part~2 and Theorem~\ref{thm:lc-diag}(2), by
cut.  (4) NEC at part~2 (Theorem~\ref{thm:lc-ent}(1)) gives
$\gPr(\bquo{\gCon}) \vdash \gPr(\bquo{\gG})$, and parts~3 and~4 of
Theorem~\ref{thm:lc-ent} continue as in part~1.  (5) $\Sigma$
contraposition of part~4, and double negation.  (6) The $\Sigma$
deduction theorem (Theorem~\ref{thm:lc-embed}(2)) at part~4, with the
matrix $\gchk(q, \bquo{\gCon})$ decided by condition~S3.  \qed

\medskip
\noindent
Every part is a derivation in the extended system, and no part uses D1
or consistency.  The layer
\textsf{g2\_ent}, in which the theorem lives, has neither D1 nor
consistency.  The step at which
the restriction of discharge to decided hypotheses shows is the
contraposition in parts~2 and~5.  The textbook proof contraposes an implication; here the
hypothesis $\gPr(\gS)$ has a decided matrix, and the decided matrix
licenses the contraposition.  Part~6 puts the decided $\gchk(q,
\bquo{\gCon})$ in antecedent position.  The theorem does not give the
implication $\gCon \limp \neg\gPr(\bquo{\gCon})$, whose antecedent
$\gCon$ is undecided; at \rga that implication is not a theorem
(Theorem~\ref{thm:g-formal}(3)).

\subsection{The theorems}
\label{sec:locale-thm}

\begin{thm}[\textsf{goedel1}, \textsf{goedel2\_textbook}]
\label{thm:lc-g2}
In every G\"odel extension of \rga:
\begin{enumerate}
\item if $\emptyset \vdash \gG$, then $\emptyset \vdash \z = \tofn 1$;
\item $\emptyset \nvdash \gG$;
\item $\emptyset \nvdash \gCon$.
\end{enumerate}
Part~1 holds already under the structural conditions and D1.
\end{thm}

\noindent
\emph{Proof.}  (1) D1 gives a number $k$ with $\emptyset \vdash
\gchk(\tofn k, \bquo{\gG})$, and condition~S1 with
Theorem~\ref{thm:lc-diag}(1) transports the verdict to $\gchk(\tofn k,
\gS)$.  The sentence $\gG$ instantiated at $\tofn k$ gives
$\neg\gchk(\tofn k, \gS)$, and \rga's rule from $\neg a$ and $a$
derives $\z = \tofn 1$.  (2) Part~1 and consistency.  (3) A
derivation of $\gCon$, cut against $\gCon \vdash \gG$
(Theorem~\ref{thm:lc-textbook}(2)), gives a derivation of $\gG$,
against part~2.  \qed

\medskip
\noindent
The theorem is classical.  The proof of part~3 uses every one of the
eight conditions.  S3 licenses the contraposition; S1 transports
verdicts along equations; S2 discharges the well-formedness side
conditions; NEC, D3 and EXPL enter through
Theorem~\ref{thm:lc-textbook}; D1 and consistency enter only in the
external step, parts~1 and~2.  The proof also uses strip, through the
$\Sigma$ deduction theorem and $\Sigma$ contraposition.  The proof
uses neither D2, nor the ex falso condition EXF, nor NEC2, nor
the checker's faithfulness, \rga's soundness, a truth predicate, cut
elimination or the disjunction property.  The theorem does not say
that the system refutes $\gCon$, nor anything about a checker other
than $\gchk$, nor that the conditions are minimal.

Section~\ref{sec:goedel}'s proof for \rga took a different route, with
explicit witnesses for EXPL and D3G in place of the transformer
conditions.  That route holds in the locale too.

\begin{thm}[\textsf{goedel2T\_\allowbreak of\_\allowbreak expl\_\allowbreak d3g}, \textsf{no\_\allowbreak loeb\_\allowbreak hypothesis\_\allowbreak 0eq1T}]
\label{thm:lc-witness}
In every G\"odel extension of \rga:
\begin{enumerate}
\item if $e$ and $\sigma$ are closed terms, the system derives
  $\forall u.\,\forall v.\,(e \cdot u \cdot v)\jnat$ and $\forall
  v.\,(\sigma \cdot v)\jnat$, and the system derives the sentences
  \begin{gather*}
    \forall u.\,\forall v.\,\bigl((\gchk(u, \bquo{\gG}) \land \gchk(v,
    \bquo{\neg\gG})) \limp \gchk(e\cdot u\cdot v, \bquo{\z = \tofn
    1})\bigr) , \\
    \forall v.\,\bigl(\gchk(v, \bquo{\gG}) \limp \gchk(\sigma\cdot v,
    \bquo{\neg\gG})\bigr) ,
  \end{gather*}
  then $\emptyset \nvdash \gCon$;
\item $\emptyset \nvdash \forall q.\,(\gchk(q, \bquo{\z = \tofn 1})
  \limp \z = \tofn 1)$.
\end{enumerate}
\end{thm}

\noindent
Both parts are classical.  Part~1 is the witness form of
Theorems~\ref{thm:g-split} and~\ref{thm:g-g2}, for every G\"odel
extension; its proof is the decided case split of
Theorem~\ref{thm:g-split} and uses the structural conditions, D1 and
consistency, with $e$ and $\sigma$ as HOL premises, and not the
transformer conditions.  At \rga, Theorems~\ref{thm:g-expl}
and~\ref{thm:g-d3g} supply $e$ and $\sigma$.  Part~2 says that the
system does not derive the decided-antecedent form of local reflection
at $\z = \tofn 1$.  From that sentence a case split on the decided
verdict would derive $\gCon$, against Theorem~\ref{thm:lc-g2}(3).

The optional layers carry the rest of Section~\ref{sec:goedel}'s
internal family.

\begin{thm}[\textsf{D2\_entT}, \textsf{formalized\_\allowbreak loeb\_\allowbreak entT}]
\label{thm:lc-loeb}
If the structural conditions and NEC, D3, EXPL and D2 hold, then:
\begin{enumerate}
\item $\gPr(\bquo{\varphi}), \gPr(\bquo{\varphi \limp \psi}), \Gamma
  \vdash \gPr(\bquo{\psi})$ for all well-formed $\varphi$, $\psi$ and
  every context $\Gamma$;
\item $\gPr(\bquo{\gPr(\bquo{\psi}) \limp \psi}) \vdash
  \gPr(\bquo{\psi})$ for every closed well-formed $\psi$;
\item $\emptyset \vdash \forall q.\,(\gchk(q, \bquo{\gPr(\bquo{\psi})
  \limp \psi}) \limp \gPr(\bquo{\psi}))$ for every closed well-formed
  $\psi$.
\end{enumerate}
\end{thm}

\noindent
All three parts are derivations in the extended system.  Part~2, formalized L\"ob's theorem as an
entailment, runs through the L\"ob sentence $L = \gPr(S_L) \limp
\psi$, with $\vdash S_L = \bquo{L}$ from the diagonal, and through
NEC, D3, D2 and $\Sigma$ contraposition; part~3 is part~2 in
decided-antecedent form.  Part~2 is the locale version of
Theorem~\ref{thm:br-loeb-ent}(1).
L\"ob's rule, ``if $\vdash \PrC(\bquo{\psi}) \limp \psi$ then $\vdash
\psi$'' (Theorem~\ref{thm:br-loeb}), is not proved in the locale.  \rga's
proof of that rule uses the disjunction property, which is not a
G\"odel condition.

\begin{thm}[\textsf{goedelGT\_\allowbreak con\_\allowbreak interderivable}, \textsf{no\_\allowbreak code\_\allowbreak unprovableT}]
\label{thm:lc-exf}
If the structural conditions and NEC, D3, EXPL and EXF hold, then:
\begin{enumerate}
\item $\gPr(\bquo{\z = \tofn 1}), \Gamma \vdash \gPr(\bquo{\varphi})$
  for every well-formed $\varphi$ and every context $\Gamma$;
\item $\gG \eqv \gCon$ and $\neg\gCon \eqv \gPr(\bquo{\gCon})$;
\item if moreover D1 and consistency hold, then for no well-formed
  $\varphi$ does $\emptyset \vdash \forall q.\,\neg\gchk(q,
  \bquo{\varphi})$.
\end{enumerate}
\end{thm}

\noindent
Parts~1 and~2 are derivations in the extended system, and part~2 is the locale version of
Theorem~\ref{thm:g-interderivable}.  Part~3 is classical.  A
derivation of $\forall q.\,\neg\gchk(q, \bquo{\varphi})$ gives one of
$\gCon$, by EXF and a case split on the decided verdict, against
Theorem~\ref{thm:lc-g2}(3).  So no well-formed sentence is derivably
unprovable; part~3 does not say that any unprovable sentence is
provable.

\begin{thm}[\textsf{bsT\_\allowbreak structure\_\allowbreak conditions}, \textsf{bsT\_\allowbreak goedel2\_\allowbreak irrefutable}]
\label{thm:lc-bs}
Read Beklemishev and Shamkanov's abstract provability
structures~\cite{beklemishev16abstract} as in \S\ref{app:goedel-ent},
with $\gPr(\bquo{x})$ for $\Box x$.
\begin{enumerate}
\item If the structural conditions, NEC, D3, EXPL, EXF and D1 hold,
  then the structure conditions C1 (its monotone half) to C5 hold
  (\S\ref{app:locale-map}).
\item If moreover consistency holds, then $\gPr(\bquo{\neg(\z =
  \z)}) \nvdash \neg(\z = \z)$.
\end{enumerate}
\end{thm}

\noindent
Here $\neg(\z = \z)$ is the structures' falsum, a classical falsum,
not \rga's divergent constant~$\bot$.  Of the seven structure
conditions of part~1 (\S\ref{app:locale-map}, in that order), the
first, fourth, sixth and seventh are entailments, and derivations in the extended system; the
second, third and fifth are rules between entailments, and classical.
Part~2 is classical.

\subsection{\rga as an instance}
\label{sec:locale-rga}

\begin{thm}[\textsf{g2\_con\_rga}, \textsf{g2\_exf\_con\_rga}]
\label{thm:lc-rga}
\rga's derivability with the tree checker $\Chk$ is a G\"odel
extension of \rga, and meets the conditions EXF, D2 and NEC2 too.
\end{thm}

\noindent
The witnesses are the machinery of Section~\ref{sec:goedel}
(\S\ref{app:locale-witness} lists them): the checker's laws at
abstract arguments for S1 to S3 (\S\ref{sec:recall-chk}
and~\S\ref{app:goedel-chk}); $T(v) = c_Q \cdot v \cdot v$ for NEC, with
$Q$ a proof of $\{A\} \vdash B$ (Theorem~\ref{thm:g-cut}); the decided
D3 at $\varphi$ (Theorem~\ref{thm:g-d3-any}); $T(u, v) = e \cdot u
\cdot v$ for EXPL (Theorem~\ref{thm:g-expl}); for D1, closed evaluation
of the checker at a proof's code, which uses neither \rga's soundness
nor the checker's faithfulness; and \rga's HOL consistency theorem.
The optional conditions take the transformers of decided D2
(Theorem~\ref{thm:g-d2}), of ex falso (\S\ref{app:goedel-cons}) and of
the two-hypothesis cut (Theorem~\ref{thm:g-cut}).

The development interprets the nine locales at \rga one at a time,
each after the layers it is built on, with rewrites that turn the
locale's $\gCon$, $\gG$ and copy of the diagonal into $\ConR$, $G$ and
$\Diag$.  So each interpreted theorem has the proposition of the
corresponding theorem of Section~\ref{sec:goedel}, and, since an
interpreted fact depends on the witnesses of its layer only, the facts
of \textsf{g2\_syn} and \textsf{g2\_ent} do not depend on the
consistency witness or, through it, on \rga's soundness.

\begin{cor}[\textsf{rga\_\allowbreak goedel2\_\allowbreak textbook}]
\label{cor:lc-rga-g2}
$\emptyset \nvdash \ConR$, as the instance at \rga of
Theorem~\ref{thm:lc-g2}(3).
\end{cor}

\noindent
The corollary is classical.  It has literally the proposition of
Theorem~\ref{thm:g-g2}, which keeps its own proof by the explicit
witness of Theorem~\ref{thm:g-d3g}, and a build-time check compares
the two propositions.  In the same way the development defines 62
facts of its G\"odel theories for \rga, many of them structural
lemmas and not G\"odel theorems, as instances of locale facts, among
them \textsf{goedel\_diag\_ent} of Theorem~\ref{thm:br-diag-ent},
\textsf{ent\_nec} of Theorem~\ref{thm:g-cut}, and the facts named in
the headers of Theorems~\ref{thm:br-dc-ent}, \ref{thm:br-loeb-ent},
\ref{thm:g-formal} and~\ref{thm:g-interderivable}.  Build-time checks
fix their names and propositions by fingerprints and keep the
instances outside the two consistency layers clear of nine named
soundness, consistency, faithfulness and truth facts of \rga
(\S\ref{app:locale-inst}, \S\ref{app:locale-checks}).  The
decider-free check of [\citen{conrga-paper}, Appendix ``The build-time check''] runs on 18
of the instances, of 17 locale facts; since an instance's dependencies
contain those of its locale theorem, no statement among the
dependencies of those 17 locale theorems mentions a forbidden
constant of the decider, a consequence that no check prints.

\subsection{Scope}
\label{sec:locale-scope}

The second theorem covers any entailment relation and checker on
\rga's language that satisfy the G\"odel conditions, whether or not
the relation is recursively enumerable or the checker faithful.  A rule
extension of \rga on \rga's terms qualifies through
Theorem~\ref{thm:lc-inst}(2), once the extension has strip and a
checker that meets the conditions;
\iga is one (\cref{sec:iga-goedel}).  Outside the scope lie systems in
another language, with new term formers or another head step, systems
in which some rule of \rga fails, and systems without a checker that
meets the conditions.  The definition does not exclude extensions of
\rga by axioms.  One family of axiom extensions is covered: \iga with
a notation claim of Kleene's $\KO$ as an axiom has a checker, \iga's
table with one axiom cell, that meets the conditions, consistency
included when the claimed code lies in $\KO$ (\S\ref{sec:adopt-claims-g2}).  For
other axiom extensions, such as \rga${}+\ConR$, the development has no
checker and no instance, and whether one meets the conditions, D1 and
consistency in particular, is not established, so the second theorem
of this section does not cover such extensions.  The transformer conditions
are per-sentence HOL existentials.  At \rga, D3 also has one internal
transformer for every code (Theorem~\ref{thm:g-d3-uniform}), and
EXPL's transformer does not depend on the sentence, but the conditions
do not ask for either.  Since they do not ask for faithfulness either,
$\gCon$ is the consistency sentence of the checker $\gchk$, not
necessarily of the system.

\section{A generic checker and \iga's checker}
\label{sec:checker}

Section~\ref{sec:locale} proves the second theorem for every extension
of \rga with a checker that meets the G\"odel conditions, so each
extension needs a checker of its own.  \rga's checker is large: the session
\texttt{RGA\_Checker}, which holds the checker's definition and its
laws, has 87,242 lines (\texttt{wc -l} over its 24 theory files); the
detectors themselves sit in \texttt{RGA\_Kit}
(Appendix~\ref{app:ck-files}).  This section writes
that checker once over a parameter, its \emph{table} of rule detectors.
The resulting \emph{generic checker} is \rga's checker when the table
is \rga's, and an extension of \rga by new rules gets its checker by
adding a tag per rule, one new detector for each new rule, once each
new detector's shift law and faithfulness are proved.  The
derivation laws of the generic checker hold at every table that
satisfies one law, the shift law of \S\ref{sec:ck-shift}.  \iga's checker is \rga's
table plus one tag, for ground induction (\S\ref{sec:ck-iga}).  The
generic layer, the session \texttt{RGA\_Checker\_Gen}, has 1,893 lines
in two theory files, and \iga's instance, the session
\texttt{IGA\_Checker}, has 1,211 lines in three.  Both sessions only
add theories on top of \rga's development, so the generic checker
changes none of the statements that the internal-reasoning paper
prints.
Section~\ref{sec:iga-goedel} uses \iga's checker to prove G\"odel's
second theorem at \iga.

Every derivation in this section is a derivation of \rga, written
$\vdash$, including the derivations about \iga's checker, which is an
\rga term.  \iga derives the same judgments, since every \rga
derivation is an \iga derivation.  Two of the kinds of statement of
\S\ref{sec:recall-forms} recur, schematic and meta.  A law at open
codes is schematic, a derived rule of \rga.  Read classically, the
law is a closure property of \rga's derivability; read as derivations
in \rga, it is a schema, each instance one \rga derivation at an open
code, in any context.  A
statement at closed codes, such as faithfulness, is a HOL fact about
which judgments \rga derives at numerals; its parts are schematic or
meta (Theorem~\ref{thm:ck-faith}).

\paragraph{\iga.}
Inductive Grounded Arithmetic~\cite{iga-paper}, \iga, adds to \rga's 42
rules one rule, which the \iga paper calls ground induction: from the
fixed point $\mu\Phi$ at a point $x$, and a derivation that a motive $c$
is closed under one unfolding of $\Phi$, conclude the motive at $x$.
Written in this paper's notation,
\[
  \frac{\Gamma \vdashI \mu\Phi\cdot x \qquad
        i\jnat,\ \Phi[\lambda c/a][i/x],\ \Gamma \vdashI c[i] \qquad
        \Gamma \vdashI x\jnat}
       {\Gamma \vdashI c[x]}
\]
where $\Phi$ is positive in its recursion variable $a$, $\Phi$ and $c$
are well-formed, and $i$ is not free in $\Gamma$, in $\lambda c$ or in
$\lambda a.\lambda x.\Phi$; $\vdashI$ is \iga's derivability.  \iga
therefore has 43 rules.  The \iga paper proves \iga
consistent~\cite{iga-paper}.

\subsection{The checker at a table}
\label{sec:ck-table}

\rga's checker (\S\ref{sec:recall-chk}) consists of 41 rule detectors,
one per rule except that the two directions of reflexivity share one,
and a validity test that requires, at every position of a coded list, a
genuine judgment code on which some detector fires against the entries
after it.  The internal-reasoning paper reads the detector list row by
row as the checker's 41-row node table
([\citen{conrga-paper}, Appendix ``The fragment contains the accepted $\forall$-free codes'']).  This paper calls the list of
detectors, and the index tree compiled from the list, a table, so that
\rga's table is that node table.  A detector, and a table, reads a
\emph{node}: a judgment together with the list of entries below it, the
node's tail.

\begin{defi}[\textsf{drvChkIdxG}, \textsf{iChkD}]
\label{def:ck-table}
A \emph{table} is an index tree $\mathcal D$ that reads a node's tail
and context and the node's conclusion, and \emph{fires} at genuine
codes when its HOL value is nonzero.  Inside \rga,
$\mathit{fires}_E(k)$ is the term $\cmp(E)\cdot
\mathit{pair}(\mathit{tail}(k), \mathit{ctx}(k))\cdot
\mathit{concl}(k) \ne \z$ (Appendix~\ref{app:ck-shift}), and \rga may
derive it at open codes.  The table $\mathcal D_L$ of a list $L$ of detectors is the
nested disjunction of $L$'s members.  At a table $\mathcal D$, the
validity test $\mathit{TreeValid}_{\mathcal D}$, the proof test
$\pi_{\mathcal D}$ (``$p$ codes a valid list whose head judgment has
code $h$'') and the checker $\Chk_{\mathcal D}$ are the
internal-reasoning paper's $\mathit{TreeValid}$, $\pi$ and $\Chk$
([\citen{conrga-paper}, definition~\textsf{iChk}] and [\citen{conrga-paper}, definition~\textsf{iHasCert}])
with $\mathcal D$ in place of \rga's rule test.  So $\Chk_{\mathcal D}(p, c)$ is the conjunction of
$\mathit{TreeValid}_{\mathcal D}(p)$, $p \ne 0$, ``the conclusion of the
list coded by $p$ has no hypotheses'' and ``that conclusion's code
equals $c$'', and only the first conjunct mentions the table.  The
certificate statement at $\mathcal D$ is $\Cert_{\mathcal D}(h) \ldef
\exists p.\ \cmp(\pi_{\mathcal D})\cdot h\cdot p \ne \z$, as in
[\citen{conrga-paper}, definition~\textsf{iHasCert}].
\end{defi}

\noindent
The existential in $\Cert_{\mathcal D}$ abbreviates a negated universal,
$\neg\forall p.\,\neg(\cdots)$, as everywhere in \rga; a derivation of
$\Cert_{\mathcal D}(h)$ need not exhibit a witness.

\begin{thm}[\textsf{chkG\_\allowbreak bridge\_\allowbreak drvChkIdx}, \textsf{chkG\_\allowbreak bridge\_\allowbreak iChk}]
\label{thm:ck-bridge}
Let $L_{\rga}$ be \rga's list of 41 detectors.  Then $\mathcal
D_{L_{\rga}}$ is \rga's rule test, and $\Chk_{\mathcal D_{L_{\rga}}}$
is $\Chk$.  Likewise the validity test, the proof test and the
certificate statement at $\mathcal D_{L_{\rga}}$ are \rga's.
\end{thm}

\noindent
The theorem's equalities are identities of HOL terms, proved by
unfolding.  The generic checker at \rga's table is therefore \rga's
checker, term for term, and adds no definition to what \rga's checker
rests on.  Thirty-eight laws of \rga's checker carry over to every
table as twin lemmas, each named after its \rga original with the
prefix \textsf{G\_}: the unfolding laws of the validity test, the dispatch on
the detectors of an arbitrary list, the append and graft of accepted
lists, and the certificate laws.  Of the 38 twins, 33 need no
hypothesis on the table, and five need the shift law of \S\ref{sec:ck-shift}
(Appendix~\ref{app:ck-twins}).

\subsection{The shift law}
\label{sec:ck-shift}

A detector satisfies the \emph{shift law} if its firing is preserved
when a node's tail is extended or rearranged, as long as every entry
that the tail had is still there.

\begin{defi}[\textsf{tblShift}]
\label{def:ck-shift}
A detector $E$ satisfies the shift law if the following holds for every
context $\Gamma$, all nodes $k$ and $k'$, and every map $\mathit{sh}$ of
positions.  Suppose that $k$ and $k'$ have the same head judgment, and
that every position $z$ below the length of $k$'s tail goes to a
position $\mathit{sh}(z)$ below the length of $k'$'s tail that holds
the same entry.  Then $\Gamma \vdash \mathit{fires}_E(k)$ gives $\Gamma \vdash
\mathit{fires}_E(k')$.  The suppositions are derivations: $\Gamma$
derives $k\jnat$, $k'\jnat$ and the equality of the two head
judgments, and the three suppositions about positions are derived at
an abstract position $z$, in every context that derives whatever
$\Gamma$ derives, where $\mathit{sh}$ also maps every natural to a
natural.
\end{defi}

\noindent
$\mathit{fires}_E(k)$ is the test that $E$ fires at the node $k$, as
defined before Theorem~\ref{thm:d-cons}.  The law is one-directional: the tail of $k'$
may be longer than the tail of $k$ and hold other entries, and
$\mathit{sh}$ need not be injective or onto.  The law is a closure
property of \rga's derivability, a HOL predicate of the detector; its
instances are HOL implications between derivability judgments at open
codes.  The law says nothing about which rule the detector recognizes,
and a detector that never fires satisfies it.
Appendix~\ref{app:ck-shift} states the law in full.

\begin{thm}[\textsf{tblShift\_\allowbreak orNest}, \textsf{tblShift\_\allowbreak RGA}]
\label{thm:ck-shift}
\leavevmode
\begin{enumerate}
\item If every member of a list $L$ satisfies the shift law, so does
  the table $\mathcal D_L$.
\item Every detector of $L_{\rga}$ satisfies the shift law, and so does
  \rga's table $\mathcal D_{L_{\rga}}$.
\end{enumerate}
\end{thm}

\noindent
Part~1 is \rga's dispatch proof run over an arbitrary list, with the
case split over \rga's 41 detectors replaced by one use of each
member's law.  Part~2 reads \rga's existing shift facts for its
detectors as instances of the predicate and then applies part~1.

\paragraph{Why one law suffices.}
In the derivations of the generic layer, a property of the table is
needed at one place only, the \emph{graft}, which lays the nodes of an accepted list on
top of another accepted list.  Each moved node then sits on a longer
tail whose old entries keep their positions, so the shift law applies
with $\mathit{sh}$ the identity.  The unfolding laws, the dispatch, the
certificate laws that do not graft, and the embedding of
\S\ref{sec:ck-embed} never use the law.  The scope of ``one law'' is
this derivation layer.  Closed faithfulness (\S\ref{sec:ck-faith})
needs a different hypothesis, that the table decides a step relation at
genuine codes, and the embedding needs only that a list contain \rga's
detectors.  The shift law and faithfulness are separate obligations,
and the development proves each one for each detector without using
the other.

\begin{thm}[\textsf{G\_\allowbreak iTreeValid\_\allowbreak appC}]
\label{thm:ck-graft}
Let the table $\mathcal D$ satisfy the shift law.  If $\Gamma$ derives
$A\jnat$, $B\jnat$, $\mathit{TreeValid}_{\mathcal D}(A)$ and
$\mathit{TreeValid}_{\mathcal D}(B)$, then $\Gamma \vdash
\mathit{TreeValid}_{\mathcal D}(\mathit{app}(A, B))$.
\end{thm}

\noindent
At \rga's table, Theorem~\ref{thm:ck-graft} is the validity part of
Theorem~\ref{thm:d-appC} on appended accepted proofs, rederived
through the shift law.  The graft
packages the theorem for the constructions that emit new nodes.
Given an accepted non-zero list $P$ with head $h$ and an accepted list
$T$, the graft hands any such construction the
appended list $\mathit{app}(P, T)$.  The theorem is a schema, and each
instance is an \rga derivation at open codes, in any context.  The
certificate laws carry over in the same way.  An accepted non-zero list
with head $h$ gives $\Cert_{\mathcal D}(h)$, and the elimination,
lifting, transport and substitution laws of $\Cert_{\mathcal D}$ hold
at every table, without the shift law.

\subsection{The embedding of \rga's trees}
\label{sec:ck-embed}

A tree that \rga's validity test accepts is accepted at the table of
every list that contains \rga's detectors.

\begin{thm}[\textsf{chkG\_\allowbreak iTreeValidD\_\allowbreak mono}, \textsf{chkG\_\allowbreak embed}]
\label{thm:ck-embed}
\leavevmode
\begin{enumerate}
\item Let $\mathcal D_1$ and $\mathcal D_2$ be tables such that, in every
  context $\Delta$ and at all $u$ and $y$ with $\Delta \vdash u\jnat$
  and $\Delta \vdash y\jnat$, $\Delta \vdash \cmp(\mathcal D_1)\cdot
  u\cdot y \ne \z$ gives $\Delta \vdash \cmp(\mathcal D_2)\cdot u\cdot
  y \ne \z$.  If $\Gamma \vdash m\jnat$ and $\Gamma \vdash
  \mathit{TreeValid}_{\mathcal D_1}(m)$, then $\Gamma \vdash
  \mathit{TreeValid}_{\mathcal D_2}(m)$.
\item Let the list $L$ contain every detector of $L_{\rga}$.  If $\Gamma
  \vdash m\jnat$ and $\Gamma \vdash \mathit{TreeValid}(m)$, then $\Gamma
  \vdash \mathit{TreeValid}_{\mathcal D_L}(m)$.
\end{enumerate}
\end{thm}

\noindent
Part~2 follows from part~1, since a fire of a member of $L_{\rga}$ is a
fire of the larger disjunction $\mathcal D_L$.  The theorem asks nothing
of the detectors that $L$ adds: no shift law, no faithfulness and no
soundness.  The build checks that the proof of part~2 reaches no shift
law (Appendix~\ref{app:ck-checks}).  Both parts are schemas, with
instances at open codes, so the embedding also works under
eigenvariables.

\begin{thm}[\textsf{chkG\_\allowbreak embed\_\allowbreak cert}, \textsf{chkG\_\allowbreak iChkD\_\allowbreak embed}]
\label{thm:ck-embed-cert}
Let the list $L$ contain every detector of $L_{\rga}$.
\begin{enumerate}
\item If $\Gamma \vdash h\jnat$ and $\Gamma \vdash \Cert(h)$, then
  $\Gamma \vdash \Cert_{\mathcal D_L}(h)$.
\item If $\Gamma \vdash p\jnat$ and $\Gamma \vdash \Chk(p, c)$, then
  $\Gamma \vdash \Chk_{\mathcal D_L}(p, c)$.
\end{enumerate}
\end{thm}

\noindent
For part~1, the proof eliminates the existential of $\Cert(h)$ at a
fresh variable, embeds the tree by Theorem~\ref{thm:ck-embed}(2), and
reintroduces the existential at the larger table.  For part~2, the
validity conjunct embeds by Theorem~\ref{thm:ck-embed}(2), and the
other three conjuncts of the checker do not mention the table.  The
embedding serves the D3, D2 and ex falso laws at \iga
(Section~\ref{sec:iga-goedel}), and the passage from \iga's consistency
sentence to \rga's in Section~\ref{sec:adopt} on adoption.

Only trees, certificates and verdicts embed.  The internal completeness
theorem of the internal-reasoning paper ([\citen{conrga-paper}, theorem~\textsf{rga\_complete\_int2}])
builds its certificates with emission laws that add one checker node
per rule.  These laws read \rga's tree acceptance, with membership in
\rga's list, so they are relative to \rga's table and do not embed.  At every table that satisfies the shift law, the generic
layer has the general emission law that emits an accepted tree through
the graft.  The generic layer has no per-rule emission laws at a
general table, neither copies of \rga's nor laws for the rules that an
extended table adds; such laws would serve internal completeness at
an extension.

\subsection{Closed faithfulness at a table}
\label{sec:ck-faith}

A table $\mathcal D$ \emph{decides} a step relation $R$ at genuine codes
when, for every proof $P$, context $\Gamma$ and term $c$, the table
fires at $P$'s code, $\Gamma$'s code and $c$'s code exactly when $R$
relates $P$, $\Gamma$ and $c$, that is, when the entry $\Gamma \vdash
c$ follows from the entries of $P$ by one step of $R$.  This hypothesis
is about one node against its tail, not about whole proofs.  An
$R$-proof is a proof list each of whose entries follows from the
entries after it by $R$.

\begin{thm}[\textsf{chkG\_\allowbreak iChkD\_\allowbreak faithful}, \textsf{chkG\_\allowbreak iChkD\_\allowbreak genuine\_\allowbreak proves}]
\label{thm:ck-faith}
Let the table $\mathcal D$ decide a step relation $R$ at genuine codes.
\begin{enumerate}
\item For every proof $P$ and term $c$: $\emptyset \vdash \Chk_{\mathcal
  D}(\bquo{P}, \bquo{c})$ if and only if $P$ is an $R$-proof of
  $\emptyset \vdash c$.
\item For every natural $p$ and term $c$: if $\emptyset \vdash
  \Chk_{\mathcal D}(\tofn{p}, \bquo{c})$, then $p = \quo{P}$ for some
  $R$-proof $P$ of $\emptyset \vdash c$.
\end{enumerate}
\end{thm}

\noindent
Theorem~\ref{thm:ck-faith} is [\citen{conrga-paper}, theorem~\textsf{iChk\_faithful}] at a table.
Its proof follows the internal-reasoning paper's.  At a deciding table,
the HOL value of the compiled validity test at $P$'s code is $1$ exactly
when $P$ is $R$-valid.  The node-code test alone, with no hypothesis on
the table, makes every number that the validity test accepts the code
of a proof list.  The two-way closed evaluation of \rga then turns these
HOL values into derivations and back.  That evaluation rests on \rga's
soundness and completeness at closed codes, so, as proved, both parts
of the theorem depend on \rga's soundness.  Part~1 is schematic and
part~2 is meta, as in the internal-reasoning paper.  No shift law enters.
The direction from an $R$-proof to its acceptance also has a proof by
one-directional closed evaluation, which avoids \rga's soundness, and
the proof of the second theorem at \iga uses that form
(Theorem~\ref{thm:ig-closed}).

\subsection{\iga's checker}
\label{sec:ck-iga}

\iga's checker adds one detector, the ground-induction tag, to \rga's
list.  The tag's design makes its shift law need no internal
reasoning about the tag's own computation.

\begin{defi}[\textsf{chkGfixIdx}, \textsf{drvChkListI}]
\label{def:ck-iga}
The \emph{ground-induction tag} $E_{\mathrm{gi}}$ is \rga's
three-premise citation former applied to a payload that reads the
node's context and conclusion and the three cited entries, and never
the node's tail.  The payload searches for the eigenvariable $i$ up to
the code of the second cited entry's context, and recomputes from the
codes the rule's three premises and its side conditions, positivity
included (Appendix~\ref{app:ck-tag}).  \iga's list is $L_{\iga} =
L_{\rga}$ with $E_{\mathrm{gi}}$ appended at the end, and \iga's
checker is $\ChkI = \Chk_{\mathcal D_{L_{\iga}}}$.
\end{defi}

\noindent
\iga's table thus has 42 detectors for \iga's 43 rules.  The payload is
\emph{table-blind}: it sees the three cited entries but never the tail
itself, so its value does not depend on where in the proof those
entries sit.  Appending the tag at the end keeps the tag numbers of \rga's
rows.  No theorem depends on the order, since every generic statement
is about the members of the list.

The positivity side condition needs a decider.  The decider is an
explicit index tree, so the tag's definition unfolds to explicit trees
all the way down (Appendix~\ref{app:ck-tcb}).

\begin{thm}[\textsf{pdecF\_\allowbreak eq}, \textsf{pr2fi\_\allowbreak pdecPIdx}]
\label{thm:ck-pos}
\leavevmode
\begin{enumerate}
\item The positivity recursion $\mathit{pdec}$ computes, for a natural
  $k$ and a term $t$, the five booleans
  \[
    \bigl[\,\mathit{posP}\,k\,t,\ \mathit{negP}\,k\,t,\ \mathit{pnexp}\,k\,t,\
    \exists b.\ t = \lambda b \wedge \mathit{posP}\,(k{+}1)\,b,\
    \exists b.\ t = \lambda b \wedge \mathit{negP}\,(k{+}1)\,b\,\bigr],
  \]
  where $\mathit{posP}$ and $\mathit{negP}$ are positivity and
  negativity in the variable $k$ and $\mathit{pnexp}$ is the auxiliary
  test of the positivity grammar.
\item The explicit index tree of the decider is an instance of the code
  library's shared course-of-values recursion over term codes, and the
  code of $\mathit{pdec}$'s result at $k$ and $t$ is read off the
  tree's HOL value at $k$ and $\quo{t}$.
\end{enumerate}
\end{thm}

\noindent
At a node with $n$ bound variables, the recursion evaluates the node's
children at $k + n$.  The theorem's statements are HOL equations about
codes.  The two sessions prove no law of the decider at open codes, and
the tag needs none.

\begin{thm}[\textsf{chkGfix\_\allowbreak fires\_\allowbreak iff}, \textsf{drvChkI\_\allowbreak fires\_\allowbreak iff}]
\label{thm:ck-iga-fires}
Let $P$ be a proof, $\Gamma$ a context and $c$ a term.
\begin{enumerate}
\item $E_{\mathrm{gi}}$ fires at the codes of $P$, $\Gamma$ and $c$ if
  and only if the entry $\Gamma \vdash c$ is the conclusion of an
  instance of ground induction whose three premises are entries of $P$.
\item $\mathcal D_{L_{\iga}}$ fires at the codes of $P$, $\Gamma$ and $c$
  if and only if the entry $\Gamma \vdash c$ follows from the entries of
  $P$ by one of \iga's rules.
\end{enumerate}
\end{thm}

\noindent
Part~1 is the correctness of the citation former with the bounded
search for $i$.  The bound holds because $i\jnat$ is a member of the
second premise's context.  Part~2 splits the fire of the disjunction
into a fire of a member of $L_{\rga}$, which \rga's rule test decides
exactly, or of the tag, and uses that \iga's rules are \rga's rules
together with ground induction.  Both parts are HOL facts about codes.
Part~2 says that \iga's table decides \iga's step relation.

\begin{thm}[\textsf{tblShift\_\allowbreak chkGfix}, \textsf{tblShift\_\allowbreak IGA}]
\label{thm:ck-iga-shift}
$E_{\mathrm{gi}}$ satisfies the shift law, and so does \iga's table
$\mathcal D_{L_{\iga}}$.
\end{thm}

\noindent
\rga's citation former already has a shift lemma with two obligations
on its payload: an internal reading equation for the payload, and
congruence of the payload in the context, the conclusion and the cited
entries.  For a table-blind payload, the reading equation drops the
tail, and the congruence is one congruence step of the compiler.
Neither obligation looks inside the payload, so the tag's shift law
needs no law of the positivity decider.  The law for the table then
follows from Theorem~\ref{thm:ck-shift}.

\begin{thm}[\textsf{iChkI\_\allowbreak faithful}, \textsf{iChkI\_\allowbreak genuine\_\allowbreak proves}]
\label{thm:ck-iga-faith}
\leavevmode
\begin{enumerate}
\item For every proof $P$ and term $c$: $\emptyset \vdash \ChkI(\bquo{P},
  \bquo{c})$ if and only if $P$ is an \iga proof of $\emptyset \vdashI
  c$.
\item For every natural $p$ and term $c$: if $\emptyset \vdash
  \ChkI(\tofn{p}, \bquo{c})$, then $p = \quo{P}$ for some \iga proof $P$
  of $\emptyset \vdashI c$.
\end{enumerate}
\end{thm}

\noindent
The theorem is Theorem~\ref{thm:ck-faith} at \iga's table, with
Theorem~\ref{thm:ck-iga-fires}(2) for the hypothesis.  The derivations
are \rga's: \rga derives that \iga's checker accepts.  Part~1 says that
\iga's checker accepts exactly the codes of \iga's proofs from no
hypotheses, and part~2 that every numeral it accepts is such a code.
The build checks that the proof goes through the tag's correctness and
the positivity decider's HOL semantics (Appendix~\ref{app:ck-checks}).
Read classically, both parts rest on \rga's soundness in HOL as
proved, as in Theorem~\ref{thm:ck-faith}; part~2, the genuineness, has
no other proof.  The theorem is not an \rga derivation of a sentence
that says \iga's checker is sound.  The theorem also says nothing
about $\ChkI$ at open codes; the decidedness of $\ChkI$'s verdict at
abstract arguments is part of Theorem~\ref{thm:ig-witness}(1).

At \iga's table the append and the graft hold with no shift premise,
by Theorem~\ref{thm:ck-iga-shift}; the embedding needs none at any
table (Theorem~\ref{thm:ck-embed}).  The build checks that \iga's graft
comes from the generic graft and the members' shift laws, not from
\rga's whole-table graft.  In particular, Theorem~\ref{thm:ck-embed-cert}(2)
at $L_{\iga}$ gives $\Gamma \vdash \ChkI(p, c)$ from $\Gamma \vdash
p\jnat$ and $\Gamma \vdash \Chk(p, c)$.

The development builds two pieces of internal machinery nowhere at
\iga's table: the per-rule emission laws of \S\ref{sec:ck-embed}, and
laws of the tag and of the positivity decider at open codes.
G\"odel's second theorem at \iga needs neither piece
(Section~\ref{sec:iga-goedel}).  \iga's checker adds to \rga's
definitions the generic checker's (Definition~\ref{def:ck-table}), the
tag with its payload and the positivity decider's trees, and its
meaning rests on \iga's rules; Appendix~\ref{app:ck-tcb} lists them.

\section{G\"odel's second theorem at \iga}
\label{sec:iga-goedel}

Inductive Grounded Arithmetic~\cite{iga-paper}, \iga, is \rga with
one more rule, ground induction (Section~\ref{sec:checker}).  Section~\ref{sec:locale} proves
G\"odel's second theorem once, for every G\"odel extension of \rga: an
entailment on \rga's term language, closed under \rga's rules, together
with a checker that meets the G\"odel conditions.
Section~\ref{sec:checker} builds the checker at a table of rule
detectors, and \iga's checker $\ChkI$, the checker at \rga's table with one
detector added for ground induction.  This section shows that \iga with
the checker $\ChkI$ is a G\"odel extension of \rga, and so obtains the
second theorem at \iga (Theorem~\ref{thm:ig-g2}): \iga does not derive
$\ConI$, its consistency sentence for the checker $\ChkI$.

Most of the work is done once, at an arbitrary table.  The transformer
laws for necessitation and EXPL hold at every table that satisfies the
shift law of \S\ref{sec:ck-shift} and contains the detectors that
the laws use; necessitation also needs a table that decides a step
relation (\S\ref{sec:ig-laws}).  A
proof of the table's step relation is accepted by one closed
evaluation, without \rga's soundness (\S\ref{sec:ig-closed}).  The
certificate of the third derivability condition is a tree of \rga's
own rules at every table, and it is accepted at every table whose list
contains \rga's (\S\ref{sec:ig-d3}).  What
is specific to \iga comes last: the witnesses and the theorems
(\S\ref{sec:ig-instance}, \S\ref{sec:ig-theorems}).  Throughout, $\vdash$
is \rga's entailment, $\vdashI$ is \iga's, and
Appendix~\ref{app:iga-goedel} lists the remaining laws, the witnesses, the
theorems that the instance yields and the build checks.

\subsection{\iga's entailment as a calculus}
\label{sec:ig-calc}

The G\"odel conditions sit on a calculus: an entailment on \rga's terms
that is closed under \rga's rules and under hypothesis stripping
(Section~\ref{sec:locale}).  \iga's entailment is such a calculus.

\begin{thm}[\textsf{gig\_\allowbreak prf\_\allowbreak strip}, \textsf{gig\_\allowbreak gcalc\_\allowbreak strip}]
\label{thm:ig-strip}
\leavevmode
\begin{enumerate}
\item If a valid \iga proof contains the judgment $H \vdash c$, then
  \iga derives $c$ from the well-formed members of $H$.
\item \iga's entailment is an instance of the calculus with hypothesis
  stripping of Section~\ref{sec:locale}.
\end{enumerate}
\end{thm}

\noindent
Both parts are classical, HOL facts about \iga's derivations, with no
hypotheses.  \iga's step relation contains \rga's, is monotone in the
proof list and preserves well-formedness; these three properties give
the calculus.  Part~1 is an induction over the proof list.  \rga's rules
are replayed as in \rga's own stripping lemma.  The one new case is
ground induction.  The rule's closure premise carries one extra
hypothesis, the motive's instance; the proof keeps that hypothesis when
it is well formed and weakens it back in when it is not, and then
applies ground induction again.  The build checks that part~2 rests on
this proof and not on \rga's stripping lemma
(Appendix~\ref{app:iga-goedel-checks}).

\subsection{The transformer laws at any table}
\label{sec:ig-laws}

Fix a list $L$ of rule detectors, and write $\mathcal D$ for the table
that $L$ defines and $\Chk_{\mathcal D}$ for the checker at that table
(Section~\ref{sec:checker}).  $\Chk_{\mathcal D}(p, c)$ is the
conjunction of the same four decided tests as $\Chk(p, c)$
(\S\ref{sec:recall-chk}), with tree validity tested at $\mathcal D$.
A table decides a step relation $R$ in the sense of
\S\ref{sec:ck-faith}, and \iga's table decides \iga's step relation
(Theorem~\ref{thm:ck-iga-fires}(2)).

The laws below are \rga derivations at an abstract proof code.  Their
premises on the table are of three kinds: the shift law for $\mathcal
D$ of \S\ref{sec:ck-shift}, which the append of two accepted
trees needs; membership in $L$ of the \rga detectors for the rules that
the law's new entries use (cut, weakening, negation elimination); and,
where an \rga-accepted tree is embedded, that $L$ contains \rga's
list.  Necessitation also needs that $\mathcal D$ decides a step
relation.

\begin{thm}[\textsf{G\_\allowbreak cut1M\_\allowbreak chk}, \textsf{G\_\allowbreak goedelE\_\allowbreak chk}]
\label{thm:ig-laws}
Let the table $\mathcal D$ of $L$ satisfy the shift law.  Let $\Gamma$
be a context and $a$ a term with $\Gamma \vdash a\jnat$.
\begin{enumerate}
\item (NEC.)  Let $\mathcal D$ decide a step relation $R$, let $L$
  contain the cut detector, and let $Q$ be an $R$-proof of $\{A\} \vdash
  B$.  If $\Gamma \vdash \Chk_{\mathcal D}(a, \bquo{A})$, then $\Gamma
  \vdash \Chk_{\mathcal D}(c_Q\cdot a\cdot a, \bquo{B})$.
\item (EXPL.)  Let $L$ contain the negation-elimination detector, and
  let $b$ be a term with $\Gamma \vdash b\jnat$.  If
  $\Gamma \vdash \Chk_{\mathcal D}(a, \bquo{g})$ and $\Gamma \vdash
  \Chk_{\mathcal D}(b, \bquo{\neg g})$, then $\Gamma \vdash \Chk_{\mathcal
  D}(e\cdot a\cdot b, \bquo{\z = \tofn 1})$.
\end{enumerate}
Here $c_Q$ and $e$ are built by the term formers of
Theorems~\ref{thm:g-cut} and~\ref{thm:g-expl}, with $Q$ now a proof of
$R$.
\end{thm}

\noindent
Each part is a rule in every context of \rga, schematic in $Q$, $A$,
$B$ and $g$; read as derivations in \rga, the parts are schemas, each instance one \rga derivation.  The proofs follow those of Theorems~\ref{thm:g-cut}
and~\ref{thm:g-expl}.  Each law appends the
accepted inputs, and for NEC the code of $Q$, by the generic graft,
which needs the shift law (Theorem~\ref{thm:ck-graft}).  The law then
puts one new entry in front, justified by cut for NEC and by negation
elimination for EXPL, through that detector's introduction law and the
law for adding an entry whose detector belongs to $L$.  A transport
along the transformer's value equation finishes the proof.  No entry of
an input tree is examined again, so the laws do not care whether an
input contains a ground-induction entry.  Closed evaluation shows that
the code of $Q$ is a valid tree at $\mathcal D$ whose head is $\{A\}
\vdash B$ (the argument of Theorem~\ref{thm:ig-closed} below, at a
head with a hypothesis).  The output
trees contain the code of $Q$ and the input codes, so they are not
\rga trees.  What the laws share with \rga is the term formers.

Three more laws of the same pattern serve the optional layers of
Section~\ref{sec:locale}: ex falso, D2 and necessitation with two
premises.  Each of these laws puts one to three new entries over an
append, and
Theorems~\ref{thm:ig-laws-exf} and~\ref{thm:ig-laws-nec2} of
Appendix~\ref{app:iga-goedel} state them with their premises.

\subsection{Closed acceptance without soundness}
\label{sec:ig-closed}

\begin{thm}[\textsf{chkG\_\allowbreak iChkD\_\allowbreak proves\_\allowbreak sf}, \textsf{chkI\_\allowbreak proves\_\allowbreak sf}]
\label{thm:ig-closed}
\leavevmode
\begin{enumerate}
\item Let the table $\mathcal D$ decide a step relation $R$.  If $P$ is
  an $R$-proof of $\emptyset \vdash c$, then $\Gamma \vdash \Chk_{\mathcal
  D}(\bquo{P}, \bquo{c})$ for every context $\Gamma$.
\item If $P$ is an \iga proof of $\emptyset \vdashI c$, then $\Gamma
  \vdash \ChkI(\bquo{P}, \bquo{c})$ for every context $\Gamma$.
\end{enumerate}
\end{thm}

\noindent
Both parts are classical schemas: HOL quantifies over the genuine
proof $P$, and the conclusion is an \rga derivation.  The compiled validity test at
$\mathcal D$ has the HOL value $1$ at the code of $P$, because $P$ is
valid and $\mathcal D$ decides $R$.  The closed-numeral law for
compiled functions ([\citen{conrga-paper}, theorem~\textsf{cmp\_eq\_num}]) derives that value
equation in \rga, and weakening moves the equation into $\Gamma$.  The
other three tests are computations on numerals.  Part~2 is part~1 at
\iga's table.  The same argument, applied to an $R$-proof whose head
has hypotheses, gives the validity of $Q$'s code at $\mathcal D$ and
its head, which Theorem~\ref{thm:ig-laws}(1) needs.

The evaluation runs in one direction only, from the value in HOL to
the equation in \rga.  The closed faithfulness of $\ChkI$
(Theorem~\ref{thm:ck-iga-faith}) also yields this acceptance, but
faithfulness is an equivalence, and the equivalence's other direction,
which reads a value off an equation that \rga derives, rests on \rga's
soundness.  A dependency check sees whole theorems, not directions, so
a proof that cites the equivalence for one direction has the soundness
theorem in its dependency cone.  The witnesses below therefore use
Theorem~\ref{thm:ig-closed}, and the build checks that the internal
theorems of \S\ref{sec:ig-theorems} avoid the soundness theorems.
Theorem~\ref{thm:ig-closed} does not say that acceptance gives a
genuine proof; that converse is the checker's genuineness, which rests
on soundness.

\subsection{D3 at every table}
\label{sec:ig-d3}

For \rga, the substantial piece of the second theorem's proof is the
certificate of the third derivability condition
(\S\ref{sec:goedel-d3g}).  At a table, the same certificate builder
serves, and the certificate is still a tree of \rga's rules.
For a table $\mathcal D$ and a code term $S$, write $\PrC_{\mathcal
D}(S) = \neg\forall q.\,\neg\Chk_{\mathcal D}(q, S)$.

\begin{thm}[\textsf{G\_\allowbreak iChk\_\allowbreak d3SigmaAt}, \textsf{chkG\_\allowbreak iChkD\_\allowbreak d3SigmaAt}]
\label{thm:ig-d3}
Let $\mathcal D$ be any table, and let $S$ be a closed code term with a
source index, as every code numeral $\bquo{\varphi}$ has.  There is a
term $\sigma_{\mathcal D, S}$ such that, for every context $\Gamma$
and every term $v$ with $\Gamma \vdash v\jnat$ and $\Gamma \vdash
\Chk_{\mathcal D}(v, S)$:
\begin{enumerate}
\item $\Gamma \vdash \Chk(\sigma_{\mathcal D, S}\cdot v,
  \bquo{\PrC_{\mathcal D}(S)})$, at \rga's own checker;
\item if $\mathcal D$ is the table of a list that contains \rga's list,
  then $\Gamma \vdash \Chk_{\mathcal D}(\sigma_{\mathcal D, S}\cdot v,
  \bquo{\PrC_{\mathcal D}(S)})$.
\end{enumerate}
\end{thm}

\noindent
Read as derivations in \rga, the theorem is a schema, each instance one \rga derivation.
Part~1 has no premise on the table.  The output $\sigma_{\mathcal D,
S}\cdot v$ is a tree of \rga's rules, accepted by \rga's own checker.
The tree certifies the value of the checker formula $\Chk_{\mathcal
D}(v, S)$ at $v$, one conjunct at a time, as in
Theorem~\ref{thm:g-d3-any}.  The validity conjunct is the compiled
validity test of $\mathcal D$.  The certificate builder of
\S\ref{sec:goedel-d3g} is proved for every index tree, by HOL
induction on the tree, so the builder certifies $\mathcal D$'s test as
it certifies \rga's.  The other three conjuncts do not depend on the
table.  The input $v$ enters only through the four conjuncts of the
checker formula, and no entry of $v$'s tree, a ground-induction entry
or any other, is checked again.  Part~2 embeds the \rga tree into the
larger table.

So the D3 law holds at every table whose list contains \rga's, for
checkers of the compiled form $\Chk_{\mathcal D}$.  Moving the law
into an extension's entailment still needs that entailment to contain
\rga's (\S\ref{sec:ig-instance}).  The generic law has no premise on
the table, and the build checks that the \iga instance goes through it
and uses no fact about the ground-induction detector or \iga's step
relation (Appendix~\ref{app:iga-goedel-checks}).  The theorem says nothing
about the size of the certificate.

\subsection{\iga as a G\"odel extension of \rga}
\label{sec:ig-instance}

\begin{defi}[\textsf{ConIGA}]
\label{def:ig-con}
The consistency sentence of \iga for the checker $\ChkI$ is
\[
  \ConI \;\ldef\; \forall p.\ \neg\ChkI(p, \bquo{\z = \tofn 1}) .
\]
\end{defi}

\noindent
$\ConI$ is $\ConR$ with $\ChkI$ in place of $\Chk$, and it is
literally the consistency sentence of Section~\ref{sec:locale} at the
checker $\ChkI$: the two terms are equal in HOL.  $\ConI$ is closed and
well formed.  Every instance of $\ConI$'s matrix is decided: from
$\Gamma \vdashI p\jnat$ follows $\Gamma \vdashI \ChkI(p,
\bquo{\z = \tofn 1})\jbool$ (Theorem~\ref{thm:ig-witness}(1)).  So
$\ConI$ is a universal over decided instances.  On the classical
reading, $\ConI$ expresses the consistency of \iga, because $\ChkI$ is
faithful at closed codes (Theorem~\ref{thm:ck-iga-faith}).

\begin{thm}[\textsf{g2\_\allowbreak syn\_\allowbreak gi}, \textsf{g2\_\allowbreak con\_\allowbreak gi}]
\label{thm:ig-witness}
\iga's entailment with the checker $\ChkI$ meets the G\"odel conditions
of the second theorem (Section~\ref{sec:locale}):
\begin{enumerate}
\item the three structural conditions: $\ChkI$ commutes with
  substitution, $\ChkI(p, c)$ is well formed exactly when $p$ and $c$
  are, and $\Gamma \vdashI p\jnat$ and $\Gamma \vdashI c\jnat$ give
  $\Gamma \vdashI \ChkI(p, c)\jbool$, over the calculus of
  Theorem~\ref{thm:ig-strip};
\item NEC, D3 and EXPL, as transformers in \iga;
\item D1, as closed acceptance: an \iga proof $P$ of $\emptyset
  \vdashI \varphi$ gives $\emptyset \vdashI \ChkI(\bquo{P},
  \bquo{\varphi})$;
\item consistency: $\emptyset \nvdashI \z = \tofn 1$.
\end{enumerate}
So \iga is a G\"odel extension of \rga.  \iga with $\ChkI$ also meets
the conditions of the optional layers: ex falso, D2 and necessitation
with two premises.
\end{thm}

\noindent
The theorem is classical, a HOL fact about \iga.  D2 is not among the conditions
of the second theorem.  The witnesses come from the laws above.  For
NEC the transformer is $v \mapsto c_Q\cdot v\cdot v$, with $Q$ an \iga
proof of $\{A\} \vdashI B$ obtained from the entailment, by
Theorem~\ref{thm:ig-laws}(1) at \iga's table.  For D3 it is
$\sigma_{\mathcal D, S}$ of Theorem~\ref{thm:ig-d3} at \iga's table and
the code numeral, with the output sentence matched to the shape that
the condition states.  For EXPL it is $e$, by
Theorem~\ref{thm:ig-laws}(2).  D1 is Theorem~\ref{thm:ig-closed}(2).
Consistency is \iga's consistency theorem, proved in HOL through
\iga's soundness for the $\omega$-model $\Mw$
([\citen{conrga-paper}, \S\,``The $\omega$-semantics''], \cite{iga-paper}).

The laws are \rga derivations from \rga hypotheses, in every context,
while the conditions are \iga entailments.  Cut transfer bridges the
gap.  A rule that holds for \rga's entailment in every context is
derived in \rga from its hypotheses, the derivation is embedded into
\iga, whose rules include \rga's, and a cut gives the \iga entailment.
The totality of each transformer transfers in the same way.  Among the
transformer witnesses, content specific to \iga enters in two places
only: NEC, and necessitation with two premises, cite a closed \iga
proof $Q$, whose code is a valid tree at the table because \iga's
table decides \iga's step relation;
and the shift law at \iga's table discharges the shift premise.
Outside the transformer witnesses, \iga enters through D1 (\iga's table
decides \iga's step relation), through consistency, through the
ground-induction case of Theorem~\ref{thm:ig-strip}, and through the
table itself, with its one added detector.

The conditions of Section~\ref{sec:locale} come in layers, and the
development interprets the layers at \iga one at a time, nine
interpretations in all.  Each interpretation rewrites the locale's
diagonal function to $\Diag$ and its consistency sentence to $\ConI$.
Every theorem of the layers thus becomes a theorem about \iga with
$\ChkI$ for the checker and $\ConI$ for the consistency sentence, and
Appendix~\ref{app:iga-goedel-reexports} lists the twenty that the
development states for \iga.

\subsection{The theorems}
\label{sec:ig-theorems}

Write $\PrI(c) = \neg\forall q.\,\neg\ChkI(q, c)$ for the provability
formula of \iga's checker.  \iga's G\"odel sentence $G_{\mathrm I}$
is $G$ of \S\ref{sec:goedel-sentences} with $\ChkI$ in place of
$\Chk$: $D \ldef \forall q.\,\neg\ChkI(q, \Diag(x, x))$, $S \ldef
\Diag(\bquo{D}, \bquo{D})$ and $G_{\mathrm I} \ldef \forall
q.\,\neg\ChkI(q, S)$.

\begin{thm}[\textsf{gi\_\allowbreak goedel2\_\allowbreak internal\_\allowbreak contra}, \textsf{gi\_\allowbreak goedel2\_\allowbreak internal}]
\label{thm:ig-internal}
\leavevmode
\begin{enumerate}
\item $\PrI(\bquo{\ConI}) \vdashI \neg\ConI$.
\item $\ConI \vdashI \neg\PrI(\bquo{\ConI})$.
\end{enumerate}
\end{thm}

\noindent
Both parts are derivations in \iga under one hypothesis.  They are the formalized second theorem for $\ChkI$:
Theorem~\ref{thm:lc-textbook}(4, 5) at \iga, the locale form of
Theorem~\ref{thm:g-formal}(1, 2).  From $\ConI \vdashI G_{\mathrm
I}$ (Theorem~\ref{thm:lc-textbook}(2) at \iga), necessitation gives
$\PrI(\bquo{\ConI}) \vdashI \PrI(\bquo{G_{\mathrm I}})$.  D3 at
$\bquo{G_{\mathrm I}}$, followed by necessitation of the entailment
$\PrI(\bquo{G_{\mathrm I}}) \vdashI \neg G_{\mathrm I}$, gives
$\PrI(\bquo{\neg G_{\mathrm I}})$ from the same hypothesis, and EXPL
gives $\PrI(\bquo{\z = \tofn 1})$, which is $\neg\ConI$.  Part~2
follows by contraposition at the $\Sigma$-shaped hypothesis
(Corollary~\ref{cor:rc-contrapos}, in the calculus) and double
negation.  The build checks that the dependency cones of both parts
avoid eleven soundness and consistency theorems, among them those of
\rga, \iga and \igaP and the genuineness of $\ChkI$
(Appendix~\ref{app:iga-goedel-checks}).

\begin{thm}[\textsf{gi\_\allowbreak goedelG\_\allowbreak refutes}, \textsf{gi\_\allowbreak goedel1}]
\label{thm:ig-g1}
\leavevmode
\begin{enumerate}
\item If $\emptyset \vdashI G_{\mathrm I}$, then $\emptyset \vdashI
  \z = \tofn 1$.
\item $\emptyset \nvdashI G_{\mathrm I}$.
\end{enumerate}
\end{thm}

\noindent
Both parts are classical.  Part~1 is a rule on theorems, a HOL fact
proved from D1 as Theorem~\ref{thm:g-g1}(2) is and as
Theorem~\ref{thm:lc-g2}(1) is in the locale, and the build checks
that its dependency cone avoids the soundness theorems.  Part~2 is a
meta theorem, part~1 with \iga's consistency.

\begin{thm}[\textsf{gi\_\allowbreak goedel2\_\allowbreak ConIGA}]
\label{thm:ig-g2}
$\emptyset \nvdashI \ConI$.
\end{thm}

\noindent
\emph{Proof.}  A derivation of $\ConI$, cut against $\ConI \vdashI
G_{\mathrm I}$ (Theorem~\ref{thm:lc-textbook}(2) at \iga), gives a derivation of $G_{\mathrm I}$, which
Theorem~\ref{thm:ig-g1}(2) excludes.  \qed

\medskip
\noindent
Theorem~\ref{thm:ig-g2} is classical, a meta theorem with no
hypotheses, the instance at \iga of Theorem~\ref{thm:lc-g2}(3).  Its
proof uses \iga's consistency, proved in HOL through \iga's soundness
for $\Mw$, as the proof of Theorem~\ref{thm:g-g2} uses \rga's.  The build checks that the theorem's dependency cone reaches
the \iga witnesses: the laws for NEC, D3 and EXPL at \iga's table,
Theorem~\ref{thm:ig-strip}, \iga's consistency and the rewriting of the
locale's consistency sentence to $\ConI$.

Theorem~\ref{thm:ig-g2} says that \iga does not derive $\ConI$, its
consistency sentence for the checker $\ChkI$.  That $\ConI$ expresses
\iga's consistency is the classical reading, through the faithfulness
of $\ChkI$.  No theorem of the development quantifies over checkers or
over consistency sentences at \iga.  The uniformity lies in
Section~\ref{sec:locale}: the locale's second theorem holds for every
pair of an entailment and a checker that meets the G\"odel conditions,
where the entailment is on \rga's language and closed under \rga's
rules.  Theorem~\ref{thm:ig-g2} does not say that \iga fails to derive
$\ConR$; that question is open (Section~\ref{sec:open}).

The adoption extension \igaP of \iga adds one rule, adoption keyed
to \iga's checker (Section~\ref{sec:adopt}).  Write $\vdashIP$ for its
entailment.

\begin{thm}[\textsf{gi\_\allowbreak separation\_\allowbreak IGA}, \textsf{gi\_\allowbreak adiU\_\allowbreak unprovable}]
\label{thm:ig-sep}
\leavevmode
\begin{enumerate}
\item $\emptyset \nvdashI \ConI$.
\item $\emptyset \vdashIP \ConI$.
\item For every context $\Gamma$ and term $c$, $\Gamma \vdashI c$ gives
  $\Gamma \vdashIP c$.
\item $\emptyset \nvdashIP \z = \tofn 1$.
\item \iga does not derive the universal sentence that adoption at \iga
  adopts.
\end{enumerate}
\end{thm}

\noindent
So \iga is a proper part of \igaP, with no hypotheses.  The theorem
is classical, a meta theorem.  Part~1 is Theorem~\ref{thm:ig-g2}, and parts~2--4
are the adoption results of Section~\ref{sec:adopt}.  Parts~1 and~4
rest on the consistency of \iga and of \igaP, both proved in HOL
through soundness for $\Mw$; the consistency of \igaP also rests on
the genuineness of $\ChkI$ (Section~\ref{sec:adopt}).  Part~5 follows from part~1, because in \rga
the adopted universal entails $\ConI$ (Section~\ref{sec:adopt}).  The
theorem does not separate \igaP from \rga through $\ConR$, which
Section~\ref{sec:adopt} does with Theorem~\ref{thm:g-g2}, and it does
not say that \iga fails to derive $\ConR$.  No checker for \igaP is
defined, so the second theorem for \igaP at its own checker is not
stated.

\paragraph{What the instance does not need.}
The laws never examine an entry of an input tree again, so the
instance needs no law about the ground-induction detector at open
codes.  It needs only the shift law at \iga's table and the fact that
\iga's table decides \iga's step relation, and D3 needs neither.  The
G\"odel conditions do not include internal completeness, the property
that the checker accepts an emitted proof of each true decided atom,
and this paper does not prove that property for $\ChkI$.  Soundness
enters the proofs of Theorems~\ref{thm:ig-g1}(2) and~\ref{thm:ig-g2}
only through consistency: every other witness lies in a dependency
cone that the build checks against the soundness theorems
(Appendix~\ref{app:iga-goedel-checks}).  The internal theorems avoid
soundness.

\section{Adoption: the reflection rule, and the consistency family}
\label{sec:adopt}

Gödel's second theorem says that \rga does not derive its consistency
sentence $\ConR$ (Theorem~\ref{thm:g-g2}).  One classical response to
that limit is reflection: a stronger system adds a principle saying
that what the weaker system proves is so.  This section adds such a
principle to \rga as a rule, called \emph{adoption}, and keys the rule
to \rga's own proof checker.  The result is a consistent system \rgaP
that derives $\ConR$, and so strictly extends \rga.  The same rule
keyed to \iga's checker gives \igaP, which derives the consistency
sentences of both \rga and \iga.  Reflection here is in the
proof-theoretic sense, a rule that extends the system.  The section
states the rule, its soundness, the derivations of the consistency
sentences, and the separations; Appendix~\ref{app:adopt} gives the
constructions in detail.  The last three subsections turn to a second
kind of adoption, the one that the IGA paper prices: \iga with a
supplied claim that a code is an ordinal notation in Kleene's
$\mathcal O$.  They state what falsifies such a claim, what G\"odel's
second theorem says about it, and why neither unrefutedness nor, for
well-foundedness sentences, consistency with the system is a closure
condition (\S\ref{sec:adopt-claims}--\S\ref{sec:adopt-ac}).

\subsection{The rule}
\label{sec:adopt-rule}

The rule's premise is a warrant, one sentence saying that, for every
number $v$, the checker accepts a computed proof of the $v$-th
instance of a universal.  An explicit primitive-recursive index tree
$F$ computes the proof, the kind of transparent witness of
\S\ref{sec:recall-transparent}.

\begin{defi}[\textsf{adrW}, \textsf{adr\_\allowbreak derives}]
\label{def:ad-rule}
Let $u$ be a closed, well-formed term, read as a predicate on numbers,
and let $F$ be an index tree.  Write $\Inst_u(v)$ for the code of $u$
applied to the numeral of $v$, a compiled function of $v$, whose value at a numeral $\tofn n$ is the code
numeral $\bquo{u\cdot\tofn n}$.  The \emph{warrant} of $u$ by $F$ is
the sentence
\[
  \WarR(F, u) \;=\; \forall v.\ \Chk\bigl(\cmp(F)\cdot v\cdot v,\ \Inst_u(v)\bigr) .
\]
\rgaP is \rga with one more rule, \adoptR: from $\Gamma \vdashRP
\WarR(F, u)$, for a closed well-formed $u$, conclude $\Gamma \vdashRP
\cA\cdot u$.  Its derivations are lists of judgments in which every
entry follows by a rule of \rga or by \adoptR; $\vdashRP$ is the
entailment of these derivations.
\end{defi}

\noindent
Since the checker accepts only proofs whose conclusion has no
hypotheses (\S\ref{sec:recall-chk}), the warrant says that each
instance $u\cdot\tofn n$ has a proof from the empty context, computed
by $F$ at $n$.  Since $u$ is closed, adoption has no parameters.  The
checker's verdict is decided at every argument
([\citen{conrga-paper}, theorem~\textsf{hq\_iChk}]), so the warrant is a universal of
decided instances about an explicit witness.  Nothing is searched, and
no ``there exists a proof'' occurs in the warrant.

\begin{thm}[\textsf{adr\_\allowbreak adopt}, \textsf{adr\_\allowbreak embed}]
\label{thm:ad-rule}
\leavevmode
\begin{enumerate}
\item If $\Gamma \vdashRP \WarR(F, u)$, $u$ is well-formed and $u$ is
  closed, then $\Gamma \vdashRP \cA\cdot u$.
\item If $\Gamma \vdash c$, then $\Gamma \vdashRP c$.
\end{enumerate}
\end{thm}

\noindent
Both parts are classical, HOL theorems about the derivability
relations, and \rgaP is closed under cut.

\paragraph{Comparison with the restricted $\omega$-rule.}
Shoenfield's restricted (recursive) $\omega$-rule~\cite{shoenfield59omega} infers
a universal from a recursive enumeration of proofs of its instances.
Iterated, that rule is complete for true arithmetic.  Adoption has the
same conclusion and a different premise, in three ways.
\begin{itemize}
\item The premise is internal: a finite derivation, inside the system,
  of the warrant, not an external fact about the enumeration.
\item The premise is decided: the warrant's matrix is a decided
  checker verdict.
\item The premise is keyed to a fixed base checker.  The rule may be
  applied any number of times, but every warrant speaks of \rga's
  checker, which accepts only \rga proofs.  The development defines no
  checker for \rgaP, so \rgaP has no rule about its own proofs, and
  adoption is not iterated.
\end{itemize}
Classically, the reflection principles of a theory assert that what
the theory proves is true.  Kreisel and L\'evy showed that uniform
reflection over Peano arithmetic is equivalent to transfinite
induction up to $\varepsilon_0$~\cite{kreisel68reflection}.  Turing's
ordinal logics~\cite{turing39ordinals} and Feferman's transfinite
progressions~\cite{feferman62transfinite} iterate reflection along
ordinal notations.  Adoption is neither a progression nor an iteration:
it is one rule, about one fixed checker.

\subsection{Soundness}
\label{sec:adopt-sound}

The semantics for soundness is the $\omega$-semantics of
[\citen{conrga-paper}, \S\,``The $\omega$-semantics''], written $\Mw$ here.  $\Mw$
is the least valuation in which a universal is true when every numeral
instance is true and false when some instance is false, so an
ungrounded term, such as $W\cdot W$ with $W = \lambda w.\ \forall x.\
w\,w$, gets no value.  A system is \emph{$\Mw$-sound} when every
derivation $\Gamma \vdash c$ with $\Gamma$ well-formed and every member
of $\Gamma$ true in $\Mw$ has a conclusion $c$ true in $\Mw$.  \rga is
$\Mw$-sound~\cite[Sect.~6]{rga-paper}.  $\Mw$-soundness is a property
relative to this one model; it is not arithmetical soundness over the
standard numbers in the classical sense.

\begin{thm}[\textsf{adr\_\allowbreak omega\_\allowbreak sound}, \textsf{adr\_\allowbreak consistent}]
\label{thm:ad-sound}
\leavevmode
\begin{enumerate}
\item If $\Gamma \vdashRP c$, $\Gamma$ is well-formed, and every
  member of $\Gamma$ is true in $\Mw$, then $c$ is true in $\Mw$.
\item $\emptyset \nvdashRP \z = \tofn 1$.
\end{enumerate}
\end{thm}

\noindent
\emph{Proof.}  The step to check is adoption.  Suppose the warrant $\WarR(F, u)$ is true
in $\Mw$, with $u$ closed and well-formed.  At each $n$, the warrant's
instance is a closed checker verdict with no quantifier, and the
instance is true in $\Mw$.
\begin{enumerate}
\item \rga derives that verdict, by its completeness for closed
  quantifier-free sentences (Appendix~\ref{app:adopt-sound}).
\item Closed evaluation turns the derived verdict into $\emptyset
  \vdash \Chk(\tofn p, \bquo{u\cdot\tofn n})$, where $p$ is the value
  of $F$ at $(n, n)$.
\item The checker is genuine ([\citen{conrga-paper}, theorem~\textsf{iChk\_faithful}]), so $p$ is the code of an \rga proof
  of $\emptyset \vdash u\cdot\tofn n$.
\item By \rga's $\Mw$-soundness, $u\cdot\tofn n$ is true in $\Mw$.
\end{enumerate}
Every instance is true, so the universal $\cA\cdot u$ is true in
$\Mw$.  Part~1 follows by induction over \rgaP derivations, with
\rga's per-step soundness for the steps by \rga's rules.  Part~2
follows, since $\z = \tofn 1$ is false in $\Mw$.  \qed

\medskip
\noindent
Theorem~\ref{thm:ad-sound} is classical, a HOL theorem about \rgaP.  The theorem
does not say that \rgaP proves its own soundness, and the theorem is not
a reflection statement inside \rgaP.  As a prediction that the
development does not prove, \rgaP and \igaP are expected to be
proof-grounded as well, that is, sound and complete for grounded truth
read with the system's own proofs as the certificates
(Theorem~\ref{thm:op-iga-pg} for \iga), since a derivation of the
warrant followed by adoption and universal elimination at a fresh
variable is a proof of the adopted universal's open instance.

\subsection{\rga derives a warrant}
\label{sec:adopt-warrant}

The rule needs a warrant that some system derives.  \rga alone
derives such a warrant, for a body whose universal is interderivable with $\ConR$.

\begin{defi}[\textsf{adrU}, \textsf{adrF}]
\label{def:ad-body}
Let $k$ be the code of the judgment $\emptyset \vdash \z = \tofn 1$,
and let $\mathit{test}(v) = \cmp(\mathit{provesChkIdx}) \cdot \tofn
k\cdot v$ be the checker's compiled test at $k$.  The body is
\[
  U \;=\; \lambda v.\ \bigl(\mathit{test}(v) = \z\bigr) ,
\]
which says that $v$ is not an accepted tree for $\emptyset \vdash \z =
\tofn 1$.  The witness $F_U$ is one index tree, a conditional on the
test: at $v$ with $\mathit{test}(v) = \z$ the tree computes a value
certificate of the instance, and at $v$ with $\mathit{test}(v) \neq
\z$ it computes an ex-falso tree (Definition~\ref{def:fa-witness}).
\end{defi}

\begin{thm}[\textsf{adr\_\allowbreak cert}, \textsf{adr\_\allowbreak warrant\_\allowbreak uC}]
\label{thm:ad-warrant}
\leavevmode
\begin{enumerate}
\item If $\Gamma \vdash v\jnat$, then $\Gamma \vdash
  \Chk\bigl(\cmp(F_U)\cdot v\cdot v,\ \Inst_U(v)\bigr)$.
\item $\emptyset \vdash \WarR(F_U, U)$.
\end{enumerate}
\end{thm}

\noindent
Part~2 is a derivation in \rga, from no hypothesis.  Part~1 is a derived rule,
whose instance at a free $v$, under the hypothesis $v\jnat$, is the
\rga derivation that $\forall$-introduction closes into part~2.
The proof of part~1 splits on the decided value of $\mathit{test}(v)$.
\begin{itemize}
\item At value $\z$, $F_U$ computes the value certificate of the
  equation $\mathit{test}(v) = \z$, which the builder of
  Theorem~\ref{thm:d-builder} supplies for every index tree.  A node of the
  rule trueEqI (from $c$ conclude $c = \tofn 1$), two backward head steps to
  the instance code and a node of the rule trueEqE (its converse) turn that certificate into an accepted
  proof of the instance $U\cdot v$.
\item At a nonzero value, $v$ is an accepted proof of $\z = \tofn 1$.
  $F_U$ appends to $v$ the fixed proof of $\neg(\z = \tofn 1)$ and one
  root that concludes the instance, an ex-falso tree at the level of
  codes (Theorem~\ref{thm:fa-exf}).
\end{itemize}
The split is on the decided value of the checker's test at the proof
code $v$, and never on the truth of the judgments inside the proof
that $v$ codes.

The body is the compiled test rather than $\neg\Chk(v, \bquo{\z =
\tofn 1})$ because the witness needs a primitive-recursive certificate
in the branch where the body holds.  The test equation is a single
compiled application, which the builder certifies directly.  For the
body $\neg\Chk(v, \bquo{\z = \tofn 1})$, a certificate of a false
checker verdict at a free $v$ would have to certify, by cases, which of
the checker's four conjuncts fails, and would need a builder for
numeral disequality; the development builds neither.

\begin{thm}[\textsf{adr\_\allowbreak U\_\allowbreak ConRGA}, \textsf{adr\_\allowbreak ConRGA\_\allowbreak U}]
\label{thm:ad-interderivable}
$\cA\cdot U \eqv \ConR$.
\end{thm}

\noindent
Both entailments are \rga derivations under one hypothesis.  From
$\cA\cdot U$: at a free $v$, decide $\Chk(v, \bquo{\z = \tofn 1})$; in
the accepting case the test is nonzero, which contradicts the
hypothesis's instance.  From $\ConR$: at a free $v$, decide whether
the test is nonzero; in that case $v$ is accepted, which contradicts
$\ConR$'s instance, and in the other case double negation gives the
instance of $U$.  The interderivability is a pair of entailments, not
an \rga proof of the implication $\cA\cdot U \limp \ConR$.

\subsection{\texorpdfstring{$\ConR$ in \rgaP}{Con in RGA+}, and the separation}
\label{sec:adopt-con}

\begin{thm}[\textsf{adr\_\allowbreak ConRGA}, \textsf{adr\_\allowbreak separation}]
\label{thm:ad-sep}
\leavevmode
\begin{enumerate}
\item $\emptyset \vdashRP \cA\cdot U$ and $\emptyset \vdashRP \ConR$.
\item $\emptyset \nvdash \ConR$; $\emptyset \vdashRP \ConR$; every
  derivation $\Gamma \vdash c$ gives $\Gamma \vdashRP c$; and
  $\emptyset \nvdashRP \z = \tofn 1$.
\end{enumerate}
\end{thm}

\noindent
\emph{Proof.}  Part~1: adoption (Theorem~\ref{thm:ad-rule}(1)) applied
to the embedded warrant of Theorem~\ref{thm:ad-warrant}(2) gives
$\cA\cdot U$, and a cut with
Theorem~\ref{thm:ad-interderivable} gives $\ConR$.  Part~2 collects
Theorem~\ref{thm:g-g2}, part~1, Theorem~\ref{thm:ad-rule}(2) and
Theorem~\ref{thm:ad-sound}(2).  \qed

\medskip
\noindent
So \rgaP is a consistent proper extension of \rga, and adoption is not
vacuous.  \rgaP derives the very sentence that G\"odel's second theorem
denies to \rga.  That sentence is true in $\Mw$
([\citen{conrga-paper}, theorem~\textsf{goedelW\_omega\_incomplete}]).  \rga does not derive the adopted
universal $\cA\cdot U$ either, by Theorem~\ref{thm:ad-interderivable}
and Theorem~\ref{thm:g-g2}.

The separation has no HOL hypothesis.  The separation rests on two HOL
theorems: G\"odel's second theorem for \rga
(Theorem~\ref{thm:g-g2}), whose proof uses \rga's consistency, and
the consistency of \rgaP (Theorem~\ref{thm:ad-sound}), whose proof uses
\rga's $\Mw$-soundness and the genuineness of \rga's checker.
The Isabelle theory behind Theorem~\ref{thm:ad-sep} checks at build
time that the seven headline results use no oracle, that the
derivation of $\ConR$ in \rgaP reaches the adoption rule, the
code-level ex falso and the warrant, that this derivation reaches no
soundness lemma, neither \rga's nor \rgaP's nor the checker's
genuineness, and, as controls, that the proof of
Theorem~\ref{thm:g-g2} reaches the checker's cons law but neither the
adoption rule nor the code-level ex falso, and that the consistency of
\rgaP does reach \rgaP's soundness
(Appendix~\ref{app:adopt-checks}).

The development states G\"odel's second theorem for \rgaP neither
way.  The theorem would need a checker for \rgaP's derivations, and
the development defines none.  Classically one expects that \rgaP
would not derive a consistency sentence built from a faithful checker
for its own derivations; that expectation is a prediction, not a
result of this paper.

\subsection{Adoption at \iga}
\label{sec:adopt-iga}

\iga (Section~\ref{sec:checker}) has \rga's 42 rules and ground
induction, and its checker $\ChkI$ is \rga's table of rule detectors
with one more tag.  The rule
\adoptI is adoption with $\ChkI$ in the warrant,
\[
  \WarI(F, u) \;=\; \forall v.\ \ChkI\bigl(\cmp(F)\cdot v\cdot v,\ \Inst_u(v)\bigr)
\]
added to \iga's rules.  \igaP is \iga with \adoptI, and $\vdashIP$
is its entailment; $\vdashI$ is \iga's.  The instance code $\Inst_u$ is the one of \rgaP, since the
instance code does not mention the checker.  \iga's consistency
sentence for the checker $\ChkI$ is
\[
  \ConI \;=\; \forall v.\ \neg\ChkI(v, \bquo{\z = \tofn 1})
\]
(Definition~\ref{def:ig-con}), $\ConR$'s shape with \iga's checker.  The instances
of $\ConI$ are decided, since $\ChkI$'s verdict is.  Read classically,
through the faithfulness of $\ChkI$ (Theorem~\ref{thm:ck-iga-faith}), $\ConI$ says
that \iga is consistent.  That reading is informal: unlike $\ConR$,
whose truth in $\Mw$ is a theorem, $\ConI$ has no named $\Mw$-truth
theorem in the development, though one follows from
Theorem~\ref{thm:ac-turing}.

\begin{thm}[\textsf{adi\_\allowbreak embed}, \textsf{adi\_\allowbreak omega\_\allowbreak sound}]
\label{thm:ad-iga-rule}
\leavevmode
\begin{enumerate}
\item From $\Gamma \vdashIP \WarI(F, u)$, with $u$ closed and
  well-formed, $\Gamma \vdashIP \cA\cdot u$ (\textsf{adi\_adopt}).
\item If $\Gamma \vdashI c$, then $\Gamma \vdashIP c$; if $\Gamma
  \vdash c$, then $\Gamma \vdashIP c$ (\textsf{adi\_embed\_rga}).
\item \igaP is $\Mw$-sound, and $\emptyset \nvdashIP \z = \tofn 1$
  (\textsf{adi\_consistent}).
\end{enumerate}
\end{thm}

\noindent
The soundness proof is the one of Theorem~\ref{thm:ad-sound}, with
the genuineness of $\ChkI$ (Theorem~\ref{thm:ck-iga-faith}(2)) in
place of \rga's checker and \iga's $\Mw$-soundness~\cite{iga-paper}
in place of \rga's.  \iga's per-step soundness is reused for \iga's 43 rules.

The body $U_{\mathrm I}$ and the witness $F_{U_{\mathrm I}}$ are those
of Definition~\ref{def:ad-body} with \iga's compiled test in place of
\rga's.

\begin{thm}[\textsf{adi\_\allowbreak warrant\_\allowbreak uC}, \textsf{adi\_\allowbreak ConRGA\_\allowbreak of\_\allowbreak ConIGA}]
\label{thm:ad-iga-warrant}
\leavevmode
\begin{enumerate}
\item $\emptyset \vdash \WarI(F_{U_{\mathrm I}}, U_{\mathrm I})$, an
  \rga derivation.
\item $\cA\cdot U_{\mathrm I} \eqv \ConI$, in \rga
  (\textsf{adi\_U\_ConIGA}, \textsf{adi\_ConIGA\_U}).
\item $\ConI \vdash \ConR$, in \rga.
\end{enumerate}
\end{thm}

\noindent
For part~1, the two branches are those of
Theorem~\ref{thm:ad-warrant}.  The value-certificate branch is an \rga
proof accepted by \rga's checker, embedded into $\ChkI$ by the law
that an \rga proof accepted by $\Chk$ is accepted by $\ChkI$
(Theorem~\ref{thm:ck-embed-cert}(2)).  The ex-falso branch is built from the
generic layer of Section~\ref{sec:checker} and uses only \rga's rule tags
(Appendix~\ref{app:adopt-iga}).  \iga's ground induction occurs
nowhere in the proof: all three parts are \rga derivations.  Part~3
uses the same embedding law at a free $v$: decide $\Chk(v, \bquo{\z =
\tofn 1})$, and in the accepting case embed the verdict into
$\ChkI(v, \bquo{\z = \tofn 1})$, which contradicts $\ConI$'s instance.

\begin{thm}[\textsf{adi\_\allowbreak ConIGA}, \textsf{adi\_\allowbreak ConRGA}]
\label{thm:ad-iga-con}
$\emptyset \vdashIP \ConI$ and $\emptyset \vdashIP \ConR$.
\end{thm}

\noindent
\emph{Proof.}  Adoption at the embedded warrant of
Theorem~\ref{thm:ad-iga-warrant}(1) gives $\cA\cdot U_{\mathrm I}$.  A cut with part~2 of that theorem gives $\ConI$,
and a further cut with part~3 gives $\ConR$.  \qed

\begin{thm}[\textsf{adi\_\allowbreak separation\_\allowbreak RGA}, \textsf{adi\_\allowbreak separation\_\allowbreak of\_\allowbreak goedel2}]
\label{thm:ad-iga-sep}
\leavevmode
\begin{enumerate}
\item $\emptyset \nvdash \ConR$; $\emptyset \vdashIP \ConR$; every
  derivation of \rga and every derivation of \iga is one of \igaP; and
  $\emptyset \nvdashIP \z = \tofn 1$.
\item If $\emptyset \nvdashI \ConI$, then: $\emptyset \nvdashI \ConI$;
  $\emptyset \vdashIP \ConI$; every derivation of \iga is one of
  \igaP; and $\emptyset \nvdashIP \z = \tofn 1$.
\end{enumerate}
\end{thm}

\noindent
Part~1 separates \igaP from \rga.  \rga is a proper part of \igaP, and
\iga is contained in \igaP.  Part~1 has no HOL hypothesis and rests on
Theorem~\ref{thm:g-g2} and on Theorem~\ref{thm:ad-iga-rule}(3).
Part~1 does not separate \igaP from \iga, since whether \iga derives
$\ConR$ is open (Section~\ref{sec:open}).  Part~2's hypothesis is
G\"odel's second theorem at \iga, which Section~\ref{sec:iga-goedel}
proves (Theorem~\ref{thm:ig-g2}).  With that theorem, \iga is a proper
part of \igaP, again with no HOL hypothesis
(Theorem~\ref{thm:ig-sep}).  That separation rests on G\"odel's second theorem at \iga, whose proof uses
\iga's consistency~\cite{iga-paper}, and on
the consistency of \igaP (Theorem~\ref{thm:ad-iga-rule}(3)).  The
statement is about one consistency sentence: \iga does not derive
$\ConI$, its consistency sentence for the checker $\ChkI$.

The Isabelle theory behind Theorem~\ref{thm:ad-iga-con} checks at
build time that seventeen results use no oracle, that the derivation
of $\ConI$ in \igaP reaches the adoption rule, the code-level ex falso
and the warrant and reaches no soundness lemma of \rga, of \iga or of
\igaP, nor the genuineness of $\ChkI$, that the consistency of \igaP reaches the
genuineness of $\ChkI$, \iga's per-step soundness and the soundness of
\adoptI, and that the code-level ex falso at \iga avoids \rga's cons
and append laws for $\Chk$ and \rga's code-level ex falso
(Appendix~\ref{app:adopt-checks}).

The development states none of three relations between the adoption
systems: that \rgaP is contained in \igaP, whether \rgaP derives
$\ConI$, and whether \iga and \rgaP are comparable.  Nor does it state
G\"odel's second theorem for \igaP, since it defines no checker for
\igaP's derivations.  Classically one expects the same as for \rgaP;
that expectation is a prediction.

\subsection{What adoption is not}
\label{sec:adopt-not}

Adoption is not an internal soundness argument.  The derivations of
the warrants split on the decided value of the checker's test at the
proof code, and never on the truth of the judgments inside a proof.
The internal-reasoning paper's no-go theorem
([\citen{conrga-paper}, theorem~\textsf{xder\_not\_fold\_residual}]) shows that no $\Mw$-sound
entailment over \rga's terms derives the guard, the motive or any
accepted residual of that paper's internal-soundness fold at one
proof, because the fold's motive is not true in $\Mw$ at a hypothesis
node whose hypothesis is ungrounded.  So the derivations of $\ConR$ above stand in no tension
with that theorem.  The no-go theorem also applies to \igaP, because
\igaP is $\Mw$-sound (Theorem~\ref{thm:ad-iga-rule}(3)), so \igaP does
not derive the fold's motive either.  The development interprets the
no-go locale at \igaP through that soundness theorem, as well as at
\iga and at \oiga (Section~\ref{sec:open-nogo} on what the no-go
theorem excludes, and Appendix~\ref{app:adopt-nogo}).

This paper does not offer adoption as a tool for proving particular
theorems.  The role of adoption in this paper is calibration.
Adoption is a sound rule with which the consistency sentences become
derivable, and it separates strictly at both levels.
The IGA paper, in its revised version~\cite{iga-paper}, computes the
accessibility ordinals of calibration systems that add reflection to
\oiga, with the adoption premise checked in HOL: $\varphi(\omega,1)$
for one reflection step and $\varphi(\omega+1,0)$ for the autonomous
progression.

\subsection{Adopting a notation claim}
\label{sec:adopt-claims}

The adoption rule of the preceding subsections concludes a universal
from a warrant that the system derives.  The IGA paper prices a second
kind of extension of \iga, by a hypothesis that is supplied rather than
derived~\cite{iga-paper}.  The development defines a form of Kleene's
$\KO$, the constructive ordinal notations~\cite{kleene38notation}, as
one positive fixed point of \iga on tagged codes: the code of $0$, the
successor $\langle 1, a\rangle$ of a code $a$, and the limit $\langle
2, e\rangle$ of a primitive-recursive sequence whose step function has
the code $e$.  The \emph{notation claim} at a code $a$ is the sentence
$\KO(\tofn a)$, the fixed point applied to the numeral of $a$.  Ground
induction at this fixed point is induction along $\KO$, so under the
claim \iga has transfinite induction along the structure of $a$, at
every closed motive (\textsf{adopt\_gives\_TI}), and the IGA paper's
one-claim theorems buy the well-foundedness of whole notation systems
in this way.

\begin{defi}[\textsf{adoptDer}, \textsf{adopt\_\allowbreak ext\_\allowbreak iff}]
\label{def:ac-ext}
The \emph{claim extension} $\Ia{a}$ is \iga with the notation claim
$\KO(\tofn a)$ as a hypothesis in every context: $\Gamma \vdashA{a} c$
when $\{\KO(\tofn a)\} \cup \Gamma \vdashI c$.
\end{defi}

\noindent
The development proves this presentation equal to \iga with one more
rule, an axiom rule that concludes the claim in every context and has
no premise.  Write $a \in \KO$, a HOL statement, when $a$ lies in the
least fixed point of $\KO$'s defining operator.  Classically, $\Ia{a}$
is $\Mw$-sound exactly when the claim $\KO(\tofn a)$ is true in $\Mw$,
which holds exactly when $a \in \KO$, and $\Ia{a}$ is then consistent
(\textsf{adopt\_sound\_iff\_lfp}, \textsf{adopt\_consistent}).  A claim
extension at a code $a \in \KO$ is a \emph{sound} claim extension.  The
claim $\KO(\tofn a)$ is a sentence of \iga; the fact $a \in \KO$ is not,
and the IGA paper adopts claims exactly where \iga does not derive
them, with the classical fact $a \in \KO$ as the warrant.  An adopted
claim is thus a hypothesis with a classical price, and the question is
what tests it.

\begin{thm}[\textsf{fu\_\allowbreak falsified}, \textsf{fu\_\allowbreak falsified\_\allowbreak by\_\allowbreak test}, \textsf{fu\_\allowbreak falsified\_\allowbreak by\_\allowbreak refuted}, \textsf{fu\_\allowbreak falsified\_\allowbreak by\_\allowbreak instance}]
\label{thm:ac-falsified}
Let $a$ be a code.  Each of the following implies $a \notin \KO$:
\begin{enumerate}
\item $\emptyset \vdashA{a} \z = \tofn 1$;
\item $\emptyset \vdashA{a} \varphi$, for a sentence $\varphi$ that is
  false in $\Mw$;
\item $\emptyset \vdashA{a} \varphi$ and $\emptyset \vdashI
  \neg\varphi$;
\item $\emptyset \vdashA{a} \cA\cdot f$ and $\emptyset \vdashI
  \neg(f\cdot\tofn n)$, for some number $n$.
\end{enumerate}
In cases~3 and~4, $\Ia{a}$ derives $\z = \tofn 1$ in every context in
which it derives $\varphi$ or $\cA\cdot f$ (\textsf{fu\_\allowbreak
refuted\_\allowbreak explodes}, \textsf{fu\_\allowbreak
instance\_\allowbreak explodes}).
\end{thm}

\noindent
The four implications are classical, HOL theorems about $\Ia{a}$ and
$\Mw$.  Parts~1 and~2 are the contrapositives of the consistency and
the $\Mw$-soundness of a sound claim extension.  The last sentence is
grounded, a derivation in $\Ia{a}$: \iga's refutation, weakened into
the context, and the derived sentence (in case~4, the derived
universal's instance at $\tofn n$) give $\z = \tofn 1$ by explosion.  A
claim is tested only through its consequences.  A consequence that is
false in $\Mw$ falsifies the claim.  When \iga refutes the consequence,
the falsification is visible inside the system, as a derivation of $\z
= \tofn 1$ in the extension.

\iga refutes the claim itself only at a \emph{bad} code
(\textsf{gi\_refutes\_only\_bad}): a code that is neither the code of
$0$ nor tagged as a successor or a limit, the successor of a bad code,
or a limit code whose sequence has no value at some argument or a bad
value at some argument.  Call a limit code \emph{well-shaped} when its
sequence has a value at every argument and no value is bad
(\textsf{oHonest}).  So \iga never refutes the claim at a well-shaped
limit code.  Can the adoption of a well-shaped code still be
inconsistent, through induction along $\KO$ rather than through a
refutation?  It can.

\begin{thm}[\textsf{fuQ\_\allowbreak falsified\_\allowbreak by\_\allowbreak consequence}, \textsf{adopt\_\allowbreak fuQ\_\allowbreak explodes}, \textsf{fuQ\_\allowbreak omv\_\allowbreak gap}]
\label{thm:ac-quine}
There is a limit code $\mathit{qn} = \langle 2, e\rangle$, the
\emph{quine code}, whose sequence has the value $\mathit{qn}$ at every
argument, such that:
\begin{enumerate}
\item $\mathit{qn}$ is well-shaped, and $\emptyset \nvdashI
  \neg\KO(\tofn{\mathit{qn}})$;
\item the claim $\KO(\tofn{\mathit{qn}})$ is neither true nor false in
  $\Mw$, and $\mathit{qn} \notin \KO$;
\item $\emptyset \vdashA{\mathit{qn}} \z = \tofn 1$.
\end{enumerate}
\end{thm}

\noindent
Parts~1 and~2 are classical.  Part~3 is grounded, a derivation in \iga
from the claim.  The code is a quine at the level of index trees.  An
explicit tree $B$ maps a number $z$ to the limit code whose step
function is the tree with the code $z$, applied to the numerals of $z$
and $0$.  The step tree of $\mathit{qn}$ is $B$ applied to the numeral
of $B$'s own code, so that tree computes, at every argument, the limit
code built from that same tree, which is $\mathit{qn}$
(\textsf{pr2fi\_fuQF}).  The derivation of part~3 is one induction along
$\KO$ at the claimed code, with the decided motive $x \neq
\tofn{\mathit{qn}}$.  A zero or successor code differs from
$\mathit{qn}$ by its tag.  In the limit case, a limit code equal to
$\mathit{qn}$ has the step code $e$.  The induction hypothesis at the
sequence's element $0$ supplies a value $r$ of the sequence with $r
\neq \tofn{\mathit{qn}}$, while the functionality of the sequence's
graph, an \rga derivation at an open argument
(\textsf{fuQ\_gr\_val}), gives $r = \tofn{\mathit{qn}}$.  The motive at
$\mathit{qn}$ itself is $\tofn{\mathit{qn}} \neq \tofn{\mathit{qn}}$,
which yields $\z = \tofn 1$.

An observation thus refutes the theory, though the theory is never
refuted.  $\Ia{\mathit{qn}}$ derives $\z = \tofn 1$, a grounded fact,
while \iga never derives the claim's negation, a classical fact.  The
step from the derived contradiction to the falsity of the claim is
modus tollens in the metatheory.  \iga has no rule that turns a
contradictory extension into a refutation of the hypothesis, and
Theorem~\ref{thm:ac-rosser}(1) below shows the same gap at a
well-foundedness sentence.  In $\Mw$ the claim has no value, as the
ungrounded term $W\cdot W$ of \S\ref{sec:adopt-sound} has none.  So
the fact that \iga does not refute a claim does not make the claim
consistent to adopt: induction along $\KO$ detects a computable
descent, here the constant sequence of $\mathit{qn}$, that no
refutation sees.  For a code outside $\KO$ with no hyperarithmetic
descending path, in the manner of Harrison's recursive
pseudo-well-orderings~\cite{harrison68recursive}, the claim extension
is not $\Mw$-sound (\textsf{adopt\_sound\_iff\_lfp}, as at every
code outside $\KO$) and is predicted to be consistent; the
prediction is not a theorem of the development.

\subsection{The second theorem at every claim}
\label{sec:adopt-claims-g2}

G\"odel's second theorem says when a stronger hypothesis is needed.  In
the currency of notation claims, the statement can be made exact, at
\iga and at every sound claim extension.

\begin{defi}[\textsf{pwcT}, \textsf{tc}]
\label{def:ac-turing}
For an explicit index tree $T$ and a number $k$, let $\tau_{T,k}$ be
the limit code whose sequence has, at $n$, the successor code of the
finite notation $n$, the notation for $n+1$, when $\cmp(T)\cdot\tofn
k\cdot\tofn n$ has the value $0$, and the ill-shaped code $\langle 3,
0\rangle$ otherwise.  Let $u_{T,k} = \lambda v.\,(\cmp(T)\cdot\tofn
k\cdot v = \z)$.  \emph{Turing's code} $\tau$ is the code $\tau_{T,k}$
for which $u_{T,k}$ is the body $U_{\mathrm I}$ of
\S\ref{sec:adopt-iga}: the sequence of $\tau$ climbs through the
finite notations as long as $n$ is not an accepted \iga proof of
$\emptyset \vdash \z = \tofn 1$.
\end{defi}

\begin{thm}[\textsf{tc\_\allowbreak ConIGA\_\allowbreak iff}, \textsf{tc\_\allowbreak in\_\allowbreak lfp}, \textsf{fu\_\allowbreak tc\_\allowbreak strict}]
\label{thm:ac-turing}
\leavevmode
\begin{enumerate}
\item $\KO(\tofn\tau) \eqv \ConI$.
\item $\tau \in \KO$.
\item $\emptyset \nvdashI \ConI$; $\emptyset \vdashA{\tau} \ConI$;
  $\emptyset \nvdashA{\tau} \z = \tofn 1$; and $\emptyset \nvdashI
  \KO(\tofn\tau)$.
\end{enumerate}
\end{thm}

\noindent
Part~1 is a pair of \rga derivations, each under one hypothesis, so
both hold in \iga and in every claim extension.  The derivations pass
through the universal $\cA\cdot U_{\mathrm I}$: for every $T$ and $k$,
the claim $\KO(\tofn{\tau_{T,k}})$ and $\cA\cdot u_{T,k}$ are
interderivable in \rga (\textsf{pwcT\_rule\_I}, \textsf{pwcT\_rule\_E}),
and $\cA\cdot U_{\mathrm I} \eqv \ConI$ is
Theorem~\ref{thm:ad-iga-warrant}(2).  Part~2 is classical: \igaP
derives $\cA\cdot U_{\mathrm I}$ (\S\ref{sec:adopt-iga}) and is
$\Mw$-sound (Theorem~\ref{thm:ad-iga-rule}(3)), so the test never
fires.  Part~3 is classical.  Its first conjunct is G\"odel's second
theorem at \iga (Theorem~\ref{thm:ig-g2}), and its last conjunct
follows from the first by part~1.

$\Ia{\tau}$ is the consistency extension of \iga, written as a notation
claim.  The sentence ``G\"odel's second theorem says when a strictly
larger claim is needed'' is exact at $\tau$: \iga does not derive the
claim, and the sound extension $\Ia{\tau}$ derives $\ConI$, because the
claim is $\ConI$ in disguise.  The theorem carries the caveat that
Turing made himself, of ``small ordinals expressed in a very complex
and artificial way''~\cite[p.~202]{turing39ordinals}, and that
Feferman's analysis of transfinite progressions
sharpened~\cite{feferman62transfinite}: the code encodes the answer.  The
sequence of $\tau$ reads \iga's proof checker.  Read classically, as an
observation that the development does not mechanize, $\tau$ is a
notation for $\omega$, the ordinal of the standard limit code for
$\omega$, the limit of the finite notations, whose claim \iga already
derives, so that adopting that claim is free (\textsf{o\_omega},
\textsf{adopt\_omega\_free}).  The ordinal is the same;
the code differs, and the code carries $\ConI$.  Classically, too, the
theories of a progression along $\KO$ depend on the notations and not
only on the ordinals they denote, and Turing's completeness theorem
for his ordinal logics is obtained at notations built from the
sentences to be proved~\cite[pp.~200--202]{turing39ordinals},
\cite{feferman62transfinite}.

The same construction at every primitive-recursive test gives the
completeness of the sound claims for the universals of such tests.

\begin{thm}[\textsf{fu\_\allowbreak pi01\_\allowbreak iff}, \textsf{fu\_\allowbreak pi01\_\allowbreak derives}, \textsf{fuU\_\allowbreak pi01\_\allowbreak iff}]
\label{thm:ac-pi01}
Let $T$ be an explicit index tree and $k$ a number, and say that the
test of $T$ at $k$ never fires when $\cmp(T)\cdot\tofn k\cdot\tofn n$
has the value $0$ at every $n$.
\begin{enumerate}
\item $\tau_{T,k} \in \KO$ if and only if the test of $T$ at $k$ never
  fires.
\item $\emptyset \vdashA{\tau_{T,k}} \cA\cdot u_{T,k}$, whether or not
  $\tau_{T,k} \in \KO$.
\item Let $\vdash_{\mathrm U}$ be the union of the sound claim
  extensions over finite sets of codes: $\Gamma \vdash_{\mathrm U} c$
  when \iga derives $c$ from $\Gamma$ together with the claims at
  finitely many codes in $\KO$.  Then $\emptyset \vdash_{\mathrm U}
  \cA\cdot u_{T,k}$ if and only if the test of $T$ at $k$ never fires.
\end{enumerate}
\end{thm}

\noindent
Parts~1 and~3 are classical; part~2 is grounded, the \rga derivation
of the universal from the claim.  The union $\vdash_{\mathrm U}$ is
$\Mw$-sound and consistent and closed under cut
(\textsf{fuU\_omega\_sound}, \textsf{fuU\_consistent},
\textsf{fuU\_cut}).  So the sound claims prove every true universal of
a primitive-recursive test and no false one.  The union is defined
from the HOL set $\KO$, and the theorem gives no way of finding the
right code; the code $\tau_{T,k}$ is built from the test itself.

The second theorem holds at every sound claim extension, with a checker
built from \iga's.

\begin{defi}[\textsf{axIdx}, \textsf{chkF}, \textsf{ConF}]
\label{def:ac-chk}
The \emph{axiom cell} $E_a$ is a rule detector that fires at a node
exactly when the node's conclusion is the code of $\KO(\tofn a)$, at
any context and whatever entries precede the node
(\textsf{ax\_fires\_iff}).  Let $L_a$ be \iga's list $L_{\iga}$ with
$E_a$ appended (Definition~\ref{def:ck-iga}).  The checker of $\Ia{a}$
is $\ChkA{a} = \Chk_{\mathcal D_{L_a}}$, its provability formula is
$\PrA{a}(c) = \neg\forall q.\,\neg\ChkA{a}(q, c)$, and its consistency
sentence is
\[
  \ConA{a} \;\ldef\; \forall p.\ \neg\ChkA{a}(p, \bquo{\z = \tofn 1}) .
\]
\end{defi}

\noindent
The axiom cell is not a citation cell, since the axiom rule has no
premise and no side condition.  The cell reads the node's conclusion
only, in the pattern of \rga's detector for the premise-free rule that
concludes $\z\jnat$.  So the cell's shift law (\textsf{tblShift\_ax})
uses no position map, bound or equality of nodes, and four premises
suffice where a citation cell needs seven.

\begin{thm}[\textsf{drvChkF\_\allowbreak fires\_\allowbreak iff}, \textsf{fu\_\allowbreak ent\_\allowbreak eq}, \textsf{fu\_\allowbreak g2\_\allowbreak ent}, \textsf{fu\_\allowbreak g2\_\allowbreak con}, \textsf{fu\_\allowbreak goedel2\_\allowbreak internal}, \textsf{fu\_\allowbreak goedel2}]
\label{thm:ac-g2}
Let $a$ be a code.
\begin{enumerate}
\item The table $\mathcal D_{L_a}$ fires at a node exactly when the
  node is a step of \iga or concludes $\KO(\tofn a)$, and the
  entailment of these steps is $\vdashA{a}$.
\item $\vdashA{a}$ with the checker $\ChkA{a}$ meets the structural
  conditions, NEC, D3, EXPL and D1, and the conditions of the optional
  layers, ex falso, D2 and necessitation with two premises.  If $a \in
  \KO$, it also meets consistency, so $\Ia{a}$ is a G\"odel extension
  of \rga.
\item $\ConA{a} \vdashA{a} \neg\PrA{a}(\bquo{\ConA{a}})$.
\item If $a \in \KO$, then $\emptyset \nvdashA{a} \ConA{a}$.
\end{enumerate}
\end{thm}

\noindent
Parts~1, 2 and~4 are classical.  Part~3 is grounded, a derivation in
$\Ia{a}$ under one hypothesis, at every code, in $\KO$ or not; the
build checks that its dependency cone contains no soundness,
consistency or $\Mw$ fact.  Part~4 is Theorem~\ref{thm:lc-g2}(3) at
$\Ia{a}$.  The witnesses of part~2 are those of \iga
(Theorem~\ref{thm:ig-witness}) run again at the claim extension's
entailment and checker, one interpretation lemma per layer
(\textsf{fu\_g2\_syn}, \textsf{fu\_g2\_ent}, \textsf{fu\_g2\_acc},
\textsf{fu\_g2\_exf}, \textsf{fu\_g2\_d2}, \textsf{fu\_g2\_nec2},
\textsf{fu\_g2\_con}), with the locale's consistency sentence
rewritten to $\ConA{a}$ (\textsf{fu\_g2\_rw\_con}).  The table-generic transformer laws of
Section~\ref{sec:checker} hold at $\mathcal D_{L_a}$ unchanged; their
premises are part~1, the shift law of the table
(\textsf{tblShift\_F}) and the presence of \rga's rule tags in $L_a$.
Each law is an \rga derivation, carried into $\Ia{a}$ by cut transfer.
The consistency of $\Ia{a}$ at $a \in \KO$ comes from \iga's
$\Mw$-soundness.  As a control, at the code $a_0$ of the finite
notation $0$, $\Ia{a_0}$ is \iga itself (\textsf{adopt\_fin\_free}),
and part~4 there is reached a second way, from G\"odel's second
theorem at \iga, because \rga derives $\ConI$ from $\ConA{a}$ at every
$a$ (\textsf{fu\_ConF\_ConIGA}, \textsf{gi\_not\_ConF},
\textsf{fu\_goedel2\_fin}); $\ConA{a_0}$ is a different sentence from
$\ConI$, since its table has the extra cell.

The claim extensions are extensions of \iga by one undecided axiom,
and Theorem~\ref{thm:ac-g2} gives each a checker that meets the
G\"odel conditions, so the scope of the conditions
(\S\ref{sec:locale-scope}) contains this family of axiom extensions.

\begin{thm}[\textsf{tcF\_\allowbreak ConF\_\allowbreak iff}, \textsf{tcF\_\allowbreak in\_\allowbreak lfp}, \textsf{fu\_\allowbreak needs\_\allowbreak larger\_\allowbreak fuel}, \textsf{fu\_\allowbreak not\_\allowbreak tcF}]
\label{thm:ac-larger}
For a code $a$, let $\tau_a$ be the code $\tau_{T,k}$ for which
$u_{T,k}$ is the body of Definition~\ref{def:ad-body} at the compiled
proof test of $\ChkA{a}$ (\textsf{tcF}).
\begin{enumerate}
\item $\KO(\tofn{\tau_a}) \eqv \ConA{a}$.
\item If $a \in \KO$, then $\tau_a \in \KO$; $\emptyset \nvdashA{a}
  \ConA{a}$; $\emptyset \vdashA{\tau_a} \ConA{a}$; and $\emptyset
  \nvdashA{a} \KO(\tofn{\tau_a})$.
\end{enumerate}
\end{thm}

\noindent
Part~1 is a pair of \rga derivations, each under one hypothesis, from
the interderivability of $\cA\cdot u_{T,k}$ with $\ConA{a}$ at any
table.  Part~2 is classical.  At $a \in \KO$, \rga refutes every
numeral instance of the matrix of $\ConA{a}$, since a closed accepting
verdict would be read back, through the faithfulness of part~1 of
Theorem~\ref{thm:ac-g2}, as a proof of $\z = \tofn 1$ in $\Ia{a}$
(\textsf{fu\_ConF\_instances\_refuted}).  So $\ConA{a}$ is true in
$\Mw$, \rga's $\Mw$-soundness carries that truth to $\cA\cdot u_{T,k}$,
and $\tau_a \in \KO$ by Theorem~\ref{thm:ac-pi01}(1).  The second
conjunct is Theorem~\ref{thm:ac-g2}(4), the third follows from part~1,
and the fourth from the second by part~1.  At $a_0$ the theorem
restates Theorem~\ref{thm:ac-turing} at a second code: $\Ia{\tau_{a_0}}$
derives $\ConA{a_0}$ and $\ConI$, is consistent, and \iga does not
derive its claim (\textsf{fu\_needs\_larger\_fuel\_fin}).

At every sound claim extension, then, G\"odel's second theorem rules
out the extension's own consistency sentence, and the sound claim
built from the extension's proof checker derives it.  The strictness
is non-containment: $\Ia{\tau_a}$ derives a sentence that $\Ia{a}$
does not, and $\Ia{a}$ does not derive the claim at $\tau_a$.  The
development does not state that $\Ia{\tau_a}$ contains $\Ia{a}$.  The
sequence $a, \tau_a, \tau_{\tau_a}, \ldots$ replaces the claim at each
step by the consistency of the previous claim extension: the claim at
$\tau_a$ and $\ConA{a}$ derive each other
(Theorem~\ref{thm:ac-larger}(1)), so $\Ia{\tau_a}$ is \iga together with
$\ConA{a}$, not $\Ia{a}$ together with $\ConA{a}$; the claim at $a$ is
dropped at each step, the sequence is not cumulative, and it is not a progression in
Turing's or Feferman's sense; a cumulative progression along an index
is left to further work.  The caveat of
Turing and Feferman applies at every step: $\tau_a$ encodes the checker of $\Ia{a}$, so the step
from $a$ to $\tau_a$ gains no information that was not put into the
code.

\paragraph{Autonomy over claims.}
Kreisel's autonomy condition admits an ordinal notation as an index of
a progression only once an earlier stage has proved it one; its
best-known application is the analysis of predicative provability,
whose limit Feferman and Sch\"utte found to be
$\Gamma_0$~\cite{kreisel58ordinal,feferman64systems,schutte65predikative}.  Over bare claims,
autonomy adds nothing.

\begin{thm}[\textsf{fuAdm\_\allowbreak iff}, \textsf{fuAut\_\allowbreak eq\_\allowbreak iga}]
\label{thm:ac-autonomy}
Let the admitted codes be the least set that contains $a_0$ and that
contains $b$ whenever $\emptyset \vdashA{a} \KO(\tofn b)$ for an
admitted $a$.  Then $b$ is admitted if and only if $\emptyset \vdashI
\KO(\tofn b)$, and $\Gamma \vdashA{a} c$ for some admitted $a$ if and
only if $\Gamma \vdashI c$.
\end{thm}

\noindent
The theorem is classical, and its proof uses cut, weakening, and
\iga's derivation of the claim at $a_0$ (\textsf{o\_fin\_num}).
If \iga derives the claim at $a$, every derivation of $\Ia{a}$ becomes
one of \iga by one cut.  So the autonomous closure of \iga under its
own notation claims is \iga itself.  The license that admits a claim is
a proof of the claim, and a claim that \iga proves is one that \iga
can already use, with induction along $\KO$ at every closed motive.  A growing
progression needs a step that a bare claim does not supply:
consistency, as the codes $\tau_a$ supply it, though $\Ia{a}$ does not
derive the claim at $\tau_a$ (Theorem~\ref{thm:ac-larger}); reflection,
as in the IGA paper's autonomous progression of reflection, whose
license is a derived accessibility sentence~\cite{iga-paper}; or truth.
Theorem~\ref{thm:ac-pi01} marks the other extreme.  The union of all
sound claims, chosen with the classical knowledge of which codes lie
in $\KO$, proves every true universal of a primitive-recursive test.
Autonomous choice of claims adds nothing, classical choice adds every
true test universal, and the content lies in the choice.

\subsection{No adequacy rule for adopted claims and well-foundedness sentences}
\label{sec:adopt-ac}

Since claims are supplied, one may look for a rule that adopts a claim,
or a well-foundedness sentence, whenever it passes a test that the
system can apply.  Two tests suggest themselves: that the system does
not refute the sentence, and that the system together with the
sentence does not derive $\z = \tofn 1$.  The first test gives no
closure condition in a grounded system, for claims or for
well-foundedness sentences.  The second gives none for
well-foundedness sentences; for claims, the development states only
the first test, and the prediction is that the second fails as well,
through a Liar code built like $\tau_{T,k}$.

For a closed, well-formed term $R$, read as an order on the numbers,
write $\WF(R)$ for the IGA paper's sentence that $R$ is well-founded,
$\forall x.\,(x\jnat \limp \Acc_R\cdot x)$, where $\Acc_R$ is the
accessible part of $R$, a positive fixed point~\cite{iga-paper}.  Let
$D$ be the \emph{successor-reversed} order, in which $n+1$ lies below
$n$.

\begin{thm}[\textsf{gi\_\allowbreak never\_\allowbreak refutes\_\allowbreak WF}, \textsf{wf\_\allowbreak stepD\_\allowbreak explodes}, \textsf{dirAC\_\allowbreak inconsistent}, \textsf{dirAC\_\allowbreak code\_\allowbreak inconsistent}]
\label{thm:ac-direct}
\leavevmode
\begin{enumerate}
\item For every closed, well-formed $R$, $\emptyset \nvdashI
  \neg\WF(R)$.
\item $\{\WF(D)\} \vdashI \z = \tofn 1$.
\item There are finitely many sentences $\WF(R)$, none refuted by
  \iga, from which \iga derives $\z = \tofn 1$, and finitely many
  notation claims, none refuted by \iga, from which \iga derives $\z
  = \tofn 1$.
\end{enumerate}
\end{thm}

\noindent
Part~1 is classical.  Accessibility is never false in $\Mw$, so
neither is $\WF(R)$, and \iga's $\Mw$-soundness excludes a refutation.
Part~2 is grounded, a derivation in \iga.  The identity sequence
descends along $D$ from $0$, and $\WF(D)$ makes $0$ accessible, so the
IGA paper's descending-sequence theorem, by which a primitive-recursive
sequence that \iga proves descending from a point it proves accessible
gives $\z = \tofn 1$~\cite{iga-paper}, applies.  Part~3 is classical,
with one sentence as the witness in each case: $\WF(D)$, by parts~1
and~2, and the claim at the quine code $\mathit{qn}$, by
Theorem~\ref{thm:ac-quine}.
In rule form, every set of sentences that is closed under \iga
consequence and adopts each unrefuted $\WF(R)$, or each unrefuted
claim, contains $\z = \tofn 1$ (\textsf{dirAC\_rule\_inconsistent},
\textsf{dirAC\_code\_rule\_inconsistent}).  In a grounded system,
unrefutedness is no test of well-foundedness, and no self-reference is
needed to see it.

The second test needs self-reference.  An entailment relation $S$ on
\rga's terms is \emph{closed under consistency-filtered adoption} when,
for every closed, well-formed $R$, if $\{\WF(R)\} \nvdash_S \z = \tofn
1$ then $\emptyset \vdash_S \WF(R)$ (\textsf{fuACclosed}); the
\emph{decided} form asks the same only for the orders $R$ whose
instances $R\cdot y\cdot x$ \iga derives decided at open $y$ and $x$
(\textsf{fuACdec}).

\begin{thm}[\textsf{fu\_\allowbreak ACdec\_\allowbreak inconsistent}, \textsf{fu\_\allowbreak AC\_\allowbreak inconsistent}, \textsf{gi\_\allowbreak AC\_\allowbreak inconsistent}, \textsf{fu\_\allowbreak AC\_\allowbreak inconsistent\_\allowbreak fuel}]
\label{thm:ac-liar}
Let $S$ be an $\Mw$-sound entailment relation that contains \iga's,
with a checker $\Chk_{\mathcal D}$ at a table $\mathcal D$ that meets
the structural conditions and D1, and that is faithful at closed
codes: whenever \rga derives $\Chk_{\mathcal D}(\tofn p,
\bquo{\varphi})$, $\emptyset \vdash_S \varphi$.  Then $S$ is not closed
under consistency-filtered adoption, not even in the decided form.  In
particular, neither \iga nor any sound claim extension $\Ia{a}$ is
closed under the decided form.
\end{thm}

\noindent
The theorem is classical, and the derivations that carry it are
grounded.  The witness is a Liar order $R_{\mathrm L}$, in which $x+1$
lies below $x$ exactly when some $p \leq x$ is accepted by
$\Chk_{\mathcal D}$ as a proof of $\WF(R_{\mathrm L})$; the
self-reference is by the diagonal function $\Diag$ of
\S\ref{sec:goedel-sentences} (in the locales' copy,
\textsf{g2DiagSub}), and \iga derives the order decided at
open arguments (\textsf{fuRL\_dec}).  If $S$ derives $\WF(R_{\mathrm
L})$, D1 gives an accepted code $p$.  Every $x \geq p$ then has $x+1$
below it, and \iga derives $\z = \tofn 1$ from $\WF(R_{\mathrm L})$ and
the acceptance of $p$, along the descent $j \mapsto j + p$
(\textsf{fuRL\_WF\_explodes\_seq}), so $S$ derives $\z = \tofn 1$,
which $\Mw$-soundness excludes.  If $S$ does not derive
$\WF(R_{\mathrm L})$, faithfulness leaves no code accepted, the order
is empty in $\Mw$, $\WF(R_{\mathrm L})$ is true in $\Mw$, and $S$
together with $\WF(R_{\mathrm L})$ is consistent by $\Mw$-soundness
(\textsf{fuRL\_ext\_consistent}); closure would then make $S$ derive
$\WF(R_{\mathrm L})$.  Either way, closure fails.

Rosser's device removes the soundness hypothesis, but the textbook
Rosser clause does not work in a grounded system.

\begin{thm}[\textsf{fu\_\allowbreak WF\_\allowbreak no\_\allowbreak neg\_\allowbreak intro}, \textsf{fu\_\allowbreak ACR\_\allowbreak inconsistent}, \textsf{gi\_\allowbreak ACR\_\allowbreak inconsistent}, \textsf{fu\_\allowbreak ACR\_\allowbreak inconsistent\_\allowbreak fuel}]
\label{thm:ac-rosser}
\leavevmode
\begin{enumerate}
\item $\{\WF(D)\} \vdashI \z = \tofn 1$, and $\emptyset \nvdashI
  \neg\WF(D)$.
\item Let $S$ be the entailment of a step relation whose table
  $\mathcal D$ fires at a node exactly when the node is a step, such
  that $S$ contains \iga's, $\Chk_{\mathcal D}$ meets the structural
  conditions and D1, and $\emptyset \nvdash_S \z = \tofn 1$.  Then $S$
  is not closed under consistency-filtered adoption.
\item In particular, neither \iga nor any consistent claim extension
  $\Ia{a}$, with $a$ in $\KO$ or not, is closed under
  consistency-filtered adoption.
\end{enumerate}
\end{thm}

\noindent
Part~1 restates Theorem~\ref{thm:ac-direct}(1, 2) at $D$: a grounded
derivation and a classical fact.  Parts~2 and~3 are classical.  The textbook Rosser clause
would let $x+1$ lie below $x$ only when no smaller code is a proof of
$\neg\WF$.  Part~1 shows that the step from ``$S$ together with $W$
derives $\z = \tofn 1$'' to ``$S$ derives $\neg W$'' is not a rule of
\iga, and over an $\Mw$-sound $S$ no $\WF$ sentence is ever refuted, so
that clause would be vacuous.  The grounded Rosser order
$R_{\mathrm{Ro}}$ cites proofs of the sequent instead: $x+1$ lies below
$x$ exactly when some $p \leq x$ is accepted as a proof of
$W = \WF(R_{\mathrm{Ro}})$ and no $q \leq p$ is accepted, by the
table's test for judgments, as a proof of $\{W\} \vdash \z = \tofn 1$.
Two derivations carry part~2.  $W$, an accepted code $p$ for $W$ and
the clause at $p$, a bounded universal over finitely many decided
instances, derive $\z = \tofn 1$ in \iga
(\textsf{fuRLR\_WF\_explodes\_of\_hit}).  An accepted proof $q$ of the
sequent, together with the refutation of every code below $q$ as a
proof of $W$, derives $W$ itself in \rga, the order being empty
(\textsf{fuRLR\_WF\_of\_refuted}).  The proofs of part~2 use the
consistency of $S$ and no soundness fact.  The decided form of the
rule at the Rosser order is not stated.

Read together with Theorem~\ref{thm:ac-falsified}, the three theorems
of this subsection say that, in a grounded system, unrefutedness is
no closure condition on adopted claims or on well-foundedness
sentences, and consistency with the system is none on well-foundedness
sentences.  That consistency with the system is none on adopted claims
either is a prediction, expected through a Liar code built like
$\tau_{T,k}$; the development does not state it.  Adoption of a claim
is a supplied hypothesis, checked only through its consequences.  The
three statements are classical; the derivations that carry them are
grounded.

\section{Does \iga alone prove \rga's consistency?}
\label{sec:open}

\igaP derives $\ConR$ (Theorem~\ref{thm:ad-iga-con}), with the help of
a reflection rule.  Whether \iga derives $\ConR$ without that rule is
open.  This section states what is known, names what each route to an
answer would have to supply, and reports the expectations, each
labeled as a prediction.  The results are HOL theorems about
derivability (the classical viewpoint); the derivations that they
build are derivations in \iga or in \rga, and the text names the
system.

\subsection{What is known}
\label{sec:open-known}

\begin{itemize}
\item \igaP derives $\ConR$ and $\ConI$
  (Theorem~\ref{thm:ad-iga-con}).
\item \iga does not derive $\ConI$ (Theorem~\ref{thm:ig-g2}), and \rga
  does not derive $\ConR$ (Theorem~\ref{thm:g-g2}).
\item In \rga, $\ConI \vdash \ConR$
  (Theorem~\ref{thm:ad-iga-warrant}(3)).
\item \oiga, \iga with \oga's rule ATI, derives $\ConR$ from no
  hypothesis, and so strictly exceeds \rga in consistency strength
  (Theorem~\ref{thm:op-oiga-con}).
\item \iga is $\Mw$-sound~\cite{iga-paper}, and
  $\ConR$ is true in $\Mw$ ([\citen{conrga-paper}, theorem~\textsf{goedelW\_omega\_incomplete}]).
\item \iga is proof-grounded: sound and complete for grounded truth
  with \iga's own proofs as the certificates, so \iga has the
  disjunction and witness properties (Theorem~\ref{thm:op-iga-pg}).
\end{itemize}

One relative-consistency fact would settle the question negatively.
If $\ConR \vdashI \ConI$, then a derivation of $\ConR$ in \iga, cut
against that entailment, would derive $\ConI$, which
Theorem~\ref{thm:ig-g2} excludes; so \iga would not derive $\ConR$.  The
development does not prove the premise.

\begin{rem}[predictions, not theorems]
\label{rem:op-relative}
The expectation is that $\ConR \vdashI \ConI$ fails.  The expectation
is also that \iga is not $\Pi_1$-conservative over \rga, that is,
that \iga derives some universal sentence of decided instances that
\rga does not derive, in view of the arithmetical theorems, such as
Goodstein's theorem, that the IGA paper derives in
\iga~\cite{iga-paper} and that are not known to be derivable in \rga.
\end{rem}

\begin{rem}[not mechanized]
\label{rem:op-decided}
A short classical diagonal argument, sketched in the development's
records and not mechanized, indicates that no decidable set of
judgments that is closed under \rga's rules excludes $\emptyset
\vdash \z = \tofn 1$, using that \rga derives every well-formed
closed $\forall$-free term that evaluates to true (open
completeness~\cite{rga-paper}).  If so, a consistency proof, in any
system, must carry undecided truth or a well-founded measure.  A
decided invariant of proofs would not suffice.
\end{rem}

\subsection{What the no-go theorem excludes}
\label{sec:open-nogo}

The internal-reasoning paper's no-go theorem
([\citen{conrga-paper}, theorem~\textsf{xder\_not\_fold\_residual}]) is stated for any entailment over
\rga's terms that is $\Mw$-sound.  \iga's entailment is an instance,
by \iga's $\Mw$-soundness.  Let
$P_{\mathrm h}$ be that theorem's accepted proof with an empty root
context, one of whose nodes is the hypothesis judgment $Y \vdash Y$ for
the ungrounded $Y = W\cdot W$.

\begin{thm}[\textsf{iga\_\allowbreak not\_\allowbreak motive}, \textsf{iga\_\allowbreak not\_\allowbreak fold\_\allowbreak residual}]
\label{thm:op-nogo}
Let $\Gamma$ be a well-formed context whose members are all true in
$\Mw$.  Then $\Gamma \nvdashI \mathit{ResGuard}(\bquo{P_{\mathrm h}})$
(\textsf{iga\_not\_iResGuard}), where $\mathit{ResGuard}$ is the guard
half of the fold's residual; \iga does not derive from $\Gamma$ the
fold's motive at the hypothesis node of $P_{\mathrm h}$, nor the
motive at every rank; and $\Gamma \nvdashI R(\bquo{P_{\mathrm h}})$
for every residual $R$ that meets the fold's arm conditions
$\mathit{arms}(R)$ and for which $R(\bquo{P_{\mathrm h}})$ is
well-formed and does not contain the fold's variable $v_1$ free, all
as in [\citen{conrga-paper}, theorem~\textsf{xder\_not\_fold\_residual}].
\end{thm}

\noindent
The theorem is classical, and the development mechanizes it for \iga.
The development also interprets the no-go locale at \igaP and at
\oiga, the union of \oga's rules~\cite{oga-paper} and \iga's (\rga's
42 rules, ATI and ground induction).  The locale assumes only an
entailment over \rga's terms that is $\Mw$-sound, and the two
interpretations use the two systems' $\Mw$-soundness theorems
(Theorem~\ref{thm:ad-iga-rule}(3) for \igaP, and the development's
soundness theorem for \oiga).  So at \igaP and at \oiga, from the same
contexts $\Gamma$, the guard half, the fold's motive and every residual
with the arm conditions are not derivable either
(Appendix~\ref{app:adopt-nogo}).  The obstruction sits at the inner node $Y \vdash Y$, so
restricting the fold to proofs with an empty root context does not
help.

The theorem rules out the internal-soundness fold with the
internal-reasoning paper's record motive, the motive built from the
value record $\Rcd_v(c)$, ``$c$ has the value $v$ at some fuel''
([\citen{conrga-paper}, definition~\textsf{iRec}]).  The theorem does not rule out
three things:
\begin{itemize}
\item motives not built from the value record, such as a motive built
  from a truth term, a stage-indexed motive or a decided motive (the
  general claim is [\citen{conrga-paper}, the remark on the limits of the fold], which is not
  mechanized);
\item arguments that are not a fold of this shape, such as Gentzen's
  (\S\ref{sec:open-gentzen});
\item a derivation of $\ConR$ itself.
\end{itemize}

\subsection{The Gentzen route}
\label{sec:open-gentzen}

A Gentzen-style consistency proof has two halves: a reduction on proof
codes, proved correct, and a proof that the reduction's descent stops.
\rga proves the shape of the argument (Appendix~\ref{app:goedel-cons}).
From a correct reduction $r$ with ordinal assignment $\mathit{od}$ and order
$\prec$, written $\mathrm{GRed}(r, \mathit{od}, \prec)$, and the
no-descent sentence
\[
  \mathrm{GTI}(F, \mathit{od}, \prec) \;=\; \forall p.\ \neg\forall n.\
  \mathit{od}(F\,p\,(\suc n)) \prec \mathit{od}(F\,p\,n) ,
\]
\rga derives $\ConR$.  Here $F$ is any closed, well-formed term with
$F\,p\,\z = p$ and $F\,p\,(\suc n) = r\cdot(F\,p\,n)$ whenever
$F\,p\,n$ lies outside the fragment of [\citen{conrga-paper}, \S\,``The fragment theorem''], all
three laws, with totality, derived in \rga.  \iga supplies the second
half when, in addition, \rga proves the assignment total and the order
inside the Veblen order, and the assignment is provably bounded below
a notation whose accessibility \iga proves
(Theorem~\ref{thm:op-gentzen}).
For Peano arithmetic the schema has an instance with every premise an \rga theorem~\cite{iga-paper}.
Write $\prec_{\mathrm v}$ for the compiled
order on Veblen notation codes of the IGA paper~\cite{iga-paper}, and
$\#\alpha$ for the code numeral of a notation $\alpha$.

\begin{thm}[\textsf{gi\_\allowbreak ConRGA\_\allowbreak of\_\allowbreak gentzen}, \textsf{gi\_\allowbreak GTI\_\allowbreak separation}]
\label{thm:op-gentzen}
Suppose $\mathrm{GRed}(r, \mathit{od}, \prec)$, that $F$ iterates $r$
as above, that \rga proves $\mathit{od}$ total and $\prec$ contained in
$\prec_{\mathrm v}$, and that \rga proves the bound $\forall q.\
\mathit{od}(q) \prec_{\mathrm v} \#\alpha$ for a valid notation
$\alpha$ below $\varphi(\omega, 0)$.  Then:
\begin{enumerate}
\item $\emptyset \vdashI \ConR$;
\item with the bound derived in \iga instead,
  $\emptyset \vdashI \ConR$ (\textsf{gi\_ConRGA\_of\_gentzen\_gen}),
  $\emptyset \vdashI \mathrm{GTI}(F, \mathit{od}, \prec)$, and
  $\emptyset \nvdash \mathrm{GTI}(F, \mathit{od}, \prec)$.
\end{enumerate}
\end{thm}

\noindent
\emph{Proof.}  The descent comes first.  \iga proves the
accessibility of $\#\alpha$, since $\alpha$ is valid and below
$\varphi(\omega,0)$~\cite{iga-paper}.  Accessibility passes down to
$\mathit{od}(F\,p\,\z)$ through the bound.  \iga's transfinite
induction along an accessible order~\cite{iga-paper} then applies at
the motive
\[
  M(z) \;=\; \forall p\, \forall k.\ \neg\bigl(\mathit{od}(F\,p\,k) = z\bigr)
    \lor \neg\forall n.\ \mathit{desc}(p, n) ,
\]
where $\mathit{desc}(p, n)$ says that step $n$ from $p$ descends.  The
progressiveness of $M$ is \rga reasoning with one decided split, on
$\mathit{od}(F\,p\,k) = z$, which is decided because $\mathit{od}$ is
total, and with no split on the descent atom.  $M$ at
$\mathit{od}(F\,p\,\z)$ and $k = 0$ gives GTI.  The \rga facts move
into \iga by embedding and cut.  For part~1, \rga's Gentzen shape
restated with GTI as a hypothesis, $\mathrm{GTI} \vdash \ConR$, is
embedded into \iga and cut against the descent.  The last statement
of part~2 is \rga's no-Gentzen-pair theorem
(Appendix~\ref{app:goedel-cons}), which rests on
Theorem~\ref{thm:g-g2}.  \qed

\medskip
\noindent
Read through Theorem~\ref{thm:g-g2}, part~2 is a separation template:
\iga derives a no-descent sentence that \rga does not derive.  The
theorem is classical, a HOL implication between derivability facts;
the derivation of $\ConR$ that it builds is a derivation in \iga.  The Isabelle theory behind Theorem~\ref{thm:op-gentzen}
checks at build time that nine results use no oracle, that part~1
reaches \iga's transfinite induction, the accessibility theorem, the
descent, \rga's Gentzen shape with GTI as a hypothesis and the chain
lemma (every iterate of an accepted proof of $\z = \tofn 1$ is one),
and, as a control, that \rga's Gentzen shape reaches the chain lemma
but not \iga's transfinite induction.

What remains for \iga on this route is exactly the premise
$\mathrm{GRed}(r, \mathit{od}, \prec)$ with $\mathit{od}$ bounded below
$\varphi(\omega, 0)$: an ordinal analysis of \rga.  No ordinal analysis
of \rga is known.  The RGA paper leaves \rga's proof-theoretic
strength open~\cite{rga-paper}.  The two further premises, that $\mathit{od}$ is
total and that $\prec$ lies inside $\prec_{\mathrm v}$, are mild for a
primitive-recursive $\mathit{od}$ and $\prec$.

\begin{rem}[predictions, not theorems]
\label{rem:op-gentzen}
The reduction is expected to be hard.  Conversion makes Gentzen's
degree measure ill-founded: $Y$ reduces to $\forall\lambda.\,Y$,
whose instance is $Y$ again.  The rank that decreases under reduction
is expected to be the grounding stage of the cut formula, which is not
primitive recursive, so a primitive-recursive substitute would have to
be read off the proof, and the development builds no such substitute.  The route
is closed for \iga if \rga has the strength of the paracomplete
Kripke--Feferman truth theory PKF~\cite{halbach06axiomatizing}.
\end{rem}

\subsection{A stage-indexed truth}
\label{sec:open-stage}

The no-go theorem leaves open an internal-soundness fold whose motive
is a \emph{staged truth} $\mathit{TS}(\beta, c, r)$, a copy of \iga's
truth term~\cite{iga-paper} in which every point carries a
\emph{stage}~$\beta$.  A stage is the code of an ordinal below
$\varepsilon_0$ in Cantor normal form, as in the IGA paper's
consistency proof for Peano arithmetic~\cite{iga-paper}, and
$\mathit{TS}(\beta, c, r)$ says that the term with code $c$ has the
value $r$ at stage $\beta$.  The clauses of the staged truth are those
of the truth term, generated from the same table of clauses, except
that every truth atom in a clause's body is read at a stage below the
point's own.  At the node $Y \vdash Y$, $\mathit{TS}(\alpha, \bquo{Y},
1)$ is false in $\Mw$ at every $\alpha$, a classical observation that
the development does not state as a theorem, so the negated hypothesis
is true there and the no-go theorem's argument does not apply.  With
this motive the fold is a theorem in \oiga
(Theorem~\ref{thm:op-oiga-con}).  The statements collected in
Remark~\ref{rem:op-stage}, about \iga alone and about \iga's own
checker, are predictions.

\noindent
A fold of this shape, without stages, is a theorem at a different
rule.  \iga with bivalence of the truth atom as a rule, which the IGA
paper~\cite{iga-paper} calls its \emph{Kripke universe} at level~$1$,
derives $\ConR$ from no hypothesis (\textsf{itb\_ConRGA}; a derivation
in the extended system).  The name follows Kripke's closing-off of
truth~\cite{kripke75outline}; a Kripke universe is a system, unrelated
to a Grothendieck universe of sets, and Kahle's universes over Frege
structures~\cite{kahle03universes} are its closest precedent.  The
fold runs through \rga's proof checker with the truth atom at the
compiled code of each node's conclusion as its motive; every step is a
derivation in \rga under four closed formulas that the system derives,
three by ground induction over the truth term and one by the bivalence
rule, and a build-time check confirms that the rule enters only there
and that no soundness theorem of any system is in the proof's cone.
Neither \iga nor \oiga has the bivalence rule.

\begin{thm}[\textsf{oiga\_\allowbreak ConRGA}, \textsf{oiga\_\allowbreak separation}]
\label{thm:op-oiga-con}
\leavevmode
\begin{enumerate}
\item $\emptyset \vdash_{\mathrm{OIGA}} \ConR$.
\item $\emptyset \nvdash \ConR$ and $\emptyset \vdash_{\mathrm{OIGA}}
  \ConR$.
\end{enumerate}
\end{thm}

\noindent
Part~1 is a derivation in \oiga from no hypothesis,
with no rule added to \oiga.  Part~2 is a classical corollary: \rga
does not derive $\ConR$ (Theorem~\ref{thm:g-g2}), so \oiga strictly
exceeds \rga in consistency strength.  A build-time check confirms
that G\"odel's second theorem is in the cone of part~2 and not in the
cone of part~1.

The proof of part~1 is the fold through \rga's checker table with the
staged truth in place of the truth atom.  Each node of an accepted
proof carries a \emph{node truth} with one stage $\delta$, chosen
before the environment: for every environment that closes the node
and every stage $\alpha$ at which the node's hypotheses are true, the
node's conclusion is true at the stage $F(\delta, \alpha) =
\omega^{e(\alpha \sqcup \omega) \oplus \delta}$, with $e(\gamma)$ the
leading exponent of $\gamma$, $\sqcup$ the maximum and $\oplus$ the
natural sum (\textsf{stFam}).  The stage comes first because the rule
that introduces a universal into the staged truth takes all the
universal's instances at one stage (\textsf{tsr\_allTI}).  Each arm of
the fold, one per rule detector of \rga's table, is a derivation in
\rga under four closed formulas (\textsf{ndT\_fold\_RGA}): the
\emph{staged bivalence}, that $\mathit{TS}(\beta, c, r)$ is decided at
every Cantor-normal-form stage $\beta$ (\textsf{tsBivF}); the staged
head step, which carries the determinacy of the truth term; the
inclusion of the staged truth in the truth term; and the invariance of
the staged truth when an environment entry is swapped for one with the
same value, at a stage raised by the stage of that value.  \oiga
derives the four formulas from no hypothesis: the staged bivalence by
one transfinite induction along the Cantor-normal-form order below
$\varepsilon_0$, in which ATI decides each clause of the staged truth
at the clause's quantifiers (\textsf{ts\_biv}), and the other three by ground induction over the staged truth term,
with the implications introduced over the staged bivalence
(\textsf{ft\_g\_kit\_hs}, \textsf{ft\_g\_kit\_sub},
\textsf{ft\_g\_kit\_sw}).  A build-time check confirms that the arms and the fold use
neither ground induction nor ATI, so that both rules enter the proof
only through the derivations of the four formulas, and that no
soundness theorem of any system, no bivalence of the truth atom and no
oracle is in the cone of part~1.

Both derivations instantiate one generic table fold at different node
truths, the arms re-proved for the staged one; each is a derivation in
\rga under four closed formulas of its own, which its system supplies.  \iga
with the bivalence rule takes the bivalence of the truth atom as a
rule, while in \oiga ATI supplies the bivalence of the staged truth,
one clause at a time along a transfinite induction.

\noindent
The fold at the bivalence rule extends to \iga's own checker.  \iga
with the bivalence rule derives $\ConI$ from no hypothesis
(\textsf{itb\_ConIGA}; a derivation in the extended system): the table gains the
ground-induction rule, whose arm is a derivation in the plain
fixed-point calculus by one ground induction with a
positivity-indexed motive (\textsf{fc\_T\_gen}), and the build-time
check extends to it.  With Theorem~\ref{thm:ig-g2}, closing off the
truth atom strictly exceeds \iga in consistency strength, a classical
corollary (\textsf{fc\_separation}): \iga does not derive $\ConI$, and
\iga with the bivalence rule does.  This separation is the analogue,
at a rule in place of an adoption, of the separation of \igaP from
\iga (Theorem~\ref{thm:ad-iga-sep}(2)).

\noindent
The closing-off iterates.  The IGA paper~\cite{iga-paper} presents the
truth term without the numeral constant $0$, so that no closed subterm
of it has a value, and builds a tower of \emph{constant-free Kripke universes}: the system
at each level closes off the constant-free truth term of the level
below.  Its level~$1$ is \iga with bivalence of the constant-free
truth term, which need not be the system above (the end of this
paragraph returns to the question).  The fold at a constant-free Kripke
universe's own checker, \iga's checker with one rule added per element
of the level's list, derives at every level the consistency sentence
of the constant-free Kripke universe below (\textsf{cf\_lvl\_Con}; a
derivation in that Kripke universe).  G\"odel's second theorem holds at
every constant-free Kripke universe, because its checker is faithful
and satisfies the shift law of its table, so that the conditions of
\S\ref{sec:locale-cond} hold at it by the same transformers
(\textsf{cf\_lvl\_goedel2}; a classical theorem, with the constant-free
Kripke universe's consistency from its soundness for the closed-off model).  Hence each
closing-off of the truth predicate strictly increases consistency
strength, by at least the consistency of what it closes off
(\textsf{cf\_separation}, \textsf{cf\_tower\_strict}).  No theorem says
that a constant-free Kripke universe adds exactly that consistency, and the tower is stated
for the constant-free truth term; whether \iga with the bivalence of
the truth term of the IGA paper and \iga with that of its constant-free
presentation are the same system, or, equivalently by a classical
theorem of the IGA paper, prove the same sentences in which neither
truth term occurs, is open there.

The staged fold is not known to extend to \iga's own checker: by the
argument below, it stops at a named obligation, and the development
does not claim $\emptyset \vdash_{\mathrm{OIGA}} \ConI$.
The obstruction rests on classical theorems of the IGA
paper~\cite{iga-paper}.  \iga and \oiga both have the accessibility
ordinal $\varphi(\omega, 0)$: each proves the accessibility of exactly
the Veblen notations below $\varphi(\omega, 0)$, and \oiga does not prove the Veblen order well-founded on any closed
domain true in $\Mw$ at the code of a notation at or above
$\varphi(\omega, 0)$ (\textsf{oiga\_\allowbreak not\_\allowbreak
WF\_\allowbreak above\_\allowbreak pw}).  In the stage-indexed omega model of the IGA paper, the accessibility
of a notation $x$ is true at no stage below $x$ (the rank lemma
\textsf{acc\_stage\_ge}); reading this as a bound on the stages of the
staged truth is part of the argument that follows, which the
development does not mechanize.  By that argument, a fold over \iga's checker would therefore need the
staged truth at stages cofinal in $\varphi(\omega, 0)$, coded in the
Veblen notation, and the staged bivalence at all of those stages, the
obligation that this paper names \textsf{ts\_biv\_pw}.  The
transfinite induction behind \textsf{ts\_biv} would then have to reach
$\varphi(\omega, 0)$, where \oiga's accessibility stops.  So Cantor
normal forms below $\varepsilon_0$ suffice for \rga's checker, and
\iga with the bivalence rule remains the result for $\ConI$.

\begin{rem}[predictions, not theorems]
\label{rem:op-stage}
\leavevmode
\begin{enumerate}
\item \iga is expected not to derive the staged bivalence.  \iga
  lacks ATI, which is what decides the clauses in the derivation of
  Theorem~\ref{thm:op-oiga-con}(1).  If \iga derived the four closed
  formulas of that derivation, then \iga would derive $\ConR$, by embedding the fold's derivation in
  \rga and cutting (a classical consequence of
  \textsf{fki\_staged\_ConRGA\_Gk}, the embedding of \rga and cut, not
  stated as one theorem).  That the other
  three formulas follow in \iga from the staged bivalence alone, as
  they do in \oiga, is a prediction.  From the staged bivalence, \iga
  would also derive that every member of the halting family behind
  [\citen{conrga-paper}, theorem~\textsf{rga\_not\_iLFI}] has a decided truth value, and
  Theorem~\ref{thm:op-iga-pg}(3) excludes both the decidedness of
  some member and the uniform sentence.  The step from the staged
  bivalence to the halting family's bivalence is the prediction; the staged truth is defined, its bivalence in \oiga is a theorem
  (\textsf{ts\_biv}), and the exclusion is a theorem.
\item \oiga is expected not to derive the obligation
  \textsf{ts\_biv\_pw}, in view of the theorems on \oiga's
  accessibility ordinal above.
\end{enumerate}
\end{rem}

\noindent
\iga has the disjunction property.  As for \rga
(Theorem~\ref{thm:br-mech}(1)), the property rests on a semantics for
which the system is sound and complete: grounded truth as in
\S\ref{sec:bg-rga}, with a proof in \iga's own proof system, in place
of one in \rga's, as the certificate at each universal.  A system that
is sound and complete for grounded truth read with its own proofs as
the certificates is \emph{proof-grounded}.  \rga is
proof-grounded~\cite{rga-paper}, and so is
\iga~\cite{iga-paper}.  Let $K$ be the set of
programs $x$ that halt on input $x$, and let $\mathit{ofam}(x)$ be the
halting family of [\citen{conrga-paper}, theorem~\textsf{iLFI\_decides\_halting}], whose instance at
$n$ says that program $x$ on input $x$ has not halted by step $n$.

\begin{thm}[\textsf{gi\_\allowbreak disj\_\allowbreak prop}, \textsf{gi\_\allowbreak not\_\allowbreak halting\_\allowbreak biv}]
\label{thm:op-iga-pg}
\leavevmode
\begin{enumerate}
\item For a well-formed $c$, $\emptyset \vdashI c$ if and only if $c$
  is grounded-true with \iga's proofs as the certificates.  In
  particular, $\emptyset \nvdashI \z = \tofn 1$.
\item If $\emptyset \vdashI a \lor b$, then $\emptyset \vdashI a$ or
  $\emptyset \vdashI b$.  If $\emptyset \vdashI \neg\forall x.\,\varphi$,
  then $\emptyset \vdashI \neg\varphi(\tofn n)$ for some numeral
  $\tofn n$.
\item \iga's derivability is recursively enumerable.  For some $x
  \notin K$, $\emptyset \nvdashI T_x\jbool$, where
  \[
    T_x \;=\; \forall n.\ \Rcd_1(\mathit{ap}(\bquo{\mathit{ofam}(x)}, \Num\cdot n))
  \]
  is the family sentence of [\citen{conrga-paper}, theorem~\textsf{iLFI\_decides\_halting}].  So $\emptyset
  \nvdashI \forall c.\ \bigl(\forall n.\ \Rcd_1(\mathit{ap}(c,
  \Num\cdot n))\bigr)\jbool$, the uniform decidedness of the family
  sentences.  Also $\emptyset \nvdashI \LF$, for the closed sentence
  $\LF$ of [\citen{conrga-paper}, definition~\textsf{iLFI}].
\item For some $x \notin K$, $\emptyset \nvdashI
  (\cA\cdot\mathit{ofam}(x))\jbool$, while \oga derives
  $(\cA\cdot\mathit{ofam}(x))\jbool$ for every $x$.  So some sentence
  that \oga derives, and \oiga with it, is not derivable in \iga.
\end{enumerate}
\end{thm}

\noindent
The theorem is classical, a set of HOL statements about derivability
in \iga.  Part~1 is the soundness and completeness of \iga for the
semantics above, and part~2 follows from part~1 as \rga's two
properties follow from \rga's semantics (\S\ref{sec:breaks-mech}).
Part~3 is the argument of [\citen{conrga-paper}, theorem~\textsf{rga\_not\_iLFI}] run in \iga:
if \iga derived the decidedness of $T_x$ at every $x \notin K$, then,
through parts~1 and~2, the complement of $K$ would be recursively
enumerable.  The proofs of the non-derivability statements of part~3,
and of the undecided member of part~4, use \iga's consistency,
disjunction property and witness property from parts~1 and~2, never
an $\Mw$-model of \iga; a build-time check confirms that these proofs
reach no $\Mw$-soundness theorem of \iga.  In part~4, \oga derives
the decidedness of every member of the halting family with its rule
ATI~\cite{oga-paper}, and \oiga has \oga's rules.  Whether \oga
contains \iga is open.

So the question of this section has two named obligations, one per
route: an ordinal analysis of \rga below $\varphi(\omega, 0)$ for the
Gentzen route, and, for the stage-indexed route, a derivation in \iga
of the four closed formulas of the fold, the staged bivalence first,
which \oiga derives with ATI (Theorem~\ref{thm:op-oiga-con}).  Whether
\iga derives the staged bivalence is open:
Theorem~\ref{thm:op-iga-pg}(3) excludes the uniform decidedness of the
halting family's sentences, which the staged bivalence is predicted to
give (Remark~\ref{rem:op-stage}(1)).  This
analysis refines the internal-reasoning paper's open item on \oga and
$\ConR$ ([\citen{conrga-paper}, \S\,``Conclusion'']) and the IGA paper's closing remarks on Gentzen's
route~\cite{iga-paper}.

\section{Mechanization}
\label{sec:mech}

Every theorem in this paper is proved in Isabelle/HOL (Isabelle2025-2).
The development builds on the RGA paper's formalization of
\rga~\cite{rga-paper}, on the internal-reasoning paper's sessions for
the coding kit, the decider and the checker~\cite{conrga-paper}, and
on the IGA paper's formalization of \iga~\cite{iga-paper}, with its
sessions \texttt{IGA\_PGSem} and \texttt{IGA\_RE} for \iga's semantics
with proof certificates (Theorem~\ref{thm:op-iga-pg}).  The
G\"odel theories for \rga, the session \texttt{RGA\_Goedel} of
Sections~\ref{sec:breaks} and~\ref{sec:goedel}, hold 11,287 lines in
ten theory files.  The sessions new in this paper hold the following
numbers of lines (\texttt{wc -l}).
\begin{itemize}
\item \texttt{GCalc\_Core}, the rule list and the embedding of the
  calculus of \S\ref{sec:locale-calc}: 665 lines in one theory.  The
  locale \textsf{gcalc} itself and its kit sit in the theory
  \texttt{RGA\_GCalc} of the \rga session (1,689 lines).
\item \texttt{G2\_Locale}, the nine locales of the G\"odel conditions
  (Section~\ref{sec:locale}): 1,093 lines in one theory.  The
  interpretations at \rga, with the 62 instance facts, are theories of
  \texttt{RGA\_Goedel}.
\item \texttt{RGA\_Checker\_Gen}, the checker at a table
  (Section~\ref{sec:checker}): 1,893 lines in two theories.
\item \texttt{IGA\_Checker}, \iga's checker: 1,211 lines in three theories.
\item \texttt{IGA\_Goedel}, the second theorem at \iga
  (Section~\ref{sec:iga-goedel}): 1,743 lines in four theories
  (the fourth a build-time audit of the decider-free cones).
\item \texttt{RGA\_Adopt} and \texttt{IGA\_Adopt}, adoption at \rga and
  at \iga (Section~\ref{sec:adopt}): 796 and 972 lines, one theory
  each.
\item \texttt{IGA\_NoGo}, the no-go theorem at \iga (Section~\ref{sec:open})
  and at \igaP and at \oiga (\S\ref{sec:open-nogo},
  Appendix~\ref{app:adopt-nogo}): 49 and 76 lines in two theories.
\item \texttt{IGA\_Gentzen}, the Gentzen conditional
  (Section~\ref{sec:open}): 578 lines in its theory
  \texttt{Gentzen\_IGA}.
\item \texttt{IGA\_Fold}, \texttt{IGA\_Fold\_Con} and
  \texttt{OIGA\_Fold}, the folds of \S\ref{sec:open-stage}: the fold
  at the bivalence rule through \rga's checker and through \iga's, and
  the staged fold in \oiga (Theorem~\ref{thm:op-oiga-con}): 19,657,
  10,882 and 24,952 lines in 33, 17 and 55 theories.
\item \texttt{IGA\_Fuel}, the adopted notation claims of
  \S\ref{sec:adopt-claims}--\S\ref{sec:adopt-ac}: the falsification
  theorems, Turing's code, G\"odel's second theorem at every claim, the
  quine code and the adequacy rules: 3,970 lines in five theories
  (\texttt{Fuel\_First}, \texttt{Fuel\_Turing},
  \texttt{Fuel\_Goedel}, \texttt{Fuel\_Quine} and
  \texttt{Fuel\_Liar}).  The parent of \texttt{IGA\_Fuel} is the
  development's session \texttt{PG\_Sets}, itself over
  \texttt{IGA\_Adopt}, and its sessions line adds \texttt{IGA},
  \texttt{IGA\_Closing} and \texttt{IGA\_Goedel}.
\end{itemize}
These sessions add theories on top of the earlier development, apart
from one change of imports in the session \texttt{PG\_Sets}: its
theory \texttt{PW\_Witness} no longer imports the theory
\texttt{PG\_Sets}, and its theorems were re-derived with unchanged
statements.  So they change
none of the statements that the internal-reasoning paper prints.  The sessions contain no \textsf{sorry}.

\paragraph{Trusted base.}
Besides Isabelle/HOL's kernel, a reader must trust the definitions in
which the theorems are stated.  For \rga these are the definitions
that the internal-reasoning paper lists ([\citen{conrga-paper}, Section ``Mechanization'']):
\rga's derivability relation with its 42 rules, the coding of terms,
the compiler $\cmp$, the checker $\Chk$ with its validity test,
$\ConR$ and the code $\bquo{\z = \tofn 1}$.  This paper adds five
groups.
\begin{itemize}
\item \iga's derivability relation, with ground induction as its 43rd
  rule~\cite{iga-paper}, and the derivability relations of \rgaP and
  \igaP, each with one adoption rule (Definition~\ref{def:ad-rule},
  \S\ref{sec:adopt-iga}); for the folds of \S\ref{sec:open-stage},
  the derivability relations of \iga with the bivalence rule of the
  truth atom and of \oiga.
\item \iga's checker $\ChkI$ and $\ConI$ (Definitions~\ref{def:ck-iga}
  and~\ref{def:ig-con}).  $\ChkI$ unfolds to \rga's definitions, the
  generic checker's definitions at a table, and the explicit index
  trees of the ground-induction tag and the positivity decider
  (Appendix~\ref{app:ck-tcb}).  As for $\Chk$, closed faithfulness and
  genuineness (Theorem~\ref{thm:ck-iga-faith}) pin $\ChkI$ down
  extensionally at closed codes.
\item The warrants $\WarR$ and $\WarI$ and the instance code $\Inst$
  (Definition~\ref{def:ad-rule}), which the adoption rules read, and
  the $\omega$-model $\Mw$ of the RGA paper, for the soundness theorems
  of Section~\ref{sec:adopt}.
\item For \S\ref{sec:adopt-claims}--\S\ref{sec:adopt-ac}, the IGA
  paper's fixed point $\KO$ and the claim extensions $\Ia{a}$
  (Definition~\ref{def:ac-ext}), the axiom cell, $\ChkA{a}$ and
  $\ConA{a}$ (Definition~\ref{def:ac-chk}), the codes $\tau_{T,k}$
  (Definition~\ref{def:ac-turing}), and the sentences $\WF(R)$ and the
  rule forms of \S\ref{sec:adopt-ac}.
\item The locales of Section~\ref{sec:locale}, for the theorems stated
  once for every G\"odel extension.  Appendix~\ref{app:locale-cond}
  gives their assumptions as the development states them.  At \rga and
  at \iga the interpretations rewrite the locale's consistency sentence
  and diagonal function to $\ConR$ or $\ConI$ and to $\Diag$, so the
  instance theorems are stated in the trusted definitions above, not
  in the locale's.
\end{itemize}

\paragraph{Assurance.}
The sections above describe the build-time checks that their theorems
pass, and each check fails the build on a violation.  In summary:
\begin{itemize}
\item \emph{No oracle.}  Each theory of the new sessions checks that
  its headline facts depend on no oracle (in \texttt{IGA\_Fuel}, four
  of its five theories; \texttt{Fuel\_First} has no check)
  (Appendices~\ref{app:locale-checks}, \ref{app:iga-goedel-checks},
  \ref{app:ck-checks} and~\ref{app:adopt-checks}).
\item \emph{Statement identity.}  The 62 instance facts at \rga keep
  the propositions that they had as \rga theorems before they became
  instances: a ground-truth table of 130 fingerprints, with a control,
  fixes them, and two instances must have literally the propositions
  of facts that keep their own \rga proofs, among them G\"odel's second
  theorem (\S\ref{sec:locale-rga}, Appendix~\ref{app:locale-checks}).
\item \emph{Soundness footprint.}  The dependency cones of the internal
  theorems avoid the soundness and consistency theorems by name: at
  \rga, the 60 instance facts outside the two consistency layers
  avoid nine named soundness, consistency, faithfulness and truth facts
  (Appendix~\ref{app:locale-checks}); at \iga, eight theorems, among
  them formalized G\"odel~II (Theorem~\ref{thm:ig-internal}), avoid
  eleven soundness and consistency theorems, among them those of \rga,
  \iga and \igaP and the genuineness of $\ChkI$,
  and through their layers' witnesses the check covers every witness
  except consistency (Appendix~\ref{app:iga-goedel-checks}); and the
  derivations of $\ConR$ in \rgaP and of $\ConI$ in \igaP avoid the
  soundness lemmas, with a control (Appendix~\ref{app:adopt-checks}).
  The proofs of the non-derivability statements of
  Theorem~\ref{thm:op-iga-pg}(3), and of the undecided member in that
  theorem's part~4, reach no $\Mw$-soundness theorem of
  \iga, with controls (\S\ref{sec:open-stage}).  The cones of the
  derivations of $\ConR$ and $\ConI$ by the fold, at the bivalence
  rule and in \oiga, contain no soundness theorem of any system, with
  controls; in \oiga's, ATI and ground induction enter only through the
  derivations of the fold's four closed formulas
  (\S\ref{sec:open-stage}).
\item \emph{Reach.}  Need lists confirm that the proofs go through the
  intended constructions, for example that \iga's graft comes from the
  generic graft and the members' shift laws and not from \rga's
  whole-table graft (Appendix~\ref{app:ck-checks}), that G\"odel's
  second theorem at \iga reaches the \iga witnesses
  (Theorem~\ref{thm:ig-g2}), and that the Gentzen conditional reaches
  \iga's transfinite induction (\S\ref{sec:open-gentzen}).
\item \emph{Decider-free statements.}  The internal-reasoning paper's
  decider check ([\citen{conrga-paper}, Appendix ``The build-time check'']) runs on 31
  facts of \rga's G\"odel theories, 18 of them instance facts of 17
  locale facts (Appendix~\ref{app:locale-checks}).  The same check
  runs on the 122 facts of the layered G\"odel locales by name, on
  \iga's G\"odel theorems and their two top interpretations, nine
  roots, and on the derivations of $\ConI$ in \igaP and of $\ConR$ in
  \rgaP.  No node of any of these 133 dependency graphs states the
  truth predicate, a decider verdict, the search or a constant of the
  decider's table.  None reaches the search witness's D3G, the decided
  search or the decider's completeness, and three controls for each
  part of the check are flagged (theory
  \texttt{Goedel\_\allowbreak IGA\_\allowbreak Audit}, with the
  shared check in the session \texttt{RGA\_\allowbreak
  Goedel\_\allowbreak Audit}).
\item \emph{Adopted claims.}  Four of the five theories of
  \texttt{IGA\_Fuel} check at build time that their headline facts
  use no oracle and that the grounded derivations, among them
  Theorem~\ref{thm:ac-turing}(1), Theorem~\ref{thm:ac-quine}(3),
  Theorem~\ref{thm:ac-g2}(3) and the derivations along the Liar and
  Rosser orders, avoid the soundness, consistency and $\Mw$ facts by
  name.  Need lists confirm that the quine's derivation goes through
  induction along $\KO$ and the Liar's through the descending-sequence
  theorem.  Avoid checks confirm, by a scan for choice-defined
  equations and a check on the pairing index, that no choice-defined
  index occurs in the dependency cones of the quine code's laws, and
  that the choice-defined index of the earlier halting witness is
  absent from the cones of the laws of the codes $\tau_{T,k}$.  At the code $a_0$, the
  second theorem at the claim extension is checked statement-identical
  to its derivation from G\"odel's second theorem at \iga, and the
  first route is checked not to use that theorem.
\end{itemize}
The checks see whole theorems, not directions of an equivalence, which
is why the proofs at \iga use one-way closed evaluation for closed
acceptance (\S\ref{sec:ig-closed}).  Statements marked as remarks or
predictions are not mechanized.

\paragraph{Artifact.}
\label{sec:mech-artifact}
The development is released at one tag of the project repository, the
same tag that the internal-reasoning paper cites, and the repository's
README gives the build command and the sessions.  The command
\begin{center}
\small\texttt{isabelle build -b -v -o threads=10 -d hol IGA\_Goedel IGA\_Gentzen IGA\_NoGo IGA\_RE}\\
\small\texttt{IGA\_Fold\_Con OIGA\_Fold IGA\_Fuel}
\end{center}
builds every session that this paper's theorems need, the earlier
sessions included.
The theorem names and theory files that this paper gives are those of
the development at that tag.

\section{Related work}
\label{sec:related}

The design genealogy of \rga itself (Fitch's basic logic, illative
combinatory logic, the recursive $\omega$-rule, Kripke's fixed points)
is discussed in the RGA paper~\cite{rga-paper}, and the internal
machines on which this paper's proofs run, with their relatives in
self-verified proof assistants and in theories of truth, in the
internal-reasoning paper~\cite{conrga-paper}.  This section places the
present paper's results.

\paragraph{Mechanized incompleteness.}
Paulson gave the first mechanized proof of G\"odel's second theorem,
in Isabelle/HOL for the theory of hereditarily finite
sets~\cite{paulson14machine,paulson15mechanised}, following
\'Swierczkowski~\cite{swierczkowski03finite}.  Popescu and Traytel
distilled abstract versions of both theorems, with the derivability
conditions as assumptions, and instantiated them using Paulson's
proofs~\cite{popescu19formally,popescu21distilling}.  Popescu and
Traytel's framework assumes only intuitionistic reasoning, but its
instances are classical, and it assumes implication introduction,
which \rga lacks at undecided hypotheses.  Saitou and Noguchi proved
the second theorem in Lean~4 for every $\Delta_1$ theory extending
$\mathrm{I}\Sigma_1$, establishing D3 by a model-theoretic argument
through the completeness theorem~\cite{saitou26mechanizing}.  Coquand
formalized in Agda the second theorem for Church's basic recursive
arithmetic, following Guard's unpublished lecture notes of 1963, which
follow Church's Princeton lectures~\cite{coquand26auto}; the report's
title has ``binary'' where its text has ``basic''.
All of these mechanizations concern classical, or free-variable and
hence everywhere decided, systems.  Carnielli and Fuenmayor analyzed
the incompleteness theorems under weak logics of formal inconsistency
in a shallow embedding in Isabelle/HOL, with the fixed point taken as
a premise and the derivability conditions cited rather than proved,
and found that the second theorem then needs extra premises asserting
the consistency of particular formulas~\cite{carnielli20goedel};
Carnielli and Fuenmayor's setting is abstract, with no arithmetic.
Gross, Gallagher and Fallenstein formalize L\"ob's theorem in Agda for
a postulated, syntactic provability predicate, and read the theorem,
through Curry--Howard, as a limit on self-interpreters for total
languages~\cite{gross16lob}.  The mechanized proofs of the second
theorem known to the author are for classical
theories~\cite{paulson15mechanised,saitou26mechanizing,hydras-goedel2}
or for quantifier-free classical arithmetic~\cite{coquand26auto}; the
non-classical ones are abstract, with the fixed point or the
derivability conditions assumed~\cite{popescu21distilling,carnielli20goedel}.
No mechanized proof of the second theorem is known for Heyting
arithmetic, and none for a concrete arithmetic with quantifiers whose
object logic is not classical.\footnote{The search covered the Lean~4
Foundation library, the Rocq developments hydra-battles and the Coq
library for first-order logic, Agda, the Isabelle Archive of Formal
Proofs and HOL Light, and, by search, Metamath, Mizar, HOL4, Nuprl and
Minlog, as of October 2026.}
Several developments mechanize the first theorem only: Shankar's in
the Boyer--Moore prover~\cite{shankar94metamathematics}, Harrison's in
HOL Light~\cite{harrison-hol-light-godel}, and the synthetic
treatments of undecidability and essential incompleteness by Kirst and
coauthors, in Coq~\cite{kirst23synthetic,kirst23godel}.  O'Connor's
development in Coq~\cite{oconnor05essential} proved the first theorem;
its maintained version also derives the second for extensions of Peano
arithmetic, from D2 and D3 taken as hypotheses~\cite{hydras-goedel2}.
The present development also differs in architecture: Isabelle/HOL
verifies statements about \rga and \iga, and \rga in turn reasons
about its own codes and about \iga's, through checkers written as
\rga terms.

\paragraph{Abstract and uniform second theorems.}
Beklemishev and Shamkanov~\cite{beklemishev16abstract} prove the
second theorem, and its formalized version, for abstract provability
structures: a consequence relation with provability and refutability
operators $\Box$ and $\boxtimes$ and no connectives, under conditions
C1--C4, with a Jeroslow fixed point $p \eqv \boxtimes p$.  Condition C1
asks that $\Box$ be monotone and $\boxtimes$ anti-monotone in the
consequence order.  Beklemishev and Shamkanov's $\top$ is read here as
$\z = \z$ and their $\bot$ as $\neg(\z = \z)$, a classical falsum, not
\rga's divergent $\bot$.  Beklemishev and Shamkanov's formalized
theorem is itself an entailment, $\boxtimes\boxtimes\top \leqslant
\boxtimes\top$, the shape of Theorem~\ref{thm:g-formal}(1).  In \rga,
with $x \leqslant y$ read as $x \vdash y$ and $\boxtimes x$ as
$\PrC(\bquo{\neg x})$, conditions C2 to C4, Beklemishev and Shamkanov's
C5 and the monotone half of C1 hold at every sentence.  The
anti-monotone half of C1, contraposition under provability, fails in
general, at the Liar, but holds at decided sentences and at $\Sigma$
sentences with a decided matrix.  That class is closed under $\Box$
and $\boxtimes$ and contains the Jeroslow fixed point $\neg G$, so on
that class \rga is an abstract provability structure, and Beklemishev
and Shamkanov's Theorem~1 holds for \rga by their proof, mechanized
(\S\ref{app:goedel-ent} gives the statements).  So the contrapositive
form of \rga's formalized second theorem, $\PrC(\bquo{\ConR}) \vdash
\neg\ConR$, is an instance of Beklemishev and Shamkanov's Theorem~1 on
that fragment: the $\boxtimes\top$ of that theorem is \rga's
$\neg\ConR$ up to interderivability, and necessitation of entailments
(Theorem~\ref{thm:g-cut}) with $\Sigma$-contraposition moves the
instance under $\PrC$.  The G\"odel conditions give Beklemishev and
Shamkanov's structure conditions too: in every G\"odel extension that
also meets the ex falso condition, C1 (its monotone half) to C5 hold
(Theorem~\ref{thm:lc-bs}).  With an implication obeying the
deduction theorem, Beklemishev and Shamkanov also show that
contraction restricted to provability formulas is what G\"odel's
argument needs, and they give a contraction-free modal system in which
formalized G\"odel~II fails.  \rga lacks that deduction theorem and
keeps contraction.  Beklemishev and Shamkanov's framework does not
provide a concrete arithmetic with gaps that meets its conditions;
\rga and \iga are two, built and checked inside the logic.

Visser discusses how to state the second theorem in a general,
coordinate-free form~\cite{visser16second}, Smory\'nski's handbook
chapter is the standard account of the derivability
conditions~\cite{smorynski77incompleteness}, and Kurahashi compares the
sufficient sets of conditions~\cite{kurahashi20note}.  Jeroslow showed
that D2 can be dispensed with for Hilbert and Bernays' consistency
statement once D3 holds for every primitive-recursive
term~\cite{jeroslow73redundancies,kurahashi20note}.  The G\"odel
conditions of Section~\ref{sec:locale} omit D2 as well: the proof
uses D2 only in its ex falso form, EXPL, which is a condition.
The G\"odel conditions differ in form from the conditions of these
frameworks.  NEC, D3 and EXPL are rules
on accepted proofs, uniform in the context and in the proof code,
rather than implications over provability formulas; D1 asserts closed
acceptance; and the entailment relation need not obey a deduction
theorem, be recursively enumerable, or come with a faithful checker.
The scope is correspondingly narrower: the conditions fix \rga's term
language and rules (\S\ref{sec:locale-scope}).  Abstract versions of
the second theorem are uniform over systems; Saitou and Noguchi's
mechanized theorem is uniform over a class of classical
theories~\cite{saitou26mechanizing}.  What the G\"odel conditions add
is a uniform theorem over the extensions of \rga's rules on \rga's
language, none of which needs implication introduction at undecided
hypotheses, with two concrete instances whose checkers share one
definition over a table of rule detectors.

\paragraph{Hilbert--Bernays and free-variable arithmetic.}
In the first detailed proof of the second
theorem~\cite{hilbert39grundlagen}, the derivability conditions are
stated on the proof predicate, with the third restricted to
primitive-recursive equations as antecedents~\cite{kurahashi20note};
\rga's decided D3, which \rga proves for every $\varphi$ by a
construction whose instance at the code term $S_G$ is D3G
(Theorem~\ref{thm:g-d3-any}), and also as one sentence over the code
(Theorem~\ref{thm:g-d3-uniform}), has a premise certificate of the
shape of that third condition.  The first two conditions and Hilbert
and Bernays' consistency statement are implications over provability
formulas, and Hilbert and Bernays' proof reasons with such
implications~\cite{kurahashi20note}.  The decided form of the route
used in Section~\ref{sec:goedel} is that of free-variable arithmetic.
In Coquand's formalization, after Guard and Church, self-reference is
a term equation, D3 and ex falso are primitive-recursive proof
transformers with decided antecedents, the G\"odel formula is derived
from the consistency formula by propositional reasoning on decidable
formulas followed by generalization, and the theorem is a
transformation of derivations~\cite{coquand26auto}; the decided case
split of Theorem~\ref{thm:g-split} is that propositional step.  In
those settings every formula used is decided, so the restriction to
decided antecedents costs nothing.  In \rga the restriction is forced,
because every implication of the textbook proof fails while the
corresponding entailments hold (Section~\ref{sec:breaks}), and the
decided forms are those entailments written as sentences.  These
decided forms must be built inside a logic with quantifiers and
undecided sentences, where each term used as a number must be shown
to be one.

Proofs of the second theorem that avoid the textbook route are not
new.  Free-variable arithmetic proves the theorem with decided proof
transformers.  Jeroslow's proof avoids D2 through a fixed point that
asserts its own refutability~\cite{jeroslow73redundancies,kurahashi20note}.
Kritchman and Raz's proof argues from Chaitin's incompleteness theorem
and uses no self-referential sentence~\cite{kritchman10surprise};
Coquand has mechanized Kritchman and Raz's proof for Basic Recursive
Arithmetic~\cite{coquand26chaitin}.  That every implication of the
textbook proof fails in \rga while every entailment holds
(Section~\ref{sec:breaks}) is the material conditional and the
disjunction property at work, not a new route.  What this paper adds
is the decided route carried out inside a gappy arithmetic: the D3G
certificate builder, the cut transformers, and formalized G\"odel~II
in decided-antecedent form with its explicit witness.  Artemov's
schematic consistency, ``$n$ is not a proof of $\z = \tofn 1$'' for
each $n$~\cite{artemov19provability}, is the classical analogue of the
instance-by-instance consistency that \rga derives on its decided
matrix (\S\ref{sec:breaks-g1}).  Heyting arithmetic does not prove its
consistency, by the textbook proof.  Heyting arithmetic does not refute
its own G\"odel sentence either, most simply because it is contained
in Peano arithmetic, which is $\Sigma_1$-sound; inside intuitionistic
arithmetic, Heyting arithmetic's existence property and closure under
Markov's rule together play the role that the witness property plays
in \rga~\cite{troelstra73metamathematical,friedman78classically}.  The
textbook proof goes through in Heyting arithmetic, which derives
$\mathit{Con} \limp G$ for its own consistency sentence and G\"odel
sentence by implication introduction over an antecedent it does not
decide.  Implication introduction over such an antecedent is the step
\rga lacks: \rga derives the entailment $\ConR \vdash G$
(Theorem~\ref{thm:br-hyp}), and formalized G\"odel~II as the
entailment $\ConR \vdash \neg\PrC(\bquo{\ConR})$
(Theorem~\ref{thm:g-formal}), but not the implications.  The two
systems also obtain their closure properties differently.  Heyting
arithmetic's closure properties come from realizability and from
Friedman's
translation~\cite{troelstra73metamathematical,friedman78classically};
\rga's witness property comes from the soundness of its
certificate-carrying semantics.

\paragraph{Self-reference by equation.}
The equation form of self-reference used in \S\ref{sec:goedel-sentences}
is standard: the strong diagonal lemma provides, for each formula
$A(x)$, a closed term $t$ with $\vdash t = \quo{A(t)}$.  The history of
this equation form runs through the Kreisel--Henkin exchange and is
surveyed by Halbach and Visser~\cite{halbach14self} and by Grabmayr
and Visser~\cite{grabmayr20self}; Kripke obtains the same effect with
constants~\cite{kripke23goedel}.  In computability terms the equation
form is Kleene's second recursion theorem~\cite{kleene38notation}, and
$S_G = \Diag(\bquo{D}, \bquo{D})$ is the usual construction of a
quine.  The device is not new.  What is specific to \rga is that the
equation form is forced: Theorem~\ref{thm:br-bicond} shows that the
biconditional form is unprovable for every sentence.

\paragraph{Forms of consequence.}
That a principle can hold as a rule or an entailment while failing as
an implication is familiar.  The admissible rules of Heyting
arithmetic outrun its derivable ones, and its admissible propositional
rules are exactly those of intuitionistic propositional
logic~\cite{visser99rules,visser02substitutions,iemhoff01admissible,rybakov97admissibility}.
In modal logic, necessitation holds for global but not for local
consequence, which is why the deduction theorem needs a restricted
necessitation rule~\cite{blackburn01modal,hakli12deduction}.  \rga's
necessitation for entailments (Theorem~\ref{thm:g-cut}) is of the
global kind.  Kleene's strong three-valued conditional is material,
and $A \limp A$ is not valid in his logic~\cite{kleene52introduction}.
Field adds a non-material conditional to recover $A \limp A$, and does
not declare all of his rules truth-preserving~\cite{field08saving}.
\rga keeps the material conditional, so it proves no implication
between independent sentences, and it recovers the textbook steps as
entailments and, for provability hypotheses, in decided-antecedent
form.  Beall and Murzi's validity Curry~\cite{beall13two} turns on a
validity predicate obeying ``if $A \vdash B$ then $\vdash
\mathrm{Val}(A, B)$'' and ``$A, \mathrm{Val}(A, B) \vdash B$'', and
they blame structural contraction.  Field replies that for validity
read as provability the second principle is local reflection, which
the second theorem already excludes~\cite{field17disarming}.  \rga
sides with Field: the first principle holds for $\PrC$ as
Theorem~\ref{thm:g-cut}, the second fails
(Theorem~\ref{thm:br-loeb-ent}), and contraction is kept.

The closest existing logic of the entailment shape is positive
provability logic: sequents $A \vdash B$ between strictly positive
modal formulas, as in Dashkov's positive fragment of
GLP~\cite{dashkov12positive} and Beklemishev's reflection
calculus~\cite{beklemishev14positive}.  The judgmental reading of
modal logic gives the same split between hypothetical and categorical
judgments: necessitation is a rule on categorical judgments, $\vdash
A$ gives $\vdash \Box A$, not an implication~\cite{pfenning01judgmental},
and \rga's necessitation for entailments (Theorem~\ref{thm:g-cut})
lifts that rule to hypothetical judgments.  The transformer conditions
of Section~\ref{sec:locale} lift it once more, to rules on accepted
proofs in every context.

Artemov's Logic of Proofs~\cite{artemov01explicit} makes the proof
operations of provability explicit.  LP's application operation
corresponds to \rga's decided D2 transformer $m_{\varphi\psi}$
(Theorem~\ref{thm:g-d2}), and LP's proof checker ${!}$, with $t{:}A
\limp {!t}{:}(t{:}A)$, to the certificate builder of D3G
(\S\ref{sec:goedel-d3g}); the witness $\sigma$ is that builder
followed by $\exists$-introduction into $\PrC$.  LP's explicit
reflection $t{:}A \limp A$ is a schema over proof terms, and its \rga
counterpart, $\vdash \Chk(\tofn n, \bquo{A}) \limp A$ for every
numeral $n$ and well-formed $A$, holds, because the checker's verdict
at a numeral is decided and genuine (Theorem~\ref{thm:d-chk-reflect}).
Like LP's explicit reflection, the \rga counterpart is a schema: one
theorem for each numeral and term, proved in HOL by cases on the
decided verdict.  The single sentence over all proofs, $\forall
p.\,(\Chk(p, \bquo{\z = \tofn 1}) \limp \z = \tofn 1)$, fails
(Theorem~\ref{thm:br-loeb}(4)), as implicit reflection does.
Artemov's selector proofs~\cite{artemov19provability,artemov25serial}
let Peano arithmetic prove each instance of its consistency scheme
uniformly without proving $\mathrm{Con}(\mathrm{PA})$.

\paragraph{Provability logic and L\"ob's rule.}
The provability logic GL axiomatizes, for Peano arithmetic, the modal
principles that the provability predicate obeys, with L\"ob's theorem
as its characteristic axiom~\cite{boolos93logic}.  This paper does not
determine the provability logic of $\PrC$.  \rga has L\"ob's rule
(Theorem~\ref{thm:br-loeb}(2)) and formalized L\"ob's theorem as an
entailment (Theorem~\ref{thm:br-loeb-ent}(1)), but not as an
implication.  With the optional D2 condition, the G\"odel conditions
give formalized L\"ob as an entailment (Theorem~\ref{thm:lc-loeb}),
but not L\"ob's rule, whose proof at \rga uses the disjunction
property.  Step-indexed models of
recursive types~\cite{appel01indexed} and Nakano's modality for
recursion~\cite{nakano00modality} underlie logics such as
LSLR~\cite{dreyer11logical} and Iris~\cite{jung18iris}.  These logics
have a ``later'' modality $\triangleright$, which LSLR traces to
G\"odel--L\"ob logic, and a L\"ob rule: from $\triangleright P \limp
P$ infer $P$, valid by induction on the step index.  \rga has no later
modality, and its fuel indexes evaluations, not the truth of
propositions.  \rga's L\"ob rule is the rule of provability logic, and
it follows from G\"odel's second theorem and the disjunction property,
not from induction on fuel.

\paragraph{Reflection principles, progressions and the $\omega$-rule.}
Kreisel and L\'evy's local reflection principle for a theory consists
of the sentences $\PrC(\quo{\varphi}) \limp \varphi$, one for each
sentence~$\varphi$, and Kreisel and L\'evy showed that uniform
reflection over Peano arithmetic is equivalent to transfinite
induction up to $\varepsilon_0$~\cite{kreisel68reflection}.  The local
principle's instance at $\z = \tofn 1$ yields the theory's
consistency, so by the second theorem a consistent theory that meets
the derivability conditions does not prove that instance.  In \rga
even the entailment form of that instance, $\PrC(\bquo{\z = \tofn 1})
\vdash \z = \tofn 1$, is equivalent to $\vdash \ConR$ and fails
(Theorem~\ref{thm:br-loeb-ent}(2)).  Turing's ordinal
logics~\cite{turing39ordinals} and Feferman's transfinite
progressions~\cite{feferman62transfinite} iterate reflection along
ordinal notations, and Beklemishev surveys reflection principles and
their strength~\cite{beklemishev05reflection}.  Shoenfield's recursive
$\omega$-rule~\cite{shoenfield59omega} infers a universal from a
recursive enumeration of proofs of its instances, and iterated, the
rule is complete for true arithmetic.  Adoption (Section~\ref{sec:adopt})
belongs to this family as a reflection rule that strengthens the
system, but it is neither a progression nor an iteration.  Adoption's premise
is a derivation inside the system, of a warrant whose matrix is a
decided verdict of a fixed base checker, and the development defines
no checker for the extended system, so no adoption step concerns the
extended system's own proofs.  Proof assistants meet the same limit.  Where a proof
assistant's logic has been verified in itself or in a close relative,
a consistency proof works in a stronger theory or under an added
assumption, as the second theorem requires: Harrison proved HOL
without the axiom of infinity consistent inside HOL Light, and full
HOL consistent under an additional axiom~\cite{harrison06towards}, and
Fallenstein and Kumar obtained a reflection principle for HOL in HOL4
by assuming a large cardinal~\cite{fallenstein15proof}.  \rgaP and
\igaP are stronger theories of this kind, with the added rule keyed
to the base system's checker.  The claim extensions of
\S\ref{sec:adopt-claims} meet Turing and Feferman in the currency of
notations.  Turing's code $\tau$ (Theorem~\ref{thm:ac-turing}) makes
the consistency extension of \iga a single notation claim, and its
sequence reads \iga's proof checker, as the notations of Turing's
completeness theorem and of Feferman's analysis are built from the
sentences to be proved~\cite[pp.~200--202]{turing39ordinals}, \cite{feferman62transfinite}.
Turing's observation that a progression's theories depend on the
notation and not on the ordinal it denotes~\cite[\S9, pp.~200
and~202--203]{turing39ordinals}, sharpened by
Feferman~\cite{feferman62transfinite} (for the account,
see~\cite[p.~261]{beklemishev05reflection}), is visible at $\tau$
directly: read classically, $\tau$ and the standard code for $\omega$
denote the same ordinal, \iga derives the claim at the standard code,
and the claim at $\tau$ is \iga's consistency.  Kreisel's autonomy
condition~\cite{kreisel58ordinal} adds nothing over bare claims
(Theorem~\ref{thm:ac-autonomy}), and the falsification theorems
(Theorems~\ref{thm:ac-falsified} and~\ref{thm:ac-quine}) test an
adopted claim only through its consequences, never by verifying it.

\paragraph{Non-classical logics and the second theorem.}
Every consistent system known to the author that proves its own
consistency does so in one of two ways.  Either it is arithmetically
weak, as Willard's self-verifying theories, which do not prove
multiplication total~\cite{willard01self},
and Pakhomov's set theory, which does not interpret Robinson
arithmetic~\cite{pakhomov19weak}; or it uses an intensional,
nonstandard provability predicate, as in Feferman's
numerations~\cite{feferman60arithmetization}.  Meyer's relevant
arithmetic $R^\#$ proves that $0 = 1$ is unprovable, but only because
$R^\#$ lacks ex falso, so that the unprovability of $0 = 1$ no longer
expresses negation consistency; Meyer notes that the classical
barriers remain for negation consistency~\cite{meyer21consistency}.
Friedman and Meyer showed that $R^\#$ is properly weaker than Peano
arithmetic: a strictly positive theorem of Peano arithmetic is
unprovable in $R^\#$~\cite{friedman92whither}.  Paraconsistent
arithmetics avoid the theorem by being inconsistent~\cite{priest06in}.
Under weak logics of formal inconsistency, the second theorem needs
extra premises asserting the consistency of particular
formulas~\cite{carnielli20goedel}, the paraconsistent counterpart of
\rga's restriction of discharge to decided hypotheses.  In Kripke's
theory of truth and its axiomatizations KF and PKF, and in Field's
theory, the gaps are confined to the truth predicate: undecided
arithmetic sentences, consistency statements included, still obey
excluded middle, and the classical proof of the second theorem
applies~\cite{kripke75outline,feferman91reflecting,halbach06axiomatizing,field08saving}.
This classical behavior holds even for PKF, whose partial, strong-Kleene logic governs
its whole language, since PKF extends Peano arithmetic and classical
logic still applies to its arithmetic sentences~\cite{halbach18costs}.

\rga takes none of these routes.  \rga proves multiplication, each
compiled primitive-recursive function and Ackermann's function total
by quantified theorems, its consistency statement uses a checker whose
structural laws it proves, it has ex falso, and it is consistent.  Its
non-classical logic changes the form in which the steps of the second
theorem's proof hold, entailments and decided-antecedent forms rather
than implications, but not the theorem.  Nor does \rga's completeness
for grounded truth~\cite{rga-paper} conflict with the theorem: the
G\"odel sentence and $\ConR$ are not grounded-true
(\S\ref{sec:bg-rga}).  Shapiro argued that a dialetheic theory
inherits incompleteness rather than escaping
it~\cite{shapiro02incompleteness}.  Arithmetics with a many-valued
truth predicate keep their arithmetic atoms crisp: Restall studies
arithmetic and truth in \L ukasiewicz's infinitely valued
logic~\cite{restall93arithmetic}, and H\'ajek, Paris and Shepherdson
add a many-valued truth predicate to crisp Peano
arithmetic~\cite{hajek00liar}.  Substructural approaches keep a
transparent truth predicate by restricting a structural rule:
contraction~\cite{zardini11truth} or the transitivity of
consequence~\cite{cobreros13reaching}.  \rga keeps both structural
rules, and its gaps reach the arithmetic sentences themselves; no
earlier proof of the second theorem is known to the author for an
arithmetic of that kind.

\paragraph{The Kleene--Rosser and Curry paradoxes.}
Kleene and Rosser showed, by a formal Richard paradox, that the early
logics of Curry and Church built on combinatory logic and the
$\lambda$-calculus are inconsistent~\cite[pp.~630--631]{kleene35inconsistency},
and Curry found a much simpler argument from other
paradoxes~\cite[p.~115]{curry42inconsistency}, the paradox behind the
restriction of \rga's discharge to decided hypotheses
(\S\ref{sec:recall-forms}).  Czajka
isolates what the Kleene--Rosser argument needs.  Every Hilbert-style
illative system that contains arithmetic, with numerals, equality,
implication introduction at propositions and an induction rule whose
motive is unrestricted, is
inconsistent~\cite[pp.~38--39 and~54]{czajka15semantic}.  The proof
defines an evaluator $T$ with $T\ulcorner X\urcorner = X$ by
conversion and an enumerator $\Theta$ of the codes of the derivable
terms, and it derives, by the unrestricted induction, the single
sentence $\forall x{:}N.\ T(\Theta x)$, that every enumerated term
holds~\cite[pp.~54--58]{czajka15semantic}.  Czajka expects that a
natural-deduction system containing arithmetic in his sense can be
shown consistent~\cite[p.~59]{czajka15semantic}, and he locates the
difference in representation: when implication introduction is
restricted, such a system may be unable to represent its own judgments
with hypotheses as terms, and so to define an enumerator for which
that sentence is derivable~\cite[pp.~59--60]{czajka15semantic}.  \rga and \iga are
consistent natural-deduction systems (\textsf{rga\_consistent},
\textsf{gi\_consistent}; classical) whose implication introduction is
restricted to decided hypotheses and whose induction rule
(\textsf{pr\_nat\_ind}) takes any well-formed motive, with the motive
itself a hypothesis in the step, as in Czajka's
(NInd)~\cite[p.~39]{czajka15semantic}.  That they contain
arithmetic in Czajka's precise sense, up to the coding of numerals and
the restriction of conversion to head steps, is a prediction.  The
argument cannot go through in them, since they are consistent, and
the predicted reason is the two limits this paper meets elsewhere.  The sentence $\forall x{:}N.\ T(\Theta x)$ is a
soundness sentence for all of the system's derivations.  The development derives
disquotation in \rga for quantifier-free formulas~\cite{conrga-paper},
the single soundness sentence over all of its proofs fails
([\citen{conrga-paper}, theorem~\textsf{rga\_no\_uniform\_soundness}]), and by the second theorem \rga
does not even derive its consistency sentence, which that sentence
would give.  Adoption (Section~\ref{sec:adopt}) does reflect \rga's
proofs into \rgaP, but it stratifies the diagonal: every warrant speaks
of \rga's checker, so the enumeration ranges over \rga's proofs, while
the Kleene--Rosser diagonal function would have to be proved total in
\rgaP, and no warrant speaks
of \rgaP's proofs, since the development defines no checker for \rgaP.
The development does not formalize the Kleene--Rosser argument.
Rosser's trick for the first theorem~\cite{rosser36extensions}, cited
in \S\ref{sec:bg}, is a different result of Rosser's and plays no part
here.

\paragraph{Consistency proofs.}
Gentzen proved the consistency of Peano arithmetic by transfinite
induction up to $\varepsilon_0$~\cite{gentzen36widerspruchsfreiheit}:
a reduction on proofs, proved correct, and a proof that the
reduction's descent stops.  \rga proves the shape of that argument,
and no reduction that \rga proves correct has a descent that \rga
proves to stop (Appendix~\ref{app:goedel-cons},
[\citen{conrga-paper}, \S\,``What the theorems show'']).  \iga supplies the descent below
$\varphi(\omega, 0)$, by its transfinite induction along accessible
notations~\cite{iga-paper}, so the Gentzen route to $\ConR$ in \iga
reduces to an ordinal analysis of \rga (Theorem~\ref{thm:op-gentzen}).
For Peano arithmetic the route is carried out in full~\cite{iga-paper}.
\iga derives, from no hypotheses, the consistency sentence of a guarded checker for Peano arithmetic, a derivation in \iga.
The reduction is proved correct in \rga, with an ordinal assignment that \rga proves bounded below $\varepsilon_0$, and \iga's transfinite induction stops the descent.
In HOL's reading the checker accepts a proof of every theorem of Peano arithmetic, so the consistency of \iga implies the consistency of Peano arithmetic, a classical claim that uses no soundness theorem of \iga.
Classically, Peano arithmetic proves the consistency of each fragment
$\mathrm{I}\Sigma_n$~\cite{hajek93metamathematics}; the
internal-reasoning paper proves \rga's counterparts, the consistency of
decided classes of \rga's proofs ([\citen{conrga-paper}, Section ``Consistency, from unconditional to conditional'']).

\section{Conclusion}
\label{sec:concl}

This paper proves G\"odel's incompleteness theorems for \rga, a
paracomplete arithmetic in which every term is a program, and for the
systems built on \rga in the same way.  \rga does not derive its
consistency sentence $\ConR$, and \rga derives its own second theorem
as the entailment $\ConR \vdash \neg\PrC(\bquo{\ConR})$ and as the
sentence $\forall q.\,(\Chk(q, \bquo{\ConR}) \limp \neg\ConR)$, with an
explicit transformer.  The second theorem holds, by one proof, for
every G\"odel extension of \rga: every entailment relation on \rga's
language, closed under \rga's rules, with a checker that meets three
structural conditions, NEC, D3 and EXPL as proof transformers, D1 as
closed acceptance, and consistency.  \rga with its tree checker is one
G\"odel extension, and \iga, with a checker that adds one rule detector
to \rga's table, is another.  Adoption, a reflection rule keyed to a
base system's checker, gives consistent extensions \rgaP and \igaP
that derive the consistency sentences, and the separations it gives
are strict.  A second kind of adoption, a supplied claim that a code
lies in Kleene's $\mathcal O$, meets the second theorem at every code:
no sound claim extension of \iga derives its own consistency sentence,
and the sound claim built from the extension's proof checker derives
that sentence.  Such a claim is tested only through its consequences,
and the absence of a refutation is no rule for adopting it; nor, for
a well-foundedness sentence, is consistency with the system (for a
claim, that consistency is no such rule either is a prediction).

The question behind this work was whether a system like \rga can prove
its own consistency.  The question had a plausible positive answer,
because the textbook proof of G\"odel's second theorem cannot be run in
\rga as written.  The textbook proof's implications are not available.
\rga's conditional is material, so by the disjunction property \rga
proves no implication between independent sentences, and these
implications provably fail, from the diagonal biconditional to
formalized G\"odel~I, and D2 and D3 at $G$; formalized L\"ob fails as
an implication too.  The answer is nevertheless negative.  Read as
entailments, every step of the proof holds, and so does every rule
the proof needs.  The decided proof-predicate route of free-variable
arithmetic is that reading of the proof, with the entailments from
provability hypotheses written as decided sentences, and that route
runs inside \rga.  The entailments that fail, D1 at $G$ and reflection
at $\z = \tofn 1$, are those for which classical theories too have
only rules.  Paracompleteness changes the form in which the proof's
steps hold, and leaves the theorem standing.  Maudlin's worry, that a
paracomplete restriction of this kind eliminates the profound results
of metalogic, does not materialize for \rga, in a precise sense, and
Section~\ref{sec:locale} shows the sense to be uniform: the conditions
that the proof needs mention no implication over an undecided
hypothesis.

How, then, does a grounded theorist come to accept \rga's
consistency?  Not from \rga: the metatheory, HOL, proves that
consistency, and on \rga's own reading $\ConR$ stays ungrounded,
although \rga derives every numeral instance of the universal inside
$\ConR$.  The $\omega$-semantics is the reading on which $\ConR$ is
true, and \rga together with $\ConR$ is sound for it
(Theorem~\ref{thm:g-formal}(4)).  Adoption gives a second answer, by a
rule rather than by a semantics.  \rga itself derives the warrant
that its checker accepts, at every number $v$, a computed proof that
$v$ is not an accepted proof of $\z = \tofn 1$
(Theorem~\ref{thm:ad-warrant}), and adoption turns that warrant into
$\ConR$ (Theorem~\ref{thm:ad-sep}).  What the grounded theorist then
accepts beyond \rga is one rule, sound for the $\omega$-model $\Mw$
(Theorem~\ref{thm:ad-sound}), whose premise is a decided internal
fact about \rga's own checker.

\paragraph{Open questions.}
The questions below remain open.

\medskip\noindent\emph{The second theorem beyond the present scope.}\par\nopagebreak
\begin{itemize}
\item \emph{Axiom extensions and other languages.}  The G\"odel
  conditions answer the question of a uniform second theorem for rule
  extensions of \rga on \rga's language with a checker that meets the
  conditions (\S\ref{sec:locale-scope}), and for one family of
  extensions by an undecided axiom, \iga with a notation claim
  (Theorem~\ref{thm:ac-g2}).  For other extensions by undecided
  axioms, such as \mbox{\rga{}${}+\ConR$}, the development has no
  checker and no instance, and whether such an extension meets the
  conditions, D1 and consistency in particular, is open.  Systems in
  another language, with new term formers or another head step, lie
  outside the conditions altogether.
\item \emph{G\"odel's second theorem for the adoption systems.}  The
  development defines no checker for the derivations of \rgaP or
  \igaP, so G\"odel's second theorem for these systems at their own
  checkers is not stated either way (\S\ref{sec:adopt-con},
  \S\ref{sec:adopt-iga}).  Classically one expects it to hold; that
  expectation is a prediction.  Three relations between the adoption
  systems are not stated either: that \rgaP is contained in \igaP,
  whether \rgaP derives $\ConI$, and whether \iga and \rgaP are
  comparable.  That the adoption systems are proof-grounded, each
  adoption certified by the warrant's derivation, is likewise a
  prediction (\S\ref{sec:adopt-sound}).
\item \emph{More instances.}  \oiga, the union of \oga's rules and
  \iga's, is expected to be a G\"odel extension by the same pattern,
  with a tag for \oga's rule ATI added to \iga's table; that
  expectation is a prediction, and the development builds no such
  checker.  Internal completeness at an
  extended table, the property that the checker accepts an emitted
  proof of each true decided atom, would need per-rule emission laws
  at a general table (\S\ref{sec:ck-embed}), which the development
  does not build.
\end{itemize}

\medskip\noindent\emph{Provability.}\par\nopagebreak
\begin{itemize}
\item \emph{Uniform transformers.}  The transformer conditions are
  per-sentence HOL existentials (\S\ref{sec:locale-scope}).  Whether D2
  and cut have transformers that take the codes as run-time arguments,
  as the uniform D3 transformer does (Table~\ref{tab:g-uniform}), is
  open.
\item \emph{The provability logic of $\PrC$.}  Which modal principles
  $\PrC$ obeys, as implications, as entailments and as rules, is open
  (Section~\ref{sec:related}).
\item \emph{Open entailments about the search.}  Whether the search for
  a proof of $\z = \tofn 1$, if the search stops, yields such a proof
  inside \rga, and whether $\ConR$ entails that the search does not
  stop, are open ([\citen{conrga-paper}, \S\,``What the theorems show'']).
\end{itemize}

\medskip\noindent\emph{Consistency.}\par\nopagebreak
\begin{itemize}
\item \emph{Does \iga derive $\ConR$?}  Section~\ref{sec:open} names one
  obligation per route.  On the Gentzen route, \iga derives $\ConR$ from
  an ordinal analysis of \rga below $\varphi(\omega, 0)$, in the form
  of the premises of Theorem~\ref{thm:op-gentzen}, and no ordinal
  analysis of \rga is known.  On the route through a stage-indexed
  truth, the fold is a theorem in \oiga, in which ATI, inside a transfinite induction below $\varepsilon_0$,
derives the bivalence of the staged truth (Theorem~\ref{thm:op-oiga-con}), and,
  with bivalence of the truth atom taken as a rule instead, at \rga's
  checker and at \iga's own (\S\ref{sec:open-stage}).  For \iga alone,
  the predicted obstacle is the staged bivalence, which \iga, lacking
  ATI, is expected not to derive (Remark~\ref{rem:op-stage}(1)).  \iga
  is proof-grounded and has the disjunction property, and \iga derives
  neither the uniform decidedness of the halting family's sentences
  nor $\LF$ (Theorem~\ref{thm:op-iga-pg}); by the prediction of
  Remark~\ref{rem:op-stage}(1), the staged bivalence would give that
  decidedness.  The step from the staged bivalence to that decidedness
  is the one unbuilt link.  The entailment $\ConR \vdashI \ConI$
  would settle the question negatively, and the development does not
  prove that entailment.
\item \emph{\oga and $\ConR$.}  Whether \oga proves $\ConR$, the
  internal-reasoning paper's open item ([\citen{conrga-paper}, \S\,``Conclusion'']), is
  still open, but narrowed.  \oga derives that $\ConR$ is decided, in
  every context (\textsf{oga\_ConRGA\_decided}; a derivation in \oga):
  ATI decides the universal of the decided instances
  $\neg\Chk(v, \bquo{\z = \tofn 1})$.  \oga does not refute $\ConR$
  (\textsf{oga\_not\_refutes\_ConRGA}; classical, by
  $\omega$-soundness).  Whether \oga derives $\ConR$ itself is the part
  that remains open.  \oiga, which has \oga's rules and \iga's, derives
  $\ConR$ and so strictly exceeds \rga in consistency strength
  (Theorem~\ref{thm:op-oiga-con}); the derivation uses \iga's ground
  induction as well as ATI, so the theorem does not settle the question
  for \oga.  The fold at the truth atom's bivalence rule derives
  $\ConR$ and $\ConI$, and it iterates along the IGA paper's tower of
  constant-free Kripke universes, the systems obtained by closing off
  the constant-free truth term level by level, each deriving the consistency sentence of the one below and, by
  G\"odel's second theorem at its own checker, not its own
  (\S\ref{sec:open-stage}).
  Whether \oiga
  derives $\ConI$ is open: by an unmechanized argument, the staged fold at \iga's own checker
stops at a named obligation, the staged bivalence at all Veblen stages
below $\varphi(\omega,0)$, where \oiga's accessibility ends
  (\S\ref{sec:open-stage}).  The no-go theorem, which the development
  interprets at \oiga (\S\ref{sec:open-nogo}), excludes there only the
  fold with the record motive, not the stage-indexed one.
\item \emph{Adopted claims.}  A code outside $\mathcal O$ with no
  hyperarithmetic descending path, in the manner of Harrison's
  pseudo-well-orderings, gives a claim extension that is not
  $\Mw$-sound, as every code outside $\mathcal O$ does
  (\textsf{adopt\_sound\_iff\_lfp}), and is predicted to give a
  consistent one (\S\ref{sec:adopt-claims}); the development
  constructs no such code and no classical model for its extension.
  Whether consistency with the system is a closure condition on
  adopted claims is likewise not stated; the prediction is that it is
  not, through a Liar code built like $\tau_{T,k}$
  (\S\ref{sec:adopt-ac}).
  Autonomy over bare claims adds nothing
  (Theorem~\ref{thm:ac-autonomy}), so a growing progression of claim
  extensions needs a further step, consistency, reflection or truth;
  no such progression over claims is mechanized.  For a consistent
  extension that is not $\Mw$-sound, only the general form of
  consistency-filtered adoption is refuted
  (Theorem~\ref{thm:ac-rosser}); whether the decided form fails there
  is not established.
\end{itemize}

\appendix

\section{Notation}
\label{app:notation}

Table~\ref{tab:notation} lists the named objects that the paper uses,
grouped by role.  For each object the table gives the symbol, a one-line meaning,
where the paper defines the object, the Isabelle constant, and how the object
stands with respect to the build-time check of
[\citen{conrga-paper}, \S\,``Machines without opaque indices''] (column ``audit'').  In that column, \emph{audited} means
that the object is a root of that check or lies in the definitional
closure of a root (the root is named); \emph{opaque} means that the
object is compiled from an index chosen by Hilbert choice, and the
check requires that the object be found; \emph{HOL} means that the object is a function, relation or
notation of the metalanguage, for which the question does not arise;
``over $\Chk$'' marks a formula with no constant of its own, built from
the audited $\Chk$; and --- marks an \rga notation outside the check's
subject (numerals, code numerals, typing abbreviations, the Liar
and $\Omega$).
Isabelle names are those of the theory text; the theory file of each constant is
in the appendix that the ``defined'' column points to.

{\footnotesize
\setlength{\tabcolsep}{3pt}
\renewcommand{\arraystretch}{1.12}
\begin{longtable}{@{}>{\raggedright\arraybackslash}p{0.15\linewidth}>{\raggedright\arraybackslash}p{0.35\linewidth}>{\raggedright\arraybackslash}p{0.12\linewidth}>{\raggedright\arraybackslash}p{0.20\linewidth}>{\raggedright\arraybackslash}p{0.12\linewidth}@{}}
\caption{The paper's named machines, predicates and classes.}
\label{tab:notation}\\
\hline
symbol & meaning & defined & Isabelle & audit \\
\hline
\endfirsthead
\hline
symbol & meaning & defined & Isabelle & audit \\
\hline
\endhead
\hline
\endfoot
\multicolumn{5}{@{}l}{\textit{Coding and compiling}} \\
$\Diag(x, k)$ & compiled diagonal substitution: the code of $t$ with
  $\tofn k$ for $\vz$, at $x = \bquo{t}$ & [\citen{conrga-paper}, \S\,``The truth predicate and the Liar''],
  \S\ref{sec:goedel-sentences} & \textsf{iDiagSub};
  \textsf{diagSubIdx} & audited (root; also in $\sigma$'s tree) \\
$\mathit{neg}$, $\mathit{or}$, $\mathit{eq}$, $\mathit{ap}$,
  $\mathit{succ}$ & compiled code constructors: the code of a negation,
  disjunction, equation, application or successor from the codes of its
  parts & App.~\ref{app:goedel} & \textsf{mkNegS}, \textsf{mkOrS},
  \textsf{mkEqS}, \textsf{mkAppS}, \textsf{mkSucS} & audited \\
$\Inst_u(v)$ & the code of $u$ applied to the numeral of $v$, a
  compiled function of $v$ & Def.~\ref{def:ad-rule} & \textsf{adrInst} &
  not checked \\
\multicolumn{5}{@{}l}{\textit{Semantics and deciders}} \\
$\Mw$ & the $\omega$-model: a universal is true when every numeral
  instance is & \S\ref{sec:adopt-sound} & \textsf{omv} & HOL \\
proof-grounded, $\omega$-grounded & a system is proof-grounded when it
  is sound and complete for grounded truth read with its own proofs as
  the certificates at the universals, and $\omega$-grounded when it is
  sound for $\Mw$ but has no such semantics.  Neither term is a
  viewpoint label: the three viewpoints are classical, a derivation in
  the named system (``a derivation in \rga''), and prediction, and the
  bare word ``grounded'' is not a label, since it names grounded truth;
  the two terms classify a system by its semantics &
  \S\ref{sec:open-stage} & \textsf{gsem}, \textsf{gsem\_sound},
  \textsf{gs\_complete\_open} & HOL \\
$\mathit{TS}(\beta, c, r)$, $F(\delta, \alpha)$ & the staged truth: a
  copy of \iga's truth term with a stage $\beta$, a Cantor-normal-form
  code below $\varepsilon_0$, at every point; the term with code $c$
  has the value $r$ at stage $\beta$.  The staged bivalence is its
  decidedness at every such stage; $F(\delta, \alpha) = \omega^{e(\alpha
  \sqcup \omega) \oplus \delta}$ is the stage family of the node truth
  & \S\ref{sec:open-stage} & \textsf{trmTS}, \textsf{tsBivF},
  \textsf{stFam}, \textsf{NdTG} & not checked \\
$\mathit{DecVal}(c, x)$ & the value that the decider reads at $c$ is
  $x$: $\exists l.\ \Dec(l, c) = \suc x$ & Table~\ref{tab:d-search} &
  \textsf{iPrDecVal} & not checked \\
\multicolumn{5}{@{}l}{\textit{Tables and checkers}} \\
$\mathcal D$, $\mathcal D_L$ & a table of rule detectors, an index
  tree; the table of a list $L$, the nested disjunction of its members &
  Def.~\ref{def:ck-table} & \textsf{Dt}; \textsf{drvChkIdxG} & not
  checked \\
$R$ & a step relation of HOL that a table decides at genuine codes
  (written $D$ in the Isabelle text) & \S\ref{sec:ck-faith} & --- &
  HOL \\
$\Chk_{\mathcal D}$, $\pi_{\mathcal D}$, $\mathit{TreeValid}_{\mathcal
  D}$, $\Cert_{\mathcal D}$ & the checker, the proof test, the validity
  test and the certificate statement at the table $\mathcal D$ &
  Def.~\ref{def:ck-table} & \textsf{iChkD}, \textsf{provesChkD},
  \textsf{iTreeValidD}, \textsf{iHasCertD} & not checked \\
$L_{\rga}$, $L_{\iga}$ & \rga's list of 41 detectors; \iga's list,
  $L_{\rga}$ with $E_{\mathrm{gi}}$ appended & Thm.~\ref{thm:ck-bridge},
  Def.~\ref{def:ck-iga} & \textsf{drvChkList}, \textsf{drvChkListI} &
  not checked \\
$E_{\mathrm{gi}}$ & the ground-induction tag, \rga's three-premise
  citation former at a table-blind payload & Def.~\ref{def:ck-iga} &
  \textsf{chkGfixIdx} & not checked \\
$\ChkI$ & \iga's checker, $\Chk_{\mathcal D_{L_{\iga}}}$ &
  Def.~\ref{def:ck-iga} & \textsf{iChkI} & not checked \\
$\vdashI$, $\vdashRP$, $\vdashIP$ & the entailments of \iga, \rgaP
  and \igaP & \S\ref{sec:checker}, Def.~\ref{def:ad-rule},
  \S\ref{sec:adopt-iga} & \textsf{gi.entails}, \textsf{adr\_ent},
  \textsf{adi\_ent} & HOL \\
$\gchk$ & the checker parameter of the G\"odel conditions, a HOL
  function on terms & \S\ref{sec:locale-cond} & \textsf{chk} & HOL \\
\multicolumn{5}{@{}l}{\textit{Provability}} \\
$\mathit{Prf}$, $\PrC$, $\ConR_T$ & classical: proof predicate,
  $\PrC(x) \ldef \exists q.\ \mathit{Prf}(q, x)$, consistency of $T$ &
  \S\ref{sec:bg-goedel} & --- & HOL \\
$\gPr$, $\gCon$ & the provability formula and the consistency sentence
  of the checker $\gchk$ & Def.~\ref{def:lc-sentences} &
  \textsf{iPrcT}, \textsf{ConT} & over $\gchk$ \\
$\PrI$, $\ConI$ & \iga's provability formula, $\neg\forall q.\
  \neg\ChkI(q, x)$, and consistency sentence, $\forall p.\
  \neg\ChkI(p, \bquo{\z = \tofn 1})$ & \S\ref{sec:ig-theorems},
  Def.~\ref{def:ig-con} & no constant; \textsf{ConIGA} & over $\ChkI$ \\
\multicolumn{5}{@{}l}{\textit{Self-reference and G\"odel's theorems}} \\
$D_G$ & the G\"odel template, $\forall q.\
  \neg\Chk(q, \Diag(x, x))$ with $x$ the free variable $\vz$ &
  \S\ref{sec:goedel-sentences} & \textsf{goedelD} & audited (in
  $\sigma$'s tree) \\
$S_G$ & $\Diag(\bquo{D_G}, \bquo{D_G})$, the code term of $G$; $\vdash S_G
  = \bquo{G}$ & \S\ref{sec:goedel-sentences} & \textsf{goedelS} &
  audited (in $\sigma$'s tree) \\
$G$ & the G\"odel sentence, $\forall q.\ \neg\Chk(q, S_G)$ &
  \S\ref{sec:goedel-sentences} & \textsf{goedelG} & audited (in
  $\sigma$'s tree) \\
D1, D2, D3 & the classical derivability conditions &
  \S\ref{sec:bg-goedel} & --- & HOL \\
$e$ & the EXPL transformer: accepted proofs $a$ of $g$ and $b$ of $\neg
  g$ give the accepted proof $e\cdot a\cdot b$ of $\z = \tofn 1$ &
  Thm.~\ref{thm:g-expl} & \textsf{goedelE}; \textsf{goedelEIdx} & audited
  (the tree) \\
(no symbol) & the ex-falso transformer: an accepted proof of $\z =
  \tofn 1$ gives an accepted proof of any well-formed $\varphi$ &
  \S\ref{sec:breaks-dc} & \textsf{exfIdx} & audited \\
$\sigma$, $\tau_\sigma$ & the D3G witness, $\sigma = \lambda v.\
  \cmp(\tau_\sigma)\cdot v\cdot v$, and its explicit tree &
  Thm.~\ref{thm:g-d3g} & \textsf{goedelSigma}; \textsf{goedelSigmaIdx}
  & audited (the witness) \\
$\sigma_\varphi$ & the transformer of the decided D3 for $\varphi$,
  $\lambda v.\ \cmp(\tau_\varphi)\cdot v\cdot v$ &
  Thm.~\ref{thm:g-d3-any} & \textsf{d3SigmaOf} & audited \\
$\sigma_{\mathrm u}$, $\mathit{pr}(c)$ & the transformer of the decided
  D3 as one sentence over the code; the code of $\PrC(x)$ with the
  numeral of $c$ for $x$ & Thm.~\ref{thm:g-d3-uniform} &
  \textsf{uSigIdx}; \textsf{prCU} & audited \\
$c_Q$ & the cut transformer of a proof $Q$ of $\{A\} \vdash B$ (or
  $\{A_1, A_2\} \vdash B$) & Thm.~\ref{thm:g-cut} & \textsf{cut1M},
  \textsf{cut1Idx}; \textsf{cut2M}, \textsf{cut2Idx} & audited \\
$T_Q$, $T_1$ & the transformers of the G\"odel sentences in
  decided-antecedent form: from
  proofs of $\ConR$ (built from a proof $Q$ of $\{\ConR\} \vdash G$),
  resp.\ of $G$, to proofs of $\z = \tofn 1$ & Thm.~\ref{thm:g-sentence}
  & \textsf{goedel2T}, \textsf{goedel2Idx}; \textsf{goedel1T},
  \textsf{goedel1Idx} & audited \\
$S_\psi$, $L_\psi$ & the L\"ob diagonal for $\psi$: $L_\psi =
  \PrC(S_\psi) \limp \psi$, with $\vdash S_\psi = \bquo{L_\psi}$ &
  \S\ref{sec:goedel-formal} & \textsf{loebS}, \textsf{loebL} & over $\Chk$
  and $\Diag$ \\
EXPL, D3G & the decided forms of D2 and D3, as properties of $e$ and
  $\sigma$: closed, provably total, and an internal sentence &
  \S\ref{sec:goedel-split} & \textsf{EXPL}, \textsf{D3G};
  \textsf{EXPLF}, \textsf{D3GF} & HOL \\
$\gD$, $\gS$, $\gG$ & the G\"odel template, code term and sentence
  at the checker $\gchk$ & Def.~\ref{def:lc-sentences} &
  \textsf{goedelDT}, \textsf{goedelST}, \textsf{goedelGT} & over
  $\gchk$ \\
$G_{\mathrm I}$ & \iga's G\"odel sentence, $G$ with $\ChkI$ for $\Chk$ &
  \S\ref{sec:ig-theorems} & \textsf{goedelGT} at \textsf{iChkI} &
  over $\ChkI$ \\
\multicolumn{5}{@{}l}{\textit{Adoption}} \\
\rgaP, \igaP & \rga with \adoptR; \iga with \adoptI &
  \S\ref{sec:adopt-rule}, \S\ref{sec:adopt-iga} & \textsf{adr\_ent},
  \textsf{adi\_ent} & HOL \\
$\WarR(F, u)$, $\WarI(F, u)$ & the warrant of $u$ by the index tree
  $F$, at $\Chk$ and at $\ChkI$ & Def.~\ref{def:ad-rule},
  \S\ref{sec:adopt-iga} & \textsf{adrW}, \textsf{adiW} & over $\Chk$;
  over $\ChkI$ \\
\multicolumn{5}{@{}l}{\textit{Adopted notation claims}} \\
$\KO$, $\KO(\tofn a)$, $a \in \KO$ & Kleene's $\mathcal O$ as the IGA
  paper's positive fixed point; the notation claim at the code $a$; the
  HOL statement that $a$ lies in the least fixed point &
  \S\ref{sec:adopt-claims} & \textsf{trmO}, \textsf{oMem}; \textsf{lfp}
  (\textsf{gfixOp} \textsf{oPhi}) & not checked \\
$\Ia{a}$, $\vdashA{a}$ & the claim extension, \iga with $\KO(\tofn a)$
  as a hypothesis in every context, and its entailment &
  Def.~\ref{def:ac-ext} & \textsf{gi.adoptDer} & HOL \\
$\ChkA{a}$, $\PrA{a}$, $\ConA{a}$ & the checker of $\Ia{a}$, at
  $L_{\iga}$ with the axiom cell $E_a$ appended; its provability
  formula and consistency sentence & Def.~\ref{def:ac-chk} &
  \textsf{iChkF}, \textsf{chkF}, \textsf{axIdx}; no constant;
  \textsf{ConF} & not checked; over $\ChkA{a}$ \\
$\tau_{T,k}$, $u_{T,k}$ & the limit code that climbs through the
  finite notations while $\cmp(T)\cdot\tofn k\cdot\tofn n$ has the
  value $0$; the universal's body $\lambda v.\,(\cmp(T)\cdot\tofn
  k\cdot v = \z)$ & Def.~\ref{def:ac-turing} & \textsf{pwcT},
  \textsf{pwF} & not checked \\
$\tau$, $\tau_a$ & Turing's code, $\tau_{T,k}$ at \iga's proof test of
  $\z = \tofn 1$; the same at the proof test of $\ChkA{a}$ &
  Def.~\ref{def:ac-turing}, Thm.~\ref{thm:ac-larger} & \textsf{tc},
  \textsf{tcF} & not checked \\
$\mathit{qn}$ & the quine code, a well-shaped limit code whose
  sequence has the value $\mathit{qn}$ at every argument &
  Thm.~\ref{thm:ac-quine} & \textsf{fuQ} & not checked \\
$\WF(R)$, $\Acc_R$ & the sentence that the order $R$ is well-founded,
  $\forall x.\,(x\jnat \limp \Acc_R\cdot x)$; the accessible part of
  $R$, a positive fixed point & \S\ref{sec:adopt-ac} &
  \textsf{trmWF} at \textsf{trmNatD}; \textsf{trmAcc} & not checked \\
$R_{\mathrm L}$, $R_{\mathrm{Ro}}$ & the Liar order and the grounded
  Rosser order at a table $\mathcal D$ & Thm.~\ref{thm:ac-liar},
  Thm.~\ref{thm:ac-rosser} & \textsf{fuRL}, \textsf{fuRLR} & not
  checked \\
\end{longtable}
}

\section{G\"odel II in full}
\label{app:goedel}
\begingroup\sloppy

This appendix proves Theorem~\ref{thm:g-g2}, $\emptyset \nvdash \ConR$,
along the dependency spine of its Isabelle endpoint
\textsf{rga\_goedel2}.  Each step is stated at the granularity of the
mechanized lemma that carries it, with every hypothesis, so that a reader
can follow the argument or re-mechanize the key steps.  Each technique
is exhibited once, repetitions are compressed into tables, and facts
below the spine are cited by theory file (relative to \texttt{hol/}).
Unless a theory is named, a name is in one of the two G\"odel theories
(\texttt{RGA\_Goedel/Goedel\_II.thy}, or, for the primitive-recursive
witness, the headline and the consequences of G\"odel's second theorem
(\S\ref{app:goedel-cons}), \texttt{RGA\_Goedel/Goedel\_II\_PR.thy}) or
is a basic \rga rule or kit lemma (under \texttt{RGA/},
\texttt{RGA\_Kit/} and \texttt{RGA\_Checker/}).  The spine is short.
\begin{itemize}\raggedright
\item \textsf{rga\_goedel2} $\leftarrow$\allowbreak\
  \textsf{rga\_goedel2\_of\_d3gcert}\,[\textsf{D3Gcert\_holds}]
  $\leftarrow$\allowbreak\ \textsf{rga\_goedel2\_of\_d3g} $\leftarrow$\allowbreak\
  \textsf{rga\_goedel2\_of\_expl\_d3g}\,[\textsf{EXPL\_goedelE}];
\item \textsf{rga\_goedel2\_of\_expl\_d3g} $\leftarrow$
  \textsf{goedelG\_of\_expl\_d3g\_con}, \textsf{goedelG\_refutes},
  \textsf{rga\_no\_proof\_0eq1}.
\end{itemize}
\S\ref{app:goedel-ent} proves the entailment results of
Sections~\ref{sec:breaks} and~\ref{sec:goedel}: the $\Sigma$
deduction theorem and its corollaries, the G\"odel ingredients as
entailments, the cut transformers, and formalized G\"odel~II and
L\"ob.  Every negative theorem below (G\"odel~I and~II, the
independence of $G$, the failing entailments, and the consequences)
passes through a soundness theorem.  Most of these theorems pass
through the companion paper's grounded soundness theorem, proved in
HOL~\cite{rga-paper}, by way of that theorem's corollary
\textsf{rga\_consistent}.  \textsf{ConRGA\_unrefutable} and the
theorems built on \textsf{ConRGA\_unrefutable} pass through
\textsf{rga\_consistent\_via\_wtrue}, the same consistency fact reached
through the truth predicate.  The $\omega$-model arguments of
\S\ref{app:goedel-ent} pass through \rga's soundness for the
$\omega$-semantics (\textsf{rga\_omega\_sound}).  That is the sense in
which the negative theorems have no hypothesis.  The entailments that hold use
no soundness theorem.

\paragraph{Notation at the level of codes.}
Besides the notation of [\citen{conrga-paper}, Section ``RGA inside RGA''], this appendix names
the following compiled functions.  Each function is a total \rga term built from
compiled index trees.
\begin{itemize}
\item \emph{Code constructors.}  $\mathit{neg}(c)$, $\mathit{or}(a,b)$,
  $\mathit{eq}(a,b)$, $\mathit{ap}(a,b)$ and $\mathit{succ}(a)$ build
  the code of a negation, a disjunction, an equation, an application and
  a successor from the codes of the parts (Isabelle \textsf{mkNegS},
  \textsf{mkOrS}, \textsf{mkEqS}, \textsf{mkAppS}, \textsf{mkSucS}).  At code numerals \rga proves, for example,
  $\emptyset \vdash \mathit{neg}(\bquo{t}) = \bquo{\neg t}$.  The code of
  an applied recursor is
  $\mathit{rap}(g,h,s) = \mathit{ap}(\mathit{ap}(\mathit{ap}(\bquo{\cR},g),h),s)$.
\item \emph{Numeral codes.}  $\Num\cdot k$ is the code of the numeral of
  $k$ (Isabelle \textsf{numC}~$k$, which is $\Num\cdot k$ after unfolding
the definition of $\Num$).
  At a numeral, $\emptyset \vdash \Num\cdot\tofn n = \bquo{\tofn n}$
  ([\citen{conrga-paper}, \S\,``Closed-numeral representability'']).
\item \emph{Proof lists.}  The empty list has code $\z$, and
  $\mathit{cons}(h,T) = \suc(\mathit{pair}(h,T))$ puts the entry $h$ in
  front of the list $T$, where $\mathit{pair}$ is the compiled Cantor
  pair (\textsf{lstCons}, \textsf{cpPair}).  $\mathit{app}(A,B)$ appends
  two lists (\textsf{appC}), $|T|$ is the length of $T$, $\mathit{nd}(T,i)$ is the
  entry at position $i$ (position $0$ is the front), and
  $\mathit{head}(T)$ is the front entry.  The code of the judgment
  $\emptyset \vdash t$ is
  $\mathit{jdg}(\bquo t) = \mathit{pair}(\bquo{\emptyset}, \bquo t)$
  (\textsf{emptyJdg}).
\item \emph{Images.}  An \emph{environment} $z$ is a coded list of
  pairs of a variable and a code.  For a fixed HOL term $t$ (a
  \emph{template}), $\mathit{img}_z(t) = \Subst(z, \bquo t)$ is the code
  of $t$ with $z$ substituted ([\citen{conrga-paper}, \S\,``Soundness, proof by proof'']).
  $\zeta_v$ is the one-entry environment that sends $\vz$ to $\Num\cdot v$
  (\textsf{envVar0}), and $\mathit{img}_v(t) = \mathit{img}_{\zeta_v}(t)$.
\item \emph{Head reduction.}  $\mathit{red}(c)$ is the code of the head
  reduct of the term coded by $c$, computed by the compiled head reducer
  (\textsf{gHR}).
\end{itemize}

\subsection{The sentences and the diagonal}
\label{app:goedel-diag}

The diagonal function substitutes a numeral for the variable $\vz$.  The
diagonal function's index tree is explicit.  The index \textsf{diagListIdx} maps $y$ to the
code of the one-entry substitution list that sends $\vz$ to $\tofn y$,
using the explicit numeral-code tree of $\Num$.  The index
\textsf{diagSubIdx} composes \textsf{diagListIdx} with the coded simultaneous substitution
\textsf{submMIdx} of $\Subst$ at depth $0$.

\begin{defi}[\textsf{diagSubIdx}, \textsf{iDiagSub}]
\label{def:d-diag}
$\Diag(x,y) \ldef \cmp(\mathit{diagSubIdx})\cdot x\cdot y$.  For a
well-formed $t$, the denotation of $\mathit{diagSubIdx}$ sends $(\quo{t}, k)$ to
$\quo{\ssub{\vz}{\tofn k}{t}}$ (\textsf{pr2fi\_diagSubIdx},
\textsf{submMX\_correct}).
\end{defi}

\begin{thm}[\textsf{iDiagSub\_closed}]
\label{thm:d-diag-closed}
For every well-formed term $t$ and every natural $k$,
$\emptyset \vdash \Diag(\bquo{t}, \tofn k) = \bquo{\ssub{\vz}{\tofn k}{t}}$.
\end{thm}

\noindent
\emph{Proof.}  The denotation computes the substitution
(\textsf{submMX\_correct}, \textsf{fmm\_0}, \textsf{trm\_sub1\_eq\_subm}),
and [\citen{conrga-paper}, theorem~\textsf{cmp\_eq\_num}] turns the HOL value into a derivation.
\qed

\begin{defi}[\textsf{goedelD}, \textsf{goedelS}, \textsf{goedelG}, \textsf{goedelF}, \textsf{ConRGA}]
\label{def:d-sentences}
With $x$ the free variable $\vz$ (which is $v_1$ under the binder),
\begin{gather*}
  D_G \ldef \forall q.\ \neg\Chk(q, \Diag(x,x)), \qquad
  S_G \ldef \Diag(\bquo{D_G}, \bquo{D_G}), \\
  F \ldef \lambda q.\ \neg\Chk(q, S_G), \qquad
  G \ldef \cA\cdot F = \forall q.\ \neg\Chk(q, S_G),
\end{gather*}
and $\ConR \ldef \forall p.\ \neg\Chk(p, \bquo{\z = \tofn 1})$.  The
term $S_G$ is closed and well-formed.
\end{defi}

\begin{thm}[\textsf{goedel\_diag}, \textsf{goedelS\_eq}]
\label{thm:d-goedelS}
$\ssub{x}{\bquo{D_G}}{D_G} = G$, an identity of HOL terms; and
$\emptyset \vdash S_G = \bquo{G}$.
\end{thm}

\noindent
\emph{Proof.}  The substitution passes the binder with a shift
(\textsf{subf\_lift}) and leaves closed subterms unchanged: the compiled
constant $\cmp(\mathit{diagSubIdx})$, the checker's components and the
numerals.  Only the two occurrences of $x$ change, and they become
$\bquo{D_G}$, which gives $G$ literally.  For the equation, apply
Theorem~\ref{thm:d-diag-closed} at $t = D_G$ (well-formed,
\textsf{wft\_goedelD}) and $k = \quo{D_G}$; the right side of that
equation is $\bquo{G}$ by the identity.  \qed

\medskip
The equation $S_G = \bquo{G}$ moves the checker's verdict between the
two codes.  The following transport lemmas do this in either argument
of $\Chk$.  The proof of these lemmas is the standard technique of the development
for moving a judgment along a provable equation.

\begin{thm}[\textsf{iChk\_conc\_transport}, \textsf{iChk\_proof\_transport}]
\label{thm:d-transport}
For every context $\Gamma$ and all terms $p, p', c, c'$: if
$\Gamma \vdash c = c'$ and $\Gamma \vdash \Chk(p,c)$, then
$\Gamma \vdash \Chk(p,c')$; and if $\Gamma \vdash p = p'$ and
$\Gamma \vdash \Chk(p,c)$, then $\Gamma \vdash \Chk(p',c)$.
\end{thm}

\noindent
\emph{Proof.}  Pick a variable $v_i$ fresh for $p$.  \rga's
equality-substitution rule derives $\ssub{v_i}{b}{\varphi}$ from
$a = b$ and $\ssub{v_i}{a}{\varphi}$ for a well-formed $\varphi$.  Take
$\varphi = \Chk(p, v_i)$.  The checker commutes with substitution
(\textsf{subf\_iChk}), and $\varphi$ is well-formed because the
derivation of $\Chk(p,c)$ makes $p$ well-formed (\textsf{entails\_wff},
\textsf{wft\_iChk}).  The proof slot is the same with $\varphi =
\Chk(v_i, c)$.  \qed

\subsection{The checker facts used}
\label{app:goedel-chk}

The proof uses the tree checker of [\citen{conrga-paper}, definition~\textsf{iChk}], its
faithfulness and genuineness ([\citen{conrga-paper}, theorem~\textsf{iChk\_faithful}],
\textsf{iChk\_faithful}, \textsf{iChk\_genuine\_proves}) and its
decidedness ([\citen{conrga-paper}, theorem~\textsf{hq\_iChk}], \textsf{hq\_iChk}), all
proved in \texttt{RGA\_Checker/Chk\_Faithful.thy}
([\citen{conrga-paper}, Appendix ``Coding and machines'']), and the list laws below.  Recall that
$\Chk(p,c)$ is the conjunction of $\mathit{TreeValid}(p)$, $p \ne \z$,
``the head of $p$ has no hypotheses'' and ``the head's conclusion is
$c$''.

\begin{thm}[\textsf{iChk\_headI}, \textsf{iChk\_headD}]
\label{thm:d-head}
Let $\Gamma \vdash P \jnat$.  If $\Gamma$ derives $c\jnat$, $P \ne \z$,
$\mathit{TreeValid}(P)$ and $\mathit{head}(P) = \mathit{jdg}(c)$, then
$\Gamma \vdash \Chk(P,c)$.  Conversely, if $\Gamma \vdash \Chk(P,c)$,
then $\Gamma \vdash \mathit{head}(P) = \mathit{jdg}(c)$.
\end{thm}

\begin{thm}[\textsf{iChk\_appC}]
\label{thm:d-appC}
If $\Gamma$ derives $A\jnat$, $B\jnat$, $\Chk(A,a)$ and $\Chk(B,b)$,
then $\Gamma$ derives $\mathit{TreeValid}(\mathit{app}(A,B))$,
$\Chk(\mathit{app}(A,B), a)$, and
$\mathit{nd}(\mathit{app}(A,B), |A|) = \mathit{jdg}(b)$ with
$|A| < |\mathit{app}(A,B)|$.
\end{thm}

\noindent
So in an appended pair of accepted proofs, $A$'s head sits at position
$0$ and $B$'s head at position $|A|$, both with the empty context.  The
validity part is \textsf{iTreeValid\_appC}
(\texttt{RGA\_Checker/Chk\_Graft.thy}).

A new entry is accepted on top of an accepted list when its rule
detector fires.  Write $\mathit{fires}_E(k)$ for the test that detector
$E$ accepts the candidate $k = \mathit{pair}(\mathit{jdg}(y), T)$: the
entry $\mathit{jdg}(y)$ placed on the list $T$ (\textsf{ndFire}).

\begin{thm}[\textsf{iChk\_lstCons}]
\label{thm:d-cons}
Let $E$ be one of the 41 detectors.  If $\Gamma$ derives $T\jnat$,
$\mathit{TreeValid}(T)$, $y\jnat$, $\IsCode(y) \ne \z$ and
$\mathit{fires}_E(\mathit{pair}(\mathit{jdg}(y), T))$, then
$\Gamma \vdash \Chk(\mathit{cons}(\mathit{jdg}(y), T), y)$.
\end{thm}

\noindent
The detectors are fired through introduction lemmas in
\texttt{RGA\_Checker/Chk\_Emission.thy}.  The proof of
Theorem~\ref{thm:g-g2} uses two of these lemmas.  Each lemma
takes as premises the readings of the cited positions $z_1, z_2$ of
$T$: both positions are below $|T|$, and both entries have the same hypothesis
set as the new entry.
\begin{itemize}
\item \textsf{ndFire\_negE\_intro} (rule negE: from $\neg a$ and $a$,
  infer any well-formed $c$).  The conclusion at $z_1$ is
  $\mathit{neg}$ of the conclusion at $z_2$, and $\Wfc(y) \ne \z$.
\item \textsf{ndFire\_notAllI\_intro} (rule notAllI: from $\neg(u\cdot
  a)$ and $a\jnat$, infer $\neg(\cA\cdot u)$).  The conclusion at $z_1$
  has the form $\mathit{neg}(\mathit{ap}(u', a'))$; the conclusion at
  $z_2$ is an equation whose two sides are equal and equal to $a'$; and
  $y = \mathit{neg}(\mathit{ap}(\bquo{\cA}, u'))$.
\end{itemize}

\subsection{EXPL: explosion as a proof transformer}
\label{app:goedel-expl}

\begin{defi}[\textsf{goedelEIdx}, \textsf{goedelE}, \textsf{EXPLF}, \textsf{EXPL}]
\label{def:d-expl}
The index $\mathit{goedelEIdx}$ is built from the explicit trees of
pairing, numerals, list append and list length, and computes
\[
  (a,b) \mapsto \mathit{cons}(\mathit{jdg}(\bquo{\z = \tofn 1}),\
  \mathit{app}(a,b)) ;
\]
$e \ldef \cmp(\mathit{goedelEIdx})$.  For a
term $e$,
\[
  \mathrm{EXPLF}(e) \ldef \forall a.\ \forall b.\
  \bigl(\Chk(a,\bquo{G}) \wedge \Chk(b,\bquo{\neg G})\bigr) \limp
  \Chk(e\cdot a\cdot b, \bquo{\z = \tofn 1}) ,
\]
and $\mathrm{EXPL}(e)$ is the HOL proposition that $e$ is closed,
$\emptyset \vdash \forall a.\ \forall b.\ (e\cdot a\cdot b)\jnat$, and
$\emptyset \vdash \mathrm{EXPLF}(e)$.
\end{defi}

\noindent
The totality conjunct is needed because \rga's $\forall$-elimination
asks for a $\jnat$ premise at the instance.

\begin{thm}[\textsf{goedelE\_eq}]
\label{thm:d-goedelE-eq}
If $\Gamma \vdash a\jnat$ and $\Gamma \vdash b\jnat$, then
$\Gamma \vdash e\cdot a\cdot b = \mathit{cons}(\mathit{jdg}(\bquo{\z =
\tofn 1}), \mathit{app}(a,b))$.
\end{thm}

\noindent
\emph{Proof.}  Unfold the index tree one node at a time with the
evaluation laws of the composition node and the initial functions
(\textsf{ir\_prSubI}, \textsf{ir\_pairI}, \textsf{ir\_sucI},
\textsf{ir\_numI}, \textsf{cmpi\_L\_eq}, \textsf{cmpi\_R\_eq}).  Each law
holds at abstract arguments because the tree is explicit
([\citen{conrga-paper}, \S\,``Compiled primitive-recursive machines'']).  \qed

\begin{thm}[\textsf{goedelE\_chk}]
\label{thm:d-goedelE-chk}
For every context $\Gamma$, every term $g$ and all terms $a$, $b$: if
$\Gamma$ derives $a\jnat$, $b\jnat$, $\Chk(a,\bquo{g})$ and
$\Chk(b,\bquo{\neg g})$, then $\Gamma \vdash \Chk(e\cdot a\cdot b,
\bquo{\z = \tofn 1})$.
\end{thm}

\noindent
\emph{Proof.}  Theorem~\ref{thm:d-appC} gives the accepted list
$\mathit{app}(a,b)$ with $\mathit{jdg}(\bquo g)$ at position $0$ and
$\mathit{jdg}(\bquo{\neg g})$ at position $|a|$.  Take the candidate
$k = \mathit{pair}(\mathit{jdg}(\bquo{\z = \tofn 1}),
\mathit{app}(a,b))$.  The conclusion at $|a|$ is $\bquo{\neg g}$, which
\rga proves equal to $\mathit{neg}(\bquo{g})$, the negation of the
conclusion at $0$.  Both entries have the empty context, and
$\bquo{\z = \tofn 1}$ is a well-formed code.  So negE fires with $z_1 =
|a|$ and $z_2 = 0$ (\textsf{ndFire\_negE\_intro}), and
Theorem~\ref{thm:d-cons} accepts
$\mathit{cons}(\mathit{jdg}(\bquo{\z = \tofn 1}), \mathit{app}(a,b))$.
Transport along Theorem~\ref{thm:d-goedelE-eq} in the proof slot
(Theorem~\ref{thm:d-transport}).  \qed

\begin{thm}[\textsf{EXPL\_antecedent\_B}, \textsf{EXPL\_goedelE}]
\label{thm:d-expl}
If $\Gamma \vdash a\jnat$ and $\Gamma \vdash b\jnat$, then
$\Gamma \vdash (\Chk(a,\bquo{G}) \wedge \Chk(b,\bquo{\neg G}))\jbool$.
And $\mathrm{EXPL}(e)$ holds.
\end{thm}

\noindent
\emph{Proof.}  The first part is [\citen{conrga-paper}, theorem~\textsf{hq\_iChk}] twice and
the conjunction of decided terms.  For the second part, $e$ is closed, and
totality is $\cmp$-totality (\textsf{cmp\_total\_app}) under two
$\forall$-introductions (\textsf{goedelE\_total}).  For
$\mathrm{EXPLF}(e)$ (\textsf{goedelE\_expl}), work in the context
$\{v_0\jnat, v_1\jnat\}$.  Under the extra hypothesis
$\Chk(v_1,\bquo G) \wedge \Chk(v_0,\bquo{\neg G})$, the two conjuncts and
Theorem~\ref{thm:d-goedelE-chk} at $g = G$ give $\Chk(e\cdot v_1\cdot
v_0, \bquo{\z = \tofn 1})$.  The antecedent is decided by the first
part, so the hypothesis may be discharged into an implication
(\textsf{ir\_impI\_bool}).  Two $\forall$-introductions close the
sentence.  \qed

\noindent
The same proof, with an arbitrary term $g$ in place of $G$, gives the
internal sentence for every $g$ (\textsf{goedelE\_expl\_any}), with the
decided antecedent from [\citen{conrga-paper}, theorem~\textsf{hq\_iChk}] twice.
Theorem~\ref{thm:d-goedelE-chk} is generic in $g$: that theorem is the
checker's internal law that explosion composes proofs.  No numerals $a$,
$b$ satisfy both checker hypotheses of that theorem in the empty context
(\textsf{iChk\_no\_contradictory\_pair}, from genuineness and
consistency), so the law is used under hypotheses, at abstract codes.

\paragraph{Worked example.}
Let $a$ code a one-entry list $[(H, c)]$ and let $b$ code any list
$Q$.  Then \rga proves that $e\cdot a\cdot b$ is the code of the list
whose head is $\emptyset \vdash \z = \tofn 1$, followed by the entry
$(H, c)$ and then by the entries of $Q$ (\textsf{goedelE\_example}).
So the head's two premises, the entry of $a$ and the head of $b$, sit at
positions $0$ and $|a| = 1$ of the list below the head, which the
head's negE detector cites.  A two-entry $b$ gives a four-entry list.

\paragraph{The decided D2 at every pair of terms.}
The names of this paragraph are in \texttt{RGA\_Goedel/}\allowbreak\texttt{Goedel\_D2}.
EXPL composes two accepted proofs into one at a
single instance.  The same technique gives the decided form of
Hilbert and Bernays's second condition, modus ponens on accepted
proofs, at every pair of well-formed terms $\varphi$, $\psi$, where
$\varphi \limp \psi$ is $\neg\varphi \lor \psi$.  The certificate
(\textsf{d2Cert}) places three new entries in front of the appended
list of a closed core $Q$, of $a$ and of $b$.  The root $\emptyset
\vdash \psi$ is a cut over the root $\emptyset \vdash \neg\varphi
\lor \psi$ of $a$ and the entry $\{\neg\varphi \lor \psi\} \vdash
\psi$.  That entry is a cut over the entry $\{\neg\varphi \lor
\psi\} \vdash \varphi$, a weakening of the root of $b$, and the root of
$Q$.  The core $Q$ (\textsf{d2Q0}) is a fixed six-entry proof of
$\{\neg\varphi \lor \psi, \varphi\} \vdash \psi$ by a disjunction
elimination whose first branch is a negE.  The core is a literal code,
so \rga accepts the core by closed evaluation of the validity test
(\textsf{d2\_valid\_accepted}).  If $\Gamma$ derives $a\jnat$,
$b\jnat$, $\Chk(a, \bquo{\varphi \limp \psi})$ and $\Chk(b,
\bquo{\varphi})$, then $\Gamma \vdash \Chk(\mathit{d2Cert}(a, b),
\bquo{\psi})$ (\textsf{d2Cert\_chk}), and the same holds for $m\cdot
a\cdot b$, where $m = \mathit{d2M}(\varphi, \psi)$ is one explicit
index tree for each pair (\textsf{d2M\_chk}).  The antecedent is
decided at naturals, so \textsf{ir\_impI\_bool} introduces the
implication over the antecedent, and two $\forall$-introductions give, for every
pair of well-formed terms $\varphi$, $\psi$,
\[
  \emptyset \vdash \forall a\, \forall b.\ \bigl(\Chk(a,
  \bquo{\varphi \limp \psi}) \wedge \Chk(b, \bquo{\varphi})\bigr)
  \limp \Chk(m\cdot a\cdot b, \bquo{\psi})
  \qquad (\textsf{D2\_decided}) .
\]
This result is a schema in $\varphi$ and $\psi$: one internal sentence,
an implication with a decided antecedent, for each pair.  The schema
uses no soundness theorem: the dependency check of
[\citen{conrga-paper}, Appendix ``The build-time check''] has \textsf{D2\_decided} and
\textsf{d2M\_chk} among its roots.  The core's acceptance cannot be
taken from the faithfulness of the validity test
(\textsf{iTreeValid\_faithful}), whose proof reaches the soundness
theorem in both directions.  Theorem~\ref{thm:g-d2} states the result
in Section~\ref{sec:goedel}.

\subsection{D3G with an explicit primitive-recursive witness}
\label{app:goedel-d3g}

\begin{defi}[\textsf{D3GF}, \textsf{D3G}]
\label{def:d-d3g}
For a term $\sigma$,
$\mathrm{D3GF}(\sigma) \ldef \forall v.\ \Chk(v,\bquo{G}) \limp
\Chk(\sigma\cdot v, \bquo{\neg G})$, and $\mathrm{D3G}(\sigma)$ is the
HOL proposition that $\sigma$ is closed,
$\emptyset \vdash \forall v.\ (\sigma\cdot v)\jnat$, and
$\emptyset \vdash \mathrm{D3GF}(\sigma)$.  The antecedent of
$\mathrm{D3GF}(\sigma)$ is decided:
$\Gamma \vdash a\jnat$ implies $\Gamma \vdash \Chk(a,\bquo{G})\jbool$
(\textsf{D3G\_antecedent\_B}).
\end{defi}

A proof of $\neg G = \neg(\cA\cdot F)$ needs, by notAllI, a refuted
instance $\neg(F\cdot a)$ and the typing $a\jnat$.  The witness built
here takes the instance at the numeral of $v$ itself.  The witness is
parametrized by an index $\mathit{cb}$, the \emph{premise builder}, that
supplies the proof of $\neg(F\cdot\tofn v)$.

\subsubsection{The witness over a premise builder}

\begin{defi}[\textsf{goedelSigmaIdx}, \textsf{goedelSigmaOf}, \textsf{D3Gcert}]
\label{def:d-sigma}
For an index $\mathit{cb}$, the index $\mathit{goedelSigmaIdx}(\mathit{cb})$ computes
\[
  (a,b) \mapsto \mathit{cons}\bigl(\mathit{jdg}(\bquo{\neg G}),\
  \mathit{app}(\cmp(\mathit{cb})\cdot a\cdot b,\ \mathit{natCert}(a,b))\bigr)
\]
from explicit trees, and $\sigma_{\mathit{cb}} \ldef \lambda v.\
\cmp(\mathit{goedelSigmaIdx}(\mathit{cb}))\cdot v\cdot v$.  The builder's
obligation $\mathrm{D3Gcert}(\mathit{cb})$ is the HOL proposition: for
every context $\Gamma$ and term $v$, if $\Gamma \vdash v\jnat$ and
$\Gamma \vdash \Chk(v, S_G)$, then
\[
  \Gamma \vdash \Chk\bigl(\cmp(\mathit{cb})\cdot v\cdot v,\
  \mathit{neg}(\mathit{ap}(\bquo{F}, \Num\cdot v))\bigr) .
\]
\end{defi}

\noindent
The code $\mathit{neg}(\mathit{ap}(\bquo{F}, \Num\cdot v))$
(\textsf{negInstC}) is the code of $\neg(F\cdot\tofn v)$ written for an
abstract $v$.  At a numeral $v = \tofn n$, \rga proves this code equal to
$\bquo{\neg(F\cdot\tofn n)}$, in every context
(\textsf{negInstC\_num}).

\begin{thm}[\textsf{natCert\_zero}, \textsf{natCert\_suc}, \textsf{iChk\_natCert}]
\label{thm:d-natcert}
$\mathit{natCertIdx}$ is the recursion index whose base is the
one-entry list $\mathit{cons}(\mathit{jdg}(\bquo{\z\jnat}), \z)$ and
whose step puts $\mathit{jdg}(\mathit{eq}(\Num\cdot\suc k,
\Num\cdot\suc k))$ in front of the previous list; $\mathit{natCert}(a,k)
\ldef \cmp(\mathit{natCertIdx})\cdot a\cdot k$.  If $\Gamma \vdash
a\jnat$ and $\Gamma \vdash k\jnat$, then
\[
  \Gamma \vdash \Chk\bigl(\mathit{natCert}(a,k),\ \mathit{eq}(\Num\cdot k,
  \Num\cdot k)\bigr) .
\]
\end{thm}

\noindent
\emph{Proof.}  One internal induction on $k$ (\textsf{ir\_nat\_ind}),
at a variable $v_i$ fresh for $a$ and not free in $\Gamma$, with motive
$\Chk(\mathit{natCert}(a,v_i), \mathit{eq}(\Num\cdot v_i, \Num\cdot
v_i))$.  The motive is a well-formed term with $v_i$ free, so the rule
applies.  For the base, the one-entry list is accepted by the rule that
$\z$ is a natural (\textsf{cert\_pr0}).  The list's conclusion $\bquo{\z\jnat}$
is provably $\mathit{eq}(\Num\cdot\z, \Num\cdot\z)$, and
$\mathit{natCert}(a,\z)$ is provably that list (\textsf{natCert\_zero},
by [\citen{conrga-paper}, theorem~\textsf{cmp\_prRec\_zero}]).  For the step, the context holds
$v_i\jnat$ and the motive at $v_i$.  The successor-typing rule (from
$c\jnat$ infer $\suc c\jnat$) accepts a new entry
$\mathit{eq}(\mathit{succ}(\Num\cdot v_i), \mathit{succ}(\Num\cdot v_i))$
on top of $\mathit{natCert}(a,v_i)$ (\textsf{cert\_link\_sucTI}).  Then
$\Num\cdot\suc v_i = \mathit{succ}(\Num\cdot v_i)$ (\textsf{numC\_suc}) and
the recursion equation of the index (\textsf{natCert\_suc}, which needs
the freshness side condition of [\citen{conrga-paper}, theorem~\textsf{cmp\_prRec\_zero}]) transport
the verdict to the motive at $\suc v_i$.  \qed

\noindent
The recursion of the certificate mirrors the induction that proves the
certificate accepted.  The value-certificate builder below has the same structure,
one level up.

\begin{thm}[\textsf{iChk\_negG\_notAllI}]
\label{thm:d-notalli}
If $\Gamma$ derives $a\jnat$, $P_1\jnat$, $P_2\jnat$,
$\Chk(P_1, \mathit{neg}(\mathit{ap}(\bquo F, a)))$ and $\Chk(P_2,
\mathit{eq}(a,a))$, then $\Gamma \vdash
\Chk(\mathit{cons}(\mathit{jdg}(\bquo{\neg G}), \mathit{app}(P_1,P_2)),
\bquo{\neg G})$.
\end{thm}

\noindent
\emph{Proof.}  As for Theorem~\ref{thm:d-goedelE-chk}, with notAllI in
place of negE.  The detector cites $z_1 = 0$ ($P_1$'s head) and $z_2 =
|P_1|$ ($P_2$'s head).  The operator tests read
$\mathit{neg}(\mathit{ap}(\bquo F, a))$ and $\mathit{eq}(a,a)$, and the
new conclusion is $\bquo{\neg G} =
\mathit{neg}(\mathit{ap}(\bquo{\cA}, \bquo F))$
(\textsf{goedelG\_neg\_code}, since $G = \cA\cdot F$).  \qed

\begin{thm}[\textsf{goedelSigmaIdx\_eq}, \textsf{iChk\_goedelSigmaOf}, \textsf{D3G\_of\_d3gcert}]
\label{thm:d-d3g-of-cert}
If $\Gamma \vdash a\jnat$ and $\Gamma \vdash b\jnat$, then $\Gamma
\vdash \cmp(\mathit{goedelSigmaIdx}(\mathit{cb}))\cdot a\cdot b =
\mathit{cons}(\mathit{jdg}(\bquo{\neg G}),
\mathit{app}(\cmp(\mathit{cb})\cdot a\cdot b, \mathit{natCert}(a,b)))$.
If $\Gamma \vdash v\jnat$ and $\Gamma \vdash
\Chk(\cmp(\mathit{cb})\cdot v\cdot v, \mathit{neg}(\mathit{ap}(\bquo F,
\Num\cdot v)))$, then $\Gamma \vdash \Chk(\sigma_{\mathit{cb}}\cdot v,
\bquo{\neg G})$.  Consequently $\mathrm{D3Gcert}(\mathit{cb})$ implies
$\mathrm{D3G}(\sigma_{\mathit{cb}})$.
\end{thm}

\noindent
\emph{Proof.}  The equation unfolds the explicit tree, as in
Theorem~\ref{thm:d-goedelE-eq}.  For the second part, a $\beta$ step
gives $\sigma_{\mathit{cb}}\cdot v =
\cmp(\mathit{goedelSigmaIdx}(\mathit{cb}))\cdot v\cdot v$
(\textsf{goedelSigmaOf\_beta}).  Apply Theorem~\ref{thm:d-notalli} with
$a = \Num\cdot v$, $P_1 = \cmp(\mathit{cb})\cdot v\cdot v$ and $P_2 =
\mathit{natCert}(v,v)$, which Theorem~\ref{thm:d-natcert} accepts for
$\mathit{eq}(\Num\cdot v, \Num\cdot v)$.  For the third part, closedness
is immediate, and totality is the $\beta$ step plus $\cmp$-totality
under $\forall$-introduction (\textsf{goedelSigmaOf\_total}).  For
$\mathrm{D3GF}(\sigma_{\mathit{cb}})$ (\textsf{goedelSigmaOf\_d3g}),
work in $\{\vz\jnat\}$ under the hypothesis $\Chk(\vz, \bquo{G})$.
Transport gives $\Chk(\vz, S_G)$, the obligation gives the premise, and
the second part gives $\Chk(\sigma_{\mathit{cb}}\cdot\vz, \bquo{\neg G})$.
The hypothesis is decided (Definition~\ref{def:d-d3g}), so
\textsf{ir\_impI\_bool} discharges it, and $\forall$-introduction
closes.  \qed

\noindent
\textsf{rga\_goedel2\_of\_d3gcert} combines the third part with
\textsf{rga\_goedel2\_of\_d3g} (\S\ref{app:goedel-assembly}):
$\mathrm{D3Gcert}(\mathit{cb})$ implies $\emptyset \nvdash \ConR$.
What remains is a premise builder that meets the obligation.  The core
of that builder is a certificate builder for compiled primitive-recursive applications.

\subsubsection{The value-certificate builder}

The premise $\neg(F\cdot\tofn v)$ is, after one $\beta$ step, the
negation of $\neg\Chk(\tofn v, S_G)$, and $\Chk$ is a conjunction of
equations between compiled applications.  So the premise builder must
certify the values of compiled applications.  The premise builder does so at
\emph{templates}: a HOL term $T$ with free variables, and an environment
$z$ that gives those variables numeral codes.  The environment's
well-formedness test $\mathit{sigWfGnChk}\cdot z$ says that every entry
of $z$ is the code of a well-formed term.

\begin{defi}[\textsf{prCertBuilder}, \textsf{prCertBuilderOK}]
\label{def:d-builder}
For each index $f$, the builder $B_f = \mathit{prCertBuilder}(f)$ is an
index, defined by HOL primitive recursion on $f$.  The builder's first argument
packs $\mathit{pair}(c_0, \mathit{pair}(x,y))$ and its second
$\mathit{pair}(P_a, P_b)$.  The proposition $\mathrm{OK}(f)$
(\textsf{prCertBuilderOK}) says: for every context $\Gamma$, every
term $z$ with $\Gamma \vdash z\jnat$ and $\Gamma \vdash
\mathit{sigWfGnChk}\cdot z$, all well-formed terms $T_1, T_2$, and all
terms $x, y, P_a, P_b$ that $\Gamma$ proves natural, if
\[
  \Gamma \vdash \Chk\bigl(P_a, \mathit{eq}(\mathit{img}_z(T_1), \Num\cdot x)\bigr)
  \quad\text{and}\quad
  \Gamma \vdash \Chk\bigl(P_b, \mathit{eq}(\mathit{img}_z(T_2), \Num\cdot y)\bigr) ,
\]
then, with $c_0 = \mathit{img}_z(\cmp(f)\cdot T_1\cdot T_2)$,
\[
  \Gamma \vdash \Chk\bigl(\cmp(B_f)\cdot \mathit{pair}(c_0, \mathit{pair}(x,y))
  \cdot \mathit{pair}(P_a,P_b),\
  \mathit{eq}(c_0, \Num\cdot(\cmp(f)\cdot x\cdot y))\bigr) .
\]
\end{defi}

\noindent
In words: if the two argument templates have certified values, the
compiled application has a certified value, and $B_f$ computes the
certificate.  The builder reads the codes it needs (argument codes,
redex codes, reducts) from $c_0$ with the checker's own kid readers and
compiled head reducer, so no environment and no template reaches the builder.

\begin{thm}[\textsf{pr\_cert\_builder}]
\label{thm:d-builder}
$\mathrm{OK}(f)$ holds for every index $f$.
\end{thm}

\noindent
The theorem is a schema over $f$: for each index $f$, it supplies a
derived rule of \rga at abstract arguments.  The theorem is not an internal
sentence that quantifies over indices.

\emph{Proof.}  HOL structural induction on $f$
(\textsf{prCertBuilderOK\_init}, \textsf{prCertBuilderOK\_sub},
\textsf{prCertBuilderOK\_rec}).  In
each case $\cmp(f)$ is a double abstraction $\lambda\lambda.\,b$, so
$c_0$ reaches $\mathit{img}_z(R)$, for the body instance $R$, in two head
steps of the code.  The top of every certificate is two
\emph{step-backs} (\textsf{iChk\_bldTop}).  A step-back turns a
certificate of $\mathit{eq}(c', w)$ into one of $\mathit{eq}(c, w)$ when
$c'$ is the head reduct of $c$ and $c, c', w$ are well-formed genuine
codes (\textsf{iChk\_certStepBack}).  The step-back's list consists of two symmetry
nodes and a transitivity node for $c'$'s self-equation, a head-step
conversion node for $c = c'$, and a transitivity node.  It remains to
certify $\mathit{img}_z(R)$.

\emph{Initial functions.}  $R$ is $\z$, $\suc T_1$, $T_1$ or $T_2$.  The
certificates are the zero list $\mathit{cons}(\mathit{jdg}(\bquo{\z\jnat}),
\z)$, a successor-congruence extension of $P_a$
(\textsf{iChk\_certSuc}), $P_a$, and $P_b$.

\emph{Composition} $f = \iSub\,g\,h_1\,h_2$.  Here $R =
\cmp(g)\cdot U_1\cdot U_2$ with $U_i = \cmp(h_i)\cdot T_1\cdot T_2$.  The
induction hypotheses for $h_1$ and $h_2$ at the templates $T_1, T_2$
certify $\mathit{img}_z(U_i)$.  The induction hypothesis for $g$ at the
nested templates $U_1, U_2$, in the same environment, certifies
$\mathit{img}_z(R)$.  The certified value is $\cmp(f)\cdot x\cdot y$ by the
composition law of $\cmp$.  No congruence node is needed, because
nesting the templates does the work.

\emph{Recursion} $f = \iRec\,g\,h\,h'$.  This case holds the only
internal induction of the builder, and this appendix exhibits that
induction in full.  The
compiled recursion is $\cmp(f) = \lambda\lambda.\ \cR\cdot(\cmp(g)\cdot
v_1\cdot\z)\cdot(\lambda\lambda.\ \cmp(h)\cdot(\cmp(h')\cdot v_3\cdot
v_1)\cdot v_0)\cdot v_0$ ([\citen{conrga-paper}, Appendix ``Coding and machines'']).  After the two
step-backs, the reduct is $R = \cR\cdot G_T\cdot H_T\cdot T_2$, with
$G_T = \cmp(g)\cdot T_1\cdot\z$ and $H_T$ the step abstraction
instantiated at $T_1$.  The builder reads $\hat G$ and $\hat H$, the
codes of $\mathit{img}_z(G_T)$ and $\mathit{img}_z(H_T)$, off
$\mathit{red}(\mathit{red}(c_0))$, and \rga proves the readings equal
to the images.  The certificate for $R$ substitutes the scrutinee's
value into the recursor spine (\textsf{iChk\_certSpine}).  From $P_b$
and a certificate of
$\mathit{eq}(\mathit{rap}(\hat G,\hat H,\Num\cdot y), \Num\cdot(\cmp(f)
\cdot x\cdot y))$ this construction builds a certificate of $\mathit{eq}(\mathit{img}_z(R),
\Num\cdot(\cmp(f)\cdot x\cdot y))$ (\textsf{prCertBuilderOK\_rec\_top}).
The missing certificate is the instance $q = y$ of the following lemma.
The lemma's builder $B_{\mathrm{rec}}$ (\textsf{bldRecIdx}) is itself a
recursion index on $q$, and the first argument of $B_{\mathrm{rec}}$, $m =
\mathit{pair}(\mathit{pair}(c_0, \mathit{pair}(x,y)),
\mathit{pair}(P_a,P_b))$, packs the parent's inputs.

\begin{thm}[\textsf{bldRec\_ind}]
\label{thm:d-bldrec}
Assume $\mathrm{OK}(g)$, $\mathrm{OK}(h)$ and $\mathrm{OK}(h')$, the
typing and well-formedness hypotheses of
Definition~\ref{def:d-builder} on $\Gamma$, $z$, $T_1$, $T_2$, $x$,
$y$, $P_a$, $P_b$, and the certificate $P_a$ of
$\mathit{img}_z(T_1)$ (no certificate of $P_b$ is needed).  If $\Gamma \vdash q\jnat$, then $\Gamma \vdash
M(q)$, where
\[
  M(q) \ldef \Chk\bigl(\cmp(B_{\mathrm{rec}})\cdot m\cdot q,\
  \mathit{eq}(\mathit{rap}(\hat G, \hat H, \Num\cdot q),\
  \Num\cdot(\cmp(f)\cdot x\cdot q))\bigr) .
\]
\end{thm}

\noindent
\emph{Proof.}  One internal induction (\textsf{ir\_nat\_ind}) on $q$, at
a variable $v_i$ fresh for $m$, $\hat G$, $\hat H$, $x$ and not free in
$\Gamma$.  The motive $M(v_i)$ is well-formed.
\begin{itemize}
\item \emph{Base.}  By the zero law of the recursion index
  ([\citen{conrga-paper}, theorem~\textsf{cmp\_prRec\_zero}]), $\cmp(B_{\mathrm{rec}})\cdot m\cdot\z$
  equals one step-back from $\mathit{rap}(\hat G,\hat H,\Num\cdot\z)$ to
  $\hat G$ over $c_g$.  Here $c_g$ is $B_g$ applied at the template
  $G_T$ with the argument certificates $P_a$ and the zero list.  The
  compiled head reducer sends $\mathit{rap}(\hat G,\hat H,\Num\cdot\z)$
  to $\hat G$ (\textsf{ir\_iHeadRedC\_R0}).  $\mathrm{OK}(g)$ certifies
  $\mathit{eq}(\mathit{img}_z(G_T), \Num\cdot(\cmp(g)\cdot x\cdot\z))$,
  and $\cmp(f)\cdot x\cdot\z = \cmp(g)\cdot x\cdot\z$ by the same zero
  law.  The step-back gives $M(\z)$.
\item \emph{Step.}  The context is $\Theta = \{v_i\jnat, M(v_i)\} \cup
  \Gamma$; write $k = v_i$.  The successor law of the recursion index
  gives $\cmp(B_{\mathrm{rec}})\cdot m\cdot\suc k =
  \cmp(\mathit{step})\cdot\mathit{pair}(m,k)\cdot c_{\mathrm{rec}}$ with
  $c_{\mathrm{rec}} = \cmp(B_{\mathrm{rec}})\cdot m\cdot k$.  So the step
  builder receives the previous certificate as data.  The step builder's output is three
  step-backs over a certificate $q_0$, along the head reductions
  \[
    e_0 = \mathit{rap}(\hat G,\hat H,\Num\cdot\suc k)
    \;\to\; e_1 = \mathit{ap}(\mathit{ap}(\hat H, \Num\cdot k),
      \mathit{rap}(\hat G,\hat H,\Num\cdot k))
    \;\to\; e_2 \;\to\; e_3 .
  \]
  The first reduction is the recursor's successor step at the code, since
  $\Num\cdot\suc k = \mathit{succ}(\Num\cdot k)$
  (\textsf{ir\_iHeadRedC\_RS\_num}).  The other two reductions are the $\beta$
  steps of $\hat H$'s double abstraction.  To certify $e_3$, extend the
  environment at a variable $v_j$ fresh for $T_1$: $z' = z$ with $v_j$
  sent to $\Num\cdot k$ (\textsf{updSubM}).  The images of $T_1$, $G_T$
  and $H_T$ do not change, since $v_j$ does not occur in them
  (\textsf{iSubM\_updSubM\_code\_nfree}), and $\mathit{img}_{z'}(v_j) =
  \Num\cdot k$.  Hence $e_1 = \mathit{img}_{z'}(H_T\cdot v_j\cdot X)$ for
  $X = \cR\cdot G_T\cdot H_T\cdot v_j$, and $e_3 =
  \mathit{img}_{z'}(\cmp(h)\cdot A_1\cdot X)$ for $A_1 =
  \cmp(h')\cdot T_1\cdot v_j$ (\textsf{gHR2\_iSubM\_beta2}).  Now three
  certificates combine.
  \begin{enumerate}
  \item $\mathrm{OK}(h')$ at the templates $T_1, v_j$, with the
    certificates $P_a$ and $\mathit{natCert}(k,k)$
    (Theorem~\ref{thm:d-natcert}), gives $q_1$ with
    $\Chk(q_1, \mathit{eq}(\mathit{img}_{z'}(A_1), \Num\cdot(\cmp(h')\cdot
    x\cdot k)))$.
  \item The induction hypothesis $M(k)$, transported along
    $\mathit{img}_{z'}(X) = \mathit{rap}(\hat G,\hat H,\Num\cdot k)$,
    says that $c_{\mathrm{rec}}$ certifies $\mathit{eq}(\mathit{img}_{z'}(X),
    \Num\cdot(\cmp(f)\cdot x\cdot k))$.
  \item $\mathrm{OK}(h)$ at the templates $A_1, X$, with the
    certificates $q_1$ and $c_{\mathrm{rec}}$, gives $q_0$ with
    $\Chk(q_0, \mathit{eq}(e_3, \Num\cdot(\cmp(h)\cdot(\cmp(h')\cdot
    x\cdot k)\cdot(\cmp(f)\cdot x\cdot k))))$.
  \end{enumerate}
  The successor law of $\cmp(f)$ turns the value into
  $\Num\cdot(\cmp(f)\cdot x\cdot\suc k)$, and the three step-backs give
  $M(\suc k)$.
\end{itemize}
The induction rule gives $M(q)$ for every $q$ with $\Gamma \vdash
q\jnat$.  \qed

\noindent
Three features make the technique work.  The HOL induction on $f$
supplies $\mathrm{OK}(g)$, $\mathrm{OK}(h)$ and $\mathrm{OK}(h')$ as
derived rules usable in any environment, while the internal induction
runs over the scrutinee's value $y$, which is a term, not a numeral, so
HOL cannot unfold the recursion $y$ times.  The step certificate is a function of
the previous certificate, so the builder of a recursion is itself a
recursion index.  And the environment is extended at a fresh variable,
not nested, so every image law is a closed-form computation on codes.
The decider-level proof of the search route has the same structure
(\S\ref{app:goedel-search}).

\subsubsection{The premise certificate}

Let $\Phi = \Chk(\vz, S_G)$, a term with $\vz$ free.  By
[\citen{conrga-paper}, definition~\textsf{iChk}], $\Phi = (A_1 \wedge A_2) \wedge (A_3 \wedge
A_4)$ (\textsf{chkPhi\_eq}), where $A_1$ says that the compiled validity
test at $\vz$ equals $\tofn 1$, $A_2$ that $\vz \ne \z$, $A_3$ that the
length of the head's hypothesis list is $\z$, and $A_4$ that the head's
conclusion equals $S_G = \cmp(\mathit{diagSubIdx})\cdot\bquo{D_G}\cdot\bquo{D_G}$.  The readers in $A_3$ and $A_4$ are compiled applications to the
predecessor of $\vz$, which is an $\cR$-term.

\begin{thm}[\textsf{iChk\_premCert}]
\label{thm:d-premcert}
If $\Gamma \vdash v\jnat$ and $\Gamma \vdash \Chk(v, S_G)$, then
\[
  \Gamma \vdash \Chk\bigl(\mathit{premCert}(v),\
  \mathit{neg}(\mathit{ap}(\bquo F, \Num\cdot v))\bigr) .
\]
\end{thm}

\noindent
\emph{Proof.}  Work with the images $\mathit{img}_v(t)$.  The image of
$\vz$ is $\Num\cdot v$, and $\mathit{natCert}(v,v)$ certifies that image.  The
environment $\zeta_v$ passes the well-formedness test
(\textsf{envVar0\_sigWfGnChk}).  The hypothesis supplies the four
conjuncts at $v$ (\textsf{iChk\_accD}, \textsf{iChk\_neD},
\textsf{iChk\_hypD}, \textsf{iChk\_concD}).  The certificate is built in
seven stages.
\begin{enumerate}
\item $A_1$: Theorem~\ref{thm:d-builder} for the validity index at the
  templates $\vz,\vz$ certifies the image of the test as
  $\Num\cdot(\text{test at } v)$, and the first conjunct turns that value
  into the image of $\tofn 1$.
\item $A_2$: since $v \ne \z$, \rga proves $v = \suc(p)$ for $p =
  \mathit{pred}(v)$.  The successor-not-zero rule, applied on top of
  $\mathit{natCert}(p,p)$ (\textsf{cert\_link\_sucNe0}), certifies
  \[
    \mathit{neg}(\mathit{eq}(\mathit{succ}(\Num\cdot p), \bquo\z)) ,
  \]
  which is provably the image of $A_2$.
\item The predecessor: a spine substitution (with
  $\mathit{natCert}(v,v)$ certifying the scrutinee) over three
  step-backs on $\mathit{natCert}(p,p)$, $p = \mathit{pred}(v)$,
  certifies the image of $\mathit{pred}(\vz)$ as
  $\Num\cdot\mathit{pred}(v)$ (\textsf{iChk\_certPred}).
\item $A_3$: Theorem~\ref{thm:d-builder} for the pair projection (twice)
  and the list length, chained through the previous certificate, certifies
  the image of the length reader.  The third conjunct turns that value into
  $\z$.
\item $A_4$: Theorem~\ref{thm:d-builder} for the conclusion projection
  certifies the image of the conclusion reader as $\Num\cdot S_G$, after
  transport along the fourth conjunct.  Theorem~\ref{thm:d-builder} for
  $\mathit{diagSubIdx}$ at the closed numeral templates $\bquo{D_G}, \bquo{D_G}$
  certifies the image of $S_G$ as $\Num\cdot S_G$.  A symmetry and a
  transitivity node (\textsf{iChk\_certEq}) give the image of $A_4$.
  $S_G$ stays a compiled application throughout: no closed numeral is
  evaluated.
\item Three conjunction nodes give $\mathit{img}_v(\Phi)$.  Since $a
  \wedge b \ldef \neg(\neg a \lor \neg b)$, each conjunction node is two
  double-negation nodes and one node of the rule that refutes a
  disjunction from refutations of both disjuncts (\textsf{iChk\_certAnd}).
\item The flips.  A double-negation node gives
  $\mathit{neg}(\mathit{neg}(\mathit{img}_v(\Phi)))$.  The false-equation
  introduction rule (from $\neg a$ infer $a = \tofn 0$) gives
  $\mathit{eq}(\mathit{neg}(\mathit{img}_v(\Phi)), \bquo{\tofn 0})$.
  Since $F\cdot\vz$ head-reduces to $\neg\Phi$ (\textsf{goedelF\_hstep}),
  a step-back at the image (\textsf{iChk\_certStepBack\_at}) gives
  $\mathit{eq}(\mathit{img}_v(F\cdot\vz), \bquo{\tofn 0})$.  The
  false-equation elimination rule gives
  $\mathit{neg}(\mathit{img}_v(F\cdot\vz))$.
\end{enumerate}
Finally $F$ is closed, so $\mathit{img}_v(F\cdot\vz) =
\mathit{ap}(\bquo F, \Num\cdot v)$, and the conclusion transports to
$\mathit{neg}(\mathit{ap}(\bquo F, \Num\cdot v))$.  \qed

\noindent
This theorem is formalized $\Sigma_1$-completeness for the decided family
$\Chk(v, S_G)$, at an abstract $v$, with the witness written out.
\S\ref{app:goedel-d3any} extends that completeness to the family $\Chk(v,
\bquo{\varphi})$ of every term $\varphi$.

\begin{thm}[\textsf{premCertIdx\_eval}, \textsf{D3Gcert\_holds}, \textsf{D3G\_goedelSigma}, \textsf{goedelSigma\_is\_pr}]
\label{thm:d-d3g}
There is an explicit index tree $\mathit{premCertIdx}$ such that
$\Gamma \vdash v\jnat$ implies $\Gamma \vdash
\cmp(\mathit{premCertIdx})\cdot v\cdot v = \mathit{premCert}(v)$.
Hence $\mathrm{D3Gcert}(\mathit{premCertIdx})$ and
$\mathrm{D3G}(\sigma)$ for $\sigma = \sigma_{\mathit{premCertIdx}}$.
Moreover $\sigma = \lambda v.\ \cmp(\tau_\sigma)\cdot v\cdot v$ with
$\tau_\sigma = \mathit{goedelSigmaIdx}(\mathit{premCertIdx})$, and
$\den{\tau_\sigma}$ is primitive recursive.
\end{thm}

\noindent
\emph{Proof.}  $\mathit{premCertIdx}$ mirrors $\mathit{premCert}$ node
for node.  The tree $\mathit{premCertIdx}$ is built from named sub-trees
($\mathit{premCertAIdx}$, \ldots, $\mathit{premCertFIdx}$) and from
combinators with evaluation laws: list cons, step-back, conjunction,
equation, predecessor, image and builder application
(\textsf{ir\_certConsIdxI}, \textsf{ir\_certStepBackIdxI},
\textsf{ir\_certAndIdxI}, \textsf{ir\_certEqIdxI},
\textsf{ir\_certPredIdxI}, \textsf{ir\_imgV0IdxI}).  Composing the laws
gives the equation.  With Theorem~\ref{thm:d-premcert} and transport in
the proof slot, the equation gives $\mathrm{D3Gcert}(\mathit{premCertIdx})$, and
Theorem~\ref{thm:d-d3g-of-cert} gives $\mathrm{D3G}(\sigma)$.  The last
statement is definitional.  \qed

\paragraph{What $\sigma$ outputs.}
At a number $v$, $\sigma\cdot v$ is the list
\[
  \mathit{cons}\bigl(\mathit{jdg}(\bquo{\neg G}),\
  \mathit{app}(\mathit{premCert}(v),\ \mathit{natCert}(v,v))\bigr) :
\]
a notAllI head over a proof of $\neg(F\cdot\tofn v)$ and a proof of
$\tofn v\jnat$, both themselves lists.  If $v$ is not a proof of $G$,
the output is still a number, since $\sigma$ has no guard, but not an
accepted one.  D3G claims acceptance only under its antecedent.  The
development computes no sizes of these lists.  The development asserts only that
$\sigma$ is compiled from one explicit index tree and that the outputs of
$\sigma$ are accepted under the antecedent.

\paragraph{No choice in the witness.}
The conjunct ``$\den{\tau_\sigma}$ is primitive recursive'' of
\textsf{goedelSigma\_is\_pr} holds of every index tree
(\textsf{pr2fi\_n2pr}), including one picked by Hilbert choice, so that
conjunct does not show that $\sigma$ is explicit.  What shows that
$\sigma$ is explicit is the
build-time audit of [\citen{conrga-paper}, \S\,``Machines without opaque indices'']: the definitional closure
of $\mathit{goedelSigmaOf}$, $\mathit{goedelSigmaIdx}$ and
$\mathit{premCertIdx}$ (1{,}182 constants) contains no
Hilbert-choice index and no other choice-defined constant outside the
check's whitelist, and the build fails otherwise
([\citen{conrga-paper}, Appendix ``Coding and machines'']).  The tree of $e$,
$\mathit{goedelEIdx}$, is a root of the same check (823 constants,
clean), as are the ex-falso tree $\mathit{exfIdx}$ and
the proof object $\mathit{neq01Prf}$ of $\neg(\z = \tofn 1)$, the tree
and the proof object both used for Theorem~\ref{thm:br-loeb}.

\subsubsection{The decided D3 for every sentence}
\label{app:goedel-d3any}

The names of this subsection are in \texttt{RGA\_Goedel/Goedel\_D3}.

The certificate of Theorem~\ref{thm:d-premcert} depends on $G$ in two
places only: in stage~5, which certifies the image of $S_G$ as
$\Num\cdot S_G$, and in the matrix $F$.  Replacing the stage-5 certificate by a
certificate for any closed code term $S$, and $F$ by $\lambda q.\,
\neg\Chk(q, S)$, the matrix of $\PrC$ at $S$, gives the decided D3 at
$S$.  The other stages are reused unchanged.  Say that $(S, \mathit{SI}, \mathit{sc})$ is a
\emph{D3 source} (\textsf{d3Src}) when $S$ is closed and well formed,
and for every context $\Gamma$ and every $v$ with $\Gamma \vdash
v\jnat$: $\Gamma \vdash S\jnat$, $\Gamma \vdash
\cmp(\mathit{SI})\cdot v\cdot v = S$, $\Gamma \vdash
(\cmp(\mathit{sc})\cdot v\cdot v)\jnat$, and $\Gamma \vdash
\Chk(\cmp(\mathit{sc})\cdot v\cdot v, \mathit{eq}(\mathit{img}_v(S),
\Num\cdot S))$: $\mathit{sc}$ builds an accepted proof that $S$
equals the numeral of its value.  Write $\sigma_{S}$
(\textsf{d3SigmaAt}) for the term $\lambda v.\ \cmp(\tau_S)\cdot v\cdot
v$ with $\tau_S = \mathit{d3SigmaIdx}(S, \mathit{SI}, \mathit{sc})$,
one explicit index tree.

\begin{thm}[\textsf{D3\_at}, \textsf{D3\_decided\_any}, \textsf{d3SigmaAt\_goedelS}, \textsf{D3GF\_goedelSigma\_inst}]
\label{thm:d-d3-any}
\leavevmode
\begin{enumerate}
\item If $(S, \mathit{SI}, \mathit{sc})$ is a D3 source, $T$ is
  closed, $\emptyset \vdash T\jnat$ and $\emptyset \vdash T = S$,
  then $\emptyset \vdash \forall v.\ \Chk(v, T) \limp
  \Chk(\sigma_S\cdot v, \bquo{\PrC(S)})$.
\item For every term $\varphi$, $\emptyset \vdash \forall v.\
  \Chk(v, \bquo{\varphi}) \limp \Chk(\sigma_\varphi\cdot v,
  \bquo{\PrC(\bquo{\varphi})})$, where $\sigma_\varphi$
  (\textsf{d3SigmaOf}) is $\sigma_S$ at $S = T = \bquo{\varphi}$.
\item At $S = S_G$, with $\mathit{SI}$ the diagonal tree of $S_G$ and
  $\mathit{sc}$ the tree of the stage-5 certificate, the term $\sigma_S$ is the
  witness $\sigma$ of Theorem~\ref{thm:d-d3g}, as an equation of HOL
  terms, and part~1 at $T = \bquo{G}$ is $\mathrm{D3GF}(\sigma)$.
\end{enumerate}
\end{thm}

\noindent
\emph{Proof.}  Part~1: the certificate of
Theorem~\ref{thm:d-premcert}, rebuilt at $S$ (\textsf{d3Cert},
\textsf{iChk\_d3Cert}), mirrored by one explicit tree
(\textsf{d3CertIdx\_eval}), under the notAllI graft of
Theorem~\ref{thm:d-notalli} at the matrix $\lambda q.\,\neg\Chk(q,
S)$ (\textsf{iChk\_notAllI\_at}).  The implication is introduced by
\textsf{ir\_impI\_bool} over the decided $\Chk(v, T)$, then
$\forall$-introduction.  Part~2: at a numeral code, the equation of
$S$ with its numeral is the numeral-typing proof $\mathit{natCert}$
(\textsf{d3Src\_num}).  Part~3: unfolding the definitions, every tree
of the general builder at $S_G$ is literally the corresponding tree of
$\mathit{premCertIdx}$, and $\lambda q.\,\neg\Chk(q, S_G)$ is $F$;
\textsf{d3Src\_goedelS} supplies the source and
$\vdash S_G = \bquo{G}$ the transport.  \qed

\noindent
So D3G is the construction of part~1 at the diagonal code term
$S = S_G$ (with $T = \bquo{G}$), and the transformer of D3G is literally
the witness $\sigma$ of G\"odel's second theorem.  D3G is not part~2 at $\varphi = G$:
that instance has the source $\bquo{G}$ in place of $S_G$, so the
instance's witness $\sigma_{G}$ and consequent code
$\bquo{\neg\forall q.\,\neg\Chk(q, \bquo{G})}$ differ from those of
D3G, whose consequent is $\bquo{\neg G} = \bquo{\neg\forall q.\,
\neg\Chk(q, S_G)}$.  G\"odel's second theorem needs only the construction at
$S_G$.  In part~2 the
transformer is a HOL function of $\varphi$, so part~2 is a schema: one
internal sentence for each $\varphi$.  The choice audit runs on
$\mathit{d3SigmaOf}$, $\mathit{d3SigmaAt}$, $\mathit{d3SigmaIdx}$
and $\mathit{d3CertIdx}$, all clean ([\citen{conrga-paper}, Appendix ``Coding and machines'']).

\subsubsection{The decided D3 as one sentence over the code}
\label{app:goedel-d3-uniform}

The names of this subsection are in
\texttt{RGA\_Goedel/}\allowbreak\texttt{Goedel\_D3\_Uniform}.  The schema of Theorem~\ref{thm:d-d3-any}(2) becomes
one internal sentence when the code term $S$ is itself a variable.  Let
$z_c$ (\textsf{envC}) be the environment that sends $v_1$ to $\Num\cdot
c$, and let $\mathit{pr}(c) = \mathit{img}_{z_c}(\PrC(v_1))$ (\textsf{prCU}):
the code of $\neg\forall q.\,\neg\Chk(q, v_1)$ with the numeral code
of $c$ substituted for $v_1$, computed inside \rga.  Let $\sigma_{\mathrm u} =
\cmp(\mathit{uSigIdx})$, one explicit index tree that takes the code $c$
and the proof code $v$ as arguments.

\begin{thm}[\textsf{D3\_uniform}, \textsf{prCU\_num}, \textsf{D3\_uniform\_num}]
\label{thm:d-d3-uniform}
\leavevmode
\begin{enumerate}
\item $\emptyset \vdash \forall c\, \forall v.\ \Chk(v, c) \limp
  \Chk(\sigma_{\mathrm u}\cdot c\cdot v, \mathit{pr}(c))$.
\item For every natural $k$, $\emptyset \vdash \mathit{pr}(\tofn k) =
  \bquo{\PrC(\tofn k)}$, and $\emptyset \vdash \forall v.\ \Chk(v,
  \tofn k) \limp \Chk(\sigma_{\mathrm u}\cdot \tofn k\cdot v, \bquo{\PrC(\tofn
  k)})$.
\end{enumerate}
\end{thm}

\noindent
\emph{Proof.}  Part~1: the premise certificate of
Theorem~\ref{thm:d-premcert} is rebuilt at the two-entry environment
that sends $\vz$ to $\Num\cdot v$ and $v_1$ to $\Num\cdot c$
(\textsf{envU}).  The variable $v_1$ plays the role of the code term
$S$ of Theorem~\ref{thm:d-d3-any}, and the stage-5 certificate at $v_1$ is the
numeral-typing proof of $c$ (\textsf{uCert}, \textsf{iChk\_uCert}).
The notAllI graft of Theorem~\ref{thm:d-notalli} is taken at the
matrix code term (\textsf{iChk\_notAllI\_code}).  The builder mirrors
the certificate as one explicit tree (\textsf{uCertIdx\_eval}), and
the transformer appends the graft (\textsf{uSigIdx\_eq}).  These
constructions give the derivation under the hypotheses $c\jnat$, $v\jnat$ and $\Chk(v,
c)$ of $\Chk(\sigma_{\mathrm u}\cdot c\cdot v, \mathit{pr}(c))$
(\textsf{iChk\_uSig}).  The checker verdict $\Chk(v, c)$ is decided at
naturals (\textsf{hq\_iChk}), so the implication is introduced over that verdict
by \textsf{ir\_impI\_bool}, and two $\forall$-introductions give the
sentence.  Part~2: at a numeral, $\mathit{pr}(\tofn k)$ is the compiled
substitution at a closed code, which \rga evaluates
(\textsf{iSubM\_closed}); the second sentence is the derivation of
part~1 at $c = \tofn k$, with this equation.  \qed

\noindent
Part~1 is one internal sentence, an implication whose antecedent is a
decided checker verdict; part~2 gives the numeral instances of that
sentence, one internal
sentence for each $k$, with the consequent read as the code of
$\PrC(\tofn k)$.  Both parts use no soundness theorem: the dependency
check of [\citen{conrga-paper}, Appendix ``The build-time check''] has \textsf{D3\_uniform}
and \textsf{D3\_uniform\_num} among its roots.  The choice audit runs
in the theory on $\mathit{uSigIdx}$, $\mathit{uCertIdx}$, $\mathit{pr}$ and
the index $\mathit{prCUIdx}$ of $\mathit{pr}$, all clean.  Theorem~\ref{thm:g-d3-uniform}
states the result in Section~\ref{sec:goedel}.

\subsection{The assembly}
\label{app:goedel-assembly}

\begin{thm}[\textsf{goedelG\_of\_expl\_d3g\_con}]
\label{thm:d-split}
If $\mathrm{EXPL}(e)$, $\mathrm{D3G}(\sigma)$ and $\emptyset \vdash
\ConR$, then $\emptyset \vdash G$.
\end{thm}

\noindent
\emph{Proof.}  Let $v = \vz$, $\Gamma_0 = \{v\jnat\}$ and $\Gamma_1 =
\{\Chk(v,S_G)\} \cup \Gamma_0$.  The derivation, as in the theory text:
\begin{enumerate}
\item $\Gamma_0 \vdash S_G\jnat$ (from Theorem~\ref{thm:d-goedelS}), so
  $\Gamma_0 \vdash \Chk(v,S_G)\jbool$ ([\citen{conrga-paper}, theorem~\textsf{hq\_iChk}]),
  that is, $\Gamma_0 \vdash \Chk(v,S_G) \lor \neg\Chk(v,S_G)$.
\item $\Gamma_1 \vdash \Chk(v,\bquo G)$: the case hypothesis,
  transported along $S_G = \bquo G$.
\item $\Gamma_1 \vdash (\sigma\cdot v)\jnat$: $\forall$-elimination of
  $\sigma$'s totality at $v$.
\item $\Gamma_1 \vdash \Chk(v,\bquo G) \limp \Chk(\sigma\cdot v,
  \bquo{\neg G})$: $\forall$-elimination of $\mathrm{D3GF}(\sigma)$ at
  $v$.  Modus ponens with step~2 gives $\Gamma_1 \vdash \Chk(\sigma\cdot
  v, \bquo{\neg G})$.
\item $\forall$-elimination of $\mathrm{EXPLF}(e)$ at $v$, then at
  $\sigma\cdot v$, and modus ponens with the conjunction of steps 2
  and~4, give $\Gamma_1 \vdash \Chk(e\cdot v\cdot(\sigma\cdot v),
  \bquo{\z = \tofn 1})$.
\item $\Gamma_1 \vdash (e\cdot v\cdot(\sigma\cdot v))\jnat$:
  $\forall$-elimination of $e$'s totality at $v$, then at $\sigma\cdot
  v$.
\item $\Gamma_1 \vdash \neg\Chk(e\cdot v\cdot(\sigma\cdot v), \bquo{\z
  = \tofn 1})$: $\forall$-elimination of $\ConR$ at that term, with
  step~6.
\item $\Gamma_1 \vdash \neg\Chk(v,S_G)$: negation elimination on steps
  7 and~5.
\item $\{\neg\Chk(v,S_G)\} \cup \Gamma_0 \vdash \neg\Chk(v,S_G)$:
  hypothesis.
\item $\Gamma_0 \vdash \neg\Chk(v,S_G)$: disjunction elimination on
  steps 1, 8 and~9.
\item $\emptyset \vdash \forall q.\ \neg\Chk(q,S_G) = G$:
  $\forall$-introduction (\textsf{ir\_allI\_open}).  \qed
\end{enumerate}

\noindent
Each $\forall$-elimination is \textsf{gd2\_allE\_beta}: from $\Gamma
\vdash \cA\cdot(\lambda b)$ and $\Gamma \vdash a\jnat$, infer the
$\beta$-instance $b[a]$.  The count is one $\forall$-introduction, one
case split on a decided sentence (steps 1 and~10), one transport along
$S_G = \bquo G$, seven $\forall$-eliminations, one conjunction
introduction, two modus ponens and one negation elimination.  The
derivation introduces no implication.  It discharges three hypotheses:
the two case hypotheses, by disjunction elimination, and the typing
$v\jnat$ of the eigenvariable, by $\forall$-introduction.  $\ConR$
enters only through $\forall$-elimination at one term.

The typing hypothesis of an eigenvariable is not a decided sentence:
\rga does not even prove $L\jnat \limp L\jnat$ for the Liar $L$
(\textsf{sat\_wall\_eigen}, \texttt{Shadows/Shadows\_Traps.thy}).  Every
$\forall$-introduction and internal induction reasons under such a
hypothesis.  What the proof avoids is an implication introduction or a
case split over an undecided sentence.  A check that runs with every
build confirms this avoidance on the theorem-dependency graph
([\citen{conrga-paper}, Appendix ``The build-time check'']).  In the cone of
Theorem~\ref{thm:g-g2} (\textsf{rga\_goedel2}), the theorems of the
two G\"odel theories that apply an implication-introduction rule are
exactly \textsf{goedelE\_expl}, the internal sentence of EXPL
(Theorem~\ref{thm:d-expl}), and \textsf{goedelSigmaOf\_d3g}, that of
D3G (Theorem~\ref{thm:d-d3g-of-cert}); each of these two theorems applies only
\textsf{ir\_impI\_bool}, whose antecedent must be decided.  The only
theorem there that applies a case-split rule is
\textsf{goedelG\_of\_expl\_d3g\_con}, the split above.  The same
holds for the search route (\textsf{rga\_goedel2\_via\_search}), with
\textsf{goedelSigmaSearch\_d3g} for D3G.  The lemmas below the two
theories that the cone reaches (the checker, the coding kit, the
$\forall$-free laws) introduce implications too.  A second check, which
also runs with every build and fails it otherwise, covers every theory
that the two cones reach (127 for \textsf{rga\_goedel2}, 139 for the
search route): no theorem there applies \textsf{ir\_impI\_provable},
whose antecedent need only be provable, or \textsf{ir\_impI\_Bw}.  The
other implication-introduction rules that the first check watches
either discharge a decided antecedent, a term $a$ with $a\jbool$ derived or $n \ne \z$
for a natural $n$, or discharge nothing (\textsf{ir\_impI1},
\textsf{ir\_impI2}, \textsf{ir\_impI\_vac}).
The induction hypotheses of the internal
inductions of \S\ref{app:goedel-d3g} are checker verdicts at natural
arguments, which [\citen{conrga-paper}, theorem~\textsf{hq\_iChk}] decides, although \rga's
induction rule does not require decidedness.

\begin{thm}[\textsf{rga\_no\_proof\_0eq1}, \textsf{goedelG\_refutes}]
\label{thm:d-refutes}
$\emptyset \nvdash \z = \tofn 1$.  If $\emptyset \vdash G$, then
$\emptyset \vdash \z = \tofn 1$.
\end{thm}

\noindent
\emph{Proof.}  \rga refutes $\z = \tofn 1$ (\textsf{ir\_num\_ne}), and
\rga does not prove both a term and its negation
(\textsf{rga\_consistent}).  For the second part, a derivation of $G$
gives a proof $P$ with $P \proves \emptyset \vdash G$
(\textsf{entails\_provesE}).  Faithfulness gives $\emptyset \vdash
\Chk(\bquo P, \bquo G)$, and transport gives $\Chk(\bquo P, S_G)$.
$\forall$-elimination of $G$ at $\bquo P$ gives $\neg\Chk(\bquo P,
S_G)$, and negation elimination gives $\z = \tofn 1$.  \qed

\noindent
The second part has no hypothesis on \rga: it is a transformation of
derivations, G\"odel~I as a recipe.

\begin{thm}[\textsf{rga\_goedel2\_of\_expl\_d3g}, \textsf{rga\_goedel2\_of\_d3g}, \textsf{rga\_goedel2}]
\label{thm:d-g2}
If $\mathrm{EXPL}(e)$ and $\mathrm{D3G}(\sigma)$, then $\emptyset \nvdash
\ConR$.  If $\mathrm{D3G}(\sigma)$, then $\emptyset \nvdash \ConR$.  And
$\emptyset \nvdash \ConR$.
\end{thm}

\noindent
\emph{Proof.}  A derivation of $\ConR$ would give a derivation of $G$
(Theorem~\ref{thm:d-split}) and then a derivation of $\z = \tofn 1$
(Theorem~\ref{thm:d-refutes}), which does not exist.  The second form
instantiates $e$ by Theorem~\ref{thm:d-expl}.  The third form instantiates
$\sigma$ by Theorem~\ref{thm:d-d3g}.  \qed

\begin{thm}[\textsf{rga\_goedel1}, \textsf{goedelG\_no\_accepted\_code}, \textsf{rga\_goedel1\_unrefutable}]
\label{thm:d-g1}
$\emptyset \nvdash G$.  For every natural $p$, $\emptyset \nvdash
\Chk(\tofn p, \bquo G)$.  And $\emptyset \nvdash \neg G$.
\end{thm}

\noindent
\emph{Proof.}  The first statement is Theorem~\ref{thm:d-refutes}.  For the
second statement, genuineness ([\citen{conrga-paper}, theorem~\textsf{iChk\_faithful}]) would make $p$ the code of a
proof of $G$.  For the third statement, suppose $\emptyset \vdash \neg(\cA\cdot
F)$.  \rga's witness property (\textsf{rga\_not\_all\_witness}) gives a
natural $n$ with $\emptyset \vdash \neg(F\cdot\tofn n)$.  The $\beta$
step $F\cdot\tofn n \to \neg\Chk(\tofn n, S_G)$, refutation transport
along a head step (\textsf{ir\_conv\_negE}) and double-negation
elimination give $\emptyset \vdash \Chk(\tofn n, S_G)$.  Transport gives
$\Chk(\tofn n, \bquo G)$, against the second statement.  \qed

\noindent
The third part uses no Rosser sentence.  The role of G\"odel's
$\omega$-consistency hypothesis is played by the witness property, a
consequence of \rga's soundness.  With consistency, the witness property
implies that \rga is $\omega$-consistent: \rga never refutes a universal
$\cA\cdot u$ all of whose numeral instances $u\cdot\tofn n$ \rga proves
(\textsf{rga\_omega\_consistent}).

\noindent
The second statement says that no numeral is provably accepted.  That
statement does not say that D3G holds for every $\sigma$; it says only that
$\mathrm{D3GF}(\sigma)$ is true in the classical reading for every
$\sigma$, because its antecedent is refuted at every numeral.  The
content of D3G is \rga's proof of $\mathrm{D3GF}(\sigma)$.

\paragraph{The internal half on its own.}
The theorems \textsf{goedelG\_of\_con} (\texttt{Goedel\_II\_PR.thy})
and \textsf{goedelG\_of\_con\_}\allowbreak\textsf{via\_search}
(\texttt{Goedel\_II.thy}) are rules: $\emptyset \vdash \ConR$ gives
$\emptyset \vdash G$, at the two witnesses.  The premise is false by
Theorem~\ref{thm:d-g2}, so as rules both theorems are vacuous.  The derivation
itself is not.  Run over an arbitrary context
(\textsf{goedelG\_of\_expl\_d3g\_ctx}: $\mathrm{EXPL}(e)$,
$\mathrm{D3G}(\sigma)$ and $\Gamma \vdash \ConR$ give $\Gamma \vdash
G$), the derivation gives the entailment $\ConR \vdash G$
(\textsf{goedelG\_under\_ConRGA} in \texttt{Goedel\_II\_PR.thy},
Theorem~\ref{thm:br-hyp}, at the explicit witness $\sigma$ of
Theorem~\ref{thm:d-d3g}, so at the checker), while the implication $\emptyset \vdash \ConR \limp
G$ fails (\textsf{no\_formalized\_goedel1}).  The entailment is the
input of formalized G\"odel~II (\S\ref{app:goedel-ent}).  The
proof of the entailment is a transformation of derivations that does not
use \rga's consistency.  The Gentzen recipe of \S\ref{app:goedel-cons}
(\textsf{rga\_con\_int\_of\_gentzen}) does not use consistency either;
\textsf{no\_gentzen\_pair} combines that recipe with consistency, through
Theorem~\ref{thm:d-g2}.

\subsection{The search route}
\label{app:goedel-search}

\textsf{rga\_goedel2\_via\_search} (\texttt{Goedel\_II.thy}) proves the
same statement $\emptyset \nvdash \ConR$ with a different D3G witness,
$\sigma_{\mathrm{s}} = \lambda v.\ \mathit{search}\cdot v\cdot\z$
(\textsf{goedelSigmaSearch}).  The search is a term-level fixed point that returns the first $P \ge 0$
passing a decided stopping test: $v$ is not an accepted proof of $G$, or
$P$ is an accepted proof of $\neg G$.  The test uses the transparent
proves-checker $\pi$ of [\citen{conrga-paper}, the table of the three proof checkers], never the decider.
Table~\ref{tab:d-search} lists the steps.  These steps replace
Definition~\ref{def:d-sigma} through Theorem~\ref{thm:d-d3g} only; EXPL
and the assembly are shared.

\begin{table}[!htb]
\centering
\footnotesize
\begin{tabular}{@{}>{\raggedright\arraybackslash}p{0.30\linewidth}>{\raggedright\arraybackslash}p{0.64\linewidth}@{}}
\hline
Isabelle name & statement or role \\
\hline
\textsf{ir\_goedelSearch\_unfold} &
  $\mathit{search}\cdot V\cdot P$ is $P$ if the test is $\z$, else
  $\mathit{search}\cdot V\cdot\suc P$ (four head steps) \\
\textsf{goedelSearch\_halts} &
  if $\Gamma$ derives $V\jnat$, $q\jnat$ and ``the test at $(V,q)$ is
  $\z$'', the search from $\z$ stops at a stopping point; one internal
  induction over the distance to $q$ \\
\textsf{D3Gdec} (interface) &
  for every $\Gamma$ with $\mathrm{wfc}(\Gamma)$ and every $v$: $\Gamma
  \vdash v\jnat$ and $\Gamma \vdash \Chk(v,S_G)$ imply $\Gamma \vdash
  \Tru(\mathit{neg}(\mathit{ap}(\bquo F, \Num\cdot v)))$ \\
\textsf{iHasCert\_negInstC}, \textsf{iHasCert\_eqReflC} &
  from that truth, $\Cert_\emptyset(\mathit{neg}(\mathit{ap}(\bquo F,
  \Num\cdot x)))$; and $\Cert_\emptyset(\mathit{eq}(\Num\cdot x,
  \Num\cdot x))$; both certificates through the rule form
  \textsf{iHasCert\_of\_iPrDecTrueAt} of [\citen{conrga-paper}, theorem~\textsf{rga\_complete\_int2}] (and
  \textsf{rga\_wtrueE}), not the internal sentence \textsf{rga\_complete\_int} \\
\textsf{iHasCert\_negG\_notAllI} &
  the two certificates combine under one notAllI head
  (Theorem~\ref{thm:d-notalli}) into $\Cert_\emptyset(\bquo{\neg G})$ \\
\textsf{goedelSearch\_halts\_of\_cert}, \textsf{goedelSigmaSearch\_total},
\textsf{goedelSigmaSearch\_d3g}, \textsf{D3G\_of\_d3gdec} &
  a certificate of $\neg G$ is a stopping point, so the search halts;
  hence totality and $\mathrm{D3GF}(\sigma_{\mathrm s})$, given
  \textsf{D3Gdec} \\
\textsf{pr\_dec\_eval} &
  the decider-level twin of Theorem~\ref{thm:d-builder}: for every index
  $f$, the decider reads the image of $\cmp(f)\cdot T_1\cdot T_2$ at
  value $\cmp(f)\cdot x\cdot y$ when it reads the images of $T_1, T_2$ at
  $x, y$ ($\forall$-free templates); HOL induction on $f$, one internal
  induction in the recursion case \\
\textsf{iPrDecVal\_iChk}, \textsf{iWtrue\_negInstC\_of\_iChk},
\textsf{D3Gdec\_holds} &
  the internal truth schema for the one term $\Chk(\vz,S_G)$: the
  decider affirms its image at $v$; then $\Tru$ of the code of
  $\neg(F\cdot\tofn v)$ (no $\mathrm{wfc}$ premise); the interface holds \\
\textsf{D3G\_goedelSigmaSearch}, \textsf{rga\_goedel2\_via\_search} &
  $\mathrm{D3G}(\sigma_{\mathrm s})$; $\emptyset \nvdash \ConR$ \\
\hline
\end{tabular}
\caption{The search route to $\emptyset \nvdash \ConR$
  (\texttt{RGA\_Goedel/Goedel\_II.thy}).  The value wrapper of
  \textsf{pr\_dec\_eval} is $\mathit{DecVal}(c,x) \ldef \exists l.\
  \Dec(l,c) = \suc x$.}
\label{tab:d-search}
\end{table}

\noindent
The two routes share everything except D3G.  The headline route needs
neither the decider nor the truth predicate nor internal completeness:
its D3G is built entirely at the checker.
The internal-reasoning paper's appendix on the build-time check
[\citen{conrga-paper}, Appendix ``The build-time check''] compares the sizes of the two
routes' dependency cones, and the decider-free check of that appendix confirms that the
headline route and the internal second theorem use neither the decider
nor the truth predicate.  Extensionally
the search is the
zero function, because no number is a proof of $G$
(Theorem~\ref{thm:d-g1}) and the test therefore stops at $P = \z$.
\rga does not prove $\forall v.\ \sigma_{\mathrm s}\cdot v = \z$
(\textsf{goedelSigmaSearch\_not\_provably\_zero}).  \emph{Proof.}
Suppose \rga did.  At an eigenvariable $v$, split on the decided verdict
$\Chk(v, \bquo G)$.  In the accepting case,
$\mathrm{D3GF}(\sigma_{\mathrm s})$ gives $\Chk(\sigma_{\mathrm s}\cdot
v, \bquo{\neg G})$, the assumption at $v$ and transport in the proof
slot give $\Chk(\z, \bquo{\neg G})$, and the checker rejects the code
$\z$ (\textsf{iChk\_neD}).  So $\forall q.\ \neg\Chk(q,\bquo G)$ follows,
which \textsf{goedelG\_code\_unprovable} excludes.  \qed

\subsection{The entailments and the cut transformers}
\label{app:goedel-ent}

This subsection proves the entailment results of
Sections~\ref{sec:breaks} and~\ref{sec:goedel}.  The names of this
subsection are mainly in five theories: \texttt{RGA/RGA\_Sound} (the $\Sigma$ kit and the
exact deduction theorem, in the companion formalization),
\texttt{RGA\_Kit/Entailment\_Kit} (the entailment lemmas over the
$\omega$-semantics), \texttt{RGA\_Goedel/Entailment\_Goedel} (the
G\"odel ingredients), \texttt{RGA\_Goedel/Goedel\_Cut} (the cut
transformers, formalized G\"odel~II and L\"ob) and
\texttt{RGA\_Liar\_Div/Entailment\_Liar} (the truth-predicate Liar and
the search for a proof of $\z = \tofn 1$).  Write $\exists D$ for
$\neg(\cA\cdot(\lambda\neg D))$, with $D$ a term in which $\vz$ is the
bound variable, and $D[a]$ for the $\beta$-instance of $D$ at $a$.  In this
notation $\PrC(c)$ is $\exists\,\Chk(\vz, \lft c)$
(\textsf{iPrc}, an abbreviation), and $\PrC(\bquo{\z = \tofn 1})$ is
literally $\neg\ConR$ (\textsf{iPrc\_cEq01}).

\paragraph{The $\Sigma$ kit.}

\begin{thm}[\textsf{sigma\_intro}, \textsf{sigma\_elim}, \textsf{sigma\_\allowbreak deduction\_\allowbreak theorem}, \textsf{sigma\_reductio}, \textsf{sigma\_contrapos}, \textsf{deduction\_\allowbreak theorem\_\allowbreak exact}]
\label{thm:d-sigma}
Let $\Gamma$ be a context and $D$ a term.
\begin{enumerate}
\item If $\Gamma \vdash D[a]$ and $\Gamma \vdash a\jnat$, then
  $\Gamma \vdash \exists D$.  If $\Gamma \vdash \exists D$ and $D,\
  \vz\jnat,\ \lft\Gamma \vdash \lft C$, then $\Gamma \vdash C$.
\end{enumerate}
Suppose further that $D$ is decided at a fresh variable:
$\vz\jnat,\ \lft\Gamma \vdash D\jbool$.
\begin{enumerate}
\setcounter{enumi}{1}
\item $\exists D,\ \Gamma \vdash C$ if and only if $\Gamma \vdash
  \forall.\,(D \limp \lft C)$.
\item If $\exists D,\ \Gamma \vdash \neg\exists D$, then $\Gamma \vdash
  \neg\exists D$.
\item If $\exists D,\ \Gamma \vdash B$, then $\neg B,\ \Gamma \vdash
  \neg\exists D$.
\end{enumerate}
Finally, for every well-formed $A$: every entailment $A,\ \Gamma
\vdash B$ gives $\Gamma \vdash A \limp B$ if and only if $\Gamma
\vdash A \lor \neg A$.
\end{thm}

\noindent
\emph{Proof.}  (1) Introduction refutes the instance
$(\lambda\neg D)\cdot a$, which head-steps to $\neg D[a]$, so a
refutation transported backward along the head step
(\textsf{ent\_conv\_negI}) and the rule \textsf{pr\_not\_allI} give
$\exists D$.  Elimination is the rule \textsf{pr\_not\_allE} with the
eigenvariable $\vz$: the rule's hypothesis, the refuted lifted instance,
becomes $D$ by a $\beta$-step under negation and double-negation
elimination.  Neither half needs $D$ decided.  (2) From the sentence:
eliminate $\exists D$, eliminate the universal at $\vz$, and apply
modus ponens.  To the sentence: open $\vz\jnat$; under the
hypothesis $D$, introduce $\exists D$ at $\vz$ and cut against the
lifted entailment; discharge $D$, which is decided
(\textsf{ir\_impI\_bool}); introduce the universal
(\textsf{sigma\_dt\_imp}, \textsf{sigma\_dt\_ent}).  (3) At $\vz$,
split on $D$.  In the case $D$, introduction gives $\exists D$, the
entailment gives $\neg\exists D$, and negation elimination gives
$\neg D$; in the case $\neg D$ there is nothing to do.  So $\forall.\,
\neg D$, and double-negation introduction gives $\neg\exists D$,
which is $\neg\neg\forall.\,\neg D$.  (4) is (3) in the context
$\{\neg B\} \cup \Gamma$, where the entailment and $\neg B$ give
$\neg\exists D$ by negation elimination.  The last statement: the
``if'' direction is \textsf{ir\_impI\_bool}; for ``only if'', take $B
= A$, so that $\Gamma \vdash \neg A \lor A$.  \qed

\noindent
Two more generic facts complete the kit.  An entailment fails if its
closed hypothesis is true in the $\omega$-semantics and its conclusion
is not (\textsf{ent\_omega\_refute}), since \rga is sound for that
semantics in every context (\textsf{rga\_omega\_sound}); and a
universal is true there when \rga derives each of its numeral
instances (\textsf{omv\_all\_of\_instances}).  A rule in every context
is an entailment (\textsf{ent\_of\_ctx\_rule},
\texttt{Entailment\_Truth}).

\paragraph{From decided forms to entailments.}
\textsf{iPrc\_transfer} turns a checker law into an entailment between
provability formulas: if, in every context $\Delta$, $\Delta \vdash
v\jnat$ and $\Delta \vdash \Chk(v, a)$ give $\Delta \vdash T(v)\jnat$
and $\Delta \vdash \Chk(T(v), b)$, then $\PrC(a),\ \Gamma \vdash
\PrC(b)$, for closed $a$, $b$ with $a$ well-formed.  The proof is
Theorem~\ref{thm:d-sigma}(1) twice: eliminate at $\vz$, apply the law,
introduce at $T(\vz)$.  \textsf{iPrc\_transfer2} is the same with two
hypotheses.  With the rule forms of the decided conditions as the
laws, this transfer gives Theorem~\ref{thm:br-dc-ent}(2) to~(4):
\textsf{D2\_ent} from \textsf{d2M\_chk}, \textsf{D3\_ent} from the
decided D3 at $\varphi$ (and \textsf{D3\_at\_ent} at any code term with
a certificate source), \textsf{D3G\_ent} at $S_G$ with transport,
\textsf{expl\_ent} from Theorem~\ref{thm:d-goedelE-chk}, and
\textsf{exf\_ent} from \textsf{exfIdx\_chk\_sf}, the ex-falso law
re-proved by closed evaluation so that the law uses no soundness
(\textsf{iChk\_accept\_sf}).  The diagonal pairs of
Theorem~\ref{thm:br-diag-ent} transport the verdict along $S_G =
\bquo{G}$ under the binder (\textsf{all\_negchk\_transport}), and
along $S_{\Tru} = \bquo{L_{\Tru}}$ for the Liar
(\textsf{iWtrue\_neg\_transport}).  \textsf{con\_notprc\_goedelG\_ent}
is a cut of \textsf{goedelG\_under\_ConRGA} with the first half of
\textsf{goedel\_diag\_ent}.  The pair $\neg G \eqv \PrC(\bquo{G})$
(\textsf{goedel\_diag\_neg\_ent}) is \textsf{iPrc\_eq\_transfer} along
$S_G = \bquo{G}$, since $\neg G$ is literally $\PrC(S_G)$; and
\textsf{not\_goedelG\_not\_con\_ent} composes that pair with
\textsf{D3G\_ent} and \textsf{expl\_ent}.  Contraposition of
entailments at a decided hypothesis is \textsf{dec\_contrapos}
(\texttt{RGA\_Sound}): from $\Gamma \vdash A \lor \neg A$ and $A,\
\Gamma \vdash B$, split on $A$ in the context $\{\neg B\} \cup
\Gamma$.

\paragraph{The cut transformers.}
For a proof object $Q$ of $\{A\} \vdash B$, let
\[
  \mathit{cut}_Q(a) \;=\; \mathit{cons}(\mathit{jdg}(\bquo{B}),\
  \mathit{app}(\bquo{Q}, a))
\]
(\textsf{cut1Cert}).  The closed list $\bquo{Q}$ is accepted, and the
head of $\bquo{Q}$ is the pair of the codes of $\{A\}$ and of $B$, both facts by closed
evaluation of the checker at the numeral $\bquo{Q}$
(\textsf{closed\_proof\_head}, through \textsf{d2\_valid\_accepted});
no soundness theorem is used.  Given $\Gamma \vdash a\jnat$ and
$\Gamma \vdash \Chk(a, \bquo{A})$, the appended list is valid
(\textsf{iTreeValid\_appC}), with $Q$'s head at position $0$ and $a$'s
head, $\emptyset \vdash A$, at position $|\bquo{Q}|$.  The cut
detector fires on the new entry $\emptyset \vdash B$
(\textsf{ndFire\_cut\_intro}): the entry at the second cited position
has the new entry's (empty) hypothesis set and conclusion $A$, and the
entry at position $0$ has hypothesis set $\{A\}$, computed as the
insertion of $A$ into the new entry's set, and the new entry's
conclusion $B$.  Theorem~\ref{thm:d-cons} accepts the list
(\textsf{cut1Cert\_chk}).  The index $\mathit{cut1Idx}(Q, B)$ builds
the list from explicit trees for a numeral, append and cons
(\textsf{certConsIdx}, \textsf{certAppIdx}), and $c_Q =
\cmp(\mathit{cut1Idx}(Q, B))$ satisfies $\Gamma \vdash c_Q\cdot a\cdot
a = \mathit{cut}_Q(a)$ (\textsf{cut1M\_eq}); transport in the proof
slot gives \textsf{cut1M\_chk}.  For $Q$ a proof of $\{A_1, A_2\}
\vdash B$, the list puts three entries above $\mathit{app}(\bquo{Q},
\mathit{app}(a, b))$: $\{A_1\} \vdash A_2$ by weakening of $b$'s
head, $\{A_1\} \vdash B$ by cut against $Q$'s head, and $\emptyset
\vdash B$ by cut against $a$'s head (\textsf{cut2Cert},
\textsf{cut2Cert\_chk}, \textsf{cut2M\_chk}).  Necessitation
(\textsf{ent\_nec}, \textsf{ent\_nec2}) obtains $Q$ from the
derivation in HOL (\textsf{entails\_provesE}) and applies
\textsf{iPrc\_transfer} with $T(v) = c_Q\cdot v\cdot v$, whose totality
is $\cmp$-totality.  The choice audit at the end of
\texttt{Goedel\_Cut} has the four transformer constants as roots
([\citen{conrga-paper}, Appendix ``The build-time check'']).

\paragraph{Formalized G\"odel~II and L\"ob.}
\textsf{rga\_goedel2\_internal\_contra} (an alias of
\textsf{prc\_con\_ent\_bot}) composes three entailments
(\textsf{ent\_trans}, \textsf{ent\_trans2}): \textsf{ent\_nec} at
\textsf{goedelG\_under\_ConRGA} gives $\PrC(\bquo{\ConR}) \vdash
\PrC(\bquo{G})$, \textsf{D3G\_ent} gives $\PrC(\bquo{\neg G})$, and
\textsf{expl\_ent} gives $\PrC(\bquo{\z = \tofn 1})$, which is
$\neg\ConR$.  \textsf{rga\_goedel2\_internal} is
Theorem~\ref{thm:d-sigma}(4) at the decided matrix of
$\PrC(\bquo{\ConR})$ (\textsf{iPrc\_dec}), followed by a cut with
double-negation introduction.  For L\"ob, let $\psi$ be closed and
well-formed, $D_\psi = \PrC(\Diag(\vz, \vz)) \limp \psi$
(\textsf{loebD}), $S_\psi = \Diag(\bquo{D_\psi}, \bquo{D_\psi})$
(\textsf{loebS}) and $L_\psi = \PrC(S_\psi) \limp \psi$
(\textsf{loebL}).  Substituting $\bquo{D_\psi}$ into $D_\psi$ gives
$L_\psi$ (\textsf{loeb\_diag}), and $\emptyset \vdash S_\psi =
\bquo{L_\psi}$ (\textsf{loebS\_eq}), so $\PrC(S_\psi) \eqv
\PrC(\bquo{L_\psi})$ (\textsf{iPrc\_eq\_transfer}).  The proof of
\textsf{formalized\_loeb\_ent} then follows the four steps of
\S\ref{sec:goedel-formal} on the formalized theorems: disjunction elimination with
\textsf{sigma\_contrapos} for $L_\psi \vdash \PrC(\bquo{L_\psi})
\limp \psi$; \textsf{ent\_nec}, \textsf{D3\_ent} at the numeral
$\bquo{L_\psi}$ and \textsf{D2\_ent} for $\PrC(\bquo{L_\psi}) \vdash
\PrC(\bquo{\psi})$; disjunction elimination with
\textsf{sigma\_contrapos} for $\PrC(\bquo{\psi}) \limp \psi \vdash
L_\psi$; and \textsf{ent\_nec} with composition.  D3 is applied at the
numeral $\bquo{L_\psi}$, so the L\"ob sentence needs no certificate
source of its own.  \textsf{formalized\_loeb\_0eq1\_ent} is the
instance at $\z = \tofn 1$.

\paragraph{The sentences in decided-antecedent form and their witnesses.}
The names of this paragraph and the next three are in
\texttt{RGA\_Goedel/Goedel\_Cut}.  \textsf{prc\_sentence\_iff} is
Theorem~\ref{thm:d-sigma} at the decided matrix of $\PrC(c)$
(\textsf{iPrc\_dec}), for a closed code term $c$ with $\Gamma \vdash
c\jnat$ and a conclusion that the lift leaves unchanged.  Each sentence
of Theorem~\ref{thm:g-sentence}(1, 3, 4) is that equivalence applied to
one entailment: \textsf{rga\_goedel2\_internal\_contra} for
\textsf{rga\_goedel2\_sentence}; \textsf{goedel\_diag\_neg\_ent}
followed by \textsf{not\_goedelG\_not\_con\_ent} for
\textsf{goedel1\_sentence}, and \textsf{not\_goedelG\_not\_con\_ent}
alone, at the hypothesis
$\neg G = \PrC(S_G)$, for \textsf{goedel1\_sentence\_diag}; and
\textsf{formalized\_loeb\_ent} for \textsf{formalized\_loeb\_sentence}.
The explicit witnesses are two index trees.  \textsf{goedel1Idx} is
$\iSub$ of EXPL's tree \textsf{goedelEIdx} over the left projection
and the tree of D3G at the diagonal term, so that at a proof $a$ of
$G$ the tree computes $e\cdot a\cdot(\sigma\cdot a)$
(\textsf{goedel1Idx\_chk}, from \textsf{goedelE\_chk} and
\textsf{iChk\_d3SigmaAt} after transport of $\Chk(a, \bquo G)$ to
$\Chk(a, S_G)$).  \textsf{goedel2Idx}~$Q$ composes \textsf{goedel1Idx} with the cut
tree \textsf{cut1Idx}~$Q$~$G$ at a HOL proof object $Q$ of $\{\ConR\}
\vdash G$ (\textsf{goedel2Idx\_chk}, from \textsf{cut1M\_chk}).  The
composition laws are \textsf{ir\_prSubI} and \textsf{cmpi\_L\_eq}.
The terms $T_1 = \lambda(\cmp(\textsf{goedel1Idx})\cdot\vz\cdot\vz)$
(\textsf{goedel1T}) and $T_Q$ (\textsf{goedel2T}) wrap the trees as
D3G's $\sigma$ does (\textsf{cmp\_diag\_beta}), and the generic
\textsf{chk\_imp\_sentence} turns a checker law at the eigenvariable
into the sentence, introducing the implication over the decided
verdict (\textsf{chk\_imp\_matrix\_B}: both sides are decided).  These
constructions give \textsf{goedel1\_sentence\_witness} and
\textsf{goedel2T\_sentence}; \textsf{rga\_goedel2\_sentence\_witness}
takes $Q$ from \textsf{goedelG\_under\_ConRGA} by
\textsf{entails\_provesE}.  The separate halves are
\textsf{goedelE\_expl} and \textsf{D3GF\_goedelSigma\_inst}; without the
witness trees the two halves are composed only inside the proof of
Theorem~\ref{thm:g-split}.  The eighth run of the choice audit,
at the end of the theory, has the four constants as roots
([\citen{conrga-paper}, Appendix ``The build-time check'']).  The failure of the implication,
\textsf{no\_formalized\_goedel2\_nonprenex}, unfolds the implication to
$\neg\PrC(\bquo{\ConR}) \lor \neg\ConR$ and applies the disjunction
property against \textsf{no\_code\_unprovable} and
\textsf{ConRGA\_unrefutable}; \textsf{prenex\_extraction\_fails\_con}
pairs that failure with \textsf{rga\_goedel2\_sentence}.

\paragraph{Explicit reflection at each numeral.}
The analogue of explicit reflection in Artemov's Logic of Proofs
(Section~\ref{sec:related}) holds at every numeral.

\begin{thm}[\textsf{chk\_reflect\_num}, \textsf{chk\_reflect\_num\_ent}]
\label{thm:d-chk-reflect}
For every number $n$ and every well-formed term $A$, $\emptyset \vdash
\Chk(\tofn n, \bquo{A}) \limp A$, and $\Chk(\tofn n, \bquo{A}), \Gamma
\vdash A$ for every context $\Gamma$.
\end{thm}

\noindent
\emph{Proof.}  At two numerals the checker's verdict is a decided
Boolean (\textsf{hq\_iChk}), so by the disjunction property
(\textsf{ir\_disj\_prop}) \rga derives $\Chk(\tofn n, \bquo{A})$ or its
negation.  In the first case the numeral is the code of a proof of $A$
(\textsf{iChk\_genuine\_proves}, [\citen{conrga-paper}, theorem~\textsf{iChk\_faithful}]), so \rga
derives $A$, and the implication by its second disjunct; in the second
case the implication follows by its first disjunct.  The entailment
follows by modus ponens.  $A$ must be well formed, or $\neg\Chk(\tofn
n, \bquo{A}) \lor A$ is not a term; $A$ may be open.  \qed

\medskip
\noindent
The proof is a case split in HOL, so the theorem is a schema, one
theorem for each numeral and term, not one internal sentence; the
sentence over all proofs at $\z = \tofn 1$ fails
(Theorem~\ref{thm:br-loeb}(4)).  The theorem uses soundness, through the
disjunction property and genuineness, but not the decider: both forms
are roots of the decider-free audit of
[\citen{conrga-paper}, Appendix ``The build-time check''].

\paragraph{G\"odel's proof read with entailments.}
\textsf{goedelG\_under\_ConRGA\_textbook} derives $\ConR \vdash G$ by
the textbook chain, independently of the decided split behind
\textsf{goedelG\_under\_ConRGA}: \textsf{not\_goedelG\_not\_con\_ent}
gives $\neg G \vdash \neg\ConR$ through D3G and EXPL; the hypothesis
of that entailment,
$\neg G = \PrC(S_G)$, has a decided matrix, so \textsf{sigma\_contrapos}
gives $\neg\neg\ConR \vdash \neg\neg G$; double-negation introduction
at $\ConR$ and elimination at $G$ give $\ConR \vdash G$.  A proof of
$\ConR$ then cuts to a proof of $G$, against \textsf{rga\_goedel1}
(\textsf{rga\_goedel2\_textbook}).

\paragraph{Interderivabilities and the theory \mbox{\rga{}${}+\ConR$}.}
\textsf{prc\_of\_derivable} is D1 without soundness: the code of an
\rga proof of $\varphi$ is accepted by closed evaluation of the checker
(\textsf{iChk\_proves\_sf}), and that code witnesses the existential, in
every context.  \textsf{goedelG\_con\_ent} and \textsf{notcon\_prc\_con\_ent}
are proved as in \S\ref{sec:goedel-formal}, and the two
interderivabilities pair these two entailments with \textsf{goedelG\_under\_ConRGA} and
\textsf{prc\_con\_ent\_bot}.  In \texttt{RGA\_Goedel/Entailment\_Goedel},
\textsf{con\_ctx\_omega\_sound} is \textsf{rga\_omega\_sound} with
\textsf{omv\_ConRGA\_true} for the members of the context that are
$\ConR$, and \textsf{con\_ent\_0eq1\_fails},
\textsf{goedelG\_ent\_0eq1\_fails} and
\textsf{con\_goedelG\_ent\_0eq1\_fails} follow because $\z = \tofn 1$ is
not $\omega$-true.

\paragraph{The derivability conditions and Beklemishev and Shamkanov's conditions.}
These names are in \texttt{RGA\_Goedel/Goedel\_Conditions}.
\textsf{chk\_derivability\_conditions} is the conjunction of D1 as a
rule, the decided D2 for every pair of well-formed formulas, the
decided D3 over the code (\textsf{D3\_uniform}), ex falso
(\textsf{exf\_int\_sf}, statement-equal to \textsf{exf\_int} and
proved without soundness) and EXPL.  For Beklemishev and Shamkanov's
abstract provability structures~\cite{beklemishev16abstract}, read $x
\leqslant y$ as $x \vdash y$, $\top$ as $\z = \z$, the structures' $\bot$ as
$\neg(\z = \z)$ (a classical falsum, not \rga's divergent constant
$\bot$), $\Box x$ as $\PrC(\bquo{x})$ and $\boxtimes x$ as
$\PrC(\bquo{\neg x})$.  \textsf{bs\_structure\_conditions} collects the
preorder, the monotone half of C1 (\textsf{ent\_nec}), C2 (D1 at
$\neg\neg(\z = \z)$), C3 (EXPL, then ex falso at $\neg(\z = \z)$), C4
(D3 at $\neg x$) and C5 (weakening of $\vdash \z = \z$), at every term.
The anti-monotone half of C1 is $\Sigma$- or decided contraposition
followed by \textsf{ent\_nec} (\textsf{bs\_C1\_anti\_sigma},
\textsf{bs\_C1\_anti\_dec}); the class of decided sentences and
$\Sigma$ sentences with a decided matrix (\textsf{bsSD}) contains
$\top$, $\bot$, every $\Box x$ and $\boxtimes x$, and $\neg G$
(\textsf{bs\_aps\_sigma}).  The Jeroslow fixed point $\neg G \eqv
\boxtimes\neg G$ uses $\neg G = \PrC(S_G)$, the code equation and
double negation under $\Box$ (\textsf{bs\_jeroslow\_fp}).
\textsf{bs\_theorem1} and \textsf{bs\_theorem1\_i\_rule} follow
Beklemishev and Shamkanov's proof for every $\Sigma$ fixed point with
a decided matrix, and \textsf{bs\_goedel2} instantiates that proof at $\neg G$:
$\boxtimes\boxtimes\top \vdash \boxtimes\top$, and $\boxtimes\top \eqv
\neg G$; \textsf{bs\_incon\_con} identifies $\boxtimes\top$ with
$\neg\ConR$ up to interderivability.  So the instance moves under
$\PrC$: Theorem~\ref{thm:d-sigma}(4) applied to \textsf{bs\_incon\_con}
gives, with double negation, $\ConR \vdash \neg\boxtimes\top$;
\textsf{ent\_nec} gives
$\PrC(\bquo{\ConR}) \vdash \boxtimes\boxtimes\top$, and with
\textsf{bs\_goedel2} and \textsf{bs\_incon\_con} this entailment is
$\PrC(\bquo{\ConR}) \vdash \neg\ConR$, the entailment that
Theorem~\ref{thm:g-formal}(1) proves directly.  The development has
each step but not the composition as one theorem.  The anti-monotone half fails in
general (\texttt{RGA\_Liar\_Div/Entailment\_Liar}):
$L_{\Tru} \wedge \neg L_{\Tru} \vdash \z = \tofn 1$, but
$\PrC(\bquo{\neg(\z = \tofn 1)}) \nvdash \PrC(\bquo{\neg(L_{\Tru}
\wedge \neg L_{\Tru})})$, since with D1 and the reflection rule that entailment
would give the contrapositive that \textsf{ent\_contrapos\_fails}
excludes (\textsf{bs\_C1\_anti\_fails}, \textsf{bs\_not\_aps}).

\paragraph{The failing entailments.}
Table~\ref{tab:d-ent-fail} lists the failing entailments with their routes.  The
$\omega$-route needs the $\omega$-values of three closed sentences:
$G$ and $\ConR$ are true there because \rga derives each of their
numeral instances (\textsf{omv\_goedelG\_true},
\textsf{omv\_ConRGA\_true}), and so is $L_{\Tru}$
(\textsf{omv\_goedelW\_true}).  A truth-predicate conclusion
$\Tru(c)$ at the code $c$ of a sentence that is never grounded-true
is false there, since \rga refutes each instance $\Aff(\tofn n, c)$
by closed evaluation of the decider and the decider's soundness
(\textsf{iWtrue\_inst\_refuted}, \textsf{omv\_iWtrue\_not\_true}).

\begin{table}[!htb]
\centering
\footnotesize
\begin{tabular}{@{}>{\raggedright\arraybackslash}p{0.30\linewidth}>{\raggedright\arraybackslash}p{0.36\linewidth}>{\raggedright\arraybackslash}p{0.28\linewidth}@{}}
\hline
Isabelle name & entailment that fails & route \\
\hline
\textsf{d1\_ent\_fails\_goedelG} & $G \vdash \PrC(\bquo G)$ & $\omega$: \textsf{goedelG\_code\_instances\_refuted}, \textsf{omv\_prc\_goedelG\_false} \\
\textsf{reflect\_ent\_fails\_0eq1}, \textsf{reflect\_ent\_iff\_con}, \textsf{reflect\_ent\_iff\_con\_ctx} & $\PrC(\bquo{\z = \tofn 1}) \vdash \z = \tofn 1$, equivalent to $\vdash \ConR$; in a context $\Gamma$, to $\Gamma \vdash \ConR$, so it holds under $\ConR$ (\textsf{con\_reflect\_0eq1\_ent}) & reductio; Theorem~\ref{thm:d-g2} \\
\textsf{reflect\_chk\_ent\_fails} & $\vz\jnat,\ \Chk(\vz, \bquo{\z = \tofn 1}) \vdash \z = \tofn 1$ & $\Sigma$ deduction theorem; the above \\
\textsf{cert\_reflect\_ent\_fails} & $\Cert_\emptyset(\bquo{\z = \tofn 1}) \vdash \z = \tofn 1$ & \textsf{sigma\_refute\_bot}, \textsf{ConRGA\_of\_cert\_refuted} \\
\textsf{prc\_true\_ent\_fails\_0eq1}, \textsf{true\_reflect\_chk\_ent\_fails}, \textsf{prc\_true\_ent\_fails} & $\PrC(\bquo{\z = \tofn 1}) \vdash \Tru(\bquo{\z = \tofn 1})$; the same from $\Chk(\vz, \cdot)$; $\PrC(\vz),\ \vz\jnat \vdash \Tru(\vz)$ & $\Sigma$ deduction theorem; [\citen{conrga-paper}, theorem~\textsf{rga\_no\_uniform\_soundness}] \\
\textsf{con\_bool\_ent\_fails} & $\ConR\jbool \vdash \ConR$; $\ConR\jbool \vdash \neg\ConR$ & reductio; $\omega$ \\
\textsf{no\_formalized\_goedel2} & (implication) $\vdash \ConR \limp \neg\PrC(\bquo\ConR)$ & disjunction property; \textsf{ConRGA\_unrefutable}, \textsf{no\_code\_unprovable} \\
\textsf{no\_formalized\_}\allowbreak\textsf{goedel2\_nonprenex} & (implication) $\vdash \PrC(\bquo\ConR) \limp \neg\ConR$ & the same \\
\textsf{con\_ent\_0eq1\_fails}, \textsf{goedelG\_ent\_0eq1\_fails}, \textsf{con\_goedelG\_ent\_0eq1\_fails} & $\ConR \vdash \z = \tofn 1$; $G \vdash \z = \tofn 1$; both hypotheses together & $\omega$ (\textsf{con\_ctx\_omega\_sound}) \\
\textsf{tbicond\_liar\_unprovable} & (implication) $\vdash \Tru(\bquo{L_{\Tru}}) \liff L_{\Tru}$ & the biconditional is explosive (\textsf{tbicond\_liar\_explosive}) \\
\textsf{tin\_ent\_fails\_liar}, \textsf{tin\_ent\_fails\_goedelG}, \textsf{tin\_ent\_fails\_con} & $t \vdash \Tru(\bquo t)$ at $L_{\Tru}$, $G$, $\ConR$ & $\omega$ \\
\textsf{tout\_ent\_fails\_liar}, \textsf{negnec\_ent\_fails\_liar} & $\Tru(\bquo{L_{\Tru}}) \vdash L_{\Tru}$; $\neg L_{\Tru} \vdash \neg\Tru(\bquo{L_{\Tru}})$ & reductio at $\Tru(S_{\Tru})$ (\textsf{liar\_true\_refl\_refuted}); \textsf{goedelW\_unprovable} \\
\textsf{stops\_ent\_fails\_0eq1} & $\mathit{Stops} \vdash \z = \tofn 1$ & cut with \textsf{prc\_stops\_ent}; the reflection failure \\
\textsf{ent\_contrapos\_fails} & $\neg(\z = \tofn 1) \vdash \neg(L_{\Tru} \wedge \neg L_{\Tru})$, the contrapositive of a derivable entailment & it would give $\vdash L_{\Tru}\jbool$; \textsf{goedelW\_undecided} \\
\textsf{complete\_imp\_fails} & (implication) $\vdash \forall c.\ \Tru(c) \limp \Cert_\emptyset(c)$, guarded or not & disjunction property at $\bquo{L_{\Tru}}$ \\
\hline
\end{tabular}
\caption{The failing entailments of the reflection pillars
  [\citen{conrga-paper}, \S\,``The reflection pillars''] and of
  Sections~\ref{sec:breaks} and~\ref{sec:goedel}, and four failing implications,
  in \texttt{Entailment\_Goedel}, \texttt{Entailment\_Truth},
  \texttt{Entailment\_Liar} and \texttt{Goedel\_Cut}.  ``$\omega$'': the hypothesis is true in
  the $\omega$-semantics and the conclusion is not
  (\textsf{ent\_omega\_refute}).  ``Reductio'': the hypothesis has the
  $\Sigma$ shape, and Theorem~\ref{thm:d-sigma}(3) would refute the hypothesis.
  The other failures of the truth-axiom table are transfers of failing
  rules on theorems by substitution and cut.}
\label{tab:d-ent-fail}
\end{table}

\paragraph{The search for a proof of $\z = \tofn 1$.}
$\mathit{Stops} = \Rcd_1(\bquo{\mathit{conSearch}})$ (\textsf{conStops},
an abbreviation).  \textsf{prc\_stops\_ent} eliminates the hypothesis
$\PrC(\bquo{\z = \tofn 1})$ at $\vz$ and applies
\textsf{con\_search\_found}, which derives $\mathit{Stops}$ from an
accepted proof of $\z = \tofn 1$ at an abstract code in any
well-formed context.  \textsf{not\_stops\_con\_ent} is
Theorem~\ref{thm:d-sigma}(4) and double-negation elimination; the
other three entailments from $\neg\mathit{Stops}$ cut
\textsf{not\_stops\_con\_ent} against
\textsf{goedelG\_of\_expl\_d3g\_ctx}, \textsf{goedel\_diag\_ent} and
\textsf{rga\_goedel2\_internal}.  \textsf{conStops\_undecided}
(\texttt{ConOGA/}\allowbreak\texttt{ConOGA\_Residual\_Frag}) combines
the record's faithfulness at the $\forall$-free search, which has no
record value (\textsf{conSearch\_no\_orec}), with
\textsf{con\_search\_unrefutable} and the disjunction property.  Two
questions are open.  $\mathit{Stops} \vdash \PrC(\bquo{\z = \tofn 1})$
would need the converse of \textsf{con\_search\_found} at an abstract
fuel, an internal inversion of the record machine through the search's
stages, which the development does not have.  The entailment $\ConR \vdash
\neg\mathit{Stops}$ would also follow from that inversion.  Neither entailment has a countermodel
in the grounded or the $\omega$-semantics.

\subsection{Consequences}
\label{app:goedel-cons}

Table~\ref{tab:d-cons} lists the main theorems that the development
derives from G\"odel's theorems.  Not listed are variants of the
Gentzen-shape theorem (\textsf{no\_gentzen\_pair\_gen}, for any
provably correct iteration of $r$, and \textsf{no\_gentzen\_pair\_ex},
for an existential form of $\mathrm{GTI}$) and
\textsf{con\_search\_natural\_unprovable} (\texttt{RGA\_Liar\_Div/}\allowbreak\texttt{Divergence}),
which says that \rga does not prove the natural closure obligation of
the proof search for $\z = \tofn 1$, because that obligation is
equivalent to $\ConR$.  Most proofs follow one pattern: open an
eigenvariable, split on a decided checker verdict, and close with
$\forall$-introduction to obtain a sentence that
Theorem~\ref{thm:d-g2} excludes.  Other proofs use the disjunction and witness
properties, or G\"odel's first theorem (Theorem~\ref{thm:d-g1}).

\paragraph{The ex-falso tree.}
Theorem~\ref{thm:br-loeb}(1) uses a second transformer of the kind of
$e$.  For each term $\varphi$, the index $\mathit{exfIdx}(\varphi)$ is
built from the same explicit trees as $\mathit{goedelEIdx}$ and computes
\[
  (a,b) \mapsto \mathit{cons}(\mathit{jdg}(\bquo{\varphi}),\
  \mathit{app}(a, \bquo{\mathit{neq01Prf}})) ,
\]
where $\mathit{neq01Prf}$ is a fixed three-node proof of $\neg(\z =
\tofn 1)$, whose code the checker accepts in every context
(\textsf{iChk\_neq01Prf}).  The second argument is unused.  Write
$x_\varphi \ldef \cmp(\mathit{exfIdx}(\varphi))$.

\begin{thm}[\textsf{exfIdx\_chk}]
\label{thm:d-exf}
For every well-formed $\varphi$, every context $\Gamma$ and every term
$a$: if $\Gamma \vdash a\jnat$ and $\Gamma \vdash \Chk(a, \bquo{\z =
\tofn 1})$, then $\Gamma \vdash \Chk(x_\varphi\cdot a\cdot a,
\bquo{\varphi})$.
\end{thm}

\noindent
\emph{Proof.}  As for Theorem~\ref{thm:d-goedelE-chk}.
Theorem~\ref{thm:d-appC} gives the accepted list $\mathit{app}(a,
\bquo{\mathit{neq01Prf}})$, with $\mathit{jdg}(\bquo{\z = \tofn 1})$ at
position $0$ and $\mathit{jdg}(\bquo{\neg(\z = \tofn 1)})$ at position
$|a|$.  The new root fires negE citing these two positions, and the
root's well-formedness check holds because $\varphi$ is well-formed
(\textsf{iWft\_code\_ne0}).  \qed

\noindent
Ex falso is also one internal implication: for every well-formed
$\varphi$, $\emptyset \vdash \forall a.\ \Chk(a, \bquo{\z = \tofn 1})
\limp \Chk(x_\varphi\cdot a\cdot a, \bquo{\varphi})$
(\textsf{exf\_int}).  The proof applies Theorem~\ref{thm:d-exf} at an
eigenvariable, introduces the implication over the decided verdict
([\citen{conrga-paper}, theorem~\textsf{hq\_iChk}], \textsf{ir\_impI\_bool}) and closes with
$\forall$-introduction, as for EXPL.

\medskip
\noindent
\textsf{ConRGA\_of\_code\_unprovable} then derives $\ConR$ from
$\emptyset \vdash \forall q.\ \neg\Chk(q, \bquo{\varphi})$, for
well-formed $\varphi$: at an eigenvariable $v$ it splits on the decided
verdict $\Chk(v, \bquo{\z = \tofn 1})$, and in the accepting case the
hypothesis at $x_\varphi\cdot v\cdot v$ contradicts
Theorem~\ref{thm:d-exf}.  No implication is introduced.

\begingroup
\footnotesize
\begin{longtable}{@{}>{\raggedright\arraybackslash}p{0.31\linewidth}>{\raggedright\arraybackslash}p{0.36\linewidth}>{\raggedright\arraybackslash}p{0.27\linewidth}@{}}
\caption{The main theorems derived from G\"odel's theorems, and \rga's
  $\omega$-consistency, in
  \texttt{RGA\_Goedel/} unless a theory is named; each theorem cites
  \textsf{rga\_goedel2} where the theorem uses G\"odel's second theorem.  ``Decided split'' is
  the case split of Theorem~\ref{thm:d-split} on a checker verdict at an
  eigenvariable.}
\label{tab:d-cons}\\
\hline
Isabelle name & statement & proof \\
\hline
\endfirsthead
\hline
Isabelle name & statement & proof \\
\hline
\endhead
\hline
\endfoot
\textsf{rga\_no\_uniform\_soundness} &
  $\emptyset \nvdash \forall p.\ \Chk(p,\bquo{\z = \tofn 1}) \limp
  \Tru(\bquo{\z = \tofn 1})$; [\citen{conrga-paper}, theorem~\textsf{rga\_no\_uniform\_soundness}] &
  decided split; \textsf{iPrDecTrueAt\_cEq01\_never}; \textsf{rga\_goedel2} \\
\textsf{rga\_no\_loeb\_hypothesis\_0eq1} &
  $\emptyset \nvdash \forall p.\ \Chk(p,\bquo{\z = \tofn 1}) \limp \z =
  \tofn 1$; Thm.~\ref{thm:br-loeb}(4) &
  decided split; \textsf{ir\_num\_ne}; \textsf{rga\_goedel2} \\
\textsf{ConRGA\_of\_code\_unprovable}, \textsf{no\_code\_unprovable},
\textsf{exfIdx\_chk} &
  for well-formed $\varphi$, $\emptyset \vdash \forall q.\
  \neg\Chk(q,\bquo{\varphi})$ implies $\emptyset \vdash \ConR$; so
  $\emptyset \nvdash \forall q.\ \neg\Chk(q,\bquo{\varphi})$;
  Thm.~\ref{thm:br-loeb}(1) &
  decided split on $\Chk(v,\bquo{\z = \tofn 1})$, ex-falso tree
  $x_\varphi$ (Thm.~\ref{thm:d-exf}); \textsf{rga\_goedel2} \\
\textsf{rga\_loeb\_rule}, \textsf{rga\_no\_formalized\_loeb\_0eq1},
\textsf{no\_D2\_schema\_goedelG}, \textsf{no\_D3\_schema\_goedelG} &
  L\"ob's rule; formalized L\"ob fails at $\z = \tofn 1$ as an implication;
  D2 and D3 fail at $G$ as implications; Thms.~\ref{thm:br-loeb}(2,3),
  \ref{thm:br-d2d3} &
  disjunction property; \textsf{no\_code\_unprovable},
  \textsf{goedelG\_code\_unprovable}, \textsf{ConRGA\_unrefutable},
  \textsf{rga\_pr\_reflect}, \textsf{rga\_goedel1} \\
\textsf{D1\_rule\_pr}, \textsf{D2\_rule\_pr}, \textsf{D3\_rule\_pr},
\textsf{rga\_pr\_reflect} &
  the rule forms of D1--D3 and the reflection rule;
  Thm.~\ref{thm:br-rules} &
  D1: faithfulness and $\neg\forall$-introduction; reflection: witness
  property, genuineness; D2 and D3 from these two \\
\textsf{no\_goedel\_biconditional}
(\texttt{RGA\_Checker/Chk\_Faithful}), \textsf{no\_liar\_biconditional}
(\texttt{RGA\_Decider/Decider\_Truth}) &
  the diagonal biconditionals fail as implications, for every term;
  Thm.~\ref{thm:br-bicond} &
  disjunction and witness properties, genuineness, faithfulness,
  \textsf{rga\_consistent} \\
\textsf{ConRGA\_unrefutable}, \textsf{ConRGA\_undecided} &
  $\emptyset \nvdash \neg\ConR$; $\emptyset \nvdash \ConR\jbool$;
  Thm.~\ref{thm:br-con} &
  witness property, genuineness; disjunction property,
  \textsf{rga\_goedel2} \\
\textsf{no\_formalized\_goedel1}, \textsf{no\_formalized\_goedel1\_code} &
  $\emptyset \nvdash \ConR \limp G$; $\emptyset \nvdash \ConR \limp
  \neg\neg\forall q.\ \neg\Chk(q,\bquo G)$; Thm.~\ref{thm:br-g1} &
  disjunction property, Thm.~\ref{thm:d-g1}, \textsf{ConRGA\_unrefutable};
  for the second statement, \textsf{goedelG\_code\_unprovable} in place of Thm.~\ref{thm:d-g1} \\
\textsf{goedelG\_code\_unprovable} &
  $\emptyset \nvdash \forall q.\ \neg\Chk(q,\bquo G)$ &
  decided split, transport, $\forall$-introduction gives $G$; Thm.~\ref{thm:d-g1} \\
\textsf{goedelG\_instances\_refuted}, \textsf{goedelG\_omega\_incomplete},
\textsf{goedelG\_true\_orec} &
  for every $n$, $\emptyset \vdash \neg\Chk(\tofn n, S_G)$, but
  $\emptyset \nvdash G$; the classical $\omega$-evaluation gives $G$ the
  value true; [\citen{conrga-paper}, theorem~\textsf{goedelW\_omega\_incomplete}] &
  decidedness at a numeral, disjunction property,
  \textsf{goedelG\_no\_accepted\_code}; Thm.~\ref{thm:d-g1} \\
\textsf{ConRGA\_instances\_refuted}, \textsf{ConRGA\_omega\_incomplete},
\textsf{ConRGA\_true\_orec}, \textsf{uniform\_soundness\_instances} &
  the same for $\ConR$: for every $n$, $\emptyset \vdash
  \neg\Chk(\tofn n, \bquo{\z = \tofn 1})$, but $\emptyset \nvdash \ConR$;
  the $\omega$-evaluation gives $\ConR$ the value true; every numeral
  instance of $U$ ([\citen{conrga-paper}, theorem~\textsf{rga\_no\_uniform\_soundness}]) is a theorem &
  decidedness at a numeral, disjunction property, genuineness and
  consistency; \textsf{rga\_goedel2}; for $U$,
  $\lor$-introduction from $\ConR$'s instance \\
\textsf{indep\_allapp\_no\_value}, \textsf{goedelG\_no\_value},
\textsf{ConRGA\_no\_value}, \textsf{ConRGA\_not\_true},
\textsf{uniform\_soundness\_unrefutable},
\textsf{uniform\_soundness\_no\_value} &
  a well-formed universal that \rga neither proves nor refutes has no
  value; so $G$, $\ConR$ and the sentence $U$ of
  [\citen{conrga-paper}, theorem~\textsf{rga\_no\_uniform\_soundness}] have no value, and $\emptyset \nvdash
  \Tru(\bquo{\ConR})$; $\emptyset \nvdash \neg U$ &
  universals take only the values true and false; open completeness;
  Thm.~\ref{thm:d-g1}, Thm.~\ref{thm:br-con}; adequacy; witness property \\
\textsf{rga\_omega\_consistent} &
  no $u$ with $\emptyset \vdash \neg(\cA\cdot u)$ and $\emptyset \vdash
  u\cdot\tofn n$ for every $n$; Thm.~\ref{thm:g-g1}(3) &
  witness property, \textsf{rga\_consistent} \\
\textsf{no\_D1\_schema\_goedelG}, \textsf{no\_self\_imp\_goedelG} &
  $\emptyset \nvdash G \limp \PrC(\bquo G)$; $\emptyset \nvdash G
  \limp G$; Thms.~\ref{thm:br-d2d3}, \ref{thm:br-mech}(3) &
  disjunction property, Thm.~\ref{thm:d-g1}, \textsf{rga\_pr\_reflect} \\
\textsf{undecided\_no\_self\_imp}, \textsf{undecided\_no\_self\_iff},
\textsf{no\_self\_imp\_ConRGA} &
  if $\emptyset \nvdash \delta$ and $\emptyset \nvdash \neg\delta$, then
  $\emptyset \nvdash \delta \limp \delta$ and $\emptyset \nvdash \delta
  \liff \delta$; so $\emptyset \nvdash \ConR \limp \ConR$;
  Thm.~\ref{thm:br-mech}(3) &
  disjunction property at $\neg\delta \lor \delta$;
  \textsf{ConRGA\_unrefutable}, \textsf{rga\_goedel2} \\
\textsf{goedelSigmaSearch\_not\_provably\_zero} &
  $\emptyset \nvdash \forall v.\ \sigma_{\mathrm s}\cdot v = \z$ &
  decided split, D3G at $\sigma_{\mathrm s}$, \textsf{iChk\_neD};
  \textsf{goedelG\_code\_unprovable} \\
\textsf{no\_internal\_afree\_normalization\_0eq1} &
  $\emptyset \nvdash \forall p.\ \Chk(p, \bquo{\z = \tofn 1}) \limp
  \neg\forall q.\ \neg\ChkQf(q, \bquo{\z = \tofn 1})$; [\citen{conrga-paper}, \S\,``The $\forall$-free fragment''] &
  disjunct refuted by [\citen{conrga-paper}, theorem~\textsf{rga\_con\_qf\_int}], $\forall$-introduction;
  \textsf{rga\_goedel2} \\
\textsf{rga\_con\_int\_of\_gentzen}, \textsf{no\_gentzen\_pair},
\textsf{gred\_nonvacuous} (\texttt{ConRGA/Gentzen\_Shape}) &
  $\mathrm{GRed}(r,\mathit{od},\prec)$ and $\mathrm{GTI}$ imply
  $\emptyset \vdash \ConR$; so $\mathrm{GRed}$ implies $\neg\mathrm{GTI}$;
  some $\mathrm{GRed}$ holds &
  decided split, one internal induction, [\citen{conrga-paper}, theorem~\textsf{rga\_con\_frag\_int}]; \textsf{rga\_goedel2} \\
\textsf{rga\_con\_int\_of\_INFF}, \textsf{rga\_not\_INF} (\texttt{ConRGA/Gentzen\_Shape}) &
  $\emptyset \vdash \mathit{INF}$ implies $\emptyset \vdash \ConR$;
  $\emptyset \nvdash \mathit{INF}$ &
  decided split, [\citen{conrga-paper}, theorem~\textsf{rga\_con\_frag\_int}]; \textsf{rga\_goedel2} \\
\textsf{con\_search\_unrefutable} (\texttt{RGA\_Liar\_Div/Divergence}) &
  $\emptyset \nvdash \neg\Rcd_1(\bquo{\mathit{conSearch}})$ &
  \textsf{rga\_con\_int\_of\_con\_search}; \textsf{rga\_goedel2} \\
\textsf{iRes\_underivable} (\texttt{ConRGA/Guarded\_Residual}) &
  $\vz\jnat \nvdash \Res(\vz)$; [\citen{conrga-paper}, corollary~\textsf{iRes\_underivable}] &
  cut with \textsf{rga\_con\_under\_res}; \textsf{rga\_goedel2} \\
\textsf{fwP3\_iRes\_underivable} (\texttt{ConOGA/ConOGA\_Residual\_Frag}) &
  $\emptyset \nvdash \Res(\bquo{P_3})$, at a proof in $\Frag$, in $\LvC$
  and in every $\LvC_N$; [\citen{conrga-paper}, theorem~\textsf{fwP3\_iRes\_underivable}] &
  disjunction property at $\mathit{conSearch}$,
  \textsf{conSearch\_no\_value}; \textsf{con\_search\_unrefutable} \\
\end{longtable}
\endgroup

\paragraph{The Gentzen shape.}
$\mathrm{GRed}(r, \mathit{od}, \prec)$ says that $r$, $\mathit{od}$ and
the term behind $\prec$ are closed and well-formed, that \rga proves $r$
total, and that \rga proves $r$ correct (\textsf{gzRedF}): for every
accepted proof $p$ of a $\forall$-free conclusion $c$, either $p$ lies in
the fragment of [\citen{conrga-paper}, \S\,``The fragment theorem''], or $r\cdot p$ is an accepted proof
of $c$ with $\mathit{od}(r\cdot p) \prec \mathit{od}(p)$.  $\mathrm{GTI}$
says that \rga proves $\forall p.\ \neg\forall n.\
\mathit{od}(\mathit{It}_p(\suc n)) \prec \mathit{od}(\mathit{It}_p(n))$,
where $\mathit{It}_p$ iterates $r$ from $p$ and stops at the fragment
(\textsf{gzIt}).  \textsf{rga\_con\_int\_of\_gentzen} uses the fragment
theorem ([\citen{conrga-paper}, theorem~\textsf{rga\_con\_frag\_int}]), and \textsf{no\_gentzen\_pair}
says: for no reduction that \rga proves correct does \rga also prove that
the reduction's descent stops.  GTI is a single no-descent sentence about
one sequence, with no predicate parameter and no progressiveness premise,
and its derivation from a transfinite-induction schema is not
formalized.  $\mathrm{GRed}$ alone is inhabited, by the identity
reduction under the total relation (\textsf{gred\_identity}), so all
well-foundedness content lies in GTI.  $\mathit{INF}$ says that every
accepted proof of a $\forall$-free conclusion has an accepted fragment
proof of that conclusion.

\paragraph{The proof search for $\z = \tofn 1$.}
$\mathit{conSearch}$ runs through the proof codes and stops at the first
one that the checker accepts for $\z = \tofn 1$.  With the value record
$\Rcd_1$ of [\citen{conrga-paper}, definition~\textsf{iRec}], the table's statement says that
\rga does not refute ``the search returns true''.
\endgroup

\section{The locale layers and the interpretations}
\label{app:locale}

This appendix gives the assumptions of the G\"odel conditions as the
development states them, the layer of each locale theorem, the
witnesses at \rga, the 62 facts of the \rga theories that are
instances of locale facts, and the build-time checks on those
instances.  The calculus is in
\texttt{RGA/RGA\_GCalc} (the locales \textsf{gcalc},
\textsf{gcalc\_strip} and \textsf{gstep}, with the derived-rule kit)
and \texttt{GCalc\_Core/GCalc\_Core} (the rule list, the embedding and
the $\Sigma$ kit), the layers in \texttt{G2\_Locale/G2\_Locale}, and the \rga instance in
the theories \texttt{Goedel\_Interp}, \texttt{Goedel\_Transformers},
\texttt{Entailment\_Goedel}, \texttt{Goedel\_Cut} and
\texttt{Goedel\_Conditions} of the session \texttt{RGA\_Goedel}.

\subsection{The calculus}
\label{app:locale-calc}

The locale \textsf{gcalc} extends \textsf{rga\_syntax}, which has no
parameters and fixes \rga's ten operators, their codes, substitution,
well-formedness and the head step.  \textsf{gcalc} fixes one
parameter, the relation $\Gamma \vdash c$ on finite sets of terms and
terms, and assumes 43 statements: \textsf{g\_wff} ($\Gamma \vdash c$
gives $\wf(c)$) and the 42 rules of \rga's derivability, each named
\textsf{ir\_}$X$ after \rga's rule $X$ (\texttt{RGA\_GCalc.thy}:481--529;
the list \textsf{g\_rules}, \texttt{GCalc\_Core.thy}:40--45).  Each rule
is \rga's rule with \rga's derivability replaced by $\vdash$; for
instance \textsf{ir\_negE} says that $\Gamma \vdash \neg a$, $\Gamma
\vdash a$ and $\wf(c)$ give $\Gamma \vdash c$.  The locale
\textsf{gcalc\_strip} adds \textsf{ir\_strip}: $\Gamma \vdash c$ gives
$\Gamma' \vdash c$ for the set $\Gamma'$ of the well-formed members of
$\Gamma$ (\texttt{RGA\_GCalc.thy}:937--938).  \textsf{gcalc\_strip} thus has 44 named
assumptions.  \textsf{gcalc\_strip} proves the $\Sigma$ rules of
Theorem~\ref{thm:lc-embed}(2) and Corollary~\ref{cor:rc-contrapos}
(\texttt{GCalc\_Core.thy}:203--381).

\subsection{The conditions as the development states them}
\label{app:locale-cond}

In the locale \textsf{g2\_syn} the checker is a parameter
$\gchk$ of HOL type $\mathit{trm} \Rightarrow \mathit{trm} \Rightarrow
\mathit{trm}$.  The provability formula is an input abbreviation,
$\gPr(c) = \neg\forall q.\,\neg\gchk(q, \lft c)$, with $\lft c$ the
lift of $c$ past the bound variable; at a code numeral the lift is
the identity, so $\gPr(\bquo\varphi)$ is the formula of
Definition~\ref{def:lc-sentences}.  The eleven named assumptions of
the layers are \textsf{chk\_subf}, \textsf{chk\_wft} and
\textsf{chk\_dec} (S1 to S3, in \textsf{g2\_syn}), \textsf{nec\_tr},
\textsf{d3\_tr} and \textsf{expl\_tr} (\textsf{g2\_ent}),
\textsf{chk\_acc} (D1, \textsf{g2\_acc}), \textsf{con0}
(\textsf{g2\_con}), \textsf{exf\_tr} (\textsf{g2\_exf}),
\textsf{d2\_tr} (\textsf{g2\_d2}) and \textsf{nec2\_tr}
(\textsf{g2\_nec2}).  They state the conditions of
Definition~\ref{def:lc-cond} and of the optional layers, with S2 as
the equivalence $\wf(\gchk(p, c)) \Leftrightarrow \wf(p) \wedge
\wf(c)$ and with HOL's implication and conjunction throughout.  Each
transformer condition has one form: a HOL implication from the
condition's HOL premise (the entailment for NEC and NEC2,
well-formedness for EXF and D2, none for D3 and EXPL) to an
existential over a HOL function $T$ on terms, followed by a universal
over the context $\Delta$ and the proof terms.  Inside, the premises
are $\Delta \vdash u\jnat$, $\Delta \vdash v\jnat$ and the stated
acceptances, and the conclusion is the conjunction of $\Delta \vdash
T(v)\jnat$ (or $T(u, v)\jnat$) and the stated acceptance of the
output.  For instance, \textsf{nec\_tr} reads: $\{A\} \vdash B$
$\Rightarrow$ there is $T$ such that for all $\Delta$ and $v$:
$\Delta \vdash v\jnat$ and $\Delta \vdash \gchk(v, \bquo{A})$
$\Rightarrow$ $\Delta \vdash T(v)\jnat$ and $\Delta \vdash
\gchk(T(v), \bquo{B})$.
The locale heads are at \texttt{G2\_Locale.thy}:85 (\textsf{g2\_syn}),
:515 (\textsf{g2\_ent}), :668 (\textsf{g2\_acc}), :700
(\textsf{g2\_con}), :756 (\textsf{g2\_exf}), :887 (\textsf{g2\_d2}),
:980 (\textsf{g2\_nec2}), :1000 (\textsf{g2\_exf\_acc}) and :1062
(\textsf{g2\_exf\_con}).  The sentences of
Definition~\ref{def:lc-sentences} are \textsf{ConT}, \textsf{goedelDT},
\textsf{goedelST} and \textsf{goedelGT} (:100--110), over the copy
\textsf{g2DiagSub} of $\Diag$ (:61), whose values
\textsf{g2DiagSub\_closed} (:64) gives as an \rga derivation, written
here with $\Diag$ for the copy:
$\emptyset \vdash \Diag(\bquo{t}, \tofn k) = \bquo{\ssub{\vz}{\tofn
k}{t}}$ for well-formed $t$.  The witness forms of
Theorem~\ref{thm:lc-witness}(1) are \textsf{EXPLT} and \textsf{D3GT}
(:380--395).

\subsection{The layer of each theorem}
\label{app:locale-map}

Table~\ref{tab:lc-map} gives the locale of each theorem of
\S\ref{sec:locale-route} and \S\ref{sec:locale-thm} and of the other
locale theorems that the \rga instances use, from the calculus up.  The structure
conditions of Theorem~\ref{thm:lc-bs}(1), as
\textsf{bsT\_structure\_conditions} states them with $x \vdash y$ for
$\{x\} \vdash y$, are: $x \vdash x$ for well-formed $x$;
$x \vdash y$ and $y \vdash z$ give $x \vdash z$; $x \vdash y$ gives
$\gPr(\bquo{x}) \vdash \gPr(\bquo{y})$; $\z = \z \vdash
\gPr(\bquo{\neg\neg(\z = \z)})$; $x \vdash \gPr(\bquo{y})$ and $x
\vdash \gPr(\bquo{\neg y})$ give $x \vdash \gPr(\bquo{\neg(\z =
\z)})$; $\gPr(\bquo{\neg x}) \vdash
\gPr(\bquo{\gPr(\bquo{\neg x})})$; and $x \vdash \z = \z$ for
well-formed $x$.  The development calls these the structure
conditions C1 (its monotone half) to C5; their proofs use NEC (third
item), D1 (fourth), EXPL and EXF (fifth) and D3 (sixth).

\begin{table}[tbp]
\centering
\small
\begin{tabular}{@{}l>{\raggedright\arraybackslash}p{0.62\linewidth}l@{}}
\hline
locale & theorems & here \\
\hline
\textsf{gcalc\_strip} & the $\Sigma$ rules (\textsf{sigma\_intro} \ldots\ \textsf{dec\_contrapos}), \textsf{ir\_ent\_trans}, \textsf{ir\_ent\_trans2} & Thm.~\ref{thm:lc-embed}(2) \\
\textsf{g2\_syn} & \textsf{goedelST\_eq}, \textsf{goedelT\_diag\_ent}, \textsf{goedelT\_diag\_ent\_all}, \textsf{goedelT\_diag\_neg\_ent}, \textsf{reflectT\_ent\_iff\_con}, \textsf{reflectT\_ent\_iff\_con\_ctx}, \textsf{conT\_reflect\_0eq1\_ent}, \textsf{goedelGT\_of\_expl\_d3g\_ctx}, \textsf{goedelGT\_of\_expl\_d3g\_con}, \textsf{prcT\_sentence\_iff}, \textsf{loebST\_eq}, \textsf{neg\_goedelGT\_sigma}, \textsf{neg\_goedelGT\_dec}, and the transport, lifting, closedness and L\"ob-diagonal lemmas & Thm.~\ref{thm:lc-diag} \\
\textsf{g2\_ent} & \textsf{ent\_necT}, \textsf{ent\_nec\_contraT}, \textsf{D3\_entT}, \textsf{D3G\_entT}, \textsf{expl\_entT}, \textsf{not\_goedelGT\_not\_con\_ent}, \textsf{goedelGT\_under\_ConT\_textbook}, \textsf{con\_notprc\_goedelGT\_ent}, \textsf{goedel2\_internal\_contra}, \textsf{goedel2\_internal}, \textsf{goedel2\_sentenceT}, \textsf{goedel1\_sentenceT}, \textsf{goedel1\_sentence\_diagT}, \textsf{bsT\_C1\_anti\_sigma}, \textsf{bsT\_C1\_anti\_dec}, \textsf{bsT\_jeroslow\_fp} & Thms.~\ref{thm:lc-ent}, \ref{thm:lc-textbook} \\
\textsf{g2\_acc} & \textsf{prcT\_of\_derivable}, \textsf{goedelGT\_refutes} & Thm.~\ref{thm:lc-g2}(1) \\
\textsf{g2\_con} & \textsf{goedel1}, \textsf{goedel2\_textbook}, \textsf{goedel2T\_of\_expl\_d3g}, \textsf{no\_loeb\_hypothesis\_0eq1T} & Thms.~\ref{thm:lc-g2}, \ref{thm:lc-witness} \\
\textsf{g2\_exf} & \textsf{exf\_entT}, \textsf{notcon\_prc\_con\_entT}, \textsf{goedelGT\_con\_ent}, \textsf{goedelGT\_con\_interderivable}, \textsf{notcon\_prc\_con\_interderivableT}, \textsf{ConT\_of\_code\_unprovable}, \textsf{bsT\_theorem1}, \textsf{bsT\_goedel2}, \textsf{bsT\_incon\_con} & Thm.~\ref{thm:lc-exf}(1, 2) \\
\textsf{g2\_d2} & \textsf{D2\_entT}, \textsf{formalized\_loeb\_entT}, \textsf{formalized\_loeb\_0eq1\_entT}, \textsf{formalized\_loeb\_sentenceT} & Thm.~\ref{thm:lc-loeb} \\
\textsf{g2\_nec2} & \textsf{ent\_nec2T} & \\
\textsf{g2\_exf\_acc} & \textsf{bsT\_structure\_conditions}, \textsf{bsT\_theorem1\_i\_rule} & Thm.~\ref{thm:lc-bs}(1) \\
\textsf{g2\_exf\_con} & \textsf{bsT\_goedel2\_irrefutable}, \textsf{no\_code\_unprovableT} & Thms.~\ref{thm:lc-exf}(3), \ref{thm:lc-bs}(2) \\
\hline
\end{tabular}
\caption{The locale in which each theorem is proved
  (\texttt{G2\_Locale}, \texttt{RGA\_GCalc}, \texttt{GCalc\_Core}).  A theorem carries all
  assumptions of its locale as premises, and its proof uses at most
  those.}
\label{tab:lc-map}
\end{table}

\subsection{The witnesses at \rga}
\label{app:locale-witness}

Table~\ref{tab:lc-witness} lists the witnesses by which \rga with
the tree checker interprets each layer (\texttt{Goedel\_Interp.thy}:146--264).
Each transformer of a condition is a HOL function of the condition's
instance: NEC's of the entailment, through a proof $Q$ of it that HOL
obtains at \rga; D3's of $\varphi$; ex falso's of
$\varphi$; D2's of $\varphi$ and $\psi$.  The D1 witness,
\textsf{iChk\_proves\_sf}, rests on a closed-evaluation lemma of the
theory \texttt{Goedel\_D2}; that lemma evaluates the checker at a
closed proof code and is not the D2 transformer law.  The second
theorem's instance uses no D2 transformer law.  \iga's witnesses are
the subject of \cref{sec:iga-goedel}.

\begin{table}[tbp]
\centering
\small
\begin{tabular}{@{}l>{\raggedright\arraybackslash}p{0.33\linewidth}>{\raggedright\arraybackslash}p{0.27\linewidth}l@{}}
\hline
condition & witness & law (theory) & here \\
\hline
S1, S2, S3 & the checker's laws at abstract arguments & \textsf{subf\_iChk}, \textsf{wft\_iChk}, \textsf{hq\_iChk} & \S\ref{sec:recall-chk}, \S\ref{app:goedel-chk} \\
NEC & $T(v) = c_Q\cdot v\cdot v$ for a proof $Q$ of $\{A\} \vdash B$ & \textsf{cut1M\_chk} (\texttt{Goedel\_Transformers}) & Thm.~\ref{thm:g-cut} \\
D3 & the decided-D3 transformer at $\bquo\varphi$, applied to $v$ & \textsf{iChk\_d3SigmaAt} (\texttt{Goedel\_D3}) & Thm.~\ref{thm:g-d3-any} \\
EXPL & $T(u, v) = e\cdot u\cdot v$ & \textsf{goedelE\_chk} (\texttt{Goedel\_II}) & Thm.~\ref{thm:g-expl} \\
D1 & closed evaluation of the checker at a proof's code & \textsf{iChk\_proves\_sf} (\texttt{Goedel\_Transformers}) & \S\ref{app:goedel-ent} \\
CON & \rga's consistency & \textsf{rga\_no\_proof\_0eq1} & \\
EXF & the ex falso transformer at $\varphi$ & \textsf{exfIdx\_chk\_sf} (\texttt{Goedel\_Transformers}) & \S\ref{app:goedel-cons} \\
D2 & the decided-D2 transformer at $\varphi$, $\psi$ & \textsf{d2M\_chk} (\texttt{Goedel\_D2}) & Thm.~\ref{thm:g-d2} \\
NEC2 & the two-hypothesis cut transformer & \textsf{cut2M\_chk} (\texttt{Goedel\_Transformers}) & Thm.~\ref{thm:g-cut} \\
\hline
\end{tabular}
\caption{The witnesses of the G\"odel conditions at \rga's tree
  checker, from \textsf{g2\_syn\_rga} to \textsf{g2\_exf\_con\_rga}.}
\label{tab:lc-witness}
\end{table}

The rewrites \textsf{g2\_rws} (\texttt{Goedel\_Interp.thy}:269--297)
turn the locale's constants into \rga's at the tree checker:
\textsf{g2DiagSub} into $\Diag$, \textsf{ConT} into $\ConR$,
\textsf{goedelDT}, \textsf{goedelST} and \textsf{goedelGT} into $D_G$,
$S_G$ and $G$, and the witness forms into \rga's EXPL and D3G
predicates; the theory \texttt{Goedel\_Cut} adds the rewrites of the
L\"ob sentence.  At the tree checker, $\gPr$ is literally \rga's
$\PrC$, since $\gPr$ is an input abbreviation.  The nine
interpretations \textsf{g2sR}, \textsf{g2eR}, \textsf{g2aR},
\textsf{g2fR}, \textsf{g2dR}, \textsf{g2nR}, \textsf{g2faR},
\textsf{g2R} and \textsf{g2fcR} (:306--362) interpret the layers one
at a time, each after the layers it is built on: \textsf{g2\_syn},
\textsf{g2\_ent}, \textsf{g2\_acc}, \textsf{g2\_exf},
\textsf{g2\_d2}, \textsf{g2\_nec2}, \textsf{g2\_exf\_acc},
\textsf{g2\_con}, \textsf{g2\_exf\_con}.

\subsection{The 62 instance facts}
\label{app:locale-inst}

The list below gives every fact of the three theories
\texttt{Entailment\_Goedel} (20 facts), \texttt{Goedel\_Cut} (31) and
\texttt{Goedel\_Conditions} (11) that the development defines as an
instance, \textsf{lemmas} $X$ \textsf{=} $\ell.Y$, of a locale fact
$Y$ through the interpretation $\ell$, grouped by interpretation.
Within a group, the facts of \texttt{Entailment\_Goedel} come first,
then those of \texttt{Goedel\_Cut}, then those of
\texttt{Goedel\_Conditions}, with semicolons between the theories.  The
locale fact's name $Y$ is the \rga name $X$ with \textsf{ConRGA},
\textsf{goedelG}, \textsf{iPrc}, \textsf{iChk} and \textsf{bs\_} read
as \textsf{ConT}, \textsf{goedelGT}, \textsf{iPrcT}, \textsf{chk} and
\textsf{bsT\_}, or, when none of these occurs in $X$, $X$ with a final
\textsf{T}; where the rule does not give $Y$, the list prints $Y$ in
parentheses.  A dagger marks the roots of the decider check of
\S\ref{app:locale-checks}.  Two \rga facts,
\textsf{prc\_con\_ent\_bot} and
\textsf{rga\_goedel2\_internal\_contra}, are instances of the same
locale fact.  Separately, \texttt{Goedel\_Interp} defines 31 locale
instances \textsf{g2x\_}$X$, each with the proposition of the \rga
theorem $X$, which keeps its own proof (a build-time check compares
the two), among them \textsf{g2x\_rga\_goedel2}, the instance of
\textsf{goedel2\_textbook} beside \textsf{rga\_goedel2}.

{\small
\begin{description}\setlength{\itemsep}{1pt}
\item[\textsf{g2sR} (27)] \textsf{lift\_\allowbreak num\_\allowbreak code},
\textsf{iPrc\_\allowbreak num}, \textsf{iPrc\_\allowbreak cEq01}, \textsf{bsub\_\allowbreak
iChk}, \textsf{lift\_\allowbreak iPrc\_\allowbreak closed}, \textsf{iPrc\_\allowbreak intro},
\textsf{iPrc\_\allowbreak dec}, \textsf{iPrc\_\allowbreak transfer}, \textsf{iPrc\_\allowbreak
transfer2}, \textsf{all\_\allowbreak negchk\_\allowbreak transport} (\textsf{all\_\allowbreak
negchkT\_\allowbreak transport}), \textsf{goedel\_\allowbreak diag\_\allowbreak
ent\_\allowbreak all} (\textsf{goedelT\_\allowbreak diag\_\allowbreak ent\_\allowbreak all}),
\textsf{goedel\_\allowbreak diag\_\allowbreak ent} (\textsf{goedelT\_\allowbreak
diag\_\allowbreak ent}), \textsf{reflect\_\allowbreak ent\_\allowbreak iff\_\allowbreak con}
(\textsf{reflectT\_\allowbreak ent\_\allowbreak iff\_\allowbreak con}),
\textsf{reflect\_\allowbreak ent\_\allowbreak iff\_\allowbreak con\_\allowbreak ctx}$^\dagger$
(\textsf{reflectT\_\allowbreak ent\_\allowbreak iff\_\allowbreak con\_\allowbreak ctx}),
\textsf{con\_\allowbreak reflect\_\allowbreak 0eq1\_\allowbreak ent}$^\dagger$
(\textsf{conT\_\allowbreak reflect\_\allowbreak 0eq1\_\allowbreak ent});
\textsf{ent\_\allowbreak trans} (\textsf{ir\_\allowbreak ent\_\allowbreak trans}),
\textsf{ent\_\allowbreak trans2} (\textsf{ir\_\allowbreak ent\_\allowbreak trans2}),
\textsf{closed\_\allowbreak loebS}, \textsf{wft\_\allowbreak loebD}, \textsf{loeb\_\allowbreak
diag} (\textsf{loebT\_\allowbreak diag}), \textsf{loebS\_\allowbreak eq}
(\textsf{loebST\_\allowbreak eq}), \textsf{iPrc\_\allowbreak eq\_\allowbreak transfer},
\textsf{goedel\_\allowbreak diag\_\allowbreak neg\_\allowbreak ent}
(\textsf{goedelT\_\allowbreak diag\_\allowbreak neg\_\allowbreak ent}),
\textsf{lift\_\allowbreak not\_\allowbreak ConRGA}, \textsf{prc\_\allowbreak
sentence\_\allowbreak iff} (\textsf{prcT\_\allowbreak sentence\_\allowbreak iff});
\textsf{neg\_\allowbreak goedelG\_\allowbreak sigma}, \textsf{neg\_\allowbreak
goedelG\_\allowbreak dec}.
\item[\textsf{g2eR} (17)] \textsf{D3\_\allowbreak ent}, \textsf{D3G\_\allowbreak
ent}$^\dagger$, \textsf{expl\_\allowbreak ent}$^\dagger$, \textsf{con\_\allowbreak
notprc\_\allowbreak goedelG\_\allowbreak ent}$^\dagger$; \textsf{ent\_\allowbreak
nec}$^\dagger$, \textsf{ent\_\allowbreak nec\_\allowbreak contra}, \textsf{prc\_\allowbreak
con\_\allowbreak ent\_\allowbreak bot}$^\dagger$ (\textsf{goedel2\_\allowbreak
internal\_\allowbreak contra}), \textsf{rga\_\allowbreak goedel2\_\allowbreak
internal}$^\dagger$ (\textsf{goedel2\_\allowbreak internal}), \textsf{not\_\allowbreak
goedelG\_\allowbreak not\_\allowbreak con\_\allowbreak ent}, \textsf{rga\_\allowbreak
goedel2\_\allowbreak internal\_\allowbreak contra}$^\dagger$ (\textsf{goedel2\_\allowbreak
internal\_\allowbreak contra}), \textsf{rga\_\allowbreak goedel2\_\allowbreak
sentence}$^\dagger$ (\textsf{goedel2\_\allowbreak sentenceT}), \textsf{goedel1\_\allowbreak
sentence}, \textsf{goedel1\_\allowbreak sentence\_\allowbreak diag},
\textsf{goedelG\_\allowbreak under\_\allowbreak ConRGA\_\allowbreak textbook}$^\dagger$;
\textsf{bs\_\allowbreak C1\_\allowbreak anti\_\allowbreak sigma}, \textsf{bs\_\allowbreak
C1\_\allowbreak anti\_\allowbreak dec}, \textsf{bs\_\allowbreak jeroslow\_\allowbreak fp}.
\item[\textsf{g2fR} (8)] \textsf{exf\_\allowbreak ent}; \textsf{notcon\_\allowbreak
prc\_\allowbreak con\_\allowbreak ent}, \textsf{goedelG\_\allowbreak con\_\allowbreak ent},
\textsf{goedelG\_\allowbreak con\_\allowbreak interderivable}$^\dagger$,
\textsf{notcon\_\allowbreak prc\_\allowbreak con\_\allowbreak interderivable}$^\dagger$;
\textsf{bs\_\allowbreak theorem1}, \textsf{bs\_\allowbreak goedel2}$^\dagger$,
\textsf{bs\_\allowbreak incon\_\allowbreak con}.
\item[\textsf{g2dR} (4)] \textsf{D2\_\allowbreak ent}$^\dagger$,
\textsf{formalized\_\allowbreak loeb\_\allowbreak ent}$^\dagger$,
\textsf{formalized\_\allowbreak loeb\_\allowbreak 0eq1\_\allowbreak ent},
\textsf{formalized\_\allowbreak loeb\_\allowbreak sentence}.
\item[\textsf{g2aR} (1)] \textsf{prc\_\allowbreak of\_\allowbreak derivable}
(\textsf{prcT\_\allowbreak of\_\allowbreak derivable}).
\item[\textsf{g2nR} (1)] \textsf{ent\_\allowbreak nec2}$^\dagger$.
\item[\textsf{g2faR} (2)] \textsf{bs\_\allowbreak structure\_\allowbreak conditions},
\textsf{bs\_\allowbreak theorem1\_\allowbreak i\_\allowbreak rule}.
\item[\textsf{g2R} (1)] \textsf{rga\_\allowbreak goedel2\_\allowbreak textbook}$^\dagger$
(\textsf{goedel2\_\allowbreak textbook}).
\item[\textsf{g2fcR} (1)] \textsf{bs\_\allowbreak goedel2\_\allowbreak irrefutable}.
\end{description}
}

\subsection{The build-time checks}
\label{app:locale-checks}

The ML checks below run in every build of the \rga instance and fail
the build when they fail.  These checks are facts about the artifact,
not mathematics.

\paragraph{Statement identity.}
Each theory compares each instance fact with its locale fact,
proposition against proposition.  Since the development defines each
instance fact as its locale fact, these comparisons hold by
construction; they
guard against a wrong name in a definition.  The substantive
comparisons are two.  First, a \emph{ground-truth table}
(\texttt{Goedel\_Conditions.thy}:211--389) records 130 named facts of
the three theories (46 of \texttt{Entailment\_Goedel}, 68 of
\texttt{Goedel\_Cut}, 16 of \texttt{Goedel\_Conditions}) with a
fingerprint of each fact's proposition as the facts stood before they
were redefined as instances (commit \texttt{07ef275f}).
A fingerprint is a SHA1 hash of the printed proposition, types
included, with every bound-variable name blanked, so that equal
fingerprints mean propositions equal up to the names of bound
variables.  For each row the check requires that the fact exist under
its full name (19 facts now in \texttt{Goedel\_Transformers} are
mapped to that theory), that the fact's short name resolve to that
full name, and that the fact's fingerprint equal the recorded one (six
constants now in \texttt{Goedel\_Transformers} are read under their
earlier names).  All 62
instance facts are rows of the table.  A control feeds the check the
fact \textsf{ent\_nec} with the fingerprint of \textsf{ent\_nec2},
and the build fails unless the check rejects it.  The build does not
check that the table lists every named fact of the three theories at
that commit, and the check does not see facts added since: the four
rewrite facts for the L\"ob sentence, \textsf{g2\_rw\_loebD},
\textsf{g2\_rw\_loebS}, \textsf{g2\_rw\_loebL} and
\textsf{g2\_rws\_loeb}, are not rows.  Second, two pairs must have equal
propositions (:393--397), each an instance against a fact that keeps
its own \rga proof: \textsf{rga\_goedel2\_textbook} and
\textsf{rga\_goedel2}, and \textsf{goedelG\_under\_ConRGA\_textbook}
and \textsf{goedelG\_under\_ConRGA}.  The same check's third pair,
\textsf{rga\_goedel2\_internal\_contra} and
\textsf{prc\_con\_ent\_bot}, holds by construction: both are
instances of \textsf{goedel2\_internal\_contra}.

\paragraph{No oracle.}
No listed fact depends on an oracle (\texttt{Goedel\_Interp.thy}:460;
\texttt{Goedel\_Conditions.thy}:434).

\paragraph{Soundness footprint.}
An audit walks the dependency cone of each listed fact, the
transitive closure of the theorems its proof cites, and checks it
against nine facts of \rga: the soundness theorems
\textsf{prf\_owsjt} and \textsf{prf\_oomjt}, the consistency theorems
\textsf{rga\_consistent} and \textsf{rga\_no\_proof\_0eq1},
$\omega$-soundness \textsf{rga\_omega\_sound}, the checker's
faithfulness \textsf{iChk\_faithful}, and the truth readings
\textsf{rga\_provable\_wosr}, \textsf{rga\_val\_iff} and
\textsf{rga\_open\_iff}.  The cones of 26 \textsf{g2x\_} facts and of
the 60 instance facts outside the two consistency layers reach none
of the nine (\texttt{Goedel\_Interp.thy}:441--463;
\texttt{Goedel\_Conditions.thy}:411--436).  For the consistency-layer
facts, five \textsf{g2x\_} facts and \textsf{rga\_goedel2\_textbook}
and \textsf{bs\_goedel2\_irrefutable}, the audit requires that the
cone reach at least one of \textsf{rga\_no\_proof\_0eq1},
\textsf{rga\_consistent} and \textsf{prf\_owsjt}, and none of the
other six.  Three further checks, run in
\texttt{Goedel\_Conditions} (:437--457) on 24 roots (21 instance
facts and three \textsf{g2x\_} facts), follow those of
\texttt{ConOGA\_Audit}.  No cone contains \textsf{prf\_owsjt} or
\textsf{prf\_oomjt} (with controls); no statement in a cone has the
shape of a soundness statement, apart from \rga's eight listed clauses
of the certificate semantics (\texttt{ConOGA\_Audit} lists thirteen
clauses, those eight and the five generic clauses of the same
evaluation stated for an arbitrary system's proofs as the certificates); and no theorem in a cone cites an induction
rule of a derivability predicate.

\paragraph{Decider-free statements.}
The decider check of the session \texttt{ConOGA}
(\texttt{ConOGA\_Audit.thy}:334--366;
[\citen{conrga-paper}, Appendix ``The build-time check'']) runs on 31 \rga facts.  It
fails the build if a theorem statement in the cone of any of them
mentions one of nine constants of the decider's verdicts, its truth
predicate, its search witness and the HOL-level decider
(\textsf{iWtrue}, \textsf{iPrDecTrueAt}, \textsf{iPrDecVal},
\textsf{goedelSigmaSearch}, \textsf{prDecV}, \textsf{prDecH},
\textsf{prDecVEnc}, \textsf{prDecHEnc}, \textsf{prDecTab}) or one of
five constants of the compiled decider's table (\textsf{iPrDecV},
\textsf{iPrDecH}, \textsf{iPrDecTab}, \textsf{iPrDecTabT},
\textsf{prDecTabIdx}), or if the cone of any of them contains one of
the theorems \textsf{D3G\_goedelSigmaSearch}, \textsf{D3Gdec\_holds}
and \textsf{iHasCert\_of\_iPrDecTrueAt}; controls that must be flagged
guard both checks, and each root's cone must have at least 1000
nodes.  The check reads the statements in the cone, not the
definitional closure of the constants that those statements mention.
Of the 31 roots, 18 are instance facts, of 17 distinct locale facts,
marked with a dagger in \S\ref{app:locale-inst}; the other 13 keep
\rga proofs.

An instance obtains its statement by instantiating the stored locale
theorem, so the instance's cone contains the locale theorem's cone,
together with the cones of the witnesses and of the rewrites.  A
clean instance therefore means that no statement in the cone of the
locale theorem mentions a forbidden constant and that none of the
three forbidden theorems is in that cone.  This consequence has three
limits.  No check prints it.  The consequence covers only the cones of
the 17 locale facts above, and so not, for example,
\textsf{bsT\_structure\_conditions}, \textsf{no\_code\_unprovableT} or
\textsf{goedel2T\_of\_expl\_d3g}.  The consequence says nothing about
\iga's instance.  The consequence is not vacuous: the session
\texttt{G2\_Locale} has the decider session as parent, and the
locale's diagonal cites two lemmas about codes of numerals from the
decider session \texttt{RGA\_Decider}, \textsf{pr2fi\_num\_codeIdx}
of \texttt{Decider\_Truth} and \textsf{num\_code\_correct} of
\texttt{RGA\_Decide\_PR}, whose statements the check finds clean.

\subsection{G\"odel II at \iga: laws, witnesses, re-exports and checks}
\label{app:iga-goedel}

This appendix completes Section~\ref{sec:iga-goedel}.  It states the
transformer laws of the optional layers, lists the witnesses with
which \iga meets each condition, lists the theorems that the nine
interpretations state for \iga, and describes the build checks.  As in
Section~\ref{sec:iga-goedel}, $\mathcal D$ is the table of a detector
list $L$, $\vdash$ is \rga's entailment and $\vdashI$ is \iga's.

\paragraph{The laws of the optional layers.}

\begin{thm}[\textsf{G\_\allowbreak exfIdx\_\allowbreak chk\_\allowbreak sf}, \textsf{G\_\allowbreak d2M\_\allowbreak chk}]
\label{thm:ig-laws-exf}
Let the table $\mathcal D$ of $L$ satisfy the shift law, and let $L$
contain \rga's list.  Let $\Gamma$ be a context and $a$ a term with
$\Gamma \vdash a\jnat$.
\begin{enumerate}
\item (Ex falso.)  Let $L$ contain the negation-elimination detector,
  and let $\varphi$ be well formed.  If $\Gamma \vdash \Chk_{\mathcal
  D}(a, \bquo{\z = \tofn 1})$, then $\Gamma \vdash \Chk_{\mathcal
  D}(x_\varphi\cdot a\cdot a, \bquo{\varphi})$.
\item (D2.)  Let $L$ contain the cut and weakening detectors, let
  $\varphi$ and $\psi$ be well formed, and let $b$ be a term with
  $\Gamma \vdash b\jnat$.  If $\Gamma \vdash \Chk_{\mathcal
  D}(a, \bquo{\varphi \limp \psi})$ and $\Gamma \vdash \Chk_{\mathcal
  D}(b, \bquo{\varphi})$, then $\Gamma \vdash \Chk_{\mathcal
  D}(d_{\varphi,\psi}\cdot a\cdot b, \bquo{\psi})$.
\end{enumerate}
Here $x_\varphi$ and $d_{\varphi,\psi}$ are built by the term formers
that the \rga instance of Section~\ref{sec:locale} uses.
\end{thm}

\noindent
The ex falso transformer puts one negation-elimination entry over the
append of $a$ and \rga's fixed proof of $\z \neq \tofn 1$, embedded
into the table.  The D2 transformer puts a weakening entry and two cut
entries over the append of a closed \rga tree, its core, and the two
inputs.  The core is valid for \rga, and its embedding into the table
is where the containment of \rga's list is needed.

\begin{thm}[\textsf{G\_\allowbreak cut2M\_\allowbreak chk}]
\label{thm:ig-laws-nec2}
Let the table $\mathcal D$ of $L$ satisfy the shift law and decide a
step relation $R$, let $L$ contain the cut and weakening detectors, and
let $Q$ be an $R$-proof of $\{A_1, A_2\} \vdash B$.  For every context
$\Gamma$ and all terms $a, b$: if $\Gamma$ derives $a\jnat$, $b\jnat$,
$\Chk_{\mathcal D}(a, \bquo{A_1})$ and $\Chk_{\mathcal D}(b,
\bquo{A_2})$, then $\Gamma \vdash \Chk_{\mathcal D}(c_Q\cdot a\cdot b,
\bquo{B})$.
\end{thm}

\noindent
The transformer is the two-hypothesis transformer of
Theorem~\ref{thm:g-cut}, with $Q$ a proof of $R$.  It puts a weakening
entry and two cut entries over the append of the code of $Q$ and the
two inputs.  Closed evaluation shows that the code of $Q$ is a valid
tree at the table whose head is $\{A_1, A_2\} \vdash B$ (the argument
of Theorem~\ref{thm:ig-closed}).  Theorems~\ref{thm:ig-laws-exf}
and~\ref{thm:ig-laws-nec2} are rules in every context of \rga.

\paragraph{The witnesses.}
Table~\ref{tab:ig-witnesses} lists, for each condition, the layer of
Section~\ref{sec:locale} that states it and the witness with which
\iga meets it.  Each witness is stated for \iga's entailment and the
checker $\ChkI$, and the laws of \rga reach \iga by cut transfer
(\S\ref{sec:ig-instance}).

\begin{table}[tbp]
\centering
\small
\setlength{\tabcolsep}{3pt}
\begin{tabular}{@{}>{\raggedright\arraybackslash}p{0.2\linewidth}>{\raggedright\arraybackslash}p{0.17\linewidth}>{\raggedright\arraybackslash}p{0.38\linewidth}>{\raggedright\arraybackslash}p{0.2\linewidth}@{}}
\hline
condition & layer & witness at \iga & law \\
\hline
calculus and the three structural conditions & \textsf{g2\_\allowbreak syn} & \iga's strip calculus; $\ChkI$ commutes with substitution and is well formed exactly at well-formed arguments; \rga's decidedness of the checker at a table, carried into \iga & Thm.~\ref{thm:ig-strip} \\
NEC & \textsf{g2\_\allowbreak ent} & $v \mapsto c_Q\cdot v\cdot v$, with $Q$ an \iga proof of $\{A\} \vdashI B$ & Thm.~\ref{thm:ig-laws}(1) \\
D3 & \textsf{g2\_\allowbreak ent} & $\sigma_{\mathcal D, S}$ at \iga's table and $S = \bquo{\varphi}$ & Thm.~\ref{thm:ig-d3} \\
EXPL & \textsf{g2\_\allowbreak ent} & $e$ & Thm.~\ref{thm:ig-laws}(2) \\
D1 & \textsf{g2\_\allowbreak acc} & the code $\bquo{P}$ of an \iga proof $P$ & Thm.~\ref{thm:ig-closed}(2) \\
consistency & \textsf{g2\_\allowbreak con} & $\emptyset \nvdashI \z = \tofn 1$ & \cite{iga-paper} \\
ex falso & \textsf{g2\_\allowbreak exf} & $v \mapsto x_\varphi\cdot v\cdot v$ & Thm.~\ref{thm:ig-laws-exf}(1) \\
D2 & \textsf{g2\_\allowbreak d2} & $d_{\varphi,\psi}$ & Thm.~\ref{thm:ig-laws-exf}(2) \\
NEC with two premises & \textsf{g2\_\allowbreak nec2} & $c_Q$, with $Q$ an \iga proof of $\{A_1, A_2\} \vdashI B$ & Thm.~\ref{thm:ig-laws-nec2} \\
\hline
\end{tabular}
\caption{The witnesses with which \iga and $\ChkI$ meet the G\"odel
  conditions (Theorem~\ref{thm:ig-witness}).  The layers
  \textsf{g2\_\allowbreak exf\_\allowbreak acc} and \textsf{g2\_\allowbreak exf\_\allowbreak con} combine layers of
  the table, and their witnesses combine the witnesses of those layers.}
\label{tab:ig-witnesses}
\end{table}

\paragraph{The re-exports.}
\label{app:iga-goedel-reexports}
Table~\ref{tab:ig-reexports} lists the twenty theorems that the
development states for \iga, each a theorem of one layer, interpreted
with $\ChkI$ for the checker and $\ConI$ for the consistency sentence.
In the table $\vdash$ is \iga's entailment, $\PrC$ is $\PrI$, $\Gamma$
is any context, and $G$ is $G_{\mathrm I}$.  The column ``reading''
follows the convention of Section~\ref{sec:locale}: a rule or a
non-derivability is classical, and an entailment or a sentence that
\iga derives is a derivation in \iga, written ``in \iga''.  The column ``checked'' marks the eight
theorems whose dependency cones the build checks against the soundness
theorems.  $\neg(\z = \z)$ plays the role of the falsum of Beklemishev
and Shamkanov's provability structures~\cite{beklemishev16abstract}, a
classical falsum, not \rga's divergent constant $\bot$.

\begin{table}[tbp]
\centering
\small
\setlength{\tabcolsep}{3pt}
\begin{tabular}{@{}>{\raggedright\arraybackslash}p{0.27\linewidth}>{\raggedright\arraybackslash}p{0.5\linewidth}>{\raggedright\arraybackslash}p{0.1\linewidth}>{\raggedright\arraybackslash}p{0.09\linewidth}@{}}
\hline
\iga name (layer) & statement & reading & checked \\
\hline
\textsf{gi\_\allowbreak goedel2\_\allowbreak internal} (ent) & $\ConI \vdash \neg\PrC(\bquo{\ConI})$ & in \iga & yes \\
\textsf{gi\_\allowbreak goedel2\_\allowbreak internal\_\allowbreak contra} (ent) & $\PrC(\bquo{\ConI}) \vdash \neg\ConI$ & in \iga & yes \\
\textsf{gi\_\allowbreak ent\_\allowbreak nec} (ent) & $A \vdash B$ gives $\PrC(\bquo{A}), \Gamma \vdash \PrC(\bquo{B})$ & classical & yes \\
\textsf{gi\_\allowbreak D3\_\allowbreak ent} (ent) & $\PrC(\bquo{\varphi}), \Gamma \vdash \PrC(\bquo{\PrC(\bquo{\varphi})})$ & in \iga & yes \\
\textsf{gi\_\allowbreak goedelG\_\allowbreak under\_\allowbreak ConIGA} (ent) & $\ConI \vdash G$ & in \iga & no \\
\textsf{gi\_\allowbreak goedelG\_\allowbreak refutes} (acc) & $\emptyset \vdash G$ gives $\emptyset \vdash \z = \tofn 1$ & classical & yes \\
\textsf{gi\_\allowbreak ConIGA\_\allowbreak of\_\allowbreak code\_\allowbreak unprovable} (exf) & for well-formed $\varphi$, $\emptyset \vdash \forall v.\,\neg\ChkI(v, \bquo{\varphi})$ gives $\emptyset \vdash \ConI$ & classical & no \\
\textsf{gi\_\allowbreak D2\_\allowbreak ent} (d2) & for well-formed $\varphi, \psi$: $\PrC(\bquo{\varphi}), \PrC(\bquo{\varphi \limp \psi}), \Gamma \vdash \PrC(\bquo{\psi})$ & in \iga & no \\
\textsf{gi\_\allowbreak formalized\_\allowbreak loeb\_\allowbreak ent} (d2) & for closed well-formed $\psi$: $\PrC(\bquo{\PrC(\bquo{\psi}) \limp \psi}) \vdash \PrC(\bquo{\psi})$ & in \iga & yes \\
\textsf{gi\_\allowbreak formalized\_\allowbreak loeb\_\allowbreak 0eq1\_\allowbreak ent} (d2) & $\PrC(\bquo{\neg\ConI \limp \z = \tofn 1}) \vdash \neg\ConI$ & in \iga & no \\
\textsf{gi\_\allowbreak formalized\_\allowbreak loeb\_\allowbreak sentence} (d2) & for closed well-formed $\psi$: $\emptyset \vdash \forall v.\,(\ChkI(v, \bquo{\PrC(\bquo{\psi}) \limp \psi}) \limp \PrC(\bquo{\psi}))$ & in \iga & no \\
\textsf{gi\_\allowbreak ent\_\allowbreak nec2} (nec2) & $A_1, A_2 \vdash B$ gives $\PrC(\bquo{A_2}), \PrC(\bquo{A_1}), \Gamma \vdash \PrC(\bquo{B})$ & classical & yes \\
\textsf{gi\_\allowbreak bs\_\allowbreak theorem1} (exf) & Beklemishev and Shamkanov's Theorem~1, for a matrix $D$ decided at a fresh variable and a sentence $\neg\forall v.\,\neg D$ that meets their two fixed-point entailments & classical & no \\
\textsf{gi\_\allowbreak bs\_\allowbreak goedel2} (exf) & three entailments between $\neg G$ and $\PrC(\bquo{\neg(\z = \z)})$ & in \iga & no \\
\textsf{gi\_\allowbreak bs\_\allowbreak structure\_\allowbreak conditions} (exf\_acc) & Beklemishev and Shamkanov's structure conditions & in \iga and classical & yes \\
\textsf{gi\_\allowbreak bs\_\allowbreak goedel2\_\allowbreak irrefutable} (exf\_con) & $\PrC(\bquo{\neg(\z = \z)}) \nvdash \neg(\z = \z)$ & classical & no \\
\textsf{gi\_\allowbreak goedel1} (con) & $\emptyset \nvdash G$ & classical & no \\
\textsf{gi\_\allowbreak goedel2} (con) & $\emptyset \nvdash \ConI$ & classical & no \\
\textsf{gi\_\allowbreak no\_\allowbreak loeb\_\allowbreak hypothesis\_\allowbreak 0eq1} (con) & $\emptyset \nvdash \forall v.\,(\ChkI(v, \bquo{\z = \tofn 1}) \limp \z = \tofn 1)$ & classical & no \\
\textsf{gi\_\allowbreak no\_\allowbreak code\_\allowbreak unprovable} (exf\_con) & for well-formed $\varphi$, $\emptyset \nvdash \forall v.\,\neg\ChkI(v, \bquo{\varphi})$ & classical & no \\
\hline
\end{tabular}
\caption{The theorems that the nine interpretations state for \iga,
  with the layer of each (\textsf{g2\_\allowbreak ent}, \textsf{g2\_\allowbreak acc}, \dots,
  abbreviated).  Theorem~\ref{thm:ig-g2}, \textsf{gi\_\allowbreak goedel2\_\allowbreak ConIGA},
  restates \textsf{gi\_\allowbreak goedel2}.}
\label{tab:ig-reexports}
\end{table}

\paragraph{The build checks.}
\label{app:iga-goedel-checks}
Each of the three theories of the session \textsf{IGA\_\allowbreak Goedel} contains
build-time checks on the dependency graphs of its theorems.

\emph{No oracle.}  No theorem among the checked heads depends on an
oracle: five heads in \textsf{GCalc\_\allowbreak IGA}, nineteen in
\textsf{Goedel\_\allowbreak Gen\_\allowbreak Transformers} and twenty in \textsf{Goedel\_\allowbreak IGA}.
The twenty heads of \textsf{Goedel\_\allowbreak IGA} do not include every
re-export: \textsf{gi\_\allowbreak D2\_\allowbreak ent}, \textsf{gi\_\allowbreak formalized\_\allowbreak loeb\_\allowbreak sentence},
\textsf{gi\_\allowbreak bs\_\allowbreak theorem1}, \textsf{gi\_\allowbreak bs\_\allowbreak goedel2},
\textsf{gi\_\allowbreak goedelG\_\allowbreak refutes}, \textsf{gi\_\allowbreak goedelG\_\allowbreak under\_\allowbreak ConIGA} and
\textsf{gi\_\allowbreak ConIGA\_\allowbreak of\_\allowbreak code\_\allowbreak unprovable} are not heads.

\emph{The calculus.}  \textsf{gig\_\allowbreak gcalc\_\allowbreak strip} reaches
\textsf{gig\_\allowbreak prf\_\allowbreak strip}, \textsf{gig\_\allowbreak gstep} and the ground-induction
lemma \textsf{ir\_\allowbreak gfix\_\allowbreak par}, and avoids \rga's \textsf{prf\_\allowbreak strip} and
\textsf{ir\_\allowbreak strip}.

\emph{The laws at \iga's table.}  Each \iga instance of a law reaches
its generic law and avoids \rga's own version of that law, by the
following avoid lists.
\begin{itemize}
\item NEC (\textsf{chkI\_\allowbreak cut1M\_\allowbreak chk}): \textsf{cut1M\_\allowbreak chk},
  \textsf{cut1Cert\_\allowbreak chk}, \textsf{iTreeValid\_\allowbreak appC},
  \textsf{iChk\_\allowbreak lstCons}, \textsf{closed\_\allowbreak proof\_\allowbreak head},
  \textsf{d2\_\allowbreak valid\_\allowbreak accepted}, \textsf{cmp\_\allowbreak num\_\allowbreak one\_\allowbreak iff},
  \textsf{prf\_\allowbreak owsjt} and \textsf{rga\_\allowbreak consistent}.  The instance
  reaches \textsf{G\_\allowbreak cut1M\_\allowbreak chk}, \textsf{G\_\allowbreak iTreeValid\_\allowbreak appC},
  \textsf{tblShift\_\allowbreak IGA}, \textsf{chkG\_\allowbreak validChkD\_\allowbreak correct},
  \textsf{cmp\_\allowbreak eq\_\allowbreak num}, \textsf{drvChkI\_\allowbreak fires\_\allowbreak iff},
  \textsf{adi\_\allowbreak iChkD\_\allowbreak lstCons} and \textsf{ndFire\_\allowbreak cut\_\allowbreak intro}.
\item EXPL (\textsf{chkI\_\allowbreak goedelE\_\allowbreak chk}): \textsf{goedelE\_\allowbreak chk},
  \textsf{iChk\_\allowbreak appC} and \textsf{iChk\_\allowbreak lstCons}.
\item Ex falso (\textsf{chkI\_\allowbreak exfIdx\_\allowbreak chk\_\allowbreak sf}): \textsf{exfIdx\_\allowbreak chk\_\allowbreak sf},
  \textsf{iChk\_\allowbreak appC} and \textsf{iChk\_\allowbreak lstCons}.
\item D2 (\textsf{chkI\_\allowbreak d2M\_\allowbreak chk}): \textsf{d2M\_\allowbreak chk},
  \textsf{d2Cert\_\allowbreak chk}, \textsf{d2\_\allowbreak cons\_\allowbreak valid} and
  \textsf{iChk\_\allowbreak lstCons}.
\item NEC with two premises (\textsf{chkI\_\allowbreak cut2M\_\allowbreak chk}):
  \textsf{cut2M\_\allowbreak chk}, \textsf{cut2Cert\_\allowbreak chk},
  \textsf{closed\_\allowbreak proof\_\allowbreak head} and \textsf{cmp\_\allowbreak num\_\allowbreak one\_\allowbreak iff}.
\item D3 (\textsf{chkI\_\allowbreak d3SigmaAt}): \rga's \textsf{iChk\_\allowbreak d3SigmaAt}
  and \textsf{iChk\_\allowbreak d3Cert}, and the four facts about \iga's table
  \textsf{tblShift\_\allowbreak IGA}, \textsf{tblShift\_\allowbreak chkGfix},
  \textsf{chkGfix\_\allowbreak correct} and \textsf{drvChkI\_\allowbreak fires\_\allowbreak iff}.  The
  instance reaches \textsf{chkG\_\allowbreak iChkD\_\allowbreak d3SigmaAt},
  \textsf{G\_\allowbreak iChk\_\allowbreak d3SigmaAt}, \textsf{G\_\allowbreak iChk\_\allowbreak d3Cert},
  \textsf{pr\_\allowbreak cert\_\allowbreak builder} and \textsf{chkG\_\allowbreak embed}.
\end{itemize}
Of these lists, only NEC's names soundness or consistency theorems.
NEC's list and the list for NEC with two premises name
\textsf{cmp\_\allowbreak num\_\allowbreak one\_\allowbreak iff}, the two-way closed-evaluation lemma whose
forward direction rests on \rga's soundness
(\S\ref{sec:ig-closed}).  The lists for EXPL, ex falso and D2 name
none of these, and D1 (\textsf{chkI\_\allowbreak proves\_\allowbreak sf}) gets only the
no-oracle check in this theory.  The soundness checks of
\textsf{Goedel\_\allowbreak IGA} cover those witnesses.

\emph{Soundness.}  The dependency cones of eight theorems avoid eleven
names: \textsf{gi\_\allowbreak consistent}, \textsf{gi\_\allowbreak step\_\allowbreak oomjt},
\textsf{omega\_\allowbreak sound}, \textsf{rga\_\allowbreak consistent},
\textsf{rga\_\allowbreak omega\_\allowbreak sound}, \textsf{prf\_\allowbreak owsjt}, \textsf{prf\_\allowbreak oomjt},
\textsf{iChkI\_\allowbreak genuine\_\allowbreak proves}, \textsf{adi\_\allowbreak omega\_\allowbreak sound},
\textsf{adi\_\allowbreak consistent} and \textsf{rga\_\allowbreak no\_\allowbreak proof\_\allowbreak 0eq1}.  The eight
theorems are \textsf{gi\_\allowbreak goedel2\_\allowbreak internal},
\textsf{gi\_\allowbreak goedel2\_\allowbreak internal\_\allowbreak contra}, \textsf{gi\_\allowbreak D3\_\allowbreak ent},
\textsf{gi\_\allowbreak ent\_\allowbreak nec}, \textsf{gi\_\allowbreak goedelG\_\allowbreak refutes},
\textsf{gi\_\allowbreak formalized\_\allowbreak loeb\_\allowbreak ent}, \textsf{gi\_\allowbreak ent\_\allowbreak nec2} and
\textsf{gi\_\allowbreak bs\_\allowbreak structure\_\allowbreak conditions}.  The check covers more than
these eight theorems' own proofs.  A theorem of a layer is exported
with the whole locale predicate of that layer as its premise, and the
interpretation discharges the premise with the layer's whole witness.
So an interpreted theorem's dependency cone holds every witness of its
layer, whether or not the theorem's proof uses that condition.  The
eight checks thus cover all ten witnesses other than consistency: the three structural ones,
NEC, D3, EXPL, D1, ex falso, D2 and NEC with two premises.  For
example, \textsf{gi\_\allowbreak bs\_\allowbreak structure\_\allowbreak conditions} belongs to the layer
\textsf{g2\_\allowbreak exf\_\allowbreak acc}, and its cone holds the witnesses for ex falso
and for D1.

Twelve of the twenty re-exports are not checked against the soundness
names.  Five of them are classical endpoints that use consistency by
design: \textsf{gi\_\allowbreak goedel1}, \textsf{gi\_\allowbreak goedel2},
\textsf{gi\_\allowbreak no\_\allowbreak loeb\_\allowbreak hypothesis\_\allowbreak 0eq1}, \textsf{gi\_\allowbreak no\_\allowbreak code\_\allowbreak unprovable}
and \textsf{gi\_\allowbreak bs\_\allowbreak goedel2\_\allowbreak irrefutable}.  The other seven are
derivations or rules: \textsf{gi\_\allowbreak goedelG\_\allowbreak under\_\allowbreak ConIGA},
\textsf{gi\_\allowbreak ConIGA\_\allowbreak of\_\allowbreak code\_\allowbreak unprovable}, \textsf{gi\_\allowbreak D2\_\allowbreak ent},
\textsf{gi\_\allowbreak formalized\_\allowbreak loeb\_\allowbreak 0eq1\_\allowbreak ent},
\textsf{gi\_\allowbreak formalized\_\allowbreak loeb\_\allowbreak sentence}, \textsf{gi\_\allowbreak bs\_\allowbreak theorem1} and
\textsf{gi\_\allowbreak bs\_\allowbreak goedel2}.  Their witnesses lie inside checked cones, so
for these seven only the proofs inside the locale are unchecked.

\emph{Reach.}  \textsf{gi\_\allowbreak goedel2\_\allowbreak ConIGA} reaches
\textsf{g2\_\allowbreak con\_\allowbreak gi}, \textsf{chkI\_\allowbreak cut1M\_\allowbreak chk},
\textsf{chkI\_\allowbreak d3SigmaAt}, \textsf{chkI\_\allowbreak goedelE\_\allowbreak chk},
\textsf{gig\_\allowbreak prf\_\allowbreak strip}, \textsf{g2i\_\allowbreak con0}, \textsf{gi\_\allowbreak consistent}
and \textsf{g2\_\allowbreak rw\_\allowbreak con\_\allowbreak gi}.  \textsf{gi\_\allowbreak separation\_\allowbreak IGA} reaches
\textsf{gi\_\allowbreak goedel2\_\allowbreak ConIGA}, \textsf{adi\_\allowbreak ConIGA} and
\textsf{adi\_\allowbreak consistent}.  These need checks also show that the
dependency computation used by the avoid checks is not blind: through
the interpretations, it reaches the witnesses.

\emph{Not checked.}  The build does not check the theorems of this
appendix against the decider.  The session of the locale has the
decider's session as its parent, for the index lemmas of the diagonal
function.

\emph{Size.}  The session \textsf{IGA\_\allowbreak Goedel} has three theories and
1,652 lines (\textsf{GCalc\_\allowbreak IGA} 329, \textsf{Goedel\_\allowbreak Gen\_\allowbreak Transformers}
996, \textsf{Goedel\_\allowbreak IGA} 327), counted by \texttt{wc -l}.

\section{The generic checker's and \iga's tables}
\label{app:checker}

This appendix gives the generic checker of Section~\ref{sec:checker}
and \iga's instance at the level of the Isabelle text, for a reader who
wants to re-mechanize them or to add a tag of their own.  Every
statement is in \texttt{context rga}, so $\vdash$ is \rga's derivability,
except \textsf{gi\_\allowbreak derives\_\allowbreak iffI}, which is in
\texttt{context gi}.  The Isabelle text writes $\vdash$, a numeral
$\tofn n$ and a code $\quo{x}$ as in the paper; \textsf{pr2fi} is the
HOL value of an index tree.

\subsection{Files and sizes}
\label{app:ck-files}

Table~\ref{tab:ck-files} lists the theory files.  The session
\texttt{RGA\_Checker\_Gen} has the parent \texttt{RGA\_Checker}, and
\texttt{IGA\_Checker} has the parent \texttt{RGA\_Checker\_Gen} and
imports the sessions \texttt{IGA} and \texttt{Ordinal}.  \texttt{RGA\_Checker\_Gen} has 76 definitions and lemmas and
\texttt{IGA\_Checker} has 111, counted as
the lines that begin with \texttt{lemma}, \texttt{theorem},
\texttt{corollary}, \texttt{definition}, \texttt{abbreviation},
\texttt{fun}, \texttt{primrec} or \texttt{lemmas}.  Both sessions only
add theories; they change no theory of \texttt{RGA\_Kit},
\texttt{RGA\_Checker}, \texttt{RGA\_Goedel} or \texttt{RGA}.

\begin{table}[ht]
\centering
\small
\begin{tabular}{@{}>{\raggedright\arraybackslash}p{0.33\linewidth}>{\raggedleft\arraybackslash}p{0.09\linewidth}>{\raggedright\arraybackslash}p{0.50\linewidth}@{}}
\hline
theory & lines & content \\
\hline
\texttt{RGA\_Checker\_Gen/\allowbreak Chk\_Gen\_Base} & 1,089 & the definitions at a table, the bridges, the unfolding and dispatch twins, the shift law and its closure, monotonicity and the embedding \\
\texttt{RGA\_Checker\_Gen/\allowbreak Chk\_Gen\_Graft} & 804 & append, graft, the certificate layer, the instances at \rga's table and at the test list, the build check \\
\texttt{IGA\_Checker/\allowbreak Chk\_IGA\_Pos} & 352 & the positivity decider and its explicit tree \\
\texttt{IGA\_Checker/\allowbreak Chk\_IGA\_Gfix} & 481 & the ground-induction tag, its payload and correctness, \iga's list and rule split, the tag's shift law \\
\texttt{IGA\_Checker/\allowbreak Chk\_IGA\_Table} & 378 & closed faithfulness at a table, $\Chk_{\mathcal D}$, \iga's checker, the instances at \iga's table, the build check \\
\hline
generic layer and \iga's instance & 3,104 & \\
\hline
\texttt{RGA\_Checker} (24 files) & 87,242 & \rga's checker and its laws \\
\texttt{RGA\_Decider} & 9,964 & the decider, and the correctness of \rga's rule test and of the citation formers \\
\texttt{RGA\_Kit} & 15,766 & the code-index library; \texttt{Code\_Index} (4,264 lines) holds the 41 detectors \\
\hline
\end{tabular}
\caption{The theory files of the generic checker and of \iga's instance,
  with \rga's checker sessions for comparison.  Lines are
  \texttt{wc -l} of the theory files.}
\label{tab:ck-files}
\end{table}

\subsection{Definitions at a table}
\label{app:ck-defs}

\rga's detector list \textsf{drvChkList} (\texttt{RGA\_Kit/\allowbreak
Code\_Index}) has 41 entries.  The table of a list is
$\textsf{drvChkIdxG}\ L \ldef \textsf{orNestIdx}\ L$, with
$\textsf{orNestIdx}\ [\,] = \textsf{zeroIdx}$ and $\textsf{orNestIdx}\ (E :: \mathit{Es}) =
\textsf{orIdx}\ E\ (\textsf{orNestIdx}\ \mathit{Es})$.  At a table
$\mathcal D$ (Isabelle \textsf{Dt}), \texttt{Chk\_Gen\_Base} defines
\textsf{validBodyD} (the node-code test conjoined with $\mathcal D$ at the
position's node), \textsf{validChkD} (the bounded universal of the body below
the list's length), \textsf{provesChkD} ($\pi_{\mathcal D}$: the sign of $p$,
the head code equal to the argument, and validity), \textsf{prfChkD}, the
abbreviation $\textsf{iTreeValidD}\ \mathcal D\ m \ldef
(\cmp(\textsf{validChkD}\ \mathcal D)\cdot m\cdot m = \tofn 1)$, and
\textsf{iHasCertD} ($\Cert_{\mathcal D}$, the negated universal
$\neg\forall p.\,\neg(\cmp(\pi_{\mathcal D})\cdot h\cdot p \ne \z)$).
\texttt{Chk\_IGA\_Table} defines
\[
  \textsf{iChkD}\ \mathcal D\ p\ c \ldef
  \bigl(\mathit{TreeValid}_{\mathcal D}(p) \wedge p \ne \z\bigr) \wedge
  \bigl(|\mathit{hyps}(p)| = \z \wedge \mathit{concl}(p) = c\bigr).
\]
Each definition is \rga's with $\mathcal D$ in place of \rga's rule
test \textsf{drvChkIdx}.  The bridges, all proved by unfolding, recover \rga's
machines (Table~\ref{tab:ck-bridges}).

\begin{table}[ht]
\centering
\small
\begin{tabular}{@{}>{\raggedright\arraybackslash}p{0.34\linewidth}>{\raggedright\arraybackslash}p{0.58\linewidth}@{}}
\hline
bridge & equation \\
\hline
\textsf{chkG\_\allowbreak bridge\_\allowbreak drvChkIdx} & $\textsf{drvChkIdxG}\ \textsf{drvChkList} = \textsf{drvChkIdx}$ \\
\textsf{chkG\_\allowbreak bridge\_\allowbreak validBody} & $\textsf{validBodyD}\ \textsf{drvChkIdx} = \textsf{validBody}$ \\
\textsf{chkG\_\allowbreak bridge\_\allowbreak validChk} & $\textsf{validChkD}\ \textsf{drvChkIdx} = \textsf{validChkIdx}$ \\
\textsf{chkG\_\allowbreak bridge\_\allowbreak provesChk}, \textsf{chkG\_\allowbreak bridge\_\allowbreak prfChk} & the same for $\pi$ and the packed proof test \\
\textsf{chkG\_\allowbreak bridge\_\allowbreak iTreeValid} & $\textsf{iTreeValidD}\ (\textsf{drvChkIdxG}\ \textsf{drvChkList})\ m = \textsf{iTreeValid}\ m$ \\
\textsf{chkG\_\allowbreak bridge\_\allowbreak iHasCert} & $\textsf{iHasCertD}\ \textsf{drvChkIdx}\ h = \textsf{iHasCert}\ h$ \\
\textsf{chkG\_\allowbreak bridge\_\allowbreak iChk} & $\textsf{iChkD}\ (\textsf{drvChkIdxG}\ \textsf{drvChkList})\ p\ c = \textsf{iChk}\ p\ c$ \\
\hline
\end{tabular}
\caption{The bridges from the generic checker at \rga's table to \rga's
  checker (Theorem~\ref{thm:ck-bridge}); the last is in
  \texttt{Chk\_IGA\_Table}, the others in \texttt{Chk\_Gen\_Base}.}
\label{tab:ck-bridges}
\end{table}

\subsection{The shift law}
\label{app:ck-shift}

Write a node as $k = \mathit{pair}(\mathit{head}(k), \mathit{tail}(k))$,
with $\mathit{head}(k)$ the judgment and $\mathit{tail}(k)$ the list
below it, and $\mathit{nd}(T, z)$ for the entry at position $z$ of the
list $T$.  The test $\mathit{fires}_E(k)$ (\textsf{ndFire},
\texttt{RGA\_Checker/\allowbreak Chk\_Formers}) is $\cmp(E)\cdot
\mathit{pair}(\mathit{tail}(k), \mathit{ctx}(k))\cdot \mathit{concl}(k)
\ne \z$.  Say that a context $\Delta$ \emph{extends} $\Gamma$ when $\Delta$
derives every term that $\Gamma$ derives.  The predicate \textsf{tblShift} of
a detector $E$ (\texttt{Chk\_Gen\_Base}) holds when, for every context
$\Gamma$, all terms $k$ and $k'$ and every map $\mathit{sh}$ from terms
to terms, the following premises give $\Gamma \vdash
\mathit{fires}_E(k')$:
\begin{enumerate}
\item $\Gamma \vdash k\jnat$ and $\Gamma \vdash k'\jnat$;
\item for every $\Delta$ that extends $\Gamma$ and every $z$: $\Delta
  \vdash z\jnat$ gives $\Delta \vdash \mathit{sh}(z)\jnat$;
\item $\Gamma \vdash \mathit{head}(k') = \mathit{head}(k)$;
\item for every $\Delta$ that extends $\Gamma$ and every $z$ with $\Delta
  \vdash z\jnat$: if $\Delta$ derives that $z$ is below
  $|\mathit{tail}(k)|$, then $\Delta$ derives that $\mathit{sh}(z)$ is
  below $|\mathit{tail}(k')|$;
\item for every such $\Delta$ and $z$: if $\Delta$ derives that $z$ is
  below $|\mathit{tail}(k)|$, then $\Delta \vdash
  \mathit{nd}(\mathit{tail}(k'), \mathit{sh}(z)) =
  \mathit{nd}(\mathit{tail}(k), z)$;
\item $\Gamma \vdash \mathit{fires}_E(k)$.
\end{enumerate}
``$z$ is below $\ell$'' is the equation $(z+1) \mathbin{\dot-} \ell = \z$
of truncated subtraction.  The predicate is the conclusion of \rga's
whole-table law \textsf{drvChk\_\allowbreak table\_\allowbreak shift} (\texttt{RGA\_Checker/\allowbreak
Chk\_Table\_Shift}), read as a property of one detector.
Table~\ref{tab:ck-shift} lists the facts about the predicate.

\begin{table}[ht]
\centering
\small
\begin{tabular}{@{}>{\raggedright\arraybackslash}p{0.30\linewidth}>{\raggedright\arraybackslash}p{0.62\linewidth}@{}}
\hline
fact & statement \\
\hline
\textsf{tblShiftD}, \textsf{tblShiftI} & elimination and introduction of the predicate \\
\textsf{G\_\allowbreak drvChk\_\allowbreak table\_\allowbreak shift} & the dispatch law over a list $L$: if every member of $L$ satisfies the law, then a fire of \textsf{orNestIdx} $L$ at $k$ gives one at $k'$, under the premises above \\
\textsf{tblShift\_\allowbreak orNest} & Theorem~\ref{thm:ck-shift}(1), from the previous fact \\
\textsf{tblShift\_\allowbreak drvChkList} & every member of \textsf{drvChkList} satisfies the law, by \rga's per-detector shift facts \textsf{ndFire\_\allowbreak tag\_\allowbreak table\_\allowbreak shifts}, \textsf{ndFire\_\allowbreak pr0\_\allowbreak table\_\allowbreak shift} and \textsf{ndFire\_\allowbreak hyp\_\allowbreak table\_\allowbreak shift} \\
\textsf{tblShift\_\allowbreak zeroIdx} & the detector that never fires satisfies the law \\
\textsf{tblShift\_\allowbreak RGA} & \rga's table satisfies the law \\
\textsf{tblShift\_\allowbreak drvChkIdx} & the same, from \rga's whole-table law, as a cross-check \\
\textsf{tblShift\_\allowbreak extList} & the table of the test list $\textsf{chkG\_\allowbreak extList} = \textsf{drvChkList} \mathbin{@} [\textsf{zeroIdx}]$ satisfies the law \\
\hline
\end{tabular}
\caption{The shift law's facts, in \texttt{Chk\_Gen\_Base}.  The test
  list adds a detector that never fires; the list is not \iga's.}
\label{tab:ck-shift}
\end{table}

\noindent
In \texttt{Chk\_Gen\_Graft} one proof,
\textsf{G\_\allowbreak iTreeValid\_\allowbreak appRun}, uses the predicate, with
$\mathit{sh}$ the identity; the append and
the graft pass the premise on.  At the test list and at \rga's table the
generic append and graft hold with no shift premise
(\textsf{chkG\_\allowbreak iTreeValid\_\allowbreak appC\_\allowbreak ext}, \textsf{chkG\_\allowbreak graft\_\allowbreak appC\_\allowbreak ext},
\textsf{chkG\_\allowbreak iTreeValid\_\allowbreak appC\_\allowbreak RGA}, \textsf{chkG\_\allowbreak graft\_\allowbreak appC\_\allowbreak RGA},
\textsf{chkG\_\allowbreak emitArm\_\allowbreak of\_\allowbreak tree\_\allowbreak RGA}).

\subsection{The twins}
\label{app:ck-twins}

Table~\ref{tab:ck-twins} lists the 38 lemmas of the generic layer that
are twins of \rga's checker laws.  Each twin's \rga original has the
same name without the prefix \textsf{G\_}, and the twin's statement is
the original's with a table $\mathcal D$ in place of \rga's rule test,
or, for the two dispatch laws and
\textsf{G\_\allowbreak drvChk\_\allowbreak table\_\allowbreak shift}, with a list
$L$ in place of \rga's detector list.
Five twins carry the shift law of $\mathcal D$ as a premise; the other
33 have no hypothesis on the table.

{\small
\setlength{\tabcolsep}{4pt}
\begin{longtable}{@{}>{\raggedright\arraybackslash}p{0.40\linewidth}>{\raggedright\arraybackslash}p{0.33\linewidth}>{\raggedright\arraybackslash}p{0.19\linewidth}@{}}
\caption{The 38 twins.  Base and Graft are \texttt{Chk\_Gen\_Base} and
  \texttt{Chk\_Gen\_Graft}; the original's theory is in
  \texttt{RGA\_Checker}.}
\label{tab:ck-twins}\\
\hline
twin (file) & original's theory & premise on $\mathcal D$ \\
\hline
\endfirsthead
\hline
twin (file) & original's theory & premise on $\mathcal D$ \\
\hline
\endhead
\hline
\endfoot
\textsf{G\_\allowbreak ir\_\allowbreak validChkIdx\_\allowbreak eq} (Base) & \texttt{Chk\_Formers} & none \\
\textsf{G\_\allowbreak ir\_\allowbreak validChk\_\allowbreak zero\_\allowbreak unfold} (Base) & \texttt{Chk\_Formers} & none \\
\textsf{G\_\allowbreak ir\_\allowbreak validChk\_\allowbreak suc\_\allowbreak unfold} (Base) & \texttt{Chk\_Formers} & none \\
\textsf{G\_\allowbreak ir\_\allowbreak validBody\_\allowbreak eq} (Base) & \texttt{Chk\_Formers} & none \\
\textsf{G\_\allowbreak ir\_\allowbreak validBody\_\allowbreak at\_\allowbreak node} (Base) & \texttt{Chk\_Formers} & none \\
\textsf{G\_\allowbreak ir\_\allowbreak validBody\_\allowbreak node\_\allowbreak D} (Base) & \texttt{Chk\_Formers} & none \\
\textsf{G\_\allowbreak ir\_\allowbreak validBody\_\allowbreak node\_\allowbreak I} (Base) & \texttt{Chk\_Formers} & none \\
\textsf{G\_\allowbreak ir\_\allowbreak validChk\_\allowbreak ne0\_\allowbreak one} (Base) & \texttt{Chk\_Formers} & none \\
\textsf{G\_\allowbreak validBody\_\allowbreak cons\_\allowbreak shift} (Base) & \texttt{Chk\_Formers} & none \\
\textsf{G\_\allowbreak iTreeValid\_\allowbreak cons\_\allowbreak split} (Base) & \texttt{Chk\_Formers} & none \\
\textsf{G\_\allowbreak iTreeValid\_\allowbreak suc\_\allowbreak split} (Base) & \texttt{Chk\_Formers} & none \\
\textsf{G\_\allowbreak iTreeValid\_\allowbreak cons\_\allowbreak intro} (Base) & \texttt{Chk\_Coder\_Bounds} & none \\
\textsf{G\_\allowbreak iTreeValid\_\allowbreak suc\_\allowbreak intro} (Base) & \texttt{Chk\_Coder\_Bounds} & none \\
\textsf{G\_\allowbreak ir\_\allowbreak drvChk\_\allowbreak fires\_\allowbreak elim} (Base) & \texttt{Chk\_Formers} & none \\
\textsf{G\_\allowbreak ir\_\allowbreak drvChk\_\allowbreak fires\_\allowbreak intro} (Base) & \texttt{Chk\_Coder\_Bounds} & none \\
\textsf{G\_\allowbreak drvChk\_\allowbreak table\_\allowbreak shift} (Base) & \texttt{Chk\_Table\_Shift} & shift law of each member \\
\textsf{G\_\allowbreak ir\_\allowbreak provesChkIdx\_\allowbreak eq} (Base) & \texttt{Chk\_Formers} & none \\
\textsf{G\_\allowbreak iTreeValid\_\allowbreak appRun} (Graft) & \texttt{Chk\_Graft} & shift law \\
\textsf{G\_\allowbreak iTreeValid\_\allowbreak appC} (Graft) & \texttt{Chk\_Graft} & shift law \\
\textsf{G\_\allowbreak graft\_\allowbreak appC} (Graft) & \texttt{Chk\_Graft} & shift law \\
\textsf{G\_\allowbreak wft\_\allowbreak iHasCert} (Graft) & \texttt{Chk\_HasCert} & none \\
\textsf{G\_\allowbreak iHasCertI} (Graft) & \texttt{Chk\_HasCert} & none \\
\textsf{G\_\allowbreak iHasCert\_\allowbreak of\_\allowbreak delivery} (Graft) & \texttt{Chk\_HasCert} & none \\
\textsf{G\_\allowbreak emitArm\_\allowbreak of\_\allowbreak tree} (Graft) & \texttt{Chk\_HasCert} & shift law \\
\textsf{G\_\allowbreak lift\_\allowbreak iHasCert\_\allowbreak f} (Graft) & \texttt{Chk\_HasCert} & none \\
\textsf{G\_\allowbreak lift\_\allowbreak iHasCert} (Graft) & \texttt{Chk\_HasCert} & none \\
\textsf{G\_\allowbreak iHasCertE\_\allowbreak sub} (Graft) & \texttt{Chk\_HasCert} & none \\
\textsf{G\_\allowbreak iHasCertE} (Graft) & \texttt{Chk\_HasCert} & none \\
\textsf{G\_\allowbreak iHasCert\_\allowbreak transport} (Graft) & \texttt{Chk\_HasCert} & none \\
\textsf{G\_\allowbreak subf\_\allowbreak iHasCert} (Graft) & \texttt{Chk\_HasCert} & none \\
\textsf{G\_\allowbreak bsub\_\allowbreak iHasCert} (Graft) & \texttt{Chk\_HasCert} & none \\
\textsf{G\_\allowbreak sub1\_\allowbreak iHasCert} (Graft) & \texttt{Chk\_HasCert} & none \\
\textsf{G\_\allowbreak ir\_\allowbreak lift\_\allowbreak iTreeValid} (Graft) & \texttt{Chk\_HasCert} & none \\
\textsf{G\_\allowbreak ir\_\allowbreak lift\_\allowbreak tree\_\allowbreak head} (Graft) & \texttt{Chk\_HasCert} & none \\
\textsf{G\_\allowbreak ir\_\allowbreak lift\_\allowbreak iHasCert} (Graft) & \texttt{Chk\_HasCert} & none \\
\textsf{G\_\allowbreak ir\_\allowbreak lift\_\allowbreak hasCert\_\allowbreak free2} (Graft) & \texttt{Chk\_HasCert} & none \\
\textsf{G\_\allowbreak iHasCertE2} (Graft) & \texttt{Chk\_HasCert} & none \\
\textsf{G\_\allowbreak iHasCertE3} (Graft) & \texttt{Chk\_HasCert} & none \\
\end{longtable}
}

\noindent
The append and graft twins also carry the shift premise, in the place
where \rga's graft uses \rga's whole-table law.
\textsf{G\_\allowbreak drvChk\_\allowbreak table\_\allowbreak shift} takes the
members' shift laws as a premise and uses one of them where \rga's
proof splits over the 41 detectors.
The dispatch laws \textsf{G\_\allowbreak ir\_\allowbreak drvChk\_\allowbreak fires\_\allowbreak elim} and \textsf{G\_\allowbreak ir\_\allowbreak drvChk\_\allowbreak fires\_\allowbreak intro}
hold over an arbitrary list.  The generic layer does not generalize
\rga's well-formedness invariant for checker trees, which reads \rga's
table.

\subsection{The ground-induction tag and the positivity decider}
\label{app:ck-tag}

\paragraph{The rule as containment.}
\textsf{gfixRuleI} $P\ \Gamma\ c'$ (\texttt{Chk\_IGA\_Gfix}) holds when there
are $\Phi$, $x$, $i$ and $c$ with $c' = c[x]$, $\Phi$ positive in its
recursion variable, the three judgments
\[
  \Gamma \vdash \mu\Phi\cdot x, \qquad
  i\jnat,\ \Phi[\lambda c/a][i/x],\ \Gamma \vdash c[i], \qquad
  \Gamma \vdash x\jnat
\]
all entries of the list $P$, $\Phi$ and $c$ well-formed, and $i$ free
in none of $\Gamma$, $\lambda c$ and $\lambda a.\lambda x.\Phi$.  This
definition is the clause of ground induction in \iga's rules
(\textsf{pr\_\allowbreak gfix}, \texttt{IGA/\allowbreak IGA\_Proof}), read as containment.
\textsf{gi\_\allowbreak derives\_\allowbreak iffI} splits \iga's step relation: \iga relates $P$,
$\Gamma$ and $c$ exactly when \rga's step relation or \textsf{gfixRuleI} does.

\paragraph{The payload.}
The tag is $\textsf{chkGfixIdx} \ldef \textsf{cite3cIdx}\ \textsf{payGfix}$, \rga's
three-premise citation former applied to the payload \textsf{payGfix}.  The
payload passes to \textsf{payGfix0} only the node's context $\Gamma$ and
conclusion $c'$ and the three cited judgments $j_1, j_2, j_3$; the
payload drops the tail component of the node's argument.  \textsf{payGfix0} is a bounded
search, \textsf{bexLeIdx}, for the eigenvariable $i$ up to the code of the
context of $j_2$, over the conjunction \textsf{gfixWitIdx}.  At the argument
$\langle\langle\langle\Gamma, c'\rangle, \langle j_1, \langle j_2,
j_3\rangle\rangle\rangle, i\rangle$, \textsf{gfixWitIdx} reads three codes off
the cited judgments: $\Phi$ from the fixed-point term in the conclusion
of $j_1$, $x$ as the argument of that conclusion, and the motive $c$ as
the conclusion of $j_2$ with the variable $i$ abstracted.
\textsf{gfixWitIdx} then conjoins fifteen checks:
\begin{enumerate}
\item the code of $\Phi$ is a term code;
\item the code of $x$ is a term code;
\item the context of $j_1$ is $\Gamma$;
\item the context of $j_3$ is $\Gamma$;
\item the conclusion of $j_1$ is the code of $\mu\Phi\cdot x$;
\item the conclusion of $j_3$ is the code of $x\jnat$;
\item the context of $j_2$ is the code of $i\jnat,\ \Phi[\lambda
  c/a][i/x],\ \Gamma$;
\item the conclusion of $j_2$ is the code of $c[i]$;
\item $\Phi$ is positive in its recursion variable (\textsf{pdecTIdx} at $1$);
\item $\Phi$ is well-formed;
\item $c$ is well-formed;
\item $i$ is not free in $\Gamma$;
\item $i+1$ is not free in $c$, that is, $i$ is not free in $\lambda c$;
\item $i+2$ is not free in $\Phi$, that is, $i$ is not free in
  $\lambda a.\lambda x.\Phi$;
\item $c' = c[x]$.
\end{enumerate}
The search ranges over $i$ only, as in \rga's detector for number
induction.  Two small trees build the codes of $\lambda$ and of the
fixed-point former (\textsf{gfxLamIdx}, \textsf{gfxGfixIdx}).
\textsf{chkGfix\_\allowbreak correct} states $\textsf{pr2fi}\ \textsf{chkGfixIdx}\ \langle\quo{P},
\quo{\Gamma}\rangle\ \quo{c'} = \quo{\textsf{try\_\allowbreak rule}\ \textsf{gfixRuleI}\
\langle\quo{P}, \quo{\Gamma}\rangle\ \quo{c'}}$, by \rga's
\textsf{cite3c\_\allowbreak rule\_\allowbreak correct} and the bounded form of the rule
(\textsf{gfixRuleI\_\allowbreak bounded}).  \textsf{ir\_\allowbreak payGfix\_\allowbreak read} is the payload's internal
reading equation, and \textsf{tblShift\_\allowbreak chkGfix} applies \rga's
\textsf{ndFire\_\allowbreak cite3c\_\allowbreak table\_\allowbreak shift\_\allowbreak read} (\texttt{RGA\_Checker/\allowbreak
Chk\_Emission}) with that equation and one congruence step
(\textsf{cmp\_\allowbreak argI} with \textsf{cpPair\_\allowbreak argI}).

\paragraph{The positivity decider.}
The recursion \textsf{pdecF} (\texttt{Chk\_IGA\_Pos}) is
\begin{align*}
  \textsf{pdecF}\ k\ (\textsf{trmVar}\ j) &= [\,j \ne k,\ j \ne k,\ j < k,\ \mathit{False},\ \mathit{False}\,],\\
  \textsf{pdecF}\ k\ (\textsf{trmOpr}\ \mathit{op}\ n\ T) &=
  \textsf{pdecStep}\ k\ \mathit{op}\ n\ T\ (\mathit{map}\ (\textsf{pdecF}\ (k+n))\ T),
\end{align*}
and \textsf{pdecF\_\allowbreak eq} is Theorem~\ref{thm:ck-pos}(1), proved by induction on
the term with the inversion lemmas for $\mathit{posP}$, $\mathit{negP}$
and $\mathit{pnexp}$.  The node step is the tree \textsf{pdecHIdx}, built from
\textsf{pdecLst}, \textsf{pdecIsOIdx} and \textsf{pdecGIdx}, and the decider's tree is
$\textsf{pdecPIdx} \ldef \textsf{covIdx}\ (\textsf{numIdx}\ \quo{[\mathit{False}, \ldots,
\mathit{False}]})\ \textsf{pdecHIdx}\ \textsf{covStepIdx}$, with five entries
$\mathit{False}$.  \textsf{pr2fi\_\allowbreak pdecPIdx} identifies its value with an
instance of the shared term recursion \textsf{primrec\_\allowbreak trmf}, whose variable
slot carries the binder depth, and \textsf{pdecF\_\allowbreak code} reads the code of
\textsf{pdecF}'s result off that value.  The tree \textsf{pdecTIdx} $K\ T$ selects
component~0, and \textsf{pr2fi\_\allowbreak pdecTIdx} gives the code of $\mathit{posP}\ k\
t$ when $K$ and $T$ evaluate to $k$ and $\quo{t}$.  Check~9 of the
payload is this tree at $1$, and \textsf{gfxPosC\_\allowbreak code} reads its value.  All of these facts are HOL equations
about codes.

\subsection{Closed faithfulness and the facts at \iga's table}
\label{app:ck-inst}

The closed-faithfulness facts of \texttt{Chk\_IGA\_Table} sit in two
unnamed contexts with the hypothesis
\[
  \mathrm{DF}:\quad \textsf{pr2fi}\ \mathcal D\ \langle\quo{P}, \quo{\Gamma}\rangle\ \quo{c} \ne 0
  \iff D\ P\ \Gamma\ c
  \qquad\text{for every proof } P \text{, context } \Gamma \text{ and term } c,
\]
which every exported fact carries as a premise; the step relation
$D$ of the display is the relation $R$ of \S\ref{sec:ck-faith}.
Table~\ref{tab:ck-inst} lists these facts and the instances at \iga's
table.

{\small
\setlength{\tabcolsep}{4pt}
\begin{longtable}{@{}>{\raggedright\arraybackslash}p{0.32\linewidth}>{\raggedright\arraybackslash}p{0.60\linewidth}@{}}
\caption{Closed faithfulness at a table and the facts at \iga's table
  (\texttt{Chk\_IGA\_Table}; \textsf{chkI\_proves\_sf} is in
  \texttt{IGA\_Goedel/\allowbreak Goedel\_Gen\_Transformers}).}
\label{tab:ck-inst}\\
\hline
fact & statement \\
\hline
\endfirsthead
\hline
fact & statement \\
\hline
\endhead
\hline
\endfoot
\textsf{chkG\_\allowbreak validChkD\_\allowbreak correct} & under DF, the HOL value of the validity fold at $\quo{P}$ is the code of ``$P$ is $R$-valid'' \\
\textsf{chkG\_\allowbreak validChkD\_\allowbreak genuine} & no DF: a number that the validity fold accepts is the code of a proof list \\
\textsf{chkG\_\allowbreak iTreeValidD\_\allowbreak faithful}, \textsf{chkG\_\allowbreak iTreeValidD\_\allowbreak genuine} & under DF, $\emptyset \vdash \mathit{TreeValid}_{\mathcal D}(\bquo{P})$ iff $P$ is $R$-valid; an accepted numeral is the code of an $R$-valid list \\
\textsf{chkG\_\allowbreak iChkD\_\allowbreak accept\_\allowbreak of\_\allowbreak parts} & under DF, an $R$-valid list $(\emptyset, c) :: Q$ is accepted for $c$ \\
\textsf{chkG\_\allowbreak iChkD\_\allowbreak faithful}, \textsf{chkG\_\allowbreak iChkD\_\allowbreak genuine\_\allowbreak proves} & Theorem~\ref{thm:ck-faith} \\
\textsf{iChkI} & the abbreviation $\ChkI = \textsf{iChkD}\ (\textsf{drvChkIdxG}\ \textsf{drvChkListI})$ \\
\textsf{iChkI\_\allowbreak validChk\_\allowbreak correct} & the validity fold at \iga's table decides \iga-validity \\
\textsf{iChkI\_\allowbreak faithful}, \textsf{iChkI\_\allowbreak genuine\_\allowbreak proves} & Theorem~\ref{thm:ck-iga-faith} \\
\textsf{chkI\_\allowbreak proves\_\allowbreak sf} & an \iga proof $P$ of $\emptyset \vdashI c$ gives $\Gamma \vdash \ChkI(\bquo{P}, \bquo{c})$ in every $\Gamma$, by one-directional closed evaluation \\
\textsf{chkI\_\allowbreak iTreeValid\_\allowbreak appC}, \textsf{chkI\_\allowbreak graft\_\allowbreak appC} & Theorem~\ref{thm:ck-graft} and the graft at \iga's table, with no shift premise \\
\textsf{chkI\_\allowbreak embed}, \textsf{chkI\_\allowbreak embed\_\allowbreak cert} & Theorems~\ref{thm:ck-embed}(2) and~\ref{thm:ck-embed-cert}(1) at \iga's list \\
\textsf{chkI\_\allowbreak iChk\_\allowbreak embed} & Theorem~\ref{thm:ck-embed-cert}(2) at \iga's list \\
\textsf{chkI\_\allowbreak emitArm\_\allowbreak of\_\allowbreak tree} & the general emission law through the graft, at \iga's table \\
\end{longtable}
}

\noindent
The bridge from derivations to values, \textsf{chkG\_\allowbreak iTreeValidD\_\allowbreak iff\_\allowbreak value}, is
\rga's \textsf{cmp\_\allowbreak num\_\allowbreak one\_\allowbreak iff}, which cites \textsf{rga\_\allowbreak val\_\allowbreak iff}; that equivalence
rests on \rga's soundness and completeness at closed codes.  As proved,
both directions of \textsf{chkG\_\allowbreak iTreeValidD\_\allowbreak faithful} therefore depend on \rga's
soundness.  The soundness-free acceptance \textsf{chkI\_\allowbreak proves\_\allowbreak sf} instead uses
\textsf{chkG\_\allowbreak validChkD\_\allowbreak correct}, the one-directional closed-numeral law
\textsf{cmp\_\allowbreak eq\_\allowbreak num} and weakening.  That property
is not checked in \texttt{Goedel\_Gen\_Transformers};
\texttt{Goedel\_IGA}'s soundness checks cover it
(Appendix~\ref{app:iga-goedel-checks}).

\subsection{Build checks}
\label{app:ck-checks}

The ML blocks below check dependencies when the sessions build.  Each block
fails the build if a listed theorem depends on an oracle, if a theorem
does not reach a fact on its need list, or if a theorem reaches a fact
on its avoid list.  A listed name matches a fact whose full name is the
name or ends in a dot followed by the name, so the avoid entry
\textsf{graft\_\allowbreak appC} matches \rga's graft and not \textsf{G\_\allowbreak graft\_\allowbreak appC}.

\begin{enumerate}
\item \texttt{Chk\_Gen\_Graft} (the block at the end of the file): no
  oracle under 32 headline facts, among them the bridges, the shift-law
  facts, the generic append and graft, their instances, monotonicity,
  the embeddings of trees and certificates, the general emission law
  and the certificate introduction, elimination and transport laws.
  \textsf{chkG\_\allowbreak graft\_\allowbreak appC\_\allowbreak ext} needs \textsf{G\_\allowbreak graft\_\allowbreak appC}, \textsf{tblShift\_\allowbreak zeroIdx},
  \textsf{tblShift\_\allowbreak drvChkList} and \textsf{G\_\allowbreak drvChk\_\allowbreak table\_\allowbreak shift}, and avoids
  \textsf{drvChk\_\allowbreak table\_\allowbreak shift} and \textsf{graft\_\allowbreak appC}.  \textsf{chkG\_\allowbreak embed\_\allowbreak ext} reaches no
  fact whose base name begins with \textsf{tblShift} or contains
  \textsf{table\_shift}.
\item \texttt{Chk\_IGA\_Table} (the block at the end of the file): no
  oracle under 20 headline facts: \textsf{pr2fi\_\allowbreak pdecTIdx}, the tag's
  correctness, the rule split, the shift laws of the tag and of \iga's
  table, closed faithfulness at a table and at \iga's table, and the
  instances at \iga's table.  \textsf{chkI\_\allowbreak graft\_\allowbreak appC} needs
  \textsf{G\_\allowbreak graft\_\allowbreak appC}, \textsf{tblShift\_\allowbreak chkGfix}, \textsf{tblShift\_\allowbreak drvChkList} and
  \textsf{G\_\allowbreak drvChk\_\allowbreak table\_\allowbreak shift}, and avoids \textsf{drvChk\_\allowbreak table\_\allowbreak shift} and
  \textsf{graft\_\allowbreak appC}.  So \iga's graft comes from the generic graft and the
  members' shift laws, not from \rga's whole-table graft.
  \textsf{iChkI\_\allowbreak faithful} needs \textsf{chkGfix\_\allowbreak correct}, \textsf{gfxPosC\_\allowbreak code}, \textsf{pdecF\_\allowbreak eq}
  and \textsf{drvChk\_\allowbreak fires\_\allowbreak iff}, so \iga's faithfulness goes through the
  explicit tag and the positivity decider's HOL semantics.
\item \texttt{IGA\_Goedel}: the checks of the three theories of
  Section~\ref{sec:iga-goedel} are described in
  Appendix~\ref{app:iga-goedel-checks}.
\end{enumerate}

\subsection{Definitions behind \iga's checker}
\label{app:ck-tcb}

$\ChkI$ unfolds to \rga's existing definitions and to new definitions
of three kinds.  The generic checker
adds \textsf{drvChkIdxG}, \textsf{validBodyD}, \textsf{validChkD} and \textsf{iChkD}
(\textsf{iTreeValidD} is an abbreviation).  The tag adds \textsf{chkGfixIdx},
\textsf{payGfix}, \textsf{payGfix0}, \textsf{gfixWitIdx}, \textsf{gfxLamIdx} and \textsf{gfxGfixIdx}, and the
positivity decider adds \textsf{pdecTIdx}, \textsf{pdecPIdx}, \textsf{pdecHIdx}, \textsf{pdecLst},
\textsf{pdecIsOIdx} and \textsf{pdecGIdx}.  These rest on definitions that \rga's
checker already uses: the citation former \textsf{cite3cIdx}, the code
destructors such as \textsf{kid1Idx} and \textsf{kid2Idx}, the shared child lookup
\textsf{liftPNFIdx}, and \rga's detectors.  The meaning of $\ChkI$ adds \iga's
rules, with ground induction \textsf{pr\_\allowbreak gfix}, through
Theorem~\ref{thm:ck-iga-faith}.

\section{Adoption in detail}
\label{app:adopt}

This appendix gives the constructions behind Section~\ref{sec:adopt}:
the witness tree of the warrant branch by branch, the ex-falso tree at
the level of codes, the soundness induction, the version at \iga, the
no-go locale at \igaP and at \oiga, and the build-time checks.  Every derivation named here is an \rga
derivation unless the text says otherwise.  Table~\ref{tab:fa-files}
lists the theories.

\begin{table}[tbp]
\centering
\small
\begin{tabular}{@{}lrl@{}}
\hline
theory & lines & content \\
\hline
\texttt{RGA\_Adopt/Adopt\_RGA} & 796 & \rgaP, the warrant at \rga, $\ConR$, the separation \\
\texttt{IGA\_Adopt/Adopt\_IGA} & 972 & \igaP, the warrant at \iga, $\ConI$, $\ConR$, the separations \\
\texttt{IGA\_Gentzen/Gentzen\_IGA} & 578 & the Gentzen conditional at \iga (\S\ref{sec:open-gentzen}) \\
\texttt{IGA\_NoGo/IGA\_NoGo} & 49 & the no-go locale at \iga (\S\ref{sec:open-nogo}) \\
\texttt{IGA\_NoGo/NoGo\_Ext} & 76 & the no-go locale at \igaP and \oiga (\S\ref{app:adopt-nogo}) \\
\hline
\end{tabular}
\caption{The theories of Sections~\ref{sec:adopt} and~\ref{sec:open},
  with their sizes in lines as counted by \texttt{wc -l}; the fold sessions of \S\ref{sec:open-stage} are
  listed in \S\ref{sec:mech}.}
\label{tab:fa-files}
\end{table}

\subsection{The witness tree}
\label{app:adopt-witness}

\begin{defi}[\textsf{adrF}, \textsf{adrGood}]
\label{def:fa-witness}
Let $k = \cpair{\quo{\emptyset}}{\quo{\z = \tofn 1}}$ be the code of
the judgment $\emptyset \vdash \z = \tofn 1$ (\textsf{adrK}), and let
$\mathit{test}(v) = \cmp(\mathit{provesChkIdx})\cdot\tofn k\cdot v$
(\textsf{adrTest}).  The witness is one index tree,
\[
  F_U \;=\; \mathit{ifIdx}(C, E, G) ,
\]
the conditional (\textsf{ifIdx}) on the test tree $C$
(\textsf{adrCIdx}) over two certificate builders: $E$
(\textsf{adrExfIdx}) at a nonzero test, $G$ (\textsf{adrGoodIdx}) at a
zero test.  At a natural $v$, \rga proves that the three trees compute
(\textsf{adr\_F\_reads}):
\begin{itemize}
\item $C$ computes $\mathit{test}(v)$;
\item $E$ computes $\mathit{Exf}(v) = \mathit{cons}\bigl(\mathit{jdg}(\Inst_U(v)),\
  \mathit{app}(v, \bquo{\mathit{neq01Prf}})\bigr)$ (\textsf{adrExf});
\item $G$ computes $\mathit{Good}(v)$ (\textsf{adrGood}): a list whose
  root is the judgment $\Inst_U(v)$, above a step-back
  (\textsf{certStepBack}) from the head reduct of $\Inst_U(v)$, above a
  judgment that the body's image at $v$ equals $1$, above the value
  certificate $\mathit{Val}(v)$ (\textsf{adrValCert}) that the builder
  of Theorem~\ref{thm:d-builder} computes for the test at $v$.
\end{itemize}
\end{defi}

\noindent
The conditional has only its HOL reading in the code-index library
(\texttt{RGA\_Kit}).  Its internal laws at a nonzero and at a zero
guard (\textsf{adr\_ifIdx\_ne0}, \textsf{adr\_ifIdx\_eq0}) follow
from the mask and addition laws of compiled machines.
$\cmp(\mathit{ifIdx}(C,E,G))\cdot a\cdot b$ is the sum of $E$'s value masked by $C$'s value and $G$'s value masked by
the zero test of $C$'s value.

\paragraph{The proof of Theorem~\ref{thm:ad-warrant}(1)
(\textsf{adr\_cert}).}  Fix $v$ with $\Gamma \vdash v\jnat$.
$\mathit{test}(v)$ is a natural, so $\mathit{test}(v) \neq \z$ is
decided, and \rga splits on it (\textsf{ir\_boolE}, \textsf{ne0\_B}).
\begin{itemize}
\item \emph{Nonzero.}  Under $\mathit{test}(v) \neq \z$,
  \textsf{adr\_ifIdx\_ne0} gives $\cmp(F_U)\cdot v\cdot v =
  \mathit{Exf}(v)$.  The nonzero test gives $\Chk(v, \bquo{\z = \tofn
  1})$ (\textsf{adr\_iChk\_of\_test}), and the instance code is a
  well-formed genuine code (\textsf{adr\_inst\_code}), so
  Theorem~\ref{thm:fa-exf} accepts $\mathit{Exf}(v)$ for $\Inst_U(v)$
  (\textsf{adr\_exf}).
\item \emph{Zero.}  Under $\neg(\mathit{test}(v) \neq \z)$, double
  negation gives $\mathit{test}(v) = \z$, and \textsf{adr\_ifIdx\_eq0}
  gives $\cmp(F_U)\cdot v\cdot v = \mathit{Good}(v)$.  The checker
  accepts $\mathit{Good}(v)$ for $\Inst_U(v)$ (\textsf{adr\_good}),
  in four links: the value certificate (\textsf{pr\_cert\_builder\_D},
  the form of Theorem~\ref{thm:d-builder} used here, which holds for
  every index tree); a trueEqI node (\textsf{adr\_link\_trueEqI}); two
  backward head steps to the instance code
  (\textsf{iChk\_certStepBack\_at}); and a trueEqE node
  (\textsf{adr\_link\_trueEqE}).
\end{itemize}
Transport along the computed equation (\textsf{iChk\_proof\_transport})
and disjunction elimination close the two cases.  The warrant,
Theorem~\ref{thm:ad-warrant}(2) (\textsf{adr\_warrant\_uC}), follows by
$\forall$-introduction at the free variable (\textsf{ir\_allI\_open}).

\subsection{Ex falso at the level of codes}
\label{app:adopt-exf}

Theorem~\ref{thm:d-exf} turns an accepted proof of $\z = \tofn 1$ into
an accepted proof of a fixed formula $\varphi$.  The warrant needs the
same construction at a target that is a code term, the instance code
$\Inst_U(v)$ at a free $v$, rather than the code numeral of a fixed
formula.

\begin{thm}[\textsf{adr\_\allowbreak exf\_\allowbreak code}]
\label{thm:fa-exf}
Suppose $\Gamma \vdash a\jnat$, $\Gamma \vdash \Chk(a, \bquo{\z =
\tofn 1})$, $\Gamma \vdash y\jnat$, and \rga proves in $\Gamma$ that
$y$ is a well-formed genuine code ($\Wfc(y) \neq \z$ and $\IsCode(y)
\neq \z$).  Then
\[
  \Gamma \vdash \Chk\bigl(\mathit{cons}(\mathit{jdg}(y),\
  \mathit{app}(a, \bquo{\mathit{neq01Prf}})),\ y\bigr) .
\]
\end{thm}

\noindent
\emph{Proof.}  As for Theorem~\ref{thm:d-exf}.  The append law
(Theorem~\ref{thm:d-appC}), with the accepted fixed proof of $\neg(\z
= \tofn 1)$ (\textsf{iChk\_proves\_sf} at \textsf{neq01Prf\_proves}),
gives the accepted list $\mathit{app}(a, \bquo{\mathit{neq01Prf}})$,
with $\z = \tofn 1$ and $\neg(\z = \tofn 1)$ at known positions.  The
new root's detector, negE, fires citing those two positions
(\textsf{ndFire\_negE\_intro}).  The format equation reads the two
conclusions off the list, and the root's well-formedness comes from the
hypothesis on $y$.  The cons law (Theorem~\ref{thm:d-cons}) accepts the
extended list.  \qed

\subsection{The soundness induction}
\label{app:adopt-sound}

The proof of the adoption step of Theorem~\ref{thm:ad-sound}
(\textsf{adr\_adopt\_omv}) runs as follows.  The warrant $\WarR(F, u)
= \cA\cdot w$, with $w$ its matrix as a $\lambda$-term, is true in
$\Mw$, so each instance $w\cdot\tofn n$ is true in $\Mw$
(\textsf{omv\_Aapp\_inv}), and so is its head reduct
(\textsf{om\_hstep\_rev\_omv}), the closed checker verdict
$\Chk(\cmp(F)\cdot\tofn n\cdot\tofn n, \Inst_u(\tofn n))$.  The verdict
has no quantifier, so \rga derives it (\textsf{om\_osr\_reflect},
\textsf{rga\_complete\_open}).  Closed evaluation of the compiled
application and of the instance code (\textsf{ir\_cmp\_num2},
\textsf{ir\_mkAppS\_code}) turns the derived verdict into $\emptyset
\vdash \Chk(\tofn p, \bquo{u\cdot\tofn n})$ with $p$ the value of $F$
at $(n, n)$.  Genuineness (\textsf{iChk\_genuine\_proves}) makes $p$
the code of an \rga proof of $\emptyset \vdash u\cdot\tofn n$, so \rga
derives $u\cdot\tofn n$ (\textsf{entails\_provesI}), which is then true
in $\Mw$ (\textsf{rga\_omega\_sound}).

For the system, each \rgaP step is $\Mw$-sound
(\textsf{adr\_step\_oomjt}): an \rga step by the per-step soundness of
\rga's rules in $\Mw$ (\textsf{rga\_step\_mjt}), the
adoption step by \textsf{adr\_adopt\_omv}.  An induction over \rgaP
proof lists (\textsf{adr\_prf\_oomjt}) gives
Theorem~\ref{thm:ad-sound}(1) (\textsf{adr\_omega\_sound}).

\subsection{The interderivability}
\label{app:adopt-interderivable}

\emph{$\cA\cdot U \vdash \ConR$} (\textsf{adr\_U\_ConRGA}).  At a free
$v$, the instance $\mathit{test}(v) = \z$ holds.  Decide $\Chk(v,
\bquo{\z = \tofn 1})$ (\textsf{hq\_iChk}).  In the accepting case the
test is nonzero (\textsf{adr\_test\_of\_iChk}), which contradicts the
instance, and negE gives $\neg\Chk(v, \bquo{\z = \tofn 1})$.  The other
case is the hypothesis $\neg\Chk(v, \bquo{\z = \tofn 1})$ itself.
$\forall$-introduction gives $\ConR$.

\emph{$\ConR \vdash \cA\cdot U$} (\textsf{adr\_ConRGA\_U}).  At a free
$v$, decide $\mathit{test}(v) \neq \z$.  In that case $\Chk(v,
\bquo{\z = \tofn 1})$ (\textsf{adr\_iChk\_of\_test}) contradicts
$\ConR$'s instance, and negE closes the case.  In the other case,
double negation (\textsf{ir\_dnegE}) gives $\mathit{test}(v) = \z$.
$\forall$-introduction gives $\cA\cdot U$.

\subsection{The version at \iga}
\label{app:adopt-iga}

At \iga the compiled proves-test is \iga's
(\textsf{adiPC}, the proves-test of the table of $\ChkI$); the body
$U_{\mathrm I}$ (\textsf{adiU}) and the witness $F_{U_{\mathrm I}} =
\mathit{ifIdx}(C_{\mathrm I}, E_{\mathrm I}, G_{\mathrm I})$
(\textsf{adiF}) have the shapes of Definition~\ref{def:fa-witness}
with that test in place of \rga's.

\paragraph{The zero branch.}  The value certificate, the trueEqI and
trueEqE links and the backward head steps are \rga's, and \rga's
checker accepts them.  The result is embedded into \iga's checker by
one law (\textsf{chkI\_iChk\_embed}): if $\Gamma \vdash p\jnat$ and
$\Gamma \vdash \Chk(p, c)$, then $\Gamma \vdash \ChkI(p, c)$.

\paragraph{The nonzero branch.}  A nonzero test gives $\ChkI(v,
\bquo{\z = \tofn 1})$ (\textsf{adi\_iChk\_of\_test}), and the code-level
ex falso at \iga applies.

\begin{thm}[\textsf{adi\_\allowbreak exf\_\allowbreak code}]
\label{thm:fa-exf-iga}
Theorem~\ref{thm:fa-exf} holds with $\ChkI$ in place of $\Chk$, both in
the hypothesis on $a$ and in the conclusion.
\end{thm}

\noindent
The construction is that of Theorem~\ref{thm:fa-exf}, on the generic
layer of Section~\ref{sec:checker}: the generic append law at \iga's table
(\textsf{G\_iTreeValid\_appC}, at the shift law \textsf{tblShift\_IGA}),
the generic cons law (\textsf{G\_iTreeValid\_suc\_intro}), the generic
firing law of a detector (\textsf{G\_ir\_drvChk\_fires\_intro}), and the
embedding \textsf{chkI\_iChk\_embed} for the fixed proof of $\neg(\z =
\tofn 1)$.  The emitted negE node uses an \rga tag, which belongs to
\iga's table.

\paragraph{Soundness.}  The adoption step of \igaP
(\textsf{adi\_adopt\_omv}) follows the proof of
Appendix~\ref{app:adopt-sound}: \rga derives the closed verdict
(\textsf{rga\_complete\_open}), the genuineness of $\ChkI$
(\textsf{iChkI\_genuine\_proves}) gives an \iga proof of the instance,
and \iga's $\Mw$-soundness (\textsf{omega\_sound}) makes the instance
true in $\Mw$.  The proof reuses \iga's per-step soundness
(\textsf{gi\_step\_oomjt}) for \iga's 43 rules.

\paragraph{$\ConI \vdash \ConR$} (\textsf{adi\_ConRGA\_of\_ConIGA}).
At a free $v$, decide $\Chk(v, \bquo{\z = \tofn 1})$.  In the accepting
case, the embedding law at the free $v$ gives $\ChkI(v, \bquo{\z =
\tofn 1})$, which contradicts $\ConI$'s instance; negE closes the case.
$\forall$-introduction gives $\ConR$.

\subsection{The no-go locale at \texorpdfstring{\igaP}{IGA+} and at \oiga}
\label{app:adopt-nogo}

The no-go locale \textsf{omega\_\allowbreak sound\_\allowbreak ext} of
the internal-reasoning paper ([\citen{conrga-paper}, theorem~\textsf{xder\_not\_fold\_residual}])
assumes one thing: an entailment over \rga's terms that is sound for
$\Mw$.  The theory \texttt{NoGo\_Ext} of the session
\texttt{IGA\_NoGo} interprets the locale twice, inside \rga's
context, as the theory \texttt{IGA\_NoGo} does at \iga:
\begin{itemize}
\item \textsf{adi\_\allowbreak nogo}, at \igaP's entailment
  \textsf{adi\_\allowbreak ent}, by \igaP's $\Mw$-soundness
  \textsf{adi\_\allowbreak omega\_\allowbreak sound}
  (Theorem~\ref{thm:ad-iga-rule}(3));
\item \textsf{oiga\_\allowbreak nogo}, at \oiga's entailment
  \textsf{oiga.\allowbreak oiga\_\allowbreak ent}, by \oiga's
  $\Mw$-soundness \textsf{oiga.\allowbreak oiga\_\allowbreak omega\_\allowbreak sound}.
\end{itemize}
Six corollaries state the three exclusions of
Theorem~\ref{thm:op-nogo} at the two systems, with that theorem's
hypotheses: a well-formed context $\Gamma$ whose members are all true
in $\Mw$ and, for a residual $R$, the arm conditions $\mathit{arms}(R)$,
a well-formed $R(\bquo{P_{\mathrm h}})$ and no free $v_1$ in
$R(\bquo{P_{\mathrm h}})$.
\begin{itemize}
\item The guard half $\mathit{ResGuard}(\bquo{P_{\mathrm h}})$ of the
  fold's residual is not derivable from $\Gamma$
  (\textsf{adi\_\allowbreak not\_\allowbreak iResGuard},
  \textsf{oiga\_\allowbreak not\_\allowbreak iResGuard}).
\item The fold's motive is not derivable from $\Gamma$, in the three
  forms of \textsf{iga\_\allowbreak not\_\allowbreak motive}: two at
  the hypothesis node of $P_{\mathrm h}$ and one at every rank
  (\textsf{adi\_\allowbreak not\_\allowbreak motive},
  \textsf{oiga\_\allowbreak not\_\allowbreak motive}).
\item No residual $R$ with these conditions has $R(\bquo{P_{\mathrm
  h}})$ derivable from $\Gamma$
  (\textsf{adi\_\allowbreak not\_\allowbreak fold\_\allowbreak residual},
  \textsf{oiga\_\allowbreak not\_\allowbreak fold\_\allowbreak residual}).
\end{itemize}
The corollaries are classical, like Theorem~\ref{thm:op-nogo}.  Like
that theorem, the corollaries exclude only the fold with the record
motive, and they leave alone the three things that
Section~\ref{sec:open-nogo} lists.  At \igaP the last of these, a
derivation of $\ConR$, exists (Theorem~\ref{thm:ad-iga-con}).

\subsection{The build-time checks}
\label{app:adopt-checks}

Each of the three theories \texttt{Adopt\_RGA}, \texttt{Adopt\_IGA}
and \texttt{NoGo\_Ext} ends with a dependency check that runs at build
time and fails the build when a condition fails.
\begin{itemize}
\item \texttt{Adopt\_RGA}: no oracle under the seven results
  \textsf{adr\_warrant\_uC}, \textsf{adr\_ConRGA},
  \textsf{adr\_separation}, \textsf{adr\_consistent},
  \textsf{adr\_omega\_sound}, \textsf{adr\_U\_ConRGA} and
  \textsf{adr\_ConRGA\_U}; \textsf{adr\_ConRGA} reaches
  \textsf{adr\_adopt}, \textsf{adr\_exf\_code} and
  \textsf{adr\_warrant\_uC}; and the control \textsf{rga\_goedel2}
  reaches \textsf{iChk\_lstCons}, so the check is not blind, and reaches
  neither \textsf{adr\_adopt} nor \textsf{adr\_exf\_code}.
\item \texttt{Adopt\_IGA}: no oracle under seventeen results;
  \textsf{adi\_ConIGA} reaches \textsf{adi\_adopt},
  \textsf{adi\_exf\_code} and \textsf{adi\_warrant\_uC};
  \textsf{adi\_consistent} reaches \textsf{iChkI\_genuine\_proves},
  \textsf{gi\_step\_oomjt} and \textsf{adi\_adopt\_omv};
  \textsf{adi\_exf\_code} reaches \textsf{G\_iTreeValid\_appC},
  \textsf{G\_iTreeValid\_suc\_intro}, \textsf{G\_ir\_drvChk\_fires\_intro},
  \textsf{tblShift\_IGA} and \textsf{chkI\_iChk\_embed}, and avoids
  \textsf{iChk\_lstCons}, \textsf{iChk\_appC} and
  \textsf{adr\_exf\_code}; and the control \textsf{gi\_consistent}
  reaches \textsf{gi\_step\_oomjt} and no fact whose name begins with
  \textsf{adi\_} or \textsf{g2i\_}.
\item \texttt{NoGo\_Ext}: no oracle under the six corollaries of
  Appendix~\ref{app:adopt-nogo}.
\end{itemize}
The checks also test that the derivations of $\ConR$ and $\ConI$
avoid the soundness lemmas.  \textsf{adr\_ConRGA} reaches none of
\textsf{omega\_sound}, \textsf{rga\_consistent},
\textsf{rga\_omega\_sound}, \textsf{prf\_owsjt}, \textsf{prf\_oomjt},
\textsf{iChk\_genuine\_proves}, \textsf{adr\_omega\_sound},
\textsf{adr\_consistent}, \textsf{adr\_adopt\_omv} and
\textsf{rga\_no\_proof\_0eq1}; the control is that
\textsf{adr\_consistent} does reach \textsf{adr\_omega\_sound} and
\textsf{adr\_adopt\_omv}.  \textsf{adi\_ConIGA} reaches none of the
eleven soundness names of the check at \iga's G\"odel theorems
(Appendix~\ref{app:iga-goedel-checks}) nor \textsf{adi\_adopt\_omv}.

\bibliographystyle{alpha}
\bibliography{logic}

\end{document}